\documentclass[12pt]{article}
\usepackage{a4wide,latexsym,amsfonts,amsmath,amssymb,epsf,bbm}
\usepackage[dvips]{graphics}

\newcounter{defthm}
\newtheorem{defthm}{\whattheorem}[section]
\newcommand\dt[1]  {\noindent\def\whattheorem{#1}\pagebreak[0]\begin{defthm}{}%
                   \samepage{$\!\!${\rm:}\nopagebreak\\[-1.91em]{}}\end{defthm}}
\newcommand\dtl[2] {\noindent\def\whattheorem{#1}\pagebreak[0]\begin{defthm}{}%
                   \samepage{$\!\!${\rm:}\label{#2}\nopagebreak\\[-1.91em]{}}%
                   \end{defthm}}

\newcommand  {\A}         {L}
\renewcommand{\AA}        {{\!A}}
\newcommand  {\Ao}        {{\A{\times}\oneHb}}
\newcommand  {\AO}        {{\A{\times}\oneHp}}
\newcommand  {\alg}       {algebra}
\newcommand  {\Alg}       {\mbox{-\sl Alg}}
\newcommand  {\apppicture}[2] {\put(#2,0) {\begin{picture}(0,0)(0,0)
                          \scalebox{.29}{\includegraphics{app#1.eps}} \end{picture}}}
\newcommand  {\bc}        {boundary condition}
\newcommand  {\be}        {\begin{equation}}
\newcommand  {\bea}       {\begin{equation}\begin{array}l}
\newcommand  {\bearl}     {\begin{array}{l}}
\newcommand  {\bearll}    {\begin{array}{ll}}
\newcommand  {\Beta}      {\Xbeta{}{}}

\newcommand  {\B}         {\varUpsilon}

\newcommand  {\CAAalm}    {\mbox{${\mathcal C}_{\!A|A}
                              ^{\alpha^-\!\mbox{\tiny-Ind}}$}}
\newcommand  {\CAAalp}    {\mbox{${\mathcal C}_{\!A|A}
                              ^{\alpha^+\!\mbox{\tiny-Ind}}$}}
\newcommand  {\CAAalpm}   {\mbox{${\mathcal C}_{\!A|A}
                              ^{\alpha^\pm\!\mbox{\tiny-Ind}}$}}
\newcommand  {\CAAm}      {\mbox{${\mathcal C}_{\!A|A}^-$}}
\newcommand  {\CAAmu}     {\mbox{${\mathcal C}_{\!A|A}^\mu$}}
\newcommand  {\CAAnu}     {\mbox{${\mathcal C}_{\!A|A}^\nu$}}
\newcommand  {\CAAo}      {\mbox{${\mathcal C}_{\!A|A}^{\,0}$}}
\newcommand  {\CAAp}      {\mbox{${\mathcal C}_{\!A|A}^+$}}
\newcommand  {\CAApm}     {\mbox{${\mathcal C}_{\!A|A}^\pm$}}

\newcommand  {\calca}     {\mbox{${\mathcal C}_{\!A}$}}
\newcommand  {\calcaa}    {\mbox{${\mathcal C}_{\!A|A}$}}
\newcommand  {\calcai}    {\mbox{${\mathcal C}_{\!A}^\ind$}}
\newcommand  {\calcal}    {\mbox{$\Ext{\mathcal C}{\!A}$}}
\newcommand  {\calcali}   {\mbox{${\mathcal C}_{\!A}^\lInd$}}

\newcommand  {\cat}       {category}
\newcommand  {\cats}      {categories}
\newcommand  {\cC}        {\mathcal{C}}
\newcommand  {\cCA}       {\cC\Alg}
\newcommand  {\cCF}       {\cC\mbox{-\sl Frob}}
\newcommand  {\cir}       {\,{\circ}\,}
\newcommand  {\cirb}      {\,{\ol\circ}\,}
\newcommand  {\circb}     {\;{\ol\circ}\;}
\newcommand  {\cD}        {\mathcal{D}}
\newcommand  {\cDA}       {\cD\Alg}

\newcommand  {\cDF}       {\cD\mbox{-\sl Frob}}
\newcommand  {\cft}       {conformal field theory}
\newcommand  {\Cft}       {Conformal field theory}
\newcommand  {\cfts}      {conformal field theories}

\newcommand  {\cG}        {\mathcal{G}}
\newcommand  {\cGG}       {{\cG\Ti\ol\cG}}
\newcommand  {\cH}        {\mathcal{H}}
\newcommand  {\cHb}       {\ol{\mathcal{H}}}
\newcommand  {\cI}        {\mathcal{I}}

\newcommand  {\class}     {classification}
\newcommand  {\cloc}      {c^A} 
\newcommand  {\cocon}     {coset construction}

\newcommand  {\con}       {conformal }
\newcommand  {\Con}       {Conformal }
\newcommand  {\convention}{declaration}
\newcommand  {\Convention}{Declaration}

\newcommand  {\complex}   {\mathbbm{C}}
\newcommand  {\corfu}     {correlation function}
\def\Corollary            {Corollary }

\newcommand  {\cQ}        {\mathcal{Q}}
\newcommand  {\csplit}    {centrally split}

\def\Definition           {Definition }
\def\Definitions          {Definitions }
\newcommand  {\df}        {\,{:=}\,}
\newcommand  {\dream}     {dream}
\newcommand  {\dsty}      {\displaystyle }
\newcommand  {\ee}        {\end{equation}}
\newcommand  {\eear}      {\end{array}}
\newcommand  {\efu}[3]    {E_{#3}^{#2}(#1)}
\newcommand  {\Efu}[2]    {E_{#2}^{}(#1)}
\newcommand  {\EFU}[2]    {E_{#2}^{#1}}
\newcommand  {\End}       {{\rm End}}
\newcommand  {\eps}       {\varepsilon}
\newcommand  {\epicture}[2] {\end{picture}\\{}\\[#1.#2em]\end{array}}
\newcommand  {\eq}        {\,{=}\,}
\newcommand  {\erf}[1]    {(\ref{#1})}

\newcommand  {\Ext}[2]    {#1_{#2}^{\sss\rmloc}}
\newcommand  {\EXt}[2]    {(#1)_{#2}^{\sss\rmloc}}
\newcommand  {\EXT}[2]    {\Llb #1 \Lrb_{#2}^{\sss\rmloc}}
\newcommand  {\F}[9]      {{\sf F}_{\!{\sss#6}#4{\sss#7},{\sss#8}#5{\sss#9}}
                          ^{\,({#1}\,{#2})\,{#3}}}
\newcommand  {\FF}        {{\sf F}}
\newcommand  {\Fi}        {\Phi_\AA}
\newcommand  {\findim}    {fi\-ni\-te-di\-men\-si\-o\-nal}
\newcommand  {\Fol}[9]    {{\ol{\sf F}}_{\!{\sss#6}#4{\sss#7},{\sss#8}#5{\sss#9}}
                            ^{\,({#1}\,{#2})\,{#3}}}
\newcommand  {\foodnode}[1] {\,\footnote{~#1}}
\newcommand  {\Frac}[2]   {\mbox{\large$\frac{#1}{#2}$}}
\newcommand  {\G}[9]      {{\sf G}_{\,{\sss#6}#4{\sss#7},{\sss#8}#5{\sss#9}}
                            ^{\,({#1}\,{#2})\,{#3}}}
\newcommand  {\Gama}      {\varGamma}

\newcommand  {\GHb}       {{\cG\Ti\ol\cH}}
\newcommand  {\GHbc}      {{\cG\Tic\ol\cH}}
\newcommand  {\GHp}       {{\cG\Ti\Hp}}

\newcommand  {\Gmat}      {{\sf G}}
\newcommand  {\Gol}[9]    {{\ol{\sf G}}_{\,{\sss#6}#4{\sss#7},{\sss#8}#5{\sss#9}}
                            ^{\,({#1}\,{#2})\,{#3}}}
\newcommand  {\haploid}   {-ha\-plo\-id}
\newcommand  {\hatota}    {{\hat\otimes}{}_{\!A}}

\newcommand  {\Hb}        {{\ol{\mathcal H}}}

\newcommand  {\HH}        {{\cH\Ti\Hp}}

\newcommand  {\Hom}       {{\rm Hom}}
\newcommand  {\HomA}      {{\rm Hom}_{\!A}}
\newcommand  {\HomAA}     {{\rm Hom}_{\!A|A}}
\newcommand  {\Hp}        {{\mathcal{H}'}}

\newcommand  {\hy}        {$\mbox{-\hspace{-.66 mm}-}$\linebreak[0]}
\newcommand  {\id}        {{\mbox{\sl id}}}
\newcommand  {\iD}        {{\mbox{\scriptsize\sl id}}}
\newcommand  {\ID}        {{\mbox{\tiny\sl id}}}
\newcommand  {\Id}        {{\mbox{\sl Id}}}
\newcommand  {\IG}        {{\mathcal I}_{\mathcal G}}
\newcommand  {\IH}        {{\mathcal I}_{\mathcal H}}

\newcommand  {\II}        {{\mathcal I}}
\renewcommand{\Im}        {{\rm Im}\,}
\newcommand  {\iN}        {\,{\in}\,}
\newcommand  {\In}        {\prec}

\newcommand  {\ind}       {{\mbox{\tiny Ind}}}
\newcommand  {\Ind}       {{\rm Ind}}
\newcommand  {\Indap}[3]  {{\rm Ind}_{#1}^{#2}(#3)}
\newcommand  {\Indk}      {{\rm Ind}^{}_{(A;\iD_\AA)}}
\newcommand  {\IndK}      {{\rm Ind}^{}_{(A;\ID_\AA)}}
\newcommand  {\J}[1]      {\Gamma_{\!#1}}
\newcommand  {\kar}[1]    {{#1}^{\rm K}}
\newcommand  {\kaR}[1]    {{#1}^{\!\rm K}}

\newcommand  {\koerper}   {\Bbbk} %% Bbbk needs amssymb
\newcommand  {\kx}        {\koerper^{\!\times}}
\newcommand  {\labl}[1] {\label{#1}\ee}
\def\Lemma                {Lemma }
\def\Lemmata              {Lemmata }
\newcommand  {\lhs}       {left hand side}

\newcommand  {\lInd}      {{\mbox{\tiny$\ell$-Ind}}}
\newcommand  {\llb}       {\mbox{\large[}}
\newcommand  {\Llb}       {\mbox{\large(}}

\newcommand  {\lrb}       {\mbox{\large]}}
\newcommand  {\Lrb}       {\mbox{\large)}}
\newcommand  {\lxt}[2]    {\mbox{$\ell$-Ind}_{#2}^{}(#1)}
\newcommand  {\Lxt}[2]    {\mbox{\scriptsize$\ell$-Ind}_{#2}^{}(#1)}
\newcommand  {\LXT}[1]    {\mbox{$\ell$-Ind}_{#1}^{}}
\newcommand  {\lxtp}[3]   {\mbox{$\ell$-Ind}_{#2}^{#3}(#1)}
\newcommand  {\Lxtp}[3]   {\mbox{\scriptsize$\ell$-Ind}_{#2}^{#3}(#1)}
\newcommand  {\LXTp}[2]   {\mbox{$\ell$-Ind}_{#1}^{#2}}
\newcommand  {\M}         {{\dot M}}
\newcommand  {\mcll}      {\multicolumn2{|l|}}
\newcommand  {\mclll}     {\multicolumn3{|l|}}
\newcommand  {\mclo}      {\multicolumn1{l|}}
\newcommand  {\mclO}      {\multicolumn1l}
\newcommand  {\modinv}    {modular invarian}
\newcommand  {\mtc}       {modular tensor category}

\newcommand  {\N}[3]      {{N_{#1#2}}^{\!\!#3}}

\newcommand  {\Obj}       { {\rm Obj} }
\newcommand  {\Objc}      {{\rm Obj}(\cC)}
\newcommand  {\ol}        {\overline}
\newcommand  {\one}       {{\bf 1}}
\newcommand  {\One}       {{0}}
\newcommand  {\oneb}      {{\ol{\one}}}
\newcommand  {\oneH}      {{\one_{\!\cH}}}
\newcommand  {\oneHb}     {{\one_{\!\ol\cH}}}
\newcommand  {\oneHp}     {{\one_{\!\Hp}}}
\newcommand  {\oneQ}      {{\one_{\!\cQ}}}
\newcommand  {\operation} {operation}
\newcommand  {\OT}        {{\oneQ{\times}T}}
\newcommand  {\ota}       {\otimes_{\!A}}
\newcommand  {\OtA}       {{\otimes}_{\!A}^{}}
\newcommand  {\otA}       {\,{\otimes}_{\!A}^{}\,}
\newcommand  {\oti}       {\,{\otimes}\,}
\newcommand  {\Oti}       {{\otimes}}
\newcommand  {\otib}      {\,{\ol\otimes}\,}
\newcommand  {\otic}      {\,\Otic\,}
\newcommand  {\Otic}      {{\otimes^{}_{\sss\!\koerper}}}

\newcommand  {\platl}     {trivialisable}
\newcommand  {\platy}     {trivialisability}

\newcommand  {\platz}     {trivialisation}
\def\Proposition          {Proposition }
\def\Propositions         {Propositions }
\newcommand  {\q}         {quantum }
\newcommand  {\Q}         {Quantum }
\newcommand  {\qed}       {\hfill$\Box$\medskip}
\newcommand  {\QH}        {{\cQ\Ti\cH}}
\newcommand  {\QHc}       {{\cQ\Tic\cH}}
\newcommand  {\QHh}       {{\cQ\Ti\cH\Ti\ol\cH}}

\newcommand  {\QHH}       {{\cQ\Ti\cH\Ti\Hp}}

\renewcommand{\r}         {{\rho}}
\newcommand  {\rholoc}[2] {\rho_{#1;#2}^{\sss\rmloc}}
\newcommand  {\Rholoc}    {\rho_{}^{\sss\rmloc}}
\newcommand  {\rmloc}     {{\ell{\rm oc}}}
\newcommand  {\rr}        {\rho_{\rm r}}
\newcommand  {\R}[5]      {{\sf R}^{(#1\,#2)#3}_{#4\,#5}}
\def\Remark               {Remark }
\def\Remarks              {Remarks }
\newcommand  {\rep}       {representation}

\newcommand  {\REP}       {{\mathcal Rep}}

\newcommand  {\retmodule} {module retract}  
\newcommand  {\rhs}       {right hand side}
\newcommand  {\Rol}[5]    {{\ol{\sf R}}^{(#1\,#2)#3}_{#4\,#5}}
\newcommand  {\Rm}[5]     {{\sf R}^{-\,(#1\,#2)#3}_{\;#4\,#5}}
\newcommand  {\Rmol}[5]   {{\ol{\sf R}}^{-\,(#1\,#2)#3}_{\;#4\,#5}}
\newcommand  {\RR}        {{\sf R}}

\newcommand  {\sect}[1]   {\section{#1}\setcounter{equation}0\setcounter{defthm}0}
\def\Section              {Section }
\def\Sections             {Sections }

\newcommand  {\sse}       {\scriptsize}
\newcommand  {\ssF}       {symmetric special Frobenius}
\newcommand  {\ssFA}      {symmetric special Frobenius algebra}

\newcommand  {\sss}       {\scriptscriptstyle}

\newcommand  {\TA}        {{F}}
\newcommand  {\tc}        {tensor category}
\newcommand  {\tcs}       {tensor categories}

\def\Theorem              {Theorem }
\def\Theorems             {Theorems }
\newcommand  {\ti}        {\,{\times}\,}
\newcommand  {\Ti}        {{\boxtimes}}
\newcommand  {\Tic}       {{\otimes^{}_\koerper}}
\newcommand  {\tildeAO}   {{\Lxt\AO\OT}}
\newcommand  {\tildeb}[2] {\tilde\gamma^{A\,#1#2}}
\newcommand  {\tildeB}    {{\Lxt B\AA}} 
\newcommand  {\TildeB}    {{\lxt B\AA}} 

\newcommand  {\tildeOT}   {{\Lxt\OT{\AO\!}}}
\newcommand  {\tr}        {{\rm tr}}
\newcommand  {\U}         {\,}

\newcommand  {\Vectc}     {\mathcal Vect_\complex}
\newcommand  {\Vectk}     {\mathcal Vect_\koerper}

\newcommand  {\Xbeta}[2]  {\tilde\Gamma^{A\,#1#2}}
\newcommand  {\Xb}[2]     {\gamma^{A\,#1#2}}
\newcommand  {\Yb}[3]     {\gamma^{C_{#3}\,#1#2}}
\newcommand  {\Xgamma}[2] {\Gamma^{A\,#1#2}}
\newcommand  {\Ygamma}[3] {\Gamma^{C_{#3}\,#1#2}}

\newcommand  {\zet}       {\mathbb{Z}}

\begin{document} 

%%%%%%%%%%%%%%%%%%%%%%%%%%%%%%%%%%%%%%%%%%%

 \begin{flushright}  {~} \\[-1cm] {\sf math.CT/0309465} \\[1mm]
 {\sf HU-EP-03/31} \\
 {\sf Hamburger Beitr\"age zur Mathematik Nr.\ 179} \\[1 mm]
 \end{flushright}

 \begin{center} \vskip 11mm
 {\Large\bf CORRESPONDENCES}\\[1.02em]
 {\Large\bf OF RIBBON CATEGORIES}
 \\[12mm]
 {\large J\"urg Fr\"ohlich}\,$^1_{}$ \ \ {\large J\"urgen Fuchs}\,$^2_{}$ \ \
 {\large Ingo Runkel}\,$^3_{}$ \ \ {\large Christoph Schweigert}\,$^4_{}$
     \\[9mm]
 $^1\;$ Institut f\"ur Theoretische Physik, \ ETH Z\"urich\\ CH\,--\,8093\,
    Z\"urich\\[2mm]
 $^2\;$ Institutionen f\"or fysik, \ Karlstads Universitet\\
    Universitetsgatan 5, \ S\,--\,651\,88\, Karlstad\\[2mm]
 $^3\;$ Institut f\"ur Physik, \ Humboldt Universit\"at Berlin\\
 Newtonstra\ss{}e 15, \ D\,--\,12\,489\, Berlin\\[2mm]
 $^4\;$ Fachbereich Mathematik, \ Universit\"at Hamburg \\
 Schwerpunkt Algebra und Zahlentheorie\\
 Bundesstra\ss e 55, \ D\,--\,20\,146\, Hamburg

  \end{center} \vskip 12mm
 \begin{quote}{\bf Abstract}\\[1mm] 
Much of algebra and representation theory can be formulated in the general 
framework of tensor categories.
The aim of this paper is to further develop this theory for braided tensor
categories. Several results are established that do not have a substantial
counterpart for symmetric tensor categories. In particular, we exhibit
various equivalences involving categories of modules over algebras in
ribbon categories. Finally we establish a correspondence of ribbon categories 
that can be applied to, and is in fact motivated by, the coset construction 
in conformal quantum field theory.  
 \end{quote} \newpage

%%%%%%%%%%%%%%%%%%%%%%%%%%%%%%%%%%%%%%%%%%%

{\small\tableofcontents}\newpage

%%%%%%%%%%%%%%%%%%%%%%%%%%%%%%%%%%%%%%%%%%%

\sect{Introduction}

In this paper we study equivalences involving categories of modules over
algebras in ribbon categories. Our main results are 
\Theorem \ref{thm:equiv} and \Theorem \ref{thm:coset}. 
To motivate these results and clarify their relevance, we start by looking 
at a classical analogue: correspondences.

\subsection{Correspondences}

Correspondences are often needed to express relations between mathematical 
objects of the same type. For instance, in algebraic geometry they enter in 
the definition of rational maps. A more recent application is the construction 
of an action of the Heisenberg algebra on the cohomology of Hilbert schemes 
of points on surfaces. In the present paper, we introduce a generalisation of 
correspondences in the setting of braided tensor categories, which turns out 
to provide a powerful tool for the study of such categories.

Correspondences deal with classes of mathematical objects for which a Cartesian
product is defined. For definiteness, let us consider finite groups.
A correspondence of two groups $G_1$ and $G_2$ is a subgroup $R$ of the 
product group $G_1\ti G_2$,
  \be  R \,\leq\,  G_1 \ti G_2 \, . \labl{corr1}
Suppose now that the representation theories of the groups $G_1$ and $R$ 
are known. One could then be tempted to formulate the following \dream:
A correspondence \erf{corr1} might allow us to express   
the category $\REP(G_2)$ of (finite-dimensional complex) representations of 
$G_2$ in terms of the representation categories $\REP(G_1)$ and $\REP(R)$.

Obviously, in this generality our \dream\ is entirely unrealistic -- just take 
$G_1$ and $R$ to be trivial. To assess the feasibility of the \dream\ in more 
general categories than representation categories of finite groups, it is 
helpful to reformulate the correspondence \erf{corr1} in the spirit of the 
Tannaka-Krein philosophy, i.e.\ to express statements about groups entirely in 
terms of their representation categories rather than in terms of the groups 
themselves. One advantage of this point of view is the following. Once the 
statements are translated to a category-theoretic setup, one can try to relax 
some of the properties of the representation category so as to arrive at 
analogous statements applying to categories that appear in other contexts, 
e.g.\ as representation categories of quantum groups, of vertex algebras, or 
of precosheaves of von Neumann algebras, and that, in turn, have important 
applications in quantum field theory.

Our starting point, i.e.\ the correspondence \erf{corr1}, is easily
reformulated in category-theoretic language. The representation category of
the product group is simply the product of the two representation categories,
$\REP(G_1{\times}G_2) \,{\cong}\, \REP(G_1)\,\Ti\,\REP(G_2)$.\,%
  \footnote{~For a precise definition of the relevant notion of product
tensor category, see \Section \ref{platy1}.}
The correspondence $R$ is, by definition, a subgroup of $G_1\ti G_2$; a 
category-theoretic analogue of the notion of subgroup is known (\cite{kios}; 
  %  theorem 2.2 
for earlier discussions compare
also \cite{w-unpublished,o-nowpublished}): There is a bijection
between subgroups $H$ of a group $G$ and commutative algebras in the
tensor category $\REP(G)$. The commutative algebra in $\REP(G)$ that is
associated to $H$ is given by the space $\complex(G/H)$ of functions on 
the homogeneous space $G/H$. The category $\REP(G)_{\complex(G/H)}$ of 
$\complex(G/H)$-modules in $\REP(G)$ is equivalent to $\REP(H)$,
  \be  \REP(G)_{\complex(G/H)} \cong \REP(H) \, . \ee
Our \dream\ can thus be stated more precisely as follows. Suppose
we are given two tensor categories $\cC_1$ and $\cC_2$ and a commutative 
semisimple algebra $A_R$
in $\cC_1\,\Ti\,\cC_2$. Denote by $\cC_R$ an appropriate tensor category of
$A_R$-modules. Then we might attempt to express $\cC_2$ in terms of $\cC_1$ 
and $\cC_R$, as the category of modules over a commutative algebra $B$
in a tensor category $\cC$ that is derived from $\cC_1$ and $\cC_R$ only.

In the particular case that $G_1$ and $R$ are trivial, our \dream\ would amount 
to the statement that $\REP(G_2)$ is equivalent to the representation category 
of a commutative semisimple algebra over $\complex$, which clearly cannot be 
true for any non-abelian group $G_2$. More explicitly, in this specific
situation the data involved in the correspondence are, in category-theoretic 
language, the tensor category $\REP(G_2)$ and the commutative algebra
$\complex(G_2)$ of functions on $G_2$, seen as an algebra in $\REP(G_2)$.  
Since all irreducible representations of $G_2$ appear as subrepresentations
of $\complex(G_2)$, this algebra has trivial representation theory:
  \be  \REP(G_2)_{\complex(G_2)} \cong \Vectc \,. \labl{l2}

It is therefore all the more remarkable that there do exist situations in
which our \dream\ {\em can\/} be realised. It involves a generalisation of 
algebra and representation theory to tensor categories that are not 
necessarily symmetric, but are still braided. Among such categories there are, 
in particular, the modular tensor categories. The interest in modular tensor 
categories comes e.g.\ from the fact that such a category contains the data 
needed for the construction of a three-dimensional topological quantum field 
theory. These categories arise in many interesting applications; for example, 
the representation categories of certain vertex algebras are modular \tcs. 

Modular tensor categories are distinguished by a 
non-degeneracy property of the braiding; in particular, the braiding is 
``maximally non-symmetric''. This makes it apprehensible that in contrast to 
the classical case above, in which all involved tensor categories are 
symmetric, such categories can indeed provide a realisation of our \dream.

%%%%%%%%%%%%%%%%%%%%%%%%%%%%%%%%%%%%%%%%%%%

\subsection{Frobenius algebras}

Many aspects of the representation theory of rings or algebras can be 
generalised to the general setting of tensor categories \cite{pare11}. In any 
tensor category one has the notions of an associative algebra with unit and 
its modules and bimodules. Similarly one can define coalgebras. A particularly 
interesting class are algebras $A$ that are also coalgebras such that the 
coproduct is a bimodule morphism from $A$ to the $A$-bimodule $A\oti A$. Such 
algebras are called {\em Frobenius algebras\/}, because Frobenius algebras in 
the modular tensor category of finite-dimensional vector spaces over some 
field are just ordinary Frobenius algebras. Frobenius algebras in 
more general tensor categories have recently attracted attention in several 
different contexts (see e.g.\ \cite{kios,fuSc16,fuRs,ostr,muge8,fuRs4}).

In contrast to bialgebras (such as Hopf algebras), Frobenius algebras can be
defined in tensor categories that are not necessarily braided.
In a braided category, however, their theory becomes much richer. One then
has the notion of a commutative algebra and, more generally, of
center(s) of an algebra. The present paper aims at developing the theory of
Frobenius algebras in such a setting. It turns out to be helpful to impose
a few additional requirements, both on the algebra and on the category.
In particular, we assume that the braided tensor category in question is 
additive, $\koerper$-linear (with $\koerper$ some field), as well as sovereign 
-- it has a left and a right duality that coincide as functors from $\cC$ to 
the opposed category; a braided sovereign tensor category is also known as a 
{\em ribbon\/} category. Other requirements imposed on the category will be 
given in the body of the paper; the setting is summarised in \convention\ 
\ref{basicprops}.

The additional properties of the algebra are that it is a {\em special\/} and
{\em symmetric\/} Frobenius algebra, see definition \ref{symm-frob-spec}. (To 
ensure the existence of various images needed in our constructions, we also 
assume that the algebra is what we call \csplit, see definition \ref{csplit} 
and \convention\ \ref{c:csplit}.) Symmetric Frobenius algebras in the category 
of vector spaces appear e.g.\ in the study of group algebras and thus play 
a central role in representation theory. It is worth noting that in a 
braided setting, a commutative Frobenius algebra is not necessarily 
symmetric. The specialness property of the Frobenius algebra $A$ implies, 
in particular, \cite{kios,fuSc16} that when the category $\cC$ is semisimple 
then the category of left $A$-modules is semisimple as well.

%%%%%%%%%%%%%%%%%%%%%%%%%%%%%%%%%%%%%%%%%%%

\subsection{Local modules and local induction}\label{intro-locind}

In this paper we study symmetric special Frobenius algebras $A$ in ribbon
categories $\cC$. Given such an algebra, there are three other categories
one should consider: The category \calca\ of left $A$-modules, the
analogous category of right $A$-modules, and the
category \calcaa\ of $A$-bimodules. The tensor product $B_1\otA B_2$
of bimodules endows $\calcaa$ with the structure of a tensor category.

The braiding of $\cC$ allows to construct two tensor functors \cite{lore,ostr}
  \be  \alpha^\pm_\AA:\quad \cC \to \calcaa \,,  \labl{Ialpha}
known as $\alpha$-induction (see \Definition \ref{alphadef}).
In \Definition \ref{def:[]-functor} we introduce two endofunctors
  \be  \EFU{l/r}\AA:\quad \cC \to \cC \,. \labl{EFU}
We show in \Proposition \ref{alphaalpha} 
that if right-adjoint functors
$(\alpha^{\pm})^\dagger$ to \erf{Ialpha} exist, then the endofunctors 
\erf{EFU} are the compositions
$\EFU{l}\AA\eq(\alpha^+_\AA)^\dagger \cir \alpha^-_\AA$ and
$\EFU{r}\AA\eq(\alpha^-_\AA)^\dagger \cir \alpha^+_\AA$. For commutative 
algebras, the two functors $\EFU{l/r}\AA$ coincide (see 
\Proposition \ref{prop:AB-alg} (iv)); 
in this case we suppress the index $l$ or $r$.

A basic principle in this paper is to try to lift a given functor
$F{:}\ \cC\,{\to}\,\cD$ between two tensor categories $\cC$ and $\cD$ to a
functor from the category $\cCA$ of algebras in $\cC$ to the category $\cDA$
of algebras in $\cD$, or even to a functor between the respective categories
$\cCF$ and $\cDF$ of Frobenius algebras. For the functors
$\EFU{l/r}\AA$ both lifts turns out to be possible. This result, established in 
\Proposition \ref{prop:AB-alg}(i), is non-trivial because $\EFU{l/r}\AA$
are not necessarily tensor functors. By abuse of notation, we use the same
symbol for the resulting endofunctors of $\cCA$ and of $\cCF$
as for the underlying endofunctors of $\cC$, i.e.\ we write
  \be  \EFU{l/r}\AA:\quad \cCA \to \cCA  \labl{EFUA}
as well as $\EFU{l/r}\AA{:}\ \cCF\;{\to}\,\cCF$.

\smallskip

The images of the endofunctors \erf{EFU} carry additional structure. To 
describe it we need two additional ingredients: a braided version of the 
concept of the center of an algebra and the concept of local modules. First, 
the braiding allows one to generalise the notion of a center of an algebra 
$A$, and for a general braiding one obtains in fact two different centers 
$C_l(A)$ and $C_r(A)$, known as the {\em left\/} and the {\em right center\/} 
of $A$, respectively. After adapting, in \Definition \ref{defLRC}, 
their description in \cite{vazh,ostr} to the present setting, we show in 
\Proposition \ref{lem:C=[1]i} 
that the centers of a symmetric special Frobenius algebra
carry the structure of commutative symmetric Frobenius algebras.
In a braided category there is also a notion of the tensor product $A\oti B$ 
of two algebras $A$ and $B$. It enters e.g.\ in the definition \cite{vazh} 
of the Brauer group of the category. Remarkably, in the braided setting the 
tensor product of two commutative algebras is not necessarily commutative. 
(Thus it is not natural to restrict one's attention to commutative algebras.) 
In \Proposition \ref{prop:tensor-center}(i) 
we compute the centers of $A\oti B$; 
they can be expressed in terms of the endofunctors \erf{EFUA}, namely
  \be  C_l(A\Oti B) \,\cong\, \efu{C_l(B)}l\AA  \qquad{\rm and}\qquad
  C_r(A\Oti B) \,\cong\, \efu{C_r(A)}r{\!B}  \labl{I0}
as Frobenius algebras.

The category of left modules over a {\em commutative\/} algebra $A$ in
$\Vectc$ is again a tensor category. In order to generalise this fact
to a braided setting, a refinement is necessary, and this refinement makes
use of the second ingredient of our construction -- the concept of
{\em dyslectic\/} \cite{pare23} or {\em local\/} module. A module $M$ over 
a commutative special Frobenius algebra $A$ in a ribbon category is local 
iff the representation morphism commutes with the twist (see 
\Proposition \ref{pro:ostr}), so that the twist on $M$ is a morphism in 
\calca. The resulting generalisation of the classical statement is given in 
\Proposition \ref{thm:mod}, 
which follows \cite{pare23} and \cite{kios}: The category 
\calca\ of left modules over a commutative \ssFA\ $A$ in a ribbon category 
$\cC$ has a natural full subcategory -- the category \calcal\ of local left
$A$-modules -- that is again a tensor category, and in fact, unlike e.g.\ 
the category of $A$-bimodules, even a ribbon category.

\medskip

With these results at hand, we proceed to show, in 
\Proposition \ref{pr:EAB-locmod}, 
that every object in the image of the endofunctors
$\EFU{l/r}\AA$ has a natural structure of a local $C_l(A)$-module,
respectively of a local $C_r(A)$-module. Thus the functors $\EFU{l/r}\AA$
give rise to two functors
  \be  \LXTp\AA{l/r}: \quad \cC\to \Ext\cC{C_{l/r}(A)} \,,  \labl{LXTp}
which we call {\em local induction\/} functors (\Definition \ref{def-lxt}).
(Again, for commutative algebras, the two functors coincide, and we shall
then suppress the index $l$ or $r$, i.e.\ just write $\LXT\AA$.)
However, in contrast to ordinary induction, local induction is not a tensor
functor. In the tensor categories $\Ext\cC{C_{l/r}(A)}$ we have the notion 
of an algebra; it turns out (\Proposition \ref{lem:[B]A-lift-i}) 
that the local induction functors can be extended to
functors between categories of algebras, too, i.e.\ (again abusing notation)
  \be  \LXTp\AA{l/r}: \quad \cCA \,\to\, \Ext\cC{C_{l/r}(A)}\Alg \,.  \ee

\medskip

All this structure enters the following result about successive local 
inductions. Let $A$ and $B$ be two commutative symmetric special Frobenius 
algebras in $\cC$. Then \calcal\ is again a tensor category, and $\lxt B\AA$ 
is a commutative algebra in that category. It thus makes sense to consider 
the tensor category of local $\lxt B\AA$-modules in \calcal. In 
\Proposition \ref{thm:[B]A-lift-ii} 
we show that this category can also be obtained as the 
category of local modules over some commutative algebra in $\cC$, and that 
this algebra is in fact just $\Efu B\AA$:
  \be  \EXt{\Ext{\cC}\AA}{\Lxt B\AA} \, \cong \, \Ext{\cC}{\Efu B\AA}
  \,.  \labl{I2}
If in addition $A$ is simple and $\Efu B\AA$ is special, then this is even an 
equivalence of ribbon categories. (An algebra is called simple iff it is
  % (absolutely)
simple as a bimodule over itself, see \Definition \ref{defsimple}.)

The next statement -- \Theorem \ref{thm:equiv} 
-- is the first main result
of this paper: Provided that the left and right centers $C_l(A)$ and $C_r(A)$ 
of a symmetric special Frobenius algebra $A$ in a ribbon category $\cC$ 
(which are symmetric Frobenius by \Proposition \ref{lem:C=[1]i}) are also 
special, the categories of local modules over $C_l(A)$ and $C_r(A)$ are 
equivalent as ribbon categories,
  \be  \Ext\cC{C_l(A)} \cong\, \Ext\cC{C_r(A)} \,. \labl{I1}
Moreover, there is in addition a ribbon equivalence of these categories
to a certain subcategory of $\alpha$-induced $A$-bimodules, the
category \CAAo\ of {\em ambichiral\/} $A$-bimodules, introduced in 
\Definition \ref{def:CAApm}.

The equivalence \erf{I1} can, in general, not be extended to an equivalence of 
the respective categories of all modules (as module categories over $\cC$) --
the left center and the right center are not necessarily Morita equivalent.

\medskip

It is instructive to see how the results \erf{I0}, \erf{I2} and \erf{I1} 
simplify for a symmetric tensor category $\cC$, in which the braiding obeys 
$c_{U,V}^{-1}\eq c_{V,U}^{}$, for all objects $U,\,V$. This includes in 
particular the `classical' situation that $\cC$ is the category $\Vectk$ of 
finite-dimensional vector spaces over a field $\koerper$, as well as the 
category of finite-dimensional super vector spaces. In this case, the notions
of left and right center coincide, there is only a single center $C(A)$.
The relations \erf{I0} then reduce to the statement that the center of the
tensor product of two algebras is 
the tensor product of the centers, $C(A\Oti B)\,{\cong}\, C(A)\oti C(B)$.

Furthermore, in a symmetric tensor category {\em all\/} modules over a 
commutative special Frobenius algebra are local. The result \erf{I2} thus 
simplifies to a simple statement about the induction with respect to the 
tensor product of two commutative algebras $A$ and $B$:
$(\calca)_{{\rm Ind}_\AA(B)}\,{\cong}\;\cC_{\AA\otimes B}.$

Finally, there is only a single $\alpha$-induction $\alpha_\AA^{}\eq
\alpha_\AA^+\eq\alpha_\AA^-$, and the two endofunctors $\EFU{l/r}\AA$ of $\cC$
just amount to tensoring objects with $C(A)$ and morphisms with $\id_{C(A)}$. 
The two functors $\LXTp\AA{l/r}$ coincide as well, and are induction
to modules over $C(A)$. Therefore, in a symmetric tensor category, our first
main result \erf{I1} becomes a tautology -- in other words, \erf{I1} is a
theorem of `braided algebra' without substantial classical analogue.

%%%%%%%%%%%%%%%%%%%%%%%%%%%%%%%%%%%%%%%%%%%

\subsection{Correspondences and the trivialisation of ribbon categories}

Before we can discuss the category-theoretic generalisation of correspondences,
we must still find an appropriate generalisation to the braided setting
of the relation \erf{l2}, i.e.\ of the fact that
the category $\REP(G)$ of representations of a group $G$
contains a commutative special symmetric Frobenius algebra $A\eq\complex(G)$
such that $\Ext{\REP(G)}\AA\,{\cong}\,\Vectc$. We call an algebra $A$ in $\cC$
with the property that $\calcal\,{\cong}\,\Vectk$ a {\em trivialising algebra\/}
for $\cC$.

Requiring the existence of a trivialising algebra is too restrictive for the 
applications we have in mind. We rather need the following more general concept 
(\Definition \ref{def:platy}): We call a
($\koerper$-linear) ribbon category $\cC$ {\em \platl\/} iff there exist
a ribbon category $\cC'$ and a commutative symmetric special Frobenius algebra 
$T$ in $\cC\,\Ti\,\cC'$ such that the category of local $T$-modules is trivial,
  \be  \Ext{(\cC\Ti\cC')}T \,\cong\, \Vectk \,.  \ee
An important class of braided tensor categories are the modular tensor
categories, which play a key role in various applications. In 
\Proposition \ref{thm:top} 
we show that every modular tensor category is \platl,
with $\cC'\eq\ol\cC$ the \tc\ dual to $\cC$.

\medskip

Combining all these results finally allows us to obtain a category-theoretic
generalisation of the correspondence \erf{corr1}. Suppose that a ribbon
category $\cC_3$ is equivalent to the category of local $A$-modules in the
product of two ribbon categories $\cC_1$ and $\cC_2$, i.e.\ that the 
correspondence takes the form
  \be  \cC_3 \cong \EXt{\cC_1\Ti\cC_2}\AA \,,  \ee
where $\cC_2$ is \platl\ with trivialising algebra $T$ in $\cC_2\Ti\cC_2'$.
The \dream\ spelt out in the beginning then amounts to expressing
$\cC_1$ as the category of local modules over a commutative special
Frobenius algebra in $\cC_3\,\Ti\,\cC_2'$. We shall indeed show 
(\Proposition \ref{prop:coset}) that, quite generally,
it is possible to express a category of local modules over a certain
commutative algebra in $\cC_1$ in terms of $\cC_3$ and $\cC_2'$:
  \be  \EXt{\cC_1}{\Lxt{A\otimes\one}{\one\otimes T}}
  \,\cong\, \EXt{\cC_3\Ti\cC_2'}{\Lxt{\one\otimes T}{A\otimes\one}} \,.
  \labl{I9}
($A\Oti\one$ and $\one\Oti T$ are algebras in $\cC_1\Ti\cC_2\Ti\cC_2'$, and
$\one$ denotes the tensor unit of the respective \cat; the product 
$\Ti$ of tensor categories is associative, see \Remark \ref{boxbox}.) 

Moreover, the situation simplifies considerably when we make the following
restrictions. First, we demand that the category $\cC_2$ is modular;
second, we require that the algebra $A$ in $\cC_1\,\Ti\,\cC_2$
has the property that the only subobject of $A$ of the form $U\ti\one$ is
$\one\ti\one$. Then the commutative algebra $\lxt{A\Oti\one}{\one\otimes T}$
in $\cC_1$ is the tensor unit, so that \erf{I9} reduces to
  \be  \cC_1 \,\cong\, \EXt{\cC_3\Ti\ol\cC_2  }{B}  \labl{I10}
with $B\eq\lxt{\one\otimes T}{A\otimes\one}$.
This result -- \Theorem \ref{thm:coset} -- is arguably the strongest
possible realisation of our \dream. We stress that only in a braided
setting such an effect can happen: It is the non-triviality of the braiding
that is responsible for getting the {\em locally\/} induced algebra
$\lxt{A\otimes\one}{\one\otimes T}$ so small.

%%%%%%%%%%%%%%%%%%%%%%%%%%%%%%%%%%%%%%%%%%%

\subsection{Applications in quantum field theory}

$\complex$-linear tensor categories have played a prominent role in quantum 
field theory, especially in connection with the general analysis of 
superselection rules and of quantum statistics \cite{dohr,doro2}. In 
two- and three-dimensional quantum field theory they have become an
indispensable tool for studying braid statistics and quantum symmetries. 
The analysis presented in this paper is primarily inspired by problems
in two-dimensional conformal field theory and string theory and
has grown out of the results presented in \cite{fuRs,fuRs4}. Concrete
applications of our results, in particular of \erf{I10}, to conformal
field theory form the subject of a forthcoming paper. Here we just
give an indication of what some of these applications consist in.

First consider the case that $\cC_3$ is trivial, $\cC_3\,{\cong}\,\Vectc$.
This case can e.g.\ be realised through certain conformal embeddings
of directs sums $\hat{\mathfrak g}_1\,{\oplus}\,\hat{\mathfrak g}_2$ of
untwisted affine Lie algebras into an untwisted affine Lie algebra
$\hat{\mathfrak g}$. The relevant tensor categories are the \cats\
$\cC_i\eq\cC(\mathfrak g_i,k_i)$ of integrable \rep s of
the affine Lie algebras $\hat{\mathfrak g}_1$ and $\hat{\mathfrak g}_2$ with 
specified values $k_{1,2}$ of the level; as \rep s for $\hat{\mathfrak g}$ 
one must take the integrable \rep s at some level $k$, and require that the 
\cat\ of those \rep s is equivalent to $\Vectc$, which
is the case for ${\mathfrak g}\eq E_8$ at level $k\eq1$.
These are modular tensor categories, and the embedding of $\hat{\mathfrak g}_1
\,{\oplus}\,\hat{\mathfrak g}_2$ into $\hat{\mathfrak g}$ provides us with a 
simple commutative \ssFA\ $A$ in their product $\cC_1\Ti\cC_2$. Our result 
\erf{I9} then asserts that a category of local modules in $\cC_1$ is equivalent
to a category of local modules in the category $\ol\cC_2$ dual to $\cC_2$. If,
in addition, the conditions are met that the only subobject of $A$ of the form
$U{\times}\one$ is $\one{\times}\one$ and the only subobject of $A$ of the form
$\one{\times}U$ is $\one{\times}\one$, then the \cats\ $\cC_1$ and $\ol\cC_2$
are equivalent; this happens e.g.\ for those conformal embeddings in
$\hat{\mathfrak g}\eq E_8^{\sss(1)}$ for which
  \be
  (\mathfrak g_1,\mathfrak g_2) \;=\; (A_2,E_6)
  \quad{\rm or}\quad(A_1,E_7)
  \quad{\rm or}\quad(F_4,G_2)  \labl{E8}
and, in each case, $k_1\eq k_2\eq 1$.
The corresponding equivalences of modular tensor categories are known.
On the other hand, the two conditions are not met for the conformal embedding
into $E_8^{\sss(1)}$ of $A_2^{\sss(1)}{\oplus}\,A_1^{\sss(1)}$ with $k_1\eq6$
and $k_2\eq16$. In this case only categories of local modules, in fact 
so-called simple current extensions, for the two categories are equivalent.

The second application we have in mind concerns coset conformal field theories.
In these theories one starts from the \rep\ categories $\cC(\mathfrak g,k)$
and $\cC(\mathfrak h,k')$ for a pair of untwisted affine Lie algebras for
which $\mathfrak h\,{\subset}\,\mathfrak g$, and desires to understand the
\rep\ \cat\ of the commutant of the conformal vertex \alg\ associated to
$(\mathfrak h,k')$ in the conformal vertex \alg\ associated to $(\mathfrak g,k)$.
The results of the present paper will form an essential ingredient of a
universal description of these representation categories,
including, in particular, the so-called maverick coset theories.
A discussion of this application is beyond the scope of this introduction;
it will appear in a separate publication.

\subsection{Relation to earlier work}

The methods and results of this paper owe much to work that has been done
within two lines of development: the study of algebras in tensor categories, 
and alpha induction for nets of subfactors. Algebras in symmetric tensor
categories already played an important role in Deligne's characterisation of
Tannakian categories (see e.g.\ \cite{SAav,Demi}). They were studied in much
detail by Pareigis (see e.g.\ \cite{pare11,pare13}), who also introduced the 
concept of local (dyslectic) modules of a commutative algebra in a braided 
tensor category \cite{pare23}. More recently, commutative algebra and local 
modules in semisimple braided tensor categories were e.g.\ studied in the 
context of conformal field theory and quantum subgroups in \cite{kios}, in 
relation to weak Hopf algebras in \cite{ostr}, and in connection with
Morita equivalence for tensor categories in \cite{muge8}.
The algebras relevant in the conformal field theory context are
symmetric special Frobenius algebras \cite{fuSc16,fuRs,fuRs4}; those
encoding properties of conformal field theory on surfaces with boundary are,
generically, non-commutative. It is also worth mentioning that while bi- or 
Hopf algebras in braided \tcs\ (for a review, see \cite{maji33}) do not play 
a role in this context, they are indeed important for other applications in 
quantum field theory, see e.g.\ \cite{KElu}.  

The concept of $\alpha$-induction (see \Definition \ref{alphadef}) was invented 
in \cite{lore} in the framework of the $C^*$-algebraic approach to quantum 
field theory (see e.g.\ \cite{dohr,doro2}). $\alpha$-induction was further 
developed in \cite{xu3} and in a series of papers by B\"ockenhauer, Evans and 
Kawahigashi (see e.g.\ \cite{boev4,boek2,boek3,boev6}), in particular
applying it to the construction of subfactors associated to modular invariants,
and it was formulated in purely categorical form (and, unlike in the quantum
field theory and subfactor context, without requiring that one deals with
a *-category) in \cite{ostr}. Also in the study of subfactors Frobenius 
algebras arise naturally, in the guise of `Q-systems' \cite{long6,loro}. 
Indeed, every Q-system is a symmetric special *-Frobenius algebra \cite{evpi},
and the product and coproduct, and unit and counit, respectively, are *'s 
of each other (the Frobenius property can then in fact be derived from the 
other properties). For instance, the trivialising algebra defined in
\Lemma \ref{lem:TG} corresponds to the Q-system that is associated to the
canonical endomorphism of a subfactor, see \Proposition 4.10 of \cite{lore}.  

\vskip2em \noindent {\bf Acknowledgement.}\\
The collaboration leading to this work was supported in part by grant
IG 2001-070 from STINT (Stiftelsen f\"or internationalisering av h\"ogre
utbildning och forskning). J.Fu.\ is
supported in part by VR under contract no.\ F\,5102\,--\,20005368/2000,
and I.R.\ is supported by the DFG project KL1070/2--1.
    
\vskip3em 

%%%%%%%%%%%%%%%%%%%%%%%%%%%%%%%%%%%%%%%%%%%%%%%%%%%%%%%%%%%%%%%%%%%%%%%%

\sect{Algebras in tensor categories}\label{sect2}

\subsection{Tensor categories}\label{sect21}

Let $\cC$ be a category. We denote the class of its objects by $\Objc$
and the morphism sets by $\Hom(U,V)$, for $U,V$ in $\Objc$; we will often 
abbreviate endomorphism sets $\Hom(U,U)$ by $\End(U)$. In this paper we
will be concerned with categories that come with the following additional 
structure. First, they are small ($\Objc$ is a set), they are additive (so 
that, in particular, they have direct sums)
and their morphism sets are 
vector spaces over the ground field $\koerper$.
Second, most often they are tensor categories. By invoking the 
coherence theorems, tensor categories will be assumed to be strict; we 
denote the associative tensor product by `$\otimes$', both for objects and for 
morphisms, and the tensor unit by $\one$. 
Third, most of the categories we will be interested in are {\em ribbon\/}
categories; this includes as a special subclass the {\em modular\/} tensor
categories.
\\[-2.1em]

\dtl{Definition}{def:ribcat}
A {\em ribbon category\/} is a tensor category with the following
additional structure. To every object $U\iN\Objc$ one assigns an object
$U^\vee{\in}\,\Objc$, called the (right-) dual of $U$,
and there are three families of morphisms,%
 \foodnode{The existence of a duality is often included in the definition of a
 tensor category. What we refer to as a \tc\ is then called a {\em monoidal\/}
 category.}
  \be\bearll
  \mbox{(Right-) Duality:}&  \quad
    b_U \in \Hom(\one,U{\otimes}U^\vee) \,, \quad
    d_U \in \Hom(U^\vee{\otimes}U,\one) \,, \\[5pt]
  {\rm Braiding:} &  \quad c_{U,V} \in \Hom(U{\otimes}V,V{\otimes}U) \,, \\[6pt]
  {\rm Twist:}    &  \quad \theta_U \in \Hom(U,U)
  \eear\ee
for all $U\iN\Objc$, respectively for all $U,V\iN\Objc$, satisfying
  \be \bearll
  (d^{}_V\oti\id^{}_{V^\vee}) \cir (\id^{}_{V^\vee}\oti b^{}_V)=\id^{}_{V^\vee}
  \,,   &
  (\id^{}_V\oti d^{}_V) \cir (b^{}_V\oti \id^{}_V)=\id^{}_V \,,
  \\{}\\[-.7em]
  c^{}_{U,V\otimes W}= (\id^{}_V\oti c^{}_{U,W}) \cir (c^{}_{U,V}\oti \id^{}_W)
  \,, \ \ &
  c^{}_{U\otimes V,W}= (c^{}_{U,W}\oti \id^{}_V) \cir (\id^{}_U\oti c^{}_{V,W})
  \,,
  \\{}\\[-.7em]
  (g\oti f) \cir c^{}_{U,W} =c^{}_{V,X} \cir (f\oti g) \,,  &
  \theta^{}_{V}\cir f=f\cir\theta^{}_{U} \,,
  \\{}\\[-.7em]
  (\theta^{}_V\oti\id^{}_{V^\vee})\cir b^{}_V
  = (\id^{}_V\oti\theta^{}_{V^\vee})\cir b^{}_V \,,  &
  \theta^{}_{V\otimes W}
  = c^{}_{W,V}\cir c^{}_{V,W}\cir (\theta^{}_V\oti\theta^{}_W)
  \eear \labl{DTB}
for all $U,V,W,X\iN\Objc$ and all $f\iN\Hom(U,V)$, $g\iN\Hom(W,X)$.

\medskip

In a tensor category with duality, one defines the morphism dual to 
$f\iN\Hom(U,V)$ by $f^\vee\,{:=}\,(d_V\oti\id_{U^\vee}) \cir 
(\id_{V^\vee}\oti f\oti\id_{U^\vee}) \cir (\id_{V^\vee}\oti b_U) 
\iN \Hom(V^\vee\!{,}\,U^\vee)$. A left-duality is an assignment of a left-dual 
object $^{\vee\!}U$ to each $U\iN\Objc$ together with a family of morphisms,
  \be\bearll  \mbox{Left-duality:}
  & \quad \tilde b_U \in \Hom(\one,{}^{\vee\!}U{\otimes}U) \,, \quad
  \tilde d_U \in \Hom(U{\otimes}^{\vee\!}U,\one)\,, \eear\ee
that obey analogous properties as a right-duality.
In a ribbon category, there is automatically also a left-duality;
it can be constructed from right-duality, braiding and twist, and in fact
coincides with the right-duality both on objects and on morphisms,
${}^{\vee\!}U\eq U^\vee$, ${}^{\vee\!\!}f\eq f^\vee$. Tensor categories
with coinciding left- and right-duality functors from $\cC$ to $\cC^{\rm opp}$
are called {\em sovereign\/}. Thus, every ribbon category is in particular
sovereign; conversely, every braided sovereign category is a ribbon category.
For a \tc\ with both a left- and a right-duality, one defines left and right
{\em traces\/} of an endomorphism $f\iN\Hom(U,U)$ as 
  \be  {\rm tr}_{\rm l}(f) := d_U \cir (\id_{U^\vee}\oti f) \cir \tilde b_U \,,
  \qquad  {\rm tr}_{\rm r}(f) := \tilde d_U \cir (f\oti\id_{U^\vee}) \cir b_U \,,
  \ee
and the left and right (quantum) dimensions of an object $U$
as $\,{\rm dim}_{\rm l/r}(U)\,{:=}\,{\rm tr}_{\rm l/r}(\id_U)$. In a ribbon 
category the left and right traces coincide, i.e.\ ribbon categories are 
{\em spherical\/}.
Accordingly, in a ribbon category we denote the trace just by \,${\rm tr}$.

The properties \erf{DTB} of the braiding, twist and duality morphisms in a 
strict \tc\ allow us to visualise them via ribbon graphs (see e.g.\ 
\cite{joSt5} and chapter XIV of \cite{KAss}). In the sequel we will make 
ample use of this graphical notation. When drawing such graphs we follow 
the conventions set up in \Section~2 
of \cite{fuRs4}. In particular, all 
diagrams are to be read from bottom to top and, for simplicity, we use the 
blackboard framing convention so that ribbons can be drawn as lines. For 
convenience, we have also summarised the basic structural data in an 
appendix; the graphs for the braiding, twist and duality morphisms are 
collected in appendix \ref{apptab1}, and the graphical transcription of 
the axioms \erf{DTB} is given in appendix \ref{apptab2}.
To give another example, the graph
  \bea \begin{picture}(100,27)(0,30)
  \put(44,0)  {\begin{picture}(0,0)(0,0)
              \scalebox{.38}{\includegraphics{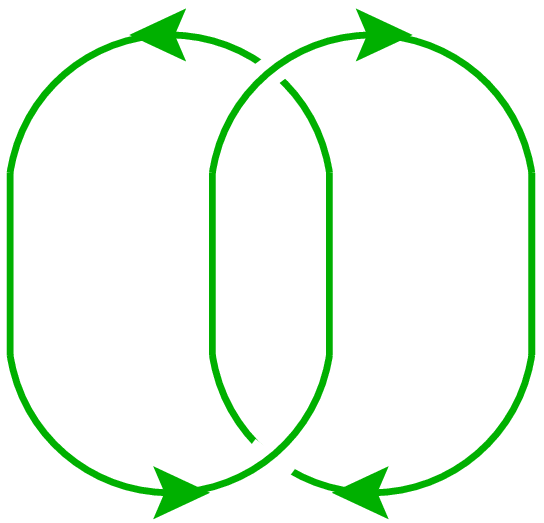}} \end{picture}}
  \put(-2,25.3)  {$s_{U,V}^{}\,:=$}
  \put(59.4,22.8){\sse$U$}
  \put(81.3,22.8){\sse$V$}
  \epicture08 \labl{sUV}
with $U,V$ any pair of objects of a braided tensor category with dualities,
is the trace
  \be  s_{U,V} = {\rm tr}(c_{U,V}^{}\cir c_{V,U}^{})
  = (d_{V}^{}\oti \tilde d_{U}^{})\cir
  [\, \id_{V^\vee}^{}\oti (c_{U,V}^{} {\circ}\, c_{V,U}^{})
  \oti\id_{U^\vee}^{} ] \cir (\tilde b_{V}^{}\oti b_{U}^{})  \ee
of the endomorphism $c_{U,V}^{}\cir c_{V,U}^{}$ of $V\oti U$.

When $\cC$ is semisimple, then we are particularly interested in simple 
objects. We denote by $\{U_i\,|\,i\iN\cI\}$ a collection of representatives 
of the isomorphism classes of (non-zero) simple objects of $\cC$, and set
  \be  \N ijk := \dim_\koerper(\Hom(U_i\oti U_j,U_k))  \labl{Nijk}
(taking values in $\zet_{\ge 0} \,{\cup}\, \{\infty\}$). Assuming that the 
tensor unit is simple, we take it as one of these representatives, so that 
$\II\,{\ni}\,0$ with $U_0\eq\one$. If an object $U$ is simple, then so is 
its dual $U^\vee$; thus in particular for every $i\iN\II$ there is a unique 
label $\bar\imath\iN\II$ such that $U_{\bar\imath}\,{\cong}\,U_i^\vee$.

\dtl{Definition}{def-modular}
A {\em modular tensor category\/} is a semisimple additive ribbon category 
for which the index set $\II$ is finite and
 % the property of semisimplicity replaced by the weaker dominance property,
 % and additiveness by the weaker "Ab" in the qualification `modular'
 % in the original definition in \cite{TUra}.
for which the $s$-{\em matrix\/} $s\eq(s_{i,j})^{}_{i,j\in\cI}$ with entries
  \be  s_{i,j}^{} := s_{U_i,U_j}^{}
  = {\rm tr}(c_{U_i,U_j}^{}{\circ}\, c_{U_j,U_i}^{})
  \labl{sij}
is non-degenerate.  

\medskip

Instead of non-degeneracy of $s$, one can equivalently \cite{brug2}
require that, up to isomorphism, the tensor unit is the only `transparent' 
simple object, i.e.\ that any simple object $U$ for which $c_{V,U}\cir c_{U,V} 
\eq \id_{U \oti V}$ holds for all $V\iN\Objc$ satisfies $U \,{\cong}\, \one$.  

\medskip

The dimension of an object $U\iN\Obj(\cC)$ is expressed
through the numbers \erf{sUV} as $\,{\rm dim}(U)\,{\equiv}\,{\rm tr}\,\id_{U}
\eq s_{U,\one} \eq s_{\one,U}$. In a modular tensor category, the square of the
matrix $s$ is, up to a multiplicative constant, a permutation matrix,
  \be  (s^2)^{}_{i,j} = \delta_{i,\bar\jmath}\sum_{k\in\II}(\dim(U_k))^2_{}
  \,. \labl{ss=c}
(In the physics literature, one usually considers the field of complex numbers, 
and instead of using $s$ it is more conventional 
to work with the unitary matrix $S$ defined as $S\,{:=}\, S_{0,0}\,s$ with 
$S_{0,0}\,{:=}\,\llb\sum_{i\in\II} (\dim(U_i))^2_{} \lrb^{-1/2}$.)

\medskip

For later reference we quote the following criterion for a functor $F$ to 
be an equivalence of categories (see e.g.\ \Theorem IV.4.1 of 
\cite{MAcl}).
\\[-2.3em]

\dtl{Proposition}{XI.1.5}
A functor $F$ is an equivalence of \cats\ if and only if
$F$ is essentially surjective (i.e.\ surjective up to
isomorphisms) and fully faithful (i.e.\ bijective on morphisms).

\medskip

Also note that when a functor $F$ is an equivalence of {\em braided\/} tensor 
categories, then, owing to the uniqueness properties of the left and right 
dualities and the fact that the twist can be expressed through the dualities 
and the braiding, $F$ is even an equivalence of {\em ribbon\/} categories.

\medskip

We will occasionally have to deal with constructions which, just like
functors, assign to each object $U$ of a \cat\ $\cC$ an object $F(U)$ of a 
\cat\ $\cD$, and to each morphism $f$ of $\cC$ a morphism $F(f)$ of $\cD$
in a manner compatible with the domain and target structure (i.e.\ such that 
$F(f)\iN\Hom(F(U),F(V))$ for $f\iN\Hom(U,V)$), but which are not, or are not 
known to be, functors. For definiteness, we will call a collection of maps 
that has these properties an {\em \operation\/} on the \cat\ $\cC$.

%%%%%%%%%%%%%%%%%%%%%%%%%%%%%%%%%%%%%%%%%%%

\subsection{Idempotents and retracts}

In order to fix our conventions and notation for subobjects and retracts
we review a few notions from category theory (for more details see e.g.\ 
\Sections I.5, V.7, VIII.1 and VIII.3 
of \cite{MAcl}). For brevity, in this 
description we often dispense with naming the source and target objects of 
a morphism explicitly; the corresponding statements are meant to hold for 
every object for which they can be formulated at all.

A morphism $e$ is called {\em monic\/} iff $e\cir f \eq e\cir g$ implies 
that $f\eq g$. A morphism $r$ is called {\em epi\/} iff $f\cir r\eq g\cir r$ 
implies that $f\eq g$. A {\em subobject\/} of an object $U$ is an equivalence 
class of monics $e\iN\Hom(\,{\cdot}\,,U)$. Here two monics $e\iN\Hom(S,U)$ 
and $e'\iN\Hom(S',U)$ are called equivalent iff there exists an isomorphism 
$\varphi\iN\Hom(S,S')$ such that $e\eq e'\cir\varphi$. A subobject $(K,e)$ of 
$U$ is a {\em kernel\/} of $f\iN\Hom(U,V)$ iff $f\cir e\eq 0$ and for every 
$h\iN\Hom(W,U)$ with $f\cir h\eq 0$ there exists a unique $h'\iN\Hom(W,K)$ 
such that $h\eq e \cir h'$. If a kernel exists, it is unique up to equivalence 
of subobjects. Cokernels are defined by reversing all arrows. The {\em image\/} 
$\Im f$ of a morphism $f$ is the kernel of the cokernel of $f$. It is often 
convenient to think of an isomorphism class of subobjects, kernels or 
cokernels as a single pair $(S,f)$. This is done by selecting a definite 
representative of the isomorphism class, invoking the axiom of choice 
(recall that all categories we consider are small).

A subobject $S$ is called {\em split\/} iff together with the monic
$e\iN\Hom(S,U)$ there also comes a morphism $r\iN\Hom(U,S)$ such that 
$r\cir e\eq \id_S$ (the letters $e$ and $r$ remind of `embedding' and 
`restriction'/`retract', respectively). We refer to the triple $(S,e,r)$ as a
{\em retract\/} of $U$ (just like for subobjects, we use the term retract both 
for the corresponding equivalence class of such triples and for an individual 
representative). We use the notations $S\,{\prec}\,U$ and $U\,{\succ}\,S$ to 
indicate that there exists a retract $(S,e,r)$ of $U$; when it is clear from 
the context what retract we are considering, we also use the abbreviations
$e\,{\equiv}\,e_S\,{\equiv}\,e_{S\prec U}$ and 
$r\,{\equiv}\,r_{\!S}\,{\equiv}\,r_{U\succ S}$. In the pictorial notation we 
will use the following shorthands for the morphisms $e,r$ specifying a retract:
  \bea \begin{picture}(95,42)(0,20)
  \put(0,0)   {\begin{picture}(0,0)(0,0)
              \scalebox{.38}{\includegraphics{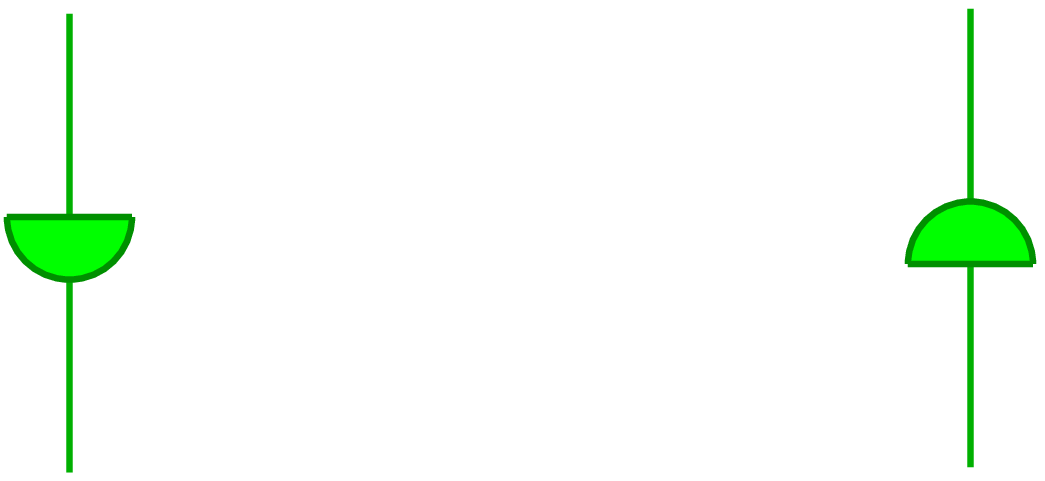}} \end{picture}}
  \put(-32,22.5) {$e\;=$}
  \put(4.8,-9.2)   {\sse$S$}
  \put(5.1,56.5)   {\sse$U$}
  \put(66,22.5)    {$r\;=$}
  \put(103.3,-9.2) {\sse$U$}
  \put(104.5,56.5) {\sse$S$}
  \epicture10 \labl{e-r}
Two retracts $S,\;S'$ are called equivalent iff $(S,e)$ and $(S',e')$
are equivalent as subobjects and $e\cir r\eq e'\cir r'$.

\smallskip

An endomorphism $p\iN \Hom(U,U)$ is called an {\em idempotent\/} (or a 
{\em projector\/}) iff $p \cir p \eq p$. To every retract $S\eq (S,e,r)$ of 
$U$ there is associated an idempotent $P_S\iN\Hom(U,U)$, namely 
$P_S\,{:=}\,e \cir r$. An idempotent $p$ is said to be {\em split\/} if, 
conversely, there exists a retract $(S,e,r)$ with $p\eq P_S \,{\equiv}\, 
e\cir r$. Thus a split idempotent has in particular an image, 
${\rm Im}(p)\eq S$, and split subobjects are precisely the images of 
split idempotents. Further, the retract $(S,e,r)$ is then unique up to 
equivalence of retracts, and
  \be
  e \circ r = p \,,   \qquad
  r \circ e = \id_S \,, \qquad
  p \circ e = e  \,,  \qquad
  r \circ p = r  \,.  \labl{eq:S-prop}
Also note that in a sovereign tensor category it follows, via the
cyclicity of the trace, that
$\tr_U(p)\eq\dim(\Im p)$, both for the left and the right trace,
for any split idempotent $p$.

\dtl{Lemma}{lem:rezI}
(i)~~For any two objects $U,V$ and any split 
idempotent $p\iN\Hom(U,U)$, there is a
natural bijection between the vector space $\Hom(\Im p,V) $ and the subspace
  \be  \Hom_{(p)}(U,V) := \{ f\iN\Hom(U,V) \,|\, f \cir p\eq f \}  \ee
of $\Hom(U,V)$.
\\[.3em]
(ii) For any two objects $U,V$ and any split
idempotent $q\iN\Hom(V,V)$, there is a
natural bijection between the vector space $\Hom(U,\Im q) $ and the subspace
  \be  \Hom^{(q)}(U,V) := \{ f\iN\Hom(U,V) \,|\, q \cir f\eq f \}  \labl{eq:dim-up-p}
of $\Hom(U,V)$.

\medskip\noindent
Proof: \\
Recall from the remarks before \erf{eq:S-prop} that $\Im p$ is in a
canonical way a retract $(\Im p,e,r)$ of $U$. With the help of the relations
\erf{eq:S-prop} one checks immediately that the map
$\Hom(\Im p,V)\,
    $\linebreak[0]$%
{\ni}\,\varphi\,{\mapsto}\,\varphi\cir r$ maps to the correct
subspace $\Hom_{(p)}(U,V)\,{\subseteq}\,\Hom(U,V)$ and that it has the map
$\Hom_{(p)}(U,V)\,{\ni}\,\psi\,{\mapsto}\,\psi\cir e$ as a two-sided inverse.
This establishes (i). Statement (ii) follows analogously, the relevant
mappings now being $\varphi\,{\mapsto}\,e\cir \varphi$ and
$\psi\,{\mapsto}\,r\cir\psi$.
\qed

\dtl{Definition}{Karoubian}
A \cat\ $\cC$ is called {\em Karoubian\/} (or {\em idempotent complete\/},
or {\em pseudo-abelian\/})
iff every idempotent is split.

\dtl{Remark}{K-rem0}
To every idempotent $p\iN\Hom(U,U)$ in an additive Karoubian \cat\ there 
corresponds
an isomorphism $U\,{\cong}\,{\rm Im}(p)\,{\oplus}\,{\rm Im}(\id_U{-}p)$.
All abelian \cats, as well as all additive semisimple \cats, are Karoubian.

\dtl{Definition}{Karoubi}
The {\em Karoubian envelope\/} (or {\em idempotent completion\/}, or 
{\em pseudo-abelian hull\/}) $\kar\cC$ of a \cat\ $\cC$ is a Karoubian \cat\ 
$\kar\cC$ together with an embedding functor $K{:}\;\cC\,{\to}\,\kar\cC$ 
which is universal in the sense that every functor $F{:}\;\cC\,{\to}\,\cD$ 
to a Karoubian \cat\ $\cD$ factors as $F\eq G\cir K$, with the functor 
$G{:}\;\kar\cC\,{\to}\,\cD$ unique up to isomorphism of functors.

\dtl{Remark}{K-rem1}
(i)~\,\,In the original definition of Karoubian envelope \cite{KAro} it is
also assumed that the \cat\ $\cC$ is additive, and the functors $K$ and $F$
are required to be additive functors. $\kar\cC$ is then an additive \cat, too.
\\[.3em]
(ii)~\,By general nonsense concerning universal properties, the Karoubian 
envelope is unique up to equivalence of categories. When $\cC$ is already 
Karoubian, then $\kar\cC\,{\cong}\,\cC$ and $K\,{\cong}\,\Id_\cC$.
\\[.3em]
(iii)~The Karoubian envelope of $\cC$ can be realised 
   \cite{KAro} 
as the \cat\ whose objects are pairs $(U;p)$ 
of objects $U\iN\Objc$ and idempotents $p\iN\Hom(U,U)$, and with morphisms
  \be  \kar\Hom((U;p),(V;q)) := \{ f\iN\Hom(U,V) \,|\, q\cir f\cir p\eq f \}
  \labl{karHom}
and $\kar\id_{(U;p)}\eq p$, so that in particular $p\iN\kar\Hom((U;p),(U;p))$.
In this realisation the embedding functor $K$ acts as $K(U)\eq(U;\id_U)$ and
$K(f)\eq f$, implying for instance that $\kar\Hom(K(U),K(V))\eq\Hom(U,V)$.
As a consequence, we may (and will) think of $\cC$ as a full sub\-\cat\ of
$\kar\cC$, and accordingly identify $U\iN\Objc$ with $(U;\id_U)\iN\Obj(\kar\cC)$.
\\
Further, when $q$ is any idempotent in $\kar\Hom((U;p),(U;p))$, we have 
$q\cir p\eq q\eq p\cir q$, implying that ${\rm Im}(q)\,{\cong}\,(U;q)$, 
independently of $p$.
\\[.3em]
(iv)~Various properties of $\cC$ are naturally
inherited by $\kar\cC$ (compare e.g.\ \cite{bklt}):\\[2pt]
a) If $\cC$ is tensor, then $\kar\cC$ becomes a \tc\ by setting
$f\,{\kar\otimes}\,g\,{:=}\,f\oti g$ and
  \be  \kaR\one := K(\one) \qquad{\rm and}\qquad
  (U;p)\,{\kar\otimes}\,(V;q) := (U\Oti V;p\Oti q) \,.  \labl{eq:K-tensor}
b) If a \tc\ $\cC$ is braided, then a braiding for the \tc\ $\kar\cC$ is
given by
  \be  \kar c_{(U;p),(V;q)} := (q\oti p)\cir c^{}_{U,V} \,.  \labl{eq:216}
c) If a \tc\ $\cC$ has a left duality, then a left duality for the
\tc\ $\kar\cC$ is given by $(U,p)^\vee\,{:=}\,(U^\vee,p^\vee)$ and
  \be  \kar d_{(U;p)} := d_U^{} \cir (\id^{}_{U^\vee} \oti p)
  \qquad{\rm and}\qquad
  \kar b_{(U;p)} := (p \oti \id^{}_{U^\vee}) \cir b_U^{} \,.  \labl{eq:217}
An analogous statement holds for a right duality.
\\[1pt]
d) If a braided \tc\ $\cC$ with duality has a twist, then a twist for
$\kar\cC$ is given by $\kar\theta_{(U;p)}\,{:=}\,p\cir\theta_U^{}$.
It follows in particular that when $\cC$ is ribbon, then $\kar\cC$
carries a natural structure of ribbon \cat\ as well.
\\
Further, dimensions in $\kar\cC$ are given by
  \be  \kar\dim((U;p)) = {\rm tr}(p) \,.  \ee
(v)~\,By the observation in \Remark \ref{K-rem0}
it thus follows that for every idempotent $p\iN\Hom(V,V)$ in an additive 
Karoubian ribbon \cat\ one has $\dim(V)\eq\dim({\rm Im}(p)) \,{+}\, 
\dim({\rm Im}(\id_V{-}p))$. In particular, if all dimensions are non-negative 
real numbers, then $\dim(U)\,{\le}\,\dim(V)$ if $U$ is a retract of $V$.

\dtl{Lemma}{le:K-funct}
For $F{:}\; \cC\,{\to}\,\cD$ a functor between \cats\ $\cC$ and $\cD$,
let $\kar F{:}\; \kar\cC\,{\to}\,\kar\cD$ be the functor between their
Karoubian envelopes given by
  \be 
  \kar F( (U;p) ) := (F(U);F(p))  \qquad{\rm and}\qquad
  \kar F(f) := F(f)  \ee
for objects $(U;p)$ and morphisms $f$ of $\kar\cC$.
\\[.2em] 
(i)~\,\,If $F$ is an equivalence functor, then so is $\kar F$.
\\[.2em] 
(ii)~\,If $\cC$ and $\cD$ are tensor categories and $F$ is a tensor functor,
  then $\kar F$ is a tensor functor, too.
\\[.2em] 
(iii)~If $\cC$ and $\cD$ are ribbon categories and $F$ is a ribbon functor,
  then $\kar F$ is a ribbon functor, too.

\medskip\noindent
Proof:\\
(i)~\,is derived easily by using the criterion of \Proposition \ref{XI.1.5}
for a functor to be an equivalence.
\\
(ii)~and~(iii)~follow by combining the respective properties of $F$ with the
prescription given in \Remark \ref{K-rem1}(iv) 
for the tensor and ribbon structure on the Karoubian envelope of a tensor 
and ribbon category, respectively.
\qed

\medskip

In the applications to rational conformal quantum field theory, the categories 
of main interest are ribbon categories that are even modular in the sense of 
\Definition \ref{def-modular}. In the present paper, also categories with much 
less structure play a role. However, a few basic properties (shared in 
particular by modular tensor categories) will generally be required below. We 
will not mention these properties repeatedly, but rather collect them in the
   \\[-2.5em]

\dtl{\Convention}{basicprops}
(i)~\,Every category $\cC$ is a small additive \cat, with all morphism
sets being vector spaces over some fixed field $\koerper$.
\\
Whenever a tensor category is not strict, we tacitly replace it by
an equivalent strict tensor category.
\\[.2em]
(ii)~\,\,Unless stated otherwise, every category is a assumed to be Karoubian. 
\\[.2em]
(iii)~Unless stated otherwise, the tensor unit $\one\iN\Objc$ of a tensor category
$\cC$ is simple, as well as absolutely simple, i.e.\ satisfies
$\End(\one)\eq\koerper\,\id_\one$.
\\[.3em]
For the categories from which our considerations start, all these properties
are {\em assumptions\/}. On the other hand, various constructions of new 
categories that we deal with in this paper -- taking the Karoubian envelope
(introduced in \Definition \ref{Karoubi}), 
the Karoubian product (see \Definition \ref{def:Ti}(ii)), 
the dual (\Definition \ref{def:dual-cat}), 
the category of modules over an algebra, and the category of local modules 
over a commutative \ssFA\ (\Definition \ref{def:ext-cat}) 
-- preserve the properties in part (i) and (ii) of the \convention; 
the procedures of taking the Karoubian envelope, the dual, or the Karoubian 
product in addition also preserve the properties stated in part (iii). 
Below this permanence will be mentioned only when it is non-trivial.

\dtl{Definition}{maxidm}
For $U$ an object in a (not necessarily Karoubian) \cat\ $\cC$,
let $H$ be a subset of the set ${\rm Idem}(U)$ of idempotents in $\End(U)$.
\\[.3em]
(i)~\,A {\em maximal idempotent in} $H$ is a morphism $P_{\!\rm max}^H\iN H$ 
such that
  \be  q \circ P_{\!\rm max}^H = q = P_{\!\rm max}^H \circ q  \ee
for all $q\iN H$.
\\[.3em]
(ii)~A {\em maximal retract with respect to} $\!H$ is a retract of $U$ such 
that $P_U$ is a maximal idempotent in $H$.

\dtl{Lemma}{lem:maxidm}
If a set $H\,{\subseteq}\,{\rm Idem}(U)$ contains a maximal idempotent,
then this maximal idempotent is unique.

\medskip\noindent
Proof:\\
Let $P_{\!\rm max}^{}$ and $P_{\!\rm max}'$ be two maximal idempotents in $H$.
Then $P_{\!\rm max}^{}\cir P_{\!\rm max}'\eq P_{\!\rm max}'$ by the maximality of
$P_{\!\rm max}^{}$ and $P_{\!\rm max}^{}\cir P_{\!\rm max}'\eq P_{\!\rm max}^{}$ 
by the maximality of $P_{\!\rm max}'$.
\qed

\dtl{Corollary}{cor:maxret}
If a maximal retract with respect to some $H\,{\subset}\,{\rm Idem}(U)$
exists, then it is unique up to isomorphism of retracts.

\dtl{Lemma}{KLEIMA}
Let $H\,{\subseteq}\,{\rm Idem}(U)$ be a set of idempotents on an object $U$
for which a maximal retract $P_{\!\rm max}^{}$ exists and is
split. Then for any split idempotent $P\iN H$, the image ${\rm Im}(P)$
is a retract of ${\rm Im}(P_{\!\rm max}^{})$.

\medskip\noindent
Proof:\\
We realise both ${\rm Im}(P)$ and ${\rm Im}(P_{\!\rm max}^{})$ as retracts 
of the object $U$, i.e.\ write $({\rm Im}(P),e,r)$ as well as 
$({\rm Im}(P_{\!\rm max}^{}),e_{\rm max}^{},r_{\!\rm max}^{})$. Then the
morphisms $\tilde e\,{:=}\, r_{\!\rm max}^{}{\circ}\,e\iN\Hom({\rm Im}(P),
{\rm Im}(P_{\!\rm max}^{}))$ and $\tilde r\,{:=}\,
r\cir e_{\rm max}^{}\iN\Hom({\rm Im}(P_{\!\rm max}^{}),{\rm Im}(P))$
obey $\tilde r \cir \tilde e\eq \id_{{\rm Im}(P)}$
owing to the maximality of $P_{\!\rm max}^{}$.
\qed

%%%%%%%%%%%%%%%%%%%%%%%%%%%%%%%%%%%%%%%%%%%

\subsection{Frobenius algebras}\label{sec:alg-fun}

The notion of an algebra over some field $\koerper$ has an analogue
in arbitrary tensor \cats. A $\koerper$-algebra is then nothing but an 
  algebra,\,%
  \footnote{~In using the term `algebra' we follow the terminology in e.g.\
  \cite{pare23,kios,muge8}. In a large part of the categorical literature
  (see e.g.\ \cite{MAcl,pare11,stre8}), the term `monoid' is used instead.}
in the \cat-theoretic sense, in the particular \tc\
$\Vectk$ of vector spaces over the field $\koerper$.
\\[-2.3em]

\dt{Definition}
An (associative) {\em algebra\/} (with unit) $A$ in a \tc\ $\cC$ is a triple 
$(A,m,\eta)$ consisting of an object $A$ of $\cC$, a multiplication morphism 
$m\iN\Hom(A\Oti A,A)$ and a unit morphism $\eta\iN\Hom(\one,A)$, satisfying
  \be
  m\circ (m\oti\id_\AA) = m \circ (\id_\AA\oti m) \quad\;{\rm and} \quad\;
  m \circ(\eta\oti \id_\AA) = \id_\AA = m \circ (\id_\AA\oti \eta) \,.
  \labl{alg}

\medskip

Other algebraic notions familiar from $\Vectk$ generalise to arbitrary
tensor \cats, too. In particular, a {\em co-algebra\/} in $\cC$ is a triple
$(A,\Delta,\eps)$ consisting of an object $A$, a comultiplication
$\Delta\iN\Hom(A,A\Oti A)$ and a counit $\eps\iN\Hom(A,\one)$ possessing
coassociativity and counit properties that amount to `reversing
all arrows' in the associativity and unit properties \erf{alg}.
Again a pictorial notation for these morphisms is helpful; we set
  \bea \begin{picture}(270,45)(0,15)
  \put(0,4.5)  {\begin{picture}(0,0)(0,0)
               \scalebox{.38}{\includegraphics{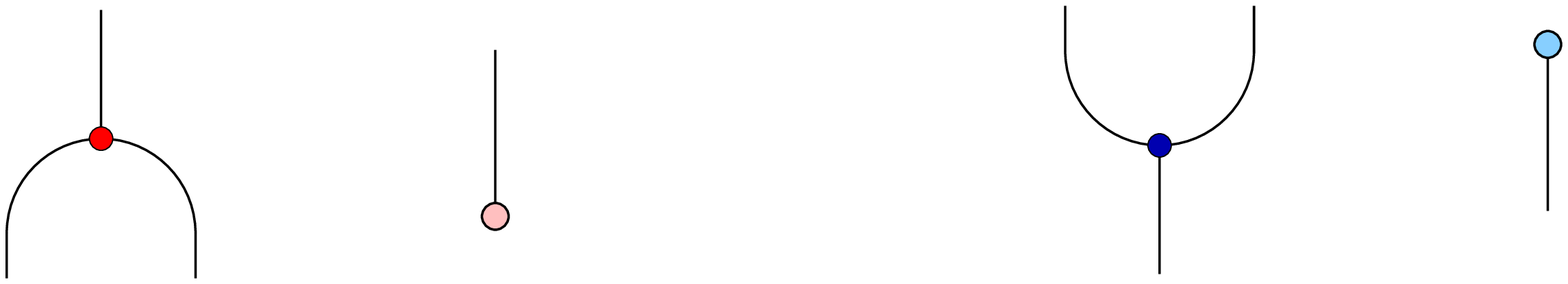}} \end{picture}}
  \put(-30,25)    {$m\,=$}
  \put(-3.5,-4.1) {\sse$A$}
  \put(12.8,52.5) {\sse$A$}
  \put(26.5,-4.1) {\sse$A$}
  \put(50.6,25)   {$\eta\,=$}
  \put(76.7,3.4)  {\sse$\one$}
  \put(76.1,45.4) {\sse$A$}
  \put(170,0) {\begin{picture}(0,0)(0,0)
  \put(-31.7,25)  {$\Delta\,=$}
  \put(-2.8,52.5) {\sse$A$}
  \put(11.9,-4.1) {\sse$A$}
  \put(28.1,52.5) {\sse$A$}
  \put(48.8,25)   {$\eps\,=$}
  \put(75.5,49.5) {\sse$\one$}
  \put(73.9,6.4)  {\sse$A$} \end{picture}}
  \epicture06 \labl{m-Delta}
Then e.g.\ the associativity of $m$ and coassociativity of $\Delta$ look like
  \bea \begin{picture}(375,45)(0,15)
  \put(25,0) {\begin{picture}(0,0)(0,0)
               \scalebox{.38}{\includegraphics{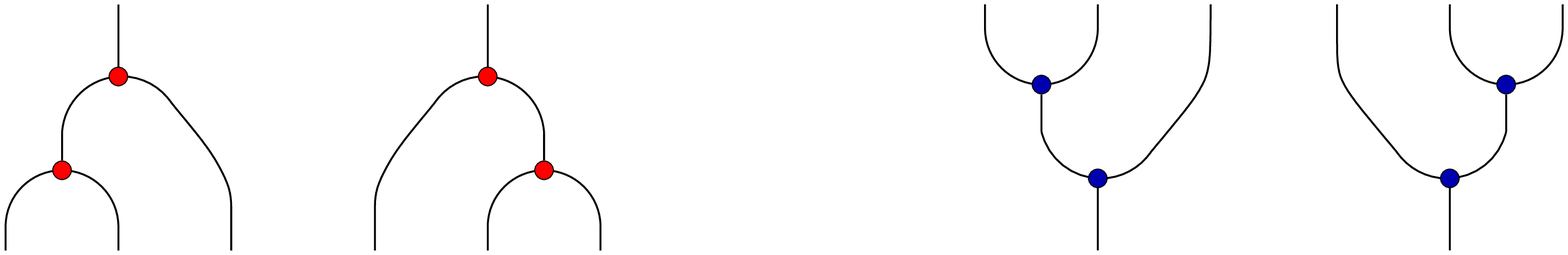}} \end{picture}}
  \put(21.1,-8.4) {\sse$A$}
  \put(43.8,-8.4) {\sse$A$}
  \put(45.2,53.5) {\sse$A$}
  \put(66.0,-8.4) {\sse$A$}
  \put(80.2,25)   {$=$}
  \put(94.6,-8.4) {\sse$A$}
  \put(117.5,-8.4){\sse$A$}
  \put(119.5,53.5){\sse$A$}
  \put(140.5,-8.4){\sse$A$}
  \put(172,25)  {and}
  \put(199.7,0) {\begin{picture}(0,0)(0,0)
  \put(19.9,54.2) {\sse$A$}
  \put(42.5,54.2) {\sse$A$}
  \put(41.8,-7.8) {\sse$A$}
  \put(65.5,54.2) {\sse$A$}
  \put(77.8,18)   {$=$}
  \put(91.5,54.2) {\sse$A$}
  \put(112.7,-7.8){\sse$A$}
  \put(113.5,54.2){\sse$A$}
  \put(135.5,54.2){\sse$A$}
  \end{picture}}
  \epicture06 \labl{ass-coass}
respectively.

\medskip

\dt{Definition}
A {\em left module\/} over an \alg\ $A\iN\Obj(\cC)$
is a pair $M\eq(\M,\r)$ consisting of an object $\M$ of $\cC$ and a
{\em \rep\ morphism\/} $\r\,{\equiv}\,\r_M^{}\iN\Hom(A\Oti\M,\M)$, satisfying
  \be  \r\circ(m\oti\id_\M) = \r \circ (\id_\AA\oti\r)
  \qquad\mbox{and}\qquad \r\circ(\eta\oti\id_\M) = \id_\M \,.  \labl{1m}

\medskip

By taking the $A$-modules as objects and the subspaces
  \be  \HomA(N,M) := \{ f\iN \Hom(\dot N,\M) \,|\,
  f \cir\r_N\eq\r_M\cir(\id_\AA\Oti f) \}   \labl{ha}
of the $\cC$-morphisms that intertwine the $A$-action as morphisms, one
gets the {\em category of left $A$-modules\/}, which we denote by \calca.
Analogously one defines right $A$-modules and their category. For
brevity we will often refer to left $A$-modules just as $A$-{\em modules\/}.
An $A$-module is called a {\em simple module\/} iff it is a simple object of
\calca. For $U\iN\Obj(\cC)$, the {\em induced\/} (left) {\em module\/}
$\Ind_\AA(U)$ is equal to $A{\otimes}U$ as an object in $\cC$, with
representation morphism $m\oti\id_U$; the full sub\cat\ of \calca\ whose 
objects are the induced $A$-modules will be denoted by \calcai.
(For more details see e.g.\ \cite{kios,fuSc16} and \Sections 4.1\,--\,3 
of \cite{fuRs4}.) When an $A$-module $N$ is a retract, as an object of \calca,
of an $A$-module $M$, we refer to it as a {\em \retmodule\/} of $M$.

\dtl{Remark}{K-rem2} 
(i)~\,If $(A,m,\eta)$ is an \alg\ in a \tc\ $\cC$, then
$((A;\id_\AA),m,\eta)$ is an \alg\ in its Karoubian envelope $\kar\cC$. 
Analogous statements hold for co\alg s, Frobenius \alg s etc.
\\[.3em]
(ii)~If $(\M,\r)$ is an $A$-module in a \tc\ $\cC$ and $p\iN\HomA(M,M)$ is
a split idempotent in $\calca$, then
  \be  ({\rm Im}(p),r{\circ}\r\,{\circ}(\id_\AA{\otimes}e))  \labl{Imp-rep}
(with $e\cir r\eq p$ as in \erf{eq:S-prop}) is an $A$-module in $\cC$, too.

\dtl{Lemma}{K-rem2iii}
For any \alg\ $A$ in a \tc\ $\cC$, the \cat\ $\kar{(\calca)}$
is equivalent to a full sub\-\cat\ of $(\kar\cC)_\AA^{}$.
\\[.2em]
In particular, if $\cC$ is Karoubian, then so is the \cat\ \calca\
of $A$-modules in $\cC$.

\medskip\noindent
Proof:\\
The first statement follows from the fact that if $M \eq (\M,\r)$ is an 
$A$-module in $\cC$ and 
$p\iN\HomA(M,M)$ is a (not necessarily split) idempotent, then
  \be  ((\M;p),p\cir\r)  \labl{Mp-rep}
is an $(A;\id_\AA)$-module in the Karoubian envelope $\kar\cC$.
\\
Since \calca\ is a full subcategory of $\kar{(\calca)}$, the second 
statement is a direct consequence of the first. More explicitly,
for any $A$-module $M\eq(\M,\r)$, every idempotent $p\iN\HomA(M,M)$ is
in particular an idempotent in $\Hom(\M,\M)$. Since $\cC$ is Karoubian,
there is thus a retract $({\rm Im}(p),e,r)$ in $\cC$. Defining
  \be  \rho_p := r\cir\rho\,{\circ}(\id_\AA{\otimes}e) \,,  \ee
we also have $e\cir\rho_p\cir(\id_\AA{\otimes}r) \eq p\cir\rho\cir
(\id_\AA{\otimes}p)$. Thus $({\rm Im}(p),\rho_p)\iN\Obj(\calca)$ is
a submodule of $M$, and hence $p$ is split as an idempotent in $\calca$.
\qed

\dtl{Remark}{K-rem2iv}
Conversely, if $((\M;p),\varrho)$ is an $(A;\id_\AA)$-module in
$\kar\cC$, with $p$ an idempotent that is already split in $\cC$, then
using the fact that $\varrho\iN\kar\Hom((A;\id_\AA)\Oti(\M;p),(\M;p))$
means (see \erf{karHom}) that
  \be  p\circ\varrho \circ (\id_\AA\oti p) = \varrho \,,  \ee
one checks that 
  \be  M_{p,\varrho} := ({\rm Im}(p),\varrho_p) \qquad{\rm with}\qquad
  \varrho_p := r\cir\varrho\cir(\id_\AA{\otimes}e)  \labl{Mpr}
is an $A$-module in $\cC$.
\\ 
Also, when the condition that the idempotent $p$ is split is not imposed
(so that $\Im(p)$ does not necessarily exist), one might be tempted to
directly interpret the pair $(\dot M,\varrho)$ as a module. But this is 
not, in general,
possible. While $(\M,\varrho)$ does satisfy the first \rep\ property
$\varrho\cir(\id_\AA\oti\varrho) \eq \varrho\cir(m\oti\id_\M)$,
for $p\,{\ne}\,\id_\M$ the second \rep\ property fails,
$\varrho \cir (\eta\oti\id_\M)\eq p$. 

\dt{Definition}
An $A$-{\em bimodule\/} is a triple $M\eq(\M,\r_{\rm l},\rr)$ such that
$(\M,\r_{\rm l})$ is a left $A$-module, $(\M,\rr)$ is a right $A$-module,
and the left and right actions of $A$ commute.

\bigskip

The category of $A$-bimodules in $\cC$ will be denoted by \calcaa. Note that
in contrast to \calca, this is always a tensor category (though
not necessarily braided).

In a {\em braided\/} tensor category, for every object $V$ the induced
{\em left\/} $A$-module $(A\Oti V,m\Oti\id_V)$ can be endowed in two obvious
ways with the structure of a {\em right\/} $A$-module $(A\Oti V,\rr^{\pm})$;
the \rep\ morphisms
$\rr^{\pm}\,{\equiv}\,\r_{V,{\rm r}}^{\pm}\iN\Hom(A\Oti V\Oti A,A\Oti V)$ are
  \be  \rr^+ := (m\oti\id_V) \cir (\id_\AA\oti c_{V,A}) \qquad{\rm and}
  \qquad \rr^- := (m\oti\id_V) \cir (\id_\AA\oti(c_{A,V})^{-1})\,, \labl{rr+-}
respectively. These are used in
\\[-2.3em]

\dtl{Definition}{alphadef}
For $A$ an \alg\ in a braided tensor category $\cC$, the functors
  \be  \alpha_\AA^{\pm}:\quad \cC \to \calcaa  \ee
of $\alpha$-{\em induction\/} are defined on objects as
  \be  \alpha_\AA^{\pm}(V) := (A\Oti V,m\Oti\id_V,\rr^\pm)  \ee
for $V\iN\Objc$, and on morphisms as
  \be  \alpha_\AA^{\pm}(f) := \id_\AA \oti f\,\in\Hom(A\Oti V,A\Oti W)  \ee
for $f\iN\Hom(V,W)$.

\medskip

The $\alpha$-inductions $\alpha_\AA^{\pm}$ are indeed functors, even tensor
functors, from $\cC$ to the category \calcaa\ of $A$-bimodules. They were first
studied in the theory of subfactors (see \cite{lore} and also e.g.\
\cite{xu3,boev6,boek3}), and were reformulated in the form used here in 
\cite{ostr}.

\medskip

We will mainly be interested in algebras with several specific additional
properties, which arise e.g.\ in applications to conformal quantum field theory
\cite{fuRs4}.
\\[-2.2em]

\dtl{Definition}{symm-frob-spec}
(i)~\,\,An algebra $A$ in a tensor category with left and right dualities
together with a morphism $\eps\iN\Hom(A,\one)$ is 
called a {\em symmetric\/} algebra iff the two morphisms
  \bea \begin{picture}(285,70)(0,12)
               \put(210,0) {\begin{picture}(0,0)(0,0)
               \scalebox{.38}{\includegraphics{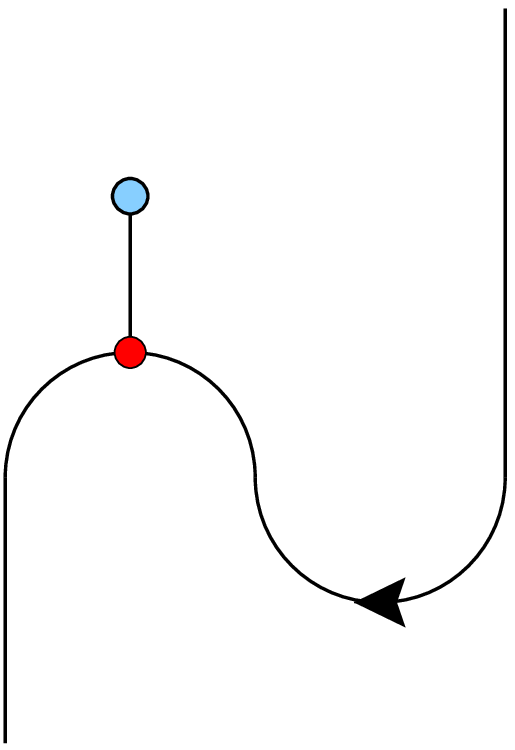}} \end{picture}}
  \put(0,30.4){$\Phi_1 \,:=\, [(\eps\cir m)\oti \id_{A^\vee}]
               \circ (\id_\AA \otimes b_A) \;=$}
  \put(206.1,-8.9){\sse$A$}
  \put(261.5,84.4){\sse$A^\vee$}
  \epicture01 \labl{eq:Phi-def}
and
  \bea \begin{picture}(285,69)(0,15)
               \put(210,0) {\begin{picture}(0,0)(0,0)
               \scalebox{.38}{\includegraphics{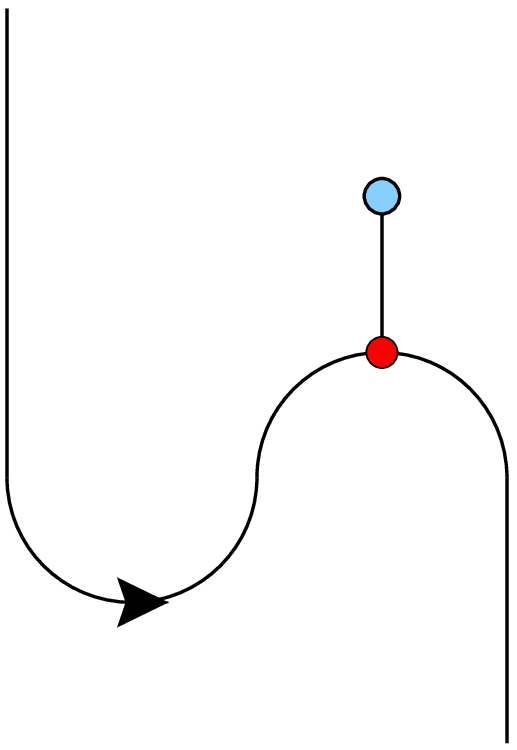}} \end{picture}}
  \put(0,38.4){$\Phi_2 \,:=\, [\id_{A^\vee}\oti (\eps\cir m)]
               \circ (\tilde b_A \oti \id_\AA) \;=$}
  \put(206.7,85.3){\sse$A^\vee$}
  \put(261.1,-9.1){\sse$A$}
  \epicture08 \labl{Phi2}
in $\Hom(A,A^\vee)$ are equal.
\\[.3em]
(ii)~\,A {\em Frobenius \alg\/} in a \tc\ $\cC$ is a quintuple
$(A,m,\eta,\Delta,\eps)$ such that $(A,m,\eta)$ is an \alg\ in $\cC$,
$(A,\Delta,\eps)$ is a co-\alg\
in $\cC$, and there is the compatibility relation
  \be  (\id_\AA\oti m) \circ (\Delta\oti\id_\AA)
  = \Delta \circ m = (m\oti\id_\AA) \circ (\id_\AA\oti\Delta)  \labl{1f}
between the two structures.
\\[.3em]
(iii)~A Frobenius \alg\ is called {\em special\/} iff
  \be  \eps\cir\eta = \beta_\one\, \id_\one
  \qquad {\rm and}\qquad
  m\cir\Delta = \beta_{\!A}\, \id_\AA \labl{eq:special-def}
for non-zero numbers $\beta_\one$ and $\beta_{\!A}$.

\medskip

Recently \cite{muge8,stre8}, in order to emphasise the analogy with classical
non-commutative ring theory (compare e.g.\ \cite{kaSt2}), the term 
``strongly separable'' was introduced for what we call ``special''.  

\medskip

For a \ssFA\ $A$ one has $\beta_\one\beta_{\!A}\eq\dim(A)$, implying in 
particular that $\dim(A)\,{\ne}\,0$. It is then convenient to normalise
$\eps$ and $\Delta$ such that $\beta_\one\eq\dim(A)$ and $\beta_{\!A}\eq1$.
We will follow this convention unless mentioned otherwise. We also set
  \be  \eps_\natural := d_A \circ (\id_{A^\vee}\oti m) \circ
  (\tilde b_A \oti\id_\AA) \,\in \Hom(A,\one)  \labl{eq:epsnat}
and write $\Phi_{1,\natural}$ for the morphism that is obtained by
replacing $\eps$ in the expression \erf{eq:Phi-def} by $\eps_\natural$.

\dtl{Remark}{prop:ssFA-unique}  
(i)~\,\,If $A$ is a special Frobenius algebra then, with the normalisation
$\beta_{\!A}\eq1$, $(A,\Delta,m)$ is a retract of $A\oti A$. The Frobenius 
property ensures that this statement holds even at the level of $A$-bimodules. 
This bears some similarity to the situation in braided tensor categories 
where the notion of a bi-algebra can be defined. In fact, the property of 
an algebra $A$ to be a bi-algebra is equivalent to the statement that the 
coproduct endows $A$ with the structure of a retract of $A\oti A$ as an 
algebra, rather than as a bimodule.
\\[.3em]
(ii)~\,When $\cC$ is semisimple and $A$ is special Frobenius, then the category
\calca\ of left $A$-mo\-dules is semisimple \cite{fuSc16}.
\\[.3em]
(iii)~For modules over any \alg\ $A$ in a tensor category $\cC$ a reciprocity
relation holds, stating that for every left $A$-module $M$ and every
object $U$ of $\cC$ there is a canonical bijection
  \be\begin{array}{rcl}
    \phi_1 :\quad \HomA(\Ind_\AA(U),M) 
    &\overset{\cong}{\longrightarrow}& \Hom(U,\M)  \\
    f &\overset{\cong}{\longmapsto}& f \cir (\eta{\otimes}\id_U)
  \eear \labl{reciUM}
between morphism spaces in $\cC$ and in \calca. If $A$ is Frobenius,
then an analogous reciprocity relation also holds when the target
of $\HomA$ is an induced module,
  \be\begin{array}{rcl}
    \phi_2 :\quad  \HomA(M,\Ind_\AA(V)) 
    &\overset{\cong}{\longrightarrow}& \Hom(\M,V)  \\
    g &\overset{\cong}{\longmapsto}& (\eps{\otimes}\id_V) \cir g \,.
  \eear \labl{reciMV}
In other words, for an arbitrary associative algebra, the induction
functor is a left adjoint functor of restriction; if the algebra carries
the additional structure of a Frobenius algebra, then induction is
a right adjoint of restriction as well.
\\
The inverses of the maps \erf{reciUM} and \erf{reciMV} can also be
given explicitly; they are 
  \bea
  \phi_1^{-1}(\tilde f) = \r_M \cir (\id_A{\otimes} \tilde f)
  \qquad {\rm and}\\{}\\[-.7em]
  \phi_2^{-1}(\tilde g) = [\id_A \oti (\tilde g\cir \rho_M)] \cir 
  [(\Delta \cir \eta) \oti \id_\M] \,.
  \eear\labl{eq:reci-inv}
\\
For a proof see e.g.\ \Propositions 4.10 and 4.11 
of \cite{fuSc16}.
\\[.3em]
(iv)~\,Given an algebra $(A,m,\eta)$ with $\dim(A)\,{\ne}\,0$ in a sovereign 
\tc\ $\cC$, morphisms $\Delta$ and $\eps$ such that $(A,m,\eta,\Delta,\eps)$ 
is a symmetric special Frobenius algebra exist iff the morphism
$\Phi_{1,\natural}$ defined after \erf{eq:epsnat} is invertible.
Further, it is always possible to normalise $\eps$ in such a way that
$\eps\eq\eps_\natural$. With this normalisation of the counit
the coproduct $\Delta$ is unique, and one has $\beta_\one\eq\dim(A)$.
\\
For a proof see lemma 3.12 of \cite{fuRs4}.
\\[.3em]
(v)~\,For every Frobenius algebra $A$ in a tensor category with left and right
duality, the morphisms $\Phi_{1,2}$ in \erf{eq:Phi-def} and \erf{Phi2} are 
invertible (see lemma 3.7 of \cite{fuRs4}), so that in particular 
$A\,{\cong}\,A^\vee$. Using $\Hom(U,V)\,{\cong}\,\Hom(V^\vee,U^\vee)$ it 
follows that for any two Frobenius algebras $A$, $B$ in a tensor category 
with left and right duality we have
  \be \Hom(A,B) \cong \Hom(B,A) \,.  \labl{eq:HomAB=HomBA}
For symmetric Frobenius algebras, there is a single distinguished isomorphism 
between the object underlying the algebra and its dual object, and as a 
consequence there is also a distinguished bijection \erf{eq:HomAB=HomBA}.
\\[.3em]
(vi)~\,If $\dim_\koerper \Hom(\one,A)\eq d$ for a Frobenius algebra $A$
in a tensor category with left and right duality, then
$I^{(d)}\,{:=}\,\one{\oplus}\one{\oplus}\cdots{\oplus}\one$ ($d$ summands)
is a retract of $A$. Indeed, as just remarked the morphisms $\Phi_{1,2}$
are then invertible, and it is not difficult to see that one can
choose $\alpha_i\iN\Hom(\one,A)$, $i\eq1,2,...\,,d$, such that
$\eps\cir m\cir(\alpha_i\oti\alpha_j)\eq\delta_{i,j}$.
(This furnishes a non-degenerate bilinear form on $\Hom(\one,A)$ and
thus endows $\Hom(\one,A)$ with the structure of a Frobenius algebra in
$\Vectk$.) Further, there are retracts $(\one,e_i,r_i)$ of $I^{(d)}$ satisfying
$r_i\cir e_j \eq\delta_{i,j}$ and $\sum_{i=1}^d e_i\cir r_i\eq\id_\AA$.
The retract in question is then given by $(I^{(d)},e,r)$ with
$e\,{:=}\,\sum_{i=1}^d \alpha_i\cir r_i$ and
$r\,{:=}\,\sum_{i=1}^d e_i\cir \eps\cir m\cir(\alpha_i\oti\id_\AA)$.  
\\[.3em]
(vii)~\,In terms of our graphical calculus, the property of an algebra to
be symmetric Frobenius in essence implies that multiplications and/or
comultiplications can be moved past each other in all possible arrangements.
Examples for such moves are provided by the defining properties
\erf{ass-coass} and \erf{1f}. Another move, which is frequently used in
our calculations below (without special mentioning), is the following:
  %% [pic~76]
  \bea
  \begin{picture}(260,97)(0,27)
  \put(0,0) {\begin{picture}(0,0)(0,0)
    \scalebox{.38}{\includegraphics{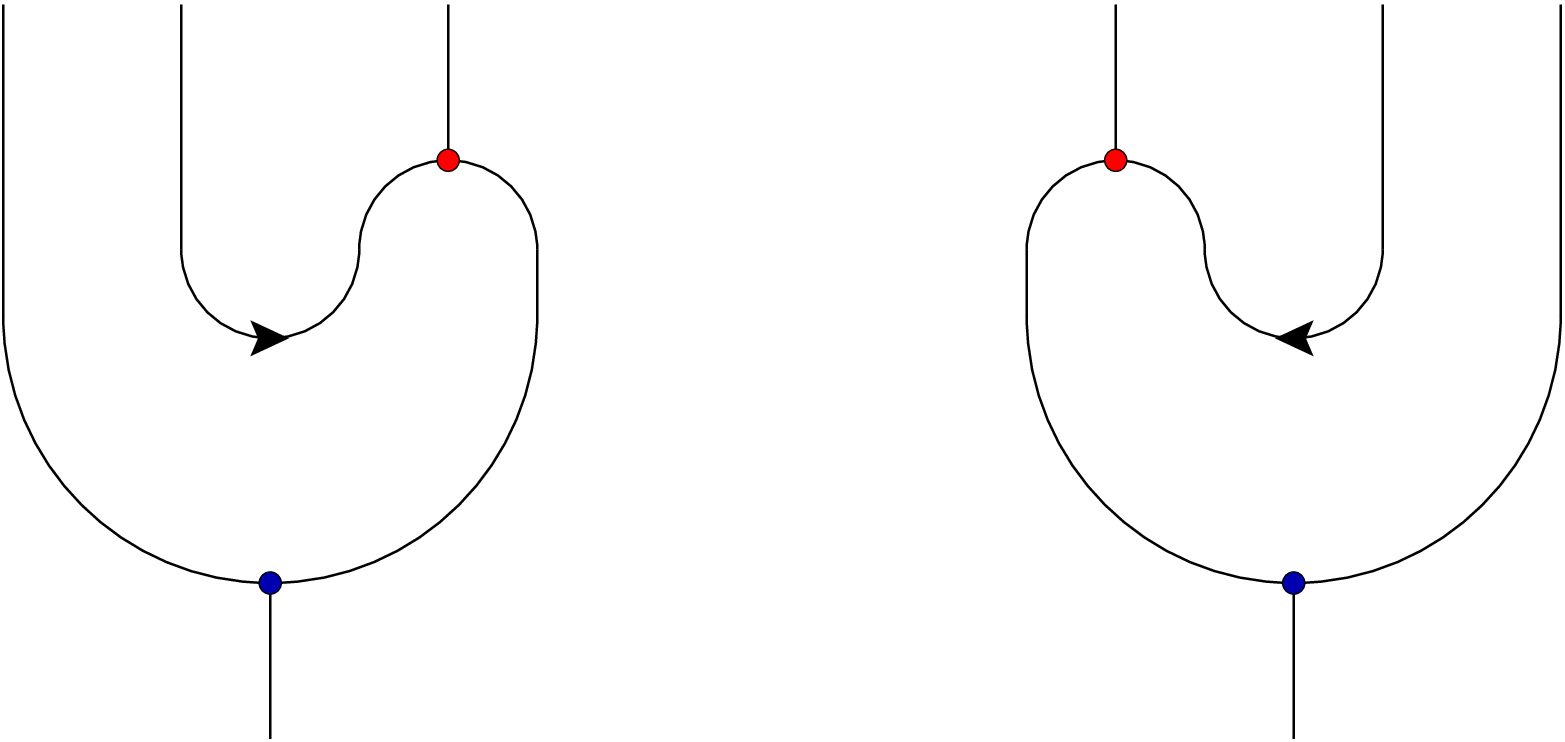}} \end{picture}
    \put(37.4,-9.2) {\sse$A$}
    \put(-3.2,118)  {\sse$A$}
    \put(22.7,118)   {\sse$A^\vee_{}$}
    \put(65.1,118)  {\sse$A$}
    \put(113,60) {=}
    \put(195.7,-9.2){\sse$A$}
    \put(168.9,118) {\sse$A$}
    \put(208,118)   {\sse$A^\vee_{}$}
    \put(237.7,118) {\sse$A$}
  }
  \epicture15 \labl{eq:more-assoc}
This identity uses both the symmetry and the Frobenius property.
First note that the latter two properties imply
  \be
  [(\eps\cir m) \oti \id_{A^\vee}] \circ
  [\id_A \oti (b_A \cir \tilde d_A)] \circ
  [(\Delta\cir\eta)\oti \id_{A^\vee}] = \id_{A^\vee} \,.  \ee
To show \erf{eq:more-assoc}, we insert this identity on the outgoing
$A^\vee$-ribbon on the left hand side. Then we use the Frobenius property 
to convert the  product $m$ on the \lhs\ of \erf{eq:more-assoc} into a 
coproduct. The latter can then be moved past the coproduct already present 
in \erf{eq:more-assoc} using coassociativity. Finally the coproduct is 
converted back to a product using the Frobenius property once again.  

\medskip

The properties of an \alg\ $A$ to be symmetric, special and Frobenius are
all indispensable in the construction of a \cft\ from $A$. Further
properties of $A$ can be important in specific applications. For us the
following is the most important one:
\\[-2.3em]

\dt{Definition}
An algebra $A$ in a braided \tc\ is said to be {\em commutative\/}
(with respect to the given braiding) iff $m\cir c_{A,A}\eq m$.

\medskip

Note that $m\cir c_{A,A}\eq m$ is equivalent to $m\cir c_{A,A}^{-1}\eq m$.
Also, while in general the category \calca\ of left $A$-modules is not a \tc, 
a sufficient condition for \calca\ to have a tensor structure is that
$A$ is commutative. The tensor structure is not canonical, though,
because in this case one can turn a left $A$-module into an $A$-bimodule in 
two different ways by using the braiding to define a right action of $A$.
However, we will see later (see also \cite{pare23,kios}) that for commutative 
$A$ there is a full subcategory of \calca, namely the category \calcal\ of 
local $A$-modules, on which the tensor structure is canonical.

\medskip

In the classical case, i.e.\ for algebras in the category of
finite-dimensional vector spaces over a field, commutative Frobenius
algebras are automatically symmetric. In a braided setting this is not true 
in general, but only if an additional condition
is satisfied. More precisely, we have
\\[-2.3em]

\dtl{Proposition}{c+s=tt}
(i)~\,\,A commutative symmetric Frobenius algebra has trivial twist, i.e.\
$\theta_\AA\eq\id_\AA$.
\\[.3em]
(ii)~\,Conversely, every commutative Frobenius algebra with trivial twist
is symmetric.
\\[.3em]
(iii)~A commutative symmetric Frobenius algebra is also cocommutative.

\medskip\noindent
Proof:\\
(i)~\,The statement follows from the equivalence
  %% [pic~61]
  \bea \begin{picture}(334,61)(0,27)
               \put(0,0) {\begin{picture}(0,0)(0,0)
               \scalebox{.38}{\includegraphics{Phi1.eps}} \end{picture}}
               \put(94,0) {\begin{picture}(0,0)(0,0)
               \scalebox{.38}{\includegraphics{Phi2.eps}} \end{picture}}
               \put(230,0) {\begin{picture}(0,0)(0,0)
               \scalebox{.38}{\includegraphics{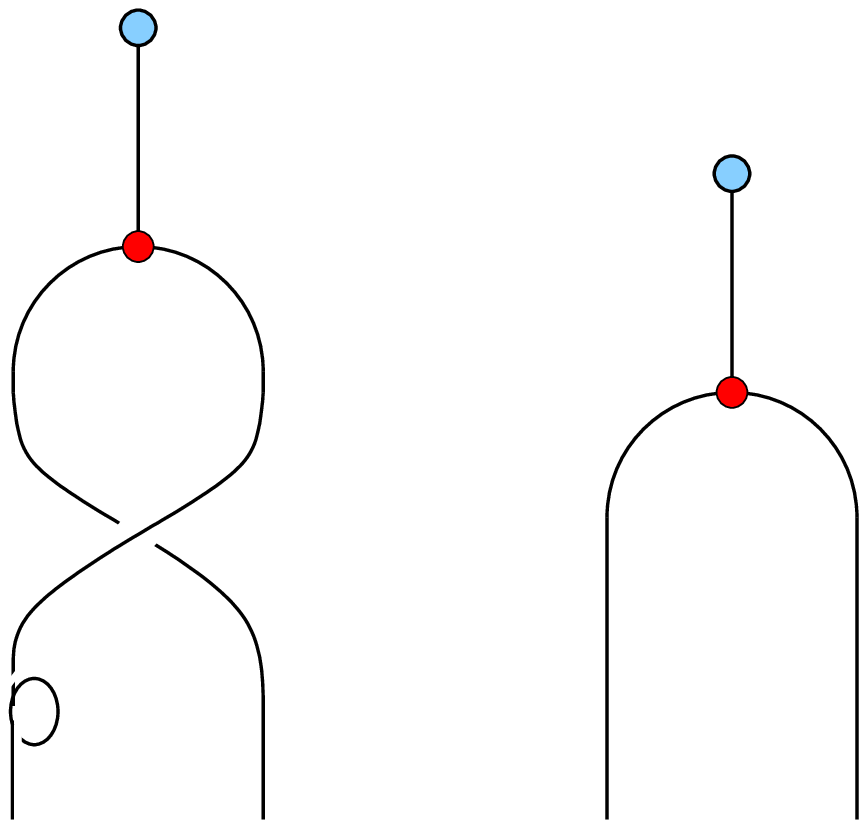}} \end{picture}}
  \put(-3.5,-8.4) {\sse$A$}
  \put(51.1,84.4) {\sse$A^\vee$}
  \put(71,37.4)   {$=$}
  \put(91.1,84.4) {\sse$A^\vee$}
  \put(145.5,-8.4){\sse$A$}
  \put(177,37.4)  {$\Longleftrightarrow$}
  \put(226.5,-8.4){\sse$A$}
  \put(255.5,-8.4){\sse$A$}
  \put(276,37.4)  {$=$}
  \put(293.9,-8.4){\sse$A$}
  \put(321.5,-8.4){\sse$A$}
  \epicture13 \labl{eq:th=1-proof}
This is obtained by bending the outgoing $A^\vee$-ribbon via a duality
morphism downwards to the left, which replaces the morphism $\Phi_{1,2}$ by
$\tilde d_{A}\cir(\id_\AA\oti\Phi_{1,2})$.
\\
Assuming symmetry, the left equality in \erf{eq:th=1-proof}
holds true. After using commutativity on the left equality,
the two sides differ only by a twist.
Tensoring this equality with $\id_\AA$ from the right and then
composing with $\id_\AA\oti(\Delta\cir\eta)$ and using the Frobenius
property finally yields $\theta_\AA\eq\id_\AA$.
Here, as well as in many arguments below, the manipulations of the
graphs also involve a process of `deforming' ribbons, making use
of the defining properties of the braiding (in particular, functoriality),
duality and twist.
\\[.3em]
(ii)~\,If $A$ is commutative and has trivial twist, the right
equality in \erf{eq:th=1-proof} holds, implying that $A$ is symmetric.
\\[.3em]
(iii) is obtained by the following moves:
  %% [pic~68]
  \bea \begin{picture}(330,59)(30,33)
  \put(0,0) {\begin{picture}(0,0)(0,0)
            \scalebox{.38}{\includegraphics{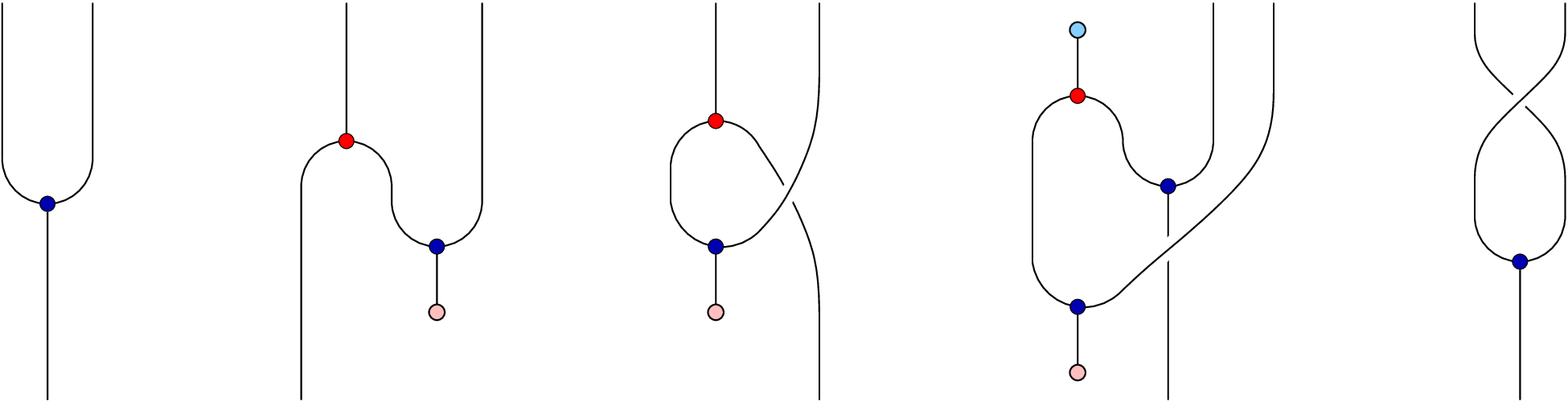}} \end{picture}}
  \put(40,44)  {$=$}
  \put(127,44) {$=$}
  \put(206,44) {$=$}
  \put(306,44) {$=$}
  \epicture13 \labl{pic68} 
Here the first equality uses the Frobenius and unit properties, the second
uses commutativity, the third the Frobenius and counit properties. The last
equality is based on the symmetry property together with the result of (i).
\qed

\dtl{Definition}{defsimple}
An algebra $A$ is called {\em simple\/}
iff all bimodule endomorphisms of $A$ as a bimodule over itself
are multiples of the identity, i.e.\ $\HomAA(A,A)\eq\koerper\,\id_\AA$.

\dtl{Lemma}{beibl}
Let $A$ be a simple symmetric special Frobenius algebra.
Then for every left $A$-module $M$ the equality
  \be
  \tilde d_\M \cir (\r_M \oti \id_{\M^\vee}) \cir (\id_A \oti b_\M)
  = \Frac{\dim(\M)}{\dim(A)}\, \eps  \labl{eq:beibl} 
holds.

\medskip\noindent
Proof:\\
The morphism $f\,{:=}\,[\id_A \oti (\tilde d_\M\cir(\r_M\Oti\id_{\M^\vee})
\cir(\id_A\Oti b_\M)) ] \cir \Delta$ is an $A$-bimodule morphism from $A$ 
to $A$. That $f$ is a left module morphism follows from the Frobenius property 
of $A$, while to show that it is also a morphism of right $A$-modules one 
needs $A$ to be symmetric. Since $A$ is simple, $f$ is thus a multiple of 
$\id_A$. Since $\eps\cir f$ gives the \lhs\ of \erf{eq:beibl}, the constant of 
proportionality is the same as on the \rhs\ of \erf{eq:beibl}. This constant, 
in turn, immediately follows from $\eps\cir f\cir\eta \eq {\rm tr}\,\id_\M$.
\qed

\dtl{Remark}{iZUV}
(i)~\,In the following an important role will be played by the dimensions of
certain spaces $\HomAA(M,N)$ of $A$-bimodule morphisms between $\alpha$-induced 
bimodules. (Recall the \Definition \ref{alphadef} 
of $\alpha$-induction.)
For any pair $U,V$ of objects and any algebra $A$ in a ribbon category we set
  \be  \tilde Z(A)_{U,V}^{} := \dim_\koerper\llb\HomAA
  (\alpha_\AA^{-}(V),\alpha_\AA^{+}(U))\lrb \,.  \labl{ZUV}
Since $A\eq\alpha_\AA^{\pm}(\one)$,
simplicity of $A$ as an algebra is thus equivalent to
  \be  \tilde Z(A)_{\one,\one} = 1 \,.   \ee
The corresponding notion with left (or right) module instead of bimodule 
endomorphisms is {\em haploidity\/}: $A$ is said to be haploid \cite{fuSc16} 
iff $\HomA(A,A)\eq\koerper\,\id_\AA$, i.e.\ iff $\,\dim\,\Hom(\one,A)\eq1$. 
Haploidity implies simplicity, but the converse is not true; however, if $A$ 
is {\em commutative\/}, then we have $\HomAA(A,A)\eq\HomA(A,A)$, so that in 
this case haploid and simple are equivalent. Moreover, every simple special 
Frobenius algebra in a semisimple \cat\ is Morita equivalent to a haploid 
algebra (see the corollary in \Section 3.3 of \cite{ostr}).  
\\[.3em]
(ii)~Every modular tensor
category gives rise to a three-dimensional topological field
theory, and thereby to invariants of ribbon graphs in three-manifolds.
(See e.g.\ \cite{KAss} or, for a brief summary, \Section 2.5 of \cite{fffs3}.)
In the three-dimensional TFT ribbons are labelled by objects of the 
underlying modular tensor category and coupons at which ribbons join by 
corresponding morphisms; our graphical notation for morphisms fits with 
the usual conventions for drawing the ribbon graphs.
For instance, as shown in \Section 5.4 of \cite{fuRs4}, when $A$ is a \ssFA\ 
in a (semisimple) modular tensor category, then the numbers \erf{ZUV} for
simple objects $U\eq U_i$ and $V\eq U_j$ coincide with the invariant
  \bea \begin{picture}(85,122)(0,34)
  \put(0,0)   {\begin{picture}(0,0)(0,0)
              \scalebox{.38}{\includegraphics{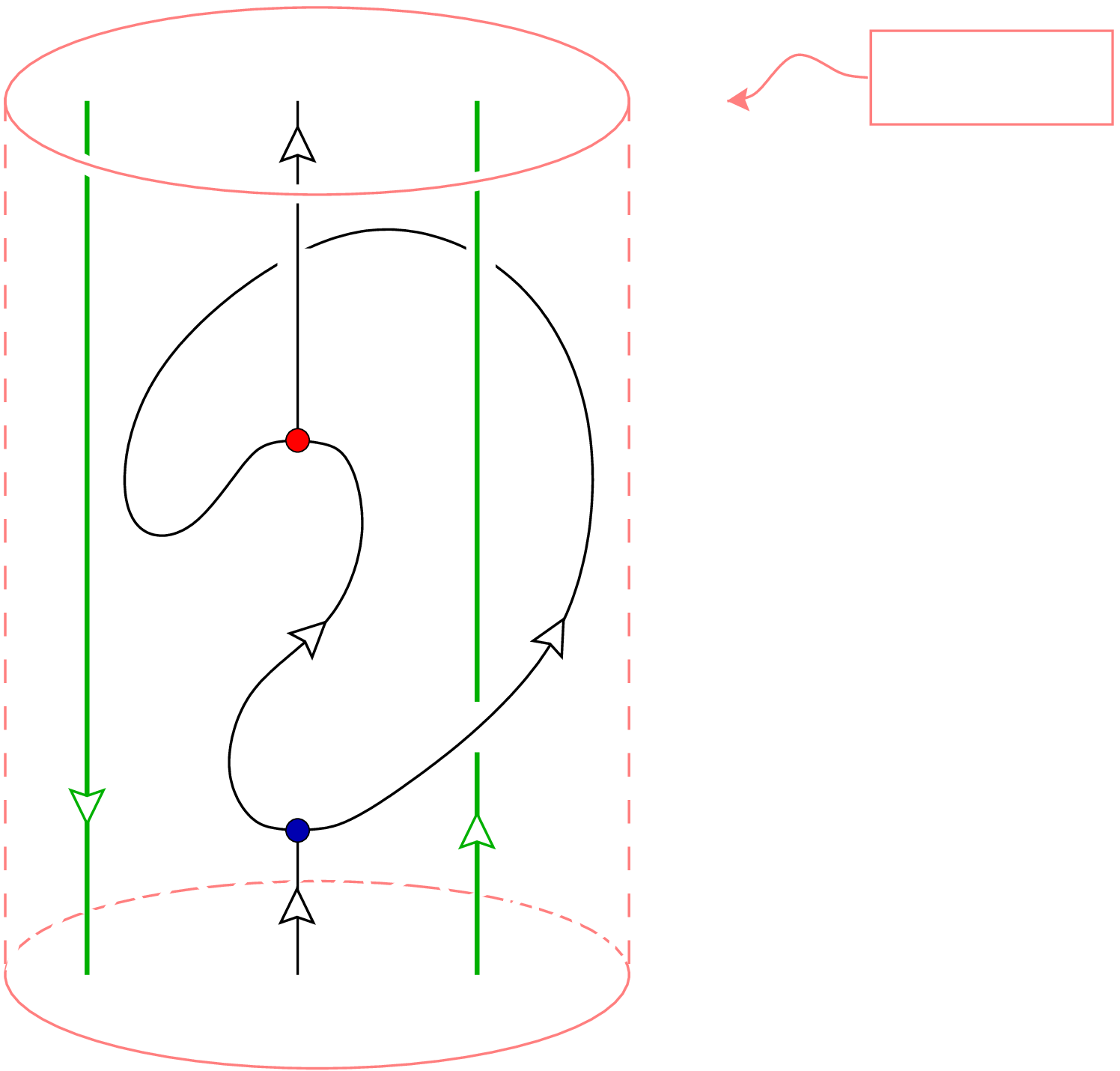}} \end{picture}}
  \put(-64.6,73.1) {$\tilde Z(A)_{ij}\;:=$}
  \put(14.8,14.9) {\sse$i$}
  \put(45.2,14.9) {\sse$A$}
  \put(71.5,14.9) {\sse$j$}
  \put(29.5,56) {\sse$A$}
  \put(79,83) {\sse$A$}
  \put(45,104) {\sse$A$}
  \put(128.5,141.6) {$S^2{\times}S^1$}
  \epicture16 \labl{tZij}
of the indicated ribbon graph in $S^2\,{\times}\,S^1$. (The $S^2$ factor is 
represented by the horizontal circles, while the $S^1$ factor is given
by the vertical direction, i.e.\ top and bottom of the figure are to be
identified.) In \Theorem 5.1 
of \cite{fuRs4} it is shown that for every
$i,j\iN\II$ the invariant $\tilde Z(A)_{ij}\,{\equiv}\,\tilde Z(A)_{U_i,U_j}$
is the trace of an idempotent and hence a
non-negative integer. The torus partition function of the conformal field
theory determined by $A$ is given by $Z(A)_{kl} \eq \tilde Z(A)_{\bar k l}$.
Furthermore (see \Proposition 5.3 
of \cite{fuRs4} and \cite{evan9}) for a modular
tensor category the integers $\tilde Z_{ij}$ defined by \erf{tZij} obey
$\tilde Z(A{\otimes}B)\eq\tilde Z(A)\, \tilde Z(B)$ as a matrix equation.

%%%%%%%%%%%%%%%%%%%%%%%%%%%%%%%%%%%%%%%%%%%
\newpage

\subsection{Left and right centers}\label{LRC}

In this subsection, $A$ denotes a symmetric special Frobenius algebra in a
ribbon category $\cC$. In a braided setting, there are two different notions
of center of an \alg, the left and right center; we establish some properties 
of the centers that we will need.  

The notion of left and right center was introduced in \cite{vazh} for
separable \alg s in abelian braided \tcs, and in \cite{ostr} for \alg s
in semisimple abelian ribbon \cats. Here we formulate it in terms of a
maximality property for idempotents; this makes it applicable even to 
non-Karoubian \cats, provided only that the particular idempotents
  %%  [pic~1b]
  \bea  \begin{picture}(300,86)(0,42)
  \put(0,0)   {\begin{picture}(0,0)(0,0)
              \scalebox{.38}{\includegraphics{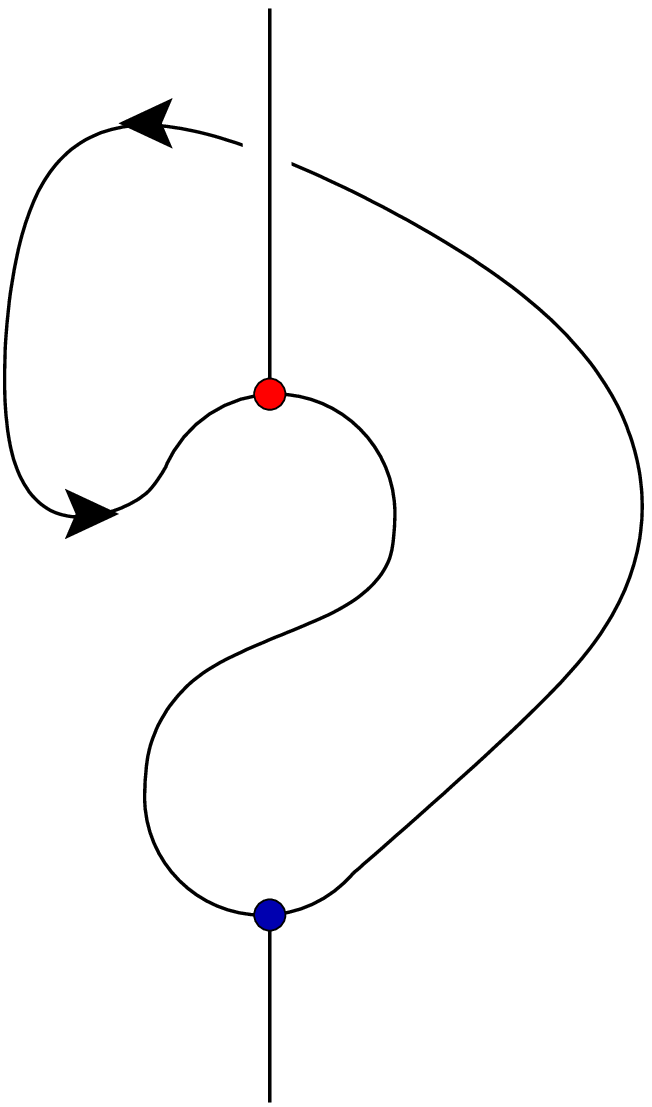}} \end{picture}}
  \put(-48.6,55.5) {$P^l_\AA\;:=$}
  \put(8.2,30)     {\sse$A$}
  \put(25.5,-8.8)  {\sse$A$}
  \put(26.1,125)   {\sse$A$}
  \put(73.2,67)    {\sse$A$}
  \put(120,55.5)   {and}
  \put(230,0)   {\begin{picture}(0,0)(0,0)
              \put(0,0)   {\begin{picture}(0,0)(0,0)
              \scalebox{.38}{\includegraphics{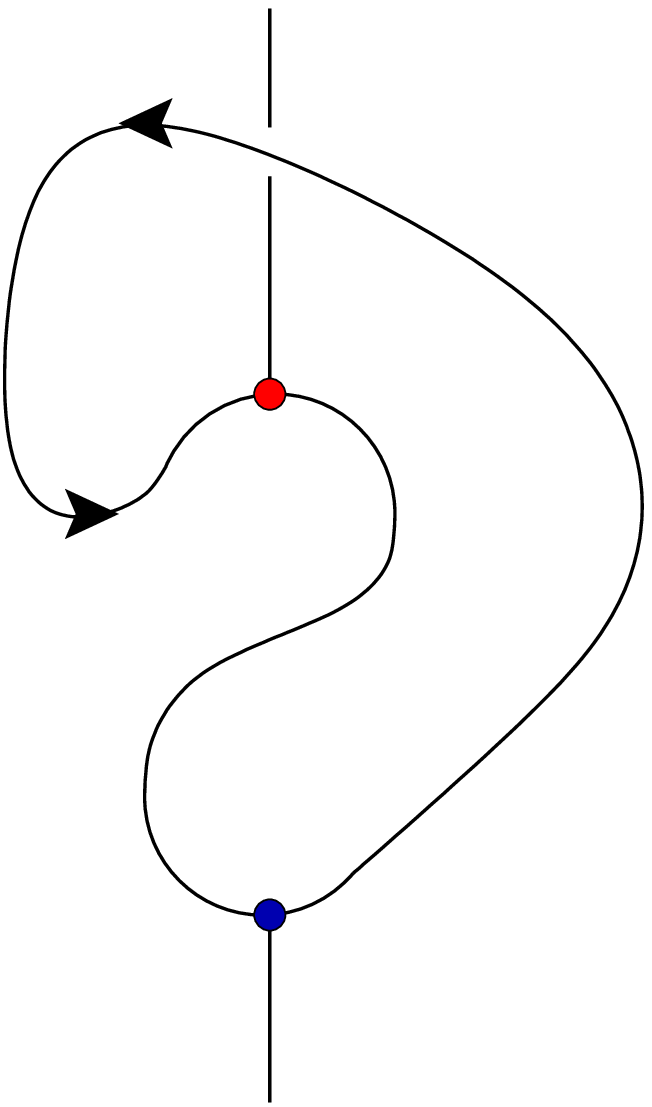}} \end{picture}}
  \put(-48.6,55.5) {$P^r_\AA\;:=$}
  \put(25.5,-8.8)  {\sse$A$}
  \put(26.1,125)   {\sse$A$}
  \put(73.2,67)    {\sse$A$}
  \put(8.2,30)     {\sse$A$} \end{picture}}
  \epicture29 \labl{PE-def}
in $\Hom(A,A)$ are split.
That $P^{l/r}_\AA$ are idempotents follows easily by using the various 
properties of $A$; for $P^l_\AA$ this is described in lemma 5.2 of \cite{fuRs4}
(setting $X\eq\one$ there), and for $P^r_\AA$ the argument is analogous.

These idempotents possess several nice properties. First, we have
\\[-2.3em]

\dtl{Lemma}{lem:C=[1]a}
The idempotents \erf{PE-def} satisfy the following relations.
\\[.2em]
(i)~~They trivialise the twist:
  \be
  \theta_\AA^{} \circ P_\AA^l = P_\AA^l
  \qquad {\rm and} \qquad
  \theta_\AA^{} \circ P_\AA^r = P_\AA^r \,.  \labl{eq:leave-theta}
(ii)~\,They are compatible with unit and counit: 
  \be  P_\AA^{l/r} \circ \eta = \eta \,, \qquad
  \eps \circ P_\AA^{l/r} = \eps \,.  \labl{oben1}
(iii)~When each of the three $A$-ribbons forming a product or coproduct is 
decorated with $P_\AA^{l/r}$, then any one of the three idempotents can be
omitted:
  \be  \bearll
  P_\AA^{l/r} \circ m \circ (P_\AA^{l/r} \oti P_\AA^{l/r})\!\!
  &= m \circ (P_\AA^{l/r} \oti P_\AA^{l/r})
  \\{}\\[-.8em]
  &=  P_\AA^{l/r} \circ m \circ (\id_\AA^{} \oti P_\AA^{l/r})
   =  P_\AA^{l/r} \circ m \circ (P_\AA^{l/r}  \oti \id_\AA^{}) \,,  
  \\{}\\[-.4em]
  (P_\AA^{l/r} \oti P_\AA^{l/r}) \circ \Delta  \circ P_\AA^{l/r} \!\!
  &= (P_\AA^{l/r} \oti P_\AA^{l/r}) \circ \Delta 
  \\{}\\[-.8em]
  &= (\id_\AA^{} \oti P_\AA^{l/r}) \circ \Delta  \circ P_\AA^{l/r}
   = (P_\AA^{l/r} \oti \id_\AA^{}) \circ \Delta  \circ P_\AA^{l/r} \,.
  \eear \labl{oben2}

\noindent
Proof:\\
(i)~\,The statement for $P_\AA^l$ follows by the moves
(the one for $P_\AA^r$ is derived analogously)
  %%  [pic~37a]
  \bea  \begin{picture}(290,69)(0,32)
  \put(0,0)  {\begin{picture}(0,0)(0,0)
             \scalebox{.38}{\includegraphics{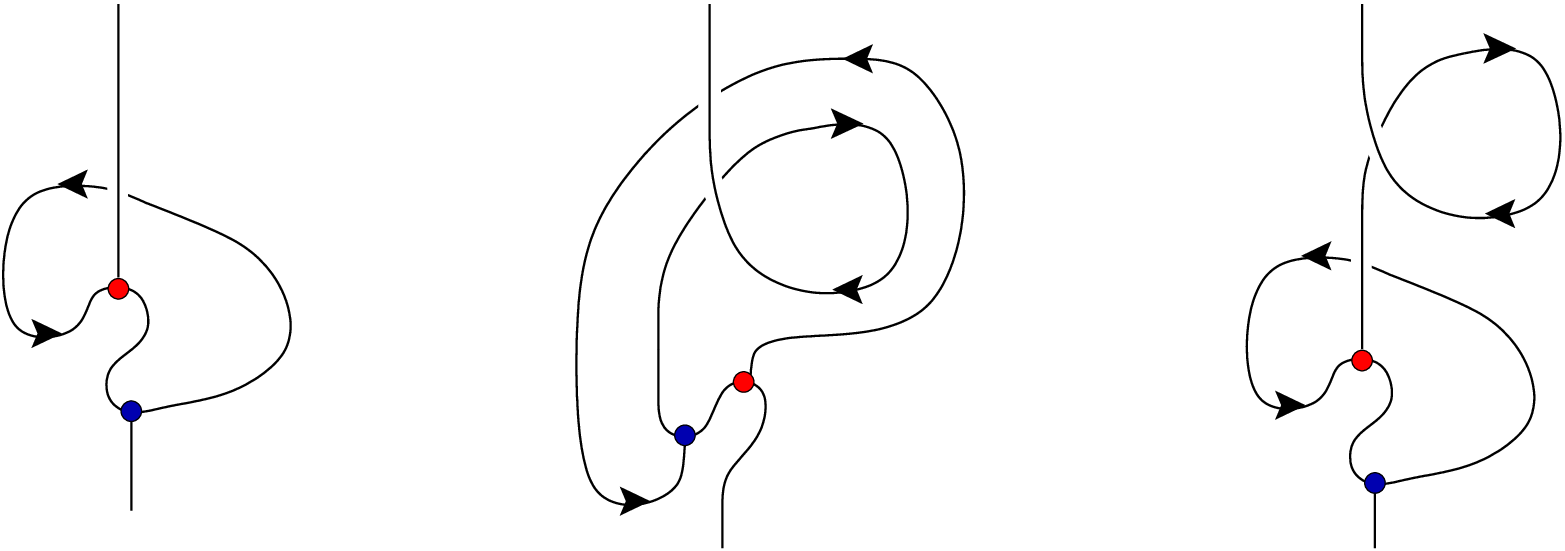}} \end{picture}}
  \put(16.2,96)    {\sse$A$}
  \put(17.5,-9.1)  {\sse$A$}
  \put(68,43.5)    {$=$}
  \put(116.1,96)   {\sse$A$}
  \put(117.2,-9.1) {\sse$A$}
  \put(180,43.5)   {$=$}
  \put(224.8,96)   {\sse$A$}
  \put(226.1,-9.1) {\sse$A$}
  \epicture17 \labl{pic37a}
To get the first equality one uses the Frobenius property and then suitably
drags the resulting coproduct along part of the $A$-ribbon. A further 
deformation and application of the Frobenius property then results in the 
second equality.
\\[.3em]
The statements in (ii)~and (iii)~are just special cases of the statements of 
\Lemma \ref{lem:remove-P} below -- they follow from those by setting $B\eq\one$.
\qed

Next recall the \Definition \ref{maxidm} 
of maximal idempotent.
It turns out that $P_\AA^l$ and $P_\AA^r$ can be characterised as 
being maximal in a subset of ${\rm Idem}(A)$ that is defined by a relation 
involving the braiding and the product of $A$:
\\[-2.3em]

\dtl{Lemma}{centralidem}
The subset $H_l\,{\subseteq}\,{\rm Idem}(A)$ consisting of those
idempotents $p$ in $\End(A)$ that satisfy
  %% [pic~8b]
  \bea  \begin{picture}(350,100)(0,1)
  \put(-5,0)  {\begin{picture}(0,0)(0,0)
              \scalebox{.38}{\includegraphics{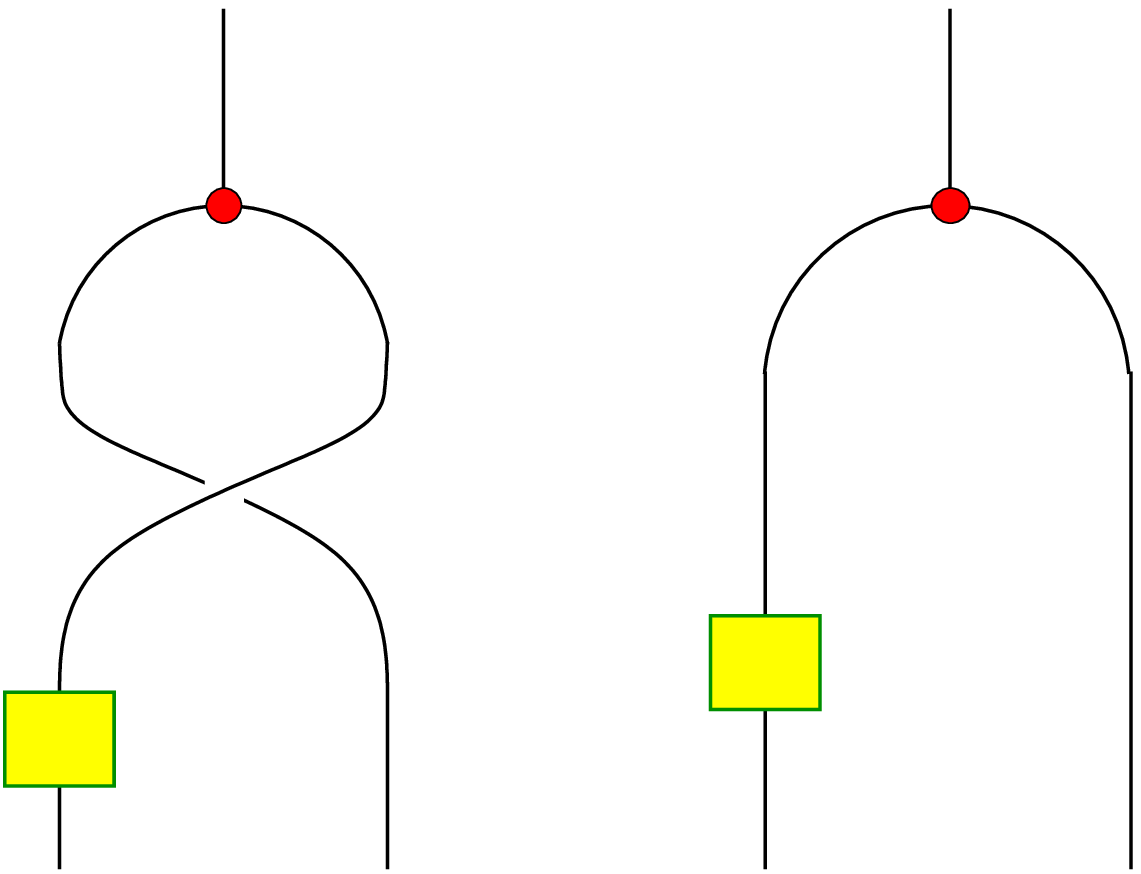}} \end{picture}}
  \put(-2.5,-8.9)  {\sse$A$}
  \put(-.8,14)     {\sse$p$}
  \put(15.5,99)    {\sse$A$}
  \put(33.2,-8.9)  {\sse$A$}
  \put(53,42.7)    {$=$}
  \put(75.5,-8.9)  {\sse$A$}
  \put(77.2,21.9)  {\sse$p$}
  \put(95.5,99)    {\sse$A$}
  \put(114.8,-8.9) {\sse$A$}
  \put(160,42.7)   {and}
  \put(215,0)  {\begin{picture}(0,0)(0,0)
              \scalebox{.38}{\includegraphics{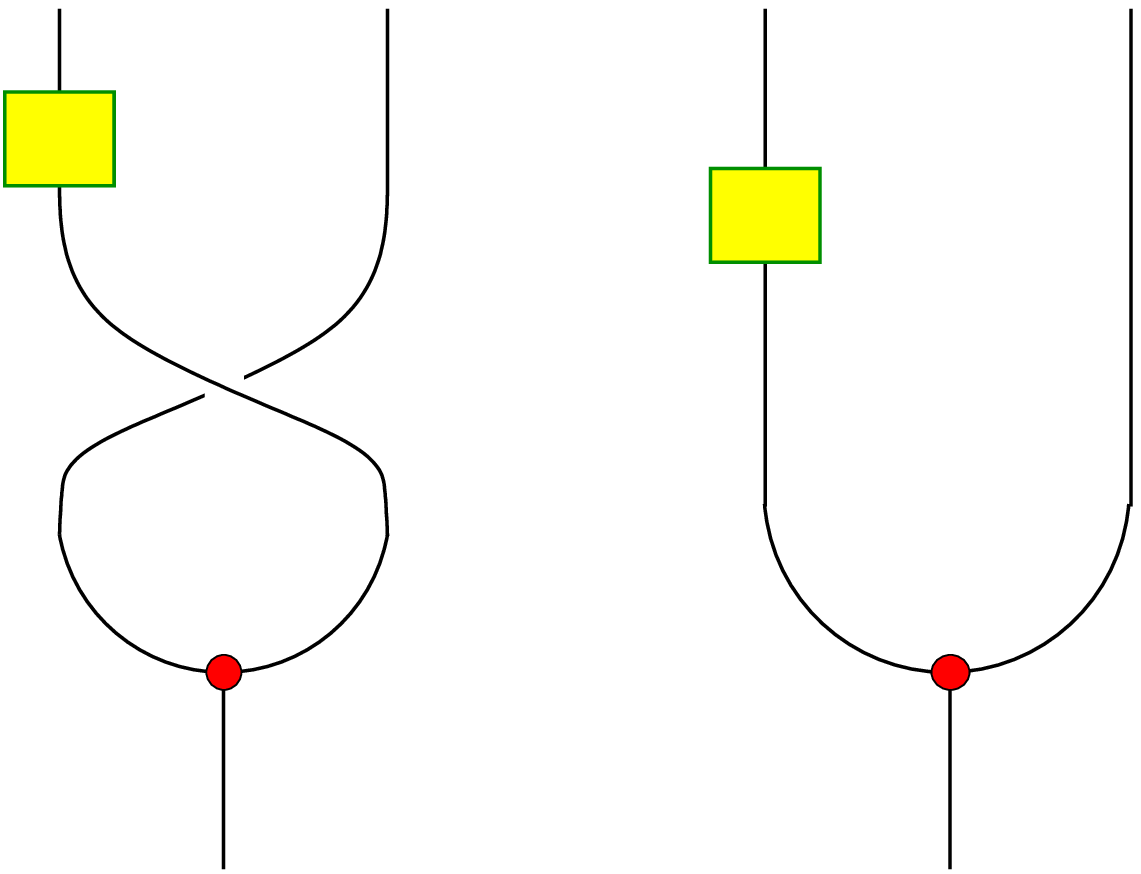}} \end{picture}}
  \put(218.1,99)   {\sse$A$}
  \put(219.4,79.8) {\sse$p$}
  \put(235.5,-8.9) {\sse$A$}
  \put(254.3,99)   {\sse$A$}
  \put(275,42.7)   {$=$}
  \put(295.8,99)   {\sse$A$}
  \put(297.1,71.1) {\sse$p$}
  \put(315.8,-8.9) {\sse$A$}
  \put(335.4,99)   {\sse$A$}
  \epicture-8 \labl{eq:Hl-def}
contains a maximal idempotent, and this is given by $P_\AA^l$ in \erf{PE-def}.
\\[.3em]
Analogously, the subset $H_r$ of those idempotents $p$ in $\End(A)$ satisfying
  %% [pic~8c]
  \bea  \begin{picture}(350,104)(0,1)
  \put(1,0)  {\begin{picture}(0,0)(0,0)
              \scalebox{.38}{\includegraphics{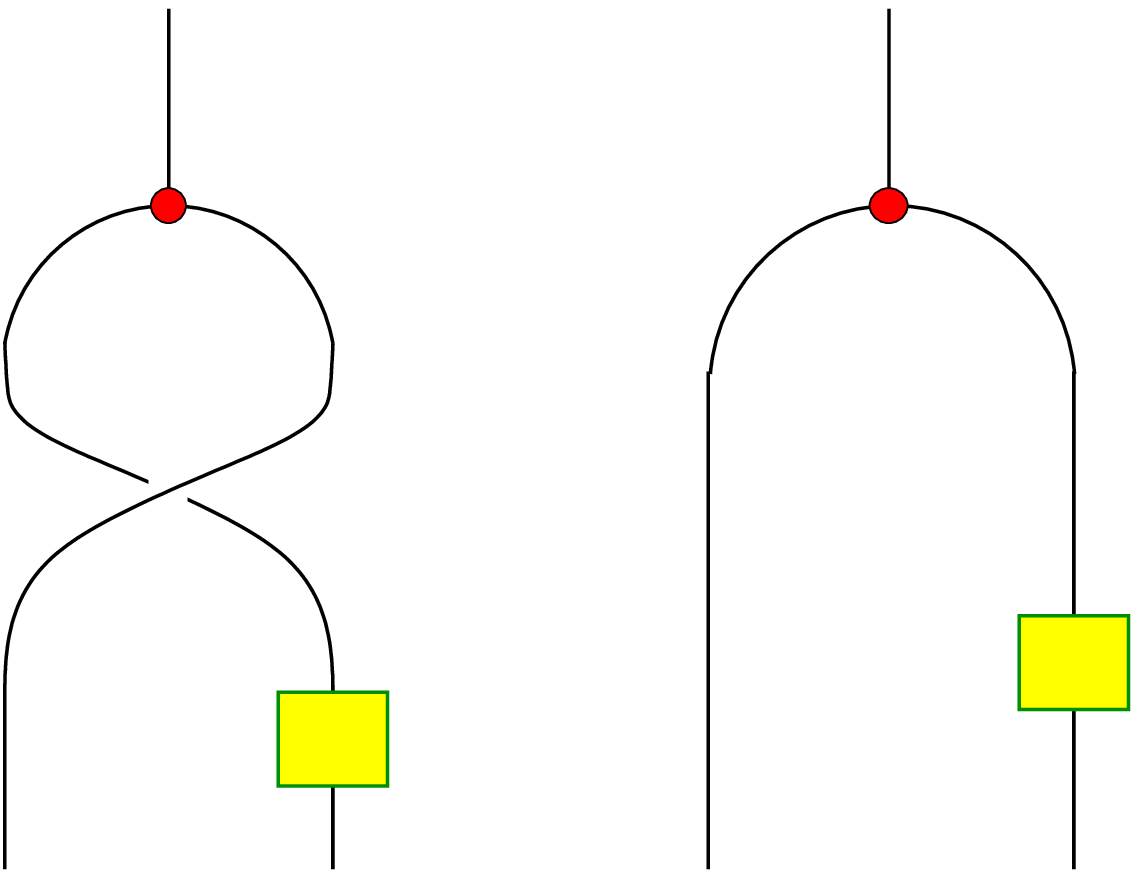}} \end{picture}}
  \put(-2.5,-8.9)  {\sse$A$}
  \put(15.5,99)    {\sse$A$}
  \put(33.2,-8.9)  {\sse$A$}
  \put(34.9,13.3)  {\sse$p$} 
  \put(53,42.7)    {$=$}
  \put(75.5,-8.9)  {\sse$A$}
  \put(95.5,99)    {\sse$A$} 
  \put(114.8,-8.9) {\sse$A$}
  \put(116.6,21.9) {\sse$p$}
  \put(164,42.7)   {and}
  \put(222,0)  {\begin{picture}(0,0)(0,0)
              \scalebox{.38}{\includegraphics{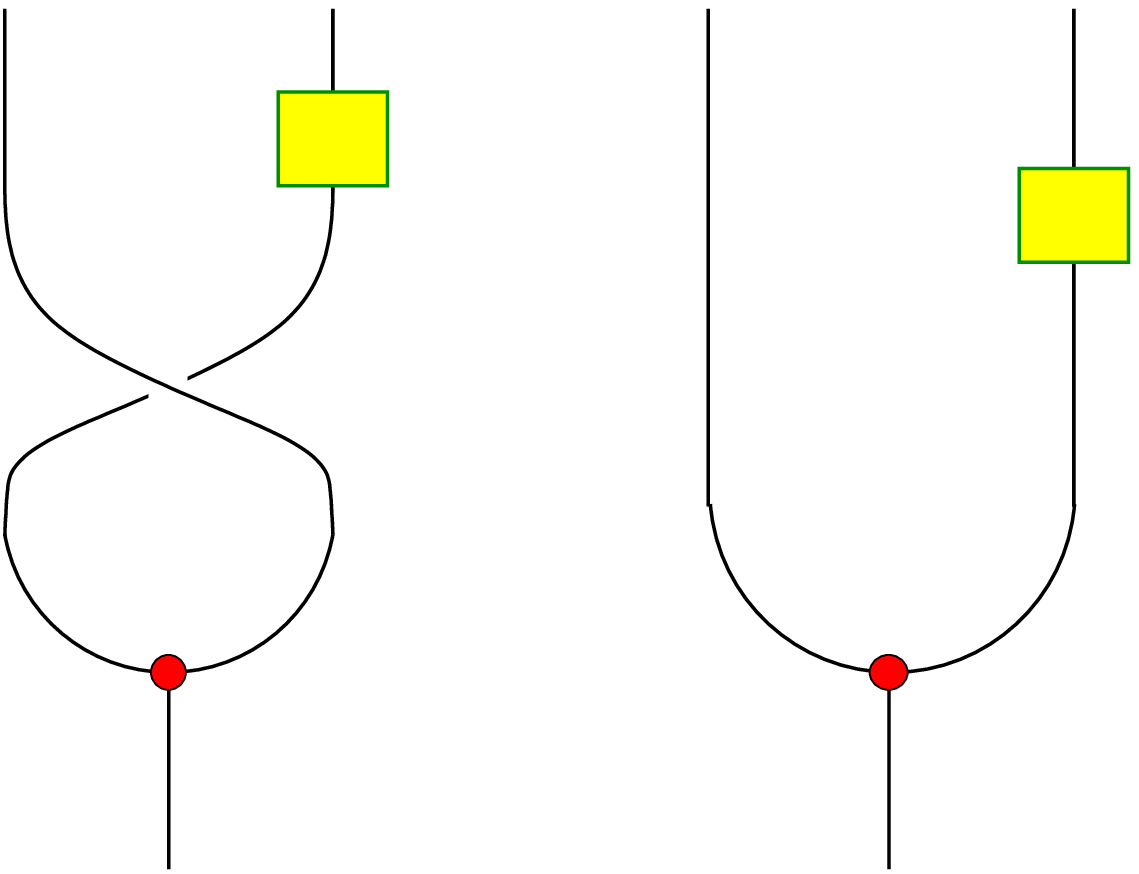}} \end{picture}}
  \put(218.1,99)   {\sse$A$}
  \put(235.5,-8.9) {\sse$A$}
  \put(254.3,99)   {\sse$A$}
  \put(255.9,79.8) {\sse$p$}
  \put(275,42.7)   {$=$}
  \put(295.8,99)   {\sse$A$}
  \put(315.8,-8.9) {\sse$A$}
  \put(335.4,99)   {\sse$A$}
  \put(337.4,70.9) {\sse$p$}
  \epicture-7 \labl{eq:Hr-def}
contains a maximal idempotent, given by $P_\AA^r$ in \erf{PE-def}.

\medskip\noindent
Proof:\\
We prove the statement for $P_\AA^l$; the statement for $P_\AA^r$
follows analogously.
\\
Consider the transformations
  %% [pic~52c]
  \begin{eqnarray}  \begin{picture}(420,172)(10,0)
  \put(0,0)  {\begin{picture}(0,0)(0,0)
              \scalebox{.38}{\includegraphics{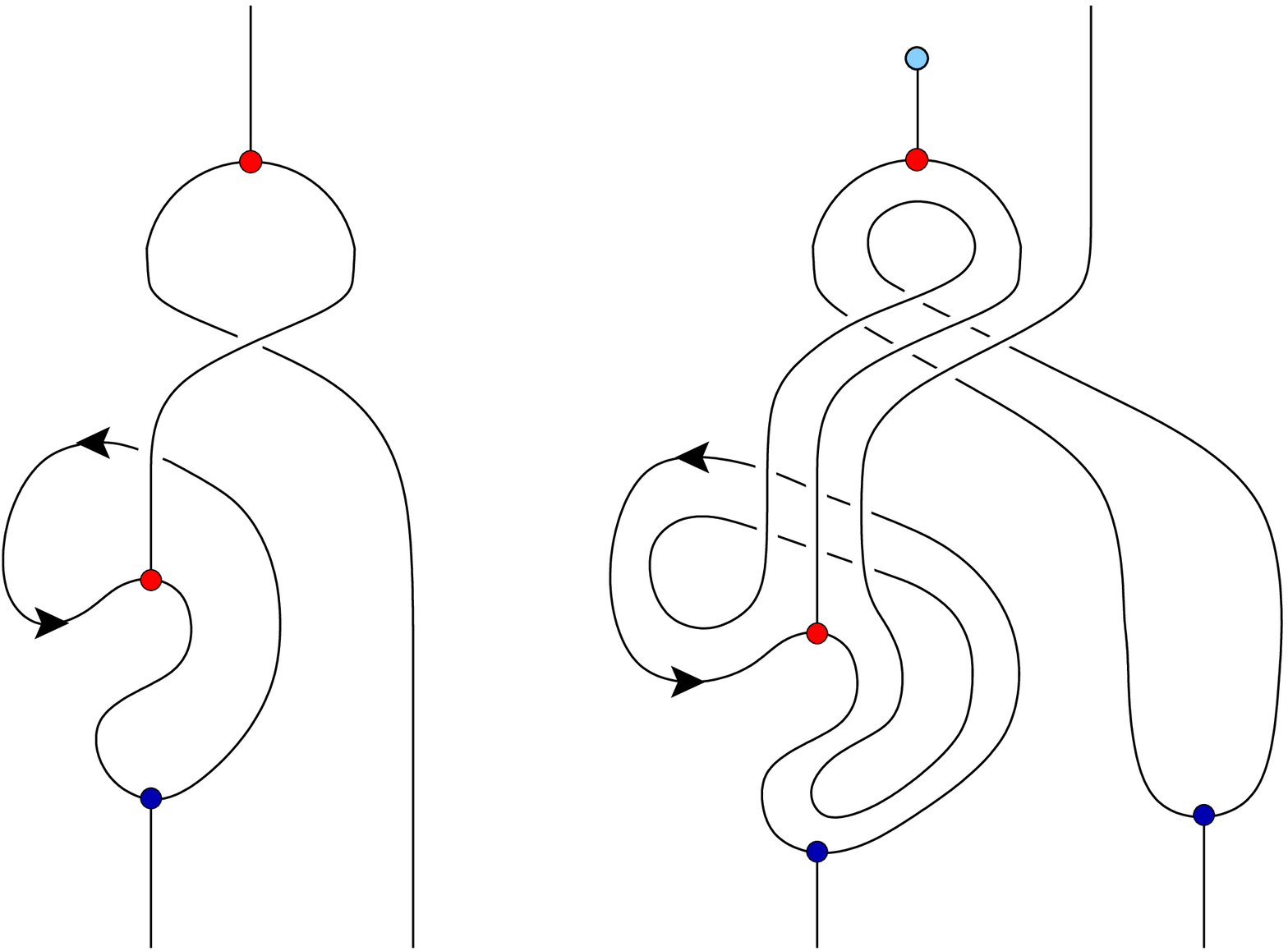}} \end{picture}}
  \put(21.5,-9.5)   {\sse$A$}
  \put(40.5,168.2)  {\sse$A$}
  \put(67.7,-9.5)   {\sse$A$}
  \put(85,81)       {$=$}
  \put(136.5,-9.5)  {\sse$A$}
  \put(185.5,168.2) {\sse$A$}
  \put(204.9,-9.5)  {\sse$A$}
  \put(238,81)      {$=$}
  \put(283.5,-9.5)  {\sse$A$}
  \put(334.5,-9.5)  {\sse$A$}
  \put(349.1,168.2) {\sse$A$}
  \put(366,81)      {$=$}
  \put(391.3,113)   {\sse$\theta_\AA^{-1}$}
  \put(408.1,-9.5)  {\sse$A$}
  \put(426.1,168.2) {\sse$A$}
  \put(454.1,-9.5)  {\sse$A$}
  \end{picture} \nonumber\\[1.1em]{} \label{pic52c}
  \\[-1.9em]{}\nonumber\end{eqnarray}
In the first step the identity $\id_\AA\eq(\eps{\otimes}\id_\AA)\cir\Delta$
is inserted in the top $A$-ribbon, and afterwards the coproduct $\Delta$
introduced this way is moved along a path that can easily be read off the
second picture, using the Frobenius and coassociativity properties of $A$. 
The next step is just a deformation of the outgoing $A$-ribbon. The third 
step uses that $A$ is symmetric, and then the Frobenius property.
\\
Afterwards, according to \Lemma \ref{lem:C=[1]a} 
the twist $\theta_\AA^{-1}$ 
can be left out. This establishes that $P_\AA^l$ satisfies the first of the
equalities \erf{eq:Hl-def}. That $P_\AA^l$ satisfies the second of those
equalities as well is deduced similarly, starting with an insertion of the 
identity $\id_\AA\eq m\cir(\eta{\otimes}\id_\AA)$ on the incoming $A$-ribbon.
Together it follows that $P_\AA^l\iN H_l$.
\\[.3em]
Furthermore, composing the first equality in \erf{eq:Hl-def} from the bottom
with $\id_A\oti b_A$ and from the top with $(\id_A\oti\tilde d_A)\cir
(\Delta\oti\id_{A^\vee})$ shows, upon using the symmetry, specialness and
Frobenius properties, that $p\eq P_\AA^l\cir p$ for $p\iN H_l$.
Similar manipulations of the second equality in \erf{eq:Hl-def} show that
also $p\eq p\cir P_\AA^l$ for $p\iN H_l$. Thus $P_\AA^l$ is maximal in $H_l$.
\qed

\dtl{Definition}{defLRC}
(i)~\,We call the morphism $P_\AA^l$ defined in \erf{PE-def} the 
{\em left central idempotent\/} of the \ssFA\ 
$A$, and $P_\AA^r$ the {\em right central idempotent\/} of $A$.
\\[.3em]
(ii)~The {\em left center\/} $C_l(A)$ of the \ssFA\ $A$
is the maximal retract of $A$ with respect to $H_l$.
\\
The {\em right center\/} $C_r(A)$ of $A$
is the maximal retract of $A$ with respect to $H_r$.

\medskip

According to \Lemma \ref{centralidem} 
the left and right central idempotents
$P_\AA^{l/r}$ are the maximal idempotents of the subsets $H_l$ and $H_r$ of
${\rm Idem}(A)$, respectively. It follows that the left (right) center of 
$A$ exists iff the left (right) central idempotent is split; if this is 
the case, then by corollary \ref{cor:maxret}, $C_{l/r}(A)$ is unique.
In the sequel we will often use the short-hand notation $C_l$ 
and $C_r$ for $C_l(A)$ and $C_r(A)$, respectively.
Also note that the definition of the centers involves both the \alg\ and 
the co\alg\ structure of $A$; in place of the term center one might 
therefore also use the term `Frobenius center'.

\dtl{Lemma}{lem:Pe=e}
Any retract $(S,e,r)$ of a \ssFA\ $A$ that obeys $\,m\cir c_{A,A}\cir 
(e^{}_{S\prec A}
    $\linebreak[0]$%
\oti\id_\AA) \eq m\cir (e^{}_{S\prec A}\oti\id_\AA)$
and $(r^{}_{\!\!A\succ S}\oti\id_\AA)\cir c_{A,A}^{-1}\cir\Delta
\eq (r^{}_{\!\!A\succ S}\oti\id_\AA)\cir\Delta$, i.e.
  %% [pic~8]
  \bea  \begin{picture}(350,113)(0,-3)
  \put(-8,0)  {\begin{picture}(0,0)(0,0)
              \scalebox{.38}{\includegraphics{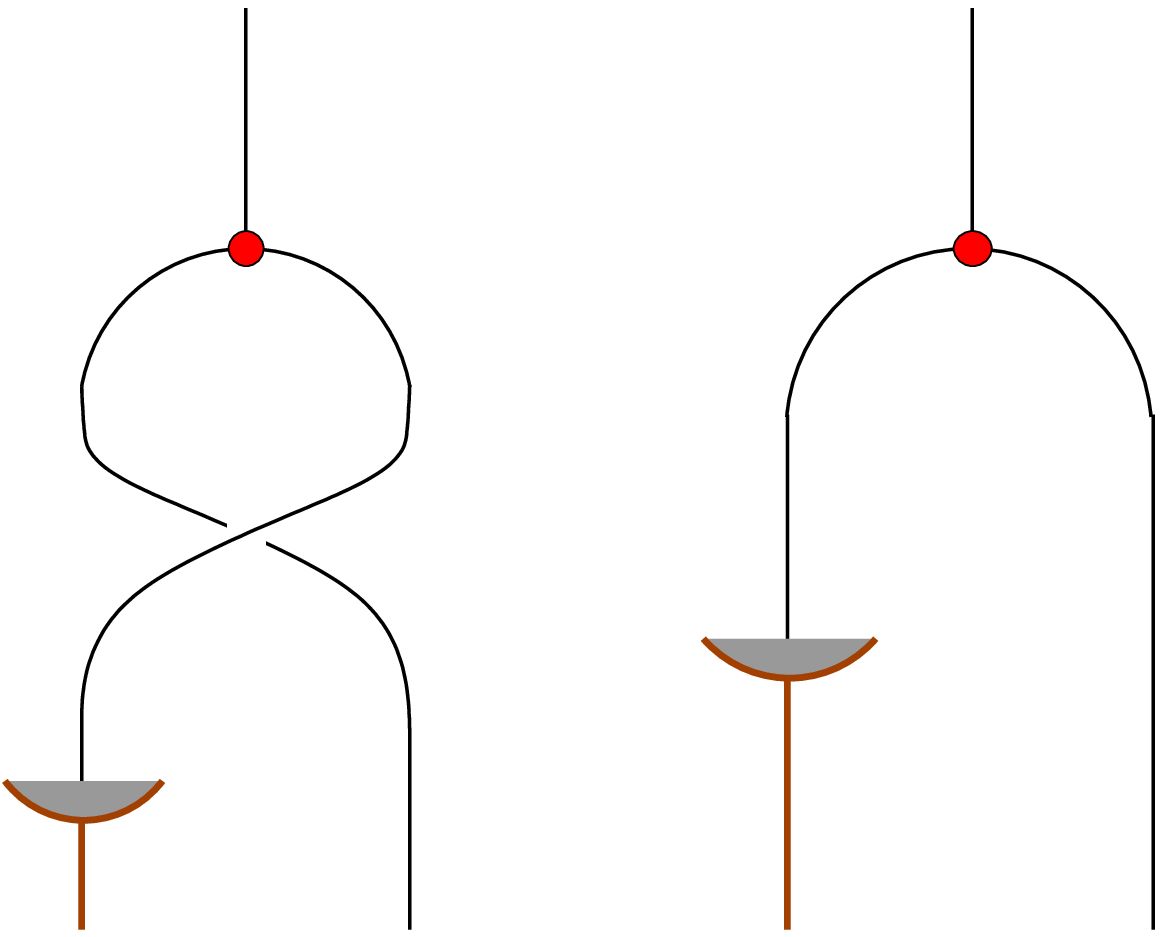}} \end{picture}}
  \put(231,0)  {\begin{picture}(0,0)(0,0)
              \scalebox{.38}{\includegraphics{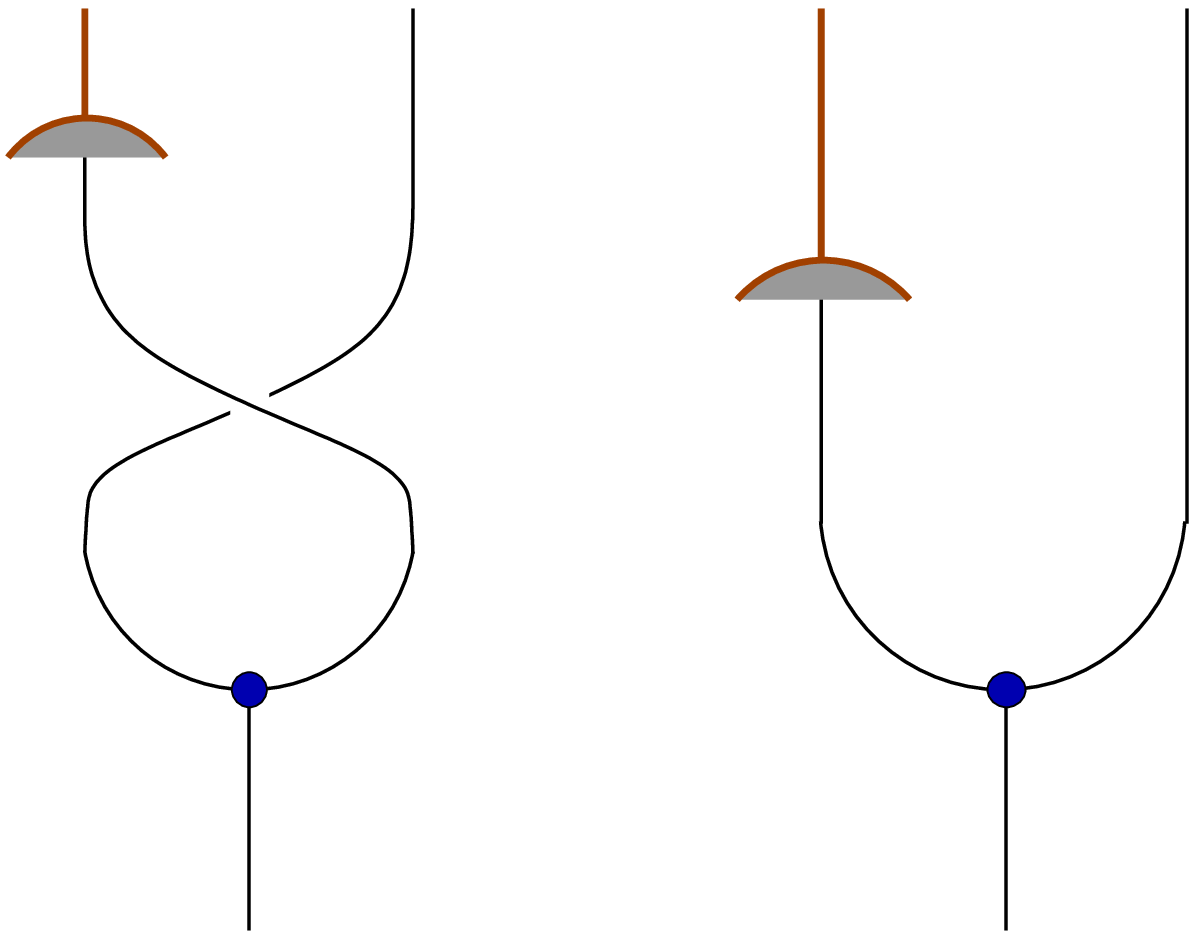}} \end{picture}} 
  \put(-3.1,-8.7)  {\sse$S$} 
  \put(15.5,105.5) {\sse$A$} 
  \put(33.2,-8.7)  {\sse$A$} 
  \put(53,37.7)    {$=$} 
  \put(75.1,-8.7)  {\sse$S$} 
  \put(95.5,105.5) {\sse$A$} 
  \put(114.8,-8.7) {\sse$A$}
  \put(170,37.7)   {and} 
  \put(238.3,105.5){\sse$S$}
  \put(254.3,-8.7) {\sse$A$}
  \put(272.3,105.5){\sse$A$}
  \put(293,37.7)   {$=$}
  \put(320.3,105.5){\sse$S$}
  \put(336.9,-8.7) {\sse$A$}
  \put(357.6,105.5){\sse$A$} 
  \epicture-5 \labl{S-Cl-prop} 
also satisfies
  \be  P_\AA^l \circ e^{}_{S\prec A} = e^{}_{S\prec A}
  \qquad{\rm and}\qquad
  r^{}_{\!\!A\succ S} \circ P_\AA^l = r^{}_{\!\!A\succ S} \,. \labl{Pe=e}
Similarly we have
  \be  \left. \begin{array}r
  m\cir c_{A,A}^{}\cir(\id_\AA\oti e^{}_{S\prec A})
    = m\cir (\id_\AA\oti e^{}_{S\prec A})
  \\{}\\[-.75em]
  (\id_\AA\oti r^{}_{\!\!A\succ S})\cir c_{A,A}^{-1}\cir\Delta
    = (\id_\AA\oti r^{}_{\!\!A\succ S})\cir\Delta
  \eear\!\right\} \, \Rightarrow \, \left\{\! \bearl 
  P_\AA^r \circ e^{}_{S\prec A} = e^{}_{S\prec A} 
  \\{}\\[-.75em]
  r^{}_{\!A\succ S} \circ P_\AA^r = r^{}_{\!A\succ S} \,.
  \end{array}\right. \labl{Pe=e2}

\medskip\noindent
Proof:\\
Composing \erf{S-Cl-prop} from the bottom with $r\oti\id_A$ shows that the 
idempotent $p\eq e^{}_{S\prec A}\cir r^{}_{\!A\succ S}$ satisfies the first 
of the equalities \erf{eq:Hl-def}. Analogously one shows that $p$ also
obeys the second of those equalities, and hence it is contained in $H_l$. 
Thus the relations \erf{Pe=e} are implied by \erf{S-Cl-prop} together with the 
maximality property of $P_\AA^l$. 
\\
The implication \erf{Pe=e2} is derived analogously.
\qed 

\dtl{Lemma}{Ctwist-etc} 
The left and right center of a \ssFA\ $A$ have trivial twist:
  \be  \theta_{C_l(A)} = \id_{C_l(A)} \,, \qquad
  \theta_{C_r(A)} = \id_{C_r(A)} \,.  \labl{eq:Ctwist}

\medskip\noindent
Proof:
\\[.3em]
The statement follows immediately from the relations \erf{eq:leave-theta}.
(Conversely, \erf{eq:leave-theta} follows from \erf{eq:Ctwist}
by functoriality of the twist.)
\qed

\dtl{Remark}{only-i}
As a consequence of \Lemma \ref{lem:Pe=e} 
the centers obey
  %% [pic~8]
  \bea  \begin{picture}(350,106)(0,1)
  \put(-8,0)  {\begin{picture}(0,0)(0,0)
              \scalebox{.38}{\includegraphics{Cl-def.eps}} \end{picture}}
  \put(-9.5,-8.7)  {\sse$C_l(A)$}
  \put(15.5,105.5) {\sse$A$}
  \put(33.2,-8.7)  {\sse$A$}
  \put(53,37.7)    {$=$}
  \put(69.5,-8.7)  {\sse$C_l(A)$}
  \put(95.5,105.5) {\sse$A$}
  \put(114.8,-8.7) {\sse$A$}
  \put(160,37.7)   {and}
  \put(220,0)  {\begin{picture}(0,0)(0,0)
              \scalebox{.38}{\includegraphics{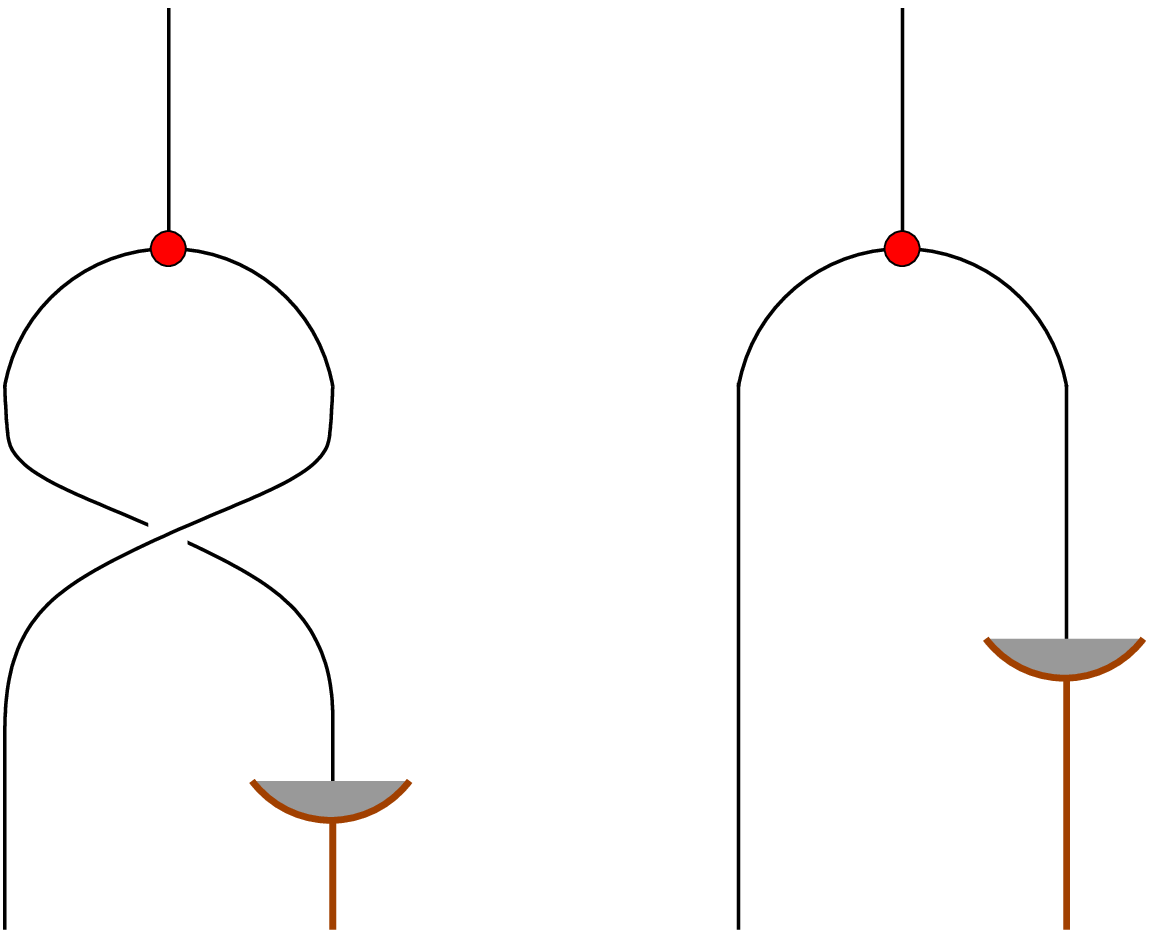}} \end{picture}}
  \put(216.1,-8.7) {\sse$A$}
  \put(235.5,105.5){\sse$A$}
  \put(246.6,-8.7) {\sse$C_r(A)$}
  \put(277,37.7)   {$=$}
  \put(296.6,-8.7) {\sse$A$}
  \put(316.2,105.5){\sse$A$}
  \put(327.1,-8.7) {\sse$C_r(A)$}
  \epicture-1 \labl{Cl-Cr-defprop}
respectively, as well as 
  %% [pic~18]
  \bea  \begin{picture}(380,376)(5,0)
  \put(0,0)  {\begin{picture}(0,0)(0,0)
              \scalebox{.38}{\includegraphics{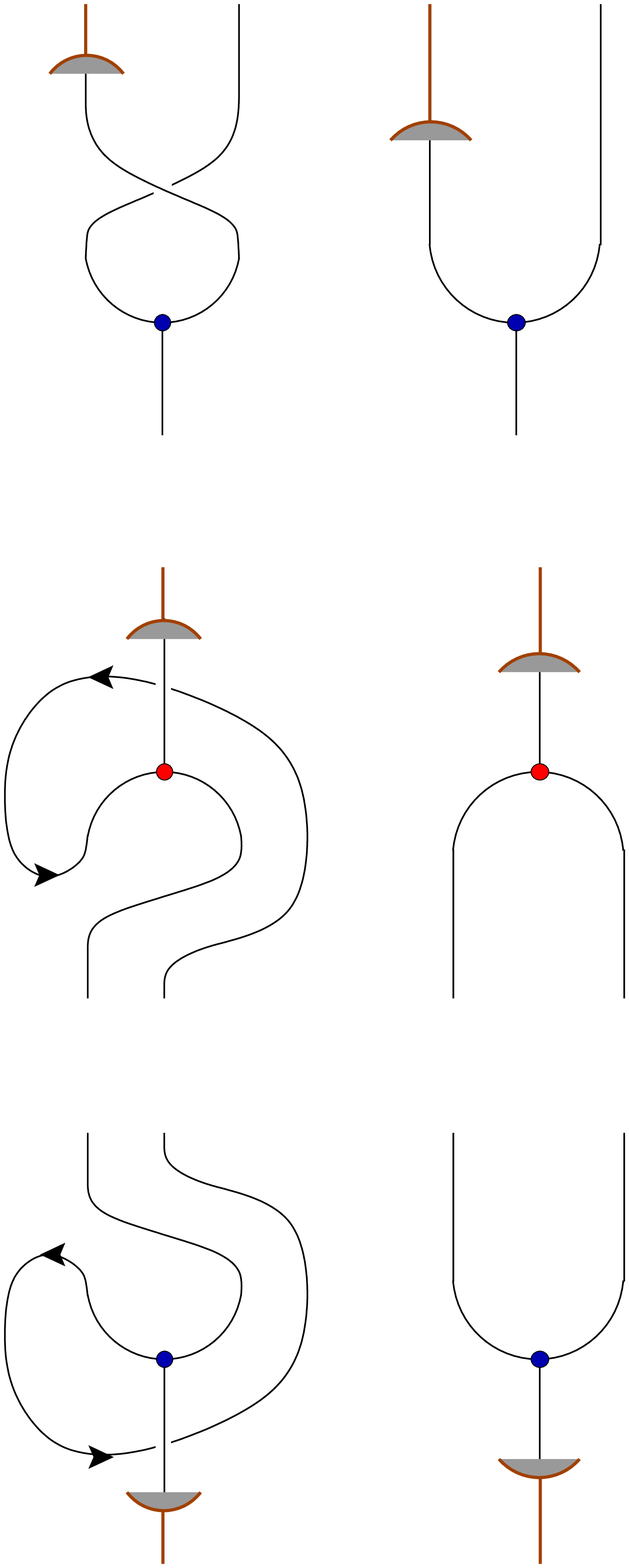}} \end{picture}}
  \put(15.3,371.2)  {\sse$C_l$}
  \put(34.3,256.2)  {\sse$A$}
  \put(52.3,371.2)  {\sse$A$}
  \put(73,309)      {$=$}
  \put(97.3,371.2)  {\sse$C_l$}
  \put(116.9,256.2) {\sse$A$}
  \put(137.6,371.2) {\sse$A$}
  \put(17.3,123.5)  {\sse$A$}
  \put(33.9,238.5)  {\sse$C_l$}
  \put(34.5,123.5)  {\sse$A$}
  \put(84,176)      {$=$}
  \put(102.1,123.5) {\sse$A$}
  \put(122.2,238.5) {\sse$C_l$}
  \put(142.3,123.5) {\sse$A$}
  \put(17.3,105.8)  {\sse$A$}
  \put(33.3,-9.2)   {\sse$C_l$}
  \put(34.7,105.8)  {\sse$A$}
  \put(84,52)       {$=$}
  \put(102.7,105.8) {\sse$A$}
  \put(121.6,-9.2)  {\sse$C_l$}
  \put(142.9,105.8) {\sse$A$}
  \put(220,0)  {\begin{picture}(0,0)(0,0)
              \scalebox{.38}{\includegraphics{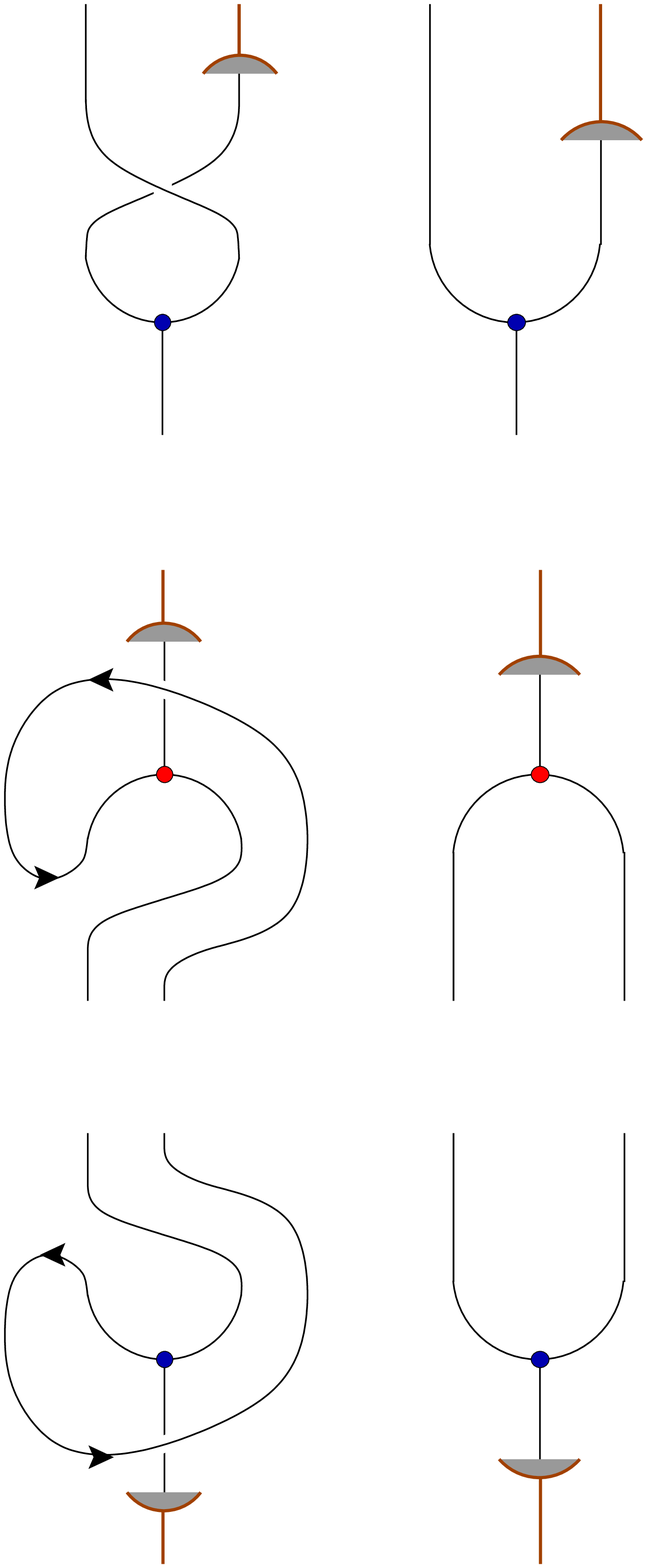}} \end{picture}}
  \put(236.1,371.2) {\sse$A$}
  \put(254.3,256.2) {\sse$A$}
  \put(271.6,371.2) {\sse$C_r$}
  \put(293,309)     {$=$}
  \put(317.4,371.2) {\sse$A$}
  \put(336.9,256.2) {\sse$A$}
  \put(356.4,371.2) {\sse$C_r$}
  \put(237.3,123.5) {\sse$A$}
  \put(253.9,238.5) {\sse$C_r$}
  \put(254.5,123.5) {\sse$A$}
  \put(304,176)     {$=$}
  \put(322.1,123.5) {\sse$A$}
  \put(342.2,238.5) {\sse$C_r$}
  \put(362.3,123.5) {\sse$A$}
  \put(237.3,105.8) {\sse$A$}
  \put(253.3,-9.2)  {\sse$C_r$}
  \put(254.7,105.8) {\sse$A$}
  \put(304,52)      {$=$}
  \put(322.7,105.8) {\sse$A$}
  \put(341.6,-9.2)  {\sse$C_r$}
  \put(362.9,105.8) {\sse$A$}
  \epicture-5 \labl{eq:remove-braiding}
together with the `mirrored' versions of these eight identities that are 
obtained by reflecting all the figures about a vertical axis. For instance, 
to establish the last of the equalities \erf{eq:remove-braiding}, one can 
start with the Frobenius relation \erf{1f} composed with 
$e_{C_r\prec A}{\otimes}\eta$, then apply the mirrored version of 
\erf{Cl-Cr-defprop} to the resulting product, and finally use the symmetry and 
Frobenius properties to remove the unit that was introduced in the first step.

\medskip

In the \Definition \erf{ZUV} of the numbers $\tilde Z(A)_{U,V}$ the two 
different $\alpha$-inductions were used. The corresponding morphism spaces 
for $\alpha$-inductions of the same type turn out to be related to the 
centers of $A$. To see this we first need
\\[-2.3em]

\dtl{Lemma}{++lemma}
Let $A$ be a symmetric Frobenius algebra in a ribbon \cat\ $\cC$ and 
  %% do not need that $A$ is special.
$U,V\iN\Objc$. Then for any $\varphi^+\iN\HomAA(\alpha_\AA^+(U),\alpha_\AA^+
(V))$ and any $\varphi^-\iN\HomAA(\alpha_\AA^-(U),\alpha_\AA^-(V))$ we have
  \be  (P_\AA^l\oti\id_V) \cir \varphi^+ = \varphi^+ \cir (P_\AA^l\oti\id_U)
  \qquad{\rm and}\qquad
  (P_\AA^r\oti\id_V) \cir \varphi^- = \varphi^- \cir (P_\AA^r\oti\id_U)
  \,.  \ee

\medskip\noindent
Proof:\\
Using functoriality of the braiding and the fact that $A$ is symmetric
Frobenius one easily rewrites the morphism $P_\AA^{l/r}\oti\id_U$ in such 
a way that it involves the left and right action of $A$ on the bimodule 
$\alpha_\AA^\pm(U)$. Since $\varphi^\pm$ is a morphism of bimodules, these
actions of $A$ can thus be passed through $\varphi^\pm$ (using again
also functoriality of the braiding). Afterwards one follows the steps used 
in rewriting $P_\AA^{l/r}\oti\id_U$ in reverse order, resulting in
$P_\AA^{l/r}\oti\id_V$.
\qed

\medskip

Using this lemma, we deduce the following relation with the centers of $A$.
\\[-2.3em]

\dtl{Proposition}{lem:[U]A-as-obj}
For any \ssFA\ $A$ in a ribbon category $\cC$ and any two objects
$U,V\iN\Objc$ there are natural bijections
  \be  \Hom(C_l(A)\Oti U,V) \,\cong\, \HomAA(\alpha_\AA^+(U),\alpha_\AA^+(V))
  \cong \Hom(U,C_l(A)\Oti V) \,\,\, \labl{neulab1}
and
  \be  \Hom(C_r(A)\Oti U,V) \,\cong\, \HomAA(\alpha_\AA^-(U),\alpha_\AA^-(V))
  \cong \Hom(U,C_r(A)\Oti V) \,.    \labl{neulab2}

\medskip\noindent
Proof:\\
We prove the first bijection in \erf{neulab1}, the proof of the others 
being analogous. 
\\
Let us abbreviate $C_l(A)\eq C$ as well as $e_{C_l(A)\prec A}\eq e$ and 
$r_{\AA\succ C_l(A)} \eq r$. Consider the mappings $\,\Phi{:}\; 
\Hom(C\Oti U,V) \,{\to}\, \Hom(A\Oti U,A\Oti V)$ and $\,\Psi{:}\; \HomAA
    $\linebreak[0]$%
(\alpha_\AA^+(U),\alpha_\AA^+(V)) \,{\to}\, \Hom(C\Oti U,V)$ defined by
  \be  \bearl
  \Phi(\varphi) := (\id_\AA\oti \varphi) \circ
  \left[\, \left( (\id_\AA\oti r)\cir\Delta) \right) {\otimes}\, \id_U \right]
  \,, \\{}\\[-.8em]
  \Psi(\psi) := (\eps\oti\id_V) \circ \psi \circ (e\oti\id_U) \,.
  \eear \labl{PhiPsi}
It is not difficult to check that for any $\varphi\iN\Hom(C\Oti U,V)$,
$\Phi(\varphi)$ intertwines both the left and the right action of $A$ on
$\alpha_\AA^+$-induced bimodules, and hence the image of $\Phi$ lies
actually in $\HomAA(\alpha_\AA^+(U),\alpha_\AA^+(V))$.
\\
Furthermore, $\Phi$ and $\Psi$ are two-sided inverses of each other.
That $\Psi\cir\Phi(\varphi)\eq\varphi$ is seen by just applying the defining
property of the counit and then using $r\cir e\eq\id_C$, while to establish 
$\Phi\cir\Psi(\psi)\eq\psi$, one invokes \Lemma \ref{++lemma} 
to move the idempotent $e\cir r\eq P_\AA^l$ arising from the composition past 
$\psi$ and then uses $\eps\cir P_\AA^l\eq\eps$ (\Lemma \ref{lem:C=[1]a}(ii)).
\qed

\noindent
It follows that in case a right-adjoint functor $(\alpha_\AA^\pm)^\dagger$
exists, the composition of $\alpha_A^\pm$ with its adjoint functor
is nothing but ordinary induction with respect to $C_{l/r}(A)$,
followed by restriction to $\cC$.

\medskip

We are now in a position to establish
\\[-2.3em] 

\dtl{Proposition}{lem:C=[1]i}
The left and right centers $C_{l/r}$ of a \ssFA\ $A$ inherit natural
structure as a retract of $A$. More precisely, we have:
\\[.2em]
(i)~\,\,$C_l$ and $C_r$ are commutative symmetric Frobenius algebras in $\cC$.
\\[.3em]
(ii)~\,If, in addition, $A$ is simple, then $C_l$ and $C_r$ are simple, too.
\\[.3em]
(iii)~If $C_{l/r}$ is simple, then it is special iff $\dim(C_{l/r})\,{\ne}\,0$.

\medskip\noindent
Proof: \\ (i)~\,We set
  \be \bearll
  m_C := r_C \cir m \cir (e_C \oti e_C) \,, \quad &
  \Delta_C := \zeta^{-1}\, (r_C \oti r_C) \cir \Delta \cir e_C \,,
  \\{}\\[-.5em]
  \eta_C := r_C \cir \eta \,,
  & \eps_C := \zeta\, \eps \cir e_C \,, \eear \labl{eq:Clr-alg}
for some $\zeta\iN\kx$, where $C\,{\equiv}\,C_{l/r}$, and with 
$e_C\,{\equiv}\,e_{C\prec A}^{}$, and $r_C\,{\equiv}\,r_{\!\!A\succ C}^{}$ the 
embedding and restriction morphisms, respectively, for $C$ as a retract of 
$A$. That is, for the product and the unit on $C$ we take the restriction of 
the product on $A$, whereas the coproduct and the counit are only fixed 
up to some invertible scalar.
\\
That $\eta_C$ and $\eps_C$ satisfy the (co-)unit properties follows 
from the corresponding properties of $A$, by \Lemma \ref{lem:C=[1]a}(ii).
The (co-)associativity of $m_C$ and $\Delta_C$ as well as the Frobenius
property are checked with the help of \Lemma \ref{lem:C=[1]a}(iii).
Thus $(C_{l/r},m_C,\eta_C,\Delta_C,\eps_C)$ are indeed Frobenius \alg s.
\\
That $C_l$ is commutative is seen by composing the first of
the equalities \erf{Cl-Cr-defprop} with $\id_{C_l}\oti e_{C_l}$ from below
and with $r_{C_l}$ from above. Commutativity of $C_r$ follows analogously.
Further, commutativity together with triviality of the twist (lemma
\ref{Ctwist-etc}) imply that $C$ is symmetric.
\\[.3em]
(ii)~It follows from \erf{neulab1}, with $U\eq V\eq \one$, and simplicity of 
$A$ that $C$ is haploid, and hence in particular simple.
\\[.3em]
(iii)~The first specialness property holds independently of the value of the
dimension of $C$: with the help of 
\erf{oben1} one finds $\eps_C\cir\eta_C\eq\zeta\dim(A)$, which is non-zero.
\\
Denote by $\eps_{C,\natural}\iN\Hom(A,\one)$ the morphism defined as in 
equation \erf{eq:epsnat}, but with the Frobenius algebra $C$ in place of $A$.
Since $C$ is commutative and simple, it is also haploid, and hence this 
morphism must be a multiple of $\eps_C$. The constant of proportionality 
can be determined by composing the equality with $\eta$; the result is
  \be
  \eps_{C,\natural} = \Frac{\dim(C)}{\dim(A)}\,\zeta^{-1}\,\eps_C \,.  \ee
It follows that $\eps_{C,\natural}$ and $\eps_C$ are non-zero multiples of each 
other iff $\dim(C)\,{\ne}\,0$. On the other hand, equality of of $\eps_{C,
\natural}$ and $\eps_C$ up to a non-zero constant is equivalent to specialness 
of the symmetric Frobenius algebra $C$; see lemma 3.11 of \cite{fuRs4}.
\\
(Note that we recover our usual normalisation convention for special
Frobenius algebras by fixing the scalar factor 
$\zeta$ in \erf{eq:Clr-alg} to $\zeta\eq\dim(C)/\dim(A)\,$.)
\qed

\dt{Remark}
(i)~\,\,Part (ii) of the proposition generalises the classic result that the 
center of a simple $\complex$-algebra is just given by $\complex$.
\\[.3em]
(ii)~\,Alternatively, symmetry of $C_{l/r}$ follows by combining symmetry of $A$
with the identity $\eps\cir P_\AA^{l/r}\eq\eps$ (\Lemma \ref{lem:C=[1]a}(ii)).
As a consequence, triviality of the twist of $C_{l/r}$ (\Lemma \ref{Ctwist-etc})
or, equivalently, \Lemma \ref{lem:C=[1]a}(i), can also be deduced by combining 
\Proposition \ref{lem:C=[1]i} with \Lemma \ref{lem:C=[1]a}(ii).
\\[.3em]
(iii)~In the proof of \Proposition \ref{lem:C=[1]i}(iii) above, as well as at
several other places below, we use conventions and results from \cite{fuRs4}. In 
\cite{fuRs4}, which builds on earlier studies in \cite{fffs3} and \cite{fuSc16}, 
the relevant categories are assumed to be abelian and semisimple. 
The proofs of those results from \cite{fuRs4} that are employed in this paper 
are, however, easily adapted to the present setting.

\medskip

We close this section with another helpful result, to be used later on, in which 
the central idempotents \erf{PE-def} arise. We present the formula with 
$P_\AA^l$; an analogous formula with $P_\AA^r$ holds in which the braiding on 
the \lhs\ is replaced by the opposite braiding.
\\[-2.2em]

\dtl{Lemma}{le4sublocind}
For $A$ a \ssFA\ in a ribbon \cat\ $\cC$, $U$ and $V$ objects of $\cC$,
and $\Phi\iN\Hom(A\Oti U,A\Oti V)$ the following identity holds:
  %%  [pic~01]
  \bea  \begin{picture}(290,155)(0,40)
  \put(0,0)  {\begin{picture}(0,0)(0,0)
             \scalebox{.38}{\includegraphics{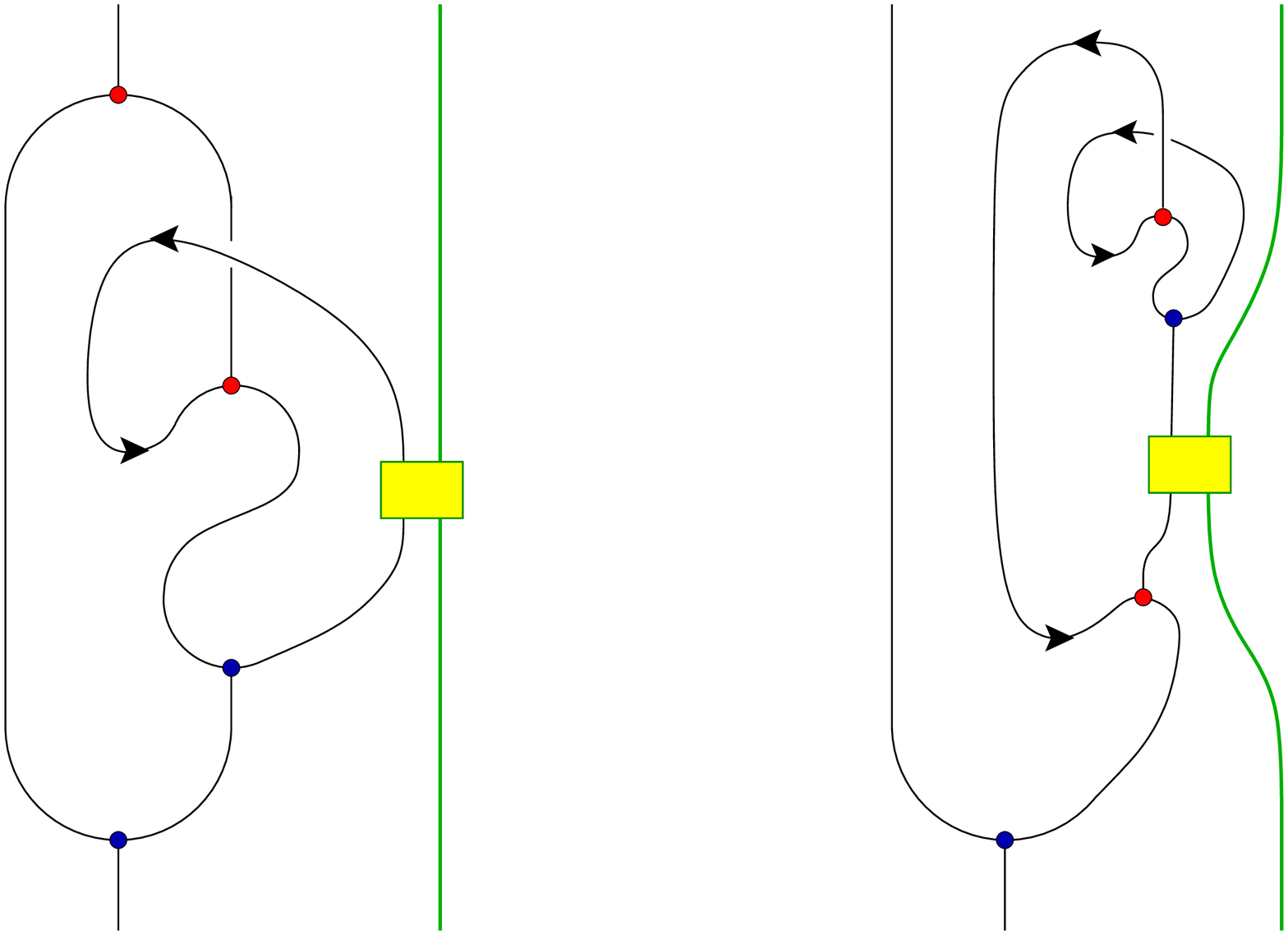}} \end{picture}}
  \put(19.5,-9.1)  {\sse$A$}
  \put(20.6,193)   {\sse$A$}
  \put(82.2,87.8)  {\sse$\Phi$}
  \put(85.5,-9.1)  {\sse$U$}
  \put(86.4,193)   {\sse$V$}
  \put(131,86.6)   {$=$}
  \put(176.8,193)  {\sse$A$}
  \put(199.1,-9.1) {\sse$A$}
  \put(237.2,92.5) {\sse$\Phi$}
  \put(255.9,-9.1) {\sse$U$}
  \put(256.9,193)  {\sse$V$}
  \epicture24 \labl{pic01}
% can also be proven without using specialness

\noindent
Proof:\\
Consider the following manipulations.
  %%  [pic~02]
  \begin{eqnarray}  \begin{picture}(390,155)(17,40)
  \put(0,0)  {\begin{picture}(0,0)(0,0)
             \scalebox{.38}{\includegraphics{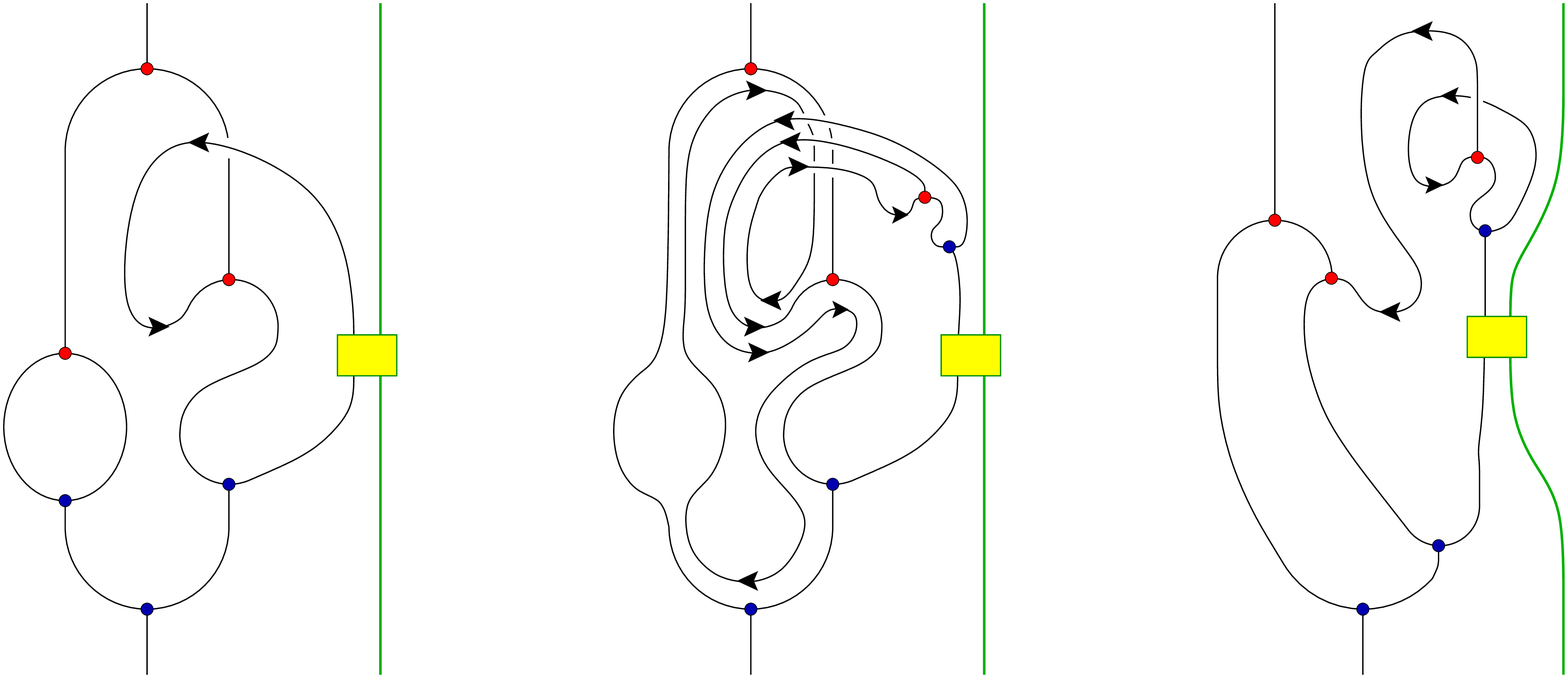}} \end{picture}}
  \put(36.5,-9.1)  {\sse$A$}
  \put(37.5,193)   {\sse$A$}
  \put(99.1,87.8)  {\sse$\Phi$}
  \put(101.5,-9.1) {\sse$U$}
  \put(102.5,193)  {\sse$V$}
  \put(137,86.6)   {$=$}
  \put(205.1,-9.1) {\sse$A$}
  \put(206.1,193)  {\sse$A$}
  \put(268.1,87.8) {\sse$\Phi$}
  \put(271.1,-9.1) {\sse$U$}
  \put(272.1,193)  {\sse$V$}
  \put(303,86.6)   {$=$}
  \put(352.8,193)  {\sse$A$}
  \put(375.8,-9.1) {\sse$A$}
  \put(414.3,93.1) {\sse$\Phi$}
  \put(433.3,-9.1) {\sse$U$}
  \put(434.3,193)  {\sse$V$}
  \end{picture} \nonumber\\[4.1em]{} \label{pic02}
  \\[-2.0em]{}\nonumber\end{eqnarray}
Here in the first step the coproduct and product in the left $A$-loop are
dragged apart, using that $A$ is Frobenius, along the $A$-ribbons 
until they result in the coproduct and product above $\Phi$ in the middle 
picture. The second step is a deformation of the $A$-ribbon that connects 
that coproduct and product, using also the properties \erf{eq:leave-theta} 
and \erf{eq:Hl-def} of the left central idempotent.
\\
The \lhs s of the equations \erf{pic02} and of \erf{pic01} are equal owing 
to specialness of $A$, and their \rhs s are equal because $A$ is special 
Frobenius. Thus \erf{pic01} follows from \erf{pic02}.
\qed

%%%%%%%%%%%%%%%%%%%%%%%%%%%%%%%%%%%%%%%%%%%%%%%%%%%%%%%%%%%%%%%%%%%%%%%%
\bigskip

\sect{Local modules}\label{sect3}

\subsection{Endofunctors related to {\boldmath$\alpha$}-induction}
\label{sec:alg-fun'}

One interesting aspect of \ssF\ algebras $A$ in a ribbon category $\cC$ is
that they allow us to construct functors to the categories of modules 
over the left and right center of $A$, respectively, which are
similar to the induction functor from $\cC$ to the category of $A$-modules.
We call these functors {\em local induction\/} functors. The construction 
makes use of certain endofunctors of $\cC$ which are associated to $A$.

For these endofunctors to exist, the \ssFA\ must have an additional
property. To motivate this property, recall from \Section \ref{LRC}
that for the left and right center of $A$ to exist, the central idempotents 
$P^{l/r}_\AA$ defined in \erf{PE-def} must be split. The construction of 
the endofunctors makes use of similar endomorphisms for each object $U$ 
of $\cC$, namely of the morphisms
  %%  [pic~1]
  \bea  \begin{picture}(205,88)(0,42)
  \put(0,0)   {\begin{picture}(0,0)(0,0)
              \scalebox{.38}{\includegraphics{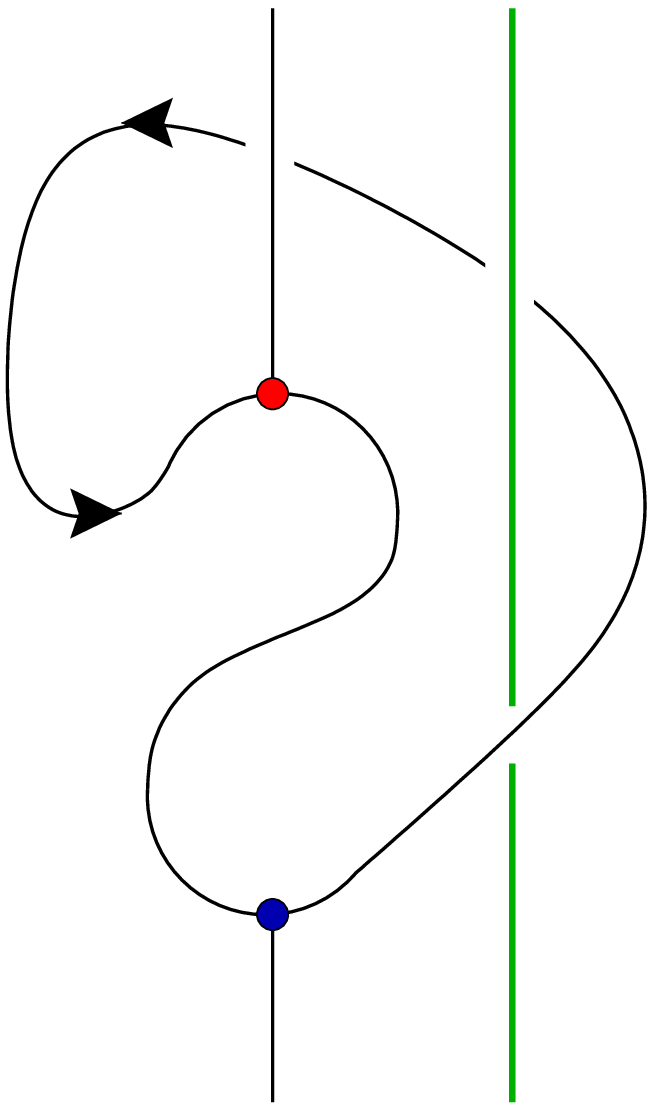}} \end{picture}}
  \put(-58.6,55.5) {$P^l_\AA(U)\;:=$}
  \put(25.5,-7.9)  {\sse$A$}
  \put(52.5,-7.9)  {\sse$U$}
  \put(26.1,125)   {\sse$A$}
  \put(53.3,125)   {\sse$U$}
  \put(73.2,67)    {\sse$A$}
  \put(8.4,32)     {\sse$A$}
  \put(170,0)   {\begin{picture}(0,0)(0,0)
              \put(0,0)   {\begin{picture}(0,0)(0,0)
              \scalebox{.38}{\includegraphics{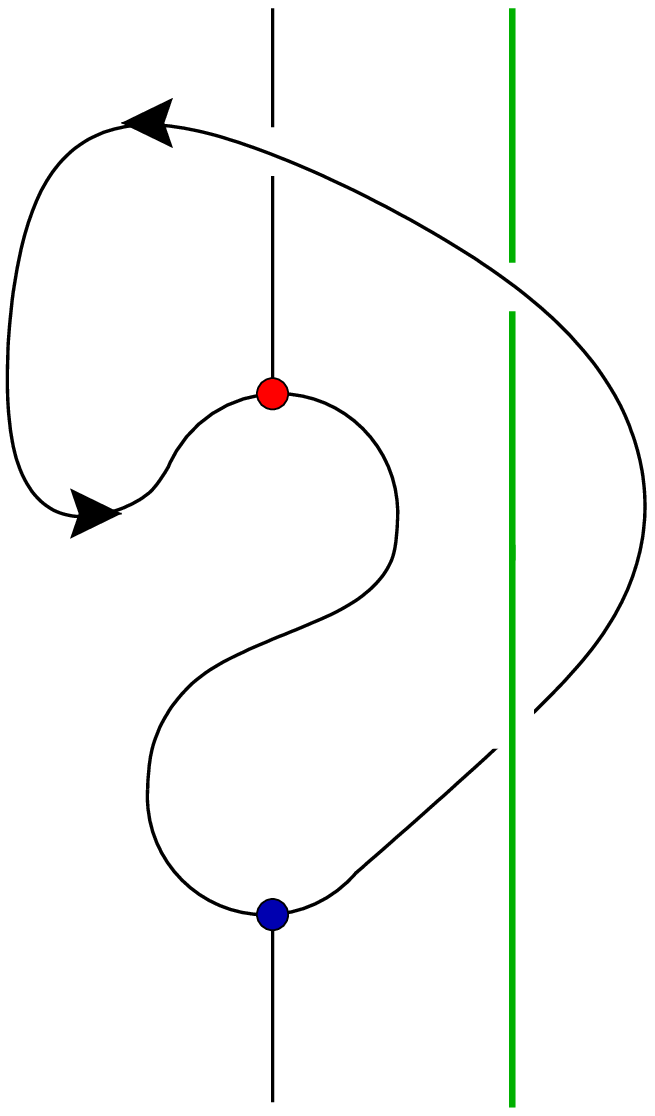}} \end{picture}}
  \put(-58.6,55.5) {$P^r_\AA(U)\;:=$}
  \put(26.1,-7.9)  {\sse$A$}
  \put(52.5,-7.9)  {\sse$U$}
  \put(25.5,125)   {\sse$A$}
  \put(53.3,125)   {\sse$U$}
  \put(73.2,67)    {\sse$A$}
  \put(8.4,32)     {\sse$A$} \end{picture}}
  \epicture30 \labl{PU-def}
in $\Hom(A{\otimes}U{,}A{\otimes}U)$. It is easily verified that, just like
$P^{l/r}_\AA\,{\equiv}\,P^{l/r}_\AA(\one)$ these endomorphisms are idempotents 
(for $P_\AA^l(U)$ this has already been shown in lemma 5.2 of \cite{fuRs4}.)

\dtl{Definition}{csplit}
A special Frobenius algebra $A$ in a ribbon category $\cC$ is called 
{\em \csplit\/} iff the idempotents \erf{PU-def} are split for every $U\iN\Objc$.

\medskip

Clearly, in a Karoubian ribbon category (and hence in particular in a modular 
tensor category) every Frobenius algebra is \csplit. Recall, however, that 
occasionally we want to allow for non-Karoubian \cats. Then
we need \csplit\ algebras in order to ensure the 
existence of the desired endofunctors. Accordingly we make the following
\\[-2.4em]

\dtl{\Convention}{c:csplit}
In the sequel every special Frobenius algebra $A$ will be assumed to be \csplit.

\bigskip 

With this agreement in mind, we can now proceed to
\\[-2.3em]

\dtl{Definition}{def:[]-functor}
For $A$ a \ssFA\ in a ribbon category $\cC$, the \operation s
$\EFU{l/r}\AA$ are defined on objects and morphisms of $\cC$ as follows.
\\
For $U\iN\Objc$, $\efu U{l/r}\AA$ are the retracts
  \be
  \efu Ul\AA := \Im P^l_\AA(U) \qquad {\rm and }\qquad
  \efu Ur\AA := \Im P^r_\AA(U)   \ee
of the induced module $A\oti U$, with the idempotents
$P^{l/r}_\AA(U)\iN\Hom(A{\otimes}U{,}A{\otimes}U)$ given by \erf{PU-def}.
\\
For $f\iN\Hom(U,V)$,
$\efu f{l/r}\AA\iN\Hom(\efu U{l/r}\AA,\efu V{l/r}\AA)$ are the morphisms
  %%  [pic~16]
  \bea  \begin{picture}(225,87)(0,19)
  \put(50,0)  {\begin{picture}(0,0)(0,0)
              \scalebox{.38}{\includegraphics{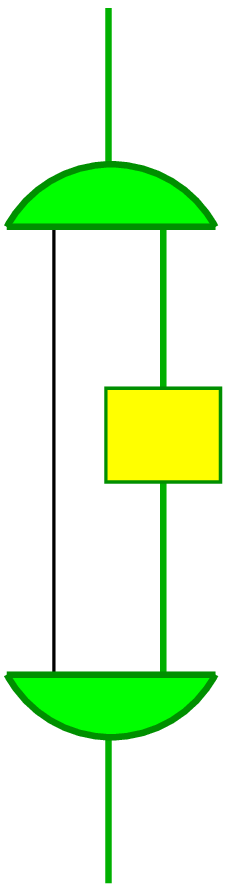}} \end{picture}}
  \put(-12,48.1)   {$\efu fl\AA\;:=$}
  \put(52.7,-9.1)  {\sse$\efu Ul\AA$}
  \put(53.3,102.1) {\sse$\efu Vl\AA$}
  \put(60.2,19.1)  {\tiny$e$}
  \put(60.2,74.9)  {\tiny$r$}
  \put(65.6,48.1)  {\sse$f$}
  \put(47.7,43)    {\sse$A$}
     \put(69,30)   {\sse$U$}
     \put(69,62)   {\sse$V$}
  \put(200,0) {\begin{picture}(0,0)(0,0)
              \scalebox{.38}{\includegraphics{Pf-def.eps}} \end{picture}}
  \put(138,48.1)   {$\efu fr\AA\;:=$}
  \put(202.7,-8.8) {\sse$\efu Ur\AA$}
  \put(203.3,102.1){\sse$\efu Vr\AA$}
  \put(210.2,19.1) {\tiny$e$}
  \put(210.2,74.9) {\tiny$r$}
  \put(215.6,48.1) {\sse$f$}
  \put(197.7,43)   {\sse$A$}
     \put(219,30)  {\sse$U$}
     \put(219,62)  {\sse$V$}
  \epicture07 \labl{Pf-def}
with $e\,{\equiv}\,e_{\efu U{l/r}\AA\!\prec A\otimes U}$
and $r\,{\equiv}\,r_{\!\!A\otimes V\succ\efu V{l/r}\AA}$.

\medskip

Let us remark that this construction is non-trivial only in a genuinely 
braided \tc. For, when $\cC$ is a {\em symmetric\/} \tc, the projection 
just amounts to considering the objects $C\oti U$, where $C$ is the center 
of the \alg\ $A$. Note that these are precisely the objects that underlie 
induced $C$-modules; as we will see later, the objects $\efu U{l/r}\AA$ 
naturally carry a module structure, too: they are modules over the left 
and right center of $A$, respectively.

\dtl{Proposition}{lem:A-functor}
The \operation s $\EFU{l/r}\AA$ are endofunctors of $\cC$.

\medskip\noindent
Proof:\\
Let $E$ stand for one of $\EFU l\AA$, $\EFU r\AA$. It follows from the 
definitions \erf{PU-def} and \erf{Pf-def} 
that for any $g\iN\Hom(U,V)$ we have $E(g) \iN \Hom(E(U),E(V))$,
i.e.\ $E(g)$ is in the correct space. It remains to check that for any
$g'\iN\Hom(V,W)$ one has $E(g'\cir g)\eq E(g') \cir E(g)$ and that
$E(\id_U)\eq \id_{E(U)}$. The second property is obvious because
$\efu{\id_U}{l/r}\AA\eq r^{l/r}{\circ}\, e^{l/r}$ is indeed nothing but
the identity morphism $\id_{E(U)}$ on the retract.
For the first property we note that, writing out the definitions for
$E(g' \cir g)$ and $E(g') \cir E(g)$, these two morphisms only differ by an 
idempotent \erf{PU-def}. By functoriality of the braiding we can shift this 
idempotent past $g$ so that it gets directly composed with the embedding 
morphism $e$, and then \erf{eq:S-prop} tells us that it can be left out.
\qed

These functors are, however, in general {\em not\/} tensor functors.
%    counter example: SU(2) level 4, simple current algebra.
%    $E(0) = 0\oplus 4 \not\cong (0)$

\medskip

The following lemma will be used in the proof of \Proposition \ref{alphaalpha}.
\\[-2.3em]

\dtl{Lemma}{lem:under-over}
(i)~For every \ssFA\ $A$ in a ribbon category $\cC$ and every $U\iN\Objc$,
and with right $A$-actions $\rr^\pm$ defined as in \erf{rr+-}, we have
   \be \bearll  P_\AA^l(U) \circ \rr^- \!\!&
   \equiv P_\AA^l(U) \circ (m\oti\id_U) \cir (\id_\AA\oti c_{A,U}^{-1})
   \\{}\\[-.8em]
   &= P_\AA^l(U) \circ
   \Llb [ m \cir c_{A,A} \cir (\id_\AA\oti \theta_\AA^{}) ] \oti\id_U \Lrb
   \cir (\id_\AA\oti c_{U,A}) \,,
   \\{}\\[-.4em]
   P_\AA^r(U) \circ \rr^+ \!\!&
   \equiv P_\AA^r(U) \circ (m\oti\id_U) \cir (\id_\AA\oti c_{U,A})
   \\{}\\[-.8em]
   &= P_\AA^r(U) \circ
   \Llb [ m \cir c_{A,A} \cir (\id_\AA\oti \theta_\AA^{-1}) ] \oti\id_U \Lrb
   \cir (\id_\AA\oti c_{A,U}^{-1}) \,.
   \eear \labl{eq70}
(ii)~If $A$ is in addition commutative, then
  \be  P_\AA(U) \circ \rr^+ = P_\AA(U) \circ \rr^-  \ee
for $P_\AA^{}(U)\,{\equiv}\,P_\AA^{l/r}(U)$.

\medskip\noindent
Proof:\\
(i)~\,The first of the formulas \erf{eq70} follows by the moves
  %%  [pic~70]
  \bea  \begin{picture}(360,138)(0,21)
  \put(0,0)  {\begin{picture}(0,0)(0,0)
              \scalebox{.38}{\includegraphics{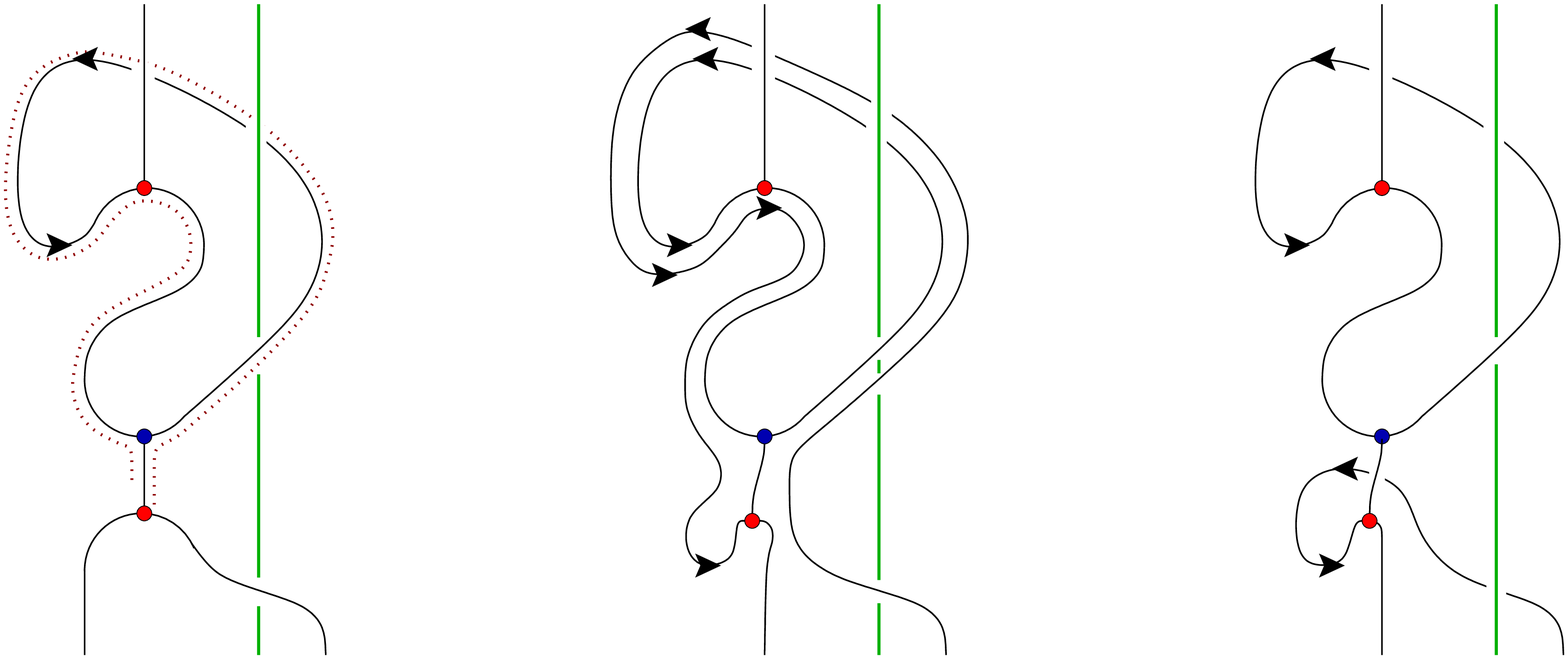}} \end{picture}}
  \put(15.1,-8.5)  {\sse$A$}
  \put(29.4,153.4) {\sse$A$}
  \put(29.5,26.8)  {\sse$m$}
  \put(55.9,-8.5)  {\sse$U$}
  \put(56.2,153.4) {\sse$U$}
  \put(70.6,-8.5)  {\sse$A$}
  \put(107,73.1)   {$=$}
  \put(171.6,-8.5) {\sse$A$}
  \put(171.6,153.4){\sse$A$}
  \put(198.6,-8.5) {\sse$U$}
  \put(198.6,153.4){\sse$U$}
  \put(213.6,-8.5) {\sse$A$}
  \put(249,73.1)   {$=$}
  \put(314.6,-8.5) {\sse$A$}
  \put(314.6,153.4){\sse$A$}
  \put(341.6,-8.5) {\sse$U$}
  \put(341.6,153.4){\sse$U$}
  \put(356.6,-8.5) {\sse$A$}
  \epicture09 \labl{pic70}
In the first picture, the dotted line is not part of the morphism, but rather
only indicates a path along which the product that is marked explicitly is
`dragged' (using functoriality of the braiding, as well as associativity and
the Frobenius property of $A$) so as to arrive at the first
equality. The second equality is obtained by deforming the $A$-ribbon that
results from this dragging.
\\
The second of the formulas \erf{eq70} is seen analogously, with under-
and overbraidings interchanged.
\\[.3em]
(ii) follows immediately form (i) by using that $A$ has trivial twist 
        (\Proposition \ref{c+s=tt}(i)) 
and the definition of commutativity.
(Also, in the commutative case we actually have $P_\AA^l(U)\eq P_\AA^r(U)$,
see the picture \erf{Pl-r} below.)
\qed

Note that, obviously, the assertions made in the lemma are non-trivial
only if the \tc\ $\cC$ is genuinely braided. The same remark applies to
several other statements below, in particular to \Theorem \ref{thm:equiv}.
(Compare also to the considerations at the end of \Section \ref{intro-locind}.)

Assume now that there exist right-adjoint functors $(\alpha_\AA^{\pm})^\dagger$ 
to the $\alpha_\AA^\pm$-induction functors. The following result 
shows that in this case the endofunctors $\EFU{l/r}\AA$ can be regarded as 
the composition of $(\alpha^{\pm}_\AA)^\dagger$ with $\alpha^\mp_\AA$.
 %  because only objects of the form \alpha^+(.) appear on rhs,
 %  but not arbitary bimodules.
(The result does not imply that such right-adjoint functors exist. They 
certainly do exist, though, if $\cC$ is semisimple with finite number of 
non-isomorphic simple objects, in particular if $\cC$ is modular.)
\\[-2.3em]

\dtl{Proposition}{alphaalpha}  
For every \ssFA\ $A$ in a ribbon category $\cC$ and any two objects
$U,V\iN\Objc$ there are natural bijections
  \be  \Hom(\efu Ul\AA,V) \,\cong\, \HomAA(\alpha_\AA^-(U),\alpha_\AA^+(V))
  \,\cong\,\Hom(U,\efu Vr\AA) \labl{-U+V}
and
  \be  \Hom(\efu Ur\AA,V) \,\cong\, \HomAA(\alpha_\AA^+(U),\alpha_\AA^-(V))
  \,\cong\,\Hom(U,\efu Vl\AA) \,. \labl{+U-V}

\medskip\noindent
Proof:\\
Let us start with the first equivalence in \erf{-U+V}. Recall that
according to the reciprocity relation \erf{reciMV} there
is a natural bijection $\Phi{:}\ \Hom(A\Oti U,V)\,{\stackrel\cong\to}\,
\HomA(\Ind_\AA(U),\Ind_\AA(V))$, and note that the target of this bijection 
contains the middle expression of \erf{-U+V} as a natural subspace.
\\
Furthermore, in view of \Lemma \ref{lem:rezI}(i), by definition of
$\efu {\cdot}l\AA$ we may identify the \lhs\ of \erf{-U+V} with the
subspace $\Hom_{(P_\AA^l(U))}(A\Oti U,V)$ of $\Hom(A\Oti U,V)$.
Thus it is sufficient to show that $\Phi$ restricts to a bijection between
this subspace and $\HomAA(\alpha_\AA^-(U),\alpha_\AA^+(V))$. The map $\Phi$ 
and its inverse are defined similarly as in formula \erf{PhiPsi}; they act as
  \be  \varphi \,\mapsto\, (\id_\AA\oti\varphi) \circ (\Delta\oti\id_U)
  \qquad {\rm and} \qquad
  \psi \,\mapsto\, (\eps\oti\id_V) \circ \psi  \labl{AUAV}
for $\varphi\iN\Hom(A\Oti U,V)$ and $\psi\iN\HomA(\Ind_\AA(U),\Ind_\AA(V))$,
respectively. The following considerations show that $\Phi$ and its inverse
restrict to linear maps between $\Hom_{(P_\AA^l(U))}(A\Oti U,V)$ and
$\HomAA(\alpha_A^-(U),\alpha_A^+(V))$.
\\
First, for $\varphi\iN\Hom_{(P_\AA^l(U))}(A\Oti U,V)$ the morphism
$\Phi(\varphi)\cir\rr^-(U)\,{\equiv}\,\Phi(\varphi\cir P_\AA^l(U))
\cir\rr^-(U)$ is given by the \lhs\ of the equality
  %%  [pic~67]
  \bea  \begin{picture}(210,156)(0,40)
  \put(0,0)  {\begin{picture}(0,0)(0,0)
             \scalebox{.38}{\includegraphics{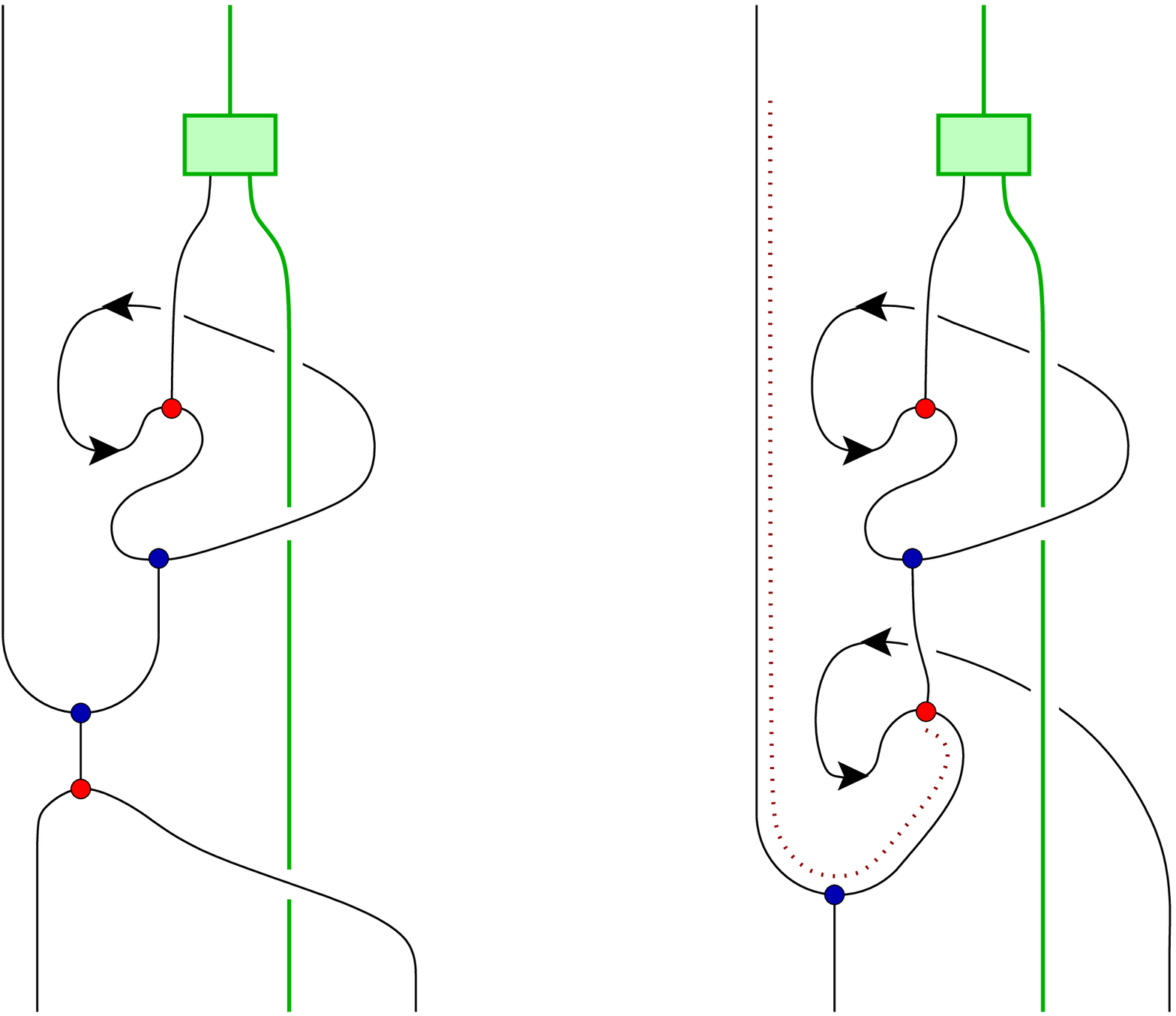}} \end{picture}}
  \put(-2.9,191.8)  {\sse$A$}
  \put(3.1,-9.1)    {\sse$A$}
  \put(39.8,160.4)  {\sse$\varphi$}
  \put(40.4,191.8)  {\sse$V$}
  \put(49.7,-9.1)   {\sse$U$}
  \put(73.3,-9.1)   {\sse$A$}
  \put(97,89.5)     {$=$}
  \put(137.1,191.8) {\sse$A$}
  \put(150.4,-9.1)  {\sse$A$}
  \put(165.5,49.5)  {\sse$m$}
  \put(178.9,160.4) {\sse$\varphi$}
  \put(179.5,191.8) {\sse$V$}
  \put(188.6,-9.1)  {\sse$U$}
  \put(211.6,-9.1)  {\sse$A$}
  \epicture26 \labl{pic67}
This equality, in turn, is a straightforward application of lemma
\ref{lem:under-over}. Further, by dragging the marked product along the path
indicated by the dashed line and deforming the resulting ribbon (using 
functoriality of the braiding) such that the braiding occurs above the 
morphism $\varphi$ and omitting again the idempotent $P_\AA^l(U)$ then yields
the graphical description of $\rr^+(U)\cir(\Phi(\varphi){\otimes}\id_\AA)$.
\\
This shows that $\Phi(\varphi)$ is indeed a morphism of $\alpha$-induced
bimodules.
\\[.3em]
The required property of $\Phi^{-1}$ is obtained by the following manipulations,
valid for every $\psi\iN\HomAA(\alpha_A^-(U),\alpha_A^+(V))$:
  %%  [pic~72]
  \bea  \begin{picture}(330,242)(0,37)
  \put(0,0)  {\begin{picture}(0,0)(0,0)
             \scalebox{.38}{\includegraphics{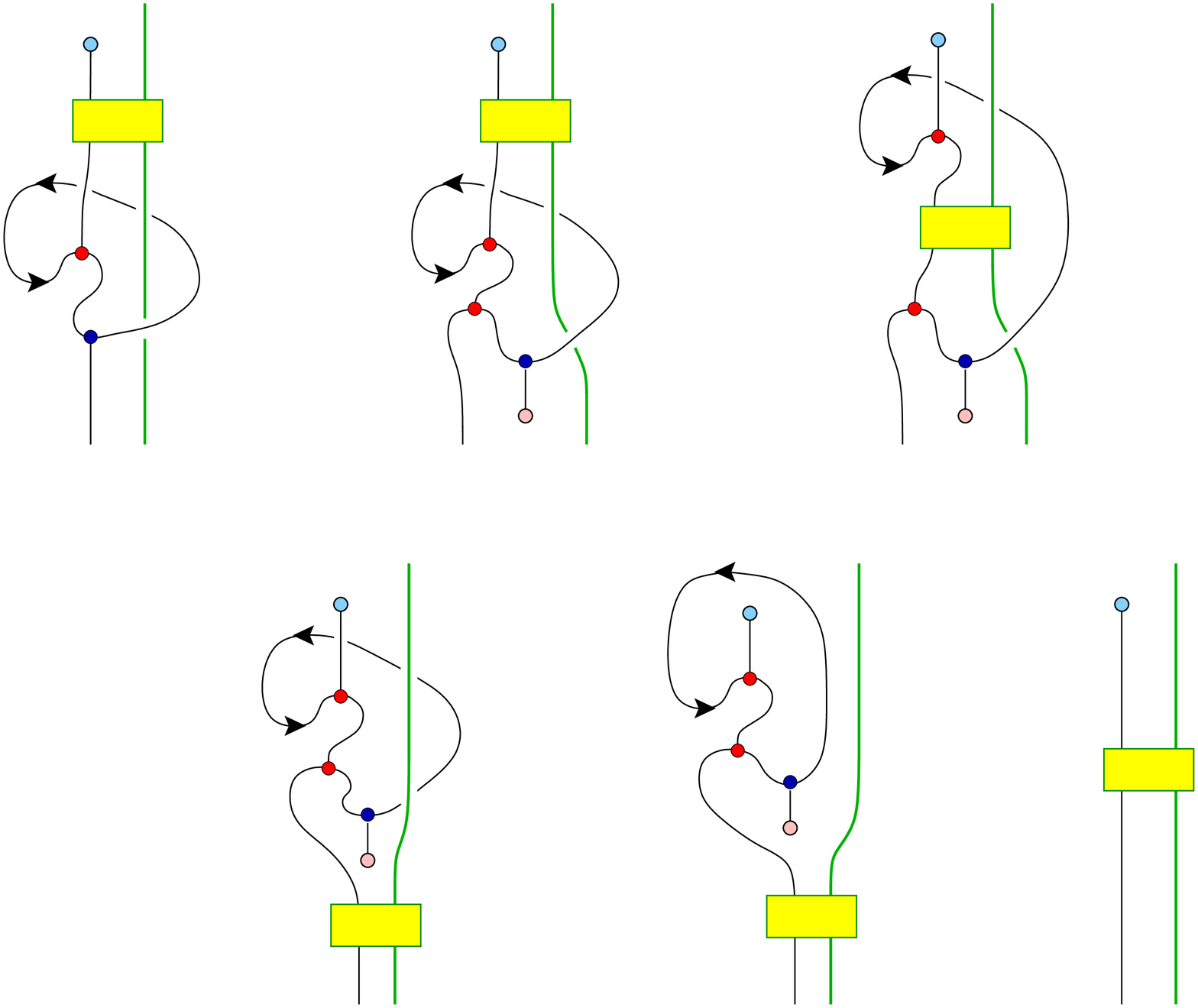}} \end{picture}}
  \put(19.7,142.1)  {\sse$A$}
  \put(28.4,237.9)  {\sse$\psi$}
  \put(35.7,142.1)  {\sse$U$}
  \put(36.8,274.9)  {\sse$V$}
  \put(78,206.5)    {$=$}
  \put(120.7,142.1) {\sse$A$}
  \put(139.2,237.9) {\sse$\psi$}
  \put(146.8,274.9) {\sse$V$}
  \put(155.5,142.1) {\sse$U$}
  \put(195,206.5)   {$=$}
  \put(238.7,142.1) {\sse$A$}
  \put(257.1,208.7) {\sse$\psi$}
  \put(264.8,274.9) {\sse$V$}
  \put(272.9,142.1) {\sse$U$}
  \put(38,59.5)     {$=$}
  \put(84.4,34.1)   {\sse$A$}
  \put(91.7,-9.1)   {\sse$A$}
  \put(98.7,20.3)   {\sse$\psi$}
  \put(103.6,-9.1)  {\sse$U$}
  \put(107.6,124.2) {\sse$V$}
  \put(150,59.5)    {$=$}
  \put(209.5,-9.1)  {\sse$A$}
  \put(216.2,22.2)  {\sse$\psi$}
  \put(221.1,-9.1)  {\sse$U$}
  \put(228.3,124.2) {\sse$V$}
  \put(261,59.5)    {$=$}
  \put(298.7,-9.1)  {\sse$A$}
  \put(307.1,62.1)  {\sse$\psi$}
  \put(313.7,-9.1)  {\sse$U$}
  \put(314.4,124.2) {\sse$V$}
  \epicture27 \labl{pic72}
The first equality uses that $A$ is Frobenius; the second and third use 
functoriality of the braiding and the fact that $\psi$ intertwines the 
$A$-bimodules $\alpha_A^-(U)$ and $\alpha_A^+(V)$ (more specifically, that 
$\psi$ is a morphism of left modules for the second, and that it is a morphism 
of right modules for the third equality); and the fourth is just a deformation 
of the $A$-loop. The last equality combines the fact that $A$ is symmetric 
Frobenius and the identification
of the counit with the morphism $\eps_\natural$ \erf{eq:epsnat}.
Thus indeed $\Phi^{-1}(\psi) \cir P^l_\AA(U)\eq\Phi^{-1}(\psi)$.
\\[.3em]
Next consider the second equivalence in \erf{+U-V}. In this case we can use
the natural bijection $\tilde\Phi{:}\ \Hom(U, A\Oti V)\,{\stackrel\cong\to}\,
\HomA(\Ind_\AA(U),\Ind_\AA(V))$ as well as the equivalence
$\Hom(U,\efu Vl\AA) \,{\cong}%\, 
    $\linebreak[0]$%
\Hom^{(P_\AA^l(V))}(U, A\Oti V)$, see equation~\erf{eq:dim-up-p}.
Explicitly, the linear map $\tilde\Phi$ and its inverse are given by
  \be  \varphi \,\mapsto\, (m\oti\id_V) \circ (\id_\AA\oti\varphi)
  \qquad {\rm and} \qquad
  \psi \,\mapsto\, \psi \circ (\eta\oti\id_U) \,. \labl{AUAV2}
Similarly to the argument above, one shows that $\tilde\Phi$ and its inverse
restrict to linear maps between $\Hom^{(P_\AA^l(V))}(U, A\Oti V)$ and
$\HomAA(\alpha_\AA^+(U),\alpha_\AA^-(V))$. For example, the pictures
occurring in the proof of $\tilde\Phi(\varphi)\cir\rr^+(U)\eq
\rr^-(U)\cir(\tilde\Phi(\varphi){\otimes}\id_A)$ look like the ones in 
\erf{pic67} except that they are `reflected' about a horizontal axis.
\\[.3em]
The remaining two equivalences are derived analogously.
\qed

\dtl{Remark}{rem:[U]A-as-obj}
When $\cC$ is in addition semisimple, then it follows that
the objects $\efu U{l/r}\AA$ decompose into simple objects as
  \be
  \efu Ul\AA \cong \bigoplus_{i\in\cI} \Big( \sum_{q\in\cI}
  \tilde Z(A)_{iq} \, n_q \Big) U_i
  \quad\ {\rm and} \quad\
  \efu Ur\AA \cong \bigoplus_{i\in\cI} \Big(\sum_{q\in\cI}
  n_q \, \tilde Z(A)_{qi} \Big) U_i \,,  \labl{eq:[U]A-as-obj}
with the non-negative integers $n_q$ defined by the decomposition
$U \,{\cong}\, \bigoplus_q n_q U_q$ of $U$.
Thus when expressing objects as direct sums of the simple objects $U_i$
with $i\iN\II$, the action of the functor $\efu{\,\cdot\,}{l/r}\AA$ on
objects amounts to multiplication
from the left and right, respectively, with the matrix $\tilde Z(A)$.

%%%%%%%%%%%%%%%%%%%%%%%%%%%%%%%%%%%%%%%%%%%%%%%%%%%%%%%%%%%%%%%%%%%%%%%%

\subsection{Endofunctors on categories of algebras}

One datum contained in a Frobenius algebra $(B,m,\Delta,\eta,\epsilon)$
is the object $B\iN\Objc$, on which we can consider the action of the 
endofunctors $\EFU{l/r}A$ associated to some \ssF\ algebra $A$. We wish 
to show that the objects $\efu B{l/r}\AA$ carry again the structure of a
Frobenius algebra. (This would be 
obvious if the functors $\EFU{l/r}\AA$ were tensor functors, because then 
we could simply take $\efu m{l/r}\AA$ as the multiplication morphism. But 
this is not the case, in general.) This will imply that $\EFU{l/r}\AA$
also provides us with endofunctors on the category of Frobenius algebras in
$\cC$. Since $\efu B{l/r\!}\AA$ is a retract of $A\Oti B$, what we first
need is the notion of a tensor product of two Frobenius algebras $A$ and $B$.

For any pair $A,B$ of Frobenius algebras in a ribbon category there are in 
fact two natural Frobenius algebra structures -- to be denoted by $A{\otimes
^\pm_{}}B\,{\equiv}\, (A{\otimes}B, m^\pm_{A{\otimes}B}, \eta^\pm_{A{\otimes}B},
\Delta^\pm_{A{\otimes}B}, \eps^\pm_{A{\otimes}B})$ -- on the tensor product
object $A\Oti B$. For the case of $\otimes^+_{}$, the structural morphisms are
  \bea \begin{picture}(255,120)(0,31)
  \put(60,88) {\begin{picture}(0,0)(0,0)
              \scalebox{.38}{\includegraphics{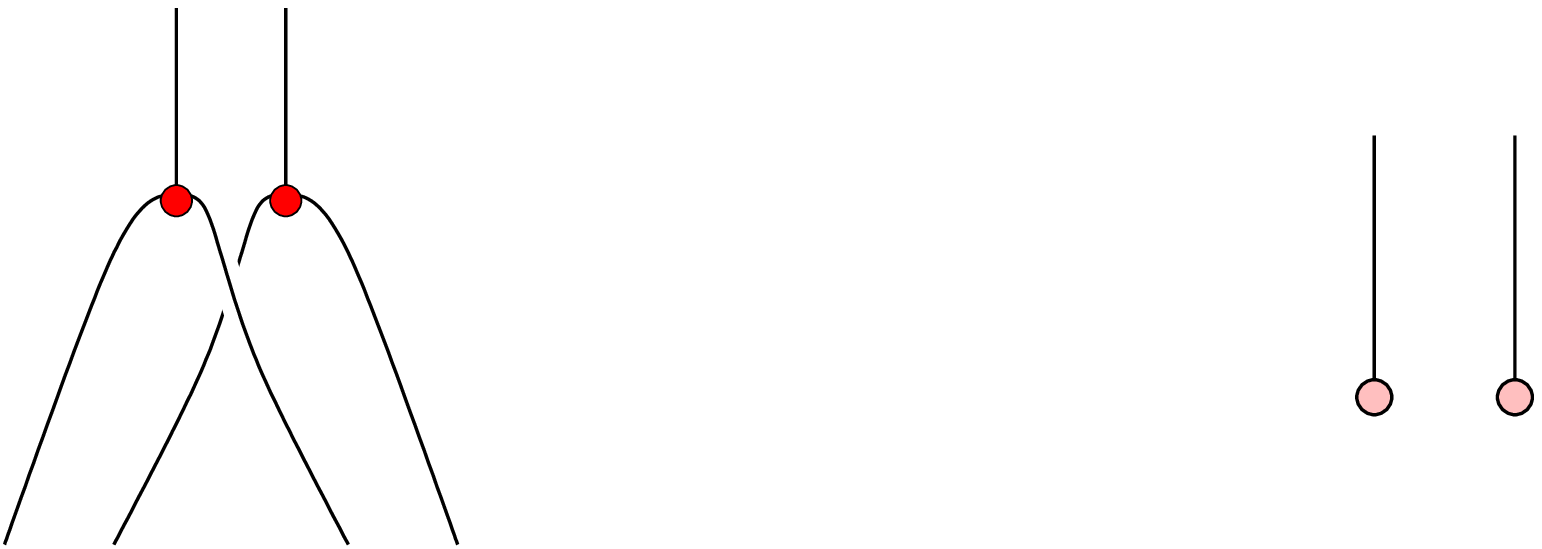}} \end{picture}}
  \put(60,0)  {\begin{picture}(0,0)(0,0)
              \scalebox{.38}{\includegraphics{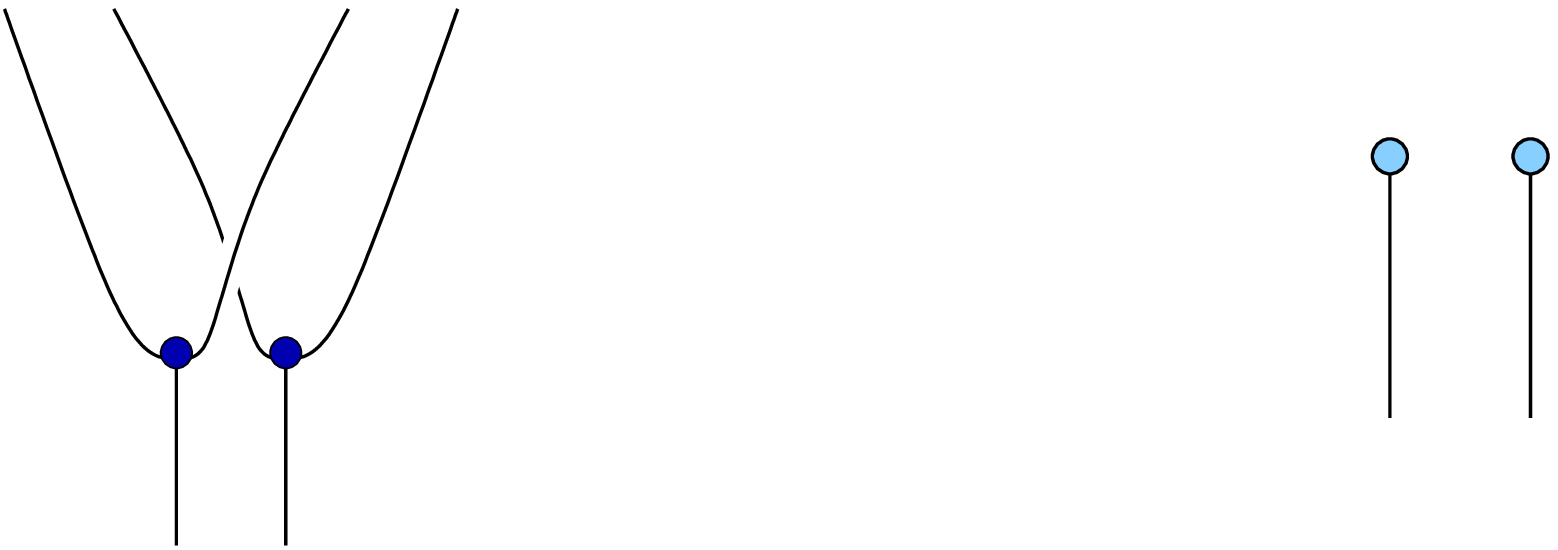}} \end{picture}}
  \put(0,119.9)    {$m^+_{A{\otimes}B} \;:=$}
  \put(55.4,79.3)  {\sse$A$}
  \put(68.4,79.3)  {\sse$B$}
  \put(75.7,150.8) {\sse$A$}
  \put(88.7,150.8) {\sse$B$}
  \put(94.8,79.3)  {\sse$A$}
  \put(108.1,79.3) {\sse$B$}
  \put(0,27.9)     {$\Delta^+_{A{\otimes}B} \;:=$}
  \put(57.2,62.9)  {\sse$A$}
  \put(69.4,62.9)  {\sse$B$}
  \put(75.5,-9.2)  {\sse$A$}
  \put(87.5,-9.2)  {\sse$B$}
  \put(95.8,62.9)  {\sse$A$}
  \put(108.1,62.9) {\sse$B$}
  \put(155,119.9)  {$\eta^+_{A{\otimes}B} \;:=$}
  \put(207.2,137.3){\sse$A$}
  \put(223.8,137.3){\sse$B$}
  \put(155,27.9)   {$\eps^+_{\AA{\otimes}B} \;:=$}
  \put(209.4,4.3)  {\sse$A$}
  \put(224.4,4.3)  {\sse$B$}
  \epicture19 \labl{mabdab}
while for $\otimes^-_{}$ over-braiding and under-braiding must be exchanged 
in the definition of both the product and the coproduct. One verifies by 
direct substitution that $A{\otimes^+}B$ is again a Frobenius algebra. 
Further, if $A,B$ are in addition symmetric and special, then so is 
$A{\otimes^+}B$. An analogous statement holds for $A{\otimes^-}B$.

In the sequel we will work with $\otimes^+_{}$; also, we slightly abuse
notation and simply write $A\Oti B$ in place of $A{\otimes^+_{}}B$ for
the tensor product of two Frobenius algebras.

\medskip

Note that even when both $A$ and $B$ are commutative, their tensor product
$A\Oti B$ is not commutative, in general. More precisely, if $A$ and $B$
are commutative, then $A\Oti B$ is commutative iff
$c_{A,B}\cir c_{B,A}\eq \id_{A\otimes B}$. While this identity holds in a
symmetric tensor category, it does not necessarily hold in a genuinely
braided \tc; in this setting it is therefore not advisable to restrict one's
attention exclusively to (braided-) commutative algebras.

\medskip

\dtl{Proposition}{prop:AB-alg}
Let $A$ be a symmetric special Frobenius algebra and $B$ a Frobenius algebra
in a ribbon category. Then the following holds:
\\[.3em]
(i)~\,$\efu B{l/r\!}\AA \,{\equiv}\, (\efu B{l/r\!}\AA, m_{l/r}, \eta_{l/r},
\Delta_{l/r}, \eps_{l/r})$, with morphisms given by
 %% [pic~2], [pic~17]
  \bea \begin{picture}(390,218)(0,71)
    \put(0,157) {\begin{picture}(0,0)(0,0)
               \put(44,35) {\begin{picture}(0,0)(0,0)
               \scalebox{.38}{\includegraphics{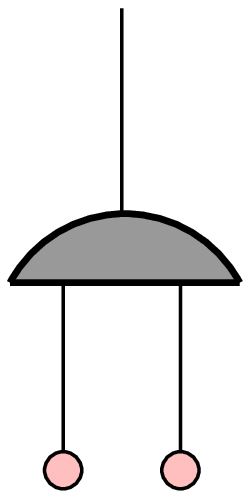}} \end{picture}}
               \put(141,0) {\begin{picture}(0,0)(0,0)
               \scalebox{.38}{\includegraphics{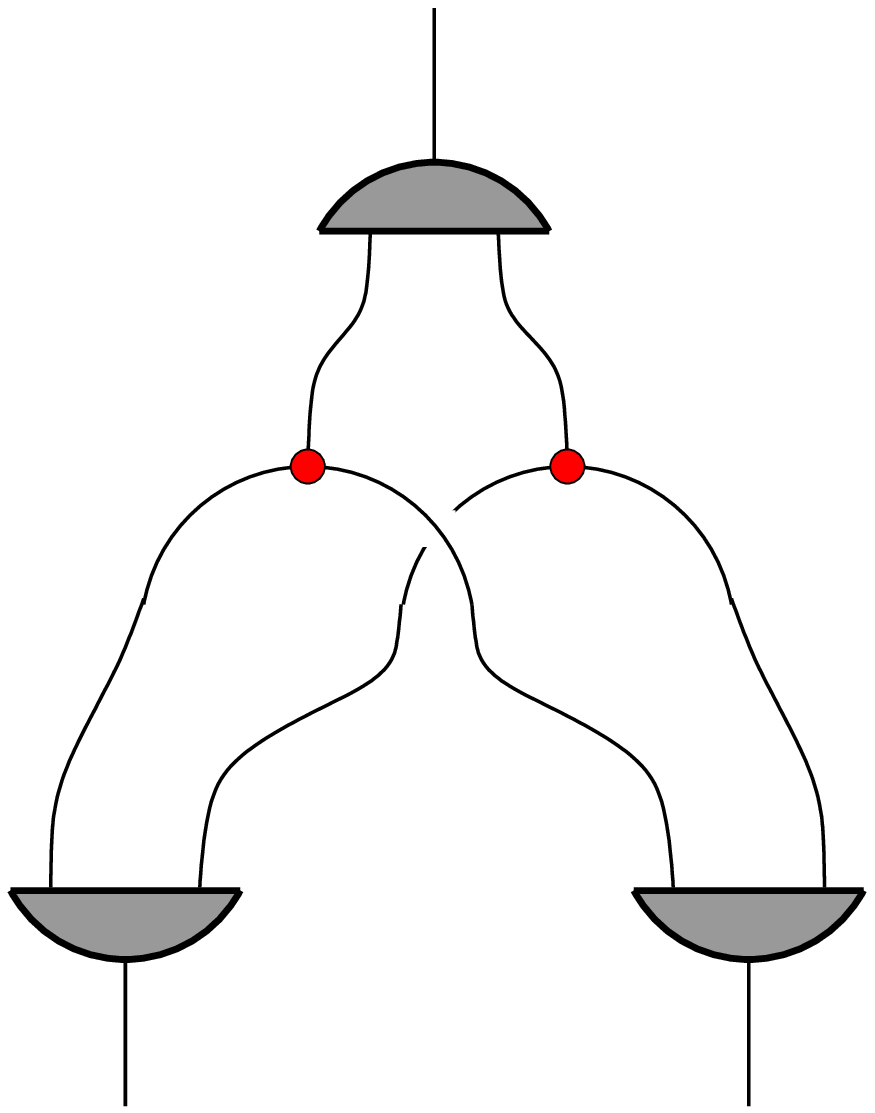}} \end{picture}}
               \put(287,0) {\begin{picture}(0,0)(0,0)
               \scalebox{.38}{\includegraphics{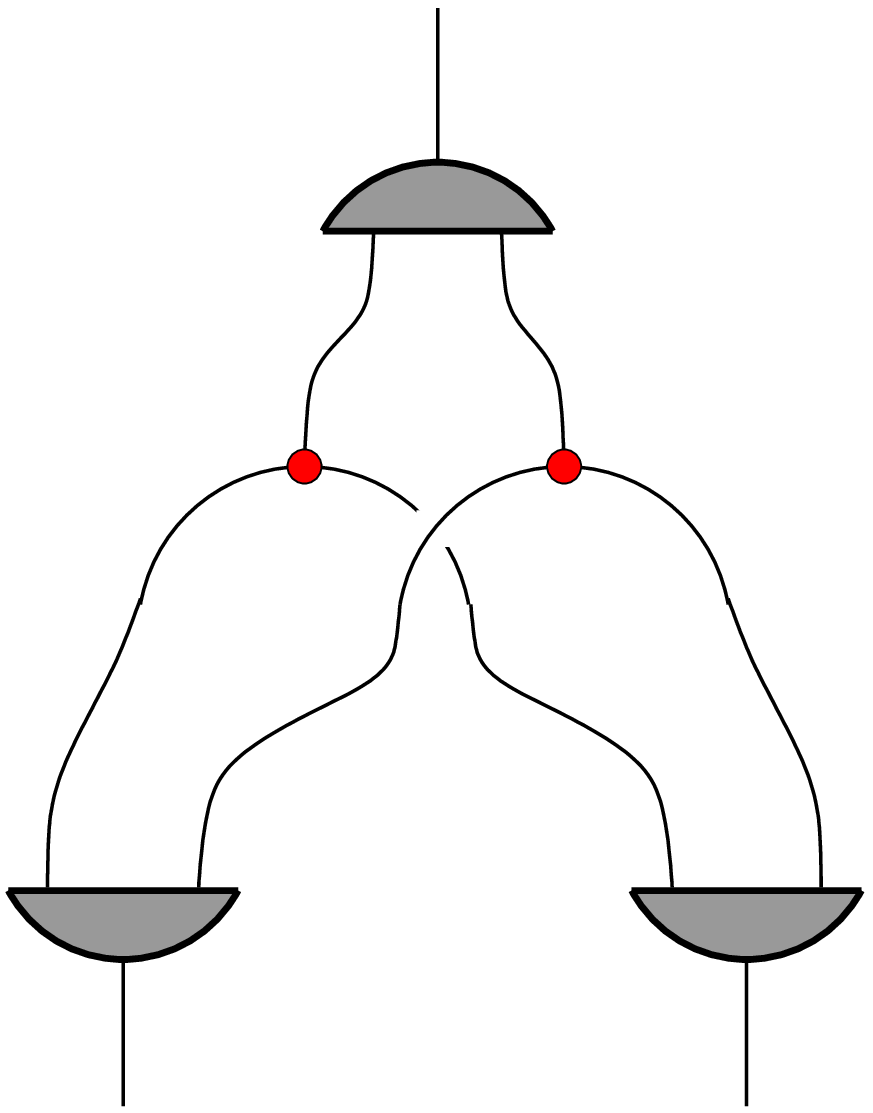}} \end{picture}}
  \put(0,58.4)     {$\eta_{l/r}\,:=$}
  \put(38.9,35.5)  {\tiny$\eta_A^{}$}
  \put(46.2,94.3)  {\sse$\efu B{l/r}\AA$}
  \put(55.5,60.2)  {\sse$r$}
  \put(67.5,35.5)  {\tiny$\eta_B^{}$}
  \put(117,58.4)   {$m_l\,:=$}
  \put(144.9,-8.5) {\sse$\efu B{l}\AA$}
  \put(153.3,18.9) {\sse$e$}
  \put(161.7,73.9) {\tiny$m_{\!A}^{}$}
  \put(177.5,127.3){\sse$\efu B{l}\AA$}
  \put(186.3,99.1) {\sse$r$}
  \put(204.7,72.9) {\tiny$m_B^{}$}
  \put(213.4,-8.5) {\sse$\efu B{l}\AA$}
  \put(221.3,18.9) {\sse$e$}
  \put(261,58.4)   {$m_r\,:=$}
  \put(290.9,-8.5) {\sse$\efu B{r}\AA$}
  \put(299.3,18.9) {\sse$e$}
  \put(307.7,73.9) {\tiny$m_{\!A}^{}$}
  \put(323.5,127.3){\sse$\efu B{r}\AA$}
  \put(332.3,99.1) {\sse$r$}
  \put(350.7,72.9) {\tiny$m_B^{}$}
  \put(359.4,-8.5) {\sse$\efu B{r}\AA$}
  \put(367.3,18.9) {\sse$e$}
    \end{picture}}
    \put(0,0) {\begin{picture}(0,0)(0,0)
               \put(44,35) {\begin{picture}(0,0)(0,0)
               \scalebox{.38}{\includegraphics{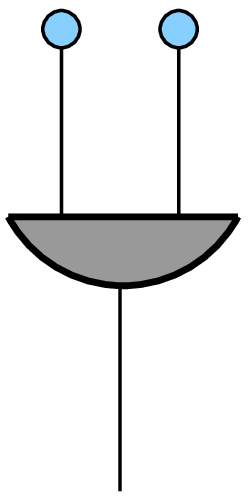}} \end{picture}}
               \put(141,0) {\begin{picture}(0,0)(0,0)
               \scalebox{.38}{\includegraphics{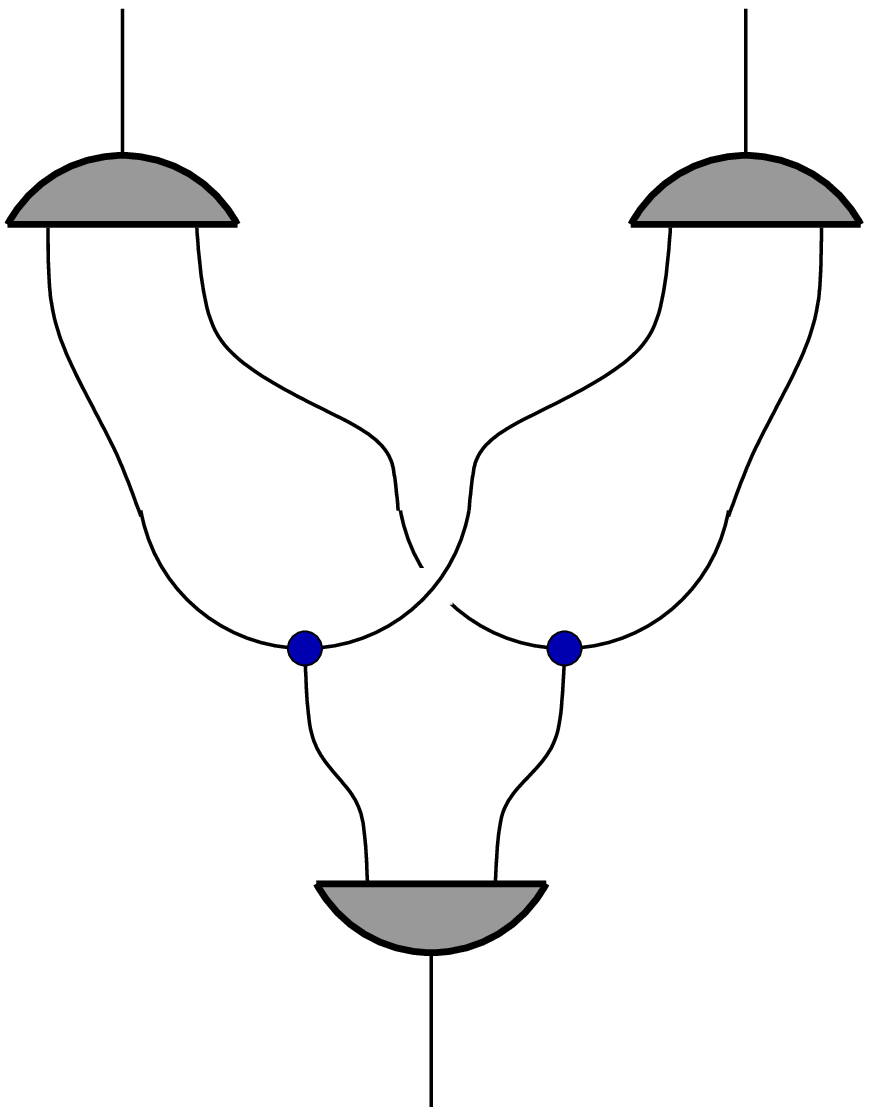}} \end{picture}}
               \put(287,0) {\begin{picture}(0,0)(0,0)
               \scalebox{.38}{\includegraphics{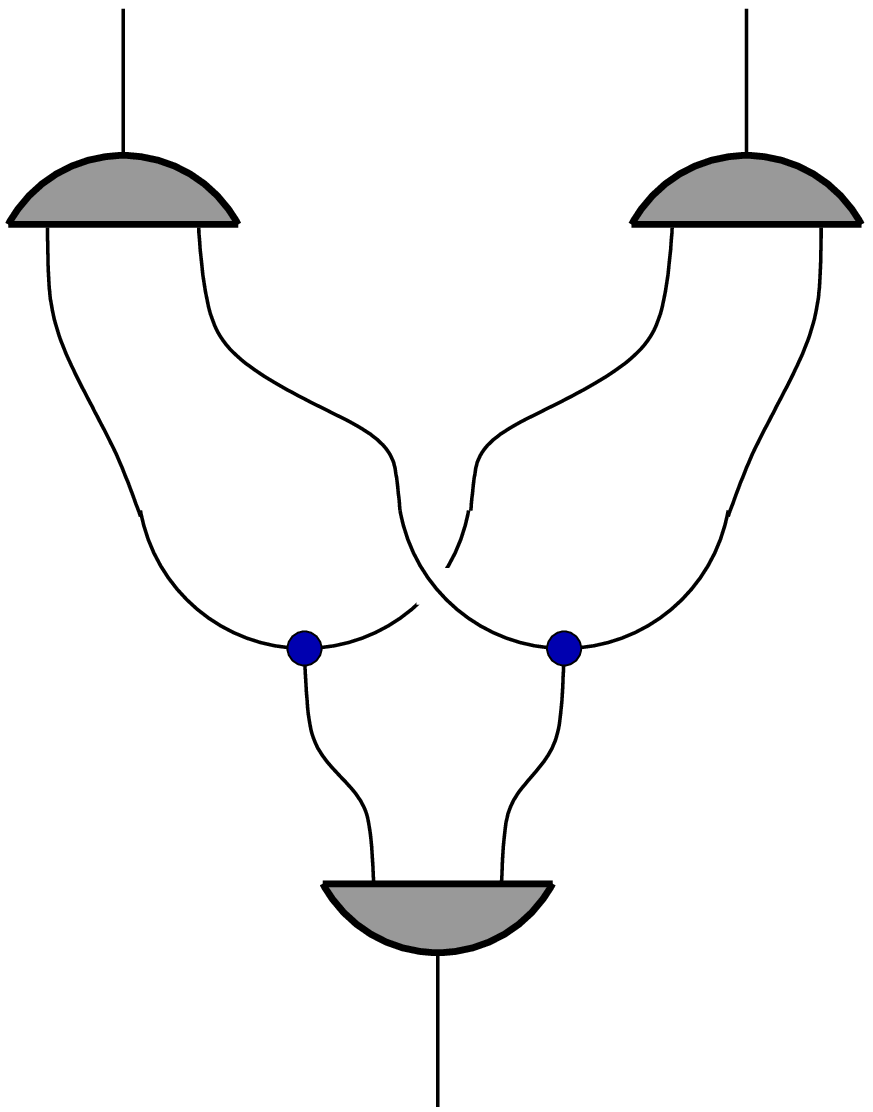}} \end{picture}}
  \put(-13,58.4)   {$\xi_{l/r}^{-1}\eps_{l/r}:=$}
  \put(39.8,88.9)  {\tiny$\eps_{\!A}^{}$}
  \put(43.2,23.7)  {\sse$\efu B{l/r}\AA$}
  \put(55.2,60.7)  {\sse$e$}
  \put(66.9,88.9)  {\tiny$\eps_B^{}$}
  \put(110,58.4)   {$\xi_l\Delta_l:=$}
  \put(144.9,126.3){\sse$\efu B{l}\AA$}
  \put(153.3,99.5) {\sse$r$}
  \put(163.1,45.5) {\tiny$\Delta_{\!A}^{}$}
  \put(178.2,-8.5) {\sse$\efu B{l}\AA$}
  \put(186.3,19.9) {\sse$e$}
  \put(203.9,44.9) {\tiny$\Delta_B^{}$}
  \put(213.4,126.3){\sse$\efu B{l}\AA$}
  \put(221.3,99.5) {\sse$r$}
  \put(254,58.4)   {$\xi_r\Delta_r:=$}
  \put(290.9,126.3){\sse$\efu B{r}\AA$}
  \put(299.3,99.5) {\sse$r$}
  \put(308.9,45.5) {\tiny$\Delta_{\!A}^{}$}
  \put(322.5,-8.5) {\sse$\efu B{r}\AA$}
  \put(332.3,19.9) {\sse$e$}
  \put(349.8,44.9) {\tiny$\Delta_B^{}$}
  \put(359.4,126.3){\sse$\efu B{r}\AA$}
  \put(367.3,99.5) {\sse$r$}
    \end{picture}}
  \epicture51 \labl{eq:[B]A-alg}
with $\xi_{l/r}\iN\kx$, is a Frobenius algebra.
\\[.3em]
(ii)~\,If $B$ is symmetric, then $\efu Bl\AA$ and $\efu Br\AA$ are symmetric.\\
   If $B$ is commutative, then $\efu Bl\AA$ and $\efu Br\AA$ are commutative.
\\[.3em]
(iii)~If $\efu Bl\AA$ is symmetric, $B$ is in addition special,
$\dim_\koerper\Hom(B,C_r(A))\eq1$, and $\dim(\efu Bl\AA)$ is non-zero,
then $\efu Bl\AA$ is in addition haploid and special.
\\
If $\efu Br\AA$ is symmetric, $B$ is in addition special,
$\dim_\koerper\Hom(B,C_l(A))\eq1$, and $\dim(\efu Br\AA)$ is non-zero,
then $\efu Br\AA$ is in addition haploid and special.
\\[.3em]
(iv)~If $A$ is commutative, then $\EFU l\AA\eq\EFU r\AA$ as functors. More 
precisely, for every $U\iN\Objc$ we have the equality $\efu Ul\AA\eq\efu Ur\AA$ 
as objects in $\cC$, and for every morphism $f$ of $\cC$ we have
$\efu fl\AA\eq\efu fr\AA$.
\\[.3em]
(v)~\,If $A$ is commutative, then $\efu Bl\AA\eq\efu Br\AA$ as Frobenius \alg s.

\bigskip

\dtl{Remark}{vazhrem}
In \cite{vazh} the notion of an Azumaya algebra in a braided tensor category 
has been introduced. The definition in \cite{vazh} can be seen to be 
equivalent to the following one: An algebra $A$ in a ribbon category $\cC$ is 
called an {\em Azumaya algebra\/} iff the functors $\alpha_\AA^+$ and 
$\alpha_\AA^-$ from $\cC$ to $\cC_{A|A}$ are equivalences of tensor 
categories. If a symmetric special Frobenius algebra $A$ is Azumaya, then 
$C_l(A) \,{\cong}\, \one \,{\cong}\, C_r(A)$. To see this note that if 
$\alpha_\AA^+$ is an equivalence functor, then it has a left and right adjoint 
$(\alpha_\AA^+)^\dagger_{}$, given by $(\alpha_\AA^+)^{-1}$. In 
Assertions (i)\,--\,(iii) of \Proposition \ref{lem:[U]A-as-obj} 
we have seen that the composition 
$(\alpha_\AA^+)^\dagger_{} \cir \alpha_\AA^+$ corresponds to tensoring with 
$C_l(A)$. This is an equivalence iff $C_l(A) \,{\cong}\, \one$. A similar
argument shows that $C_r(A) \,{\cong}\, \one$.
\\
Assertions (i)\,--\,(iii) of \Proposition \ref{prop:AB-alg} thus imply in 
particular that every Azumaya algebra defines two endofunctors of the full 
subcategory of haploid commutative symmetric special Frobenius algebras in 
a given ribbon category $\cC$. Algebras of the latter type can be used to 
construct new ribbon categories starting from $\cC$, see \Proposition 
\ref{thm:mod} below.

\bigskip

The proof of \Proposition \ref{prop:AB-alg} 
will fill the remainder of this section. We need the following three lemmata.
\\[-2.4em]

\dtl{Lemma}{lem:remove-P}
Let $A$ and $B$ be as in \Proposition \ref{prop:AB-alg}, 
and $\tilde m_l\iN\Hom((A{\otimes}B)\oti (A{\otimes}B),A{\otimes}B)$ denote the 
morphism obtained from $m_l$ of \erf{eq:[B]A-alg} by omitting the embedding 
and restriction morphisms $e,r$; define $\tilde m_r$ and $\tilde \Delta_{r/l}$ 
similarly. Further let $\tilde\eta\df \eta_A\oti\eta_B$ and $\tilde\eps\df 
\eps_A\oti\eps_B$. The idempotent $P^l \,{\equiv}\, P^l_\AA(B)$ fulfills
  \be  \bearll
   P^l \circ \tilde\eta = \tilde\eta \,, &\!\!\! 
   \tilde\eps \circ P^l = \tilde\eps \,,
  \\{}\\[-.7em]
   P^l \circ \tilde m_l \circ (P^l\oti P^l)
  &\!\!\!  = \id_{A{\otimes}B} \circ \tilde m_l \circ (P^l\oti P^l)
  \\{}\\[-.7em]
  &\!\!\!  = P^l \circ \tilde m_l \circ (\id_{A{\otimes}B}\oti P^l)
   = P^l \circ \tilde m_l \circ (P^l\oti \id_{A{\otimes}B}) \,.
  \eear\labl{eq:remove-P}
Analogous relations hold for
$P^r \,{\equiv}\, P^r_\AA(B)$ and $\tilde m_r, \tilde \eta, \tilde \eps$,
as well as for $P^{l/r}$ and $\tilde\Delta_{l/r}$.
\\
In terms of the graphical calculus, this means that at any product or coproduct 
vertex for which each of the three attached ribbons carries an idempotent 
$P^l$, or each a $P^r$, any one out of the three idempotents can be omitted.

\medskip\noindent
Proof:\\
The proof is similar for all relations. As examples we present it for 
$\tilde\eps \cir P^l \eq \tilde\eps$ and for $(P^r\oti\id_{A{\otimes}B}) \cir 
\tilde\Delta^r \cir P^r$. The first of these relations is easily seen from
  %% [pic~33]
  \bea  \begin{picture}(260,83)(0,20)
  \put(0,0)   {\begin{picture}(0,0)(0,0)
              \scalebox{.38}{\includegraphics{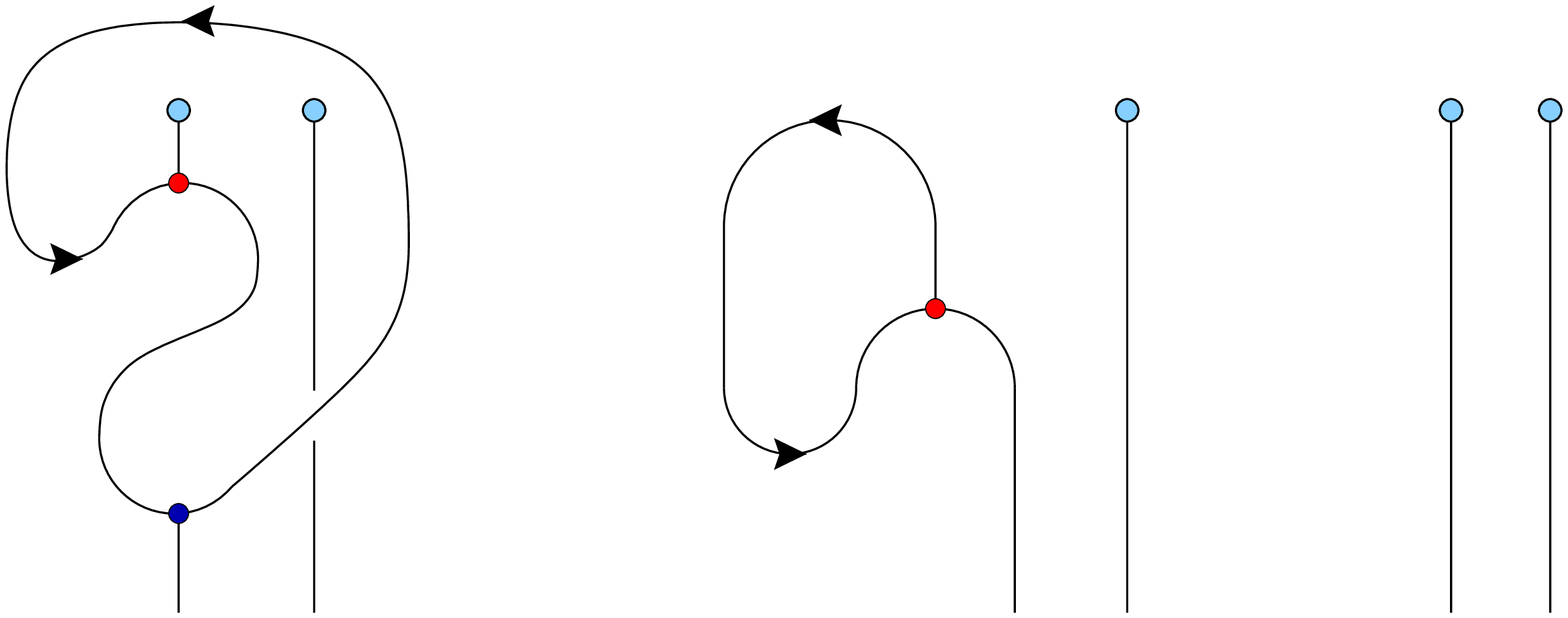}} \end{picture}}
  \put(26.2,-9.2) {\sse$A$}
  \put(49.5,-9.2) {\sse$B$}
  \put(171.1,-9.2){\sse$A$}
  \put(191.1,-9.2){\sse$B$}
  \put(246.5,-9.2){\sse$A$}
  \put(264.7,-9.2){\sse$B$}
  \put(93,44.1)  {$=$}
  \put(219,44.1) {$=$}
  \epicture09 \labl{epst-Pl}
In the first step one substitutes the definition of $P^l$ and deforms the 
graph slightly; then one uses the Frobenius and counit properties to get rid 
of the counit of $A$. The final step re-introduces this counit by using the 
fact that it obeys $\eps\eq\eps_\natural$ with $\eps_\natural$ given by 
\erf{eq:epsnat} (and also that $A$ is symmetric Frobenius).
\\[.3em]
To obtain the second relation one considers the following series of
transformations, for which all defining properties of the \ssFA\ $A$
are needed:
  %% [pic~34]
  \bea \mbox{\hspace{-1.9em}}
  \xi_r\,(P^r\oti P^r)\circ\tilde\Delta^r\circ P^r 
  \\
        \begin{picture}(380,216)(5,-24)
  \put(0,0)   {\begin{picture}(0,0)(0,0)
              \scalebox{.38}{\includegraphics{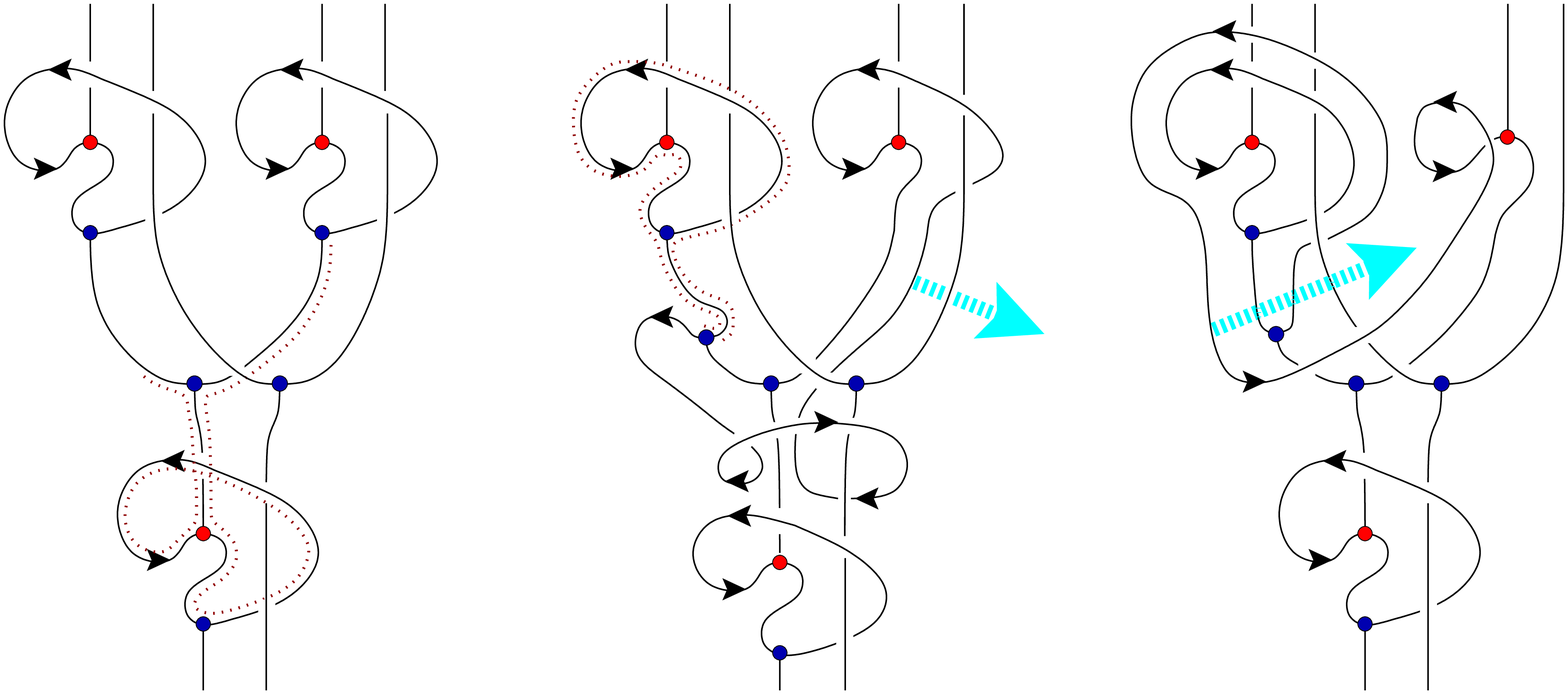}} \end{picture}}
  \put(-15,81.1)   {$=$}
  \put(18.4,176.6) {\sse$A$}
  \put(34.1,176.6) {\sse$B$}
  \put(46.2,-9.2)  {\sse$A$}
  \put(61.4,-9.2)  {\sse$B$}
  \put(67.9,111.6) {\tiny$\Delta_{\!A}^{}$}
  \put(77.1,176.6) {\sse$A$}
  \put(92.8,176.6) {\sse$B$}
  \put(123,81.1)   {$=$}
  \put(163.2,176.6){\sse$A$}
  \put(166.9,83.1) {\tiny$\Delta_{\!A}^{}$}
  \put(178.4,176.6){\sse$B$}
  \put(189.4,-9.2) {\sse$A$}
  \put(206.9,-9.2) {\sse$B$}
  \put(221.2,176.6){\sse$A$}
  \put(236.9,176.6){\sse$B$}
  \put(273,81.1)   {$=$}
  \put(309.3,176.6){\sse$A$}
  \put(307.9,85.6) {\tiny$\Delta_{\!A}^{}$}
  \put(325.5,176.6){\sse$B$}
  \put(327.9,71.6) {\tiny$\Delta_{\!A}^{}$}
  \put(337.2,-9.2) {\sse$A$}
  \put(352.4,-9.2) {\sse$B$}
  \put(372.2,176.6){\sse$A$}
  \put(387.3,176.6){\sse$B$}
  \epicture-1 \labl{pic34}
{}\\[-3.5em]
The first equality just consists of writing out the \Definition 
\erf{eq:[B]A-alg}
of the coproducts and the idempotents. To arrive at the second equality, one
drags the coproduct that is marked explicitly in the first graph along the 
path that is drawn as a dotted line, so that its new location is the one 
marked in the second graph. The third equality is obtained by first pulling 
an $A$-ribbon under the right $B$-ribbon, which is indicated by the big 
shaded arrow, and then moving it back in the opposite direction, but this 
time {\em over\/} the $B$-ribbon (as well as over another $A$-ribbon). 
In addition, one continues to drag the coproduct that was already moved 
during the previous step, now along the dotted path in the second graph;
this way it returns to the same location, but is now attached from the
opposite side.
\\
Starting from the third graph, one can now pull the left-most $A$-ribbon in 
the direction of the shaded arrow, over various $A$-ribbons as well as over 
the left $B$-ribbon; the twists on this ribbon then cancel. Afterwards one can 
use co-associativity (on the two coproducts that are marked explicitly in the
graph) and then the specialness of $A$ so as to arrive at the desired result.
\qed

\dtl{Lemma}{lem:add-twist}
For every \ssFA\ $A$ we have
  %%  [pic~36]
  \bea  \begin{picture}(290,47)(0,34)
  \put(0,0)  {\begin{picture}(0,0)(0,0)
             \scalebox{.38}{\includegraphics{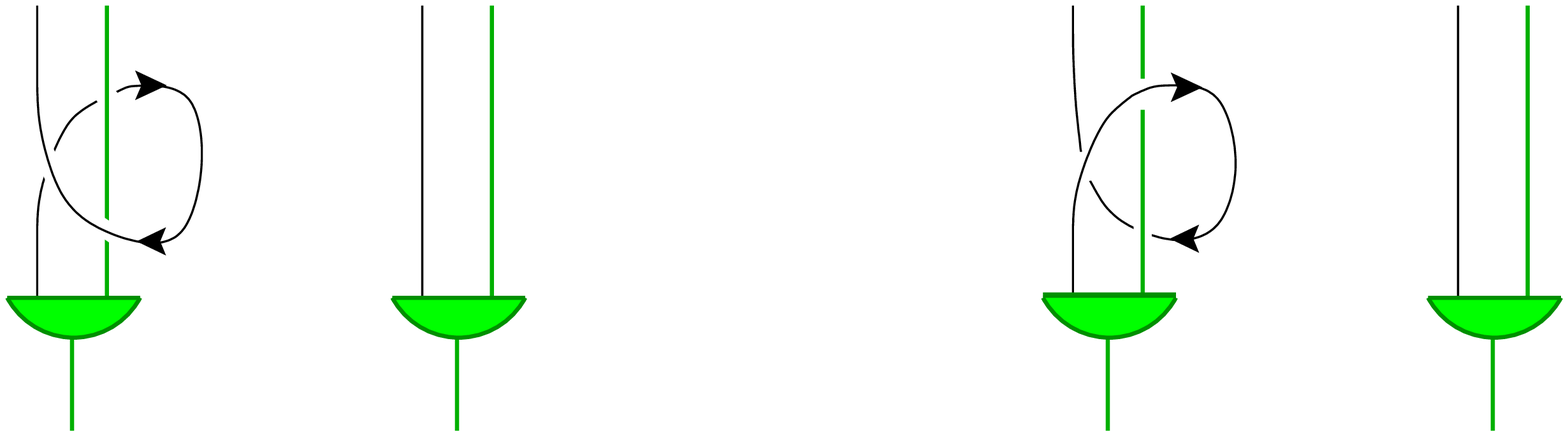}} \end{picture}}
  \put(2.7,78.5)   {\sse$A$}
  \put(2.1,-9.1)   {\sse$\efu Ul\AA$}
  \put(15.7,78.5)  {\sse$U$}
  \put(47,37.5)    {$=$}
  \put(69.3,78.5)  {\sse$A$}
  \put(68.7,-9.1)  {\sse$\efu Ul\AA$}
  \put(82.3,78.5)  {\sse$U$}
  \put(124,37.5)   {and}
  \put(181.3,78.5) {\sse$A$}
  \put(180.7,-9.1) {\sse$\efu Ur\AA$}
  \put(194.3,78.5) {\sse$U$}
  \put(226,37.5)   {$=$}
  \put(248.3,78.5) {\sse$A$}
  \put(247.2,-9.1) {\sse$\efu Ur\AA$}
  \put(260.6,78.5) {\sse$U$}
  \epicture21 \labl{eq:add-twist}
as well as the analogous relations for $r_{\!\!\!A\otimes U \succ
\efu U{l/r}\AA}$ instead of $e_{\!\efu U{l/r}\AA\prec A\otimes U}$.

\medskip\noindent
Proof:\\
We show the moves needed to derive the left equality -- the right one
and the relations for $r$ follow analogously:
  %%  [pic~37]
  \bea  \begin{picture}(280,75)(0,47)
  \put(0,0)  {\begin{picture}(0,0)(0,0)
             \scalebox{.38}{\includegraphics{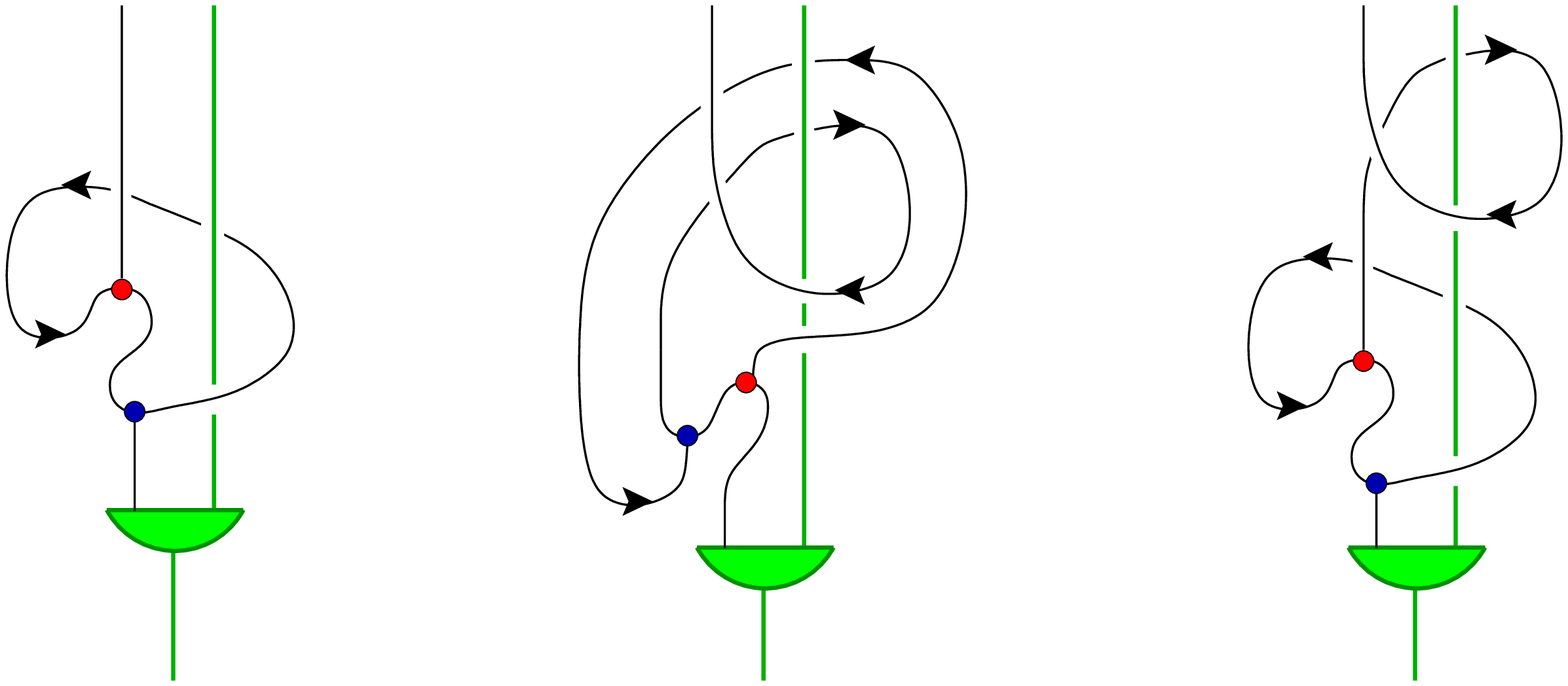}} \end{picture}}
  \put(16.2,119.1) {\sse$A$}
  \put(17.5,-9.1)  {\sse$\efu Ul\AA$}
  \put(32.7,119.1) {\sse$U$}
  \put(68,55.5)    {$=$}
  \put(116.1,119.1){\sse$A$}
  \put(117.8,-9.1) {\sse$\efu Ul\AA$}
  \put(133.1,119.1){\sse$U$}
  \put(179,55.5)   {$=$}
  \put(224.8,119.1){\sse$A$}
  \put(226.1,-9.1) {\sse$\efu Ul\AA$}
  \put(241.3,119.1){\sse$U$}
  \epicture32 \labl{pic37}
The first expression is the \rhs\ of the first equality in
\erf{eq:add-twist}, with a redundant idempotent $P^l_\AA(U)$ inserted. To 
arrive at the second graph one uses the Frobenius property and suitably drags 
the resulting coproduct along part of the $A$-ribbon. A further deformation 
and application of the Frobenius property results in the graph on the \rhs. 
In this last expression the idempotent $P^l_\AA(U)$ is again redundant; 
removing it yields the \lhs\ of the first equality in \erf{eq:add-twist}.
\qed

\dtl{Lemma}{le:mE0=mE1}
Let $A$ be a symmetric special Frobenius algebra and $B$ a Frobenius
algebra in a ribbon \cat. Denote by $m_E$ the multiplication morphism
of $E \eq \efu B{l/r\!}\AA$ as defined in \erf{eq:[B]A-alg}. Then
  \be
  m_E \circ c_{E,E} = m_{E'}  \labl{mEE'}
with $E' \eq \efu {B'}{l/r\!}\AA$, where
$B' \eq (B,m_B\cir c_{B,B},\eta_B,\Delta_B\cir c_{B,B}^{-1},\eps_B)$
(i.e., the opposite algebra of $B'$ is $B$).

\medskip\noindent
Proof:\\
We prove the relation for $\efu Bl\AA$, the case of $\efu Br\AA$ being 
analogous. Consider the following moves:
  %% [pic~38]
  \begin{eqnarray}  \begin{picture}(305,120)(-30,14)
  \put(-80,0)   {\begin{picture}(0,0)(0,0)
              \scalebox{.38}{\includegraphics{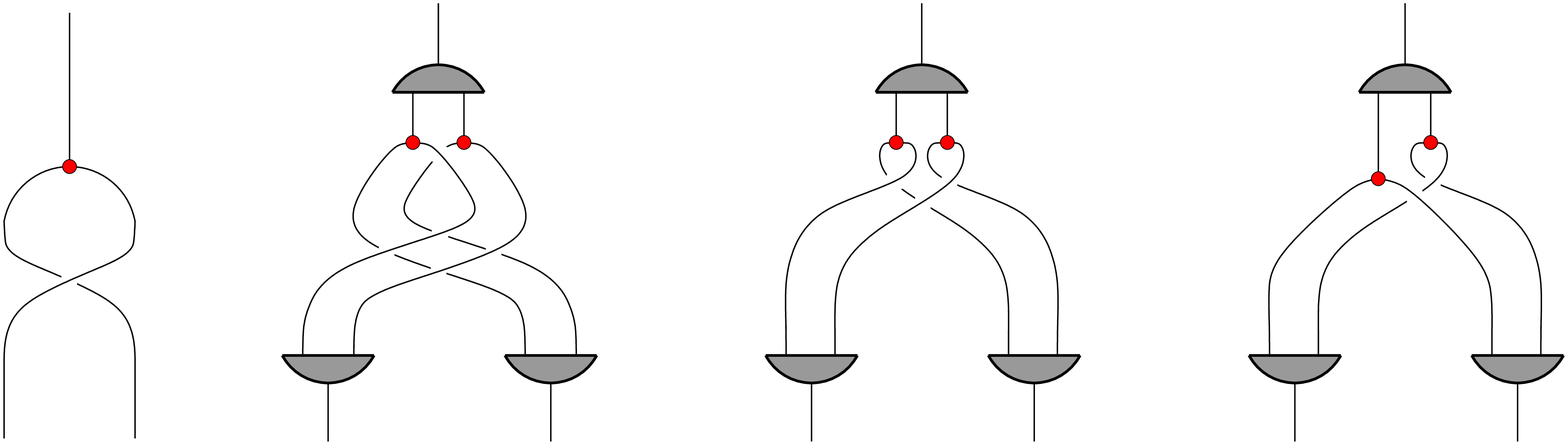}} \end{picture}
  \put(58,52.5)    {$=$}
  \put(-10.6,-9.2) {\sse$\efu Bl\AA$}
  \put(8.6,127.2)  {\sse$\efu Bl\AA$}
  \put(25.6,-9.2)  {\sse$\efu Bl\AA$}
  \put(79.4,-9.2)  {\sse$\efu Bl\AA$}
  \put(109.6,127.2){\sse$\efu Bl\AA$}
  \put(138.6,-9.2) {\sse$\efu Bl\AA$}
  \put(212.2,-9.2) {\sse$\efu Bl\AA$}
  \put(240.2,127.2){\sse$\efu Bl\AA$}
  \put(273.3,-9.2) {\sse$\efu Bl\AA$}
  \put(181,52.5)   {$=$}
  \put(311,52.5)   {$=$}
  \put(344.6,-9.2) {\sse$\efu Bl\AA$}
  \put(374.2,127.2){\sse$\efu Bl\AA$}
  \put(406.3,-9.2) {\sse$\efu Bl\AA$}
  }
  \end{picture} \nonumber\\[2.3em]{} \label{pic38}
  \\[-1.8em]{}\nonumber\end{eqnarray} 
The first step implements the definition \erf{eq:[B]A-alg} 
of the product on $\efu Bl\AA$, while in the second step the resulting ribbons
are deformed slightly. The third expression in \erf{pic38} is already
almost the multiplication of $E'$, except that the braiding 
$c_{A,A}$ must be removed and the braiding $c_{B,A}$ must be replaced by 
$c_{A,B}^{-1}$. This is achieved in two steps. First we use the equality
$r\eq r\cir P^l_\AA(B)$ to insert an idempotent $P^l_\AA(B)$ before the
restriction morphism and then carry out the moves displayed in figure 
\erf{pic70} backwards. This replaces $m \cir c_{A,A}$ by 
$m \cir (\id_A \oti \theta_\AA^{-1})$. After a further slight deformation
of ribbons one arrives at a graph for which the right ingoing leg is just 
given by the leftmost graph in \erf{eq:add-twist}.  Using the first 
equality in \erf{eq:add-twist} we then arrive at the last expression in 
\erf{pic38}, which is precisely the multiplication of $E'$.
\qed

\medskip\noindent
Proof of \Proposition {\bf \ref{prop:AB-alg}}:
\\[.2em]
We restrict our attention to the case of $\efu Bl\AA$.
For $\efu Br\AA$ the reasoning works in the same way.
\\[.3em]
(i)~\,\,The checks of the (co)associativity, (co)unit and Frobenius properties
all work by direct computation: After writing out the definition, one uses 
\Lemma \ref{lem:remove-P} to remove the projector on the `internal' 
$A$-ribbon. The (co)associativity, Frobenius and (co)unit relations then 
follow directly from the corresponding properties of $A$ and $B$.
\\[.3em]
(ii)~\,The check that symmetry of $B$ implies symmetry of $\efu Bl\AA$ can be 
performed by the same method as in (i). To see that commutativity of $B$ implies 
commutativity of $\efu Bl\AA$, first note that from \Lemma \ref{le:mE0=mE1}, 
$m_E\cir c_{E,E}\eq m_{E'}$. If $B$ is commutative, then $B'\eq B$ as a 
Frobenius \alg, so that \erf{mEE'} reduces to $m_E \cir c_{E,E}\eq m_{E}$.
\\[.3em]
(iii)~That $\efu Bl\AA$ is haploid follows from \Proposition \ref{alphaalpha}
by specialising to $U\eq B$ and $V\eq\one$, together with \erf{eq:HomAB=HomBA} 
and \Lemma \ref{lem:C=E} 
below, by which we have the bijections $\Hom(B,C_r(A)) 
\,{\cong}\, \Hom(B,\efu\one r\AA) \,{\cong}\, \Hom(\efu Bl\AA,\one)$.
(Note that here the assumption about the non-vanishing of the 
dimension of $\efu Bl\AA$ is not yet needed.)
\\
To see that $\efu Bl\AA$ is special, we can use lemma 3.11 of \cite{fuRs4},
according to which it suffices to show that the counit $\eps$ given in
\erf{eq:[B]A-alg} is a non-zero multiple of $\eps_\natural$ as defined in
\erf{eq:epsnat} (evaluated for the algebra $\efu Bl\AA$). Since $\efu Bl\AA$ 
is haploid, it is guaranteed that $\eps_\natural\eq\gamma\,\eps$ with 
$\gamma\iN\koerper$. The proportionality constant $\gamma$ can be determined 
by composing the equality with $\eta$; the result is $\gamma\eq\xi^{-1} 
\dim(\efu Bl\AA){/}\!\dim(A)\dim(B)$, which is non-zero by assumption.
\\[.3em]
(iv)~It is sufficient to check that the projectors on $A{\otimes}U$ are equal, 
i.e.\ $P^l_\AA(U) \eq P^r_\AA(U)$. Since $A$ is commutative and symmetric, 
and thus also has trivial twist, the desired equality can be rewritten as
   % [pic~35]
  \bea  \begin{picture}(180,99)(0,29)
  \put(0,0)   {\begin{picture}(0,0)(0,0)
              \scalebox{.38}{\includegraphics{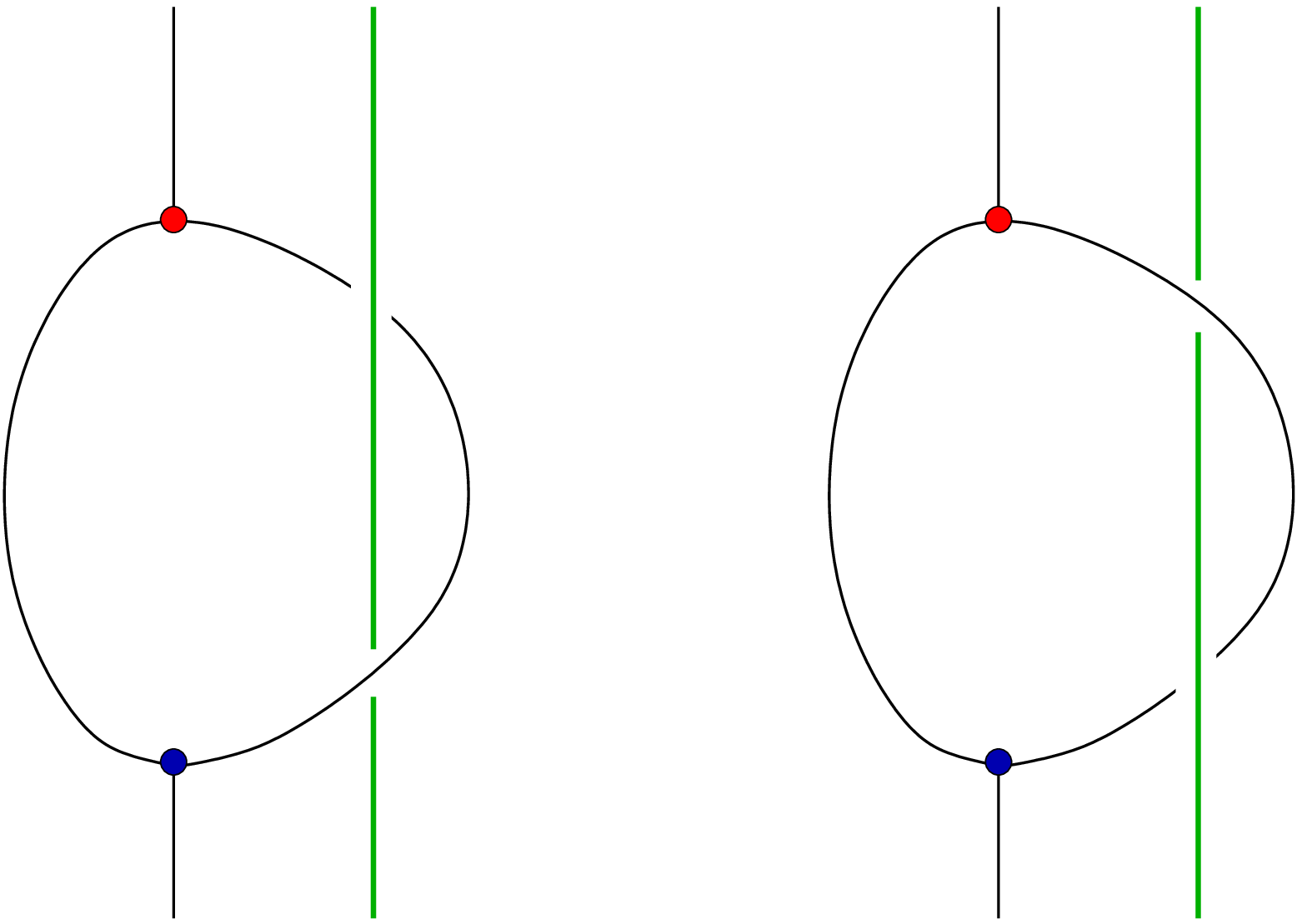}} \end{picture}}
  \put(18.8,-9.9)  {\sse$A$}
  \put(19.8,125.5) {\sse$A$}
  \put(45.5,-9.9)  {\sse$U$}
  \put(46.3,125.5) {\sse$U$}
  \put(82,54.1)    {$=$}
  \put(127.7,-9.9) {\sse$A$}
  \put(128.7,125.5){\sse$A$}
  \put(154.4,-9.9) {\sse$U$}
  \put(155.2,125.5){\sse$U$}
  \epicture20 \labl{Pl-r}
This latter equality, in turn, can be verified as follows.
First one deforms the $A$-loop on the \lhs\ of \erf{Pl-r} in such a
manner that the order of the braidings $c_{A,U}^{}$ and $c_{U,A}^{-1}$
gets interchanged, and then uses, consecutively, commutativity,
the Frobenius property, again commutativity, symmetry to obtain
   % [pic~69]
  \begin{eqnarray}  \begin{picture}(330,99)(48,29)
  \put(0,0)   {\begin{picture}(0,0)(0,0)
              \scalebox{.38}{\includegraphics{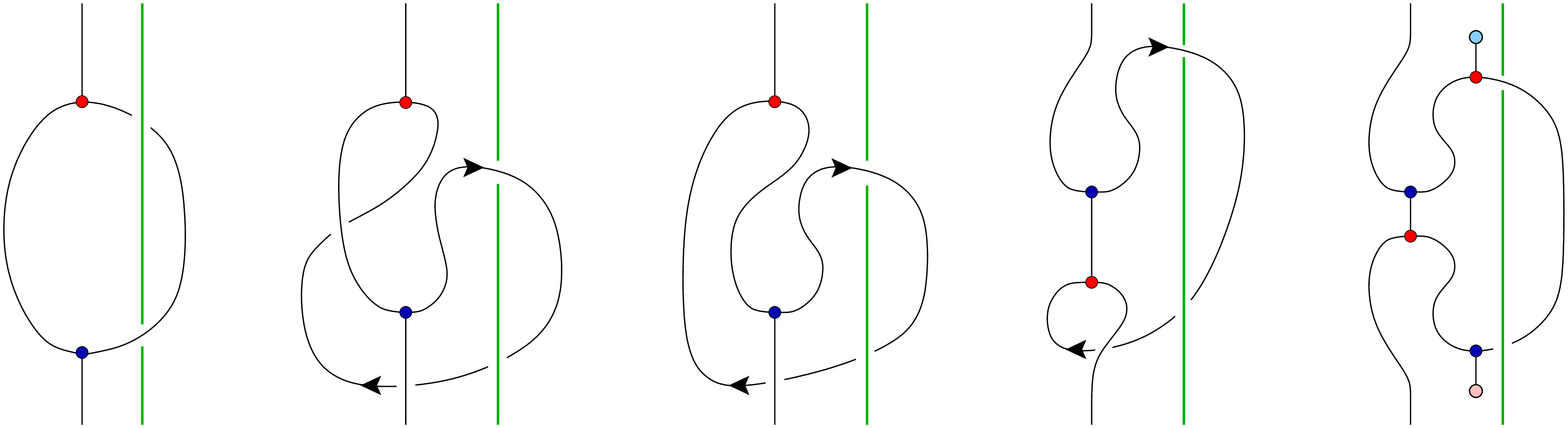}} \end{picture}}
  \put(18.4,-9.2)  {\sse$A$}
  \put(19.8,125.2) {\sse$A$}
  \put(36.4,-9.2)  {\sse$U$}
  \put(37.6,125.2) {\sse$U$}
  \put(67,56.5)    {$=$}
  \put(110.8,-9.2) {\sse$A$}
  \put(112.2,125.2){\sse$A$}
  \put(137.8,-9.2) {\sse$U$}
  \put(139.0,125.2){\sse$U$}
  \put(172,56.5)   {$=$}
  \put(215.8,-9.2) {\sse$A$}
  \put(217.2,125.2){\sse$A$}
  \put(243.1,-9.2) {\sse$U$}
  \put(244.3,125.2){\sse$U$}
  \put(278,56.5)   {$=$}
  \put(306.8,-9.2) {\sse$A$}
  \put(308.2,125.2){\sse$A$}
  \put(333.7,-9.2) {\sse$U$}
  \put(334.9,125.2){\sse$U$}
  \put(366,56.5)   {$=$}
  \put(397.5,-9.2) {\sse$A$}
  \put(398.9,125.2){\sse$A$}
  \put(425.2,-9.2) {\sse$U$}
  \put(426.4,125.2){\sse$U$}
  \end{picture} \nonumber\\[3.2em]{} \label{pic69}
  \\[-2.2em]{}\nonumber\end{eqnarray}
from which one arrives at the \rhs\ of \erf{Pl-r} by another (twofold) use of
the Frobenius property.
\\[.3em]
(v)~\,In addition to having $\efu Bl\AA\,{\cong}\,\efu Br\AA$ as objects
in $\cC$, one further verifies that $m_l\eq m_r$, $\eta_l\eq\eta_r$,
$\Delta_l\eq\Delta_r$ and $\eps_l\eq\eps_r$. Let us only show how to
check equality of $m_l$ and $m_r$; for the other morphisms similar
arguments apply. One considers the transformations (for better
readability we suppress the arrows indicating the duality morphisms)
  %% [pic~60]
  \bea  \begin{picture}(390,115)(15,17)
  \put(38,0)   {\begin{picture}(0,0)(0,0)
              \scalebox{.38}{\includegraphics{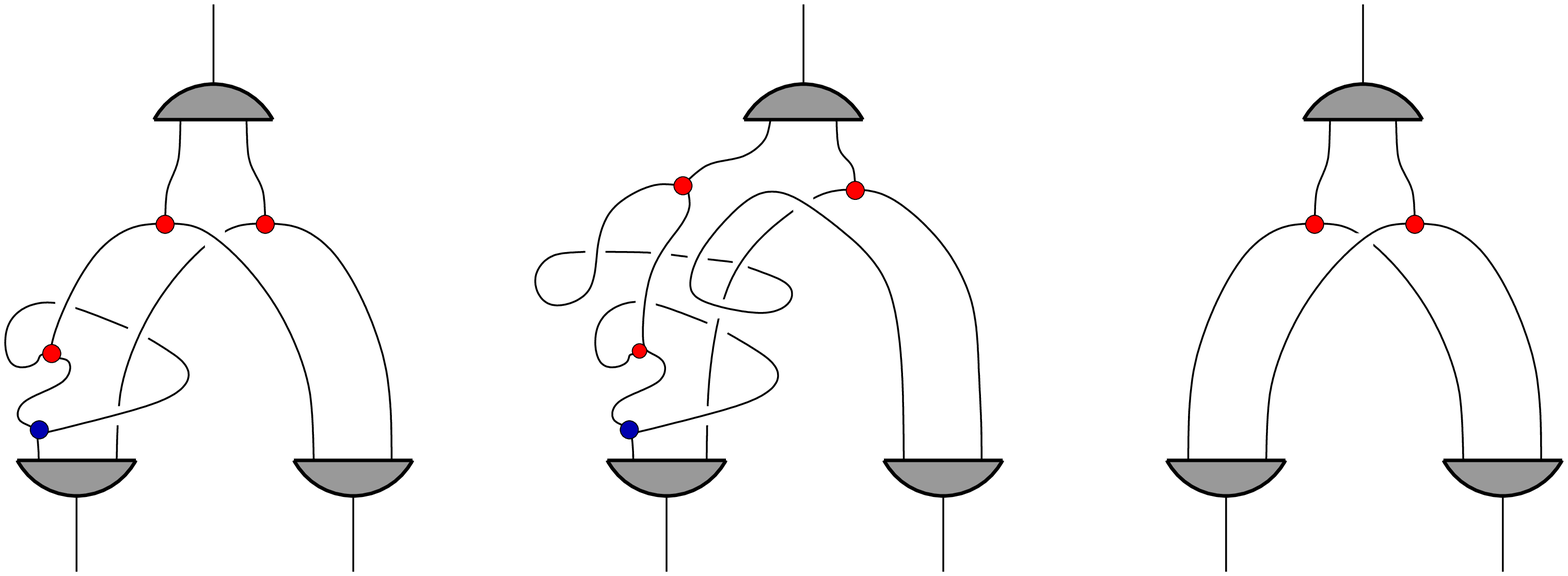}} \end{picture}}
  \put(0,51.5)     {$m_l\,=$}
  \put(42.8,-9.2)  {\sse$\efu Bl\AA$}
  \put(95.8,28.7)  {\sse$A$}
  \put(72.6,126.7) {\sse$\efu Bl\AA$}
  \put(101.9,-9.2) {\sse$\efu Bl\AA$}
  \put(113.6,62.2) {\sse$B$}
  \put(134,51.5)   {$=$}
  \put(168.8,-9.2) {\sse$\efu Bl\AA$}
  \put(197.6,126.7){\sse$\efu Bl\AA$}
  \put(227.3,-9.2) {\sse$\efu Bl\AA$}
  \put(263,51.5)   {$=$}
  \put(287.8,-9.2) {\sse$\efu Bl\AA$}
  \put(316.6,126.7){\sse$\efu Bl\AA$}
  \put(346.5,-9.2) {\sse$\efu Bl\AA$}
  \put(377,51.5)   {$=\,m_r$}
  \epicture13\labl{pic60}
In the first step the definition of $m_l$ is written out and
an idempotent $P^l_\AA(B)$ is inserted on top of an embedding
$\efu Bl\AA\,{\prec}\,A{\otimes}B$. Afterwards the multiplication
on $A$ is moved along the idempotent. In the third step the
$A$-ribbon is rearranged and the commutativity of $A$ is used.
\qed

%%%%%%%%%%%%%%%%%%%%%%%%%%%%%%%%%%%%%%%%%%%

\subsection{Centers and endofunctors}

{}From the \Definitions \ref{defLRC} and \ref{def:[]-functor} 
it is clear that the centers of an \alg\ can be interpreted as images of the
endofunctors $\EFU{l/r}\AA$, i.e.\ $C_{l/r}(A)\,{\cong}\,\efu\one{l/r}\AA$
as objects of $\cC$. We will now see that, upon endowing $C_{l/r}(A)$
with the structure of a Frobenius algebra inherited from $A$, and
$\efu\one{l/r}\AA$ with the Frobenius structure described in 
\Proposition \ref{prop:AB-alg}, 
this is even an isomorphism of Frobenius algebras.

\dtl{Lemma}{lem:C=E} 
For every \ssFA\ $A$ we have isomorphisms
  \be  C_l(A) \cong \efu\one l\AA
  \qquad {\rm and} \qquad C_r(A) \cong \efu\one r\AA  \labl{C1A}
as Frobenius algebras.
 % If $A$ is in addition simple, then $C_l$ and $C_r$ are
 % haploid commutative symmetric special Frobenius algebras.

\medskip\noindent
Proof: \\
{}From \Lemma \ref{centralidem} we know that $C_{l/r}$ is the image of the 
split idempotent $P^{l/r}_A$ defined in equation \erf{PE-def}. Also, comparison 
with the idempotents \erf{PU-def} shows that $P^{l/r}_A\eq P^{l/r}_A(\one)$. 
Thus $C_{l/r} \,{\cong}\, \efu\one{l/r}\AA$ as an object in $\cC$. Further, 
the definition of the algebra structure on $\efu B{l/r}\AA$ in 
\Proposition \ref{prop:AB-alg}(i) 
reduces to the one of $C_{l/r}$ (as given in equation 
\erf{eq:Clr-alg}) when $B\eq\one$.
\qed

This lemma can be used to establish the following more general result
for tensor product \alg s:
\\[-2.3em]

\dtl{Proposition}{prop:tensor-center}
(i)~\,For any pair $A,\,B$ of \ssFA s in a ribbon category $\cC$ one has
  \be
  C_l(A{\otimes}B) \,\cong\, \efu{C_l(B)}l\AA
  \qquad{\rm and}\qquad
  C_r(A{\otimes}B) \,\cong\, \efu{C_r(A)}r{\!B} \labl{eq:Cl-Cr}
as symmetric Frobenius algebras. 
\\[.3em]
(ii)~If in addition $\dim(C_r(A))\,{\ne}\,0$, $\dim(C_l(B))\,{\ne}\,0$ and
$\dim(C_{l/r}(A{\otimes}B))\,{\ne}\,0$, as well as
$\dim_\koerper\Hom(C_r(A),C_l(B))\eq1$, then 
$\efu{C_l(B)}l\AA$ and $\efu{C_r(A)}r{\!B}$ are haploid and special.

\medskip\noindent
Proof:\\
(i)~\,Let us start with the second relation in \erf{eq:Cl-Cr}. The following
series of equalities shows that the braiding $(c_{A,B})^{-1}$ relates
the idempotents \erf{PU-def} for $C_r(A{\otimes}B) \,{\cong}\,
\efu\one r{A{\otimes}B}$ and for $\efu{C_r(A)}rB$:
  %% [pic~15]
  \bea  \begin{picture}(350,105)(0,38)
  \put(0,0)  {\begin{picture}(0,0)(0,0)
              \scalebox{.38}{\includegraphics{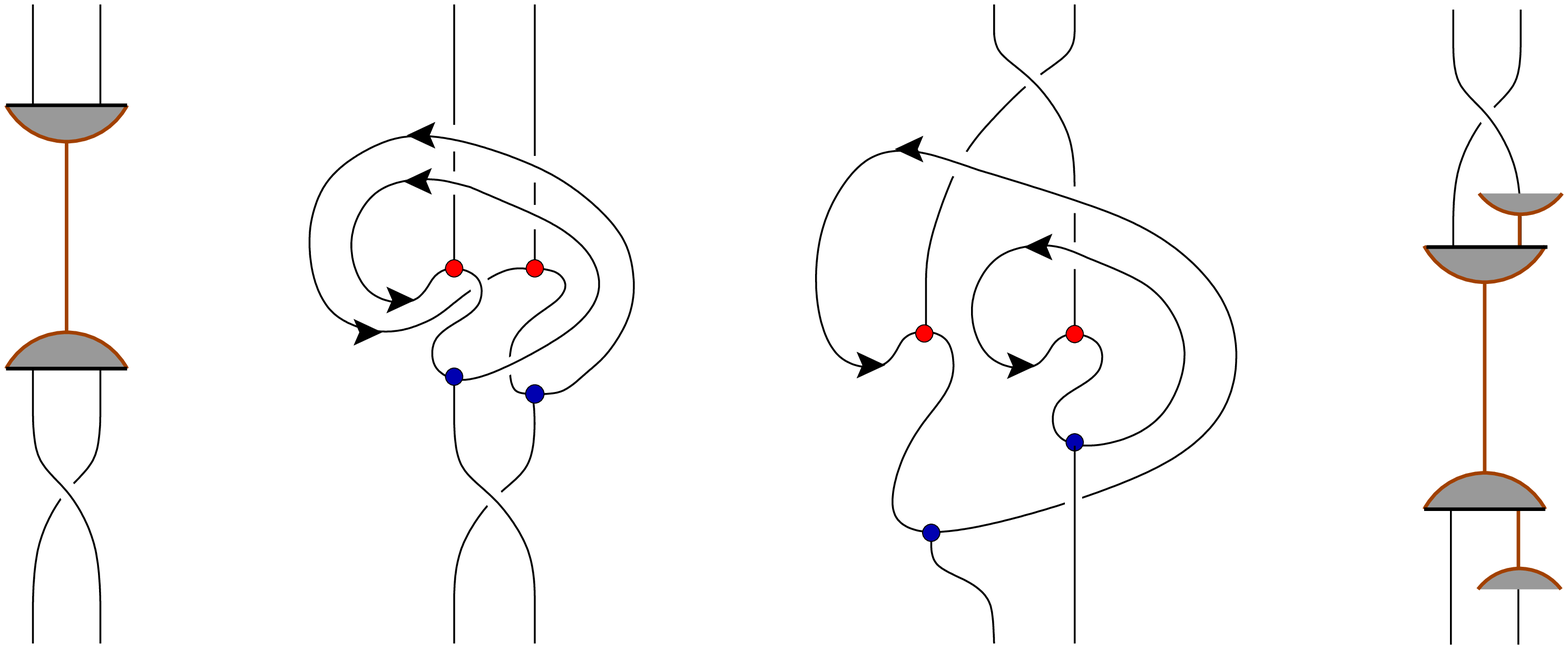}} \end{picture}}
  \put(2.2,-9.2)    {\sse$B$}
  \put(2.2,137.7)   {\sse$A$}
  \put(-24.1,82.8)  {\sse$C_r(A{\otimes}B)$}
  \put(16.8,-9.2)   {\sse$A$}
  \put(17.2,137.7)  {\sse$B$}
  \put(41.1,66.8)   {$=$}
  \put(90.2,-9.2)   {\sse$B$}
  \put(90.9,137.7)  {\sse$A$}
  \put(107.2,-9.2)  {\sse$A$}
  \put(107.6,137.7) {\sse$B$}
  \put(146.1,66.8)  {$=$}
  \put(202.2,-9.2)  {\sse$B$}
  \put(202.2,137.7) {\sse$A$}
  \put(219.5,-9.2)  {\sse$A$}
  \put(219.9,137.7) {\sse$B$}
  \put(273.1,66.8)  {$=$}
  \put(298.1,-9.2)  {\sse$B$}
  \put(298.1,137.7) {\sse$A$}
  \put(310.4,54.4)  {\sse$\efu{C_r(A)}rB$}
  \put(312.6,-9.2)  {\sse$A$}
  \put(312.6,137.7) {\sse$B$}
  \put(317.4,19.9)  {\sse$C_r(A)$}
  \put(319.5,84.8)  {\sse$C_r(A)$}
  \epicture26 \labl{eq:pic15}
Define the morphisms $\varphi\iN\Hom(\efu{C_r(A)}r{\!B},C_r(A{\otimes}B))$
and $\psi\iN\Hom(C_r(A{\otimes}B),\efu{C_r(A)}r{\!B})$ by
  %% [pic~12]
  \bea  \begin{picture}(210,80)(0,35)
  \put(35,0)  {\begin{picture}(0,0)(0,0)
              \scalebox{.38}{\includegraphics{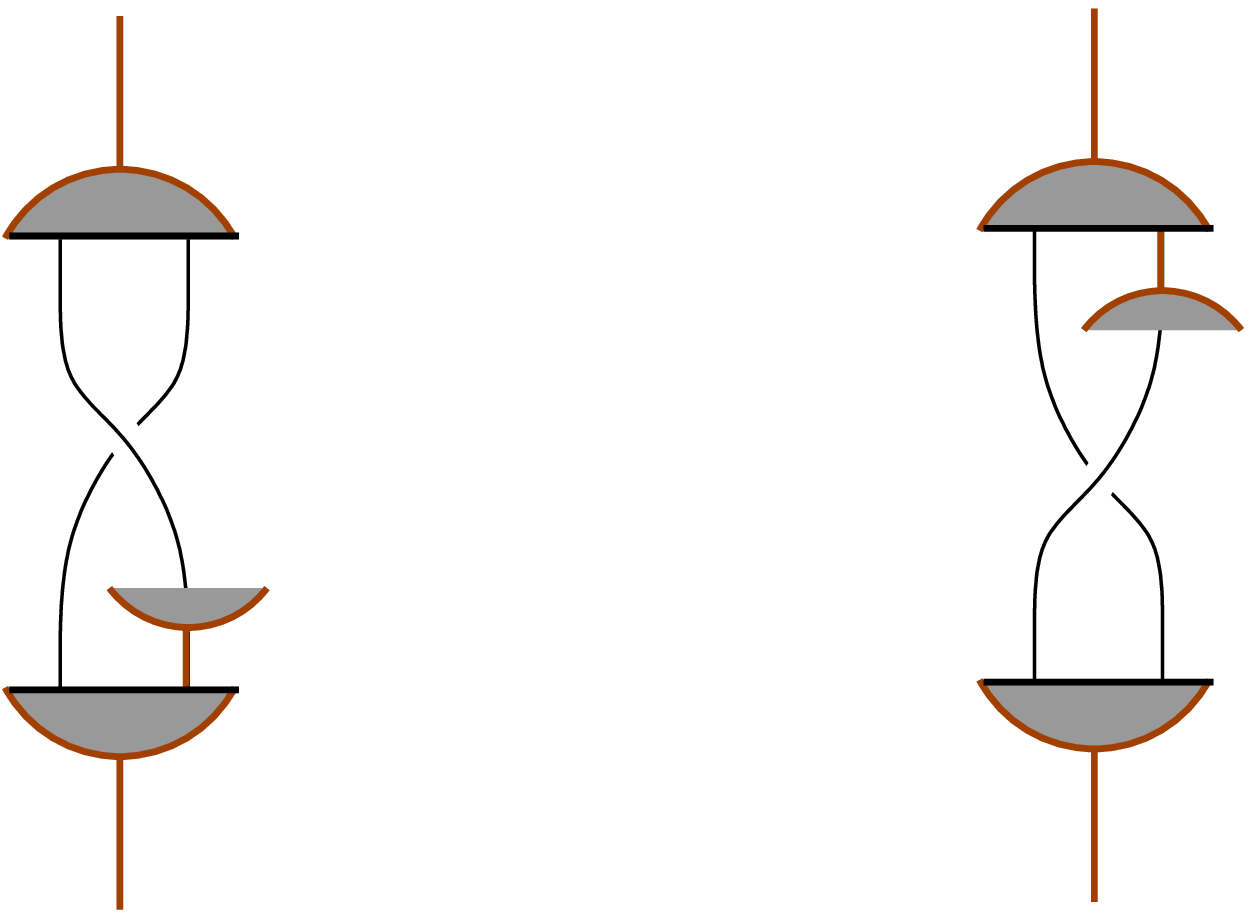}} \end{picture}}
  \put(0,47.8)      {$\varphi\;:=$}
  \put(30.4,-9.2)   {\sse$\efu{C_r(A)}rB$}
  \put(29.4,103.3)  {\sse$C_r(A{\otimes}B)$}
  \put(33.4,30.6)   {\sse$B$}
  \put(54.9,43.4)   {\sse$A$}
  \put(58.1,26.7)   {\sse$C_r(A)$}
  \put(110,47.8)    {$\psi\;:=$}
  \put(137.4,-9.2)  {\sse$C_r(A{\otimes}B)$}
  \put(137.4,103.3) {\sse$\efu{C_r(A)}rB$}
  \put(141.1,30.6)  {\sse$A$}
  \put(163.4,30.6)  {\sse$B$}
  \put(167.1,69.2)  {\sse$C_r(A)$}
  \epicture27 \labl{pic12}
Using \erf{eq:pic15} one can verify that
  \be
  \varphi \circ \psi = \id_{C_r(A{\otimes}B)} \qquad {\rm and} \qquad
  \psi \circ \varphi = \id_{\efu{C_r(A)}rB} \,.  \ee
Next we would like to see that $\varphi$ is compatible with the symmetric 
special Frobenius structure of the two algebras. We need to check that
  \be\bearll
  \varphi^{-1} {\circ}\, \eta_{C_r(A{\otimes}B)}^{}
  = \eta_{\efu{C_r(A)}rB}^{} \,,  \qquad &
  \varphi^{-1} {\circ}\, m_{C_r(A{\otimes}B)}^{} \cir
  (\varphi{\otimes}\varphi) = m_{\efu{C_r(A)}rB}^{} \,,
  \\{}\\[-.7em]
  \eps_{C_r(A{\otimes}B)}^{} {\circ}\, \varphi
  = \eps_{\efu{C_r(A)}rB}^{}  \,,  &
  (\varphi^{-1}{\otimes}\varphi^{-1}) \cir 
  \Delta_{C_r(A{\otimes}B)}^{} \cir \varphi = \Delta_{\efu{C_r(A)}rB}^{} \,.
  \eear\labl{eq:phi-is-hom}
The relations for $\eta$ and $\eps$ are immediate when inserting the 
definitions \erf{eq:[B]A-alg} 
and \erf{eq:Clr-alg}. Using again \erf{eq:pic15}, for the multiplication we find
  %% [pic~14]
  \bea  \begin{picture}(320,115)(0,44)
  \put(0,0)  {\begin{picture}(0,0)(0,0)
              \scalebox{.38}{\includegraphics{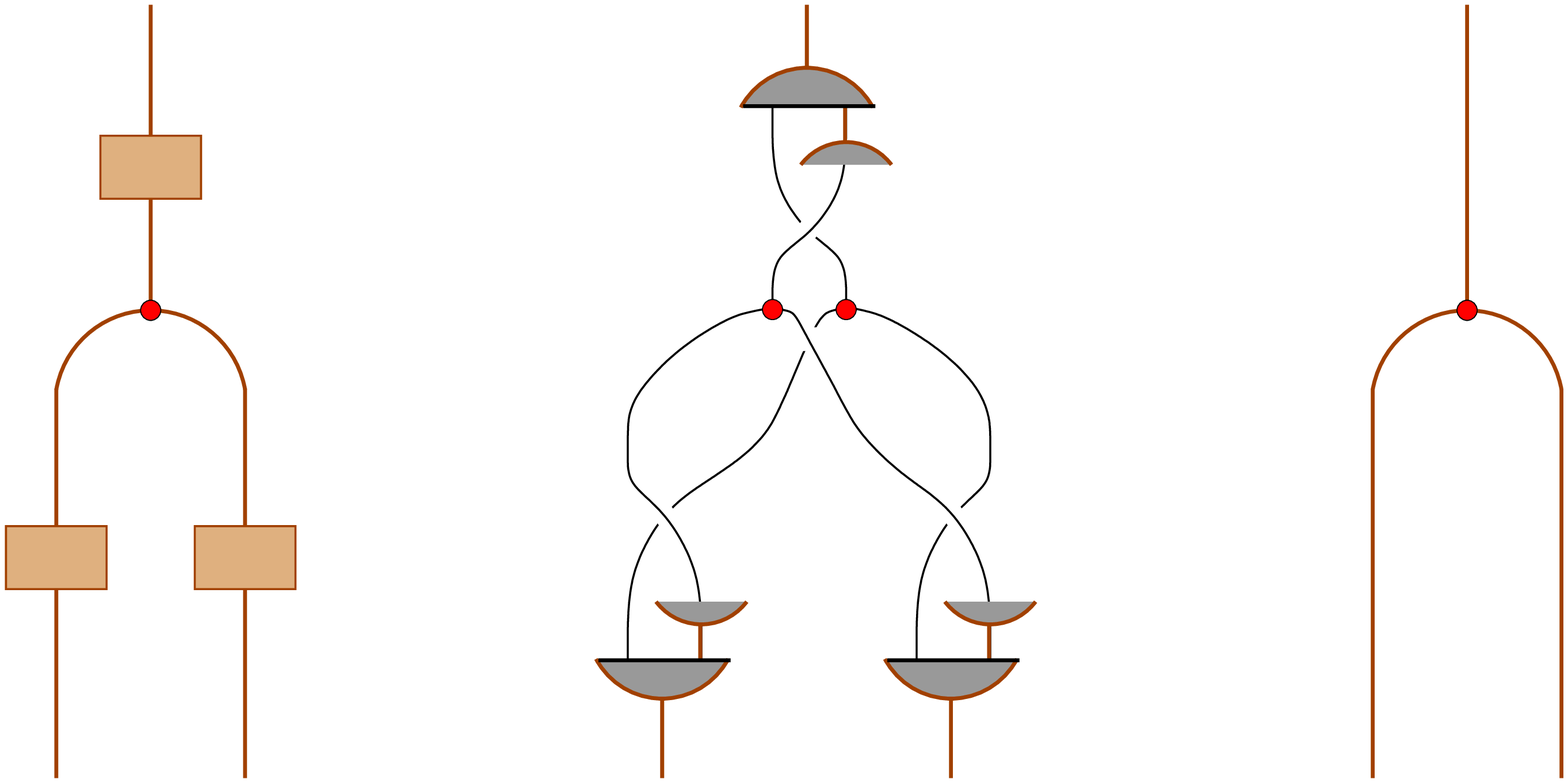}} \end{picture}}
  \put(-13.4,-9.2)  {\sse$\efu{C_r(A)}rB$}
  \put(7.5,41.8)    {\sse$\varphi$}
  \put(14.4,154.5)  {\sse$\efu{C_r(A)}rB$}
  \put(21.1,115.2)  {\sse$\varphi^{-1}_{}$}
  \put(33.4,-9.2)   {\sse$\efu{C_r(A)}rB$}
  \put(43.4,83.2)   {\sse$C_r(A{\otimes}B)$}
  \put(43.1,41.8)   {\sse$\varphi$}
  \put(79.1,55.8)   {$=$}
  \put(108.4,-9.2)  {\sse$\efu{C_r(A)}rB$}
  \put(112.4,39.2)  {\sse$B$}
  \put(133.3,39.2)  {\sse$A$}
  \put(136.1,154.5) {\sse$\efu{C_r(A)}rB$}
  \put(136.6,24.9)  {\sse$C_r(A)$}
  \put(139.6,96.6)  {\sse$A$}
  \put(161.8,96.6)  {\sse$B$}
  \put(163.6,121.9) {\sse$C_r(A)$}
  \put(164.6,-9.2)  {\sse$\efu{C_r(A)}rB$}
  \put(168.4,39.2)  {\sse$B$}
  \put(189.3,39.2)  {\sse$A$}
  \put(192.6,24.9)  {\sse$C_r(A)$}
  \put(219.9,55.8)  {$=$}
  \put(237.6,-9.2)  {\sse$\efu{C_r(A)}rB$}
  \put(262.4,154.5) {\sse$\efu{C_r(A)}rB$}
  \put(287.4,-9.2)  {\sse$\efu{C_r(A)}rB$}
  \epicture34 \labl{pic14}
The corresponding relation for the comultiplication
in \erf{eq:phi-is-hom} is demonstrated similarly.
\\[.2em]       
The proof of the first relation in \erf{eq:Cl-Cr} works along
the same lines, but this time $\varphi$ and $\psi$ take the easier form
  %% [pic~13]
  \bea  \begin{picture}(210,78)(0,34)
  \put(35,0)  {\begin{picture}(0,0)(0,0)
              \scalebox{.38}{\includegraphics{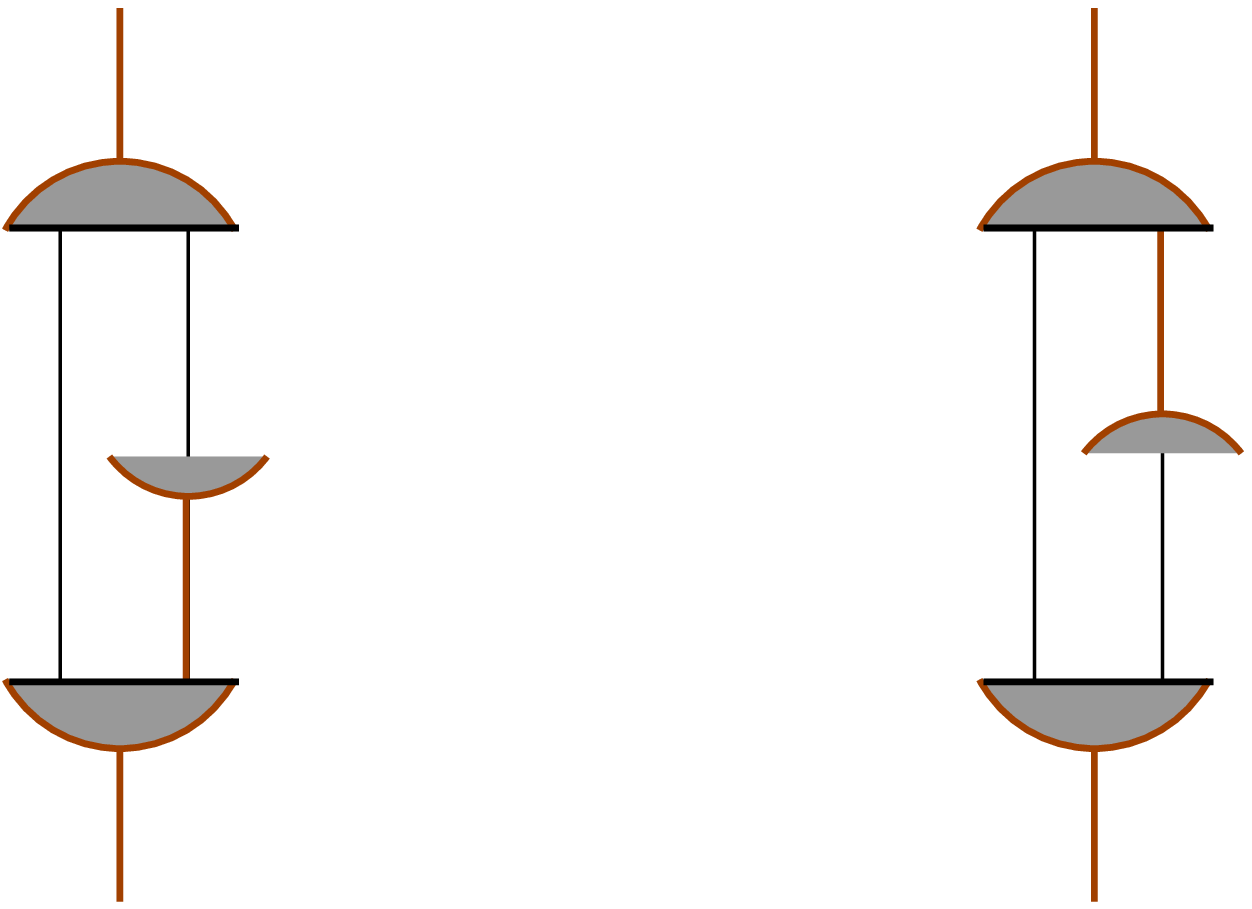}} \end{picture}}
  \put(1,47.8)      {$\varphi\;=$}
  \put(30.4,-9.2)   {\sse$\efu{C_l(B)}l\AA$}
  \put(29.4,103.3)  {\sse$C_l(A{\otimes}B)$}
  \put(33.4,32.5)   {\sse$A$}
  \put(58.6,61.4)   {\sse$B$}
  \put(58.6,32.5)   {\sse$C_l(B)$}
  \put(111,47.8)    {$\psi\;=$}
  \put(137.4,-9.2)  {\sse$C_l(A{\otimes}B)$}
  \put(137.4,103.3) {\sse$\efu{C_l(B)}l\AA$}
  \put(140.8,32.5)  {\sse$A$}
  \put(165.5,32.5)  {\sse$B$}
  \put(165.5,61.4)  {\sse$C_l(B)$}
  \epicture22 \labl{pic13}
Correspondingly there is no braiding in the analogue of \erf{eq:pic15}.
\\[.2em]
By \Propositions \ref{lem:C=[1]i}(i) and \ref{prop:AB-alg}(ii),
$E^l_A(C_l(B))$ and $E^r_B(C_r(A))$ are symmetric. 
\\[.3em]
(ii)~By \Proposition \ref{lem:C=[1]i}(i), 
$C_r(A)$ and $C_l(B)$ are commutative symmetric Frobenius \alg s. Further, 
since by \Remark \ref{prop:ssFA-unique}(vi)
the tensor unit is a retract of every Frobenius algebra, the condition 
$\dim_\koerper\Hom( C_r(A),C_l(B))\eq1$ implies in particular that $C_r(A)$ 
and $C_l(B)$ are haploid and thus simple. Since their dimensions are non-zero 
by assumption, \Proposition \ref{lem:C=[1]i}(iii) 
then tells us that the two 
centers are also special. Together with the assumptions $\dim(C_{l/r}
(A{\otimes}B)) \,{\ne}\,0$ and $\dim_\koerper\Hom(C_r(A),C_l(B))\eq1$, 
as well as the isomorphisms of part (i), we can finally apply 
\Proposition \ref{prop:AB-alg}(ii) and \ref{prop:AB-alg}(iii) to see that 
$\efu{C_r(A)}r{\!B}$ is haploid and special. 
\\
Similarly, again by \Proposition \ref{prop:AB-alg}(ii) 
and \ref{prop:AB-alg}(iii), this time together with the bijection 
\erf{eq:HomAB=HomBA}, $\efu{C_l(B)}l\AA$ is haploid and special as well.
\qed

%%%%%%%%%%%%%%%%%%%%%%%%%%%%%%%%%%%%%%%%%%%

\subsection{The ribbon subcategory of local modules}\label{sec:local-modules}

As noticed after \erf{Pf-def}, in symmetric tensor categories the objects 
$\efu U{l/r}A$ are closely related to induced modules over the {\em center\/}
of the algebra $A$. We therefore now consider categories of modules over
{\em commutative\/} symmetric Frobenius algebras. As it turns out, this is 
still appropriate in the genuinely braided case. Note that according to 
\Proposition \ref{prop:AB-alg}(iv), for commutative $A$,
there is only one endofunctor $\Efu\cdot A\,{\equiv}\,\efu\cdot{l/r}A$.

The relevant class of modules is introduced in
\\[-2.3em]

\dtl{Definition}{def:loc-mod}
A left module $M\eq(\M,\r_M)$ over a commutative \ssFA{} $A$
in a ribbon category $\cC$ is called {\em local\/} iff
  \be  \r_M^{} \cir P_\AA(\M) = \r_M^{} \,, \labl{eq:def-loc}
where $P_\AA(\M)\,{\equiv}\,P_\AA^{l/r}(\M)$ is the idempotent defined
in \erf{PU-def}.

\medskip

As we will see in \Proposition \ref{pro:ostr} below, our concept of locality 
is equivalent to the one of \Definition 3.2 of \cite{kios}. The latter, which
says that $M$ is local iff $\r_M^{}\cir c_{\M,A}^{}\cir c_{A,\M}^{}\eq\r_M^{}$,
had been introduced earlier for modules over an algebra in a general braided
\tc\ in \cite{pare23}, where such modules were termed {\em dyslectic\/}.
A main motivation for the introduction of dyslectic modules in \cite{pare23} was
the fact that they form a full subcategory that can be naturally endowed with a
  tensor structure and a braiding.
This property will be crucial in the present context, too.
Here we prefer the qualification ``local'' to the term ``dyslectic'' because
it agrees with the standard use \cite{scya5} in conformal quantum field theory 
in the context of so-called simple current extensions (compare 
\Remark \ref{+sicu}(ii) below).

To show the equivalence between the two characterisations (as well as a
third one) we first give the
\\[-2.3em]

\dtl{Lemma}{ab(c)lemma}
For $M\eq(\M,\r)$ a left module over a commutative \ssFA\ in a ribbon category
the morphism
  %% [pic~39]
  \bea  \begin{picture}(175,88)(0,22)
  \put(50,0)  {\begin{picture}(0,0)(0,0)
              \scalebox{.38}{\includegraphics{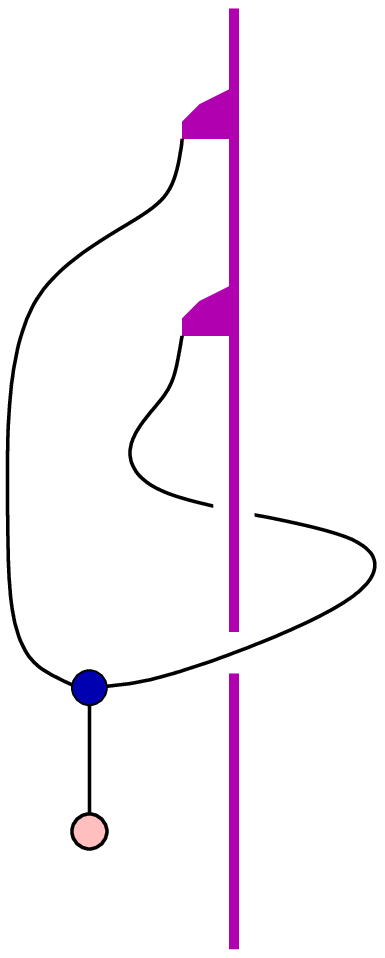}} \end{picture}}
  \put(0,51.5)     {$Q_M\,:=$}
  \put(43.6,40.2)  {\sse$A$}
  \put(70.6,-8.9)  {\sse$\M$}
  \put(71.6,107.7) {\sse$\M$}
  \put(77.7,69.7)  {\tiny$\r_{\!M}^{}$}
  \put(77.7,91.7)  {\tiny$\r_{\!M}^{}$}
  \put(109.5,51.5) {$\in\,\Hom(\M,\M)$}
  \epicture07 \labl{eq:QM-def}
satisfies\\[-.7em]
  \be \mbox{\hspace{-5.5em}} {\rm (i)} \quad\,
  \r_M^{}\cir P_{\!A}(\M) = \r_M^{}\cir(\id_\AA^{}{\otimes}Q_M^{})
  \qquad\ {\rm and}\qquad\ {\rm (ii)} \quad
  Q_M^{}\cir Q_M^{} = Q_M^{} \,. \labl{abc}
\mbox{$\ $}\\[-1.4em]
\mbox{$\ $}\hspace{.8em}(iii)\,\, $M$ is local iff $\,Q_M\eq\id_\M$.

\bigskip
\noindent
Proof:\\
To see (i) we note that
  %% [pic~40]
  \bea  \begin{picture}(265,86)(0,22)
  \put(90,0)  {\begin{picture}(0,0)(0,0)
              \scalebox{.38}{\includegraphics{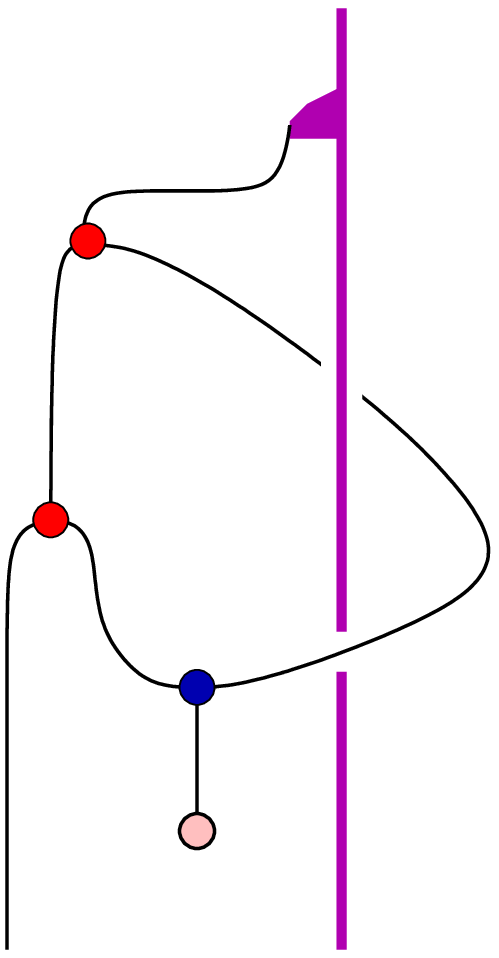}} \end{picture}}
  \put(0,51.5)     {$\r_M \cir P_\AA(\M) \,=$}
  \put(85.6,-8.2)  {\sse$A$}
  \put(121.6,-8.9) {\sse$\M$}
  \put(122.6,107.7){\sse$\M$}
  \put(160,51.5)   {$=\, \r_M \cir (\id_\AA \oti Q_M) \,.$}
  \epicture07 \labl{pic40}
The equality on the left uses the fact that $A$ is commutative and has
trivial twist (recall the corresponding comment in \Section \ref{sec:alg-fun}),
as well as the Frobenius property.
The equality on the right is obtained by applying the representation
property first for the upper and then on the lower of the two products.
\\[.3em]
To deduce (ii) we combine (i) with the fact that $P_\AA(\M)$ is an idempotent,
and insert (i) twice (using also that $\id_\AA \oti Q_M$ commutes
with $P_\AA(\M)$), so as to arrive at
  \be \bearll
  \rho_M \circ (\id_\AA \oti Q_M) \!\! & = \rho_M \circ P_\AA(\M)
  \\{}\\[-.8em] &
  = \rho_M \circ P_\AA(\M) \circ P_\AA(\M) = \rho_M \circ (\id_\AA \oti Q_M)
  \circ (\id_\AA \oti Q_M) \,.  \eear \labl{qqq}
Property (ii) now follows by composing both sides of \erf{qqq} with
$\eta\oti\id_M$.
\\[.3em]
(iii)~For local $M$ (\ref{abc}(i)) reads $\r_M \cir(\id_\AA\oti Q_M)\eq\r_M$,
which when composed with $\eta\oti\id_\M$ yields $Q_M\eq\id_\M$. Conversely, 
inserting $Q_M\eq\id_\M$ into (i) yields the defining property of locality.
\qed

\medskip

We are now in a position to present
\\[-2.5em]

\dtl{Proposition}{pro:ostr}
For a left module $M\eq(\M,\r)$ over a commutative \ssFA{} in a
ribbon category $\cC$ the following conditions are equivalent: \\[.3em]
(i)~~~$M$ is local.
\\[4pt]
(ii)\,~~$\theta_\M\iN\HomA(M,M)$\,.
\\[4pt]
(iii)~~$\r\cir c_{\M,A}\cir c_{A,\M} \eq\r$\,.

     \newpage
\medskip\noindent
Proof:\\
(i)\,$\Rightarrow$\,(ii)\,: We start with the equalities
  %% [pic~74]
  \bea  \begin{picture}(330,188)(57,4)
  \put(0,0)  {\begin{picture}(0,0)(0,0)
             \scalebox{.38}{\includegraphics{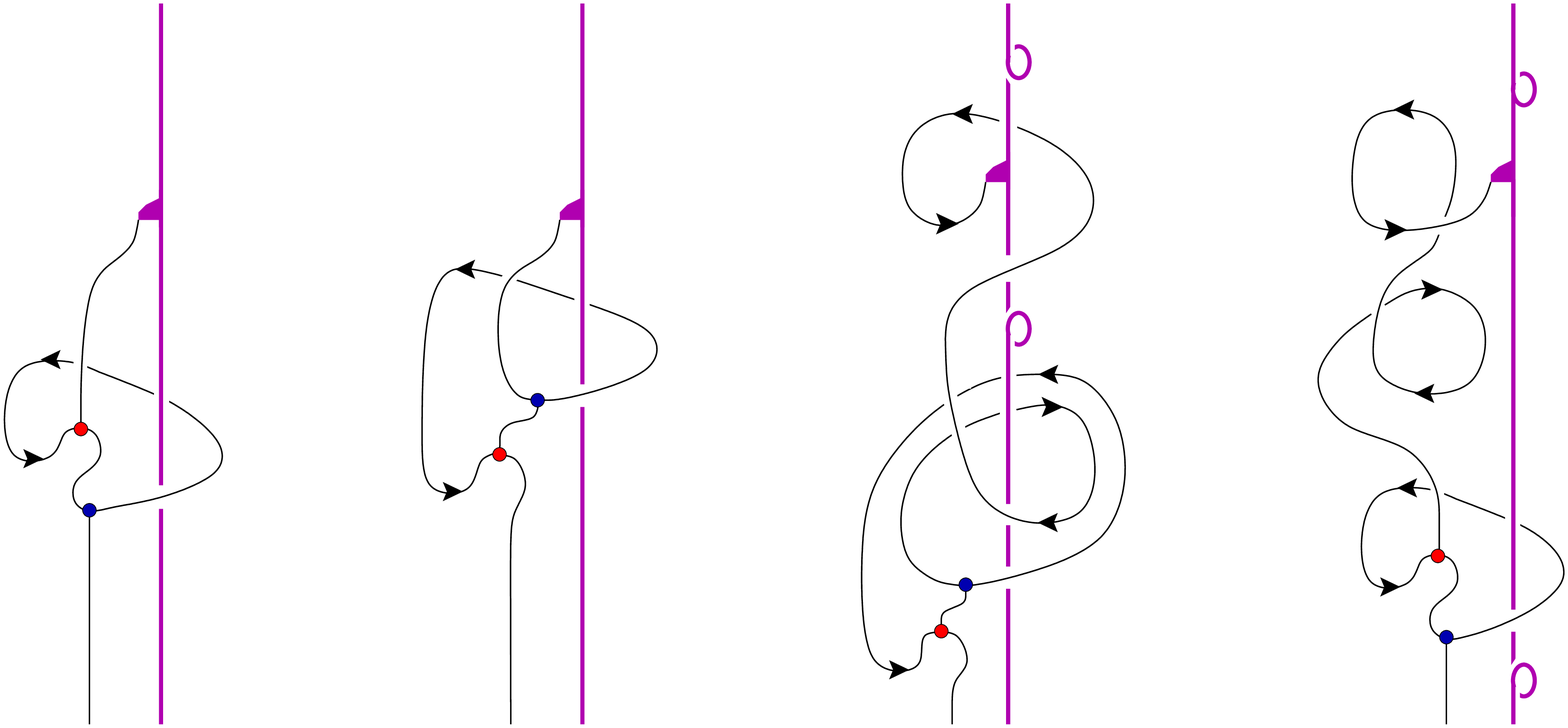}} \end{picture}}
  \put(17.8,-8.9)  {\sse$A$}
  \put(35.1,-9.5)  {\sse$\M$}
  \put(35.7,186.6) {\sse$\M$}
  \put(74,86.5)    {$=$}
  \put(123.8,-8.9) {\sse$A$}
  \put(141.2,-9.5) {\sse$\M$}
  \put(141.8,186.6){\sse$\M$}
  \put(191,86.5)   {$=$}
  \put(234.7,-8.9) {\sse$A$}
  \put(249.4,-9.5) {\sse$\M$}
  \put(250.1,186.6){\sse$\M$}
  \put(262.1,100.5){\sse$\theta_{\!\M}$}
  \put(262.1,166.2){\sse$\theta_{\!\M}^{-1}$}
  \put(304,86.5)   {$=$}
  \put(359.7,-8.9) {\sse$A$}
  \put(375.9,-9.5) {\sse$\M$}
  \put(376.6,186.6){\sse$\M$}
  \put(389.9,14.5) {\sse$\theta_{\!\M}$}
  \put(389.9,159.2){\sse$\theta_{\!\M}^{-1}$}
  \epicture-1 \labl{pic74}
The first equality uses that $A$ is Frobenius, the second combines the identity
  %% [pic~73]
  \bea  \begin{picture}(135,94)(0,37)
  \put(0,0)  {\begin{picture}(0,0)(0,0)
             \scalebox{.38}{\includegraphics{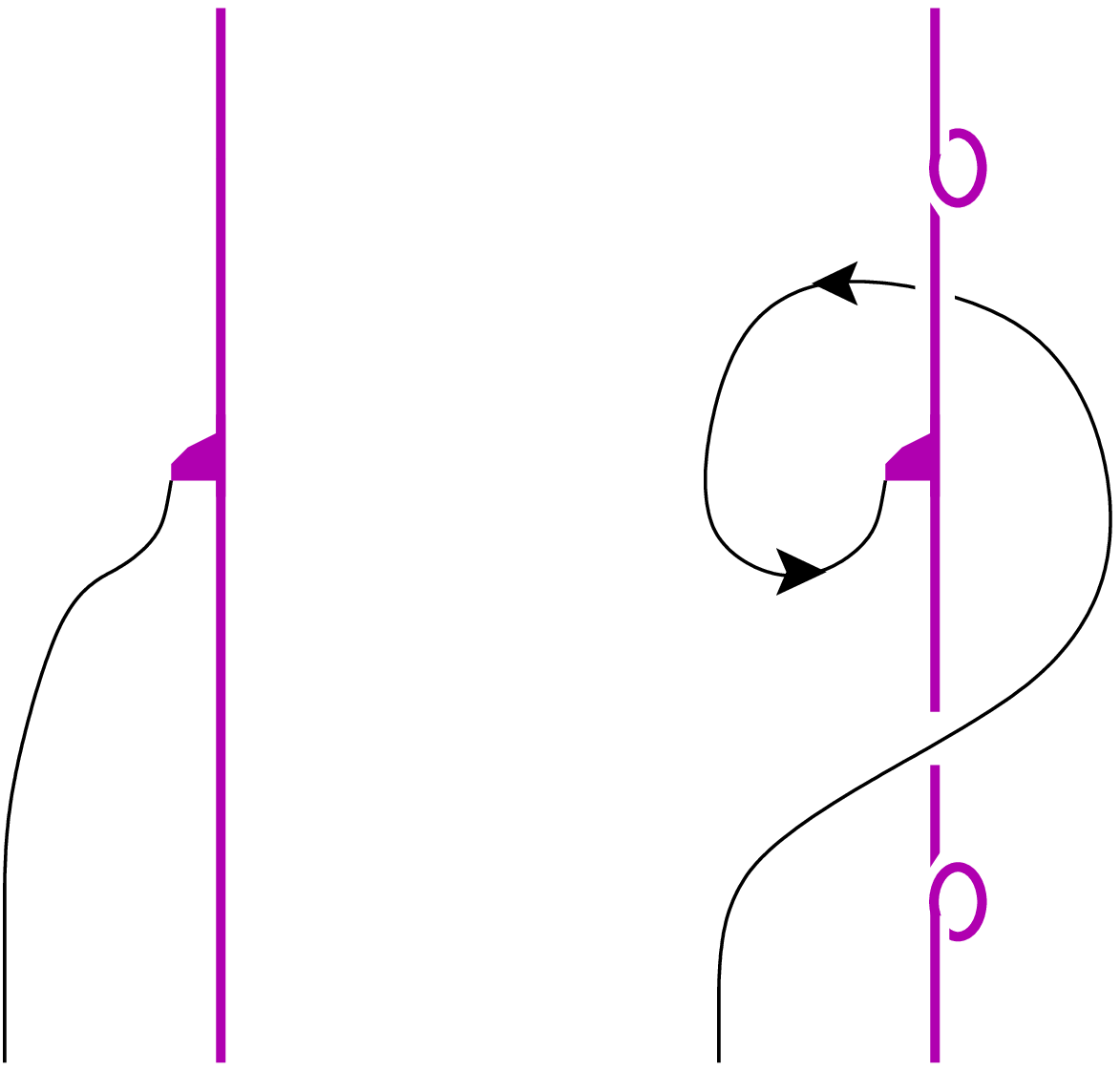}} \end{picture}}
  \put(-3.5,-8.9)  {\sse$A$}
  \put(20.3,-9.5)  {\sse$\M$}
  \put(21.5,126.3) {\sse$\M$}
  \put(27.7,70.1)  {\sse$\r_{\!M}^{}$}
  \put(52,58.5)    {$=$}
  \put(78.8,-8.9)  {\sse$A$}
  \put(102.9,-9.5) {\sse$\M$}
  \put(104.1,126.3){\sse$\M$}
  \put(109.4,70.1) {\sse$\r_{\!M}^{}$}
  \put(115.5,18.2) {\sse$\theta_{\!\M}^{}$}
  \put(115.5,102.8){\sse$\theta_{\!\M}^{-1}$}
  \epicture22 \labl{pic73}
with moves similar to those in figure \erf{pic37}, and the third combines a 
deformation of the upper $A$-ribbon with an application of the Frobenius 
property in the lower part of the graph.  On the \rhs\ of \erf{pic74}, we 
can in addition straighten the upper $A$-ribbon. Afterwards, by composing
the left and \rhs s with $\theta_\M$ from the top and removing the idempotent 
$P_\AA(\M)$ (as allowed by locality) we arrive at the statement 
that $\theta_\M \iN \End(\dot M,\dot M)$ is actually in $\End_A(M,M)$.
\\[.3em]
(ii)\,$\Rightarrow$\,(i)\,:
By $\theta_\AA\eq\id_\AA$ and the
compatibility between braiding and twist we have
  \be
  \r_M \circ c_{\M,A}^{} \circ c_{A,\M}^{} =
  \theta_\M \circ \r_M  \circ (\id_\AA \oti \theta_\M^{-1}) \,.
  \labl{eq:ostr-1-2}
To show that $Q_M$ is the identity morphism, we insert this relation into the 
definition \erf{eq:QM-def} of the morphism $Q_M$. Then by (ii)
we can take the lower \rep\ morphism $\r_M$ past $\theta_\M$ without
introducing any braiding or twist. Using that $A$ is  special, the
$A$-ribbon can then be removed, resulting in $Q_M\eq\id_\M$. By 
\Lemma \ref{ab(c)lemma}(iii) it follows that $M$ is local.
\\[.3em]
(ii)\,$\Leftrightarrow$\,(iii)\,
follows immediately from relation \erf{eq:ostr-1-2}.
\\
(For semisimple $\cC$ this equivalence is \Theorem 3.4.1 of \cite{kios}.)
\qed

\medskip

In applications one is often interested in {\em simple\/} modules.
Therefore we separately state the following result which makes it easy
to test if a simple module is local.
\\[-2.3em]

\dtl{Corollary}{cor:local}
For a simple module $M$, with $\dim(M)\,{\ne}\,0$, over a commutative 
\ssFA{} in a ribbon category $\cC$ the following conditions are equivalent: 
\\[.3em]
(i)~~~$M$ is local. \\[.3em]
(ii)\,~~$\tr\, \theta_\M \ne 0\,$.\\[.3em]
(iii)~~$\theta_\M = \xi_M\,\id_\M\ $ for some $\ \xi_M\iN\kx$.

\medskip\noindent
Proof:\\
We first note that when $M$ is simple, then the morphism $Q_M\iN\Hom(\M,\M)$
given by \erf{eq:QM-def} satisfies
  %% [pic~41]
  \bea  \begin{picture}(195,85)(0,26)
  \put(50,0)  {\begin{picture}(0,0)(0,0)
              \scalebox{.38}{\includegraphics{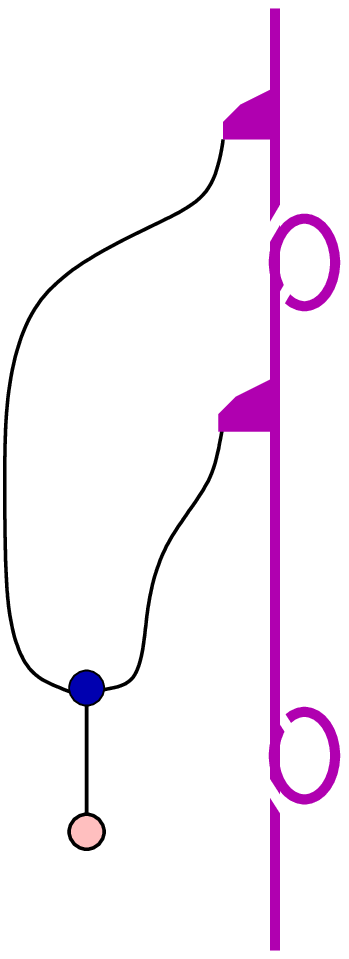}} \end{picture}}
  \put(-3,51.5)    {$Q_M \,=$}
  \put(43.6,40.2)  {\sse$A$}
  \put(75.1,-8.9)  {\sse$\M$}
  \put(76.1,107.7) {\sse$\M$}
  \put(104,51.5)   {$=\, \dsty\frac{\tr\,\theta_\M}{\dim(\M)}
                     \;\theta_\M^{-1}\,.$}
  \epicture13 \labl{(ab)c}
To get the first equality, the $\M$-ribbon is twisted so as to remove
the braidings; because of $\theta_\AA\eq\id_\AA$ the resulting twist of the
$A$-ribbon can be left out. The second equality is a consequence of 
\Lemma 4.4 of \cite{fuRs4}. By definition (see \Definition 4.3 of \cite{fuRs4})
of the $A$-{\em averaged morphism\/}
$\theta^{\rm av}_\M$, which is an element of
$\Hom_A(M,M)$, the graph is just $\theta^{\rm av}_\M\cir\theta_\M^{-1}$.
Since $M$ is simple, $\Hom_A(M,M)$ is one-dimensional, so that
$\theta^{\rm av}_\M$ is proportional to $\id_\M$. The constant of
proportionality is determined by comparing the traces, resulting in the
final expression in \erf{(ab)c}.
\\[.3em]
(i)\,$\Rightarrow$\,(ii)\,: By \Lemma \ref{ab(c)lemma}(iii), 
for local $M$ we have $Q_M\eq\id_\M$. By \erf{(ab)c} this, in turn, means that
  \be  \Frac{\tr\,\theta_\M}{\dim(\M)}\, \theta_\M^{-1} = \id_\M \,.  \ee
Since $\id_\M$ is invertible, this requires $\tr\,\theta_\M$ to be non-zero.
\\[.3em]
(ii)\,$\Rightarrow$\,(iii)\,:
The equality obtained by inserting \erf{(ab)c} into the projection property
(\ref{abc}(ii)) can hold only if $\tr\,\theta_\M\eq0$ or if
  \be
  \Frac{\tr\,\theta_\M}{\dim(\M)} \, \id_\M = \theta_\M \,.  \ee
Since by (ii) the first possibility is excluded, we arrive at (iii) with
$\xi_M\eq\tr\,\theta_\M{/}\!\dim(\M)$ (which is non-zero).
\\[.3em]
(iii)\,$\Rightarrow$\,(i)\,:
When combined with \erf{(ab)c}, the statement (iii) implies $Q_M\eq\id_\M$.
Together with (\ref{abc}(iii)) it then follows that $M$ is local.
\qed

\dtl{Remark}{+sicu}
(i)~\,\,When $\cC$ is semisimple, it follows immediately from the definition
that in the decomposition of a local module $M$ into simple modules
$M_\kappa$ all the $M_\kappa$ are local as well.
\\[.3em]
(ii)~\,The case when the commutative algebra $A$ is a direct sum of
{\em invertible\/} simple objects is known in the physics literature as a
{\em simple current extension\/}. Then the local modules $M$
are those for which the `monodromy charge' with respect to $A$ vanishes,
which means that for all simple subobjects $J$ of $A$ and all simple
subobjects $U_i$ of $M$ one has the equality $s_{J,U_i}\eq s_{\one,U_i}$,
where $s$ is the $s$-matrix defined in \erf{sUV}. Precisely these modules 
appear in the chiral conformal field theory obtained by a simple current 
extension \cite{scya6}.
\\[.3em]
(iii)~For $\cC\eq\REP_{\rm DHR}(\mathfrak C)$ the category of DHR 
superselection sectors \cite{dohr,doro2} of a local rational quantum field 
theory $\mathfrak C$, there is a bijection between finite index extensions 
$\mathfrak C_{\rm ext}\,{\supseteq}\,\mathfrak C$ and symmetric special 
Frobenius algebras $A$ in $\REP_{\rm DHR}(\mathfrak C)$, and 
$\mathfrak C_{\rm ext}$ is again a local quantum field theory iff 
$A$ is commutative \cite{lore}.
\\
For the case that $\cC\eq\REP(\mathfrak V)$ is the category of modules
over a rational vertex algebra $\mathfrak V$ with certain nice properties, 
the fact that $\Ext{\mathcal C}{\REP(\mathfrak V)}$ is equivalent to 
$\REP(\mathfrak V_{\rm ext})$, with $\mathfrak V_{\rm ext}$ the vertex 
algebra for the extended \cft, has been observed in \cite{kios} (\Theorem 5.2).

\dtl{Definition}{def:ext-cat}
Let $A$ be a commutative \ssFA{} in a ribbon category $\cC$. The
{\em category of local $A$-modules\/}, denoted by \calcal, is
the full subcategory of \calca\ whose objects are local $A$-modules.

\medskip

Under suitable conditions on $\cC$ and $A$, the category \calcal\
inherits various structural properties from $\cC$, such as being braided tensor 
(\Theorem 2.5 of \cite{pare23}) or modular (\Theorem 4.5
of \cite{kios}). We collect some of these properties in
\\[-2.3em]

\dtl{Proposition}{thm:mod}
For every commutative symmetric special 
Frobenius \alg\ $A$ in a ribbon category $\cC$ the following holds:
\\[.3em]
(i)~~\,\,\calcal\ is a ribbon category.
\\[.3em]
(ii)~\,If $\cC$ is semisimple, then \calcal\ is semisimple.
If $\cC$ is closed under direct sums and subobjects, then \calcal\ is
closed under direct sums and subobjects.
\\[.3em]
(iii)~If $\cC$ is modular and if $A$ is in addition simple,
then \calcal\ is modular.

\medskip\noindent
The proof is a straightforward combination of the results contained in the
proofs of \Theorems 1.10, 1.17 and 4.5 of \cite{kios} 
(which are derived in a semisimple setting and with $A$ assumed to be haploid,
but are easily adapted to the present framework, using in particular the fact 
that simple commutative algebras are also haploid) and the permanence properties
established in \Section 5 of \cite{fuSc16}. 
(For the simple current case that was mentioned in \Remark \ref{+sicu}(ii) 
above, see also \cite{fusS6,brug2,muge7}.) 

\medskip
   
Let us describe the tensor structure of \calcal\ in some detail.
For any algebra $A$, one defines the tensor product $M\OtA N$ of a right 
$A$-module $M$ and a left $A$-module $N$ as the cokernel of the morphism
$\r_M\oti\id_N \,{-}\, \id_M\oti\r_N$, provided that the cokernel exists.
In the present context, i.e.\ for $A$ a commutative \ssFA\ and
$M$ and $N$ two local left $A$-modules, the tensor product can 
conveniently be described as the image
  \be
  M \otA N := \Im P_{M{\otimes}N}  \labl{eq:2-tensor}
of a suitable idempotent in $\End(\M\Oti\dot N)$, provided that this idempotent
is split. The idempotent in question is given by 
(compare lemma 1.21 of \cite{kios})
  %% [pic~46]
  \bea  \begin{picture}(220,50)(0,28)
  \put(63,0)  {\begin{picture}(0,0)(0,0)
              \scalebox{.38}{\includegraphics{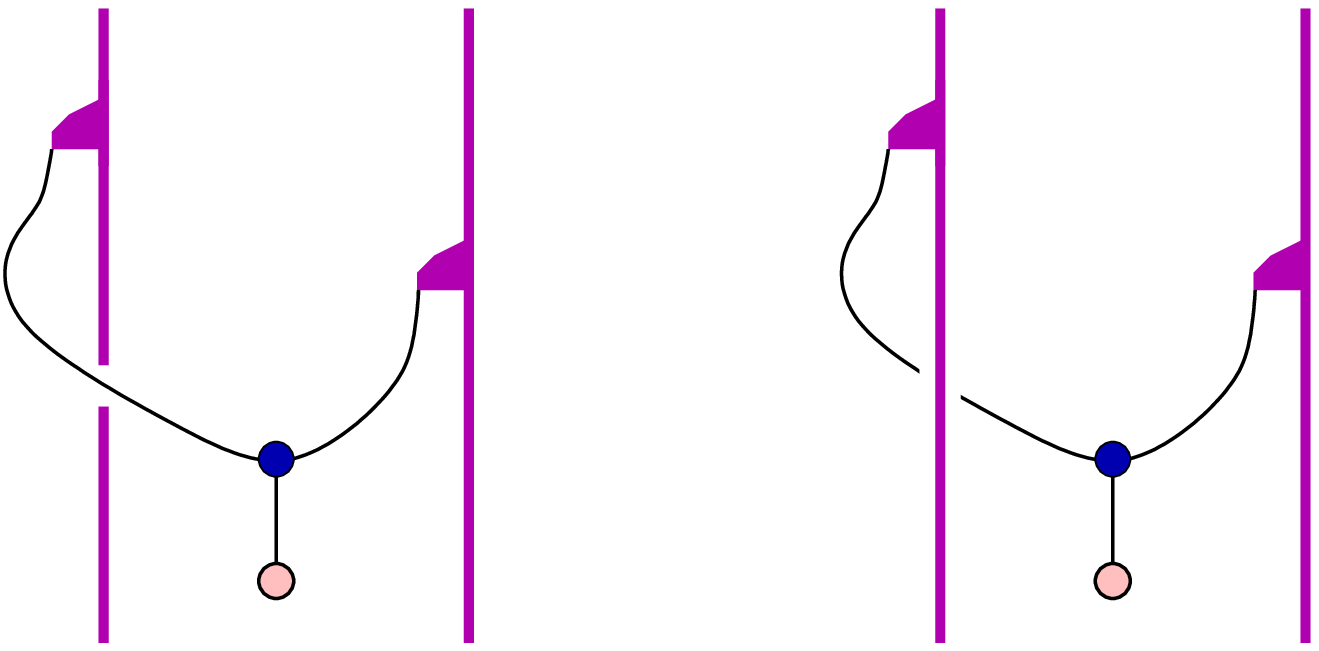}} \end{picture}}
  \put(0,30.5)    {$P_{M{\otimes}N} \,= $}
  \put(69.5,-8.9) {\sse$\M$}
  \put(69.9,74.9) {\sse$\M$}
  \put(85.1,27.3) {\sse$A$}
  \put(110.6,-8.9){\sse$\dot N$}
  \put(110.9,74.9){\sse$\dot N$}
  \put(133,30.5)  {$=$}
  \put(160.7,-8.9){\sse$\M$}
  \put(161.2,74.9){\sse$\M$}
  \put(176.1,27.3){\sse$A$}
  \put(201.9,-8.9){\sse$\dot N$}
  \put(202.2,74.9){\sse$\dot N$}
  \epicture17 \labl{eq:P-2-tensor}
(Owing to \Proposition \ref{pro:ostr}(iii), 
applied to the \rep\ morphism $\r_M$
for the local module $M$, the morphisms given by the left and right pictures
are equal.) Similarly, multiple tensor products can then be described as 
images of the idempotents
  %% [pic~47]
  \bea  \begin{picture}(250,50)(0,29)
  \put(82,0)  {\begin{picture}(0,0)(0,0)
              \scalebox{.38}{\includegraphics{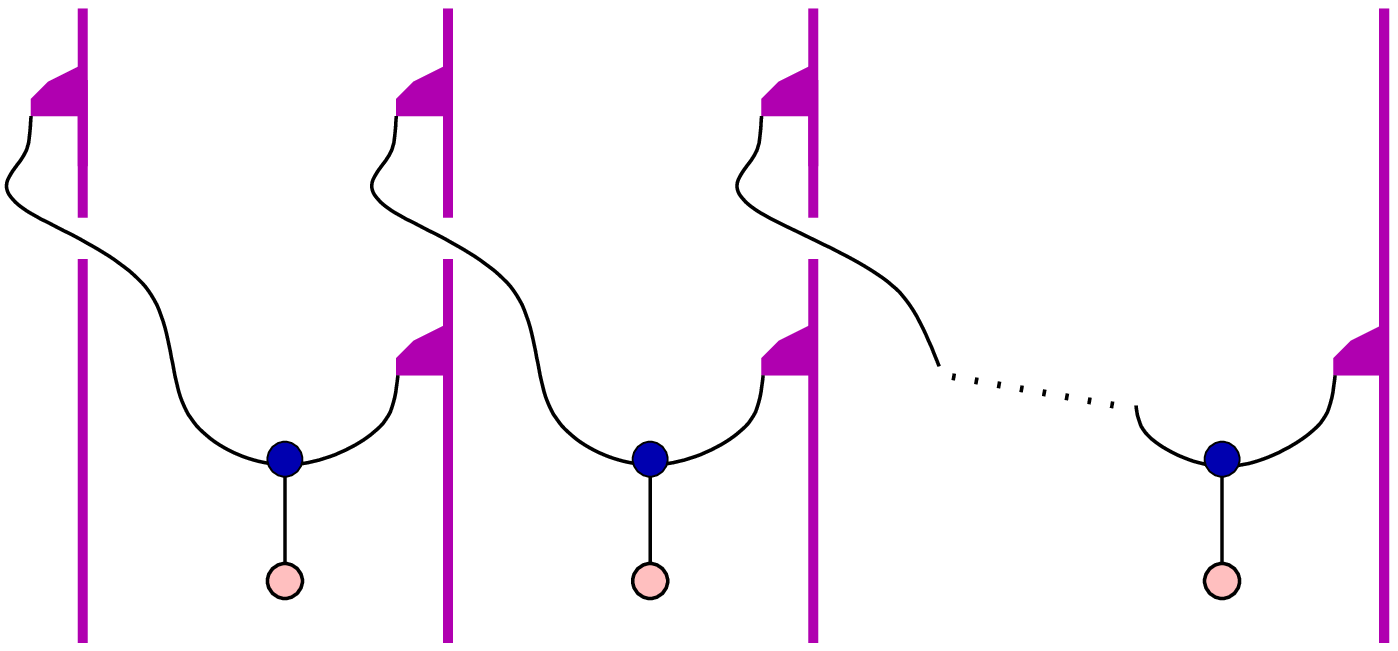}} \end{picture}}
  \put(0,30.5)    {$P_{M_1\otimes\cdots\otimes M_k} \,=$}
  \put(85.4,-8.9) {\sse$\M_1$}
  \put(85.7,74.9) {\sse$\M_1$}
  \put(104.1,27.3){\sse$A$}
  \put(125.7,-8.9){\sse$\M_2$}
  \put(125.9,74.9){\sse$\M_2$}
  \put(144.6,27.3){\sse$A$}
  \put(228.7,-8.9){\sse$\M_k$}
  \put(228.9,74.9){\sse$\M_k$}
  \epicture16 \labl{eq:mult-tens}
Note that this way of defining multiple tensor products is consistent with
the iterative application of \erf{eq:2-tensor}. Indeed one easily verifies 
that the idempotents $P_{(M\otimes_\AA N)\otimes K}$ and $P_{M\otimes
(N\otimes_\AA K)}$ are both equal to $P_{M\otimes N \otimes K}$.

Finally, denoting by $e_{M_1\otimes\cdots\otimes M_k}$ and $r_{\!M_1\otimes
\cdots\otimes M_k}$ the embedding and restriction morphisms for the idempotent
\erf{eq:mult-tens}, the tensor product of morphisms $f_i\iN\HomA(M_i,N_i)$
($i\eq1,2,...\,,k$) takes the form
  \be
  f_1 \otA \cdots \otA f_k =
  r_{\!N_1\otimes\cdots\otimes N_k}^{} \circ (f_1\,\Oti\cdots\Oti\,f_k )
  \circ e_{M_1\otimes\cdots\otimes M_k}^{} \,.  \labl{eq:morph-Atensor}

The definition \erf{eq:2-tensor} of the tensor product is based on the 
assumption that the idempotents $P_{M_1\otimes\cdots\otimes M_k}$ are split, 
for which it is sufficient that $\cC$ is Karoubian. If we do not impose 
Karoubianness of $\cC$, it can happen that $P_{M_1\otimes\cdots\otimes M_k}$ is 
not split; then we must work with the Karoubian envelope of \calcal\ and define
  \be
  M_1 \otA \cdots \otA M_k
  := (M_1\,\Oti\cdots\Oti\, M_k ; P_{M_1\oti\cdots\oti M_k}) \,.
  \labl{eq:multi-tensor-obj}
If $\cC$ {\em is\/} Karoubian so that we can define the tensor product as an 
image, we still must select $M\OtA N$ as a specific object in its isomorphism 
class (recall that we use the axiom of choice to regard images as objects).
We make this choice in a way compatible with \erf{eq:multi-tensor-obj}. 
With this definition of the tensor product the associativity constraints of
the category \calcal\ are, just as the ones of $\cC$, identities. However,
in general $A\otA M$ and $M$ are different objects of $\cC_\AA$
so that the unit constraints are non-trivial. In particular, the module
category is in general {\em not\/} a strict tensor category.

The ribbon structure of \calcal\ is inherited in a rather obvious manner
from $\cC$. Concretely, the {\em braiding\/} on \calcal\ is given by the family
  \be
  \cloc_{M,N} := r \circ c_{\M,\dot N}^{} \circ e
  \;\in \HomA(M\OtA N, N\OtA M)  \labl{eq:Cloc-braid}
of morphisms, where $e$ is the embedding morphism for the retract $M\OtA N 
\,{\prec}\, \M{\otimes}\dot N$, $c_{\M,\dot N}$ is the braiding in $\cC$, and
$r$ the restriction morphism for $\dot N{\otimes}\M \,{\succ}\, N\OtA M$.
The {\em twist\/} on \calcal\ just coincides with the one of $\cC$, i.e.\
$\theta^A_M\eq \theta^{}_{\M}$ (see \Proposition \ref{pro:ostr}), 
and the {\em duality\/} of \calcal\ is the assignment
$M\,{\mapsto}\,M^\vee \eq (\M^\vee,(d_\M\Oti\id_{\M^\vee}) \cir
(\id_{\M^\vee}\Oti\r_M\Oti\id_{\M^\vee}) \cir (c^{-1}_{\M^\vee,A}\Oti b_\M))$ 
together with the morphisms 
  \be \!\! \bearll
  b^A_M := r^{\phantom|}_{\!M\OtA M^\vee} \cir (\r_M\oti\id_{\M^\vee})
    \cir (\id_\AA\oti b_\M)
  = \r_{M\OtA M^\vee} {\circ}\,
    [\id_A \oti (r^{}_{\!M\OtA M^\vee}{\circ}\, b_\M )]
  &{\rm and} \\{}\\[-.5em]
  d^A_M := [\id_A \oti (d_\M \cir (\id_{\M^\vee} \oti \r_M) 
    \cir (c^{-1}_{\M^\vee,A} \oti \id_M) ]
    \cir [ (\Delta\cir\eta) \oti e^{}_{M^\vee\OtA M} ]
  \\{}\\[-.85em] \quad\;\ \
  = [\id_A \oti (d_\M\cir e^{}_{M^\vee\OtA M}\cir\r^{\phantom|}_{M^\vee\OtA M})]
    \cir [ (\Delta\cir\eta) \oti \id_{M^\vee\OtA M} ]
   \eear \labl{bA,dA}
(compare \Theorem 1.15 of \cite{kios} and 
   section 5.3 
of \cite{fuSc16}).

\dtl{Lemma}{neuir}
For $A$ a simple commutative special Frobenius algebra 
in a ribbon \cat\ $\cC$ and $A$-modules $M,N\iN\Obj(\cC_\AA)$ one has
  \be  \dim(M\otA N) = \frac{\dim(\M)\,\dim(\dot N)}{\dim(A)} \,.  \ee

\medskip\noindent
Proof:\\
We have
  \bea  \begin{picture}(330,40)(15,27)
  \put(225,0) {\begin{picture}(0,0)(0,0)
              \scalebox{.38}{\includegraphics{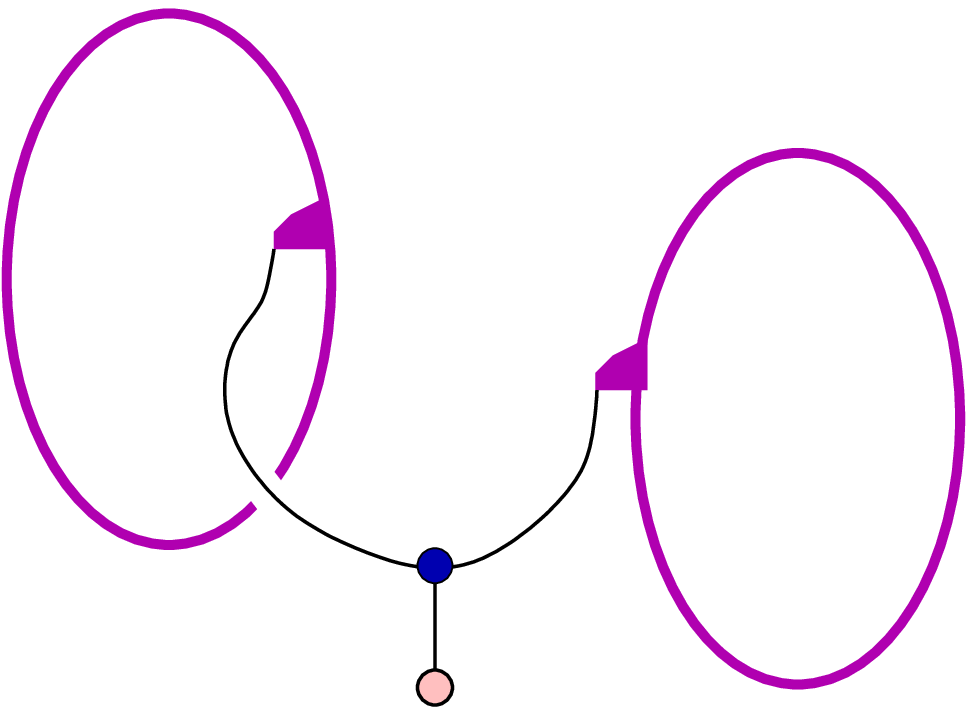}} \end{picture}}
  \put(-12,35.5)   {$\dim(M{\otimes_\AA} N) \,=\, \tr(\id_{M\otimes_\AA N})
                  \,=\, \tr(P_{M\otimes N}) \,=$}
  \put(227.3,55) {\sse$\M$}
  \put(320.2,16) {\sse$\dot N$}
  \epicture08 \labl{trPMN}
Now since $A$ is haploid, for every $\varphi\in\Hom(\one,A)$ we have
$\varphi\eq\beta_\one^{-1}(\eps{\circ}\varphi)\,\eta\eq(\eps{\circ}\varphi)\,
\eta/\dim(A)$. It follows that removing the $A$-lines from the graph on the
\rhs\ of \erf{trPMN} just amounts to a factor of $1{/}\!\dim(A)$; but
removing the $A$-ribbons leaves us just with an $\M$- and an $\dot N$-loop,
i.e.\ with $\,\dim(\M)\,\dim(\dot N)$.
\qed

\medskip

When $A$ is symmetric, this result is also implied by \Lemma \ref{beibl},
and for the case that in addition $\cC$ is semisimple, it has already been
established in \cite{kios} (corollary to \Theorem 1.18).

\dtl{Remark}{rem:Dk-1}
(i)\,To a modular tensor category $\cC$ one associates a {\em dimension\/} 
${\rm Dim}(\cC)$ and the (unnormalised) {\em charges\/} $p^\pm(\cC)$ by
  \be
  {\rm Dim}(\cC) := \sum_{i \in \II} {\dim(U_i)}^2_{} 
  \qquad {\rm and} \qquad
  p^\pm(\cC) := 
  \sum_{i\in\II} \theta_i^{\pm 1}\, {\dim(U_i)}^2_{} \,, \labl{eq:dim-kappa-def}
where $\{U_i\,|\,i\iN\cI\}$ are representatives of the isomorphism classes
of simple objects of $\cC$. The numbers ${\rm Dim}(\cC)$ and $p^\pm(\cC)$ are 
non-zero (see e.g.\ \Corollary 3.1.8 of \cite{BAki}) and satisfy
$p^+(\cC)\, p^-(\cC) \eq {\rm Dim}(\cC)$.
\\
Let $A$ be a haploid commutative symmetric special Frobenius algebra in $\cC$.
Combining \Theorem 4.1 of \cite{kios} with \Theorem 3.1.7 of
\cite{BAki}, one sees that the dimension and charge obey
  \be
  {\rm Dim}(\calcal) = \frac{{\rm Dim}(\cC)}{{(\dim^\cC(A))}^2}               
  \qquad {\rm and} \qquad
  p^\pm(\calcal) = \frac{p^\pm(\cC)}{\dim^\cC(A)} \,.  \labl{eq:dim-kappa-CAloc}
Suppose now that $\koerper\,{=}\,\complex$ and that $\dim(U)\,{\ge}\,0$ for all
$U$ (as is e.g.\ the case if $\cC$ is a *-category \cite{loro}). Then one has 
in fact $\dim(U)\,{\ge}\,1$ for all non-zero objects, as well as 
${\rm Dim}(\cC)\,{\ge}\,1$ and $|p^+/p^-|\eq1$. It also follows 
that either $\dim(A)\eq1$ or else $\dim(A)\,{\ge}\, 2$, so 
that for any non-trivial $A$ the dimension of \calcal\ is at most one quarter of 
the dimension of $\cC$. The relation ``being a category of local $A$-modules'' 
(with $A$ a haploid commutative symmetric special Frobenius algebra in another 
category) thus induces a partial ordering `$>$' on modular tensor categories,
given by $\cC \,{>}\, \cD$ iff $\cD \,{\cong}\, \calcal$ for some 
$A\,{\not\cong}\,\one$. Also note that owing to ${\rm Dim}(\cC)\,{\ge}\,1$ one 
can repeat the procedure of ``going to the category of local modules'' only a 
finite number of times. Conversely, it follows that the dimension of a haploid 
commutative special Frobenius \alg\ in a modular tensor category $\cC$ is 
bounded by the square root of the dimension of $\cC$.  
\\[.3em]
(ii)~In case $A$ is a commutative simple \ssFA, the numbers $s^A$ that are the 
analogs of the numbers \erf{sUV} in the \cat\ 
\calcal\ can be expressed in terms of morphisms of $\cC$ as
  \bea \begin{picture}(70,42)(0,28)
  \put(54,0)  {\begin{picture}(0,0)(0,0)
              \scalebox{.38}{\includegraphics{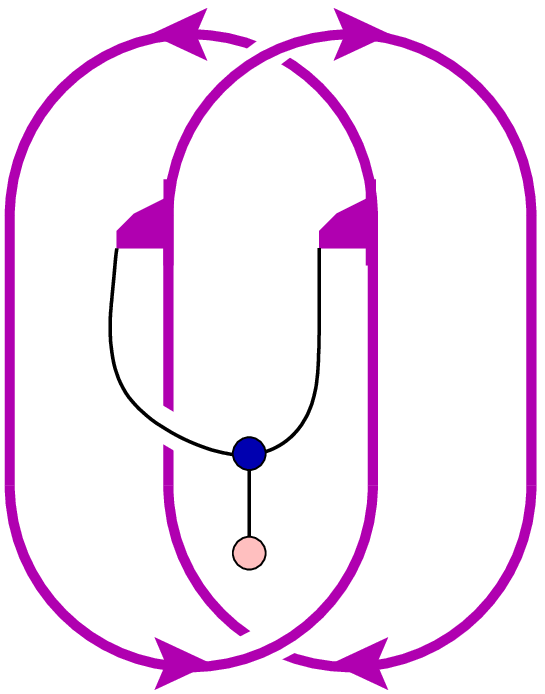}} \end{picture}}
  \put(-45,32.3)  {$s^A_{M,M'}\,=\;\dsty\frac1{\dim(A)}$}
  \put(59.9,31.3){\sse$A$}
  \put(73.7,39.3){\sse$\M$}
  \put(96.2,29.3){\sse$\M'$}
  \epicture13 \labl{sA}
It follows e.g.\ that
  \be  \dim_A(M) \equiv s^A_{M,0} = \dim(\M)\,/\dim(A)   \labl{didi}
(see \Theorem 1.18 of \cite{kios}). Note that the label 0 on $s^A$ refers 
to the tensor unit of \calcal, which is the simple local module $A$ itself. 
In the application to conformal field theory, $s^A$ is also closely related 
to the modular S-transformation of conformal one-point blocks on a torus
with insertion $A$ (see \cite{bant6} and \Section 5.7 of \cite{fuRs4}).

\bigskip

Next we study what can be said about Karoubianness of categories of local 
modules. Recall the statements about $A$-modules in \Remarks \ref{K-rem2} 
and \ref{K-rem2iv}. It follows immediately with the help of the functoriality 
of the braiding that if the $A$-module $(\M,\r)$ is in addition local, then so 
are the $A$-module $(\Im(p),r{\circ}\r\,{\circ}(\id_\AA{\otimes}e))$ 
\erf{Imp-rep} in $\cC$ and the $A$-module $((\M;p),p\cir\r)$ \erf{Mp-rep} 
in the Karoubian envelope $\kar\cC$.

According to \Remark \ref{K-rem2iv}, non-split idempotents in $\cC$ can be 
used to build $(A;\id_\AA)$-modules in $\kar\cC$ which do not come from
an $A$-module in $\cC$. Thus in general the category $\kar{(\calca)}$ is a 
{\em proper\/} sub\cat\ of $(\kar\cC)_\AA^{}$. On the other hand, we still 
have the following results, which later on will allow us to establish, 
in corollary \ref{cor:CAK-CKA-mod}, equivalence of these two \cats\ if 
$A$ is not just an algebra but even a special Frobenius \alg.
\\[-2.3em]

\dtl{Lemma}{lem:CAK-CKA}
(i)~\,If $A$ is a commutative \ssFA\ in a Karoubian ribbon \cat\ $\cC$, then
the \cat\ $\calcal$ of local $A$-modules in $\cC$ is Karoubian as well.
\\[.3em]
(ii)~For any \alg\ $A$ in a (not necessarily Karoubian) \tc\ $\cC$ the
\cat\ $(\kar\cC)_{(A;\iD_\AA)}$ is Karoubian, i.e.\ one has the equivalence
  \be  \kar{ \Llb (\kar\cC)^{}_{(A;\iD_\AA)} \Lrb }
  \cong (\kar\cC)^{}_{(A;\iD_\AA)}  \labl{CKA1}
of categories.
If $\cC$ is ribbon and $A$ is commutative symmetric special Frobenius,
then also the \cat\ $\Ext{(\kar\cC)}{(A;\iD_\AA)}$ is Karoubian,
and one has the equivalence
  \be  \kar{ \Llb \Ext{(\kar\cC)}{(A;\iD_\AA)} \Lrb }
  \cong \Ext{(\kar\cC)}{(A;\iD_\AA)} \,.  \labl{CKA2}
of ribbon categories.

\medskip\noindent
Proof:\\
(i)~\,Since $\calcal$ is a full sub\cat\ of $\calca$, the assertion
follows from immediately from the analogous statement about $\calca$ in 
\Lemma \ref{K-rem2iii}.  \\[.3em]
(ii)~Since $\kar\cC$ is Karoubian, the two equivalences
are directly implied by \Lemma \ref{K-rem2iii} 
and by (i), respectively.  That the second equivalence preserves the 
ribbon structure is easily seen by writing out the equivalence explicitly. 
\qed

%%%%%%%%%%%%%%%%%%%%%%%%%%%%%%%%%%%%%%%%%%%%%%%%%%%%%%%%%%%%%%%%%%%%%%%%
\newpage

\sect{Local induction}\label{sect3komma5}

\subsection{The local induction functors}\label{locind}

We have already announced above that the endofunctors $\EFU{l/r}A$
with respect to a \ssFA\ $A$
are related to {\em local induction\/}, i.e.\ functors from $\cC$ to a full
subcategory of the category $\cC_{C_{l/r}(A)}$ of modules over the left
and right center of $A$, respectively, that share many properties of 
induction. As shown in \Proposition \ref{pr:EAB-locmod} below, the objects 
$E^{l/r}_A(U)$ in the image of these endofunctors possess an additional 
property: they are {\em local\/} ${C_{l/r}(A)}$-modules. Accordingly, 
the relevant full subcategories are the categories $\Ext{\cC}{C_{l/r}(A)}$ of
local $C_{l/r}(A)$-modules. The corresponding local induction functors, to
be denoted by $\LXTp\AA{l/r}$, from $\cC$ to $\Ext{\cC}{C_{l/r}(A)}$ will be
introduced in \Definition \ref{def-lxt} below. In the special case that
already $A$ itself is commutative, the centers coincide with $A$, and
accordingly there is only a single local induction procedure, which is
a functor from $\cC$ to the category $\Ext{\cC}\AA$ of local $A$-modules.
\\[-2.3em]

\dtl{Proposition}{pr:EAB-locmod}
Let $A$ be a symmetric special Frobenius algebra in a ribbon category
$\cC$. Then for any object $U$ of $\cC$, $\efu Ul\AA$ is
a local $C_l(A)$\,-module and $\efu Ur\AA$ is a local $C_r(A)$\,-module.
The representation morphisms are given by
  %% [pic~75]
  \bea  \begin{picture}(180,90)(0,14)
  \put(0,0)  {\begin{picture}(0,0)(0,0)
              \scalebox{.38}{\includegraphics{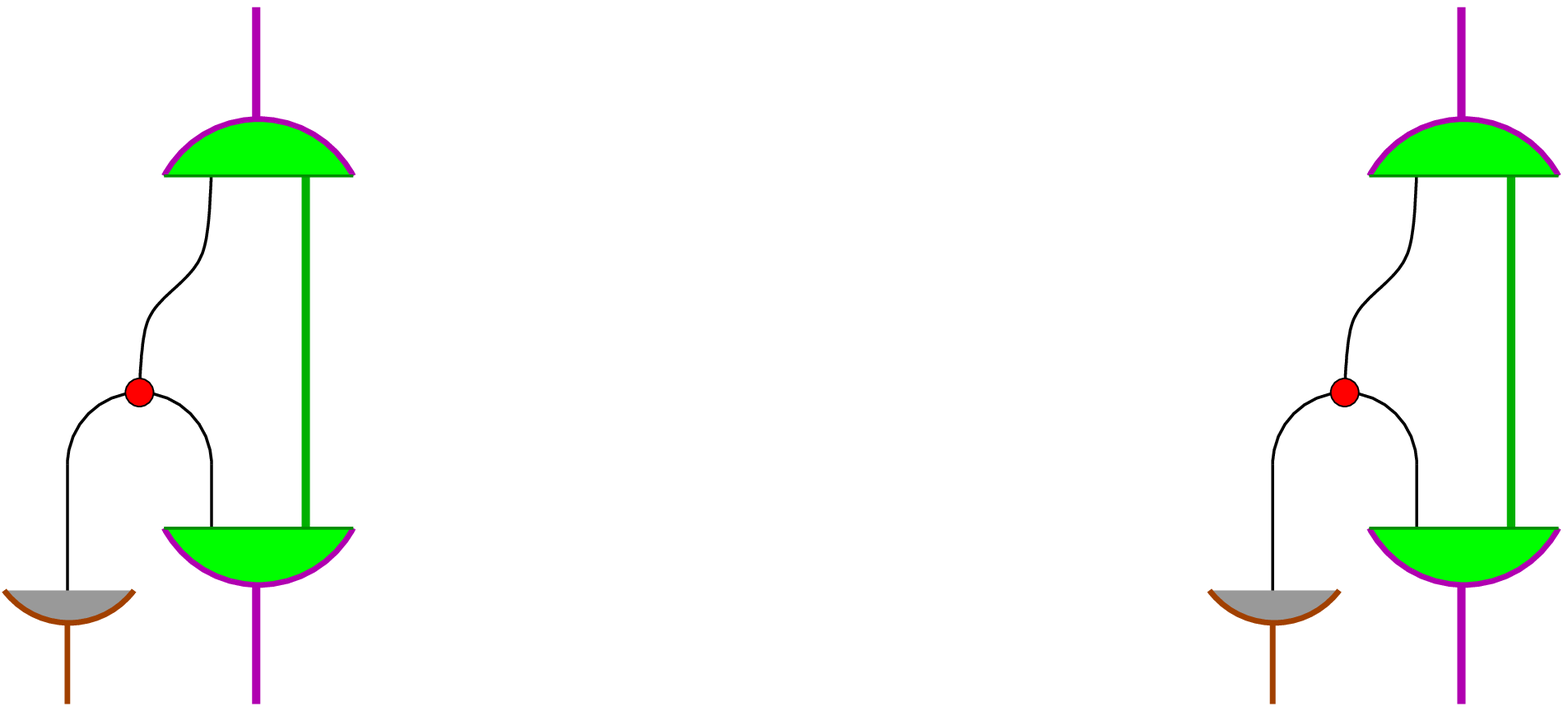}} \end{picture}}
  \put(-60,44.2)    {$\rholoc {C_{l}(A)}U\;=:$}
  \put(4.4,-9.2)    {\sse$C_l$}
  \put(10.9,32.2)   {\sse$A$}
  \put(22.4,100.6)  {\sse$\efu Ul\AA$}
  \put(24.1,-9.2)   {\sse$\efu Ul\AA$}
  \put(42.1,45.2)   {\sse$U$}
  \put(102,44.2)    {$\rholoc {C_{r}(A)}U\;=:$}
  \put(166.1,-9.2)  {\sse$C_r$}
  \put(172.5,32.2)  {\sse$A$}
  \put(184.1,100.6) {\sse$\efu Ur\AA$}
  \put(185.7,-9.2)  {\sse$\efu Ur\AA$}
  \put(203.7,45.2)  {\sse$U$}
  \epicture05 \labl{eq:Elr-rep}

\medskip\noindent
Proof:\\
Using the properties \erf{Cl-Cr-defprop} it is easily verified that 
$\rholoc{C_{l/r}(A)}U$ as defined in \erf{eq:Elr-rep} possess the properties 
of a representation morphism for $C_l(A)$ and $C_r(A)$, respectively. To 
establish locality we must check that $\rholoc{C_{l/r}(A)}U \cir P_{C_{l/r}}
(U)\eq\rholoc{C_{l/r}(A)}U$. This can be seen by inserting an idempotent 
$P^{l/r}_\AA(U)$ in front of the embedding morphism $e$ of $\efu U{l/r}\AA$; 
afterwards this idempotent can be used to remove $P_{C_{l/r}}(U)$. For the 
case of $E^r_\AA(U)$, the corresponding moves look as follows.
  %% [pic~56]
  \bea  \begin{picture}(330,395)(0,4)
  \put(0,0)  {\begin{picture}(0,0)(0,0)
              \scalebox{.38}{\includegraphics{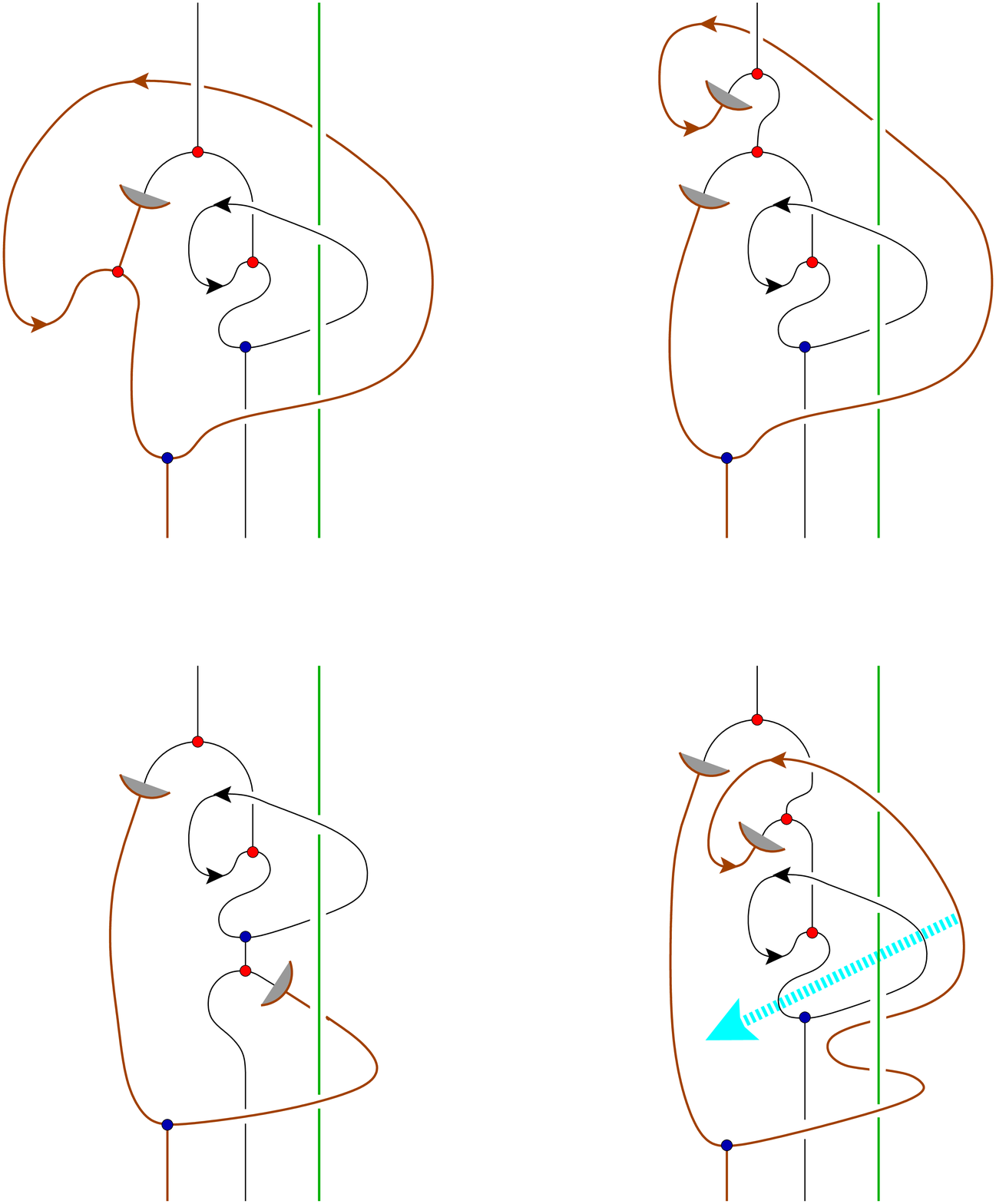}} \end{picture}}
  \put(50.2,208.2)  {\sse$C_r$}
  \put(60.8,396.9)  {\sse$A$}
  \put(75.8,208.2)  {\sse$A$}
  \put(100.7,208.2) {\sse$U$}
  \put(101.1,396.9) {\sse$U$}
  \put(172,298)     {$=$}
  \put(233.2,208.2) {\sse$C_r$}
  \put(243.8,396.9) {\sse$A$}
  \put(259.1,208.2) {\sse$A$}
  \put(284.3,208.2) {\sse$U$}
  \put(284.7,396.9) {\sse$U$}
  \put(-5,80)       {$=$}
  \put(50.2,-9.5)   {\sse$C_r$}
  \put(60.8,179.1)  {\sse$A$}
  \put(76.1,-9.5)   {\sse$A$}
  \put(101.3,-9.5)  {\sse$U$}
  \put(101.7,179.1) {\sse$U$}
  \put(165,80)      {$=$}
  \put(233.2,-9.5)  {\sse$C_r$}
  \put(243.8,179.1) {\sse$A$}
  \put(259.1,-9.5)  {\sse$A$}
  \put(284.1,-9.5)  {\sse$U$}
  \put(284.5,179.1) {\sse$U$}
  \epicture-4 \labl{eq:local-aux}
Here the embedding and restriction morphisms for $E^r_\AA(U)\,{\prec}\,A\oti U$ are 
omitted. To establish these equalities one needs in particular \erf{eq:remove-P}
and the properties \erf{Cl-Cr-defprop} and \erf{eq:remove-braiding} of $C_r$.  
\qed

\dtl{Corollary}{lem:[U]A-module-i}
Let $A$ be a commutative \ssFA\ in a ribbon category $\cC$ and
$U\iN\Obj(\cC)$. Then the object $\Efu U\AA\,{:=}\,\efu Ul\AA\eq\efu Ur\AA$
carries a natural structure of local $A$-module with representation morphism
  %% [pic~42]
  \bea  \begin{picture}(100,64)(0,32)
  \put(50,0)  {\begin{picture}(0,0)(0,0)
              \scalebox{.38}{\includegraphics{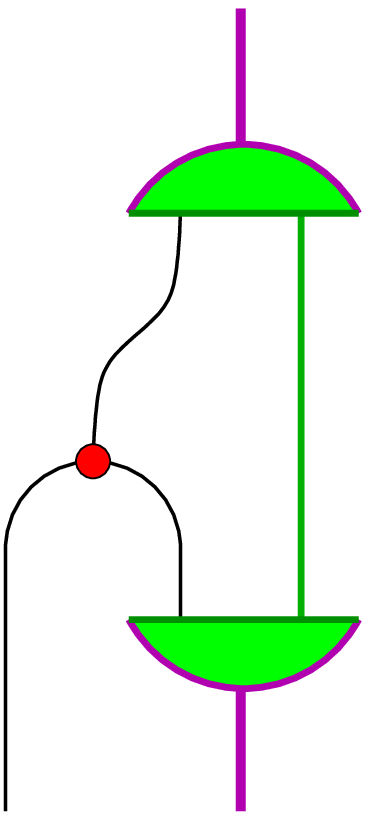}} \end{picture}}
  \put(-1,43.5)   {$\rholoc \AA U\,:=$}
  \put(46.3,-8.1) {\sse$A$}
  \put(67.5,-8.1) {\sse$\Efu U\AA$}
  \put(67.5,95.5) {\sse$\Efu U\AA$}
  \put(85.1,38.7) {\sse$U$}
  \epicture19 \labl{eq:[]rho-def}

It follows that given any \ssFA\ $A$ in a ribbon \cat, by regarding
$\efu U{l/r}\AA$ as an object of the category $\Ext\cC{C_{l/r}(A)}$ of local
$C_{l/r}$-modules we have a functor from $\cC$ to $\Ext\cC{C_{l/r}(A)}$.
\\[-2.3em]

\dtl{Definition}{def-lxt}
The functors $\LXTp\AA{l/r}$, called (left, respectively right) {\em local
induction functors\/}, from $\cC$ to $\Ext\cC{C_{l/r}(A)}$ are defined by
  \be
  \lxtp U\AA{l/r} := (\efu U{l/r}\AA,\rholoc {C_{l/r}(A)} U) \,, \qquad
  \lxtp f\AA{l/r} := E_\AA^{l/r}(f)  \,.  \ee
When $A$ is commutative, we write $\LXT\AA$ for $\LXTp\AA l\eq\LXTp\AA r$.

\medskip

The qualification `local' used here indicates that the $A$-module $\lxt U\AA$
is local; that we speak of local {\em induction\/} is justified by the
observation that there exists an embedding
of $\lxt U\AA$ into the induced module $\Ind_\AA(U)$.
More precisely, we have the following result, which allows us to use
reciprocity theorems of ordinary induction when working with local induction.
\\[-2.3em]

\dtl{Proposition}{lem:[U]A-module-ii}
For $A$ a commutative \ssFA\ in a ribbon category $\cC$ and
$\lxt U\AA$ endowed with the $A$-module structure given in corollary
\ref{lem:[U]A-module-i}, for every local $A$-module $M$ one has
  \be \bearl
  \HomA(M,\lxt U\AA) \,\cong\, \HomA(M,\Ind_\AA(U))  \qquad{\rm and}
  \\{}\\[-.6em]
  \HomA(\lxt U\AA,M) \,\cong\, \HomA(\Ind_\AA(U),M)  \,. \eear \labl{locreci}

\medskip\noindent
Proof:\\
Consider the first isomorphism in \erf{locreci}. Apply \Lemma \ref{lem:rezI}
to the objects $M$ and $\Ind_\AA(U)$ of \calca\ to see that there is a
natural bijection
  \be
  \HomA(M,\lxt U\AA) \cong \{ \varphi\iN\HomA(M,\Ind_\AA(U)) \,|\,
  P_\AA(U) \cir \varphi\eq\varphi \} \,.  \labl{xcvb}
Further, observe that for 
every $A$-module $M$ and every $\varphi\iN\HomA(M,\Ind_\AA(U))$ we have
  %% [pic~71]
  \bea  \begin{picture}(335,96)(22,28)
  \put(0,0)  {\begin{picture}(0,0)(0,0)
             \scalebox{.38}{\includegraphics{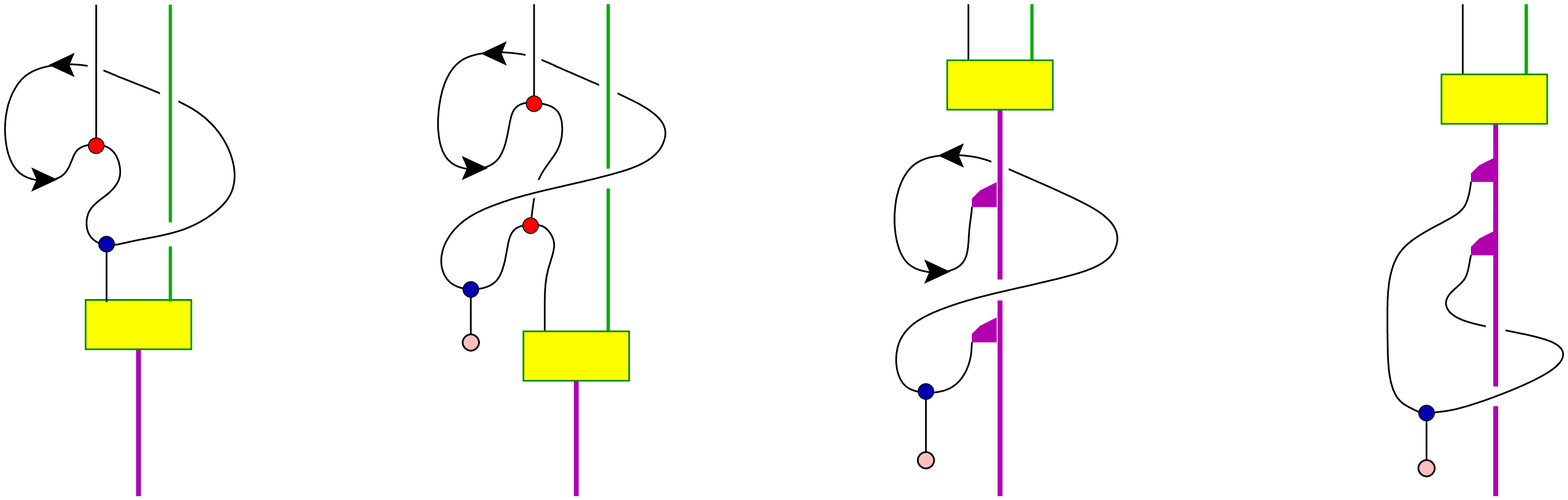}} \end{picture}}
  \put(18.1,118.1) {\sse$A$}
  \put(26.5,-9.2)  {\sse$\M$}
  \put(29.3,39.2)  {\sse$\varphi$}
  \put(35.9,118.1) {\sse$U$}
  \put(71,52.5)    {$=$}
  \put(118.6,118.1){\sse$A$}
  \put(126.6,-9.2) {\sse$\M$}
  \put(129.9,32.6) {\sse$\varphi$}
  \put(136.8,118.1){\sse$U$}
  \put(170,52.5)   {$=$}
  \put(218.4,118.1){\sse$A$}
  \put(224.1,-9.2) {\sse$\M$}
  \put(226.4,94.1) {\sse$\varphi$}
  \put(233.8,118.1){\sse$U$}
  \put(281,52.5)   {$=$}
  \put(331.7,118.1){\sse$A$}
  \put(337.9,-9.2) {\sse$\M$}
  \put(340.6,90.9) {\sse$\varphi$}
  \put(347.3,118.1){\sse$U$}
  \epicture15 \labl{pic71}
Here the first equality uses that $A$ is commutative and symmetric Frobenius, 
the second that $\varphi$ is an $A$-module morphism, and the third is a
rearrangement of the lower $A$-ribbon that uses that $A$ is commutative 
and symmetric and that (since it is also Frobenius) it has trivial twist.
\\
When $M$ is local, then by \Lemma \ref{ab(c)lemma}(iii) 
the \rhs\ of \erf{pic71} equals $\varphi$. Further, the left hand side of 
\erf{pic71} is nothing but $P_\AA(U) \cir \varphi$. Thus if $M$ is local and 
$\varphi$ a morphism in $\HomA(M,\Ind_\AA(U))$, then $P_\AA(U)\cir\varphi
\eq\varphi$ holds automatically. Together with \erf{xcvb} this implies 
the first bijection in \erf{locreci}.
\\[.3em]
The second of the bijections \erf{locreci} follows analogously by an identity 
between morphisms that looks like figure \erf{pic71} turned upside down.
\qed

\dtl{Lemma}{lem:CAK-CKA-ind}
Let $A$ be an \alg\ in a (not necessarily Karoubian) \tc\ $\cC$.
\\[.2em]
(i)~\,There is an equivalence
  \be  \kar{ \Llb (\kar\cC)_{(A;\iD_\AA)}^\ind \Lrb }
  \cong \kar{(\calcai)}  \labl{CKA3}
between Karoubian envelopes of \cats\ of induced modules.
\\[.2em]
(ii)~If $\cC$ is ribbon and $A$ is commutative symmetric special Frobenius,
then there is an equivalence
  \be  \kar{ \Llb (\kar\cC)_{(A;\iD_\AA)}^\lInd \Lrb }
  \cong \kar{(\calcali)}  \labl{CKA5}
between Karoubian envelopes of \cats\ of locally induced modules.

\medskip\noindent
Proof:\\
(i)~\,We will construct a functor $F$ from $\kar{(\calcai)}$ to
$\cD\,{:=}\,\kar{((\kar\cC)_{(A;\iD_\AA)}^\ind)}$ that satisfies the
criterion of \Proposition \ref{XI.1.5}.  \\
But first we consider the category $\cD$ in more detail.
Objects of $\,\cD$ are of the form%
  \foodnode{%
  We slightly abuse notation by writing just $\Ind_{(A;\iD_\AA)}(U;p)$
  in place of $\Ind_{(A;\iD_\AA)}((U;p))$.}
$(\Ind_{(A;\iD_\AA)}(U;p);\pi)$ with $U\iN\Objc$, and with 
$p\iN\End(U)$ and $\pi\iN\End(A\Oti U)$ idempotents satisfying
  \be  (\id_\AA\oti p) \cir \pi \cir (\id_\AA\oti p) = \pi
  \qquad{\rm and}\qquad
  \pi \cir (m\oti p) = (m\oti p) \cir (\id_\AA\oti\pi) \,.  \labl{pippi}
The latter properties imply that
  \be  \pi \cir (m\oti\id_U) = \pi \cir (m\oti p)
  = (m\oti\id_U) \cir (\id_\AA\oti\pi)  \,,  \labl{pippk}
which in turn allows us to regard $\pi$ as an idempotent in
$\End_{(A;\iD_\AA)}(\Indk(U;\id_U))$, i.e.\ in the space of endomorphisms of
an induced $(A;\id_\AA)$-module for which $p$ is replaced by $\id_U$. As a
consequence, $(\Indk(U;\id_U);\pi)$ is an object of $\cD$, and we have
  \be  \id^{}_{(\IndK(U;\iD_U);\pi)} = \pi = \id^{}_{(\IndK(U;p);\pi)}  \,. \ee
(All morphism spaces are regarded as subspaces of the corresponding spaces
of morphisms in $\cC$.)
\\
Furthermore, again using \erf{pippi}, it follows that the morphism spaces
of $\cD$ of our interest are of the form
  \be \bearl  \Hom^\cD((\Indk(U;q);\varpi),(\Indk(U;q');\varpi'))
  \\{}\\[-.7em]\mbox{\hspace{6em}}
  = \{\, f\iN\End(A\Oti U) \,|\, \varpi'{\circ}f{\circ}\,\varpi \eq f
  \eq (\id_\AA\Oti q')\,{\circ}f{\circ}\,(\id_\AA\Oti q)
  \\{}\\[-.97em]\mbox{\hspace{16em}}
  \;{\rm and}\;
  f{\circ}\,(m\Oti q) \eq (m\Oti q'){\circ}(\id_\AA\Oti f) \,\} \,. \eear \ee
By similar calculations as in \erf{pippk} one can then check that
  \be \bearl  \pi \in \Hom^\cD((\Indk(U;\id_U);\pi),(\Indk(U;p);\pi))
  \qquad{\rm and}  \\{}\\[-.7em]
  \pi \in \Hom^\cD((\Indk(U;p);\pi),(\Indk(U;\id_U);\pi)) \,,  \eear \ee
so that $(\Indk(U;p);\pi)$ and $(\Indk(U;\id_U);\pi)$ are isomorphic
as objects of $\cD$,
  \be  (\Indk(U;p);\pi) \cong (\Indk(U;\id_U);\pi) \,.  \labl{pippj}
Finally we observe that objects of $\kar{(\calcai)}$ are of the form
$(\Ind_\AA(U);\pi)$ with $U\iN\Objc$ and $\pi\iN\End_\AA(\Ind_\AA(U))$
an idempotent. Therefore by setting
  \be  F:\quad (\Ind_\AA(U);\pi) \,\mapsto\, (\Indk(U;\id_U);\pi)  \labl{iF}
for objects and defining $F$ to be the identity map on morphisms
provides us with a functor $F{:}\;\kar{(\calcai)}{\to}\,\cD$.
Because of \erf{pippj}, $F$ is essentially surjective, and it is
bijective on morphisms. By \Proposition \ref{XI.1.5}, 
$F$ thus furnishes an equivalence of \cats.
\\[.3em]
(ii)~The proof works along the same lines as for part (i). First note that 
objects of the \cat\ $\cD^{\sss\rmloc}\,{:=}\,\kar
{ \Llb (\kar\cC)_{(A;\iD_\AA)}^\lInd \Lrb }$ are of the form 
$(\lxt{(U;p)}{(A;\iD_A)}; \pi)$. On the other
hand, by definition we have $\lxt U\AA\eq (\Ind_\AA(U); P_\AA(U))$, so that
  \be
  (\lxt{(U;p)}{(A;\iD_A)}; \pi) = (\Ind_{(A;\iD_A)}(U;p); \pi ) \ee
with $P_\AA(U) \cir \pi \cir P_\AA(U)\eq\pi$. The rest of the arguments
in (i) go through unmodified, telling us that
  \be
  (\lxt{(U;p)}{(A;\iD_A)}; \pi) \,\cong\, (\lxt{(U;\id_U)}{(A;\iD_A)}; \pi) \,.
  \ee
Therefore the functor $F^{\sss\rmloc}$, defined as $F$ in \erf{iF} with
$\LXT{(A;\iD_A)}$ in place of $\Ind_{(A;\iD_A)}$, is essentially surjective
and bijective on morphisms, and hence furnishes an equivalence of categories.  
\qed
 
\dtl{Remark}{lem:[U]A-module-iii}
For any commutative \ssFA\ $A$ and any object $U$ of $\cC$
the dimension of $\Efu U\AA\iN\Objc$ is given by
  \be  \dim(\Efu U\AA) = s_{U,A}^{} \,.  \labl{sUA}
(The dimension of $\lxt U\AA$ as an object of $\Ext{\cC}A$ then follows 
via \erf{didi}.) The equality \erf{sUA} is easily verified by drawing the 
corresponding ribbon graphs:
  %% [pic~45]
  \bea  \begin{picture}(280,84)(0,34)
  \put(50,0)  {\begin{picture}(0,0)(0,0)
              \scalebox{.38}{\includegraphics{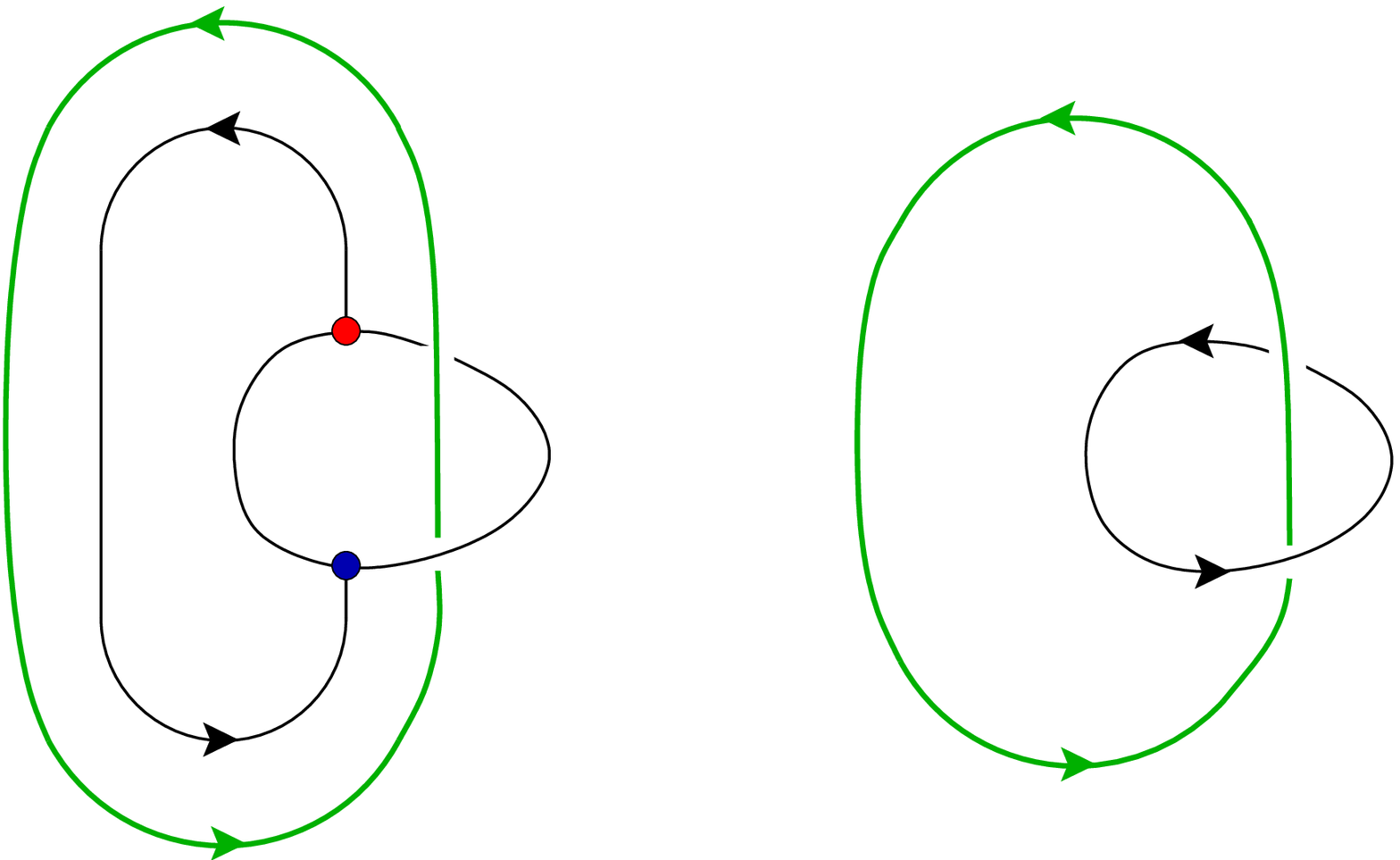}} \end{picture}}
  \put(-44,54.5)  {$\dim(\Efu U\AA)\;=$}
  \put(74.4,50.3) {\sse$A$}
  \put(111.4,77.7){\sse$U$}
  \put(121.2,43.7){\sse$A$}
  \put(141,54.5)  {$=$}
  \put(189.3,57.8){\sse$A^{\!\vee}_{}$}
  \put(225.5,79.3){\sse$U$}
  \put(237.5,43.5){\sse$A$}
  \put(258,54.5)  {$=\;s_{U,A}^{}\,.$}
  \epicture13 \labl{pic45}
The first equality uses the fact that for any retract $(S,e,r)$ of $U$
one has $\dim(S)\eq{\rm tr}_S\,\id_S\eq
    $\linebreak[0]$%
{\rm tr}_S\,r\cir e \eq{\rm tr}_U\, e\cir r\eq{\rm tr}_U\,P$,
applied to the idempotent $P\eq P_\AA$. In the second step the $A$-loop
that does not intersect the $U$-ribbon is omitted, using in particular 
the Frobenius property and specialness of $A$. The resulting graph is equal 
to $s_{U,A^\vee}^{}$; but $A \,{\cong}\, A^\vee$, since $A$ is Frobenius.

\dt{Remark}
When $\cC$ is modular, one may obtain \erf{locreci} also as follows. 
Proposition 5.22 of \cite{fuRs4} expresses the dimension $\,\dim\Hom_A
(M{\otimes}{U_k},N)$ as the invariant of a ribbon graph in $S^2{\times}S^1$:
  \bea  \begin{picture}(0,92)(25,28)
  \put(0,0)   {\begin{picture}(0,0)(0,0)
              \scalebox{.38}{\includegraphics{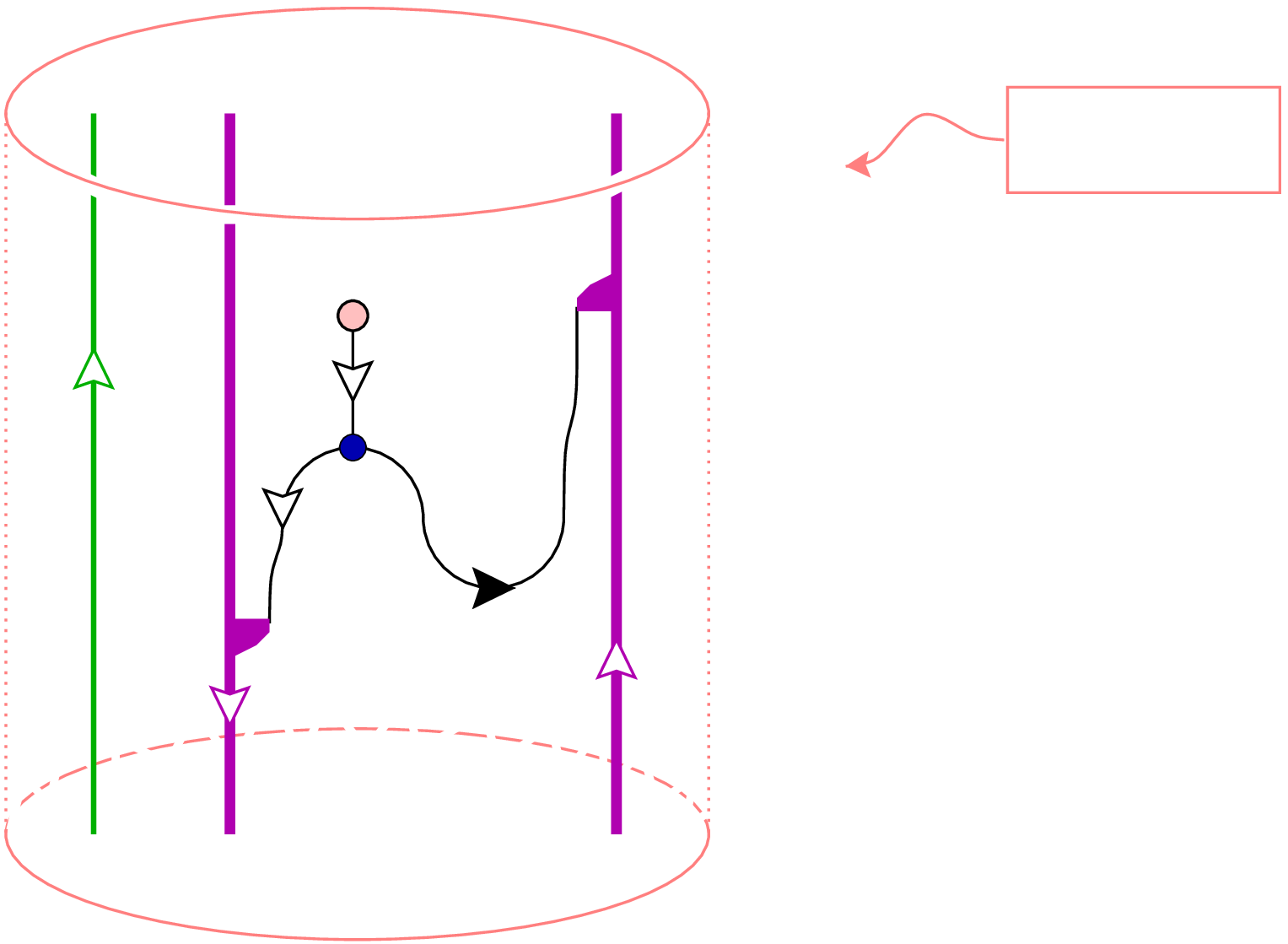}} \end{picture}}
  \put(-147,59.5) {$\dim\,\Hom_A(M{\otimes}{U_k},N)\;=$}
  \put(13.5,14.9) {\sse$k$}
  \put(31.1,14.9) {\sse$\dot N$}
  \put(70.9,14.9) {\sse$\M$}
  \put(132,100.6) {$S^2{\times}S^1$}
  \put(34,35.2) {\sse$\r_{\!N}^{}$}
  \put(64,86) {\sse$\r_{\!M}^{}$}
  \put(36.4,48) {\sse$A$}
  \put(57,53) {\sse$A$}
  \epicture07 \labl{Akmn}
Let us consider the case that $U_k\eq\one$, $M\eq\lxt U\AA$ and $N$ a
local module. Inserting the definition \erf{eq:[]rho-def} 
of $\r^{\Efu U\AA}$ and moving the restriction morphism $r$ around the (vertical)
$S^1$-direction so as to combine with the embedding $e$ to a projector,
then yields for $\dim\,\Hom_A(\lxt U\AA,N)$ the graph on the \lhs\ of
  %% [pic~44]
  \bea  \begin{picture}(330,81)(60,41)
  \put(3,0)  {\begin{picture}(0,0)(0,0)
              \scalebox{.38}{\includegraphics{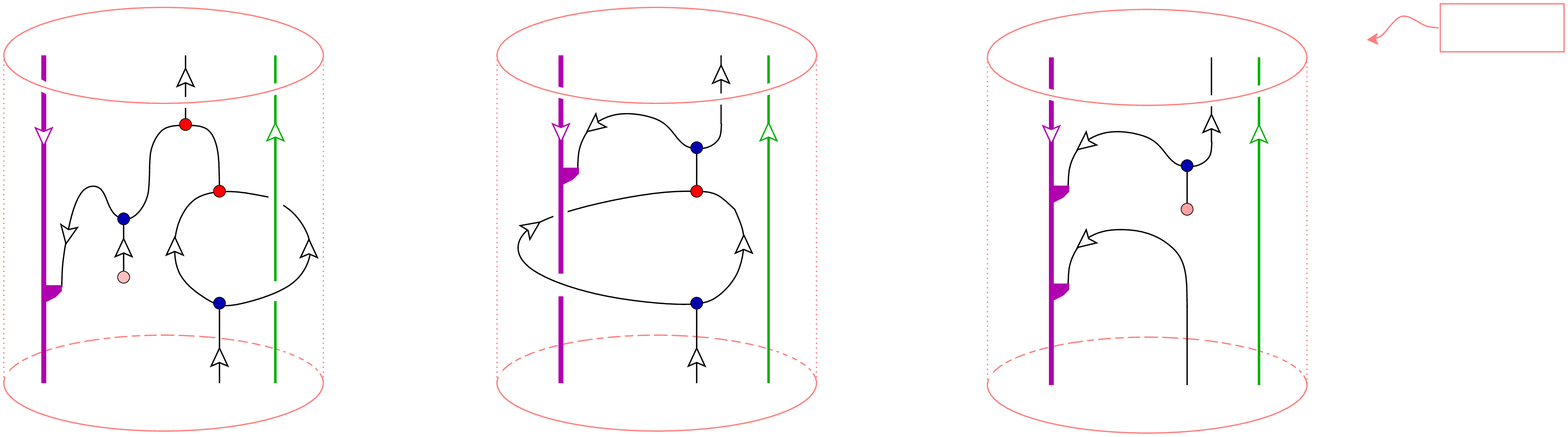}} \end{picture}}
  \put(5.9,90.7)    {\sse$\dot N$}
  \put(54.9,79.9)   {\sse$A$}
  \put(83.1,90.7)   {\sse$U$}
  \put(115.2,56.8)  {$=$}
  \put(153.3,89.2)  {\sse$\dot N$}
  \put(183.5,83.3)  {\sse$A$}
  \put(223.9,90.7)  {\sse$U$}
  \put(256.1,56.8)  {$=$}
  \put(293.7,89.2)  {\sse$\dot N$}
  \put(322.2,50.7)  {\sse$A$}
  \put(322.2,79.4)  {\sse$A$}
  \put(364.9,89.2)  {\sse$U$}
  \put(416.8,111.9) {\small$S^2{\times}S^1$}
  \epicture23 \labl{pic44}
The equalities shown here are obtained as follows. In the first step the 
$A$-ribbon of the projector is taken around the (horizontal) $S^2$-direction 
until it wraps around the $\dot N$-ribbon. This can be transformed into 
a locality projector for $N$ and thus -- as $N$ is local by assumption -- 
be left out. The second step is then completed by using the representation 
property for $N$. In the graph on the \rhs\ one can now move one of the
representation morphisms around the $S^1$-direction, and then use the 
representation property again; afterwards the $A$-ribbon can be removed, 
using that $A$ is special. The invariant of the resulting graph in 
$S^2\times S^1$ is $\dim\,\Hom(U,\dot N)$.

%%%%%%%%%%%%%%%%%%%%%%%%%%%%%%%%%%%%%%%%%%%%%%%%%%%%%%%%%%%%%%%%%%%%%%%%
\newpage

\subsection{Local modules from local induction}

In the sequel it will be very helpful to express categories of (local)
modules in terms of the corresponding categories of (locally) induced
modules. A crucial ingredient is the
\\[-2.3em]

\dtl{Lemma}{lem:sub-of-ind}
Let $A$ be a special Frobenius algebra in a (not necessarily Karoubian)
tensor category $\cC$.
\\[.2em]
(i)~\,For every module $M$ over $A$ the object $\M$ is a retract of $A\oti\M$.
\\[.2em]
(ii)~Every module over $A$ is a \retmodule\ of an induced $A$-mo\-dule.

\medskip\noindent
Proof: \\
(i)~\,The retract is given by $(\M,e_M,\r_M)$ with $\r_M$ the \rep\ morphism
of $M$ and $e_M\,{:=}\,(\id_\AA\oti\r_M)\cir((\Delta{\circ}\eta)\oti\id_\M)$.
That $\r_M\cir e_M\eq\id_\M$ is verified by first using the \rep\ property of
$\r_M$, then specialness of $A$, and then the unit property of $\eta$.
\\
Note that the Frobenius property \erf{1f} of $A$ is not used in this argument.
\\[.3em]
(ii)~We show that any $A$-module $M$ is a \retmodule\ of $\Ind_\AA(\M)$. In
view of (i), all that needs to be checked is that the morphisms $\r_M$ and
$e_M$ are module morphisms. That $\r_M\iN\HomA(\Ind_\AA(\M),M)$ follows
directly from the \rep\ property of $\r_M$, while $e_M\iN
\HomA(M,\Ind_\AA(\M))$ is a consequence of the Frobenius property of $A$.
\qed

This result has already been established in lemma 4.15 of \cite{fuSc16}.
(There the assumption was made that the \cat\ $\cC$ of which $A$
is an object is abelian, but the proof does not require this property.)

\dtl{Proposition}{sub-of-ind-ii}
Let $A$ be a special Frobenius \alg\ in a (not necessarily Karoubian) \tc\
$\cC$. Then, while the module category \calca\ is not necessarily Karoubian,
still the Karoubian envelopes of \calca\
and of its full sub\cat\ \calcai\ of induced $A$-modules coincide:
  \be  \kar{(\calca)} \cong \kar{(\calcai)}  \,.  \labl{calcai}
It follows in particular that in case that $\cC$ {\em is\/} Karoubian
(so that by \Lemma \ref{K-rem2iii}
\calca\ is Karoubian, too), then $\calca\,{\cong}\,\kar{(\calcai)}$.

\medskip\noindent
Proof:\\
Lemma \ref{lem:sub-of-ind} implies in particular that
every object of the category \calca\ of $A$-modules in $\cC$ is of the form
  \be  \Indap\AA pU := ({\rm Im}(p),r\cir(m\oti\id_U)\cir (\id_\AA\oti e)
   )  \labl{Indap}
with a suitable object $U\iN\Objc$ and $p$ a split idempotent such that
  \be  p\iN\HomA(\Indap\AA{}U,\Indap\AA{}U) \,, \qquad
  p\cir p = p\,, \quad\; e \cir r = p \,, \quad\; r \cir e = \id_{{\rm Im}(p)}
  \,.  \ee
This implies the equivalence \erf{calcai}.
\qed

\medskip

Not surprisingly, \Lemma \ref{lem:sub-of-ind} and 
\Proposition \ref{sub-of-ind-ii} have analogues for local modules. Indeed,
when combined with the previous result \erf{locreci}, they imply:
\\[-2.3em]

\dtl{Corollary}{cor:loc-lind}
Let $A$ be a \csplit\ commutative \ssFA\ in a (not necessarily Karoubian) 
ribbon category $\cC$. Then every local module over $A$ is a module retract 
of a locally induced $A$-module, and we have
  \be  \kar{(\calcal)} \cong \kar{(\calcali)}  \,.  \labl{calcali}

\smallskip

The equivalence \erf{calcai} can be combined with previously established
equivalences, in particular \Lemma \ref{lem:CAK-CKA-ind}, 
to establish the following properties of module \cats\ over special 
Frobenius \alg s. They are much stronger than \Lemma \ref{lem:CAK-CKA-ind},
and they do not hold, in general, for \alg s that are not special Frobenius.
\\[-2.3em]

\dtl{Corollary}{cor:CAK-CKA-mod}
(i)~\,For any special Frobenius \alg\ $A$ in a (not necessarily Karoubian)
tensor \cat\ $\cC$ there is an equivalence
  \be  (\kar\cC)_{(A;\iD_\AA)}^{} \cong \kar{(\calca)} \,,  \labl{CKA4}
i.e.\ the operations of taking the Karoubian envelope and of
forming the module \cat\ commute.
\\[.3em]
(ii)~For any commutative \ssFA\ $A$ in a (not necessarily Karoubian) ribbon
\cat\ $\cC$ there is an equivalence
  \be  \Ext{(\kar\cC)}{(A;\iD_\AA)} \cong \kar{(\calcal)} \,,  \labl{CKA7}
    % could also show: as ribbon
i.e.\ the operations of taking the Karoubian envelope and of
forming the \cat\ of local modules commute.

\medskip\noindent
Proof:\\
(i)~\,We have
  \be  \kar{(\calca)} \cong \kar{(\calcai)}
  \cong \kar{ \Llb (\kar\cC)_{(A;\iD_\AA)}^\ind \Lrb }
  \cong \kar{ \Llb (\kar\cC)_{(A;\iD_\AA)}^{} \Lrb }
  \cong (\kar\cC)_{(A;\iD_\AA)}^{} \,.  \ee
The last equivalence follows by \Lemma \ref{K-rem2iii},
the second equivalence is the one of \Lemma \ref{lem:CAK-CKA-ind}(i),
and the other two equivalences hold by \Proposition \ref{sub-of-ind-ii}.
\\[.3em]
(ii)~Analogously,
  \be  \kar{(\calcal)} \cong \kar{(\calcali)}
  \cong \kar{ \Llb (\kar\cC)_{(A;\iD_\AA)}^\lInd \Lrb }
  \cong \kar{ \Llb \Ext{(\kar\cC)}{(A;\iD_\AA)} \Lrb }
  \cong \Ext{(\kar\cC)}{(A;\iD_\AA)} \,.  \ee
The last equivalence follows by \Lemma \ref{lem:CAK-CKA}(i),
the second equivalence is the one of \Lemma \ref{lem:CAK-CKA-ind}(ii)
and the other two equivalences hold by corollary \ref{cor:loc-lind}.
\qed

\medskip

The statements of \Proposition \ref{lem:[U]A-module-ii} 
and the results above about commutative Frobenius \alg s that are based on 
that proposition do not directly generalise to the non-commutative case.
However, there is the following substitute:
\\[-2.2em]

\dtl{Proposition}{sub-of-locind} 
Let $A$ be a \ssFA\ in a ribbon \cat\ $\cC$, and assume that the 
commutative symmetric Frobenius \alg\ $C_l(A)$ is special.
\\[.2em]
Then every local $C_l(A)$-module $M$ is a \retmodule\ of a locally induced
$A$-module, $M\,{\prec}\,\lxtp U\AA l$ with suitable $U\iN\Objc$.
\\[.2em]
Similarly, if $C_r(A)$ is special, then every local $C_r(A)$-module is a
\retmodule\ of $\lxtp U\AA r$ with suitable $U\iN\Objc$.

\medskip\noindent
Proof:\\
We establish the statement for $C_l\,{\equiv}\,C_l(A)$.
\\
Let $M$ be a local $C_l$-module. Choose $U\eq\Im \efu \M r\AA$ and define
morphisms $\tilde e$ and $\tilde r$ as
  %% [pic~f6/f7]
  \bea  \begin{picture}(330,94)(0,34)
  \put(30,0)  {\begin{picture}(0,0)(0,0)
              \scalebox{.38}{\includegraphics{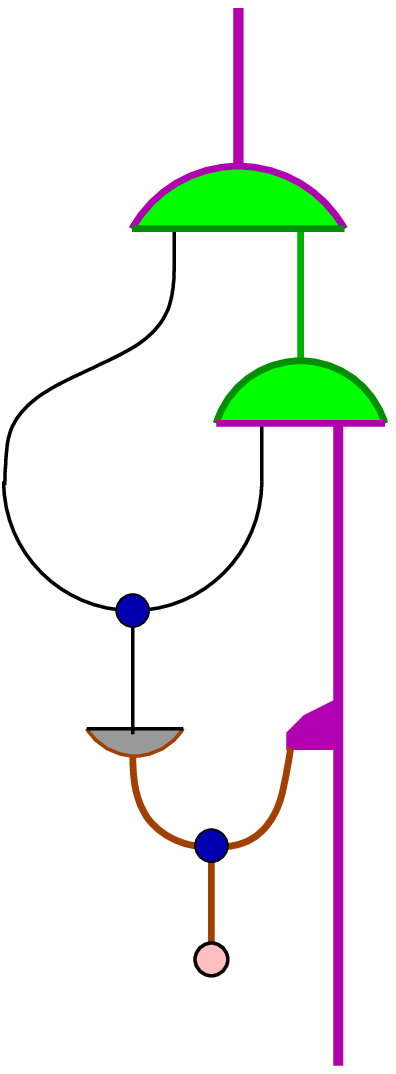}} \end{picture}}
  \put(270,0) {\begin{picture}(0,0)(0,0)
              \scalebox{.38}{\includegraphics{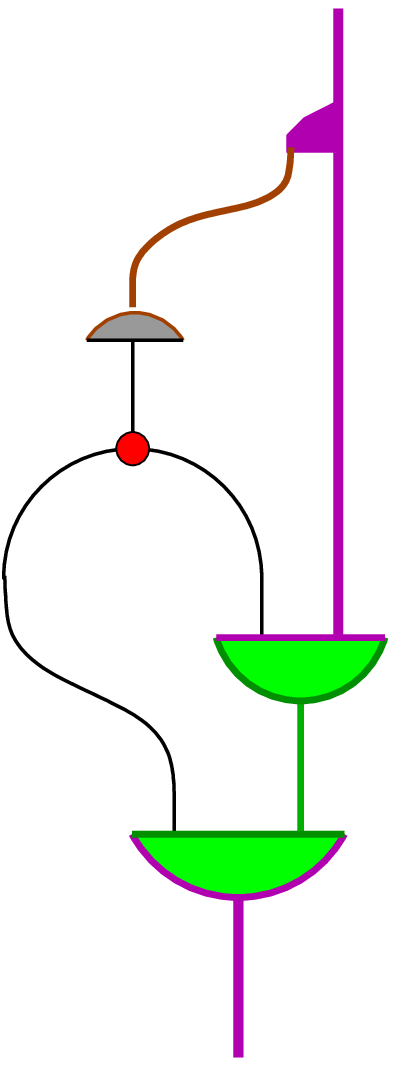}} \end{picture}}
  \put(-10,50)      {$\tilde e\,:=$}
  \put(36.8,55.6)   {\sse$A$}
  \put(41.7,121.5)  {\sse$\lxtp U\AA l$}
  \put(43.7,18.1)   {\sse$C_l$}
  \put(61.7,-9.6)   {\sse$\M$} 
  \put(64.8,82.8)   {\sse$U$} 
  \put(130,50)      {and $\quad\qquad\tilde r\,:= \Frac{\dim(A)}{\dim(C_l)} $}
  \put(273.6,58.5)  {\sse$A$}
  \put(282.4,-9.6)  {\sse$\lxtp U\AA l$} 
  \put(287.2,85.1)  {\sse$C_l$}
  \put(302.8,121.2) {\sse$\M$}
  \put(304.4,30.4)  {\sse$U$}
  \epicture20 \labl{picf6/f7}
These are $C_l$-intertwiners, i.e.\ $\tilde e\iN \Hom_{C_l}(M,\lxtp U\AA{l})$
and $\tilde r\iN \Hom_{C_l}(\lxtp U\AA{l},M)$. To establish
that $(M,\tilde e,\tilde r)$ is a $C_l$-module retract of
$\lxtp U\AA{l}$ we must show that $\tilde r \cir \tilde e\eq\id_M$.
This is seen by the following series of moves.
  %% [pic~f8]
  \begin{eqnarray}  \begin{picture}(400,240)(14,0)
  \put(80,0)  {\begin{picture}(0,0)(0,0)
              \scalebox{.38}{\includegraphics{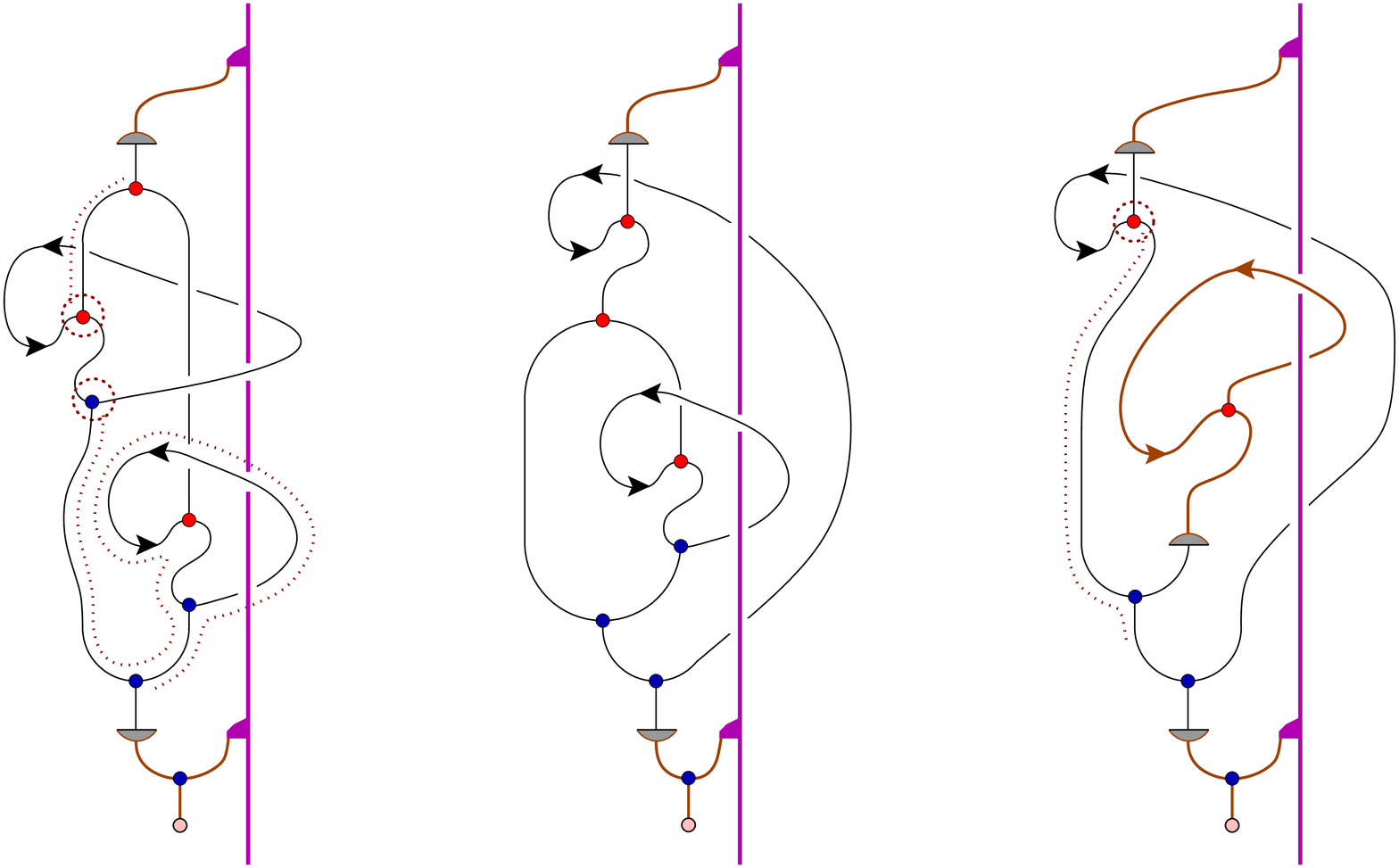}} \end{picture}}
  \put(-5,111)      {$\Frac{\dim(C_l)}{\dim(A)}\,\tilde r \cir \tilde e \;=$}
  \put(108.5,117.2) {\sse$A$}
  \put(117.2,16.2)  {\sse$C_l$}
  \put(119.7,199)   {\sse$C_l$}
  \put(139.7,-9.6)  {\sse$\M$}
  \put(140.7,234.5) {\sse$\M$}
  \put(185,111)     {$=$}
  \put(220.8,117.2) {\sse$A$}
  \put(251.7,199)   {\sse$C_l$}
  \put(252.8,16.2)  {\sse$C_l$}
  \put(271.7,-9.6)  {\sse$\M$}
  \put(272.7,234.5) {\sse$\M$}
  \put(332,111)     {$=$}
  \put(363.5,71.2)  {\sse$A$}
  \put(397.4,16.2)  {\sse$C_l$}
  \put(399.7,199)   {\sse$C_l$}
  \put(409,123.3)   {\sse$C_l$}
  \put(420.2,-9.6)  {\sse$\M$}
  \put(421.2,234.5) {\sse$\M$}
  \end{picture} \nonumber\\ \begin{picture}(160,210)(0,0)
  \put(80,0)  {\begin{picture}(0,0)(0,0)
              \scalebox{.38}{\includegraphics{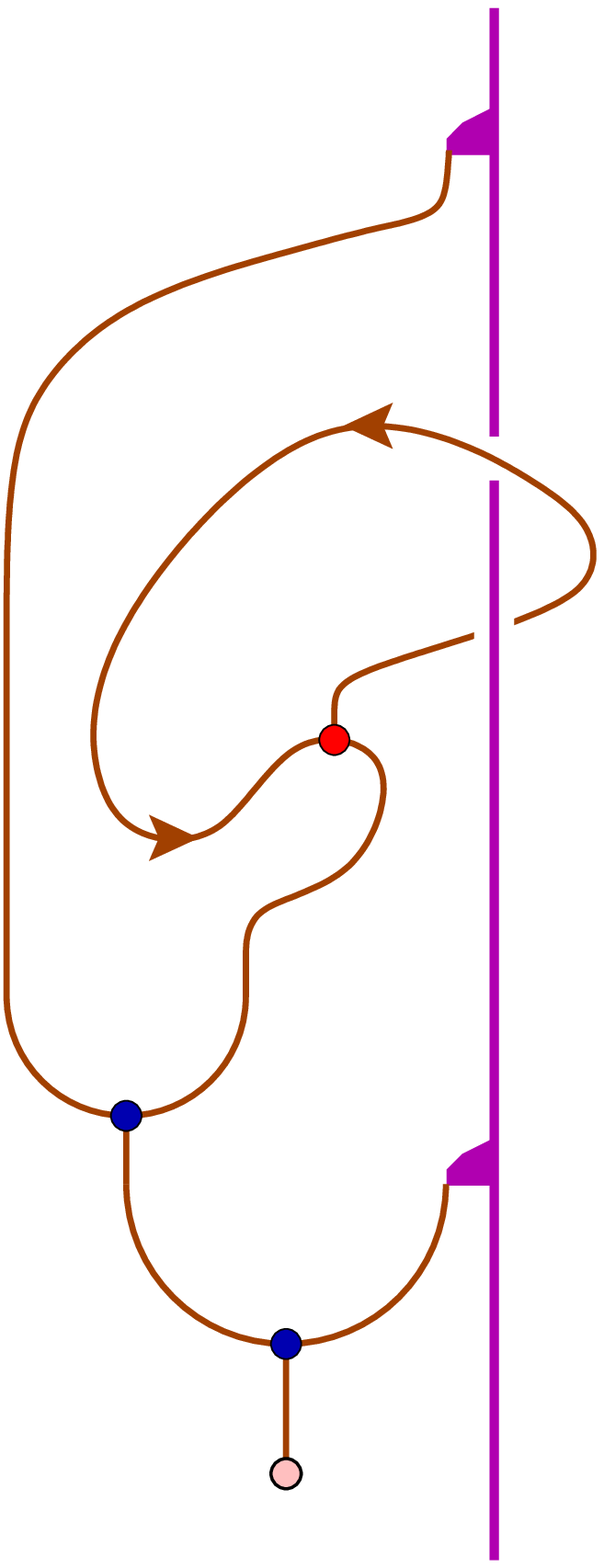}} \end{picture}} 
  \put(0,90)        {$ =\; \Frac{\dim(C_l)}{\dim(A)} $}
  \put(113.4,88.8)  {\sse$C_l$}
  \put(133.7,-9.6)  {\sse$\M$}
  \put(134.7,191.5) {\sse$\M$}
  \end{picture} \nonumber\\[-.7em]{} \label{picf8}
  \\[-1.7em]{}\nonumber\end{eqnarray}
In the first step the idempotents resulting from the composition are drawn
explicitly. Then the multiplication and comultiplication are moved along
the paths indicated. To the resulting morphism in the second
picture one can apply \Lemma \ref{le4sublocind} 
with $U\eq V\eq \M$ and 
$\Phi\eq c_{\AA,\M}^{-1} \cir c_{\M,A}^{-1}$. This results in the
insertion of an idempotent $P^l_\AA$. Using \Lemma \ref{lem:C=[1]a}(iii)
and the definition of the multiplication on $C_l$ in \erf{eq:Clr-alg}
one arrives at the third morphism. In the final step the marked coproduct 
is moved along the path indicated, resulting in another idempotent $P^l_\AA$, 
which can be omitted against the embedding morphism $e_{C_l}$. Inserting the 
definition of the comultiplication on $C_l$ in \erf{eq:Clr-alg} one finally 
arrives at the morphism on the \rhs.
\\
There, the $C_l$-loop can be rearranged to be equal to $P_{C_l}(M)$, using 
the fact that $C_l$ is a commutative symmetric Frobenius algebra. 
Afterwards, by the \Definition \ref{def:loc-mod} of a local module, the
idempotent $P_{C_l}(M)$ can be omitted. The representation property together
with specialness of $C_l$ imply that the resulting morphism it is equal to
$\dim(C_l)/{\dim}(A)\,\id_M$. Altogether we thus
have $\tilde r \cir \tilde e\eq\id_M$, showing that
$M$ is indeed a retract of $\lxtp U\AA{l}$.
\qed

Note that specialness of $C_{l/r}(A)$, which is assumed in the proposition, 
is guaranteed e.g.\ if $A$ is simple and $\dim(C_{l/r}(A))$ is non-zero, see 
\Proposition \ref{lem:C=[1]i}, and also if $A$ is commutative, because 
then $C_{l/r}(A)\eq A$ and $A$ is special by assumption.

%%%%%%%%%%%%%%%%%%%%%%%%%%%%%%%%%%%%%%%%%%%%%%%%%%%%%%%%%%%%%%%%%%%%%%%%

\subsection{Local induction of algebras}

Since for any \ssFA\ $A$ the \cats\ $\Ext{\cC}{C_{l/r}(A)}$ of local
modules over the left and right center of $A$ are \tcs, one can
study algebras in these \cats\ and, in particular, ask whether for an
\alg\ $B$ in $\cC$ the locally induced module $\lxtp B\AA{l/r}$ inherits an
algebra structure from $B$. We shall show that indeed the algebra
$\efu B{l/r\!}\AA$ as defined by \Proposition \ref{prop:AB-alg}(i) 
lifts to an algebra in $\Ext{\cC}{C_{l/r}(A)}$ and inherits further structural
properties. As a consequence, $\LXTp\AA{l/r}$ furnishes a functor from the
\cat\ of (symmetric special Frobenius) \alg s in $\cC$ to the
\cat\ of (symmetric special Frobenius) \alg s in $\Ext{\cC}{C_{l/r}(A)}$.

We start by formulating conditions that allow
an algebra $B$ in $\cC$ to be `lifted' to an algebra in $\Ext{\cC}A$:
\\[-2.3em]

\dtl{Lemma}{lem:alglift}
Let $A$ be a commutative symmetric special Frobenius algebra in a
ribbon category $\cC$. Let $B\,{\equiv}\,(B,m_B,\eta_B,\Delta_B,\eps_B)$ be a
Frobenius algebra. Let $(B,\r_B)$ carry the structure of a local
$A$-module, and the product $m_B$ on $B$ satisfy 
  \be
  m_B\in\HomA(B{\otimes}B,B)  \qquad{\rm and}\qquad
  m_B\cir P_{B{\otimes}B} = m_B \,.  \ee
(i)~\,$\tilde B\,{\equiv}\,(B,\tilde m_B,\tilde \eta_B,\tilde \Delta_B,
\tilde \eps_B)$ with
  \be
  \tilde m_B^{} := m_B^{} \circ e_{B\otimes B}^{} \,, \qquad
  \tilde\eta_B^{} := \r_B^{} \circ (\id_\AA^{} \oti \eta_B^{})
  \labl{eq:lift-mult}
and
  \be
  \tilde\Delta_B^{} := r_{\!B\otimes B}^{} \circ \Delta_B^{} \,, \qquad
  \tilde\eps_B^{} := (\id_\AA^{}\oti\eps_B^{}) \circ (\id_\AA^{}\oti\r_B^{})
  \circ ([\Delta_A^{}\cir\eta_A^{}]\oti\id_B^{})  \ee
is a Frobenius algebra
in \calcal.
\\[.3em]
(ii)~Let $A$ in addition be simple.
If $B$ has in addition any of the properties of being
commutative, haploid, simple, special, or symmetric, then so has $\tilde B$.

\medskip\noindent
Proof:\\
(i)~\,We start by showing that $P_{B\otimes B}\cir\Delta_B\eq\Delta_B$
is implied by $m_B\cir P_{B{\otimes}B}\eq m_B$. The ultimate reason is
that the coproduct can be expressed in terms of the product as
  \be
  \Delta_B = (\id_B \oti m_B) \circ (\id_B \oti \Phi_1^{-1}
  \oti \id_B) \circ (b_B \oti \id_B)  \labl{eq:Delta-m}
with the morphism $\Phi_1$, defined as in \erf{eq:Phi-def}, being invertible
because $B$ is a Frobenius algebra (see formula (3.36) of \cite{fuRs4} and,
for the proof, lemma 3.7 of \cite{fuRs4}). Consider the equivalences
  %% [pic~49]
  \bea  \begin{picture}(390,51)(7,28)
  \put(0,0)  {\begin{picture}(0,0)(0,0)
              \scalebox{.38}{\includegraphics{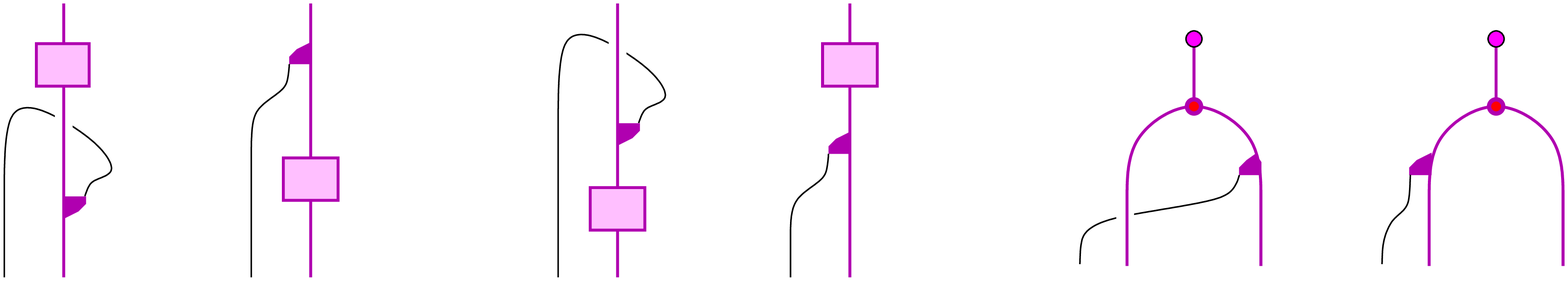}} \end{picture}}
  \put(-3.6,-8.7) {\sse$A$}
  \put(10.0,53.1) {\tiny$\Phi_1^{\!-\!1}$}
  \put(10.8,-8.7) {\sse$B^\vee$}
  \put(13.2,74.5) {\sse$B$}
  \put(41.6,31.5) {$=$}
  \put(59.2,-8.7) {\sse$A$}
  \put(73.2,23.8) {\tiny$\Phi_1^{\!-\!1}$}
  \put(73.6,-8.7) {\sse$B^\vee$}
  \put(76.0,74.5) {\sse$B$}
  \put(102.9,31.5){$\Longleftrightarrow$}
  \put(137.8,-8.7){\sse$A$}
  \put(153.2,17.3){\tiny$\Phi_1^{}$}
  \put(152.2,74.5){\sse$B^\vee$}
  \put(154.6,-8.7){\sse$B$}
  \put(180.8,31.9){$=$}
  \put(196.8,-8.7){\sse$A$}
  \put(212.2,54.1){\tiny$\Phi_1^{}$}
  \put(211.2,74.5){\sse$B^\vee$}
  \put(213.6,-8.7){\sse$B$}
  \put(242.0,31.5){$\Longleftrightarrow$}
  \put(270.8,-6.1){\sse$A$}
  \put(283.7,-6.1){\sse$B$}
  \put(316.8,-6.1){\sse$B$}
  \put(335.6,31.5){$=$}
  \put(347.9,-6.1){\sse$A$}
  \put(360.8,-6.1){\sse$B$}
  \put(393.9,-6.1){\sse$B$}
  \epicture16 \labl{eq:Phi1-rho}
The first equivalence follows by composing both sides of the first equality
with $\Phi_1$ both from the top and from the bottom. The second equivalence
is obtained by composing the middle equality with the duality morphism $d_B$
and writing out the definition \erf{eq:Phi-def} of $\Phi_1$. Now
the last equality in \erf{eq:Phi1-rho} indeed holds true. This can be seen
by replacing $m_B$ with $m_B\cir P_{B{\otimes}B}$ and using
commutativity and the Frobenius property of $A$
to move the action of $A$ along the resulting $A$-ribbon
from the right $B$-factor to the left. We can therefore write
  %% [pic~50]
  \bea  \begin{picture}(330,76)(0,47)
  \put(80,0)  {\begin{picture}(0,0)(0,0)
              \scalebox{.38}{\includegraphics{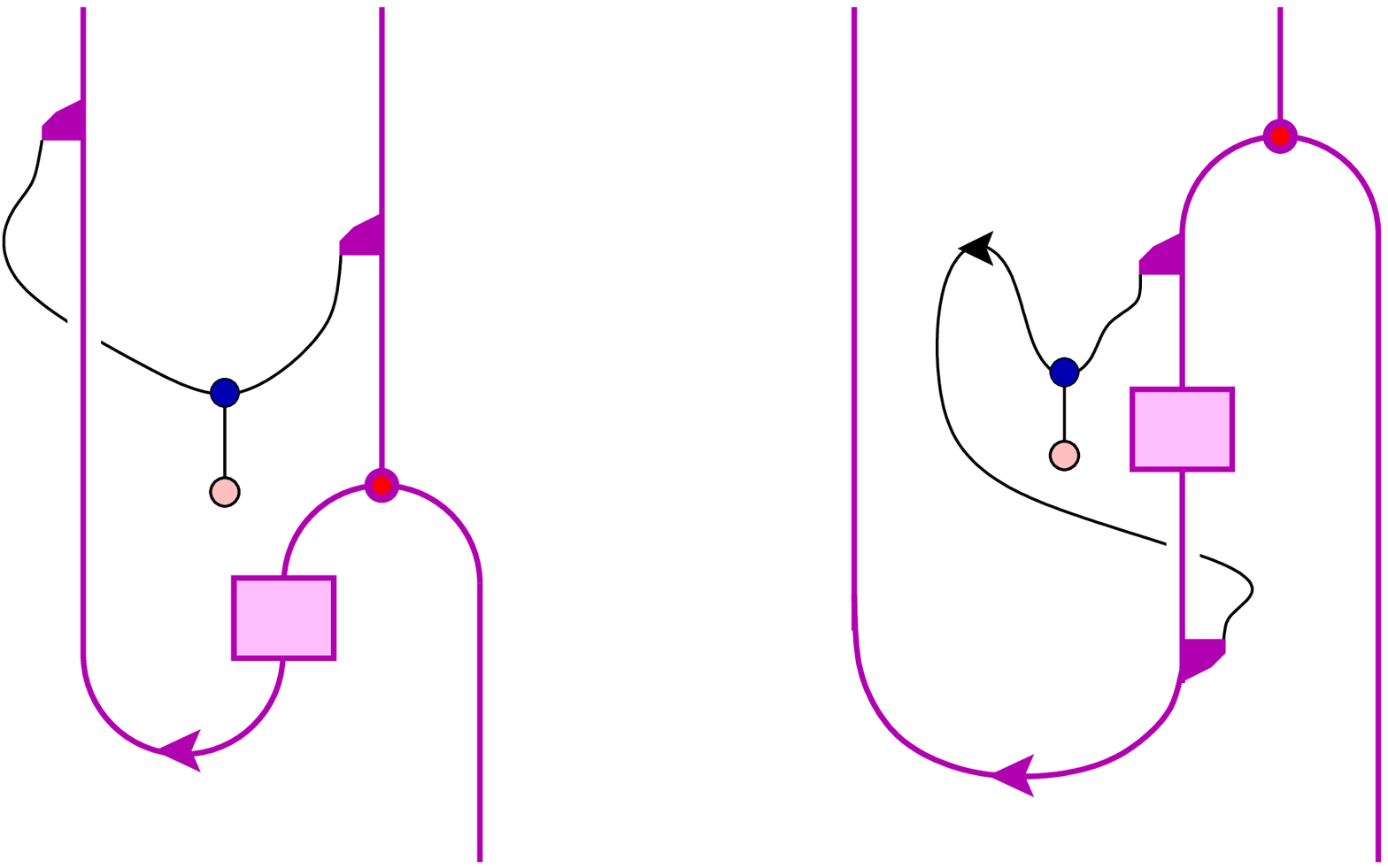}} \end{picture}}
  \put(-7,56.5)    {$P_{B\otimes B} \circ \Delta_B\,=$}
  \put(87.6,119.5) {\sse$B$}
  \put(99.8,69.7)  {\sse$A$}
  \put(112.8,31.7) {\tiny$\Phi_1^{\!-\!1}$}
  \put(128.1,119.5){\sse$B$}
  \put(141.1,-8.7) {\sse$B$}
  \put(162.6,56.5) {$=$}
  \put(191.6,119.5){\sse$B$}
  \put(214.8,39.4) {\sse$A$}
  \put(234.1,57.1) {\tiny$\Phi_1^{\!-\!1}$}
  \put(249.1,119.5){\sse$B$}
  \put(262.1,-8.7) {\sse$B$}
  \put(286,56.5)   {$=\,\Delta_B\,.$}
  \epicture33 \labl{PD=D}
The left-most graph is obtained by writing out the definition of
$P_{B{\otimes}B}$ and inserting relation \erf{eq:Delta-m} for $\Delta_B$.
The next step uses in particular that $m_B\iN\HomA(B{\otimes}B,B)$. The
final step follows from the first equality in \erf{eq:Phi1-rho} together
with the properties of $A$ to be symmetric and special.
\\[.3em]
It is easy to check that the morphisms defined in \erf{eq:lift-mult}
are elements of the relevant $\HomA$-spaces, i.e.\
$\tilde m_B\iN\HomA(B{\otimes_\AA}B,B)$ and $\tilde\eta_B\iN\HomA(A,B)$,
and analogously for $\tilde\Delta_B$ and $\tilde\eps_B$. Of the defining 
properties for $\tilde B$ to be a Frobenius algebra we will verify 
explicitly only associativity -- the other properties are checked analogously.
\\
Associativity is deduced as follows:
  \be \bearll
  \tilde m_B \circ ( \tilde m_B \otA \id_B ) \!\!
  &= m_B \circ e_{B\otimes B}^{} \circ r_{B\otimes B}^{} \circ (m_B\oti\id_B)
   \circ e_{B\otimes B\otimes B}^{}
  \\{}\\[-.6em]
  &= m_B \circ (m_B\oti\id_B) \circ e_{B\otimes B\otimes B}^{}
  = \,\cdots\, =
  \tilde m_B \circ (\id_B\otA\tilde m_B) \,.  \eear\ee
In the first step the definitions \erf{eq:morph-Atensor} 
and \erf{eq:lift-mult} are inserted; afterwards the idempotent
$e_{B\otimes B} \cir r_{B\otimes B} \eq P_{B\otimes B}$ is omitted,
which is allowed by assumption. Afterwards one can apply associativity
of $B$, and then the previous steps are followed in reverse order.
\\[.3em]
(ii)~Note that since $A$ is commutative and simple, by \Remark \ref{iZUV}(i) 
it is also haploid.
\\  
Out of the list of properties, let us look at specialness, commutativity
and haploidity as examples; the remaining cases are analysed similarly.
\\[.3em]
{\em Specialness:}
The first specialness relation for $\tilde B$ follows as
  %% [pic~48]
  \bea  \begin{picture}(140,52)(7,35)
  \put(0,0)  {\begin{picture}(0,0)(0,0)
              \scalebox{.38}{\includegraphics{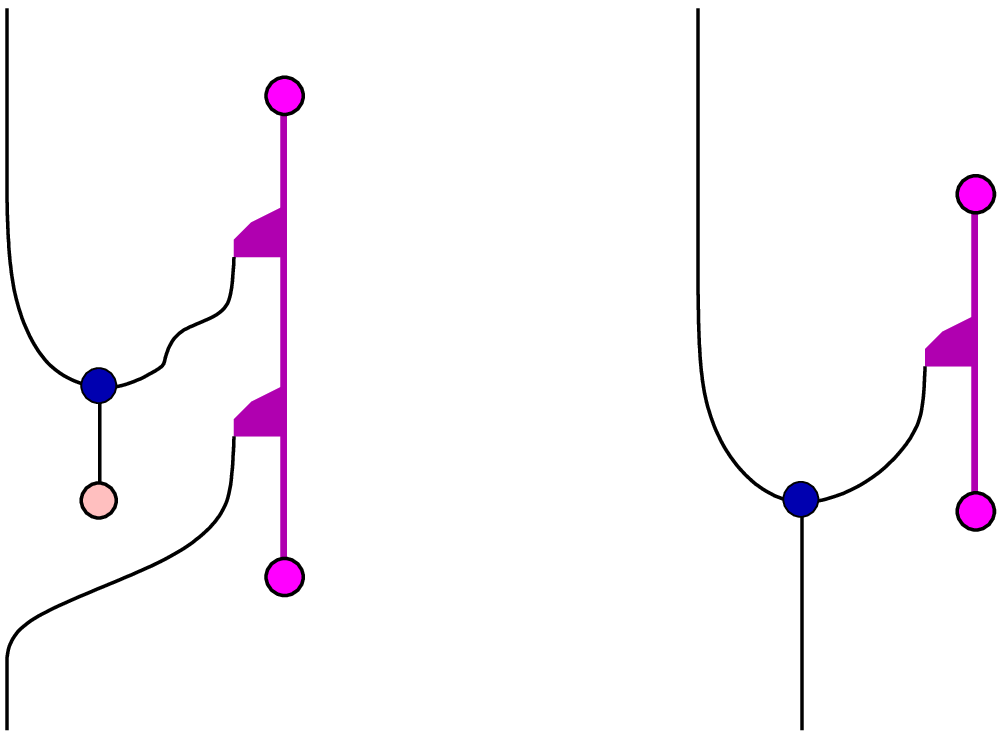}} \end{picture}}
  \put(-68,37.7)  {$\tilde\eps_B\circ\tilde\eta_B \,=$}
  \put(-2.9,-8.7) {\sse$A$}
  \put(-2.3,82.7) {\sse$A$}
  \put(32.6,24.7) {\sse$B$}
  \put(53,37.7)   {$=$}
  \put(73.8,82.7) {\sse$A$}
  \put(84.1,-8.7) {\sse$A$}
  \put(108.2,29.7){\sse$B$}
  \put(130,37.7)  {$=\;\Frac{\dim(B)}{\dim(A)}\,\id_\AA\,.$}
  \epicture19 \labl{eq:spec-eps-eta}
In the first step the definitions are substituted, while the second step
uses the representation property of $\r_B$ and the Frobenius property of
$A$. The resulting morphism is an element of $\Hom_A(A,A)$. Since $A$ is
haploid, this space is one-dimensional, so that the morphism must be
proportional to $\id_\AA$; comparing the traces determines the constant.
\\
The second specialness condition is implied by
  \be
  \tilde m_B^{} \circ \tilde \Delta_B^{} = m_B^{} \circ e_{B\otimes B}^{}
  \circ r_{\!B\otimes B}^{} \circ \Delta_B^{} = m_B^{} \circ \Delta_B^{}
  = \id_B^{} \,.  \ee
Here in the next to last step we used again that $m_B\cir 
P_{B\otimes B}\eq m_B$; the last equality holds because $B$ is special.
\\[.3em]
{\em Commutativity:}
When $B$ is commutative it follows directly from the form of the braiding in
$\Ext{\cC}{A}$ -- i.e.\ $c^A\eq r\cir c\cir e$ -- and from the definition 
\erf{eq:lift-mult} of $\tilde m_B$ that $\tilde B$ is commutative as well.
\\[.3em]
{\em Haploidity} of
$\tilde B$ is equivalent to $\dim\HomA(A,B)\eq1$. Since $A\eq \Ind_\AA(\one)$,
the reciprocity \erf{reciUM} implies $\,\dim\HomA(A,B)\eq\dim \Hom(\one,B)$.
If $B$ is haploid, then this equals 1, so that $\tilde B$ is haploid as well.
\qed

\medskip

The following assertion shows that for any simple \ssFA\ $A$, local induction
also supplies us with a functor from the category of Frobenius algebras in
$\cC$ to the category of Frobenius algebras in $\Ext{\cC}{C_{l/r}(A)}$.
\\[-2.3em]

\dtl{Proposition}{lem:[B]A-lift-i}
Let $A$ be a \ssFA\ and $B$ a Frobenius algebra in a ribbon category $\cC$,
and assume that the symmetric Frobenius algebras $C_l(A)$ and $C_r(A)$ are 
also special.
\\[.2em] 
(i)~\,The local $C_l(A)$\,-module $\lxtp B\AA l\eq(\efu Bl\AA,
\rholoc{C_l(A)}B)$ can be endowed with the structure of a Frobenius algebra 
in the category $\Ext{\cC}{C_l(A)}$ of local $C_l(A)$\,-modules.
\\[.2em]
(ii)~Let $A$ be in addition simple. If the Frobenius algebra
$\efu Bl\AA \iN \Obj(\cC)$ has any of the properties of being commutative, 
haploid, simple, symmetric, or special, then so has the Frobenius algebra
$\lxtp B\AA l \iN \Obj(\Ext{\cC}{C_l(A)})$.
\\[.2em]
Analogous statements apply to $C_r(A)$ and $\efu Br\AA$.

\medskip\noindent
Proof:\\
We show the claims for $C_l(A)$ and $\efu Bl\AA$; the corresponding
statements for $C_r(A)$ and $\efu Br\AA$ can be seen similarly.
The statements follow by applying \Lemma \ref{lem:alglift} 
to the Frobenius algebra $\efu Bl\AA$. Accordingly we must check that the 
requirements of that lemma are satisfied.  Abbreviate $C_l(A)$ by $C$.
Recall the definition \erf{eq:Elr-rep} of $\Rholoc$, which according to 
\Proposition \ref{pr:EAB-locmod} 
gives a local $C$-module structure on $\efu Bl\AA$. Furthermore, we have
  \be
  \rholoc CB \cir (\id_C\oti m) = m \cir (\rholoc CB\oti\id_{\efu Bl\AA}) \,,
  \labl{eq:m-HomA}
i.e.\ the multiplication $m$ of
$\efu Bl\AA$ is indeed in $\Hom_C(\efu Bl\AA\oti\efu Bl\AA,\efu Bl\AA)$.
To see this, we write out the definitions \erf{eq:Elr-rep} 
and \erf{eq:[B]A-alg} for the action of $C$ and the multiplication
on $\efu Bl\AA$, after which we can replace the resulting combination
$e \cir r$ by $P_\AA^l(B)$; then associativity of $A$ as well as the
properties \erf{oben2} and \erf{Cl-Cr-defprop} of the center $C$
relate the two sides of \erf{eq:m-HomA}.
\\
The equality $m\cir P_{\!\efu Bl\AA\otimes\efu Bl\AA}\eq m$ can be
verified in the same way, using in addition that $C$ is special.
\qed

Let us reformulate the statement of \Proposition \ref{lem:[B]A-lift-i}(i)
for later reference:
\\[-2.3em]

\dtl{Corollary}{deKac}
Let $A$ a be \ssFA\ such that $C_l(A)$ and $C_r(A)$ are special, and $B$ a 
Frobenius algebra, in a ribbon category $\cC$. 
Then there is a Frobenius algebra
  \be
  \LXTp\AA{l/r}(B) \;\in \Obj(\Ext{\cC}{C_{l/r}(A)}) \,.  \labl{nosep}
in the category of local $C_{l/r}(A)$-modules. 
The underlying object of the module $\LXTp\AA{l/r}(B)$ is $\efu B{l/r}\AA$.

\medskip

Note that we do not introduce a separate notation to indicate
the Frobenius algebra structure of the module \erf{nosep}.

For the following statement we take $A$ to be commutative, so that \erf{nosep} 
is now an algebra in the category of local $A$-modules, denoted by $\TildeB$.
\\[-2.3em]

\dtl{Proposition}{thm:[B]A-lift-ii}
Let $A$ and $B$ be commutative \ssFA s in a ribbon category $\cC$. Suppose in 
addition that $A$ is simple and that the Frobenius algebra $\Efu B\AA$ is 
special. Then $\TildeB$ is special, too, and we have an equivalence
  \be
  \EXt{\Ext{\cC}\AA}{\tildeB} \, \cong \, \Ext{\cC}{\Efu B\AA}  \labl{yW}
of ribbon categories.

\medskip\noindent
Proof:\\
By \Proposition \ref{prop:AB-alg}, $E_A(B)$ is a commutative symmetric 
Frobenius algebra. By assumption it is also special. Since $A$ is simple, by 
\Proposition \ref{lem:[B]A-lift-i}(ii) all properties of $E_A(B)$ get 
transported to $\TildeB$. In particular the three algebras $A$, $\TildeB$ 
and $E_A(B)$ are commutative \ssFA s, and by \Proposition \ref{thm:mod} 
all categories of local modules in \erf{yW} are ribbon categories.
\\[.2em]
The equivalence \erf{yW} is established by specifying two functors
  \be
  F :\ \ \EXt{\Ext{\cC}\AA}{\tildeB} \rightarrow \Ext{\cC}{\Efu B\AA}
  \qquad {\rm and} \qquad
  G :\ \ \Ext{\cC}{\Efu B\AA} \rightarrow \EXt{\Ext{\cC}\AA}{\tildeB}  \ee
and showing that they are each other's inverse and that they are ribbon.
\\[.3em]
{\em The functor $F$}:
An object $M$ in $\EXt{\Ext{\cC}\AA}{\tildeB}$ can be regarded as a
triple $(\M,\r^A,\r^{\tildeB})$ consisting of an object $\M$ in $\cC$,
a representation morphism $\r^A\,{\equiv}\,\r^A_M \iN \Hom(A{\otimes}\M,\M)$
that endows $(\M,\r^A)$ with the structure of a local $A$-module, and
a morphism $\r^{\tildeB}\,{\equiv}\,\r^{\tildeB}_M
       $\linebreak[0]$
\iN \HomA(\TildeB{\otimes_A}M,M)$ such that $(M,\r^{\tildeB})$ is a local
$\TildeB$-module in \calcal.
To define $F$ on objects we turn $M$ into a local $\Efu B\AA$-module by 
providing a morphism $\rho^{\Efu BA} \iN \Hom(\Efu B\AA{\otimes}M,M)$ which
has the appropriate properties; we set
  \be
  \rho^{\Efu B\AA} := \rho^{\tildeB} \circ r_{\!\Efu B\AA \otimes M}^{} \,.
  \labl{eq:CAB-rho[]}
(Recall from formula \erf{eq:morph-Atensor} that $r_{\!\Efu B\AA\otimes M}^{}$ 
is a short hand for $r_{\!\Efu B\AA \otimes M \succ \Efu B\AA \OtA M}^{}$; 
analogous abbreviations are implicit in $e_2$ and $e_3$ below.) To check 
the first representation property in \erf{1m} one computes -- abbreviating 
$\r \,{\equiv}\, \r^{\Efu B\AA}$, $m \,{\equiv}\, m^{\Efu B\AA}$, $\tilde\r 
\,{\equiv}\, \r^{\tildeB}$, $\tilde m \,{\equiv}\, m^{\tildeB}$  as well as
  $e_2 \,{\equiv}\, e^{}_{\Efu B\AA \otimes M}$,
  $e_3 \,{\equiv}\, e^{}_{\Efu B\AA \otimes
\Efu B\AA \otimes M}$ and similarly for $r_2$, $r_3$ -- as follows:
  \be \bearll
    \r \circ (\id_{\Efu B\AA}\oti\r) \!\!
  & \stackrel{(\rm a)}{=}
    \tilde\r \circ r_2 \circ P_{\Efu B\AA \otimes M} \circ
    ( \id_{\Efu B\AA}\oti\tilde\r ) \circ ( \id_{\Efu B\AA}\oti r_2 )
  \\{}\\[-.6em]
  & \stackrel{(\rm b)}{=}
    \tilde\r \circ r_2 \circ (\id_{\Efu B\AA}\oti\tilde\r)
    \circ (\id_{\Efu B\AA}\oti r_2) \circ e_3 \circ r_3
  \\{}\\[-.6em]
  & \stackrel{(\rm c)}{=}
    \tilde\r \circ (\id_{\Efu B\AA}\otA\tilde\r) \circ r_3
  \, \stackrel{(\rm d)}{=}
    \tilde\r \circ (\tilde m\otA\id_M) \circ r_3
  \\{}\\[-.6em]
  & \stackrel{(\rm e)}{=}
    \r \circ (m\oti\id_M) \,.
  \eear \labl{eq:CAB-rho[]-rep}
In step (a) definition \erf{eq:CAB-rho[]} of $\rho$ is substituted and
the idempotent $P_{\Efu B\AA \otimes M}\,{\equiv}\, P_2 \eq e_2\cir r_2
       $\linebreak[0]$
\iN\End(\Efu B\AA\oti M)$ is inserted before the second restriction
morphism $r_2$. Substituting the definition \erf{eq:P-2-tensor} 
for this idempotent, we see that it can be commuted past the first representation
and restriction morphisms $\tilde \rho$ and $r_2$, both these morphisms
being in $\HomA$, and afterwards due to the presence of
$r_2\eq r_2\cir P_2$ it can be replaced by
$P_{\Efu B\AA\otimes\Efu B\AA\otimes M} \,{\equiv}\, P_3 \eq e_3\cir r_3
       $\linebreak[0]$
\iN\End(\Efu B\AA\oti \Efu B\AA\oti M)$; this has been done in (b). In (c) the 
definition \erf{eq:morph-Atensor} of the tensor product over
$A$ for morphisms is substituted, while step (d) is the representation 
property of $\r^{\tildeB}$. Finally in (e) the tensor product over $A$ 
is replaced by \erf{eq:morph-Atensor}, the multiplication $\tilde m$ of 
$\TildeB$ expressed through \erf{eq:lift-mult} and the definition 
\erf{eq:mult-tens} substituted for the resulting $e_3\cir r_3$; then all 
$A$-ribbons can be removed, yielding the final expression in 
\erf{eq:CAB-rho[]-rep}. The second property in \erf{1m} can be checked 
similarly.
\\
Locality of the module $(M,\rho^{\Efu B\AA})$ is most easily verified
with the help of the condition (ii) in \Proposition \ref{pro:ostr}.
Indeed we have $\theta_M\cir\r\eq\theta_M\cir\tilde\r\cir r_2\eq\tilde\r
\cir (\id_{\Efu B\AA}\otA\theta_M)\cir r_2$, where the second step uses
locality of $M$ with respect to $\TildeB$. As a consequence,
$\theta_M\cir\r\eq \tilde\r\cir r_2\cir(\id_{\Efu B\AA}\oti\theta_M)\cir
P_{\Efu B\AA\otimes M} \eq \r\cir (\id_{\Efu B\AA}\oti\theta_M)$, where in the
first equality the morphism $\id_{\Efu B\AA}\otA\theta_M$ is substituted,
giving rise to the appearance of the idempotent $e_2 \cir r_2\eq
P_{\Efu B\AA \otimes M}$, while the second step uses locality of $M$ with
respect to $A$ to commute $\theta_M$ with the idempotent, which is then
omitted against $r_2$.
\\[.3em]
A morphism $f$ from $M$ to $N$ in $\EXt{\Ext{\cC}\AA}{\tildeB}$ is
an element of $\Hom(\M,\dot N)$ that commutes with the two actions
$\r^A$ and $\r^{\tildeB}$. The functor $F$ is defined to act as the
identity on morphisms: $F(f)\,{:=}\,f$. If $f$ commutes with $\r^A$ 
and $\r^{\tildeB}$, then it commutes with $\r^{\Efu B\AA}$ as well, 
because (using abbreviations similar to those in \erf{eq:CAB-rho[]-rep})
  \be \bearll
  f \cir \r_M^{}  \!\!
  &= f \cir \tilde \r_M^{} \cir r_{2,M}^{} \!\!
   = \tilde\r_N^{} \cir (\id_{\Efu B\AA}{\otimes_{\!A}}f) \cir r_{2,M}^{}
  \\{}\\[-.7em]
  &= \tilde\r_N^{} \cir r_{2,M}^{} \cir(\id_{\Efu B\AA}{\otimes}f)
    \cir P_{\Efu B\AA\otimes M}
  \\{}\\[-.7em]
  &= \tilde\r_N^{}\cir r_{2,M}^{}\cir P_{\Efu B\AA\otimes M}\cir
    (\id_{\Efu B\AA}{\otimes}f)
  = \r_N^{} \cir (\id_{\Efu B\AA} \oti f) \,.
  \eear \ee
In the second step the $\TildeB$-intertwiner property of $f$ is used.
The fact that $f$ is also in $\HomA$ allows one to commute it, in the fourth
step, with $P_{\Efu B\AA \otimes M}$.
\\[.3em]
{\em The functor $G$}:
We will be still more sketchy in the definition of $G$. On morphisms it
acts as the identity, $G(f)\,{:=}\,f$, just like $F$. To a local
$\Efu B\AA$-module $(M,\rho^{\Efu B\AA})$ it assigns the object
$G(M,\rho^{\Efu B\AA}) \,{:=}\, (M,\rho^A,\rho^{\tildeB})$
of $\EXt{\Ext{\cC}\AA}\tildeB$ as follows:
  \be\bearll
  \r^A := \r^{\Efu B\AA} \circ (e_{\Efu B\AA}\oti\id_M)
  \circ (\id_\AA \oti \eta^B \oti \id_M) \!\!
  &\in \Hom(A\oti M,M) \,, \\{}\\[-.6em]
  \r^{\tildeB} := \r^{\Efu B\AA} \circ e_{\tildeB\otimes M}
  &\in \Hom(\TildeB \otA M, M) \,.
  \eear\labl{eq:C[]-rhoAB}
To verify the representation property of $\r^{\tildeB}$ one needs
the relation
  \be
  \rho^{\Efu B\AA} \circ P_{\Efu B\AA\otimes M} = \r^{\Efu B\AA} \,,
  \labl{eq:rho-P=rho}
which can be seen by combining the definition \erf{eq:P-2-tensor}
of $P_{\Efu B\AA\otimes M}$ and of $\r^A$ in \erf{eq:C[]-rhoAB} with
the representation property of $\r^{\Efu B\AA}$ and the definition 
\erf{eq:[B]A-alg} of the product on $\Efu B\AA$. Using the condition of 
\Proposition \ref{pro:ostr}(ii) 
one can further convince oneself that $\r^A$ and $\r^{\tildeB}$ are 
local; we omit the calculation.
\\[.3em]
{\em $F$ and $G$ as inverse functors}:
$F$ and $G$ are clearly inverse to each other on morphisms. That
$F \cir G$ is the identity on objects follows from \erf{eq:rho-P=rho}.
To see $G\cir F\eq\Id$ on objects one must verify that
  \bea
  \r^A = \r^{\tildeB}
  \circ r_{\!\Efu B\AA \otimes M \succ \Efu B\AA \OtA M}^{}
  \circ (e^{}_{\Efu B\AA\prec A\otimes B} \oti\id_M)
  \circ (\id_\AA \oti\eta^B \oti\id_M) \,,
  \\{}\\[-.7em]
  \r^{\tildeB} = \r^{\tildeB}
  \circ r_{\!\Efu B\AA \otimes M \succ \Efu B\AA \OtA M}^{}
  \circ e_{\Efu B\AA \OtA M \prec \Efu B\AA \otimes M}^{} \,.  \eear\ee
The second equality is obvious. To see the first equality one replaces
$\id_\AA$ by $m \cir(\id_\AA \oti\eta^A)$ and uses the fact that all morphisms
are in $\HomA$ to trade the multiplication first for the representation
of $A$ on $\Efu B\AA$, then on $\Efu B\AA \otA M$, and finally on $M$.
The morphism $\rho^{\tildeB}$ is now applied to the unit of $\TildeB$
and can be left out. The remaining morphism is precisely $\rho^A$, the
action of $A$ on $M$.
\\[.3em]
{\em F as tensor functor:}
Denote by $\otimes_1$ the tensor product in $\EXt{\Ext{\cC}\AA}{\tildeB}$ 
and by $\otimes_2$ the tensor product in $\Ext{\cC}{\Efu B\AA}$. We need 
to show that $F(M\,{\otimes_1}\,N) \,{\cong}\, F(M)\,{\otimes_2}\,F(N)$; as 
we will see, the two objects are in fact equal. Since $F$ only changes the 
representation morphisms of $M$ and $N$, but not the underlying objects 
$\M$ and $\dot N$ we have (working with the Karoubian envelope, see
formula \erf{eq:multi-tensor-obj})
  \be
  F(M \,{\otimes_1} N) = \big( (\M\Oti\dot N ; P_1) , \rho_1 \big)
  \qquad {\rm and} \qquad
  F(M) \,{\otimes_2}\, F(N) = \big( (\M\Oti\dot N ; P_2) , \rho_2 \big)
  \,,  \ee
where $\M$ and $\dot N$ are objects in $\cC$ and $\r_{1,2}$ are 
representation morphisms for the algebra $E_A(B)$. Further, $P_1$ is the 
idempotent in $\End(\M\Oti\dot N)$ whose retract is $M\,{\otimes_1}\,N$, 
while $P_2$ gives the retract $F(M)\,{\otimes_2}\,F(N)$, i.e.
  \be
  P_1 = e \cir e' \cir r' \cir r \qquad {\rm and} \qquad
  P_2 = e'' \cir r'' \,, \ee
where the abbreviations $e\eq e_{M\OtA N \prec M\otimes N}$, 
$e'\eq e_{M \otimes_{\tildeB} N \prec M\OtA N}$,
$e''\eq e_{M\otimes_{E_A(B)}N \prec M\otimes N}$, as well as an analogous
notation for $r$, $r'$, $r''$ are used.
By direct substitution of the definitions one verifies that $P_1\eq P_2$.
It then  remains to compare the representation morphisms $\r_1$ and $\r_2$.
Again by substituting the definitions one finds that they are 
  \be
  \r_1 = \r_2 = (\r^{\tildeB}_M {\circ}\, r_{\!E_A(B)\oti M})
  \oti \id_N ~~\in \Hom(E_A(B) \oti M \oti N , M \oti N) \,.  \ee
{\em F as a ribbon functor:}
The duality and braiding are defined as in \erf{eq:216} and 
\erf{eq:217}, with the idem\-potents given by the idempotents \erf{eq:mult-tens} 
for the corresponding tensor products. But since the idempotents $P_{1,2}$ 
which define the retracts $M \,{\otimes_1} N \,{\prec}\, \dot M\oti \dot N$ and
$F(M)\,{\otimes_2}F(N) 
  $\linebreak[0]$ %%\,
{\prec}\,\M\oti\dot N$ are equal and $F$ acts as the 
identity on morphisms, duality and braiding of $\EXt{\Ext{\cC}\AA}{\tildeB}$
get mapped to duality and braiding of $\Ext{\cC}{\Efu B\AA}$.
\qed

\newpage  

\sect{Local modules and a subcategory of bimodules} \label{sect5n} 

The aim of this section is to establish -- in \Theorem \ref{thm:equiv} 
-- an equivalence between the three ribbon categories
$\cC_{C_{l}(A)}^{\sss\rmloc}$, $\cC_{C_{r}(A)}^{\sss\rmloc}$ and $\CAAo$.
Here $\CAAo$ denotes the full subcategory of \calcaa\ whose objects are 
those $A$-bimodules which are at the same time a sub-bimodule of an 
$\alpha_\AA^+$-induced and of an $\alpha_\AA^-$-induced bimodule.

To obtain this equivalence we introduce families of morphisms
in the category of left modules and in the category of bimodules
over a \ssFA. These families will be called pre-braidings.
The terminology derives from the fact that for left modules the
pre-braiding restricts to the braiding defined in \erf{eq:Cloc-braid}
if the algebra is commutative and the modules are local, while
for bimodules it gives rise to a braiding when restricted to $\CAAo$ 
(\Propositions \ref{pr:cmn=cloc} and \ref{pr:beta-braid}).

After discussing these preparatory concepts, a tensor functor from
$\cC_{C_{l/r}(A)}^{\sss\rmloc}$ to $\CAAo$ is constructed. Then it is first 
shown that this functor respects the braiding, and finally that it provides 
an equivalence, thus establishing the theorem.

\subsection{Braiding and left modules} \label{sect5n1}

Let $A$ be a \ssFA\ in a ribbon category $\cC$. If $A$ is in addition 
{\em commutative\/}, then one can define two tensor products $\ota^\pm$ on the 
\cat\ \calca\ of left $A$-modules, by extending the tensor product on its full 
sub\cat\ \calcal\ of local $A$-modules (see \Section \ref{sec:local-modules}) 
in two different ways. The basic ingredients are the idempotents introduced 
in \erf{eq:P-2-tensor}, i.e.
  \bea  \begin{picture}(320,50)(0,28)
  \put(63,0)  {\begin{picture}(0,0)(0,0)
              \scalebox{.38}{\includegraphics{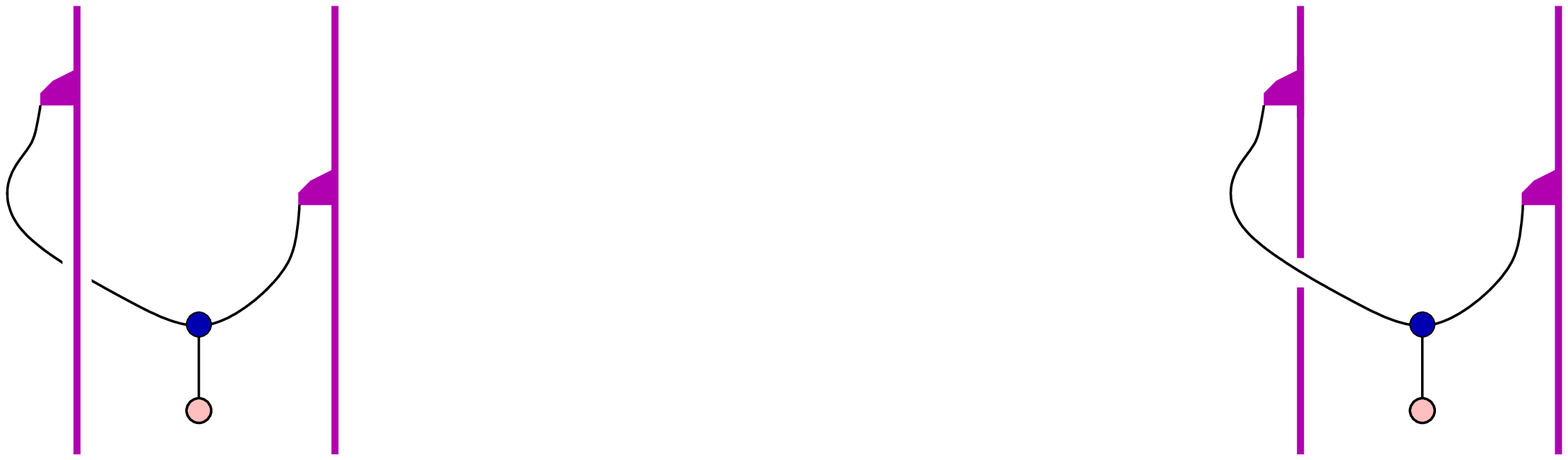}} \end{picture}}
  \put(0,30.5)    {$P_{\!M\ota^+ \!N} \,:= $}
  \put(69.5,-9.2) {\sse$\M$}
  \put(70.3,74.5) {\sse$\M$}
  \put(85.1,27.3) {\sse$A$}
  \put(110.2,-9.2){\sse$\dot N$}
  \put(110.9,74.5){\sse$\dot N$}
  \put(200,30.5)  {$P_{\!M \ota^- \!N} \,:= $}
  \put(99,0){
    \put(160.4,-9.2){\sse$\M$}
    \put(161.2,74.5){\sse$\M$}
    \put(176.1,27.3){\sse$A$}
    \put(201.3,-9.2){\sse$\dot N$}
    \put(201.8,74.5){\sse$\dot N$}
  }
  \epicture17 \labl{P-ota+-}
for any pair $M$, $N$ of $A$-modules. 

If $M$ is local, then $P_{\!M\ota^+\!N}\eq P_{\!M\ota^-\!N}\eq P_{\!M\OtA\!N}$ 
as defined in \erf{eq:P-2-tensor}, and one deals with tensor product $\OtA$ on
\calcal\ described in \Section \ref{sec:local-modules}.  In contrast, for 
general $A$-modules we get two distinct tensor products $\ota^\pm$. If, for 
$\nu\iN\{\pm\}$, the idempotent $P_{\!M\ota^\nu\!N}$ is split, we denote the 
associated retract by $(\Im P_{\!M\ota^\nu \!N},e^\nu_{M\otimes N},
r^\nu_{\!M\otimes N} )$, and thus the tensor product $\ota^\nu$ is given by
  \be
  M \,{\ota^\nu}\, N = \Im P_{M \ota^\nu N} 
  \qquad {\rm and} \qquad
  f \,{\ota^\nu}\, g = r^\nu_{\!M'\otimes N'} \cir
  (f \oti g) \cir e^\nu_{M\otimes N}  \labl{ota-nu}
for $M,M',N,N'\iN\Obj(\calca)$ and $f\iN\HomA(M,M')$, $g\iN\HomA(N,N')$.
If $P_{\!M\ota^\nu \!N}$ is not split, we must instead work with the
Karoubian envelope; then the same comments apply as in the case of 
\calcal\ that was discussed in \Section \ref{sec:local-modules}.

\medskip

When the \ssFA\ $A$ is {\em not\/} commutative, \calca\ is, in general, not a 
tensor category. However, we can still perform an operation that has
some similarity with a tensor product. This then allows us in particular to
introduce a `pre-braiding' on \calca\ that shares some properties of a genuine
braiding. To this end we restrict, for the moment, our attention to induced 
modules. For any pair $U,V$ of objects of $\cC$ we introduce the endomorphisms
  \be \! \bearl
  P_{\!\hatota^+}(U,V) := \big[ 
  (m\oti \id_U \oti \id_A) \cir (\id_A \oti c_{U,A}^{} \oti \id_A) \cir
  (\id_A \oti \id_U \oti \Delta) \big] \oti \id_V \quad\ {\rm and}
  \\{}\\[-.8em]
  P_{\!\hatota^-}(U,V) := \big[ 
  (m\oti \id_U \oti \id_A) \cir (\id_A \oti c_{\!A,U}^{-1} \oti \id_A) \cir
  (\id_A \oti \id_U \oti \Delta) \big] \oti \id_V 
  \eear \labl{hatota-nu}
in  $\End_\AA(\Ind_\AA(U\Oti A\Oti V))$, with $c$ the braiding on $\cC$.

\dt{Lemma}
The morphisms $P_{\!\hatota^\pm}(U,V)$ are split idempotents, with
image $\,\Ind_\AA(U\Oti V)$.

\medskip\noindent
Proof:\\
That $P_{\!\hatota^\pm}(U,V)$ are idempotents follows easily by using 
(co)associativity and specialness of $A$. To show that they are split, we 
just give explicitly the corresponding embedding and restriction morphisms 
$e_{UV}^\pm \eq e_{\hatota^\pm}(U,V) \iN \HomA(\Ind_\AA(U{\otimes}V),
\Ind_\AA(U{\otimes}A{\otimes}V))$ and $r_{\!UV}^\pm \eq r_{\!\hatota^\pm}(U,V)
\iN \HomA(\Ind_\AA(U{\otimes}A{\otimes}V),\Ind_\AA(U{\otimes}V))$: 
  \be\bearll
  e_{UV}^+ = \big[(\id_A \Oti c_{U,A}^{-1}) \cir (\Delta \Oti \id_U)
    \big] \oti \id_V \,, \quad &
  r_{\!UV}^+ = \big[(m\Oti \id_U)\cir (\id_A\Oti c_{U,A}^{}) 
    \big] \oti \id_V \,,
  \\[5pt]
  e_{UV}^- = \big[(\id_A \Oti c_{\!A,U}^{}) \cir (\Delta \Oti \id_U)
    \big] \oti \id_V 
  \,,  \quad &
  r_{\!UV}^- = \big[(m\Oti \id_U)\cir (\id_A\Oti c_{\!A,U}^{-1}) 
    \big] \oti \id_V \,.
  \eear\labl{eq:er-pm-def}
That $e_{UV}^\nu\cir r_{\!UV}^\nu\eq P_{\!\hatota^\nu}(U,V)$ is an immediate 
consequence of the Frobenius property of $A$. Further, as a result of 
specialness of $A$ the composition $r_{\!UV}^\nu \cir e_{UV}^\nu$ is equal to 
$\id_A{\otimes}\id_U{\otimes}\id_V$, hence the statement about the image.
\qed 

\medskip

The \retmodule s associated to the idempotents $P_{\!\hatota^\pm}(U,V)$
are used in
\\[-2.3em]

\dt{Definition}
The \operation s $\hatota^\nu{:}\; \calcai{\times}\calcai \,{\to}\, \calcai$
($\nu\iN\{\pm\}$) are given by
  \be
  \Ind_\AA(U) \,{\hatota^\nu}\, \Ind_\AA(V) := \Im P_{\!\hatota^\nu}(U,V)
  = ( \Ind_\AA(U\Oti V), e_{UV}^\nu, r_{UV}^\nu )  \ee
and
  \be
  f \,{\hatota^\nu}\, g := r_{\!U'V'}^\nu \cir (f\oti g) \cir e_{UV}^\nu  \ee
for $U,V,U',V'\iN\Objc$ and $f \iN \HomA(\Ind_\AA(U),\Ind_\AA(U'))$,
$g \iN \HomA(\Ind_\AA(V),\Ind_\AA(V'))$.

\bigskip 

In general, $(f_1{\hatota^\nu}g_1) \cir (f_2{\hatota^\nu}g_2)$ is not equal 
to $(f_1{\circ}f_2) \,{\hatota^\nu}\, (g_1{\circ}g_2)$, so that $\hatota^\nu$ 
is not a functor from $\calcai{\times}\calcai$ to \calcai, and hence in 
particular it is not a tensor product. However, for {\em commutative\/} 
algebras $\hatota^\nu$ does constitute a tensor product on \calcai. Indeed, 
the following statement can be verified by direct substitution of the 
respective definitions:
\\[-2.3em]

\dtl{Lemma}{lem:tensor=hattensor}
For every commutative \ssFA\ $A$ the \operation s $\hatota^\nu$ and $\ota^\nu$ 
coincide on $\calcai{\times}\calcai$, i.e.\ $\Ind_\AA(U)\,{\hatota^\nu}\,
\Ind_\AA(V)\eq\Ind_\AA(U)\,{\ota^\nu}\,\Ind_\AA(V)$ and
$f{\hatota^\nu}\, g \eq f {\ota^\nu}\,g$ for all $U,V,U'\!,V'\iN\Objc$ and all 
$f\iN\HomA(\Ind_\AA(U),\Ind_\AA(U'))$, $g\iN\HomA(\Ind_\AA(V),\Ind_\AA(V'))$.

\dtl{Definition}{def:b-fam}
Let $A$ be a (not necessarily commutative) \ssFA\ in a ribbon category $\cC$.
For $\mu,\nu\iN\{\pm\}$, we denote by $\Xb\mu\nu$ the family of morphisms
  \be
  \Xb\mu\nu_{UV} := \id_A \oti c_{U,V}^{} \qquad{\rm for}\quad U,V\iN\Objc
  \labl{Xb}
in $\HomA( \Ind_\AA(U){\hatota^\mu}\Ind_\AA(V) ,
\Ind_\AA(V) {\hatota^\nu} \Ind_\AA(U) )$. 

\bigskip

We will refer to the family $\Xb\mu\nu$, and likewise to
similar structures occurring below, as a {\em pre-braiding\/} on
\calcai. While $\Xb\mu\nu$ is itself not a braiding, it will 
give rise to one when restricted to a suitable sub\cat.

For the rest of this subsection we suppose that the \ssFA\ $A$ is
commutative. Then $\Xb\mu\nu$ can indeed be used to obtain a braiding
on the \cat\ \calcal\ of local $A$-modules, and this braiding coincides with 
the one already described in \erf{eq:Cloc-braid}. To obtain a statement about 
\calcal\ we must, however, get rid of the restriction to induced modules. To 
this end we recall from \Lemma \ref{lem:sub-of-ind}(ii)
that every $A$-module, and hence in particular every local $A$-module, is a
\retmodule\ of an induced module. Accordingly for each local $A$-module 
$M$ we select an object $U_M\iN\Objc$ such that $(M,e_M,r_M)$ is a module 
retract of $\,\Ind_\AA(U_M)$. Then for $\mu,\nu\iN\{\pm\}$ we define a family
$\Xgamma\mu\nu_{MN}$ of morphisms of \calcal\ by
  \be
  \Xgamma\mu\nu_{MN} := ( r_N \,{\ota^\nu}\, r_M ) \cir
  \Xb\mu\nu_{U_M\,U_N} \cir (e_M \,{\ota^\mu}\, e_N)  \labl{eq:gamdef}
for $M,N\iN\Obj(\calcal)$. Note that even though $\ota^\pm\eq\OtA$ for local 
modules, here we still must use the \operation\ $\ota^\pm$, because 
the induced module $\Ind_\AA(U_M)$ is not necessarily local, 
  % but compare  corollary \ref{cor:loc-lind} for the commutative case
so that e.g.\ the morphism $e_M \iN \HomA(M,\Ind_\AA(U_M))$ 
is, in general, only a morphism in \calca, but not in \calcal.

The following result implies that $\Xgamma\mu\nu$ does not depend on the 
particular choice of the triple $(U_M,e_M,r_M)$. It also establishes that 
$\Xgamma\mu\nu$ is actually independent of $\mu$ and $\nu$,
that it furnishes a braiding on \calcal, and that this braiding 
coincides with the braiding $c^A$ defined in \erf{eq:Cloc-braid}.
\\[-2.2em]

\dtl{Proposition}{pr:cmn=cloc}
Let $A$ be a commutative \ssFA\ and $M,N$ be local $A$-modules. Then 
  \be  \Xgamma\mu\nu_{MN} = c_{MN}^A  \ee
for $\mu,\nu\iN\{\pm\}$.

\medskip\noindent
Proof:\\
Writing out the definition of $\Xgamma\mu\nu_{MN}$ gives
  \be
  \Xgamma\mu\nu_{MN} = r_{\!N \otimes M} \cir ( r_{\!N} \oti r_{\!M} )
  \cir e_{\hatota^\nu} \cir (\id_A \oti c_{U_M,U_N}^{})
  \cir r_{\hatota^\mu} \cir ( e_M \oti e_N) \cir e_{M\otimes N} \,.
  \labl{eq:gamma-expl}
In the sequel we consider the case $\mu\eq{-}$, $\nu\eq{+}$ as an example.
(The other cases are verified similarly.)
In pictorial notation, formula \erf{eq:gamma-expl} is the first
equality in the following series of transformations:
  \begin{eqnarray} \begin{picture}(390,302)(9,0)
  \put(70,0)  {\begin{picture}(0,0)(0,0)
              \scalebox{.38}{\includegraphics{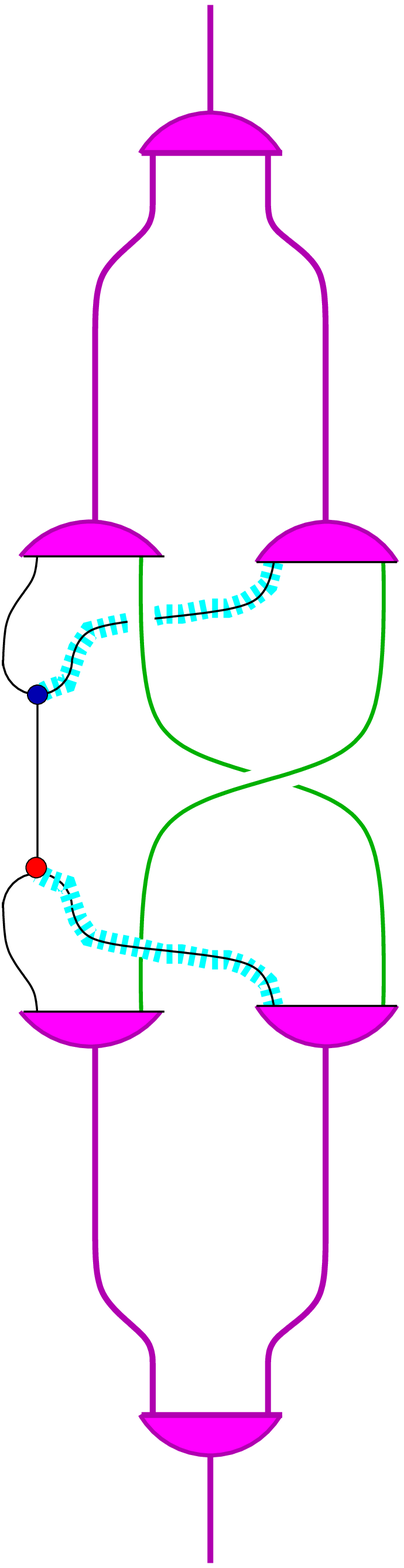}} \end{picture}}
  \put(210,0)  {\begin{picture}(0,0)(0,0)
              \scalebox{.38}{\includegraphics{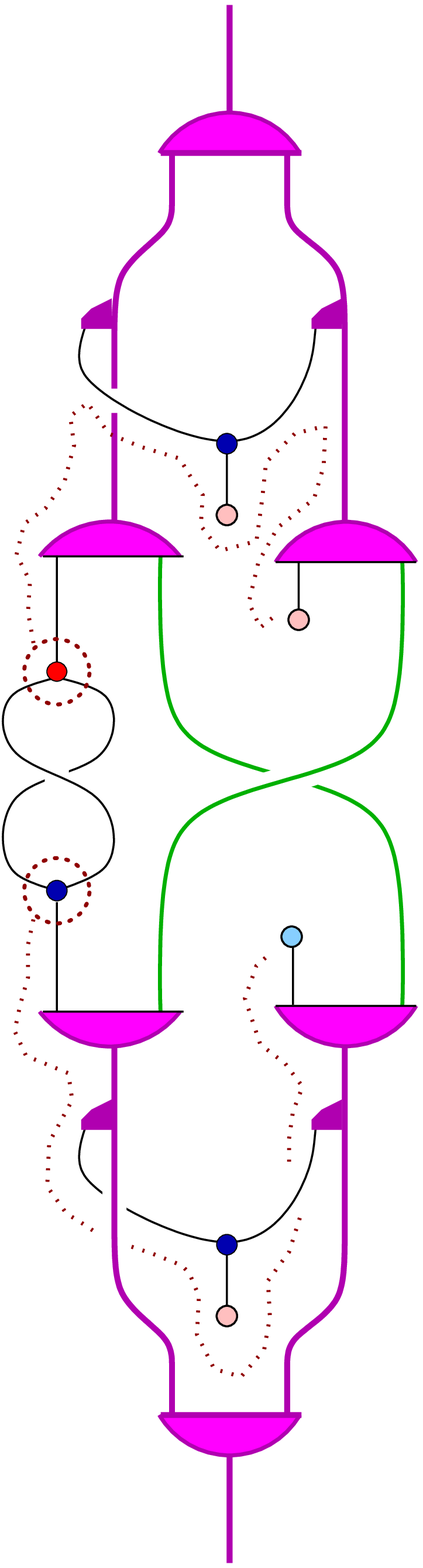}} \end{picture}}
  \put(350,0)  {\begin{picture}(0,0)(0,0)
              \scalebox{.38}{\includegraphics{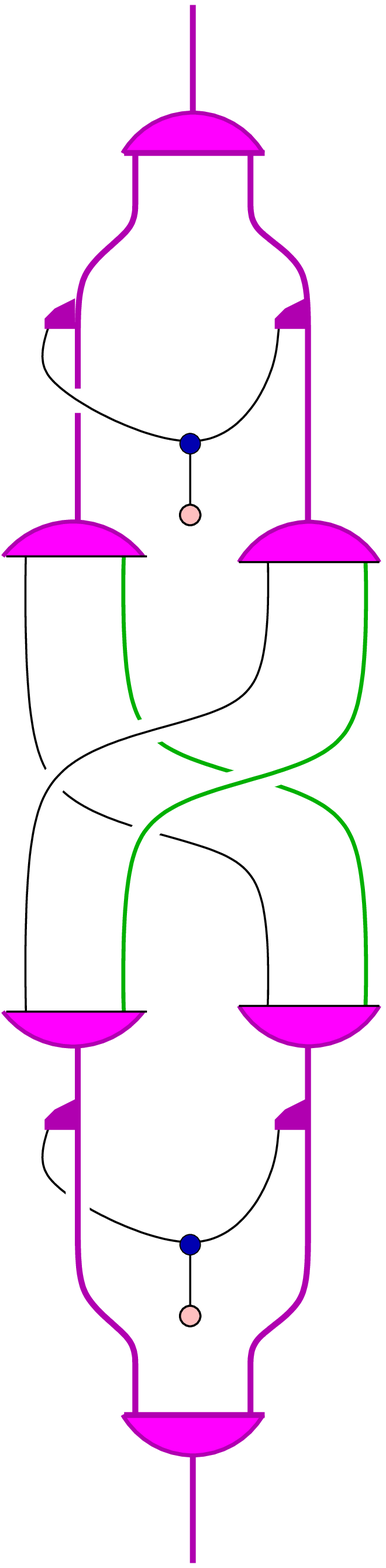}} \end{picture}}
  \put(0,145)      {$\Xgamma-+_{MN} \;= $}
  \put(69.3,152)   {\sse$A$}
  \put(77.4,66)    {\sse$M$}
  \put(78.5,224)   {\sse$N$}
  \put(97,-9.2)    {\sse$M{\OtA}N$}
  \put(98.8,298)   {\sse$N{\OtA}M$}
  \put(132.4,66)   {\sse$N$}
  \put(132.6,224)  {\sse$M$}
  \put(143.1,121)  {\sse$U^{}_{\!M}$}
  \put(143.1,169)  {\sse$U^{}_{\!N}$}
  \put(173,145)    {$=$}
  \put(207.1,146.6){\sse$A$}
  \put(232.1,33)   {\sse$M$}
  \put(233.2,256)  {\sse$N$}
  \put(241,-9.2)   {\sse$M{\OtA}N$}
  \put(242.8,298)  {\sse$N{\OtA}M$}
  \put(249.3,64)   {\sse$A$}
  \put(249.3,214)  {\sse$A$}
  \put(265.4,33)   {\sse$N$}
  \put(265.4,256)  {\sse$M$}
  \put(287.2,121)  {\sse$U^{}_{\!M}$}
  \put(287.2,169)  {\sse$U^{}_{\!N}$}
  \put(313,145)    {$=$}
  \put(347.2,178)  {\sse$A$}
  \put(347.2,114)  {\sse$A$}
  \put(365.1,33)   {\sse$M$}
  \put(366.3,256)  {\sse$N$}
  \put(373.3,-9.2) {\sse$M{\OtA}N$}
  \put(375.1,298)  {\sse$N{\OtA}M$}
  \put(382.9,64)   {\sse$A$}
  \put(382.9,214)  {\sse$A$}
  \put(398.7,33)   {\sse$N$}
  \put(398.7,256)  {\sse$M$}
  \put(420.3,121)  {\sse$U^{}_{\!M}$}
  \put(420.3,169)  {\sse$U^{}_{\!N}$}
  \end{picture} \nonumber\\[-.6em]{} \label{picf1}
  \\[-2.0em]{}\nonumber\end{eqnarray}
The second step of these manipulations involves a rewriting of
the marked $A$-ribbons as idempotents $P_{\!M \ota^\pm \!N}$, which
uses in particular that $A$ is commutative and that $M$ and $N$ are local. 
Furthermore, the identity $\id_A\eq m \cir c_{A,A}^{-1} \cir \Delta$,
which holds because $A$ is special and commutative, is inserted. In the 
last step, the marked multiplication and comultiplication morphisms
are dragged along the paths indicated (becoming representation morphisms 
for part of the way); this relies again on $A$ being commutative.
\\
In the final picture, the idempotents $P_{\!M \ota^\pm \!N}$ 
can be removed, while the morphisms $e_{M/N}$ and $r_{\!M/N}$
combine to the identity morphism on $M$ and $N$, respectively. Comparison
with \erf{eq:Cloc-braid} then shows that $\Xgamma-+_{MN}\eq c^A_{MN}$,
as claimed. 
\qed

%%%%%%%%%%%%%%%%%%%%%%%%%%%%%%%%%%%%%%%%%%%%%%%%%%%%%%%%%%%%%%%%%%%%%%%% 
\newpage

\subsection{Braiding and bimodules}\label{sect5n2}

{}From now on $A$ is again a general \ssFA, not necessarily commutative. 

The \cat\ \calcaa\ of $A$-bimodules contains interesting full sub\cats\
which were studied in \cite{boek2} and \cite{ostr}.

\dtl{Definition}{def:CAApm}
The full subcategories of \calcaa\ whose objects are the $\alpha_\AA^+$-induced 
and the $\alpha_\AA^-$-induced bimodules, respectively, are denoted by 
\CAAalp\ and \CAAalm, and their Karoubian envelopes by
  \be  \CAApm := \kar{(\CAAalpm)} \,.  \ee
The \cat\ \CAAo\ of {\em ambichiral\/} $A$-bimodules is the full subcategory 
of \calcaa\ whose objects are both in \CAAp\ and in \CAAm, i.e.
  \be  \CAAo := \CAAp \cap \CAAm \,.  \ee

\medskip

One can wonder whether the pre-braiding $\Xb\mu\nu$ on \calcai\ 
can be lifted to the bimodule \cat\ \calcaa. We will see that
this is indeed possible, by constructing families $\tildeb\mu\nu$ of
morphisms satisfying $R_A(\tildeb\mu\nu_{UV})\eq \Xb\mu\nu_{UV}$, where 
  \be  R_A:\quad \calcaa\,{\to}\;\calca  \labl{RAcc}
is the restriction functor whose action on objects consists in forgetting 
the right-action of $A$ on a bimodule. To do so first note that, as follows 
again by a straightforward application of the definitions, we have 
  \be  \alpha^\mu(U) \otA \alpha^\nu(V) 
  = (\alpha^\nu(U\Oti V), e_{UV}^\mu, r_{UV}^\mu) \,,  \ee
with $e_{UV}^\pm$ and $r_{UV}^\pm$ defined as in \erf{eq:er-pm-def},  
as a bimodule retract of $\alpha^\mu(U) \oti \alpha^\nu(V)$. 
To proceed we set
  \be
  \tildeb\mu\nu_{UV} := \id_A \oti c_{U,V}^{}  \labl{btilde}
for $\mu,\nu\iN\{\pm\}$ and $U,V\iN\Objc$ as in formula \erf{Xb}, but now 
regarded as morphisms from $\alpha_\AA^\mu(U) {\OtA} \alpha_\AA^\nu(V)$ to
$\alpha_\AA^\nu(V) {\OtA} \alpha_\AA^\mu(U)$. 
These families will again be called pre-braidings.

\dtl{Lemma}{lem:tildeb}
The pre-braidings $\tildeb\mu\nu_{UV}$ defined by \erf{btilde} have the
following properties.
\\[.2em]
(i)~\,For $(\mu\nu) \iN\{ (++), (+-), (--) \}$ they are bimodule morphisms,
i.e.
  \be
  \tildeb\mu\nu_{UV} \in \HomAA(
  \alpha_\AA^\mu(U) {\OtA} \alpha_\AA^\nu(V) , 
  \alpha_\AA^\nu(V) {\OtA} \alpha_\AA^\mu(U) ) \,.  \ee
(ii)~They fulfill
  \be  R_A(\tildeb\mu\nu_{UV}) = \Xb\mu\nu_{UV} \,,  \ee
with $R_A$ the restriction functor \erf{RAcc}. 

\medskip\noindent
Proof:\\
(i)~\,Compatibility of $\tildeb\mu\nu$ with the left action of $A$ is clear.
In the case of the right action $\rr^\pm$, given for $\alpha$-induced 
bimodules in \erf{rr+-}, we must show that
  \be
  \tildeb\mu\nu_{UV} \cir (\id_{\alpha_\AA^\mu(U)} \otA \rr^\nu(V))
  = (\id_{\alpha_\AA^\nu(V)} \otA \rr^\mu(U)) \cir 
  (\tildeb\mu\nu_{UV} \oti \id_A) \,.  \ee
Writing out the definitions, this amounts to 
  \bea
  (\id_A \oti c_{U,V}^{}) \cir r_{\!UV}^\mu \cir
  (\id_A \oti \id_U \oti \rr^\nu(V)) \cir ( e_{UV}^\mu \oti \id_A)
  \\[5pt] \quad\qquad =
  r_{\!VU}^\nu \cir (\id_A \oti \id_V \oti \rr^\mu(U)) \cir
  (e_{VU}^\nu \oti \id_A) \cir (\id_A \oti c_{U,V}^{} \oti \id_A) \,.
  \eear\labl{eq:R-al-al}
Inserting also the definitions of $\rr^\pm$, $e$ and $r$ 
one verifies, separately for each choice of 
$(\mu\nu) \iN\{ (++), (+-), (--) \}$, that this equality
follows from the properties of $A$ and of the braiding in $\cC$.
\\[.3em]
(ii)~For $\alpha$-induced bimodules we have
$R_A(\alpha^\pm(U))\eq\Ind_\AA(U)$, so that
  \be
  R_A(\alpha_\AA^\mu(U) {\OtA} \alpha_\AA^\nu(V)) = 
  R_A(\alpha_\AA^\nu(U\Oti V)) = \Ind_\AA(U\Oti V) \,.  \ee
Thus $R_A$ maps the source and target objects of $\tildeb\mu\nu_{UV}$ to 
those of $\Xb\mu\nu_{UV}$. As a consequence, the equality
  \be
  R_A(\tildeb\mu\nu_{UV}) = R_A(\id_A \Oti c_{U,V}^{}) = 
  \id_A \oti c_{U,V}^{} = \Xb\mu\nu_{UV}  \ee
follows immediately.
\qed

The morphisms $\tildeb\mu\nu_{UV}$ are not all functorial, as would be 
required for a braiding. But still we have the following properties.
\\[-2.2em]

\dtl{Lemma}{le:tb-prop}
For any $U,V,R,S\iN\Objc$ the following identities hold in \calcaa.
\\[.4em]
\begin{tabular}{lll}
(i) & $\tildeb++_{UV} \cir (\id {\OtA} g) = 
     (g{\OtA} \id) \cir \tildeb++_{US}$
   & for\, $g\iN\HomAA(\alpha_\AA^+(S), \alpha_\AA^+(V))$\,.
   \\{}\\[-.7em]
(ii) & $\tildeb--_{UV} \cir (f {\OtA} \id) = 
     (\id{\OtA} f) \cir \tildeb--_{RV}$
   & for\, $f\iN\HomAA(\alpha_\AA^-(R), \alpha_\AA^-(U))$\,.
   \\{}\\[-.7em]
(iii) & $\tildeb+-_{UV} \cir (\id {\OtA} g) = 
     (g{\OtA} \id) \cir \tildeb++_{US}$
   & for\, $g\iN\HomAA(\alpha_\AA^+(S), \alpha_\AA^-(V))$\,.
   \\{}\\[-.7em]
(iv) & $\tildeb+-_{UV} \cir (f {\OtA} \id) = 
     (\id{\OtA} f) \cir \tildeb--_{RV}$
   & for\, $f\in\HomAA(\alpha_\AA^-(R), \alpha_\AA^+(U))$\,.
   \\{}\\[-.7em]
(v) & $\tildeb+-_{UV} \cir (f {\OtA} g) = 
     (g{\OtA} f) \cir \tildeb+-_{RS}$
   & for\, $f\iN\HomAA(\alpha_\AA^+(R), \alpha_\AA^+(U))$
   \\{}\\[-.97em]
  && and\, $g\iN\HomAA(\alpha_\AA^-(S), \alpha_\AA^-(V))$\,.
\end{tabular}

\medskip\noindent
Proof:\\
The statements are all verified in a similar manner; we present the proof
of (iv) as an example. Substituting the definitions we find
  \bea  \begin{picture}(380,153)(0,55)
  \put(130,0)  {\begin{picture}(0,0)(0,0)
              \scalebox{.38}{\includegraphics{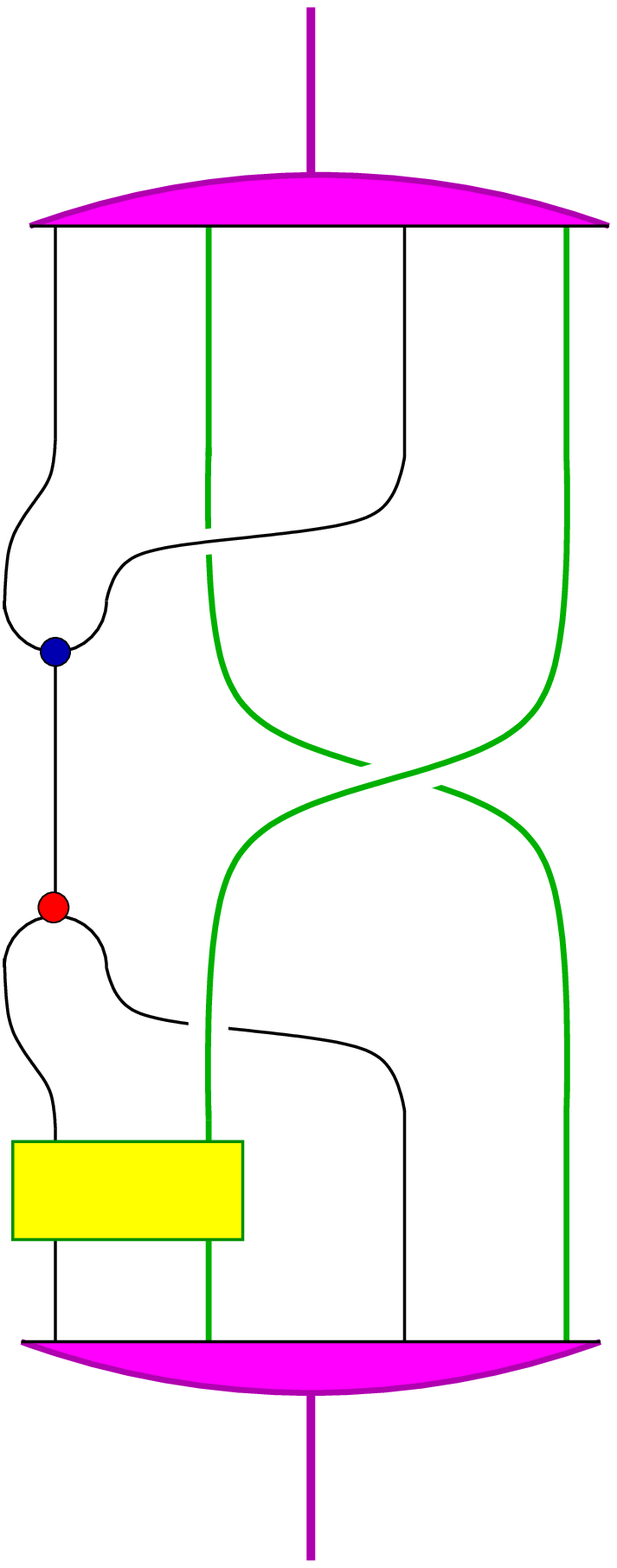}} \end{picture}}
  \put(270,0)  {\begin{picture}(0,0)(0,0)
              \scalebox{.38}{\includegraphics{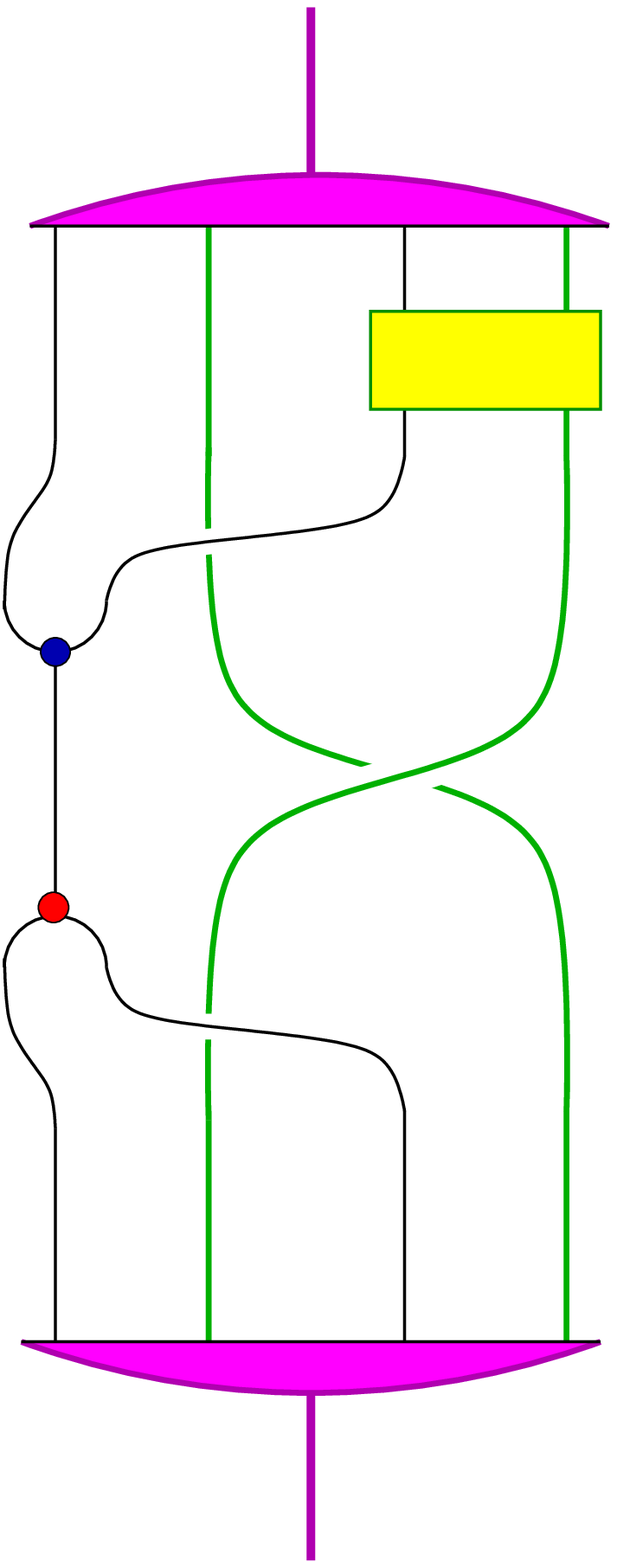}} \end{picture}}
  \put(5,97)       {$\tildeb+-_{UV} \cir (f {\OtA} \id) \,= $}
  \put(129.2,105)  {\sse$A$}
  \put(129.2,32.5) {\sse$A$}
  \put(141.4,-9.2) {\sse$\alpha_\AA^-(R){\OtA}\alpha_\AA^-(V)$}
  \put(141.4,205)  {\sse$\alpha_\AA^-(V){\OtA}\alpha_\AA^+(U)$}
  \put(143.8,46.1) {\sse$f$}
  \put(157.5,32.5) {\sse$R$}
  \put(203.2,32.5) {\sse$V$}
  \put(203.7,152)  {\sse$U$}
  \put(230,97)     {$=$}
  \put(268.7,100)  {\sse$A$}
  \put(281.4,-9.2) {\sse$\alpha_\AA^-(R){\OtA}\alpha_\AA^-(V)$}
  \put(281.4,205)  {\sse$\alpha_\AA^-(V){\OtA}\alpha_\AA^+(U)$}
  \put(328.8,151.3){\sse$f$}
  \put(297.7,32.5) {\sse$R$}
  \put(322.8,162.3){\sse$A$}
  \put(343.4,162.3){\sse$U$}
  \put(343.4,32.5) {\sse$V$}
  \epicture40 \labl{picf2}
In the first step the definition of $\tildeb+-$ and of the tensor product of 
morphisms is inserted. The second step uses first that the morphism $f$ 
intertwines the right action of $\alpha_\AA^-(R)$ and $\alpha_\AA^+(U)$ so 
as to take it past the multiplication, and next that it intertwines the 
left action (and hence, by the Frobenius property, the left co-action as
well) to commute it past the comultiplication. The resulting morphism
on the \rhs\ is equal to $(\id{\OtA} f) \cir \tildeb--_{RV}$.
\qed

So far we have a pre-braiding on the \cats\ \CAAalpm\
of $\alpha$-induced bimodules. We proceed to construct pre-braidings
$\Xbeta\mu\nu$ for \CAApm. 

\dtl{Definition}{def:beta-mu-nu}
Select, for each bimodule $X\iN\Obj(\CAAmu)$ and $\mu\iN\{\pm\}$,
an object $U^\mu_X\iN\Objc$ and morphisms $e^\mu_X$, $r^\mu_{\!X}$ such
that $(X,e^\mu_X,r^\mu_X)$ is a bimodule retract of 
$\alpha_\AA^\mu(U^\mu_X)$. 
Then for $X \iN \Obj(\CAAmu)$, $Y \iN \Obj(\CAAnu)$ and 
$(\mu\nu) \iN \{ (++), (+-), (--) \}$ the morphism $\Xbeta\mu\nu_{XY}$ 
is defined as
  \be
  \Xbeta\mu\nu_{XY} := (r^\nu_{\!Y} \,{\OtA}\, r^\mu_{\!X}) \cir 
  \tildeb\mu\nu_{U^\mu_X\,U^\nu_Y} \cir (e^\mu_X \,{\OtA}\, e^\nu_Y) \,.  \ee

\medskip

We will now show that
the families $\Xbeta\mu\nu$ of morphisms have similar properties as those 
of the pre-braidings $\tildeb\mu\nu$ that were listed in lemma 
\ref{le:tb-prop}. In particular, the morphisms $\Xbeta+-_{XY}$ turn out to 
be functorial and thus furnish a {\em relative braiding\/} between 
\CAAp\ and \CAAm, which coincides with the relative braiding introduced in 
\Proposition 4 of \cite{ostr}. Indeed we have
\\[-2.3em]

\dtl{Lemma}{le:beta-prop}
For any $X^\mu,Y^\mu,R^\mu,S^\mu\iN\Obj(\CAAmu)$ ($\mu\iN\{\pm\}$) the 
following identities hold in \calcaa.
\\[.4em]
\begin{tabular}{lll}
(i) & $\Xbeta ++_{XY} \cir (\id {\OtA} g) =
     (g{\OtA} \id) \cir \Xbeta ++_{XS}$  & for\, $g\iN\HomAA(S^+,Y^+)$\,.
   \\{}\\[-.7em]
(ii) & $\Xbeta --_{XY} \cir (f {\OtA} \id) =
     (\id{\OtA} f) \cir \Xbeta --_{RY}$  & for\, $f\iN\HomAA(R^-,X^-)$\,.
   \\{}\\[-.7em]
(iii) & $\Xbeta +-_{XY} \cir (\id {\OtA} g) =
     (g{\OtA} \id) \cir \Xbeta ++_{XS}$  & for\, $g\iN\HomAA(X^+,Y^-)$\,.
   \\{}\\[-.7em]
(iv) & $\Xbeta +-_{XY} \cir (f {\OtA} \id) =
     (\id{\OtA} f) \cir \Xbeta --_{RY}$  & for\, $f\in\HomAA(R^-,X^+)$\,.
   \\{}\\[-.7em]
(v) & $\Xbeta +-_{XY} \cir (f {\OtA} g) =
     (g{\OtA} f) \cir \Xbeta +-_{RS}$    & for\, $f\iN\HomAA(R^+,X^+)$
   \\{}\\[-.97em]
  && and\, $g\iN\HomAA(S^-, Y^-)$\,.
\end{tabular}
\\[.2em]
Here the abbreviations $\Xbeta ++_{XY}\eq\Xbeta ++_{\!\!X^+Y^+}$ etc.\ are used.

\medskip\noindent
Proof:\\
These properties of $\Xbeta\mu\nu$ are easily reduced to the corresponding
properties of $\tildeb\mu\nu$ in \Lemma \ref{le:tb-prop}.
Let us treat (i) as an example. Writing out the definition of $\Xbeta++$ on 
the \lhs\ of (i) gives (abbreviating also $r^+_{\!X}\eq r^+_{\!X^+}$ etc.)
  \be
  \Xbeta++_{XY} \cir (\id_{X^+} \otA g)
  = (r^+_{\!Y} \otA r^+_{\!X}) \cir 
  \tildeb++_{U_XU_Y} \cir (e^+_X \otA (e^-_Y\cir g)) \,,
  \labl{eq:lem-beta-1}
while for the \rhs\ we have
  \be
  (g \otA \id_{X^+}) \cir \Xbeta++_{XS} 
  = ((g\cir r^+_{\!S}) \otA r^+_{\!X}) \cir 
  \tildeb++_{U_X\,U_S} \cir (e^+_X \otA e^+_S) \,.  \labl{eq:lem-beta-2}
Since $e^+_Y {\OtA} e^+_X$ is monic and $r^+_{\!X}{\OtA} r^+_{\!S}$ is
epi, it is sufficient to show equality after composing the two expressions
\erf{eq:lem-beta-1} and \erf{eq:lem-beta-2} with $e^+_Y {\OtA} e^+_X$ from 
the left and with $r^+_{\!X} {\OtA} r^+_{\!S}$ from the right. The resulting
expressions are indeed equal, as is seen by using \Lemma \ref{le:tb-prop}(i)
with $\id_{X^+}\otA(e^-_Y\cir g\cir r^+_{\!S})$ in place of $\id\otA g$.
\qed

\medskip

The pre-braiding $\Xbeta\mu\nu$ gives rise to a braiding on \CAAo. 
The following observations will be instrumental to establish this result.
\\[-2.2em]

\dtl{Lemma}{le:be-pp-functor}
(i)~\,The morphisms $\Xbeta++$ satisfy
  \be  \Xbeta++_{XY} \cir (f \,{\OtA}\, g) = (g\,{\OtA}\, f) \cir \Xbeta++_{RS}
  \ee
for $X,R,S \iN \Obj(\CAAp)$, $Y \iN \Obj(\CAAo)$,
and $f\iN\HomAA(R,X)$, $g\iN\HomAA(S,Y)$.
\\
(ii)~The morphisms $\Xbeta--$ satisfy
  \be  \Xbeta--_{XY} \cir (f \,{\OtA}\, g) = (g\,{\OtA}\, f) \cir \Xbeta--_{RS}
  \ee
for $X \iN \Obj(\CAAo)$, $Y,R,S \iN \Obj(\CAAm)$,
and $f\iN\HomAA(R,X)$, $g\iN\HomAA(S,Y)$.
\\[.3em]
(iii)~When restricted to $\CAAp{\times}\,\CAAo$,
   the morphisms $\Xbeta++$ are functorial;
   when restricted to $\CAAm{\times}\,\CAAo$,
   the morphisms $\Xbeta--$ are functorial.

\medskip\noindent
Proof:\\
We establish (i); the proof of (ii) works analogously, while (iii) is an immediate
consequence of (i) and (ii). 
\\
By assumption on $Y$ 
there are bimodule retracts $(Y,e^+_Y, r^+_Y)$ of $\alpha_\AA^+(U^+_Y)$ and 
$(Y,e^-_Y, r^-_Y)$ of $\alpha_\AA^-(U^-_Y)$. Since 
$e^-_Y \iN \HomAA(Y,\alpha_\AA^-(U^-_Y))$ is a monic, it is
sufficient to verify that
  \be
  (e^-_Y \otA \id_X) \cir \Xbeta++_{XY} \cir (f \otA g) 
  = (e^-_Y \otA \id_X) \cir (g\otA f) \cir \Xbeta++_{RS} \,.  \ee
That this equality holds can be seen by using the properties of 
$\Xbeta\mu\nu$ established in \Lemma \ref{le:beta-prop}:
  \bea
  (e^-_Y \otA \id_X) \cir \Xbeta++_{XY} \cir (f \otA g) 
  \,\overset{\rm(iii)}{=}\,
  \Xbeta+-_{X\,\alpha^-(U^-_Y)} \cir (\id_X \otA e^-_Y) \cir (f \otA g) 
  \\[3pt] \qquad\qquad  \overset{\rm(v)}{=}\,
  (\id_{\alpha^-(U^-_Y)} \otA f) \cir \Xbeta+-_{R\,\alpha^-(U^-_Y)} 
    \cir (\id_R \otA e^-_Y) \cir (\id_R \otA g) 
  \\[3pt] \qquad\qquad  \overset{\rm(iii)}{=}\,
  (\id_{\alpha^-(U^-_Y)} \otA f) \cir (e^-_Y \otA \id_X)
    \cir \Xbeta++_{RY} \cir (\id_R \otA g) 
  \\[3pt] \qquad\qquad  \overset{\rm(i)}{=}\,
  (e^-_Y \otA \id_X) \cir (g\otA f) \cir \Xbeta++_{RS}  \eear\ee
(above the equality signs it is indicated which part of \Lemma 
\ref{le:beta-prop} is used).
\qed

\dtl{Proposition}{pr:beta-braid}
When restricting $\Xbeta\mu\nu$ with $(\mu\nu) \iN \{ (++), (+-), (--) \}$
to $\CAAo{\times}\,\CAAo$, we have:
\\[.2em]
(i)~\,\,The three families $\Xbeta\mu\nu$ coincide. Thus we can set
  \be  \Beta_{XY} :=  \Xbeta++_{XY} = \Xbeta+-_{XY} = \Xbeta--_{XY}  \ee
for all $X,Y\iN\Obj(\CAAo)$.
\\[.3em]
(ii)~\,The morphism $\Xbeta{}{}_{XY}$ is independent of the choices 
$e^\pm_{X,Y}$, $r^\pm_{\!X,Y}$ and $U^\pm_{X,Y}$ that are used in its definition.
\\[.3em]
(iii)~The family $\Beta$ of morphisms furnishes a braiding on \CAAo.

\medskip\noindent
Proof: \\
(i)~\,\,We demonstrate explicitly only the case $\Xbeta++_{XY}\eq\Xbeta+-_{XY}$;
the case $\Xbeta--_{XY}\eq\Xbeta+-_{XY}$ can be shown in the same way.
\\
We have $X,Y \iN \Obj(\CAAp)$, so there are bimodule retracts 
$(X,e^+_X, r^+_{\!X})$ of $\alpha_\AA^+(U^+_Y)$ and 
$(Y,e^+_Y, r^+_{\!Y})$ of $\alpha_\AA^+(U^+_Y)$. Furthermore 
$r^+_{\!X}{\OtA}r^+_{\!Y}$ is epi, so that it is sufficient to establish that
  \be
  \Xbeta++_{XY} \cir (r^+_{\!X}\otA r^+_{\!Y}) 
  = \Xbeta+-_{XY} \cir (r^+_{\!X}\otA r^+_{\!Y}) \,.
  \labl{eq:beta-braid-1}
Because of $Y\iN\Obj(\CAAo)$ we can apply lemma 
\ref{le:be-pp-functor}(i) to the \lhs, yielding
  \be
  \Xbeta++_{XY} \cir (r^+_{\!X}\otA r^+_{\!Y}) 
  = (r^+_{\!Y}\otA r^+_{\!X}) \cir 
  \Xbeta++_{\!\!\alpha_\AA^+(U^+_X)\,\alpha_\AA^+(U^+_Y)} \,.
  \labl{eq:beta-braid-2}
For the \rhs\ of \erf{eq:beta-braid-1} we get
  \be \bearll
  \Xbeta+-_{XY} \cir (r^+_{\!X}\otA r^+_{\!Y}) \!\!\!
  &= (\id_Y \otA r^+_X) \cir \Xbeta+-_{\alpha_\AA^+(U^+_X)\,Y} 
  \cir (\id_X \otA r^+_{\!Y})
  \\{}\\[-.7em]
  &= (r^+_{\!Y} \otA r^+_{\!X}) \cir 
  \Xbeta++_{\!\!\alpha_\AA^+(U^+_X)\,\alpha_\AA^+(U^+_Y)} \,,
  \eear\labl{eq:beta-braid-3}
where the first step amounts to \Lemma \ref{le:beta-prop}(v), 
while in the second step \Lemma \ref{le:beta-prop}(iii) 
is used, which is allowed because 
the source of the morphism $r^+_{\!Y} \iN \HomAA(\alpha_\AA^+(U^+_Y),Y)$ is 
in \CAAp\ and its target is in \CAAo\ and thus in particular in \CAAm.
\\
Comparing \erf{eq:beta-braid-2} and \erf{eq:beta-braid-3}
we see that \erf{eq:beta-braid-1} indeed holds true.
\\[.4em]
(ii)~\,is implied by (i). Indeed, $\Xbeta++_{X,Y}$ cannot depend on the 
choices of $e^+_{X/Y}$, $r^+_{\!X/Y}$ or $U^+_{X/Y}$, because $\Xbeta--_{X,Y}$
manifestly does not. Conversely, $\Xbeta--_{XY}$ must be independent of 
$e^-_{X/Y}$, $r^-_{\!X/Y}$ and $U^-_{X/Y}$. Likewise, 
since $\Xbeta+-_{XY}$ equals $\Xbeta++_{XY}$, it is independent
of the choices for $e^+_{X}$, $r^+_{\!X}$ and $U^+_{X}$,
and since it equals $\Xbeta--_{XY}$, it is independent
of the choices for $e^-_{Y}$, $r^-_{\!Y}$ and $U^-_{Y}$.
\\[.4em]
For the proof of (iii)
the tensoriality of the braiding -- the second line of formula \erf{DTB} --
must be verified. This can be done by direct computation. We do not present 
this calculation, but rather prefer to use a different argument later on, 
as part of the proof of \Theorem \ref{thm:equiv} in \Section \ref{sect5n3}.
\qed

%%%%%%%%%%%%%%%%%%%%%%%%%%%%%%%%%%%%%%%%%%%%%%%%%%%%%%%%%%%%%%%%%%%%%%%%

\subsection{A ribbon equivalence between local modules and ambichiral
            bimodules} \label{sect5n3}

Given a \ssFA\ $A$ and any pair $U$, $V$ of objects of a ribbon \cat\ $\cC$, 
define the linear maps $\Phi^{l/r}_{\!A;UV}$ by
  \be \begin{array}{lr}
  \Phi^{l/r}_{\!A;UV}:\ & \HomA(\Ind_\AA(U),\Ind_\AA(V)) \,\to\,
  \Hom(C_{l/r}\Oti U,C_{l/r}\Oti V)\;\, \\{}\\[-.8em]&
  f \,\mapsto\, (r_{\!C_{l/r}}\oti\id_V) \cir f \cir (e_{C_{l/r}}\oti\id_U)
  \,, \eear  \labl{eq:phi-map-def}
where $C_{l/r}$ stands for $C_l(A)$ and $C_r(A)$, respectively, and 
$r_{\!C_{l/r}}$ and $e_{C_{l/r}}$ are the restriction and embedding morphisms 
for the retract $C_{l/r}\,{\prec}\,A$. One checks that 
$\Phi^{l/r}_{\!A;UV}(f)$ commutes with the action of $C_{l/r}$, i.e.\ we have
  \be
  \Phi^{l/r}_{\!A;UV} :\quad \HomA(\Ind_\AA(U),\Ind_\AA(V))
  \longrightarrow \Hom_{C_{l/r}}(\Ind_{C_{l/r}}(U),\Ind_{C_{l/r}}(V)) \,.  \ee

\dt{Definition}
For $x\iN\{l,r\}$, the \operation s
  \be
  \Fi^{x}:\quad \calcai \to \cC_{C_{x}}^{\ind}  \labl{eq:Phi-fun-def}
are defined on objects as
  \be  \Fi^{x}(\Ind_\AA(U)) := \Ind_{C_x}(U)  \ee
for $U\iN\Objc$, and on morphisms as
  \be  \Fi^{x}(f) := \Phi^x_{\!A;UV}(f) \,,  \ee
with $\Phi^{x}_{\!A;UV}$ defined by \erf{eq:phi-map-def},
for $f\iN\HomA(\Ind_\AA(U),\Ind_\AA(V))$.

\medskip

The following properties of the maps $\Phi^{l/r}_{\!A;UV}$ are immediate 
consequences of the definitions. 
\\[-2.3em]

\dtl{Lemma}{lem:phi-prop}
The maps $\Phi^{x}_{\!A;UV}$ defined in \erf{eq:phi-map-def} fulfill
  \be  \Phi^{x}_{\!A;UU}(\id_A\Oti\id_U) = \id_{C_{x}} \oti \id_U  \ee
as well as
  \be
  \Phi^{x}_{\!A;VW}(g) \circ \Phi^{x}_{\!A;UV}(f) = 
  \Phi^{x}_{\!A;UW}(g{\circ}(P^x_\AA \Oti \id_V){\circ}f)  \labl{eq:phi-prop2} 
for $f\iN\HomA(\Ind_\AA(U),\Ind_\AA(V))$ and
$g\iN\HomA(\Ind_\AA(V),\Ind_\AA(W))$.

\medskip

As indicated by the appearance of the idempotent $P^{l/r}_\AA$ on the \rhs\ of 
\erf{eq:phi-prop2}, the \operation\ $\Fi^{l/r}$ is not a functor. However,
as will be seen below, $\Fi^{l/r}$ can be used to define a functor from
\CAApm\ to $\cC_{C_{l/r}}$.

\dtl{Lemma}{lem:F-b-to-b}
The \operation s $\Fi^{l/r}$ are compatible with the pre-braiding $\Xb{}{}$
in the sense that
  \be
  \Fi^l(\Xb++_{UV}) = \Yb++l_{UV} \qquad{\rm and}\qquad
  \Fi^r(\Xb--_{UV}) = \Yb--r_{UV}  \ee
for all $U,V\iN\Objc$.

\noindent
Proof:\\
As a straightforward application of the definitions, we have
  \be
  \Fi^l(\Xb++_{UV}) = 
  \Phi^l_{\!A;U{\otimes}V,V{\otimes}U}(\id_A \Oti c_{U,V})
  = (r_{C_l} \cir \id_A \cir e_{C_l}) \oti c_{U,V} = \Yb++l_{UV} \,.  \ee
and similarly for $\Fi^r(\Xb--_{UV})$.
\qed

\dtl{Lemma}{le:equiv-prep-1}
The map $\Phi^{l/r}_{\!A;UV}$ restricts as follows to bijections between 
spaces of bimodule morphisms of $\alpha$-induced $A$-bimodules and module
morphisms of (locally) induced $C_{l/r}$-modules:
  \bea
  \Phi^{l\,++}_{\!A;UV} :\quad \HomAA( \alpha_\AA^+(U), \alpha_\AA^+(V) ) 
    \,\overset{\cong}{\longrightarrow}\,
    \Hom_{C_l}(\Ind_{C_l}(U),\Ind_{C_l}(V)) \,, \\{}\\[-.7em]
  \Phi^{r\,--}_{\!A;UV} :\quad \HomAA( \alpha_\AA^-(U), \alpha_\AA^-(V) ) 
    \,\overset{\cong}{\longrightarrow}\,
    \Hom_{C_r}(\Ind_{C_r}(U),\Ind_{C_r}(V)) \,.
  \eear\labl{eq:equiv-prep-1}

\noindent
Proof:\\
This is a consequence of the \Proposition \ref{lem:[U]A-as-obj} 
together with the reciprocity relations (see 
\Remark \ref{prop:ssFA-unique}(iii))
  \be
  \Hom(C_{l/r}{\otimes}U,V) \,\cong\,
  \Hom_{C_{l/r}}(\Ind_{C_{l/r}}(U),\Ind_{C_{l/r}}(V)) \,.  \ee
Using the explicit form \erf{eq:reci-inv} and \erf{PhiPsi} of 
these isomorphisms, one can check that they are indeed given by
restrictions of the maps $\Phi^{l/r}_{\!A;UV}$.
\qed 

We now compose the \operation s $\Fi^{l/r}$ with the restriction functor $R_A$
\erf{RAcc}.
\\[-2.2em]

\dt{Definition}
For $A$ a \ssFA\ in a ribbon \cat, the \operation s
  \be
  G^\ind_+ :\quad \CAAalp \to \cC_{C_l}^\ind 
  \qquad {\rm and} \qquad
  G^\ind_- :\quad \CAAalm \to \cC_{C_r}^\ind  \ee 
are defined as the compositions $G^\ind_+ \,{:=}\, \Fi^l\cir R_A$ and
$G^\ind_- \,{:=}\, \Fi^r\cir R_A$ of the \operation s \erf{eq:Phi-fun-def}
with the restriction functor \erf{RAcc}.  

\dtl{Lemma}{lem:G-ind-tens}
Let $A$ be a \ssFA\ in a ribbon category $\cC$ such that the symmetric 
Frobenius algebras $C_l(A)$ and $C_r(A)$ are special. Then we have:
\\[.2em]
(i)~\,\,The \operation s $G^\ind_\pm$ are functors.
\\[.2em]
(ii)~\,They constitute tensor equivalences between 
the categories \CAAalpm\ and $\cC_{C_{l/r}}^\ind$.
\\[.2em]
(iii)~They satisfy
  \be
  G_+^\ind(\Xb++_{UV}) = \Yb++l_{UV} \qquad{\rm and}\qquad
  G_-^\ind(\Xb--_{UV}) = \Yb--r_{UV} \,.  \ee

\noindent
Proof:\\
We establish the properties for $G^\ind_+$; the proofs
for $G^\ind_-$ work analogously.
\\[.2em]
(i)~\,$G_+^{\ind}$ {\em is a functor\/}: Recall from the comment before 
\Lemma \ref{lem:F-b-to-b} that $G^\ind_+$ is not
a priori a functor, since $\Fi^l$ is not. However. after composition with
$R_A$ we have, for $f\iN\HomAA(\alpha_\AA^+(U),\alpha_\AA^+(V)$ and
$g\iN\HomAA(\alpha_\AA^+(V),\alpha_\AA^+(W)$,
  \be
  \Fi^l \cir R_A(g{\circ}f) = \Fi^l (g{\circ}f) = \Phi^l_{\!A;UW}(g{\circ}f)
  \labl{eq:G-ind-tens-1}
as well as
  \be \bearll
  \Fi^l(R_A(g)) \cir \Fi^l(R_A(f)) \!\!\!
  &= \Fi^l(g) \cir \Fi^l(f) =
  \Phi^l_{\!A;VW}(g) \cir \Phi^l_{\!A;UV}(f) 
  \\{}\\[-.7em]
  &= \Phi^l_{\!A;UW}(g \cir (P^l_\AA{\otimes}\id_V) \cir f)
  = \Phi^l_{\!A;UW}((P^l_\AA{\otimes}\id_W) \cir g \cir f) \,.
  \eear\ee
Here in the third step \Lemma \ref{lem:phi-prop} 
is used, and in the
last step the idempotent $P^l_\AA$ is moved past $g$, which is allowed by 
\Lemma \ref{++lemma}. Finally, when inserting \erf{eq:phi-map-def} for 
$\Phi^l_{\!A;UW}(\cdot)$, the idempotent $P^l_\AA$ can be left out because
of the presence of $r_{C_l}$, thereby yielding the \rhs\ of 
\erf{eq:G-ind-tens-1}.
Thus $G_+^{\ind}(g \cir f)\eq G_+^{\ind}(g) \cir G_+^{\ind}(f)$.
\\ 
That $G_+^{\ind}(\id)\eq\id\,$ follows again from \Lemma \ref{lem:phi-prop}.
\\[.3em]
(ii)~\,$G_+^{\ind}$ {\em is an equivalence functor\/}: 
Clearly $G_+^{\ind}$ is essentially surjective on objects. Further,
by the first equivalence in \Lemma \ref{le:equiv-prep-1},
$G_+^{\ind}$ is an isomorphism on morphisms,
    \be
    G_+^{\ind} :\quad \HomAA(\alpha_\AA^+(U),\alpha_\AA^+(V)) 
    \overset{\cong}{\longrightarrow}
    \Hom_{C_l} (\Ind_{C_l}(U),\Ind_{C_l}(V)) \,.  \ee
Thus by the criterion of \Proposition \ref{XI.1.5}, 
$G_+^{\ind}$ is an equivalence functor.
\\[.3em]
$G_+^{\ind}$ {\em is a tensor functor\/}:
Using equation \erf{eq:R-al-al} we have
  \be \bearll
  G^{\ind}_+( \alpha_\AA^+(U) \otA \alpha_\AA^+(V) ) \!\!
  &= \Fi^l( \Ind_\AA(U{\otimes}V) ) 
  \\{}\\[-.7em]
  &= \Ind_{C_l}(U{\otimes}V) = 
  \Ind_{C_l}(U) \,{\otimes_{C_l}}\, \Ind_{C_l}(V) \,.  \eear \labl{GFIII}
The \rhs\ of \erf{GFIII} is equal to 
$G^{\ind}_+(\alpha_\AA^+(U)) \,{\otimes_{C_l}}\, G^{\ind}_+(\alpha_\AA^+(V))$. 
Thus $G^{\ind}_+$ is tensorial on objects. 
\\
For $f \iN \HomAA(\alpha_\AA^+(U), \alpha_\AA^+(U'))$ and
$g \iN \HomAA(\alpha_\AA^+(V), \alpha_\AA^+(V'))$ the morphisms
$G^{\ind}_+(f{\OtA}g)$ and $G^{\ind}_+(f) \,{\otimes_{C_l}}\, G^{\ind}_+(g)$
read
  \bea  \begin{picture}(260,152)(0,50)
  \put(105,0) {\begin{picture}(0,0)(0,0)
              \scalebox{.38}{\includegraphics{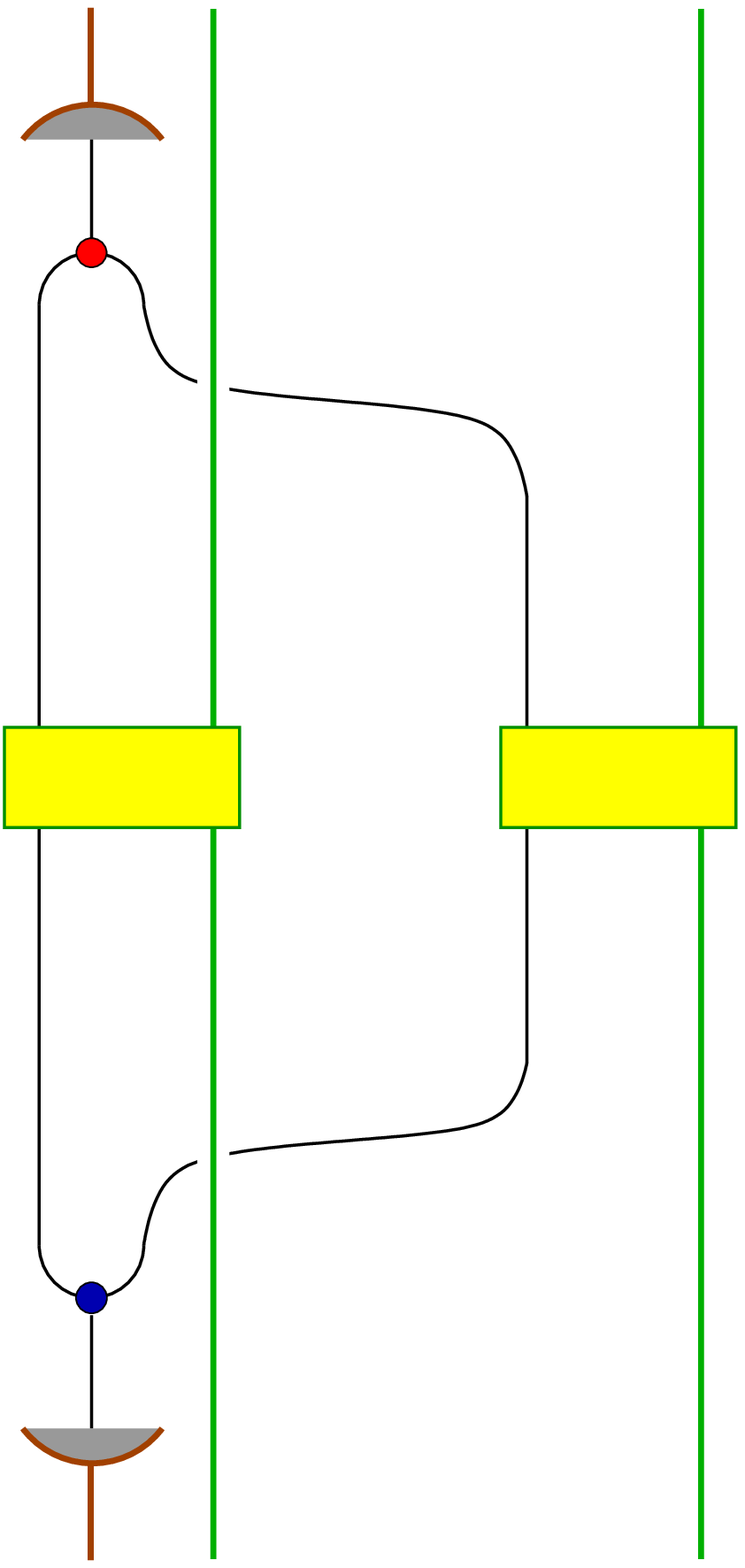}} \end{picture}}
  \put(0,97)       {$G^{\ind}_+(f \otA g) \,= $}
  \put(112.2,-8.3) {\sse$C_l$}
  \put(112.2,200)  {\sse$C_l$}
  \put(112.8,156.2){\sse$A$}
  \put(113.2,39.2) {\sse$A$}
  \put(116.7,97.6) {\sse$f$}
  \put(128.2,-8.3) {\sse$U$}
  \put(128.2,200)  {\sse$U'$}
  \put(179.2,98.7) {\sse$g$}
  \put(188.2,-8.3) {\sse$V$}
  \put(188.2,200)  {\sse$V'$}
  \epicture29 \labl{eq:G-ind-tens-2}
and
  \bea  \begin{picture}(320,156)(0,53)
  \put(165,0) {\begin{picture}(0,0)(0,0)
              \scalebox{.38}{\includegraphics{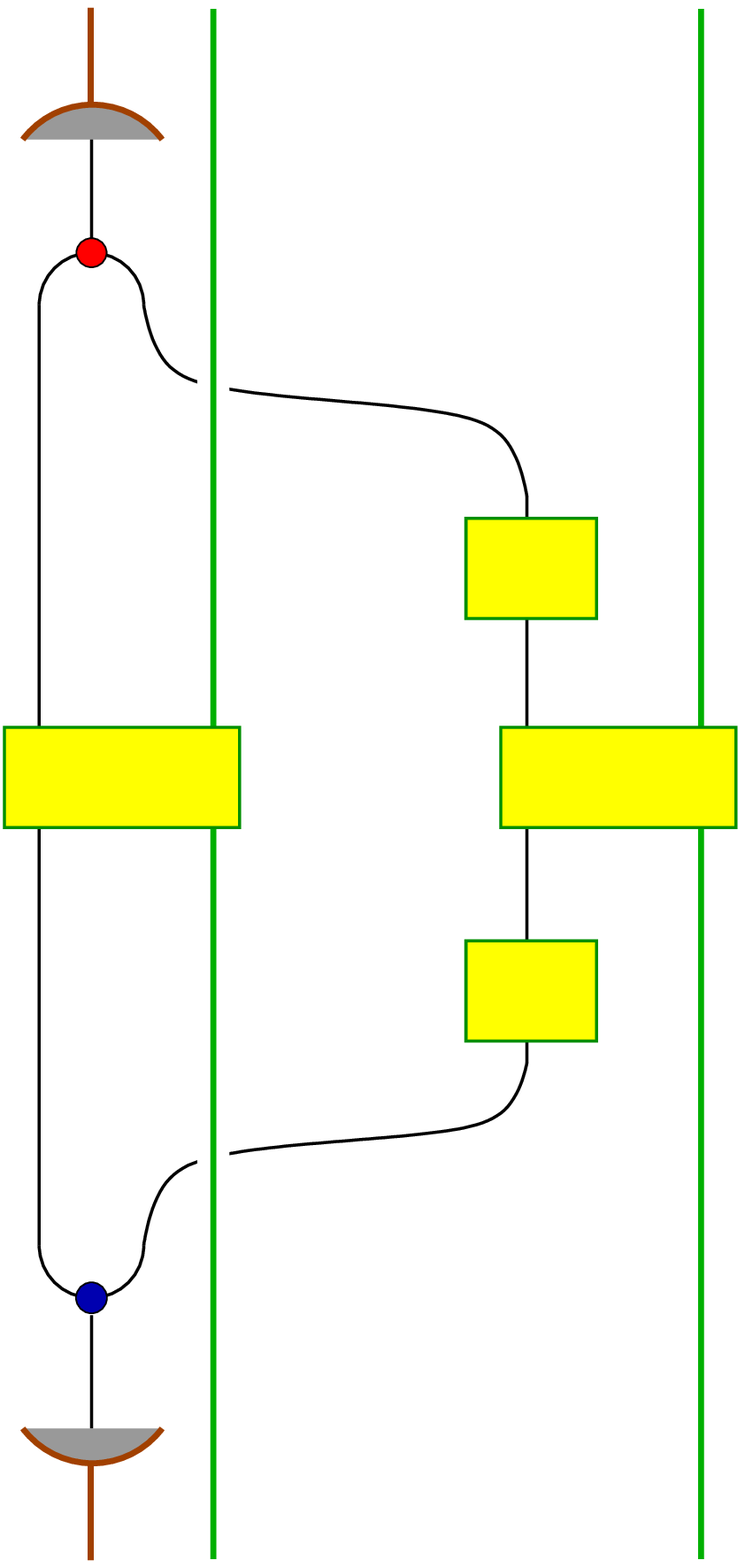}} \end{picture}}
  \put(0,95)       {$G^{\ind}_+ \,{\otimes_{C_l}}\, G^{\ind}_+(g) \,=\;
                   \Frac{\dim(A)}{\dim(C_l)} $}
  \put(172.2,-8.3) {\sse$C_l$}
  \put(172.2,200)  {\sse$C_l$}
  \put(172.8,156.2){\sse$A$}
  \put(173.2,39.2) {\sse$A$}
  \put(176.7,97.6) {\sse$f$}
  \put(188.2,-8.3) {\sse$U$}
  \put(188.2,199)  {\sse$U'$}
  \put(225.9,71.4) {\sse$P^l_\AA$} 
  \put(225.9,124)  {\sse$P^l_\AA$}
  \put(239.2,98.7) {\sse$g$}
  \put(248.9,-8.3) {\sse$V$}
  \put(249.2,200)  {\sse$V'$}
  \epicture40 \labl{eq:G-ind-tens-3}
In \erf{eq:G-ind-tens-3} the definition of the (co)multiplication on $C_l$ 
has been substituted and \Lemma \ref{lem:C=[1]a}(iii) has been used to omit 
one of the two resulting idempotents $P^l_\AA$ at $m$ and $\Delta$.
\\
To see that \erf{eq:G-ind-tens-2} and \erf{eq:G-ind-tens-3} are equal we 
consider the following identity, which can be obtained by dragging the marked 
multiplication along the path indicated. In order to do so one first uses that
$f\iN \HomAA(\alpha_\AA^+(U),\alpha_\AA^+(U'))$ (applied to the right action
of $A$) and next that $g \iN \HomAA(\alpha_\AA^+(V), \alpha_\AA^+(V'))$
(applied to the left action of $A$).
  \bea  \begin{picture}(390,127)(0,42)
  \put(70,0)  {\begin{picture}(0,0)(0,0)
              \scalebox{.38}{\includegraphics{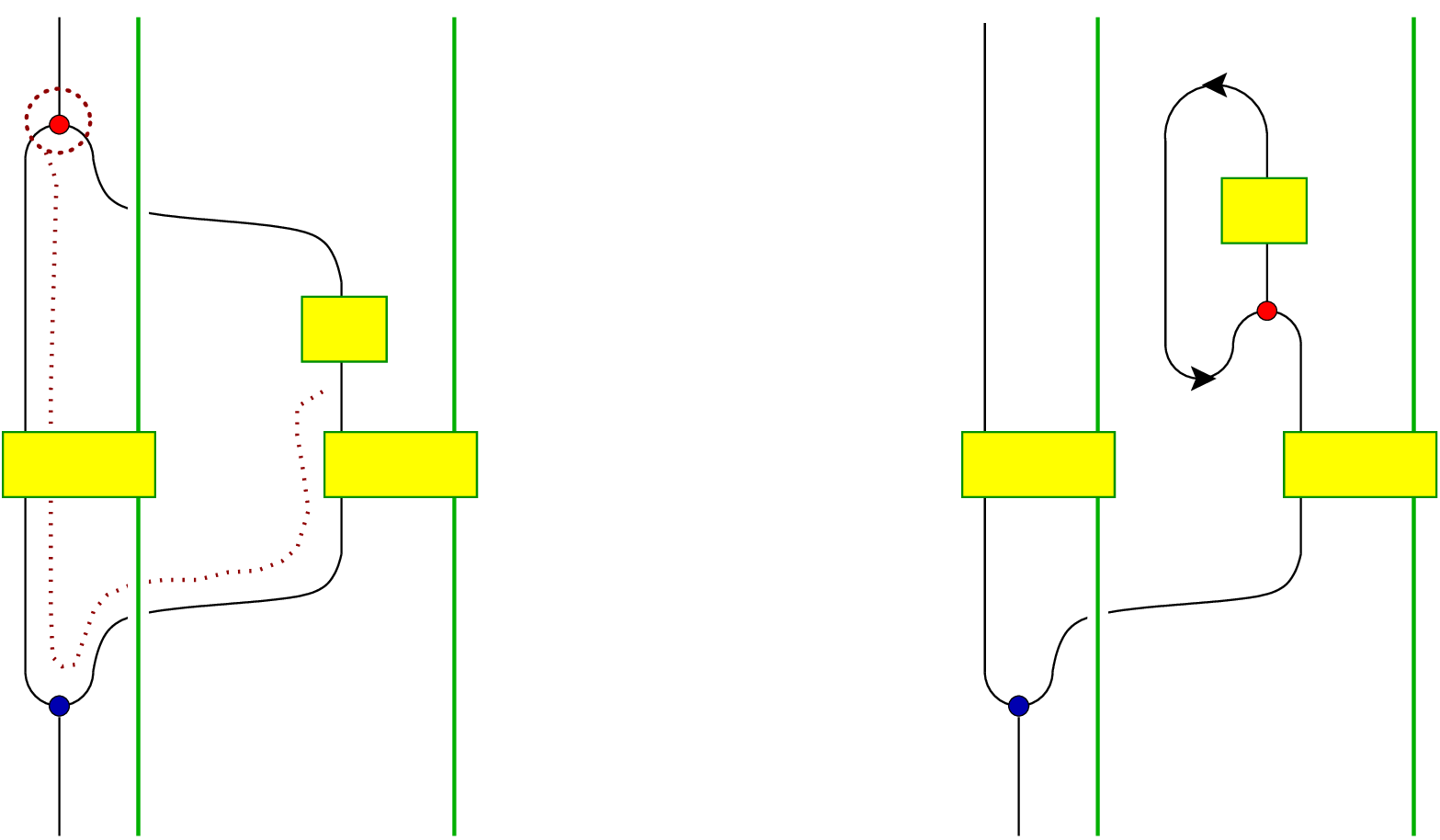}} \end{picture}}
  \put(0,76)       {$h(q) \,:= $}
  \put(77.2,-9.2)  {\sse$A$}
  \put(77.9,161.4) {\sse$A$}
  \put(82.2,69.6)  {\sse$f$}
  \put(93.2,-9.2)  {\sse$U$}
  \put(93.2,161.4) {\sse$U'$}
  \put(133.4,96.7) {\sse$q$}
  \put(143.2,70.7) {\sse$g$}
  \put(153.6,-9.2) {\sse$V$}
  \put(153.9,161.4){\sse$V'$}
  \put(203,76)     {$=$}
  \put(255.4,161.4){\sse$A$}
  \put(261.2,-9.2) {\sse$A$}
  \put(266.7,69.6) {\sse$f$}
  \put(277.7,-9.2) {\sse$U$}
  \put(277.7,161.4){\sse$U'$}
  \put(310.3,119.2){\sse$q$}
  \put(308.4,92.3) {\sse$A$}
  \put(327.4,70.7) {\sse$g$}
  \put(337.6,-9.2) {\sse$V$}
  \put(337.9,161.4){\sse$V'$}
  \epicture30 \labl{picf5}
For $q\eq\id_A$, the $A$-loop on the \rhs\ is equal to the counit $\eps_\AA$. 
On the other hand, for $q\eq P^l_\AA$, \Lemma \ref{lem:C=[1]a}(iii) 
allows us to replace the $A$-loop by a $C_l$-loop, which by specialness of $C_l$
is equal to the counit of $C_l$ and a restriction to
$C_l$, i.e.\ to replace the $A$-loop by $\eps_{C_l} \cir r_{\!C_l}$.
The latter, in turn, is equal to $\dim(C_l)/{\dim}(A)\,\eps_\AA$. Thus 
  \be  h(\id_A) = \Frac{\dim(C_l)}{\dim(A)} \, h(P^l_\AA) \,.  \labl{h-h}
Now the \rhs\ of \erf{eq:G-ind-tens-2} is equal to
$(r_{C_l}{\otimes}\id_{U'}{\otimes}\id_{V'}) \cir h(\id_A)
\cir (e_{C_l}{\otimes}\id_U{\otimes}\id_V)$, while the \rhs\ of 
\erf{eq:G-ind-tens-3} equals -- after eliminating one of the two
idempotents with the help of \Lemma \ref{++lemma} 
-- $\dim(A)/{\dim}(C_l)\, (r_{C_l}{\otimes}\id_{U'}{\otimes}\id_{V'}) 
\cir h(P^l_\AA) \cir (e_{C_l}{\otimes}\id_U{\otimes}\id_V)$. 
Hence the equality \erf{h-h} implies that
$G^{\ind}_+(f{\OtA}g)\eq G^{\ind}_+(f)\,{\otimes_{C_l}}\,G^{\ind}_+(g)$.
\\[.3em]
(iii)~$G_+^{\ind}$ {\em is compatible with\/} $\Xb{}{}$: 
The equality
  \be
  G_+^{\ind}( \tildeb\mu\nu_{UV} )
  = \Fi^l(R_A(\tildeb\mu\nu_{UV}))
  = \Fi^l(\Xb\mu\nu_{UV}) = \Xb\mu\nu_{UV}  \ee
follows by just combining \Lemma \ref{lem:tildeb}(ii) and \Lemma 
\ref{lem:F-b-to-b}.
\qed

Via Karoubification the functors $G^\ind_\pm$ induce functors
  \be
  G_+ :\quad \CAAp \to \cC_{C_l} \qquad {\rm and} \qquad
  G_- :\quad \CAAm \to \cC_{C_r} \,.  \ee 

\dtl{Proposition}{pr:G-tens}
The functors $G_\pm$ are tensor equivalences and satisfy
  \be \bearll
  G_+(\Xbeta++_{XY}) = \Ygamma++l_{\!G_+(X)\,G_+(Y)}
  \quad {\rm for}\ X,Y \iN \Obj(\CAAp)
  &\quad {\rm and} \\{}\\[-.5em]
  G_-(\Xbeta--_{XY}) = \Ygamma--r_{\!G_-(X)\,G_-(Y)}
  \quad {\rm for}\ X,Y \iN \Obj(\CAAm) \,.  \eear \labl{GXY}

\medskip\noindent
Proof:\\
By \Proposition \ref{sub-of-ind-ii} 
we have $\cC_{C_l}\,{\cong}\,\kar{(\cC_{C_l}^{\ind})}$. That $G_\pm$ is a tensor 
equivalence then follows from the corresponding property of $G^{\ind}_\pm$ 
established in \Lemma \ref{lem:G-ind-tens} by invoking \Lemma \ref{le:K-funct}.
\\[.2em]
The proof of the property \erf{GXY} will be given for $G_+$ only, the
one for $G_-$ being analogous.
Using the realisation of the Karoubian envelope via idempotents, let 
$X\eq(\alpha_\AA^+(U_X^+);p^+_X)$ and $Y\eq(\alpha_\AA^+(U_Y^+);p^+_Y)$. Then
  \be
  \Xbeta++_{XY} = (p^+_Y \otA p^+_X) \cir \Xb++_{U_X^+\,U_Y^+}
  \cir (p^+_X \otA p^+_Y) \,.  \ee
Also, if $M\eq(\Ind_{C_l}(U);p)$ and $N\eq(\Ind_{C_l}(V);q)$ are
objects in $\kar{(\cC_{C_l}^{\ind})}$, then
  \be
  \Xgamma++_{MN} = (q \,{\otimes_{C_l}}\, p) \cir \Xb++_{UV}
  \cir (p \,{\otimes_{C_l}}\, q) \,.  \ee
By definition, $G_+(X)\eq(G_+^{\ind}(\alpha_\AA^+(U_X^+));G_+^{\ind}(p_X^+))$;
the desired property of $G_+$ thus follows from the equalities
  \be \bearll
  G_+(\Xbeta++_{XY}) \!\!
  &= G_+^{\ind}((p^+_Y{\OtA}p^+_X) \cir \Xb++_{U_X^+\,U_Y^+}
  \cir (p^+_X{\OtA}p^+_Y))
  \\{}\\[-.7em]
  &= [G_+^{\ind}(p^+_Y) \,{\otimes_{C_l}}\, G_+^{\ind}(p^+_X)] \cir
  \Ygamma++l_{U_X^+\,U_Y^+} \cir
  [G_+^{\ind}(p^+_X) \,{\otimes_{C_l}}\, G_+^{\ind}(p^+_Y)]
  \\{}\\[-.7em]
  &= \Ygamma++l_{G_+(X)\,G_+(Y)} \,, \eear\ee
where we also used the compatibility of $\Xb{}{}$ with $G_+^{\ind}$ from 
\Lemma \ref{lem:G-ind-tens}(iii).  \qed

We are now in a position to present our first main result, the ribbon
equivalences between local $C_{l/r}(A)$-modules and ambichiral
$A$-bimodules; based on results of \cite{boek2}, these equivalences
have been conjectured in `claim 5' of \cite{ostr}.

\dtl{Theorem}{thm:equiv}
Let $A$ be a \ssFA\ in a ribbon category $\cC$ such that the symmetric 
Frobenius algebras $C_{l/r}(A)$ are special as well. Then there are equivalences
  \be
  \cC_{C_l(A)}^{\sss\rmloc} \,\cong\,
  \CAAo \,\cong\, \cC_{C_r(A)}^{\sss\rmloc}  \ee
of ribbon categories.

\medskip

We will only present the proof of the equivalence 
$\cC_{C_l(A)}^{\sss\rmloc}\,{\cong}\,\CAAo$ explicitly; the second equivalence 
can be shown by similar means.%
  \foodnode{Recall also \convention s \ref{basicprops} and \ref{c:csplit}.} 
As a preparation we need the following two lemmata.
\\[-2.3em]

\dtl{Lemma}{le:equiv-prep-3}
We have the following bijections
between spaces of bimodule morphisms of $\alpha$-induced $A$-bimodules and 
module morphisms of (locally) induced $C_l$-modules:
  \bea
  \Psi^{l\,+-}_{\!A;UV} :\quad \HomAA( \alpha_\AA^+(U), \alpha_\AA^-(V) ) 
    \,\overset{\cong}{\longrightarrow}\,
    \Hom_{C_l}(\Ind_{C_l}(U),\lxtp V\AA{l}) \\{}\\[-.7em]
  \Psi^{l\,-+}_{\!A;UV} :\quad \HomAA( \alpha_\AA^-(U), \alpha_\AA^+(V) ) 
    \,\overset{\cong}{\longrightarrow}\,
    \Hom_{C_l}(\lxtp U\AA{l},\Ind_{C_l}(V)) \,.
  \eear\labl{eq:equiv-prep-3a}
The maps $\Psi^{l\,+-}_{\!A;UV}$ and $\Psi^{l\,-+}_{\!A;UV}$ are given by
  \be
  \Psi^{l\,+-}_{\!A;UV}(f) = f \cir (e_{C_l} \otimes \id_U)
  \qquad {\rm and} \qquad
  \Psi^{l\,-+}_{\!A;UV}(g) = (r_{C_l} \otimes \id_V) \cir g \,.
  \labl{eq:equiv-prep-3b}
In the definition of $\Psi^{l\,+-}_{A;UV}$, the realisation of $\lxtp V\AA{l}$ 
as $(\Ind_\AA(V);P^l_\AA(V))$ is implied for obtaining the relevant subspace of 
$\Hom_{C_l}(\Ind_\AA(U),\Ind_\AA(V))$, and similar implications hold for the
definition of $\Psi^{l\,-+}_{A;UV}$.
The bijections \erf{eq:equiv-prep-3a} satisfy
  \be
  \Psi^{l\,-+}_{\!A;VU}(g)\cir\Psi^{l\,+-}_{\!A;UV}(f)
  = \Phi^{l\,++}_{\!A;UU}(g \cir f) \,.  \labl{eq:equiv-prep-3c}

\noindent
Proof:\\
This is a consequence of \Proposition \ref{alphaalpha} 
together with the reciprocity relations (see 
\Remark \ref{prop:ssFA-unique}(iii))
  \bea
  \Hom(U,\efu Vl\AA)\cong\Hom_{C_l}(\Ind_{C_l}(U),\lxtp V\AA{l})
  \qquad{\rm and} \\{}\\[-.7em]
  \Hom(\efu Ul\AA,V)\cong\Hom_{C_l}(\lxtp U\AA{l},\Ind_{C_l}(V)) \,.  \eear\ee
Using the explicit form \erf{eq:reci-inv} and \erf{AUAV} of these bijections,
one checks that $\Psi^{l\,+-}_{\!A;UV}$ and $\Psi^{l\,-+}_{\!A;UV}$ are 
indeed given by the maps \erf{eq:equiv-prep-3b}. Furthermore, substituting 
\erf{eq:equiv-prep-3b} and the definition \erf{eq:phi-map-def}, 
it is immediate that \erf{eq:equiv-prep-3c} holds true.
\qed

\dtl{Lemma}{le:equiv-prep-2}
The following two statements are equivalent:
\\[.2em]
(i)~\,$((\alpha_\AA^+(U);p),e,r)$ is a $A$-bimodule retract of 
  $\alpha_\AA^-(V)$,
\\[.2em]
(ii)~$(  (\Ind_{C_l}(U); \Phi^l_{\!A;UU}(p)) , \Psi^{l\,+-}_{\!A;UV}(e),
  \Psi^{l\,-+}_{\!A;VU}(r) )$ is a $C_l$-module retract of $\lxtp V\AA{l}$.
  
\medskip\noindent
Proof:\\
We will need two series of identities, both of which hold for any choice of
morphisms $p \iN \HomAA(\alpha_\AA^+(U),\alpha_\AA^+(U))$, 
$e\iN\HomAA(\alpha_\AA^+(U),\alpha_\AA^-(V))$ and
$r\iN\HomAA(\alpha_\AA^-(V),\alpha_\AA^+(U))$. 
\\
The first series of identities is
  \be \bearll
  \Phi^l_{\!A;UU}(p) \cir \Phi^l_{\!A;UU}(p) \!\!
  &= \Phi^l_{\!A;UU}(p \cir (P^l_\AA \cir \id_U) \cir p)
  \\{}\\[-.7em]
  &= \Phi^l_{\!A;UU}( (P^l_\AA \cir \id_U) \cir p \cir p)
  = \Phi^l_{\!A;UU}(p \cir p) \,.  \eear \labl{eq:equiv-prep-2a}
Here the first step holds by \Lemma \ref{lem:phi-prop}, 
in the second step \Lemma \ref{++lemma} 
is used, and in the third step the idempotent
$P^l_\AA$ is omitted against the restriction morphism $r_{C_l}$ contained, 
by definition, in $\Phi^l_{\!A;UU}$. The second series of identities is
  \be 
  \Psi^{l\,-+}_{\!A;VU}(r) \cir \Psi^{l\,+-}_{\!A;UV}(e) 
  = \Phi^{l\,++}_{\!A;UU}(r \cir e) = \Phi^l_{\!A;UU}(r \cir e) \,,
  \labl{eq:equiv-prep-2b}
where the first equality uses \erf{eq:equiv-prep-3c} and in the second 
equality holds because $\Phi^{l\,++}_{\!A;UU}$ is just a restriction 
of $\Phi^l_{\!A;UU}$ to a subspace.
\\[3pt]
(i)\,$\Rightarrow$\,(ii):
By assumption (i), $p$ is an idempotent and we have $r\cir e\eq p$. By 
\erf{eq:equiv-prep-2a} this implies that $\Phi^l_{\!A;UU}(p)$ is an 
idempotent, too. Furthermore, by the equality \erf{eq:equiv-prep-2b} we have
$\Psi^{l\,-+}_{\!A;VU}(r) \cir \Psi^{l\,+-}_{\!A;UV}(e)
\eq \Phi^l_{\!A;UU}(p)$, which is equal to the identity morphism of
$(\Ind_{C_l}(U); \Phi^l_{\!A;UU}
   $\linebreak[0]$
(p))$, thus establishing that
we are indeed dealing with a $C_l$-module retract.
\\[3pt]
(ii)\,$\Rightarrow$\,(i): Conversely, suppose that
$\Psi^{l\,-+}_{\!A;VU}(r) \cir \Psi^{l\,+-}_{\!A;UV}(e)\eq \Phi^l_{\!A;UU}(p)$ 
and that $\Phi^l_{\!A;UU}(p)$ is an idempotent. Then equations 
\erf{eq:equiv-prep-2a} and \erf{eq:equiv-prep-2b} tell us
that also $\Phi^l_{\!A;UU}(p \cir p)\eq \Phi^l_{\!A;UU}(p)$ and
$\Phi^l_{\!A;UU}(r \cir e)\eq \Phi^l_{\!A;UU}(p)$.
Since, by the first isomorphism in \Lemma \ref{le:equiv-prep-1},
$\Phi^l_{\!A;UU}$ is injective on $\End_{A|A}(\alpha_\AA^+(U))$, it follows
that $p$ is an idempotent, and that $e \cir r\eq p$, 
which is the identity morphism in $\End_{A|A}( (\alpha_\AA^+(U);p) )$.
\qed

\medskip\noindent
Proof of \Theorem {\bf\ref{thm:equiv}}: \\
Denote by $G{:}\;\CAAo\,{\to}\,\cC_{C_l}$ the restriction of $G_+$ to $\CAAo$. 
We will show that $G$ is a ribbon equivalence between $\CAAo$ and
$\cC_{C_l(A)}^{\sss\rmloc}$.
\\[3pt]
(i)~{\em The image of $G$ consists of local modules\/}:
Objects in \CAAp\ are of the form $B\eq(\alpha_\AA^+(U);p)$.
If $B$ is also in \CAAm, then there exist morphisms $e,r$ such that
$((\alpha_\AA^+(U);p),e,r)$ is a bimodule retract of $\alpha_\AA^-(V)$ 
for some $V{\in}\,\Objc$. By \Lemma \ref{le:equiv-prep-2} 
it follows that $( (\Ind_{C_l}(U); \Phi^l_{\!A;UU}(p)) , 
    $\linebreak[0]$%
\Psi^{l\,+-}_{\!A;UV}(e), \Psi^{l\,-+}_{\!A;VU}(r) )$ is a $C_l$-module
retract of the local module $\lxtp V\AA{l}$.
\\
Thus $G(B)\eq(\Ind_{C_l}(U); \Phi^l_{\!A;UU}(p))$ is a
retract of a local module, and hence local itself.
\\[3pt]
(ii)~\,{\em $G$ is essentially surjective on the category of local modules\/}:
{}From (i) we know that $G$ is a functor from $\CAAo$ to $\cC_{C_l(A)}
^{\sss\rmloc}$. By \Proposition \ref{sub-of-locind} 
every local module $M$ is
isomorphic to a retract of a locally induced module $\lxtp V\AA{l}$ for some 
$V\iN\Objc$. We can write $M \,{\cong}\, (\lxtp V\AA{l}; q)$ for some 
idempotent $q \iN \End_{C_l}(\lxtp V\AA{l})$. 
However, we want to make a statement involving $C_l$-modules rather than
locally induced $A$-modules. To this end we introduce the morphisms
  \bea
  e' := q \cir r_{\!\!A\otimes V\succ\Lxtp V{l}\AA} 
  \cir (m \oti \id_V) \cir (e_{C_l} \oti \id_A \oti \id_V)
  \;\in \Hom(C_l\Oti A\Oti V,\lxtp V\AA l)
  \quad\ {\rm and} \\{}\\[-.7em]
  r' := (r_{C_l} \oti \id_A \oti \id_V ) \cir( \Delta \oti \id_V )\cir 
  e_{\Lxtp V{l}\AA\prec A\otimes V} \cir q
  \,\in \Hom(\lxtp V\AA l,C_l\Oti A\Oti V) \,.  \eear\ee 
\vskip.2em\noindent
These morphisms fulfill $e' \cir r'\eq q$, as can be seen as follows. First note 
that \Lemma \ref{le4sublocind}, specialised to $U\eq V\eq\one$ and 
$\Phi\eq\id_A$, together with \Lemma \ref{lem:C=[1]a}(iii) 
and \ref{lem:C=[1]a}(ii) 
as well as specialness of $C_l$, implies that
$m \cir (\id_A \cir P^r_\AA) \cir \Delta\eq \id_A$. It is then easy
to convince oneself that an appropriately modified version of 
\Lemma \ref{le4sublocind} gives rise to the analogous identity
$m \cir (P^l_\AA \cir \id_A) \cir \Delta\eq \id_A$. This, in turn, implies
$e' \cir r'\eq q$.
\\
Next define $p'\,{:=}\, r' \cir e'$. Because of $q \cir e'\eq e'$,
$p'$ is an idempotent. Thus by construction we have an isomorphism
  \be
  ( (\Ind_{C_l}(A\oti V);p') , e', r' ) \cong (\lxtp V\AA{l};q) \ee
of $C_l$-modules. Thus $( (\Ind_{C_l}(U);p'),e',r' )$ with $U\,{:=}\,A\oti V$
is a module retract of $\lxtp V\AA{l}$.
\\
By the \Lemmata \ref{le:equiv-prep-1} and \ref{le:equiv-prep-3}  
we can now find morphisms $p\iN\End_{\!A|A}(\alpha_\AA^+(U))$, 
$e \iN \HomAA
    $\linebreak[0]$%
(\alpha_\AA^+(U),\alpha_\AA^-(U) )$ and
$r \iN \HomAA( \alpha_\AA^-(U), \alpha_\AA^+(U) )$ such that
$\Phi^{l\,++}_{\!A;UU}(p)\eq p'$, $\Psi^{l\,+-}_{\!A;UU}(e)\eq e'$
and $\Psi^{l\,-+}_{\!A;UU}(r)\eq r'$. Then we can use lemma 
\ref{le:equiv-prep-2} to conclude that $((\alpha_\AA^+(U);p),e,r)$ is an 
$A$-bimodule retract of $\alpha_\AA^-(V)$. Thus we have found an object 
$B\eq (\alpha_\AA^+(U);p)$ in $\CAAo$ such that $G(B)\,{\cong}\,M$.
\\[3pt]
(iii)~{\em $G$ is an equivalence of ribbon categories\/}:
Note that $G{:}\;\CAAo\,{\to}\,\cC_{C_l}^{\sss\rmloc}$ is an equivalence 
functor because first, it is essentially surjective on objects, and second, it 
is a restriction of $G_+$, which is bijective on morphisms. Since $G_+$ is 
a tensor functor, so is $G$. Furthermore, for the family $\Xbeta{}{}_{XY}$ 
of morphisms we have
  \be
  G(\Xbeta{}{}_{\!XY}) = G_+(\Xbeta++_{XY}) = \Ygamma++l_{\!G_+(X)\,G_+(Y)}
   = c^{C_l}_{G(X)\,G(Y)} \,, \ee
where we first used \Proposition \ref{pr:beta-braid}~(ii), 
then \Proposition \ref{pr:G-tens} and finally 
\Proposition \ref{pr:cmn=cloc}. 
Since $\Xbeta{}{}_{XY}$ is mapped to the braiding $c^{C_l}$ on $\cC_{C_l}^
{\sss\rmloc}$ by an equivalence functor, it follows that $\Xbeta{}{}_{XY}$ 
defines a braiding on \CAAo. -- This completes the proof of 
\Proposition \ref{pr:beta-braid} by also establishing part (iii) of the 
proposition.
\\
Hence the tensor equivalence $G$ is compatible with the braiding. Thus $G$ is an 
equivalence of braided tensor categories, and thereby also of ribbon categories.
\qed

\dt{Remark}
Denote by $G_{lr}{:}\;\cC_{C_l(A)}^{\sss\rmloc}\,{\to}\,\cC_{C_r(A)}
^{\sss\rmloc}$ and $G_{rl}{:}\;\cC_{C_r(A)}^{\sss\rmloc}\,{\to}\,
\cC_{C_l(A)}^{\sss\rmloc}$ the functorial equivalences
of the ribbon categories $\cC_{C_l(A)}^{\sss\rmloc}$
and $\cC_{C_r(A)}^{\sss\rmloc}$ constructed in \Theorem \ref{thm:equiv}. 
\\
One can give an explicit representation of $G_{l/r}$ using retracts. Consider 
the two morphisms $Q_{lr}(M_l) \iN \Hom(A{\otimes}\M_l,A{\otimes}\M_l)$ and
$Q_{rl}(M_r) \iN \Hom(A{\otimes}\M_r,A{\otimes}\M_r)$ given by
  %% [pic~19]
  \bea  \begin{picture}(190,85)(0,34)
  \put(0,0)  {\begin{picture}(0,0)(0,0)
              \scalebox{.38}{\includegraphics{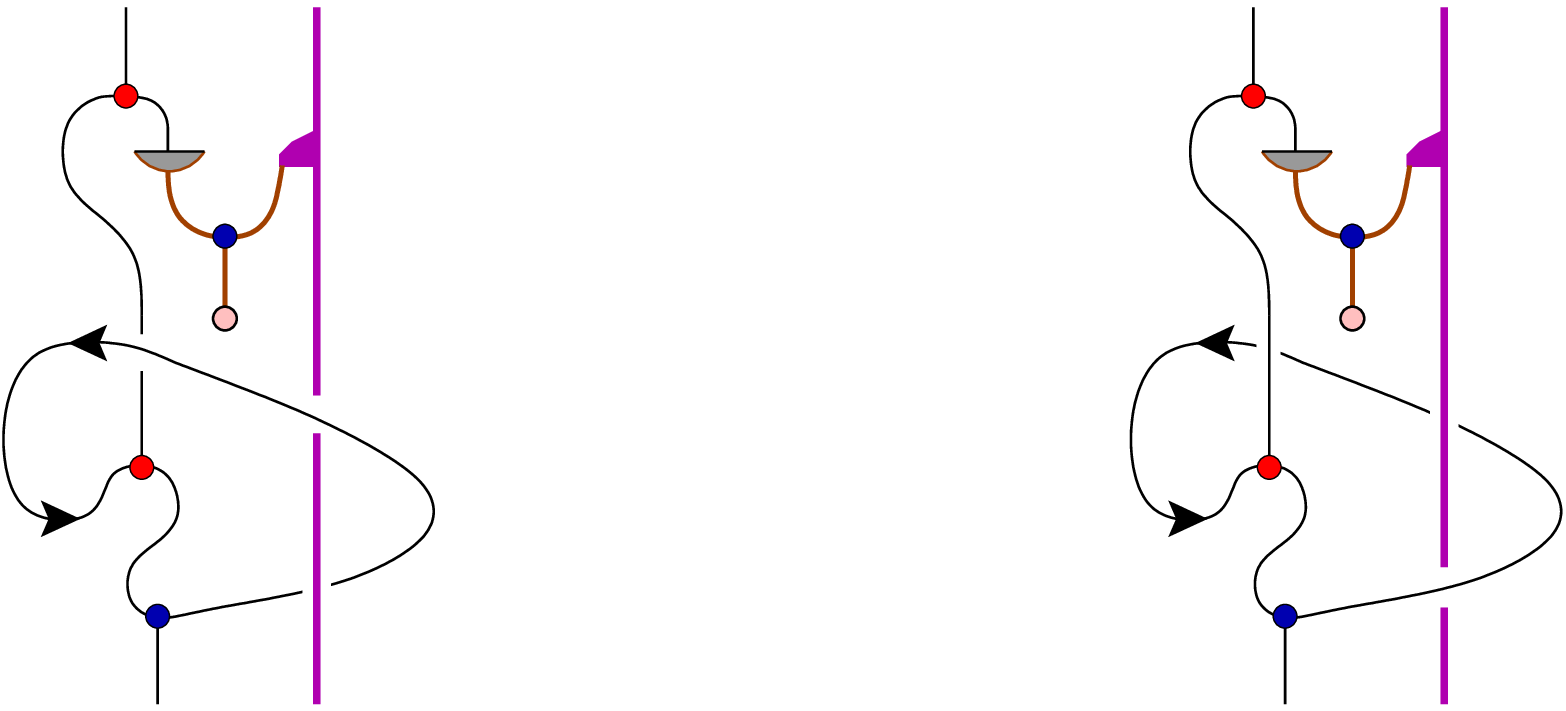}} \end{picture}}
  \put(-68,48.2)    {$Q_{lr}(M_l)\,:=$}
  \put(16.1,110.5)  {\sse$A$}
  \put(19.7,-9.2)   {\sse$A$}
  \put(35.4,64.4)   {\sse$C_l$}
  \put(42.6,-9.2)   {\sse$\M_l$}
  \put(42.6,110.9)  {\sse$\M_l$}
  \put(103,48.2)    {$Q_{rl}(M_r)\,:=$}
  \put(187.1,110.5) {\sse$A$}
  \put(190.7,-9.2)  {\sse$A$}
  \put(206.4,64.4)  {\sse$C_r$}
  \put(213.6,-9.2)  {\sse$\M_r$}
  \put(213.6,110.9) {\sse$\M_r$}
  \epicture23 \labl{eq:pic19}
By combining several previous results one sees that $Q_{lr/rl}(M_{l/r})$ are
idempotents: the morphisms $P_\AA^{l/r}(M_{r/l})$ from \erf{PU-def} are 
idempotents, \erf{eq:remove-braiding} can be used to commute the ribbon 
connecting $A$ to $M_{r/l}$ past the $A$-loop, and finally one can use 
\erf{eq:remove-P} together with specialness of $C_{l/r}$, 
which holds by the assumptions in \Theorem \ref{thm:equiv}.
\\
For local $C_{l/r}$-modules $M_{l/r}$ one has
  \be
  G_{lr}(M_l)= \Im Q_{lr}(M_l) \qquad {\rm and} \qquad
  G_{rl}(M_r)= \Im Q_{rl}(M_r)  \ee
(recall that we work with Karoubian categories, so that all idempotents are
split), and the action of the functors $G_{lr}$ and $G_{rl}$ on morphisms reads
  %% [pic~20]
  \bea  \begin{picture}(130,80)(0,24)
  \put(0,0)  {\begin{picture}(0,0)(0,0)
              \scalebox{.38}{\includegraphics{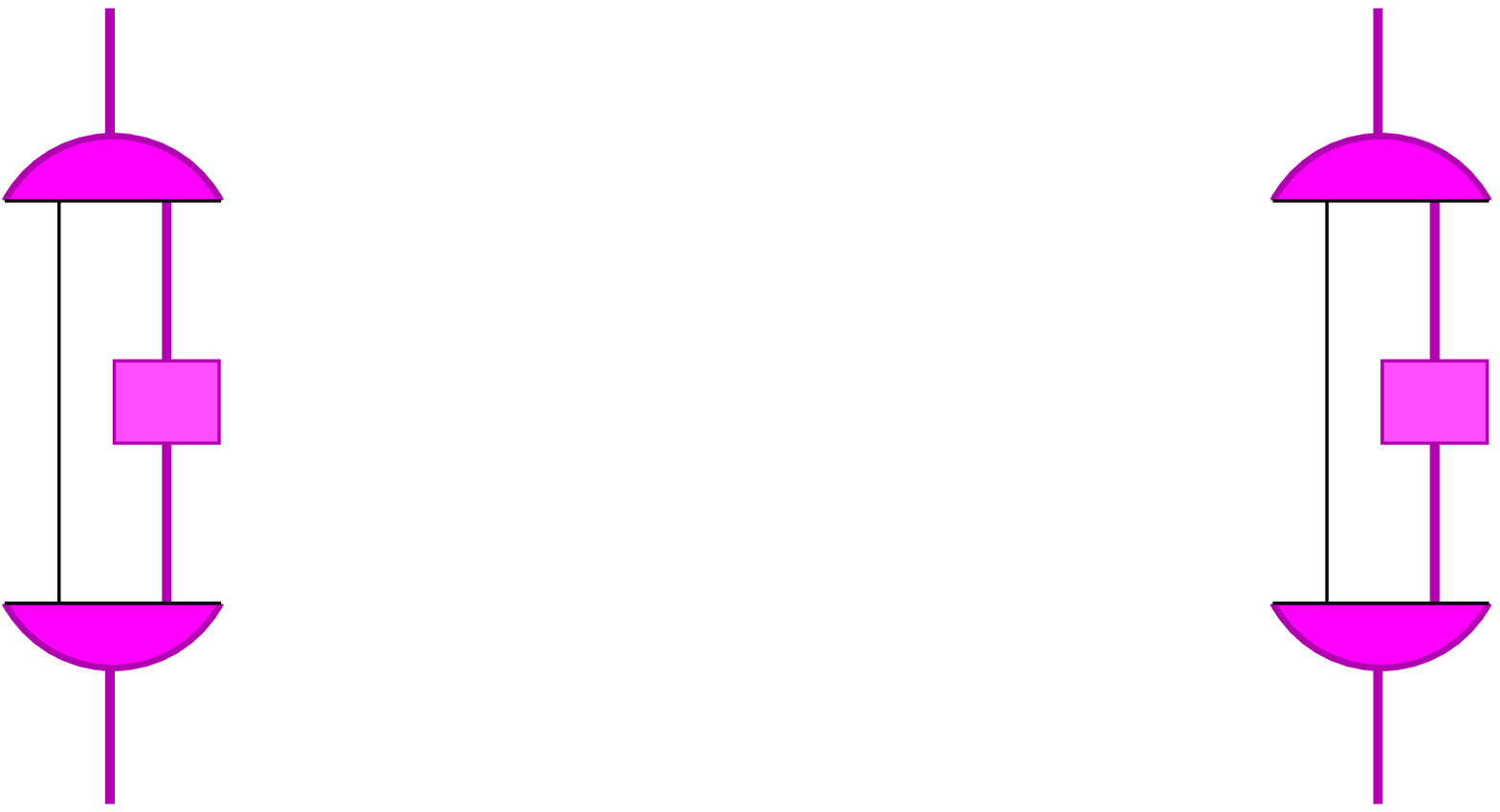}} \end{picture}}
  \put(-61,45.2)    {$G_{lr}(f_l)\;:=$}
  \put(-.9,33.5)    {\sse$A$}
  \put(1.1,-9.2)    {\sse$G_{lr}(M_l)$}  
  \put(2.1,99.5)    {\sse$G_{lr}(N_l)$}
  \put(11.1,19.5)   {\sse$e$}
  \put(10.5,73.3)   {\sse$r$}
  \put(15.9,46.1)   {\sse$f_l$}
  \put(20.8,31.9)   {\sse$M_l$}
  \put(20.8,59.9)   {\sse$N_l$}
  \put(87,45.2)     {$G_{rl}(f_r)\;:=$}
  \put(147.7,33.5)  {\sse$A$}
  \put(148.6,-9.2)  {\sse$G_{rl}(M_r)$}
  \put(149.6,99.5)  {\sse$G_{rl}(N_r)$}
  \put(160.9,19.5)  {\sse$e$}
  \put(160.3,73.3)  {\sse$r$}
  \put(164.1,46.1)  {\sse$f_r$}
  \put(169.1,31.9)  {\sse$M_r$}
  \put(169.1,59.9)  {\sse$N_r$}
  \epicture12 \labl{eq:pic20}
for $f_r\iN\Hom_{C_r}(M_r,N_r)$ and $f_l\iN\Hom_{C_l}(M_l,N_l)$.

\dt{Remark}  
(i)~\,The equivalence of the categories of local modules over the left and right
centers given in \Theorem \ref{thm:equiv} 
is a category theoretic analogue of \Theorem 5.5 of \cite{boev4}, 
which was obtained in the study of relations between 
nets of braided subfactors and modular invariants. In the context of module 
categories, the equivalence, including the relation to the category of ambichiral
bimodules, has been formulated, as a conjecture, in \Section 5.4 of \cite{ostr}.
\\[.3em]
(ii)~It is known \cite{mose2} that in conformal quantum field theory, every 
modular invariant torus partition function can be described in terms of 
extensions of the chiral algebras for left movers and right movers. The two 
extensions need not be the same, but they should lead to extended theories 
with isomorphic fusion rules. The additional information in a modular 
invariant partition function is the choice of an isomorphism of these fusion 
rules. This structure is sometimes summarised by saying that the torus partition 
function of every full conformal field theory has the form of `a fusion rule 
isomorphism on top of (maximal) extensions of the chiral algebras'.
\\
This statement has been obtained in \cite{mose2} using the action of the
(cover of the) modular group $SL(2,\zet)$ on the characters of a chiral CFT, 
the invariance of the torus partition function under this action, and the
non-negativity of its coefficients.
\\[.2em]
The connection between this description of partition functions and our study 
of algebras in tensor categories is supplied by the insight \cite{fuRs,fuRs4} 
that, given a chiral rational conformal field theory, a full rational CFT,
including in particular its torus and annulus partition functions,
can be constructed from 
a symmetric special Frobenius algebra $A$ in the modular tensor category $\cC$ 
that describes the chiral data of the CFT. (But not every modular invariant 
bilinear combination of characters of the chiral CFT is the torus partition 
functions of some full CFT.) The structure of partition functions described 
above can be obtained from \Theorem \ref{thm:equiv} as follows. The procedure 
of `extending the chiral algebra for left movers and right movers' corresponds 
to passing to the modular tensor categories $\Ext{\cC}{C_{l/r}(A)}$ of local 
modules of the left center and the right center, respectively, of $A$. By 
\Theorem \ref{thm:equiv} these two categories are equivalent, so that in 
particular they have isomorphic fusion rules,
  \be  K_0\!\left(\Ext{\cC}{C_{l}(A)}\right) \cong
  K_0\!\left(\Ext{\cC}{C_{r}(A)}\right) .  \ee 
We may lift the algebra $A$ to algebras in $\Ext{\cC}{C_{l/r}(A)}$ via lemma 
\ref{lem:alglift} to obtain algebras with trivial center. In this sense, the 
two extensions are `maximal' and the isomorphism of the fusion rules
is encoded in the `non-commutative part' of the algebra $A$.

%%%%%%%%%%%%%%%%%%%%%%%%%%%%%%%%%%%%%%%%%%%%%%%%%%%%%%%%%%%%%%%%%%%%%%%%

\sect{Product \cats\ and \platy}\label{platy}

In many respects the simplest tensor \cats\ are the \cats\ of \findim\ vector
spaces over some field $\koerper$; we denote
the latter \cat\ by $\Vectk$. It is therefore interesting to find 
commutative \ssFA s $A$ in ribbon \cats\ which are `trivialising' in the 
sense that the \cat\ \calcal\ of local $A$-modules is equivalent to $\Vectk$. 
For a generic ribbon \cat\ $\cC$ such a trivialising \alg\ need not exist.
A class of \cats\ for which a trivialising \alg\ does exist is provided by
the \rep\ \cats\ for so-called holomorphic orbifolds \cite{dipr,bant4}: for 
these, the trivialising \alg\ affords the extension of the corresponding 
orbifold conformal field theory to the underlying un-orbifolded theory.

We may, however, relax the requirement and instead look, for given $\cC$,
for some `compensating' ribbon \cat\ $\cC'$ and a trivialising \alg\ $T$ in 
the (suitably defined) product of $\cC$ with $\cC'$ -- for the precise 
formulation of this concept of \platy, see \Definition \ref{def:platy} 
below.
The main purpose of this section is to establish that such a \cat\ $\cC'$ and 
\alg\ $T$ always exist when $\cC$ is a {\em modular\/} \tc. In that case, for 
$\cC'$ we can take the \cat\ dual to $\cC$, a concept that will be discussed in 
\Section \ref{platy2}.

%%%%%%%%%%%%%%%%%%%%%%%%%%%%%%%%%%%%%%%%%%%

\subsection{Product categories and the notion of \platy}\label{platy1}
     
But first we must introduce a suitable concept of product which to any pair of
$\koerper$-linear categories $\cC$ and $\cD$ associates a product
category that shares with $\cC$ and $\cD$ all the relevant properties, such as
the basic properties listed in the \convention\ 
\ref{basicprops}. This is done in 
\\[-2.3em]

\dtl{Definition}{def:Ti}
Let $\cC$ and $\cD$ be $\koerper$-linear categories.
\\[.3em]
(i)~\,The category $\cC\Tic\cD$ is the category whose objects are pairs 
$U\ti X$ with $U\iN\Obj(\cC)$ and $X\iN\Obj(\cD)$ and whose morphism spaces 
are tensor products (over $\koerper$)
  \be  \Hom^{\cC\Tic\cD}(U{\times}X,V{\times}Y)
  := \Hom^\cC(U,V)\,{\otimes^{}_\koerper}\,\Hom^\cD(X,Y) \labl{ctimes}
of those of $\cC$ and $\cD$.
\\[.3em]
(ii)~The {\em Karoubian product\/} $\cC\Ti\cD$ is the Karoubian envelope
of $\cC\Tic\cD$,
  \be  \cC\Ti\cD := \kar{(\cC\Tic\cD)}  \,.  \labl{boxtimes}

\dtl{Remark}{rem1box}
(i)~~Taking the tensor product over $\koerper$ rather than the
Kronecker product of the morphism sets accounts for the fact that the
\cats\ of our interest are enriched over $\Vectk$. The price to
pay is that $\cC\Tic\cD$ has idempotents that are not tensor products
of idempotents in $\cC$ and $\cD$, so that even when $\cC$ and $\cD$ 
are Karoubian we get, in general, a Karoubian product 
\cat\ only after taking the Karoubian envelope.
\\[.3em]
(ii)~\,In accordance with \Remark \ref{K-rem1}(iii) we regard the category
$\cC\Tic\cD$ as a full sub\cat\ of $\cC\Ti\cD$, i.e.\ in particular identify
$U{\times}X\iN\Obj(\cC\Tic\cD)$ 
with $(U{\times}X,\id_U\Otic\id_X)\iN\Obj(\cC\Ti\cD)$.
\\[.3em]
(iii)~When $\cC$ and $\cD$ are small \cats, then so are $\cC\Tic\cD$ and
$\cC\Ti\cD$. When $\cC$ and $\cD$ are additive, then so is $\cC\Ti\cD$.
When $\cC$ and $\cD$ are semisimple, then so is $\cC\Ti\cD$.
\\[.3em]
(iv)~When $\cC$ and $\cD$ are modular tensor categories, then so is their
Karoubian product $\cC\Ti\cD$, see \Proposition \ref{prop:Ti-mod}(iii) 
below. It is easy to verify that the dimension and charge of modular tensor 
categories, as defined in \erf{eq:dim-kappa-def}, are multiplicative, i.e.
  \be
  {\rm Dim}(\cC\Ti\cD) = {\rm Dim}(\cC) \, {\rm Dim}(\cD)
  \qquad {\rm and} \qquad
  p^\pm(\cC\Ti\cD) = p^\pm(\cC) \, p^\pm(\cD) \,.  \ee

\dtl{Proposition}{prop:Ti-mod}
(i)~~When $\cC$ and $\cD$ are tensor \cats, then $\cC\Tic\cD$ can be
naturally equipped with the structure of a \tc, by setting
  \be  \bearl
  (U{\times}X) \oti^{\!\cC\Tic\cD} (V{\times}Y)
  := (U{\otimes^\cC}V) \ti (X{\otimes^\cD}Y) \,,
  \qquad \one^{\cC\Tic\cD} := \one^\cC{\times}\one^{\!\cD} \qquad{\rm and}
  \\{}\\[-.6em]
  (f\Otic\, g) \oti^{\!\cC\Tic\cD} (f'\Otic\, g')
  := (f{\otimes^\cC}f') \otic (g\,{\otimes^\cD}g')
  \,. \eear \labl{Ti-tensor}
%\\[.3em]
(ii)~\,Similarly, $\cC\Tic\cD$ inherits from $\cC$ and $\cD$ the properties
of having a (left or right) duality, a braiding, and a twist, by setting
  \be  \bearl  d_{U\times X}^{\cC\Tic\cD} := d_U^\cC \otic d_V^\cD
  \qquad {\rm etc.}\,,
  \\{}\\[-.6em]
  c_{U\times X,V\times Y}^{\cC\Tic\cD} := c_{U,V}^\cC \otic c_{X,Y}^\cD\,,
  \\{}\\[-.6em]
  \theta_{U\times X}^{\cC\Tic\cD} := \theta_U^\cC \otic \theta_V^\cD
  \,. \eear \labl{Ti-ribbon}
In particular, when $\cC$ and $\cD$ are ribbon \cats, then $\cC\Tic\cD$
is naturally equipped with the structure of a ribbon \cat. Moreover,
  \be  s_{U\times X,V\times Y}^{\cC\Tic\cD} = s_{U,V}^\cC \, s_{X,Y}^\cD
  \,;  \labl{eq:boxtimes-smat}
in particular, the dimensions in $\cC\Tic\cD$ are given by
  \be  \dim^{\cC\Tic\cD\!}_{}(U{\times}X)
  = \dim^{\cC\!}_{}(U)\, \dim^{\!\cD\!}_{}(X) \,. \ee
\\[.3em]
(iii)~Analogous statements as in (i) and (ii) apply to the Karoubian product
$\cC\Ti\cD$. In addition, if $\cC$ and $\cD$ are modular tensor \cats,
then the category $\cC\Ti\cD$ has a natural structure of modular \tc.

\bigskip\noindent
Proof:\\
(i), (ii)~\,Using the relevant properties of $\cC$ and $\cD$, it is 
straightforward to check that with the definitions \erf{Ti-tensor} and 
\erf{Ti-ribbon}, all required relations for morphisms in $\cC\Tic\cD$
are satisfied.
\\[.3em]
(iii)~then holds by combining these results with the properties of the 
Karoubian envelope listed in \Remark \ref{K-rem1}(iv). 
For modular $\cC$ and $\cD$, $\cC\Ti\cD$ is additive and semisimple by 
\Remark \ref{rem1box}(iii),
and the $s$-matrix \erf{eq:boxtimes-smat} is non-degenerate because those 
of $\cC$ and $\cD$ are. Thus $\cC\Ti\cD$ is indeed modular.
\qed

\medskip

We are now in a position to introduce the concept of {\em \platy\/} of $\cC$:
\\[-2.3em]

\dtl{Definition}{def:platy}
A ribbon \cat\ $\cC$ is called {\em \platl\/} iff there exist a 
ribbon \cat\ $\cC'$ and a commutative \ssFA\ $T$ in 
$\cC\Ti\cC'$ such that the \cat\ of local $T$-modules is equivalent
to the \cat\ of \findim\ vector spaces over $\koerper$,
  \be  \Ext{(\cC\Ti\cC')}T \,\cong\, \Vectk \,.  \ee
The data $\cC'$ and $T$ are then called a {\em \platz\/} of $\cC$.

\medskip

The rest of this subsection is devoted to the study of the Karoubian
product of tensor categories and its behaviour in the context of module 
categories.

\dtl{Lemma}{boxtimes2}
The Karoubian product of two \cats\ is equivalent to the
Karoubian product of their Karoubian envelopes,
  \be  \cC\,\Ti\,\cD  \,\cong\, \kar\cC\,\Ti\,\kar\cD .  \labl{eq:boxtimes2}
If $\cC$ and $\cD$ are ribbon, then this is an equivalence of ribbon categories.

\medskip\noindent
Proof:\\
According to \Proposition \ref{XI.1.5} to show the equivalence it is 
sufficient to construct a functor $F{:}\;\kar\cC\Ti\kar\cD\,{\to}\,\cC\Ti\cD$
that is essentially surjective on objects and bijective on morphisms.
\\[.2em]
The objects of $\cC\,\Ti\,\cD\,{\equiv}\,\kar{(\cC\Tic\cD)}$ are triples
$(U\ti X;\pi)$ with $\pi$ an idempotent in $\End(U{\times}X)\,
    $\linebreak[0]$%
{\cong}\,\End(U)\otic\End(X)$, while
the objects of $\kar\cC\,\Ti\,\kar\cD$ are quintuples
$((U;p)\ti(X;q);\hat\pi)$, where $U\iN\Objc$, $X\iN\Obj(\cD)$, 
$p\iN\End(U)$ and $q\iN\End(X)$ are idempotents in $\cC$ and $\cD$,
respectively, and $\hat\pi\iN\End(U{\times}X)$
is an idempotent obeying the Karoubi condition
  \be  (p\otic q)\circ \hat\pi = \hat\pi = \hat\pi \circ (p\otic q)
  \,.  \labl{pqpi}
We define the functor $F$ on objects as
  \be  F\Llb ((U;p)\ti(X;q);\pi) \Lrb := (U\ti X;\pi) \,. \labl{fff}
It then follows that we get every object $(U{\times}X;\pi)$ of $\cC\Ti\cD$
as the image under $F$ of the object $((U;\id_U)\ti(X;\id_X);\pi)$.
Hence $F$ is surjective on objects.
\\[.2em]
To define $F$ on morphisms, we first introduce, for any
two objects $((U;p){\times}(X;q);\pi)$ and $((V;p'){\times}(Y;q');\pi')$
of $\kar\cC\,\Ti\,\kar\cD$, certain endomorphisms $P$, $Q$ and $\Pi$ of 
vector spaces:
  \be \!
  \begin{array}{rcl} P:\;\  \Hom^{\cC}(U,V) &\!\!\to\!\!& \Hom^{\cC}(U,V) \\
                            f &\!\!\mapsto\!\!& p'\cir f\cir p  \end{array}
  \quad\,{\rm and}\quad\,
  \begin{array}{rcl} Q:\;\  \Hom^{\cD}(X,Y) &\!\!\to\!\!& \Hom^{\cD}(X,Y)\\
                            g &\!\!\mapsto\!\!& q'\cir g\cir q  \end{array}
  \ee
as well as
  \be
  \begin{array}{rcl} \Pi:\quad \Hom^{\cC\Tic\cD}(U{\times}X,V{\times}Y)
                     &\!\to\!& \Hom^{\cC\Tic\cD}(U{\times}X,V{\times}Y)\\
                     \psi &\!\mapsto\!& \pi'\cir\psi\cir\pi \,;  \end{array} 
  \ee
$P$, $\Pi$ and $Q$ are idempotents of vector spaces.
One checks that, by definition of the Ka\-rou\-bian envelope, 
  \be  
  \Hom^{\kar\cC\Ti\kar\cD\!}(((U;p)\ti(X;q);\pi),((V;p')\ti(Y;q');\pi'))
  \cong\, 
  {\rm Im}(P)\otic{\rm Im}(Q) \,\cap\, {\rm Im}(\Pi) \,,  
  \labl{ImPQPi}
while
  \be  \Hom^{\cC\Ti\cD}((U{\times}X;\pi),(V{\times}Y;\pi'))
  \,\cong\,
  {\rm Im}(\Pi) \,.
  \labl{ImPi}
In addition, from \erf{pqpi} it follows that 
$(P\otic Q)\cir\Pi\eq\Pi\eq\Pi\cir(P\otic Q)$, which in turn implies that
  \be  {\rm Im}(\Pi) \,\subseteq\,
  {\rm Im}(P\Otic Q) = {\rm Im}(P)\otic{\rm Im}(Q) \,.  \ee
We can thus conclude that the morphism spaces \erf{ImPQPi} and \erf{ImPi} 
are actually identical subspaces of $\Hom^{\cC\Tic\cD}(U{\times}X,V{\times}Y)
\eq \Hom^\cC(U,V)\,{\Tic}\,\Hom^\cD(X,Y)$. 
\\ 
We now simply define $F$ to be the identity map on morphisms, so that $F$ is 
in particular bijective on morphisms. It is easy to check that together with 
\erf{fff} this yields a functor from $\kar\cC\Ti\kar\cD$ to $\cC\Ti\cD$.
\\[.2em]
Thus $F$ is an equivalence functor from $\kar\cC\,\Ti\,\kar\cD$
to $\cC\,\Ti\,\cD$. Suppose now that $\cC$ and $\cD$ are ribbon. Instead of
directly verifying that $F$ is a ribbon equivalence, it is slightly more 
convenient to work with its functorial inverse, to be denoted by $G$. On 
objects $R\eq(U{\times}X;\pi)$ of $\cC\Ti\cD$ we have 
$G(R)\eq( (U;\id_U) \ti (X;\id_X) ; \pi )$, while on morphisms $G$ acts as 
the identity map. Using the definition of the ribbon structure on the Karoubian
envelope of a category and on the Karoubian product of \cats, as given in 
\Remark \ref{K-rem1}(iv) and in \Proposition \ref{prop:Ti-mod}, 
respectively, one verifies by direct substitution that $G$ is an 
equivalence of ribbon categories. We present details of the
calculation only for the tensor product and for the braiding.
\\[.1em]
Let $R\eq(U{\times}X;\pi)$ and $S\eq(V{\times}Y;\varpi)$ be objects
of $\cC\Ti\cD$. Using \erf{eq:K-tensor} and \erf{Ti-tensor} we get
  \be \bearll
  G(R \Oti^{\cC\Ti\cD} S) \!\!
  &= G\Llb ( (U{\otimes^{\cC}}V) \ti (X{\otimes^{\cD}}Y) ;
    \pi{\otimes^{\cC\Tic\cD}}\varpi ) \Lrb  \\{}\\[-.7em]
  &= ( (U{\otimes^{\cC}}V ; \id_{U{\otimes^{\cC}}V}) \ti
       (X{\otimes^{\cD}}Y ; \id_{X{\otimes^{\cD}}Y}) ;
       \pi{\otimes^{\cC\Tic\cD}}\varpi )
  \eear \ee
as well as
  \be
  G(R) \,{\otimes^{\kar\cC\Ti\kar\cD}}\, G(S)
  = ((U;\id_U)\ti(X;\id_X);\pi) \otimes^{\kar\cC\Tic\kar\cD}
  ((V;\id_V)\ti(Y;\id_Y);\varpi) \,,  \ee
so that indeed $G(R\Oti^{\cC\Ti\cD} S)\eq G(R)\,{\otimes^{\kar\cC\Ti\kar\cD}}
G(S)$. For morphisms, equality of $G(f\Oti^{\cC\Ti\cD} g)$ and 
$G(f)\,{\otimes^{\kar\cC\Ti\kar\cD}}G(g)$ is immediate
because $G$ is the identity on morphisms.
\\[.1em]
Concerning the braiding note that, using \erf{eq:216} and \erf{Ti-ribbon},
  \be
  G(c_{R,S}^{}) = G( c_{(U\times X;\pi),(V\times Y;\varpi)}^{} )
  = G( (\varpi {\otimes^{\cC\Tic\cD}} \pi) \circ (c_{U,V}^{}\Otic c_{X,Y}^{}) )
  \ee
and
  \be
  c_{G(R),G(S)}^{} =
  c^{}_{( (U;\iD_U)\times(X;\iD_X) ;\pi ),( (V;\iD_V)\times(Y;\iD_Y) ;\varpi )}
  = (\varpi{\otimes^{\cC\Tic\cD}}\pi) \circ (c_{U,V}^{} \otic c_{X,Y}^{}) )
  \,.  \ee
Since $G$ is the identity on morphisms, this implies that
$G(c_{R,S}^{})\eq c_{G(R),G(S)}^{}$.  
\qed

\dtl{Remark}{boxbox}
The product $\Tic$ of categories is associative. Together with lemma
\ref{boxtimes2}, this implies in particular that the Karoubian product of
\cats\ is associative as well, i.e.\ we have
  \be  
  (\cC\Ti\cD) \,\Ti\, \mathcal{E} 
  \,\cong\, (\cC\Tic\cD) \,\Ti\, \mathcal{E} 
  \,\cong\, \kar{(\cC\Tic\cD\Tic\mathcal{E})}
  \,\cong\, \cC \,\Ti\, (\cD\Tic\mathcal{E})
  \,\cong\, \cC \,\Ti\, (\cD\Ti\mathcal{E})  \labl{eq:boxbox}
for any triple $\cC,\;\cD,\;\mathcal{E}$ of \cats.
If $\cC,\,\cD,$ and $\mathcal{E}$ are ribbon, then these are
equivalences of ribbon categories.  

\dtl{Lemma}{prodvect}
For any (additive, $\koerper$-linear) \cat\ $\cC$, taking the product, 
in the sense of \erf{ctimes}, with the \cat\ $\Vectk$ of \findim\
vector spaces yields a \cat\ equivalent to $\cC$,
  \be  \cC\,\Tic\,\Vectk \,\cong\, \cC \,,  \labl{516a}
while taking the Karoubian product with $\Vectk$ yields the Karoubian 
envelope of $\cC$,
  \be  \cC\,\Ti\,\Vectk \,\cong\, \kar\cC \,.  \labl{516b}
If $\cC$ is ribbon, then these are equivalences of ribbon categories.

\medskip\noindent
Proof:\\
Consider the functor $F{:}\; \cC\,{\to}\,\cC\Tic\Vectk$ defined by
$F(U)\,{:=}\,U{\times}\,\koerper$ on objects and by $F(f)\,{:=}%\,
    $\linebreak[0]$%
f\Otic\id_\koerper$ on morphisms. Clearly, $F$ 
is bijective on morphisms. Next, note that every object $X\iN\Obj(\Vectk)$ is 
isomorphic to a direct sum $X\,{\cong}\;\koerper\,{\oplus}\cdots{\oplus}\,
\koerper$. Furthermore we have an isomorphism $(U{\oplus}\cdots{\oplus}\,U) 
\ti \koerper \;{\cong}\; U \ti (\koerper{\oplus}\cdots{\oplus}\,\koerper)$.
Thus every object $U\ti X$ of $\cC\Tic\Vectk$ is isomorphic
to an object of the form $U'\ti\koerper$, implying in particular
that $F$ is essentially surjective, and hence provides an equivalence
of categories by \Proposition \ref{XI.1.5}. This establishes \erf{516a}.
\\
Suppose now that $\cC$ is ribbon. Using the definition of the ribbon
structure on $\cC\,\Ti\,\Vectk$ as given in \Proposition \ref{prop:Ti-mod},
one immediately verifies that in this case $F$ is a ribbon functor.
\\[.2em]
The equivalence \erf{516b} is obtained from \erf{516a} by taking the
Karoubian envelope on both sides, using \Lemma \ref{le:K-funct}.  
\qed

\dt{Lemma}
(i)~\,\,When $A$ and $B$ are algebras in tensor categories $\cC$ and $\cD$,
respectively, then setting
  \be  m^{\cC\Tic\cD}_{\!A\times B} := m_\AA^\cC \otic m_B^\cD
  \qquad{\rm and}\qquad
  \eta^{\cC\Tic\cD}_{A\times B} := \eta_A^\cC \otic \eta_B^\cD
  \labl{meta-prod}
endows $A{\times}B\iN\Obj(\cC\Tic\cD)$ 
with the structure of an algebra in $\cC\Tic\cD$.
\\[.3em]
(ii)~\,An analogous statement holds for coalgebras, with
  \be  \Delta^{\cC\Tic\cD}_{\!A\times B} := \Delta_\AA^\cC \otic \Delta_B^\cD
  \qquad{\rm and}\qquad
  \eps^{\cC\Tic\cD}_{A\times B} := \eps_A^\cC \otic \eps_B^\cD
  \,. \labl{Deps-prod} 
\vskip.3em\noindent
(iii)~If $A$ and $B$ are haploid, then so is $A\ti B$.
\\[.3em]
(iv)~If in addition $\cC$ and $\cD$ are braided and $A$ and $B$ are
(co-)\,commutative, then $A\ti B$ is (co-)\,commutative as well.
\\[.3em]
(v)~\,When $A$ and $B$ are Frobenius \alg s in ribbon
categories $\cC$ and $\cD$, respectively, then \erf{meta-prod} and
\erf{Deps-prod} equip $A\ti B\iN\Obj(\cC\Tic\cD)$ with the structure of a
Frobenius \alg\ in $\cC\Tic\cD$. If in addition both $A$ and $B$ are
symmetric and/or special, then so is $A\ti B$.

\medskip\noindent
Proof:\\
All required relations of the structural morphisms
$m^{\cC\Tic\cD}_{\!A\times B}$, $\eta^{\cC\Tic\cD}_{A\times B}$ etc.\
easily follow from the corresponding ones of $A$ and $B$.
\qed

\medskip
Just like in many other respects, special Frobenius \alg s are especially 
well-behaved also with respect to taking product \cats. In particular, we have
\\[-2.3em] 

\dtl{Lemma}{CKA-KCA-box}
For $A$ and $B$ special Frobenius \alg s in (not necessarily Karoubian)
ribbon \cats\ $\cC$ and $\cD$, respectively, there is an equivalence
  \be  (\cC\Ti\cD)_{(A\times B;\iD_\AA\Otic\iD_B)}^{}
  \cong \kar{\Llb (\cC\Tic\cD)_{\!A\times B}^{} \Lrb} \,.  \labl{CKA6}
If $A$ and $B$ are in addition symmetric and commutative, then 
there is also an equivalence
  \be  \Ext{(\cC\Ti\cD)}{(A\times B;\iD_\AA\Otic\iD_B)}
  \cong \kar{\Llb \Ext{(\cC\Tic\cD)}{\!A\times B} \Lrb} \,.  \labl{CKA8}
     % could also show: ribbon
involving \cats\ of local modules.

\medskip\noindent
Proof:\\
The assertions follow immediately by applying corollary
\ref{cor:CAK-CKA-mod}(i) and (ii), respectively, to the special Frobenius 
\alg\ $A\ti B$ in the ribbon \cat\ $\cC\Tic\cD$.
\qed

\medskip

In the sequel we will often identify $\Obj(\cC\Tic\cD)$ with the
corresponding full sub\cat\ of $\Obj(\cC\Ti\cD)$, and accordingly identify
the \alg\ $(A{\times} B;\id_\AA\Otic\id_B)$ with the \alg\ 
$A{\times}B \iN \Obj(\cC\Tic\cD)\,{\subseteq}\,\Obj(\cC\Ti\cD)$.  

\medskip

A natural question is to which extent the modules over $A\ti B$ can be
understood in terms of $A$- and $B$-modules.
We first note
\\[-2.3em]

\dtl{Lemma}{boxtimes3}
(i)~\,\,For $A$ and $B$ algebras in
tensor categories $\cC$ and $\cD$, and $A{\times}B\iN\Obj(\cC\Tic\cD)$ 
endowed with the algebra structure \erf{meta-prod}, we have the equivalence
  \be  \calcai \,\Tic\, \cD_{\!B}^\ind \,\cong\,
  (\cC\Tic\cD)_{\!A\times B}^\ind  \labl{eq:boxtimes3}
of \cats\ of induced modules.
\\[.3em]
(ii) If in addition $\cC$ and $\cD$ are (not necessarily Karoubian) ribbon 
\cats\ and $A$ and $B$ are \csplit\ commutative \ssFA s then we have the 
equivalence
  \be  \calcali \,\Tic\, \cD_{\!B}^\lInd
  \,\cong\,(\cC\Tic\cD)_{\!A\times B}^\lInd \labl{eq:boxtimes3ii}
of \cats\ of locally induced modules.

\medskip\noindent
Proof:\\
(i)~\,\,The induced $A{\times}B$-modules in $\cC\Tic\cD$ are pairs 
consisting of objects $(A\Oti U)\ti(B\Oti X)$ and the $A{\times}B$-action
$(m_A\Oti\id_U)\otic(m_B\Oti\id_X)$. They are thus in natural bijection
with the objects $(A\Oti U,m_A\Oti\id_U)\ti(B\Oti X,m_B\Oti\id_X)$ of
$\calcai\Tic\cD_{\!B}^{\ind}$. Analogously there are natural
isomorphisms between the respective morphism spaces.
\\[.3em]
(ii) follows from (i) because also the idempotents \erf{PU-def} in the 
two categories that define the locally induced modules coincide.  
\qed

\medskip

The following is yet another result for which it is essential that
the algebras are special Frobenius:
\\[-2.3em]

\dtl{Proposition}{prop:boxtimes}
(i)~\,\,For $A$ and $B$ special Frobenius algebras in (not necessarily 
Karoubian) ribbon categories $\cC$ and $\cD$, there is an equivalence
  \be  \calca \,\Ti\, \cD_{\!B} \,\cong\, (\cC\Ti\cD)_{\!A\times B}^{}
  \labl{eq:boxtimes} 
of \cats.
\\[.3em]
(ii)~\,If in addition $A$ and $B$ are \csplit, symmetric and commutative, 
then there is an equivalence
  \be  \calcal \,\Ti\, \Ext\cD{\!B} \,\cong\, 
  \Ext{(\cC\Ti\cD)}{\AA\times B}   \labl{ClocBoxDloc}
of ribbon categories.

\medskip\noindent
Proof:\\
We combine the \Lemmata \ref{boxtimes2}, \ref{CKA-KCA-box} and \ref{boxtimes3}, 
\Proposition \ref{sub-of-ind-ii} and corollary \ref{cor:loc-lind}.
\\[.3em]
(i)~\,We have
  \be  \bearll \calca \,\Ti\, \cD_{\!B} \!\!
   &\cong\, \kar{(\calca)} \,\Ti\, \kar{(\cD_{\!B})}
  \,\cong\, \kar{(\calcai)} \,\Ti\, \kar{(\cD_{\!B}^\ind)}
  \\{}\\[-.6em]
   &\cong\, \calcai \,\Ti\, \cD_{\!B}^\ind
  \,\equiv\, \kar{(\calcai\Tic\cD_{\!B}^\ind)}
  \\{}\\[-.6em]
   &\cong\, \kar{((\cC\Tic\cD)_{\AA\times B}^\ind)}
  \,\cong\, \kar{((\cC\Tic\cD)_{\AA\times B}^{})}
  \,\cong\, (\cC\Ti\cD)_{(A\times B;\iD_\AA\Otic\iD_B)}^{}
  \,, \eear \ee
where in the first line we use first \erf{eq:boxtimes2} and then \erf{calcai},
in the second line again \erf{eq:boxtimes2}, and in the last line
\erf{eq:boxtimes3}, \erf{calcai} and finally \erf{CKA6}.
\\[.3em]
(ii)~Analogously,
  \be  \bearll \calcal \,\Ti\, \Ext\cD{\!B} \!\!
   &\cong\, \kar{(\calcal)} \,\Ti\, \kar{(\Ext\cD{\!B})}
  \,\cong\, \kar{(\calcali)} \,\Ti\, \kar{(\cD_{\!B}^\lInd)}
  \\{}\\[-.6em]
   &\cong\, \calcali \,\Ti\, \cD_{\!B}^\lInd
  \,\equiv\, \kar{(\calcali\Tic\cD_{\!B}^\lInd)}
  \\{}\\[-.6em]
   &\cong\, \kar{((\cC\Tic\cD)_{\AA\times B}^\lInd)}
  \,\cong\, \kar{(\Ext{(\cC\Tic\cD)}{\AA\times B})}
  \,\cong\, \Ext{(\cC\Ti\cD)}{\AA\times B}
  \,, \eear \labl{633}
where in the first line we use first \erf{eq:boxtimes2} and then \erf{calcali},
in the second line again \erf{eq:boxtimes2}, and in the last line
\erf{eq:boxtimes3ii}, \erf{calcali} and finally \erf{CKA8}.  
\\[.2em]
Next we note that, by corollary \ref{cor:loc-lind}, objects of 
$\calcal \,\Ti\, \Ext\cD{\!B}$ can be written as $( (\lxt U\AA;p){\times}
    $\linebreak[0]$%
(\lxt X B;q) ; \pi )$ with $U\iN\Objc$, $X\iN\Obj(\cD)$,
$p$ and $q$ the respective idempotents that describe a local module
as \retmodule\ of a locally induced module, and $\pi$ the idempotent
that arises in taking the Karoubian envelope of $\calcal\Tic\Ext\cD{\!B}$.
Similarly, objects of $\Ext{(\cC\Ti\cD)}{\AA\times B}$ can be written as
$( \lxt {(V{\times}Y;\varpi)} {\AA \times B} ; \hat\pi )$ with
$V\iN\Objc$, $Y\iN\Obj(\cD)$, $\varpi$ the idempotent arising in
taking the Karoubian envelope of $\cC\Tic\cD$, and $\hat\pi$ the idempotent
describing a local $A{\times}B$-module as \retmodule\ of a locally induced 
$A{\times}B$-module.
\\
With this description of the objects, the functor 
$F{:}\;\calcal\Ti\Ext\cD{\!B}\,{\stackrel\cong\to}\, \Ext{(\cC\Ti\cD)}
{\AA\times B}$ that maps the \lhs\ of \erf{633} to the \rhs\ is given by
  \be
  F:\quad ( (\lxt U\AA;p) \ti (\lxt X B;q) ; \pi ) \,\mapsto\,
  ( \lxt {(U{\times X};\id_{U{\times}X})} {\AA \times B} ; \pi )
  \labl{634}
on objects, and is the identity map on morphisms, with the latter regarded
as elements in (a subspace of)
$\Hom^{\cC\Tic\cD}((A \Oti U) \ti (B \Oti X) , (A \Oti V) \ti (B \Oti Y))$.
(That the idempotents $p$ and $q$ do not appear on the \rhs\ of \erf{634}
is seen by the same reasoning as in the proof of \Lemma \ref{boxtimes2}.)
\\
Now one checks by inserting the relevant definitions -- formula
\erf{Ti-tensor} for the tensor product on products of \cats, formula 
\erf{eq:K-tensor} for the tensor product on the Karoubian envelope of a \cat, 
as well as formula \erf{eq:multi-tensor-obj} for the tensor product of local 
modules -- that the prescription \erf{634} respects the tensor product, i.e.\
$R\,{\otimes^{\Ext\cC{\!A}\Ti\Ext\cD{\!B}}}S \,{\stackrel F\mapsto}\,
F(R)\,{\otimes^{\Ext{(\cC\Ti\cD)}{\AA\times B}}}F(S)$
(together with an analogous equality for the tensor product of morphisms,
which follows trivially). Thus $F$ is a tensor functor.
\\
Similarly, using the formulas \erf{Ti-ribbon} for the braiding on products of 
\cats, \erf{eq:216} for the braiding on the Karoubian envelope, and 
\erf{eq:Cloc-braid} for the braiding of local modules, one verifies that the
braidings on $\calcal\Ti\Ext\cD{\!B}$ and on $\Ext{(\cC\Ti\cD)}{\AA\times B}$ 
are compatible in the sense that $c^{\Ext\cC{\!A}\Ti\Ext\cD{\!B}}_{R,S}{=}\;
c^{\Ext{(\cC\Ti\cD)}{\AA\times B}}_{F(R),F(S)}$. Since $F$ is the identity 
on morphisms, this means that $F$ is braided, and hence 
  % by the uniqueness properties of the left and right dualities and the
  % fact that the twist can be expressed through the dualities and the braiding
that $F$ is a ribbon functor.
\qed 

\dtl{Corollary}{CAxD0}
If $\cC$ and $\cD$ are (not necessarily Karoubian) ribbon categories and $A$ 
is a \csplit\ commutative \ssFA\ in $\cC$, then there are equivalences
  \be
  (\cC\Tic\cD)_{\!A\times\one_\cD}^\lInd \,\cong\, \calcali \,\Tic\, \cD 
  \qquad {\rm and} \qquad
  \EXt{\cC\Ti\cD}{\!A\times\one_\cD} \,\cong\, \calcal \,\Ti\, \cD ~.
  \labl{eq:CAxD1}
The first is an equivalence of categories, the second an equivalence of 
ribbon categories.

\medskip\noindent
Proof:\\
These equivalences follow by setting $B\eq\one_\cD$ in the equivalences
\erf{eq:boxtimes3ii} and \erf{ClocBoxDloc}, respectively.
\qed

\medskip

Before we specialise to a special situation of particular interest --
$\cC$ a \mtc\ and $\cC'$ being dual to $\cC$ --  
let us mention that another large class of \platl\ pairs $\cC$ and $\cC'$ 
is provided by conformal embeddings similar to those listed in \erf{E8}.

%%%%%%%%%%%%%%%%%%%%%%%%%%%%%%%%%%%%%%%%%%%

\subsection{The dual of a tensor category}\label{platy2}

As already mentioned above, an important class of \platl\ \cats\
is given by modular \tcs, and for these $\cC'$ is the dual of 
$\cC$. We therefore turn to the discussion of the concept of dual \tc.
\\[-2.3em]

\dtl{Definition}{def:dual-cat}
The {\em dual category\/} $\ol\cC$ of a tensor category $(\cC, \otimes)$ is
the \tc\ $(\cC^{\rm opp}, \otimes)$.

\medskip\noindent
More concretely, when marking quantities in $\ol\cC$ by an overline, we have
  \be\bearll
  {\rm Objects:}  &  \quad
    \Obj(\ol\cC) = \Obj(\cC)\,, \ {\rm i.e.} \ \
    \ol U \iN \Obj(\ol\cC) \ \ {\rm iff}\ \ U \iN \Obj(\cC)\,, \\[5pt]
  {\rm Morphisms:} &  \quad
    \ol\Hom(\ol U,\ol V) = \Hom(V,U)\,, \\[5pt]
  {\rm Composition:} &  \quad
    \ol f \cirb \ol g = \ol{g{\circ}f}\,,   \\[5pt]
  {\rm Tensor\ product:} &  \quad
    \ol U \otib \ol V = \ol{U{\otimes}V}\,, \quad
    \ol f \otib \ol g = \ol{f{\otimes}g}\,, \\[5pt]
  {\rm Tensor\ unit:} &  \quad
    \oneb = \one \,.
  \eear\ee

\dt{Remark}
(i)~\,Since $\cC$ is strict, $\ol\cC$ is indeed again a (strict) \tc.
If the tensor category $\cC$ is small, then so is $\ol\cC$.
If $\cC$ is additive, then so is $\ol\cC$.
If $\cC$ is semisimple, then so is $\ol\cC$.
\\[.3em]
(ii)~If the tensor category $\cC$ is Karoubian, then so is $\ol\cC$.
More generally, since the idempotents in $\cC$ coincide with the idempotents 
in $\ol\cC$, for any tensor category $\cC$ the Karoubian envelope of 
$\ol\cC$ is the dual category of the Karoubian envelope of $\cC$, 
i.e.\ $\kar{\ol\cC}{=}\,\ol{\kar\cC}$.

\medskip

The following result is analogous to lemma 2.9 of \cite{muge8}:
\\[-2.3em]

\dtl{Lemma}{lem:dual-mod}
(i)~\,If the tensor category $\cC$ has a left (right) duality, then its
dual category $\ol\cC$  has a right (left) duality. If $\cC$ has a braiding, 
then so has $\ol\cC$, and if $\cC$ has a twist, then so has $\ol\cC$.
\\
In particular, the dual $\ol\cC$ of a ribbon category $\cC$ is
naturally a ribbon category, too.
\\
The values of $s$ for $\cC$ and $\ol\cC$ are related via
  \be  \ol s_{\ol U,\ol V}^{} = s_{U,V^\vee}^{}
  \quad (\; = s_{U^\vee,V}^{} \;) \,,  \labl{ol-s}
so that in particular
  \be  \ol{\rm dim}(\ol U) = \dim(U) \,.  \ee
\vskip.3em\noindent
(ii)~The dual category $\ol\cC$ of a modular tensor category $\cC$
carries a natural structure of a modular tensor category.

\medskip\noindent
Proof:\\
(i)~\,We set
  \be  \ol U^\vee := \ol{{}^{\vee\!}_{\phantom|}U} \,, \qquad
  {}^\vee_{\phantom i}\ol U:= \ol{U^\vee_{\phantom|}}  \ee
and
  \be\bearll
  {\rm Dualities:}&  \quad
    {\ol b}_{\ol U} := \ol{(\tilde d_U)}
    \in \ol\Hom(\ol\one,\ol U\otib\ol U^\vee) \,, \quad \
    {\ol d}_{\ol U} := \ol{(\tilde b_U)}
    \in \ol\Hom(\ol U^\vee\otib\ol U,\ol\one) \,, \\{}\\[-.5em]
  &\quad
    {\ol{\tilde b}}_{\ol U} := \ol{(d_U)}
    \in \ol\Hom(\ol\one,{}^\vee_{\phantom i}\ol U\otib\ol U) \,, \quad \
    {\ol{\tilde d}}_{\ol U} := \ol{(b_U)}
    \in \ol\Hom(\ol U\otib{}^\vee_{\phantom i}\ol U,\ol\one) \,, \\{}\\[-.5em]
  {\rm Braiding:} &  \quad
    {\ol c}_{\ol U,\ol V} := \ol{(c_{U,V})^{-1}}
    \in \ol\Hom(\ol U\otib\ol V,\ol V\otib\ol U) \,, \\{}\\[-.5em]
  {\rm Twist:}    &  \quad
    \ol\theta_{\ol U} := \ol{(\theta_U^{-1})} \in \ol\Hom(\ol U,\ol U) \,.
  \eear\ee
By direct substitution one verifies that these morphisms
satisfy all properties of dualities, braiding and twist.
\\
For $s$ as defined by \erf{sUV} one computes
  \be\bearll \ol s_{\ol U,\ol V}^{} \!\!
  &= (\ol d_{\ol V}^{}\otib {\ol{\tilde d}}_{\ol U}) \circb
     [\, \id_{\ol V^\vee}^{}\otib (\ol c_{\ol U,\ol V}^{} \circb
     \ol c_{\ol V,\ol U}^{}) \otib \id_{\ol U^\vee}^{} \,]
     \circb ({\ol{\tilde b}}_{\ol V}\otib \ol b_{\ol U}^{})
  \\{}\\[-.7em]
  &= (\ol{(\tilde b_V)}\otib\ol{(b_U)}) \circb
     [\, \id_{\ol{V^\vee}}^{}\otib (\ol{(c_{U,V})^{-1}} \circb
     \ol{(c_{V,U})^{-1}}) \otib \id_{\ol{U^\vee}}^{} ]
     \circb (\ol{(d_V)}\otib\ol{(\tilde d_U)})
  \\{}\\[-.7em]
  &= \ol{ (d_V\oti\tilde d_U) \circ [\, \id_{U^\vee}^{}\oti
     ((c_{V,U})^{-1} \circ (c_{U,V})^{-1}) \oti \id_{V^\vee}^{} ]
     \circ (\tilde b_V\oti b_U) }
  \\{}\\[-.7em]
  &= s_{U,V^\vee}^{} = s_{U^\vee,V}^{} \,.
  \eear\ee
The manipulations leading to the last two equalities may be summarised
in the language of ribbon graphs, analogously as in \erf{sUV}:
The second-to-last corresponds to a $180^\circ$ rotation of the $V$-ribbon, 
and the last to a $180^\circ$ rotation of the $U$-ribbon.
\\[.3em]
(ii)~The simple objects of $\ol\cC$ are $\ol V$ with $V$ a simple
object of $\cC$; in particular, $\ol\cC$ has as many isomorphism
classes of simple objects as $\cC$ has. Finally, owing to \erf{ol-s}
invertibility of the matrix $\ol s\,{\equiv}\,(\ol s_{i,j}^{})$ follows
immediately from invertibility of $s$.
\qed

\dtl{Remark}{rem:duals}
As in \Remarks \ref{rem:Dk-1}(i) and \ref{rem1box}(iv) we may consider 
the behaviour of the dimension and charge of a modular tensor category. One
verifies that under taking duals one has
  \be
  {\rm Dim}(\ol\cC) = {\rm Dim}(\cC) \qquad{\rm and}\qquad
  p^\pm(\ol\cC) = p^\mp(\cC) \,.  \ee 

\dt{Lemma}
(i)~~If $(A,m,\eta)$ is an \alg\ in a \tc\ $\cC$, then $(\ol A,\ol m,\ol\eta)$ 
is a co\alg\ in $\ol\cC$, and if $(A,\Delta,\eps)$ is a co\alg\ in $\cC$, then
$(\ol A,\ol\Delta,\ol\eps)$ is an \alg\ in $\ol\cC$.
\\[.3em]
(ii)~If $(A,m,\eta,\Delta,\eps)$ is a (commutative) \ssFA\ in a
ribbon \cat\ $\cC$, then $(\ol A,\ol\Delta,\ol\eps,\ol m,\ol\eta)$ is a 
(commutative) \ssFA\ in $\ol\cC$.

\medskip\noindent
Proof:\\
The relevant properties in the dual \cat\ are nothing but the
corresponding properties of the dual morphisms.
\qed

\medskip

For the rest of this subsection we assume that $\cC$ is a \tc\ with a
finite number of isomorphism classes of simple objects, i.e.\ that
the index set $\II$ (see \Section \ref{sect21}) 
is finite. Then for every triple of simple objects $U_i$, $U_j$, $U_k$ with
$i,j,k\iN\II$ we fix once and for all a basis
$\{\alpha\}\,{\subset}\,\Hom(U_i{\otimes}U_j,U_k)$ and a dual basis
$\{\alpha\}\,{\subset}\,\Hom(U_k,U_i{\otimes}U_j)$.%
 \foodnode {See \Section 2.2 of \cite{fuRs4} for more details. There the
 notation $\bar\alpha$ was used for the second type of basis elements;
 here the overbar is suppressed to avoid confusion with quantities
 referring to the dual \cat\ $\ol\cC$.}
Then the 6j-symbols, or fusing matrices, \FF, of $\cC$ and their inverses
\Gmat\ are defined by
(in the figures we abbreviate the simple objects $U_i$ by their labels $i$)
  \bea \begin{picture}(240,60)(0,38)
  \put(0,0)   {\begin{picture}(0,0)(0,0)
              \scalebox{.38}{\includegraphics{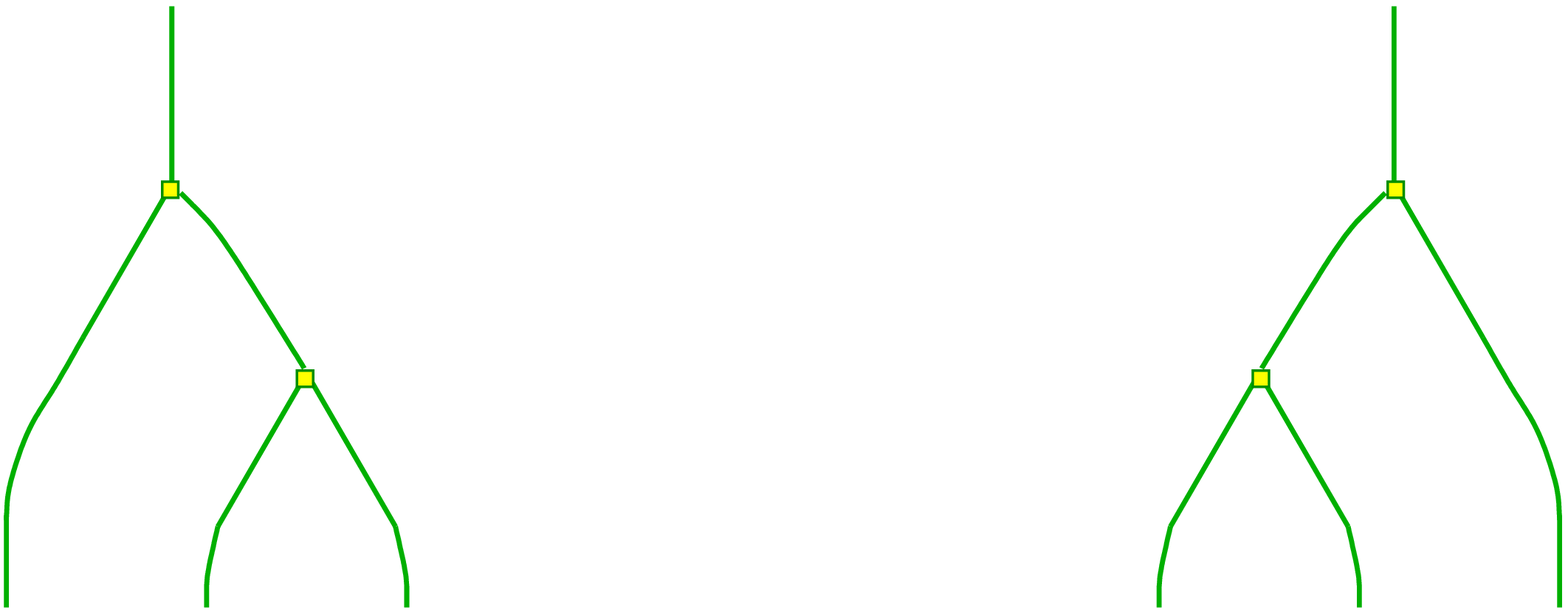}} \end{picture}}
  \put(-0.7,-8)  {\sse$i$}
  \put(18.2,61.7){\tiny$\alpha$}
  \put(24.8,93.5){\sse$l$}
  \put(28.1,-8)  {\sse$j$}
  \put(42.1,48)  {\sse$p$}
  \put(38.6,33.1){\tiny$\beta$}
  \put(58.1,-8)  {\sse$k$}
  \put(76,42)    {$=\ \dsty\sum_{q\in\II}\sum_{\gamma,\delta}\;
                 \F{i\U j}klpq\alpha\beta\gamma\delta$}
  \put(170.7,-8) {\sse$i$}
  \put(181.1,34.5){\tiny$\gamma$}
  \put(190.1,50) {\sse$q$}
  \put(199.5,-8) {\sse$j$}
  \put(210.4,61.6){\tiny$\delta$}
  \put(206.6,93.5){\sse$l$}
  \put(229.5,-8) {\sse$k$}
  \epicture26 \labl{def-fmat}
  \bea \begin{picture}(220,60)(0,38)
  \put(0,10)   {\begin{picture}(0,0)(0,0)
              \scalebox{.38}{\includegraphics{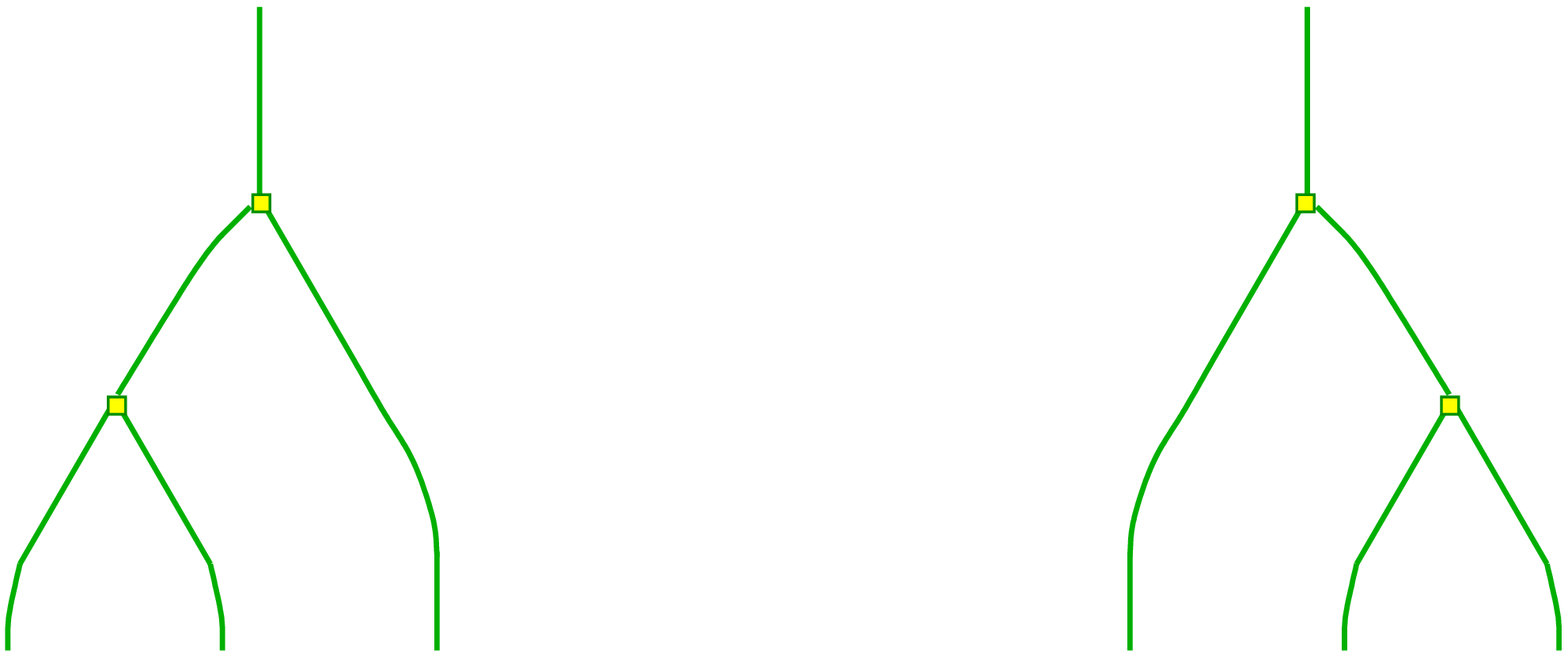}} \end{picture}}
  \put(34.8,103.8) {\sse$l$}
  \put(29.3,73)    {\tiny$\beta$}
  \put(17,58)      {\sse$p$}
  \put(8.8,45)     {\tiny$\alpha$}
  \put(-1.5,2.8)   {\sse$i$}
  \put(27.3,2.8)   {\sse$j$}
  \put(58,2.8)     {\sse$k$}
  \put(70,42)      {$=\ \dsty\sum_{q,\gamma,\delta}\;
                     \G{i\U j}klpq\alpha\beta\gamma\delta$}
  \put(181.2,103.8){\sse$l$}
  \put(174.8,71)   {\tiny$\gamma$}
  \put(194,62)     {\sse$q$}
  \put(195.8,43)   {\tiny$\delta$}
  \put(156,2.8)    {\sse$i$}
  \put(185,2.8)    {\sse$j$}
  \put(214,2.8)    {\sse$k$}
  \epicture28 \labl{gmat}
Furthermore, when $\cC$ is braided, then the braiding matrices \RR\ of
$\cC$ are defined by
  \bea \begin{picture}(150,60)(0,29)
  \put(0,0)   {\begin{picture}(0,0)(0,0)
              \scalebox{.38}{\includegraphics{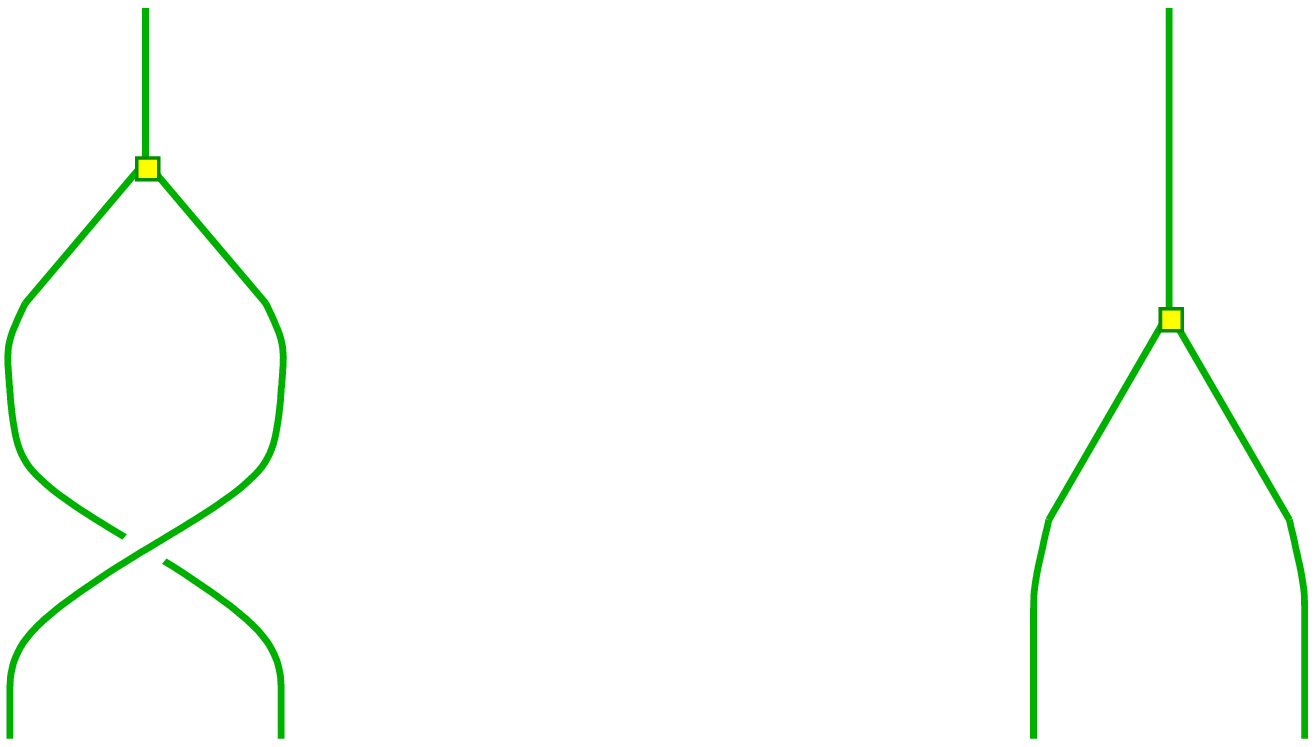}} \end{picture}}
  \put(45,36)      {$=:\ \dsty\sum_{\beta}\R ijk\alpha\beta$}
  \put(-0.7,-7)    {\sse$i$}
  \put(13.9,84)    {\sse$k$}
  \put(17.7,66)    {\tiny$\alpha$}
  \put(28.7,-7)    {\sse$j$}
  \put(112.3,-7)   {\sse$i$}
  \put(126.6,84)   {\sse$k$}
  \put(129.7,50.7) {\tiny$\beta$}
  \put(140.7,-7)   {\sse$j$}
  \epicture19 \labl{rmat}
$\R ijk{}{}$ is a square matrix with rows and columns labelled by the basis 
$\{\alpha\}$ of $\Hom(U_i{\otimes}U_j,U_k)$; its inverse with respect to 
this matrix structure is $\Rm jik{}{}$, which is defined analogously as 
$\R jik{}{}$, but with an under-braiding instead of an over-braiding.

The choice of bases in the spaces $\Hom(U_i{\otimes}U_j,U_k)$
and $\Hom(U_k,U_i{\otimes}U_j)$ of $\cC$ allow us to choose a correlated
basis in $\ol\cC$. For example to pick a basis
$\{\ol\alpha\}\,{\subset}\,\ol{\Hom}(\ol U_i\,{\ol\otimes}\,\ol U_j,\ol U_k)$
we use that by definition
$\ol{\Hom}(\ol U_i\,{\ol\otimes}\,\ol U_j,\ol U_k)\eq\Hom(U_k,U_i{\otimes}U_j)$
and take the basis we have already chosen in the latter.

To simplify notation, in the remainder of the paper we
will omit the overlines on quantities of the dual \cat\ $\ol\cC$ whenever
from the context it is so obvious that $\ol\cC$-quantities are meant
that no confusion can arise. For instance, we write the fusing matrices of
$\ol\cC$ as $\Fol{i\,j}klpq\alpha\beta\gamma\delta$ instead of
$\Fol{\ol i\,\ol j}{\ol k}{\ol l}{\ol p}{\ol q}{\ol\alpha}{\ol\beta}
{\ol\gamma}{\ol\delta}$.

\dtl{Lemma}{lem:dual-basis}
The fusing and braiding matrices of the dual $\ol\cC$ of a
braided tensor category $\cC$ with finite index set $\II$ are given by
  \be
  \Fol{i\,j}klpq\alpha\beta\gamma\delta
    = \G {i\,j}klqp\gamma\delta\alpha\beta \, , \quad
  \Gol{i\,j}klpq\alpha\beta\gamma\delta
    = \F {i\,j}klqp\gamma\delta\alpha\beta \, , \quad
  \Rol ijk\alpha\beta = \Rm jik\beta\alpha \, , \quad
  \Rmol ijk\alpha\beta = \R jik\beta\alpha
  \,.  \ee

\noindent
Proof:\\
It follows from the definition of dual bases that the fusing matrices
also appear in the relation
  \bea \begin{picture}(250,55)(0,40)
  \put(0,0)   {\begin{picture}(0,0)(0,0)
              \scalebox{.38}{\includegraphics{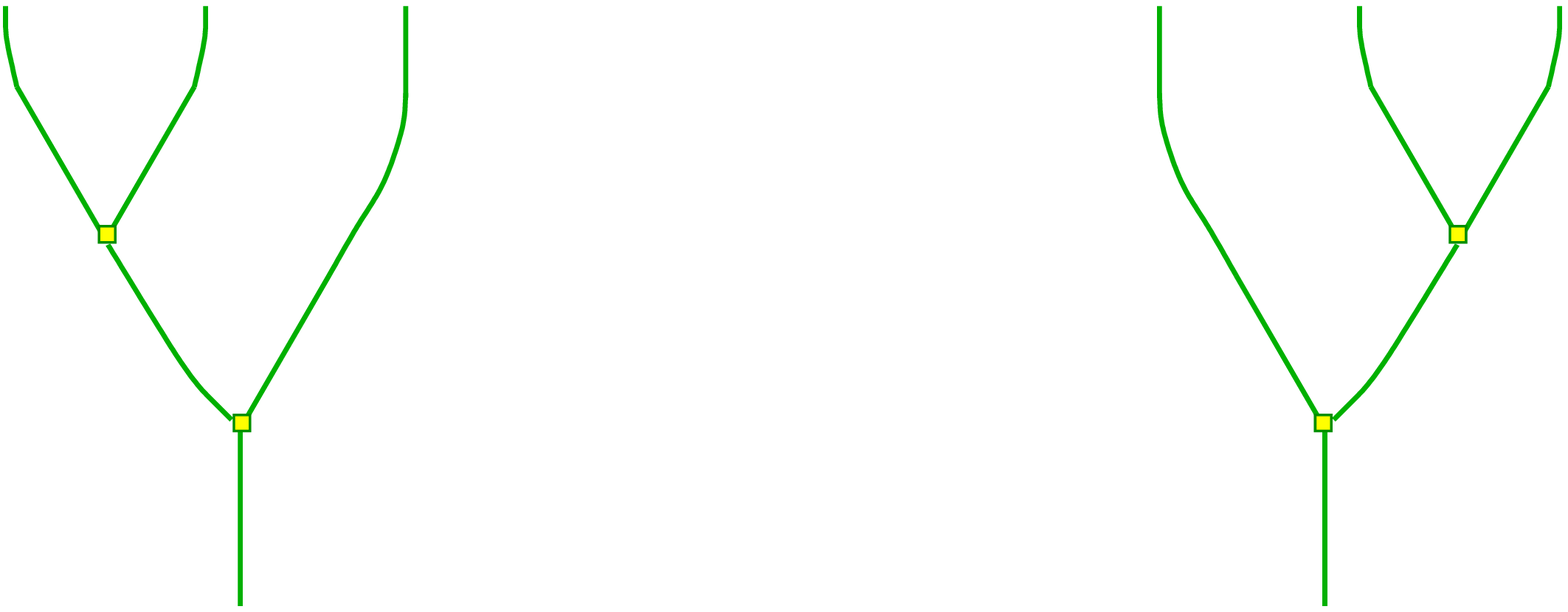}} \end{picture}}
  \put(-0.5,93.5){\sse$i$}
  \put(18.6,54.9){\tiny$\gamma$}
  \put(34.5,-8)  {\sse$l$}
  \put(18.5,38.9){\sse$q$}
  \put(29.4,94.9){\sse$j$}
  \put(38.8,25.7){\tiny$\delta$}
  \put(58.3,93.5){\sse$k$}
  \put(76,42)    {$= \ \dsty\sum_{p\in\II}\sum_{\alpha,\beta}\;
                 \F{i\U j}klpq\alpha\beta\gamma\delta$}
  \put(171.7,93.5){\sse$i$}
  \put(195.9,-8) {\sse$l$}
  \put(200.1,25.7){\tiny$\alpha$}
  \put(201.2,94.9){\sse$j$}
  \put(210.5,38.2){\sse$p$}
  \put(220.2,54.6){\tiny$\beta$}
  \put(230.7,93.5){\sse$k$}
  \epicture25 \labl{fmat-dual}
Combining this result for the category $\cC$ with the definition of the
morphisms $\ol\Hom$ and their composition $\ol\circ$ in $\ol\cC$ one arrives
at the first equality. The other relations follow by an analogous reasoning.
\qed

%%%%%%%%%%%%%%%%%%%%%%%%%%%%%%%%%%%%%%%%%%%

\subsection{The trivialising algebra \boldmath{$T_\cG$}} \label{noname-sec}

Recall that we denote by $\II$ the index set such that $\{U_i\,|\,i\iN\II\}$ 
is a collection of representatives for the equivalence classes of simple 
objects in a \cat. In this subsection we consider ribbon \cats\ $\cG$ 
which are semisimple and have finite index set $\IG$.

We start by introducing an interesting algebra $T\,{\equiv}\,T_\cG$ in
the Karoubian product $\cGG$ of $\cG$ with its dual.
This is done in the following lemma, which is essentially 
\Proposition 4.1 of \cite{muge9}:
\\[-2.3em]

\dtl{Lemma}{lem:TG}
Let $\cG$ be a semisimple ribbon category with a finite number of
equivalence classes of simple objects. 
\\[.2em]
(i)~~The triple $T_\cG\,{\equiv}\, (T_\cG,m,\eta)$ with
  \bea
  T_\cG := \displaystyle\bigoplus_{k\in\IG} U_k {\times} \ol{U_k}\;
       \in \Obj(\cGG) \,,
  \\{}\\[-.95em]
  \eta \,:= e^{}_{\one\times\ol\one \prec T_\cG} \;
  \in \Hom^{\cGG\!}(\one {\times} \ol\one, T_\cG) \,,
  \\{}\\[-1.1em]
  m \,:= \displaystyle\sum_{i,j,k\in\IG}\, \sum_{\alpha}
  %% [pic~6]
  \begin{picture}(120,53)(9,43)
  \put(17,0)  {\begin{picture}(0,0)(0,0)
              \scalebox{.38}{\includegraphics{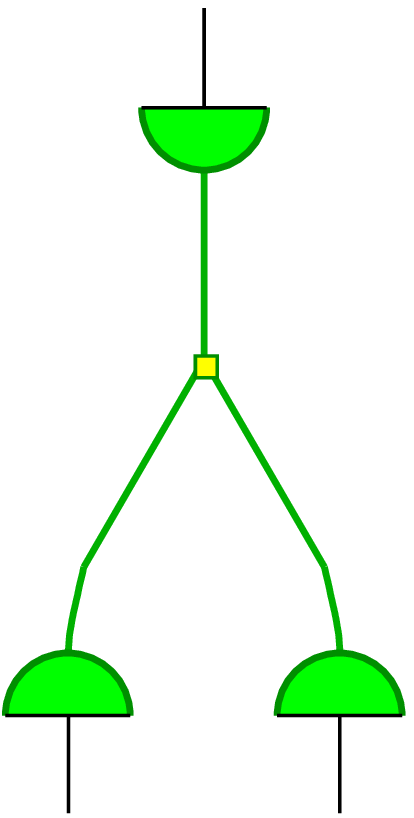}} \end{picture}}
  \put(77,0)  {\begin{picture}(0,0)(0,0)
              \scalebox{.38}{\includegraphics{xijk2.eps}} \end{picture}}
  \put(22.4,-8.5) {\sse$ $}
  \put(25.1,33.5) {\sse$i$}
  \put(32.8,58.5) {\sse$k$}
  \put(34.9,92.5) {\sse$ $}
  \put(42.7,48)   {\tiny$\alpha$}
  \put(50.2,-8.5) {\sse$ $}
  \put(51.9,33.5) {\sse$j$}
  \put(64,43)     {$\Otic$}
  \put(81.8,-8.5) {\sse$ $}
  \put(85.1,33.5) {\sse$\ol\imath$}
  \put(92.6,58.5) {\sse$\ol k$}
  \put(95.3,92.5) {\sse$ $}
  \put(102.8,47.9){\tiny$\ol\alpha$}
  \put(110.7,-8.5){\sse$ $}
  \put(112.5,33.5){\sse$\ol\jmath$}
  \end{picture}
      \in \Hom^{\cGG\!}(T_\cG\Oti T_\cG,T_\cG)
  \\[2.6em] {}
  \eear\labl{TG}
is an \alg\ in $\cGG$.
\\[.3em]
(ii)~\,The \alg\ $(T_\cG,m,\eta)$ extends to a haploid commutative 
\ssFA{} in $\cGG$.

\bigskip\noindent
Proof:
\\[.1em]
(i)~~The unit property of the multiplication $m$ follows from the
normalisation of the morphisms that was chosen in (2.33) of \cite{fuRs4}, 
which states that the basis vector chosen in 
$\Hom(U_i{\otimes}\one,U_i)$ and $\Hom(\one{\otimes}U_i,U_i)$ is $\id_{U_i}$.
\\
To see associativity one notes that
  \bea \begin{picture}(350,300)(20,16)
  \put(24,220) {\begin{picture}(0,0)(0,0)
              \scalebox{.38}{\includegraphics{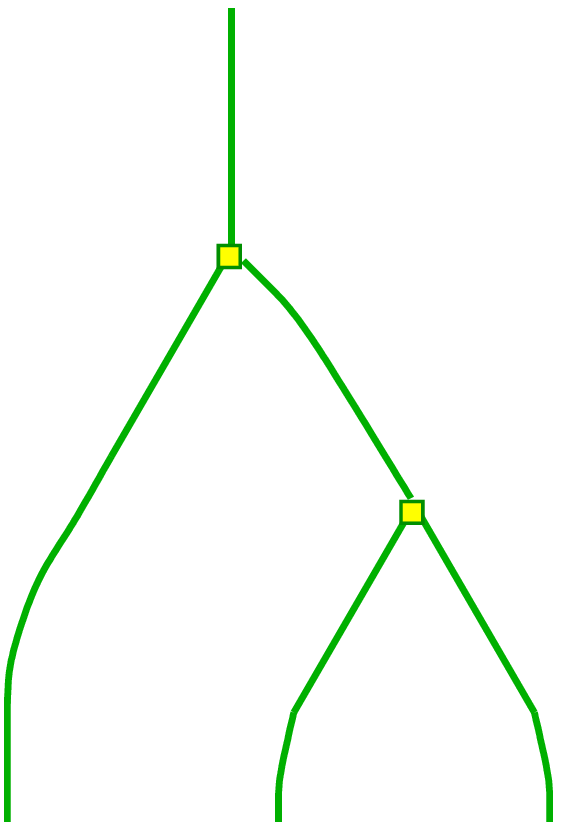}} \end{picture}}
  \put(120,220) {\begin{picture}(0,0)(0,0)
              \scalebox{.38}{\includegraphics{yijkl.eps}} \end{picture}}
  \put(206,110)  {\begin{picture}(0,0)(0,0)
              \scalebox{.38}{\includegraphics{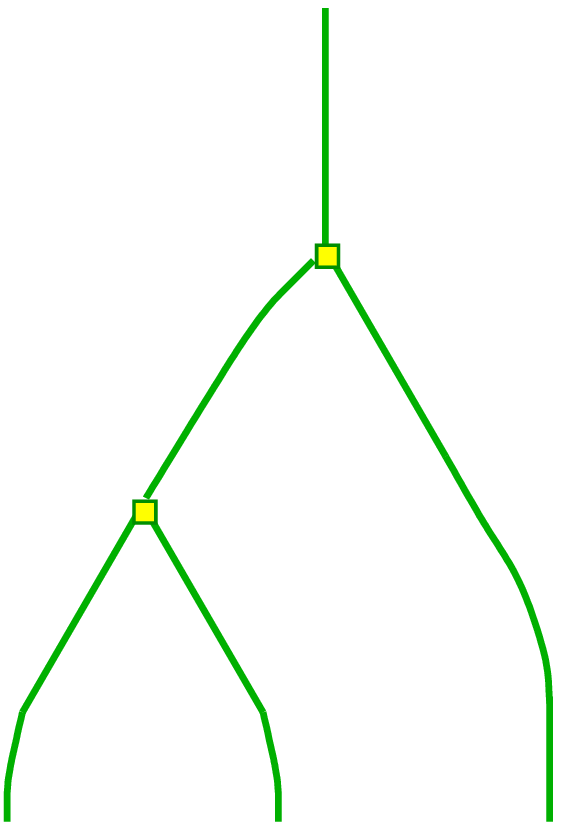}} \end{picture}}
  \put(302,110)  {\begin{picture}(0,0)(0,0)
              \scalebox{.38}{\includegraphics{xijkl.eps}} \end{picture}}
  \put(64,0)   {\begin{picture}(0,0)(0,0)
              \scalebox{.38}{\includegraphics{xijkl.eps}} \end{picture}}
  \put(160,0)  {\begin{picture}(0,0)(0,0)
              \scalebox{.38}{\includegraphics{xijkl.eps}} \end{picture}}
  \put(-10,260)    {$\dsty\sum_{p,\alpha,\beta}$}
  \put(23.2,212.7) {\sse$i$}
  \put(42.2,281.8) {\tiny$\alpha$}
  \put(48.8,314.2) {\sse$l$}
  \put(52.2,212.7) {\sse$j$}
  \put(62.2,254.2) {\tiny$\beta$}
  \put(64.2,271.9) {\sse$p$}
  \put(82.2,212.1) {\sse$k$}
  \put(97,260)     {$\Otic$}
  \put(119.2,212.8){\sse$\ol\imath$}
  \put(138.2,281.8){\tiny$\ol\alpha$}
  \put(144.8,314.2){\sse$\ol l$}
  \put(148.2,212.8){\sse$\ol\jmath$}
  \put(158.2,254.2){\tiny$\ol\beta$}
  \put(160.2,271.9){\sse$\ol p$}
  \put(178.2,211.8){\sse$\ol k$}
  \put(20,150)     {$=\ \dsty\sum_{r,\rho,\rho'}\ \sum_{s,\sigma,\sigma'}\
                    \sum_{p,\alpha,\beta}\ \F{i\,j}klpr\alpha\beta\rho{\rho'}\
                    \Fol{i\,j}klps\alpha\beta\sigma{\sigma'}$}
  \put(205.2,102.7){\sse$i$}
  \put(224.2,144.6){\tiny$\rho$}
  \put(226.2,162.5){\sse$r$}
  \put(234.2,102.7){\sse$j$}
  \put(240.8,204.2){\sse$l$}
  \put(245.3,171.9){\tiny$\rho'$}
  \put(264.2,102.1){\sse$k$}
  \put(279,150)    {$\Otic$}
  \put(301.2,102.7){\sse$\ol\imath$}
  \put(320.9,144.1){\tiny$\ol\sigma$}
  \put(322.2,162.5){\sse$\ol s$}
  \put(330.2,102.7){\sse$\ol\jmath$}
  \put(336.8,204.2){\sse$\ol l$}
  \put(341.3,171.9){\tiny$\ol\sigma'$}
  \put(360.2,101.4){\sse$\ol k$}
  \put(20,40)      {$=\; \dsty\sum_{q,\gamma,\delta}$}
  \put(63.2,-7.9)  {\sse$i$}
  \put(82.2,34.8)  {\tiny$\gamma$}
  \put(83.2,51.9)  {\sse$q$}
  \put(92.2,-7.9)  {\sse$j$}
  \put(98.8,94.2)  {\sse$l$}
  \put(103.3,61.3) {\tiny$\delta$}
  \put(122.2,-8.5) {\sse$k$}
  \put(137,41)     {$\Otic$}
  \put(159.2,-8.2) {\sse$\ol\imath$}
  \put(178.2,34.2) {\tiny$\ol\gamma$}
  \put(179.2,51.9) {\sse$\ol q$}
  \put(188.2,-8.2) {\sse$\ol\jmath$}
  \put(194.8,93.8) {\sse$\ol l$}
  \put(199.3,60.6) {\tiny$\ol\delta$}
  \put(218.2,-9.2) {\sse$\ol k$}
  \epicture10 \labl{pic64} 
The second step uses \Lemma \ref{lem:dual-basis} 
to relate $\ol\FF$ to the inverse of $\FF$.
\\[.3em]
(ii)~\,Thus $T_\cG$ is an algebra. It is clearly haploid. Commutativity
follows from
  %% [pic~65]
  \bea \begin{picture}(340,66)(15,14)
  \put(30,0)  {\begin{picture}(0,0)(0,0)
              \scalebox{.38}{\includegraphics{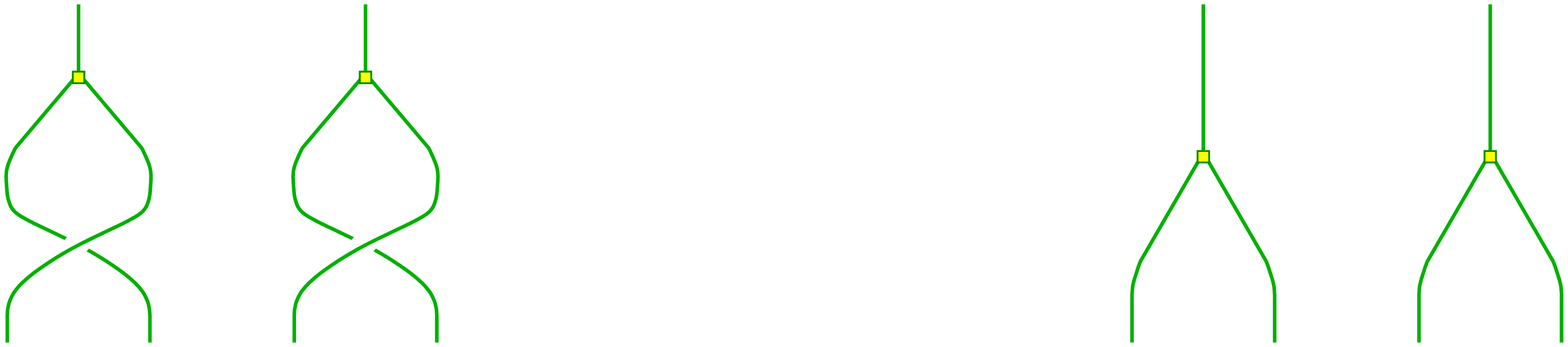}} \end{picture}}
  \put(0,31)       {$\dsty\sum_{\alpha}$}
  \put(29.2,-7.9)  {\sse$i$}
  \put(44.1,74.4)  {\sse$k$}
  \put(48.9,55.4)  {\tiny$\alpha$}
  \put(58.2,-7.9)  {\sse$j$}
  \put(68,31)      {$\Otic$}
  \put(89.2,-8.1)  {\sse$\ol\imath$}
  \put(104.1,74.4) {\sse$\ol k$}
  \put(108.9,55.4) {\tiny$\ol\alpha$}
  \put(118.2,-8.1) {\sse$\ol\jmath$}
  \put(142,31)     {$=\ \dsty\sum_{\beta,\gamma}\,\sum_{\alpha}\,
                    \R ijk\alpha\beta\;\Rol ijk\alpha\gamma $}
  \put(264.2,-7.9) {\sse$i$}
  \put(279.1,74.4) {\sse$k$}
  \put(284.2,38.3) {\tiny$\beta$}
  \put(293.2,-7.9) {\sse$j$}
  \put(305,31)     {$\Otic$}
  \put(324.2,-8.1) {\sse$\ol\imath$}
  \put(339.1,74.4) {\sse$\ol k$}
  \put(343.2,38.3) {\tiny$\ol\gamma$}
  \put(353.2,-8.1) {\sse$\ol\jmath$}
  \epicture07 \labl{pic65}
together with \Lemma \ref{lem:dual-basis}.  \\[.1em]
To show that $T_\cG$ extends to a \ssFA, by \Remark \ref{prop:ssFA-unique}(iv) 
it is sufficient to verify that the morphism $\Phi_{1,\natural}$, which was 
defined after \erf{eq:epsnat}, is invertible. 
  % It is also easily checked that $T_\cG$ has trivial twist, in agreement
  % with proposition} \ref{c+s=tt}.    
Now for every $i\iN\IG$ we have
  %% [pic~66]
  \bea \begin{picture}(370,90)(0,16)
  \put(0,0)  {\begin{picture}(0,0)(0,0)
              \scalebox{.38}{\includegraphics{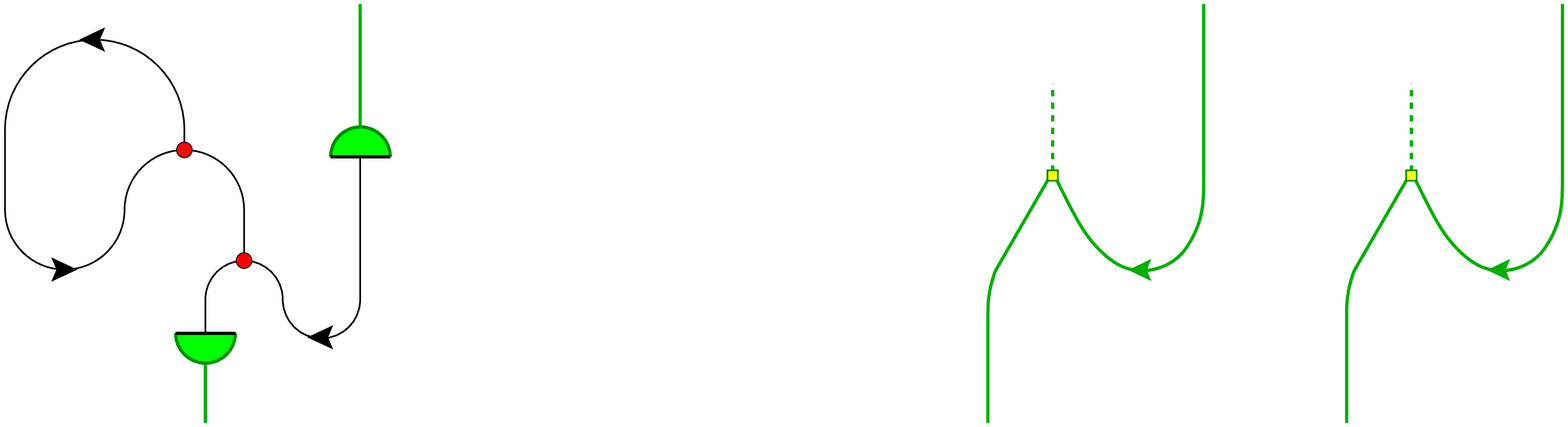}} \end{picture}}
  \put(-12,58)     {\sse$T_{}^\vee$}
  \put(33,-8.5)    {\sse$U_i\ti\ol U_i$}
  \put(42.7,71)    {\sse$T$}
  \put(57,45)      {\sse$T$}
  \put(65.3,101.5) {\sse$U_i^\vee\ti\ol U_i^\vee$}
  \put(84,45)      {\sse$T_{}^\vee$}
  \put(118,43)     {$=\ \dsty\sum_{p\in\IG}\;(\dim(U_p))^2_{}$}
  \put(223,-8.5)   {\sse$U_i$}
  \put(271,101.5)  {\sse$U_i^\vee$}
  \put(286,43)     {$\Otic$}
  \put(305,-8.5)   {\sse$\ol U_i$}
  \put(353,101.5)  {\sse$\ol U_i^\vee$}
  \epicture04 \labl{eq:TG-phinat}
because only the tensor unit of $\cGG$ contributes
in the $T_\cG$-ribbon that is connected to the $T_\cG$-loop
and the resulting isolated $T_\cG$-loop amounts to a factor $\dim(T_\cG)$.
Substituting the definition of $m$ then gives the \rhs\ of \erf{eq:TG-phinat}.
Since the morphism on the \rhs\ is
invertible for every $i\iN\IG$, so is $\Phi_{1,\natural}$.
\qed

\dtl{Lemma}{lem:TG'}
With $T_\cG$ defined by \erf{TG}, we have:
\\[.3em]
(i)~\,The induced $T_\cG$-modules
  \be
  M_k := \Ind_{T_\cG}( \one \ti \ol U_k ) \ee
($k\iN\IG$) are mutually distinct and simple.
\\[.3em]
(ii)~The induced modules $\Ind_{T_\cG}( U_k{\times}\ol U_l)$
decompose into a direct sum of simple $T_\cG$-modules according to
  \be
  \Ind_{T_\cG}( U_k{\times}\ol U_l)
  \,\cong\, \bigoplus_{r \in \IG} \N krl\, M_r \,,  \labl{eq:ind-decomp} 
with $\N ijk$ the dimension of $\Hom(U_i\Oti U_j,U_k)$, as introduced in
\erf{Nijk}. 

\medskip\noindent
Proof:
\\[.3em]
(i)~\,Since $\cG$ is semisimple, $\cGG$ is semisimple as well, and 
hence the object $\M_k$ underlying induced module $M_k$ is a direct sum of 
simple objects of $\cGG$. The decomposition into simple objects reads
  \be
  \M_k = T_\cG \otimes (\one\ti\ol U_k) \,\cong\,
  \bigoplus_{r,s\in\IG} \N rks \, U_r \ti \ol U_{\!s} \,, \labl{eq:Mk-decomp}
with $\N ijk\eq \dim\Hom(U_i\Oti U_j,U_k)$.
When combined with the reciprocity relation \erf{reciUM}, this implies
  \be \bearll
  \Hom_{T_\cG}(M_k,M_l) \!\!
  & \cong \bigoplus_{r,s\in\IG} \Hom^\cG(U_r\oti U_l,U_s)
  \oti \Hom^{\cGG}(\one \ti\ol U_k,U_r\ti\ol U_{\!s})
  \\{}\\[-.7em]
  & \cong \Hom^\cG(U_l,U_k) \,\cong\, \delta_{k,l}\,\koerper \,, \eear \ee
which proves the claim.
\\[.3em]
(ii)~We first check that the simple modules $M_r$ appear
in $\Ind_{T_\cG}( U_k{\times}\ol U_l )$ with
multiplicity $\N krl$. To this end we we use again reciprocity:
  \be
  \Hom^{\cGG}_{T_\cG}(M_r,\Ind_{T_\cG}( U_r{\times}\ol U_l ))
  \,\cong\, \Hom^{\cGG}_{}(\M_r, U_k{\times}\ol U_l ) \,\cong\,
  \koerper_{}^{\N krl} \,. \ee
The last equality follows from the decomposition of $\M_r$ into
simple objects given in \erf{eq:Mk-decomp}.
\\[.2em]
We now know that the \rhs\ of \erf{eq:ind-decomp} is a submodule of 
$\Ind_{T_\cG}( U_k{\times}\ol U_l )$. Next we check that $\Ind_{T_\cG}( 
U_k{\times}\ol U_l )$ does not contain any further submodules. It is 
sufficient to verify that \erf{eq:ind-decomp} is correct as a relation 
for objects in $\cGG$.  For the two sides of \erf{eq:ind-decomp} we find
  \bea
  \Ind_{T_\cG}( U_k{\times}\ol U_l )
  \,\cong\, \dsty \bigoplus_{r,u,v\in\IG} \N rku \N rlv \; U_u \ti \ol U_v
  \qquad {\rm and} \\{}\\[-.8em] \dsty
  \bigoplus_{r\in\IG} \N krl M_r
  \,\cong\, \bigoplus_{r,u,v\in\IG} \N krl \N urv \; U_u \ti \ol U_v \,,
  \eear\ee
respectively. Using the identities $\N rku\eq\N u{\bar k}r$ and
$\N krl\eq\N{\bar k}lr$, we see that the two expressions
coincide owing to associativity of the tensor product.
\qed

%%%%%%%%%%%%%%%%%%%%%%%%%%%%%%%%%%%%%%%%%%%

\subsection{Modularity implies \platy} \label{tococo}

We will now apply some of the results above in the particular case that 
the tensor category under consideration is even modular. We are going to
show that such categories are \platl, with the compensating \cat\ given by the 
dual and the trivialising \alg\ of the form given in \Lemma \ref{lem:TG}.

\medskip

In this subsection $\cG$ always denotes a modular \tc.
As a preparation we need
\\[-2.3em]

\dtl{Lemma}{lem:s-inv}
(i)~\,Let $U_k$ be a simple object in a modular tensor category $\cC$.
If the relation $\theta_s / (\theta_k \theta_r)\eq1$ holds
for all simple objects $U_r,\,U_s$ ($r,s\iN\II$) such that 
$\N rks \,{\ne}\, 0$, then $U_k\eq\one$.
\\[.3em]
(ii)~Conversely, let $\cC$ be a semisimple additive ribbon category with
ground field $\koerper$ and with finite index set $\II$. If the equality
$\theta_s/ \theta_k \theta_r\eq 1$ for all $r,s\iN\II$ such that
$\N rks \,{\ne}\, 0$ implies that $k\eq 0$, then $\cC$ is modular.

\medskip\noindent
Proof:\\
(i)~\,Fix a basis $\{\lambda_{kr,\alpha}^{\,s}\}\,{\subset}\,
\Hom(U_k{\otimes}U_r,U_s)$. Then one has
  \be
  \lambda_{kr,\alpha}^{\,s} \circ c_{r,k} \circ c_{k,r}
  = \frac{\theta_s}{\theta_k \theta_r}\, \lambda_{kr,\alpha}^{\,s} \ee
(see e.g.\ \Section 2.2 
of \cite{fuRs4} for more details). By assumption, all
the factors $\theta_s/(\theta_k \theta_r)$ in this expression
are equal to one. Since $s$ and $\alpha$ run over a basis, this implies that
  \be
  c_{r,k} \circ c_{k,r} = \id_{U_k{\otimes}U_r} \ee
for all $r\iN\II$.
Taking the trace of this formula yields $s_{r,k}\eq s_{k,\One} s_{r,\One}$.
Thus the $k$th column of the $s$-matrix \erf{sij} is proportional to the
$\one$-column, with a factor of proportionality equal to $s_{k,0}$.
Since the $s$-matrix is invertible, this is only possible if $k\eq\One$.
\\[.3em]
(ii)~The same calculations show that the conditions are equivalent to the
statement that the equality
$c_{U_r,U_k}c_{U_k,U_r}\eq\id_{U_k\otimes U_r} $ for all $r\iN\II$ implies that
$k\eq0$. Taking the trace, we learn that $k\eq0$ is the only element of $\II$
such that $s_{U_r,U_k}\eq\dim (U_k) \dim (U_r)$ for all $r\iN\II$. According to
\Proposition 1.1 of \cite{brug2}, this property in turn implies that the 
ribbon category $\cC$ is modular.
\qed 

\dtl{Lemma}{lem:TG-local}
For $\cG$ a modular tensor category and $T_\cG$ as defined in lemma
\ref{lem:TG}, up to isomorphism
the only local simple $T_\cG$-module is $M_\one\eq T_\cG$ itself.

\medskip\noindent
Proof:\\
By corollary \ref{cor:local} it is enough to compute the twist on the simple 
modules $M_k$ and check whether it is of the form $\xi_k \id_{M_k}$ for 
some $\xi_k \iN \koerper$.  Since $\one \ti\ol U_k$ is always a subobject of 
$M_k$, if it exists $\xi_k$ must be equal to $\theta_k^{-1}$. Evaluating the 
twist for all other subobjects of $M_k$ we find the following condition: $M_k$ 
is local iff $\theta_r \theta_s^{-1}\eq\theta_k^{-1}$ for all $r,s$ such that 
$\N rks \,{\ne}\, 0$. By \Lemma \ref{lem:s-inv} this implies that $k\eq\One$.
\qed

\medskip

\dtl{Proposition}{thm:top}
For $\cG$ a modular tensor category and $T_\cG$ as defined in lemma
\ref{lem:TG}, there is an equivalence
  \be 
  \EXt{\cG\Ti\ol\cG}{T_\cG} \;{\cong}\; \Vectk  \ee
of modular tensor categories.

\medskip\noindent
Proof:\\
Combining the \Lemmata \ref{lem:TG'}\,--\,\ref{lem:TG-local} above, we
conclude that $\EXt{\cG\Ti\ol\cG}{T_\cG}$ is a modular \tc\ that, up to
isomorphism, has the tensor unit $\one$ as its single simple object.
Any such \cat\ is equivalent to $\Vectk$.

%%%%%%%%%%%%%%%%%%%%%%%%%%%%%%%%%%%%%%%%%%%%%%%%%%%%%%%%%%%%%%%%%%%%%%%%
\newpage

\sect{Correspondences of tensor categories}\label{sec6}

\subsection{Ribbon categories} \label{s6s1} 

We are now finally in a position to establish correspondences between certain 
ribbon categories $\cQ$ and $\cG$. They make use of another ribbon category 
$\cH$, which must be \platl. The strongest result, to be derived in 
\Section \ref{s6s2}, is obtained when $\cH$ is even a modular \tc. In the 
present 
subsection, this special property of $\cH$ is not required. Also, $\cQ$ and 
$\cH$ are not assumed to be Karoubian. Given $\cQ$ and $\cH$, we consider a 
ribbon category $\cG$ that is obtained as the \cat\ of local modules over a 
suitable \alg\ $\A$ in the Karoubian product of $\cQ$ and $\cH$. 

\dtl{Proposition}{prop:coset}
Let $\cQ$ be a ribbon category, $\cH$ a \platl\ ribbon category, with
\platz\ data $\Hp$ and $T$,
and let $\A$ be a haploid commutative \ssFA{} in the \cat\ $\QH$ satisfying 
$\,\dim_\koerper\Hom(\one_\cQ{\times}T,L{\times}\one_{\Hp}) \eq 1$.
Denote by $\cG$ the ribbon category of local $\A$-modules,
  \be
  \cG:= \EXt\QH{\!\A} \,.  \labl{defG}
Further, let $\B$ be the object
  \be  \B := \lxt\OT{\AO\!}  \ee
in $\GHp$, endowed with the structure of Frobenius \alg\ in $\GHp$ via the
prescription given in the proof of \Proposition \ref{lem:[B]A-lift-i};
similarly, let $\Gama$ be the Frobenius \alg\
  \be  \Gama := \lxt\AO\OT  \ee
in $\EXt\QHH\OT$. We have
  \be  \kar\cQ \cong \EXt\QHH\OT  \,.  \labl{QHHOT}
Furthermore, if $\B$ and $\Gama$ have non-zero dimension, then they are 
haploid commutative \ssFA s, and there is an equivalence
  \be
  \EXt{\kar\cQ}\Gama \,\cong\, \EXt{\GHp}\B  \labl{QCB}
of (Karoubian) ribbon \cats.  

\bigskip\noindent
Proof:\\
(i)\,\,~To verify the equivalence \erf{QHHOT}, we first apply 
           \Lemma \ref{prodvect}, then the fact that, by assumption, 
$\Hp$ and $T$ provide a \platz\ for $\cH$, and then corollary \ref{CAxD0}:
  \be
  \kar\cQ \,\cong\, \cQ \,\Ti\, \Vectk \,\cong\, \cQ \,\Ti\, \EXt{\HH}{T}
  \,\cong\, \EXt\QHH\OT \,.  \labl{QQV}
\vskip.3em\noindent
(ii)~\,That $\B$ and $\Gama$ are haploid commutative symmetric special 
Frobenius algebras can be seen by combining \Proposition \ref{prop:AB-alg} 
  % and \ref{lem:[B]A-lift-i} 
and corollary \ref{deKac} as well as \Proposition \ref{lem:[B]A-lift-i}(ii).
Note in particular that we can apply \Proposition \ref{prop:AB-alg}(iii), 
because both $L$ and $T$ are symmetric and special, 
the dimensions of $\B$ and $\Gama$ are non-vanishing, and the condition on 
the centers is implied by $\dim_\koerper\Hom(\one_\cQ{\times}T, L{\times}
\one_{\Hp}) \eq 1$ together with the commutativity of $L$ and $T$.
\\[.3em]
(iii)~For the next two preparatory calculations, we invoke successively 
       \Proposition \ref{thm:[B]A-lift-ii}, 
corollary \ref{CAxD0} and the definition \erf{defG} 
of $\cG$ (as well as the associativity of the Karoubian product 
$\boxtimes$ from \Remark \ref{boxbox}) to write
  \be \bearll
  \EXt\QHH{\!\Efu\OT\AO} \!\!&
  \cong \EXT{ \EXt{(\QH)\Ti\Hp}{\AO} }{\tildeOT}
  \\{}\\[-.6em]&
  \cong \EXT{ \EXt{\QH}{\A} \,\Ti\, \Hp }{\tildeOT}
  \\{}\\[-.6em]&
  \cong \EXt{ \GHp }{\tildeOT}  \eear \labl{QHH-GHp}
and similarly, using \erf{QQV} in the second step,
  \be \bearll
  \EXt\QHH{\!\Efu\AO\OT} \!\!&
  \cong \EXT{ \EXt{\cQ\Ti(\HH)}{\OT} }{\tildeAO}
  \\{}\\[-.6em]&
  \cong \EXt{ \kar\cQ }{\tildeAO} \,.  \eear \labl{QHH-Q}
         % also used \erf{516b} 
(Recall from \Lemma \ref{lem:CAK-CKA}(i) that the \cat\ of local
modules over any commutative \ssFA\ in a Karoubian ribbon \cat\ is
again Karoubian. Thus all the module \cats\ appearing here are Karoubian.)
\\[.3em]
(iv)~Consider now the tensor product \alg\ 
  \be  \TA := (\OT) \oti (\AO)  \labl{TA}
in $\QHH$. Recall that in a braided setting the tensor product of two commutative
\alg s is not commutative, in general. Concretely, applying \Proposition 
\ref{prop:tensor-center} we learn that the left and right centers of $\TA$ are
  \be
  C_l(\TA) \cong \Efu\AO\OT  \qquad {\rm and} \qquad
  C_r(\TA) \cong \Efu\OT{\AO\!}  \,,  \labl{ClrTA}
respectively. Further, by \Theorem \ref{thm:equiv} 
the \cats\ of local $C_l(\TA)$- and local $C_r(\TA)$-modules are equivalent,
  \be  \EXt\QHH{C_l(\TA)} \cong \EXt\QHH{C_r(\TA)}
  \,.  \labl{QHHlr}
Combining this information with 
the results in step (iii) and \erf{ClrTA}, we finally obtain
  \be  \EXt{ \GHp }{\tildeOT} \cong \EXt{ \kar\cQ }{\tildeAO} \,,
  \labl{QHHlr1}
thus establishing the equivalence \erf{QCB}. This is a ribbon equivalence 
because all the intermediate equivalences
     we used are ribbon.
\qed

%%%%%%%%%%%%%%%%%%%%%%%%%%%%%%%%%%%%%%%%%%%%%%%%%%%%%%%%%%%%%%%%%%%%%%%%

\subsection{Modular tensor categories} \label{s6s2}

It is desirable to find also a description of the \cat\ $\kar\cQ$ itself, not
just of some module \cat\ over $\kar\cQ$, in terms of $\cG$ and $\Hp$. As it 
turns out, this can be achieved if we assume that $\cH$ is {\em modular\/}
such that it has a trivialisation of the form described in 
\Proposition \ref{thm:top}, i.e.\
  \be  \Hp = \ol\cH \qquad{\rm and}\qquad  T = T_\cH  \ee
with $T_\cH$ as given in \Lemma \ref{lem:TG}. In addition, also one further 
condition on the \alg\ $\A$ and one further condition on the category
$\cQ$ must be imposed; these properties are the following.
\\[-2.3em]

\dtl{Definition}{Cnice}
An algebra $A$ in the Karoubian product $\cC\Ti\cD$ of two tensor \cats\ 
$\cC$ and $\cD$ is called {\em $\cC$\haploid\/} iff
  \be  \Obj(\cC\Ti\cD) \ni\, U\ti\one_\cD \In A \
  \,\Rightarrow\,\  U \cong \one_\cC \,, \labl{U=1-prop}
i.e.\ iff up to isomorphism the only retract
of $A$ of the form $U{\times}\one_\cD\,$ is $\one_\cC{\times}\one_\cD$.  

\dtl{Definition}{brimfull}
A sovereign tensor category $\cC$ is called {\em separable\/}
if every idempotent $p$ with $\tr(p)\eq0$ is the zero morphism.

\dtl{Remark}{rem:C-hapl}
(i)~\,It follows from \Remark \ref{prop:ssFA-unique}(vi) 
that if $\dim_\koerper \Hom(\one,A)\eq d$ for a Frobenius algebra $A$ in 
$\cC\Ti\cD$, then $I^{(d)}_\cC\ti\one_\cD$ with $I^{(d)}_\cC\eq 
\one_\cC{\oplus}\one_\cC{\oplus}\cdots{\oplus}\one_\cC$ ($d$ summands) 
is a retract of $A$, and hence in particular $A$ is not $\cC$\haploid. 
Conversely, if $A$ is $\cC$\haploid, then it is in particular haploid.
\\[.2em]
Also, when $\cC\,{\cong}\,\Vectk$, for Frobenius algebras the notions 
of haploidity in $\cD$ and of $\cC$-haploidity coincide upon identifying 
$\cC\Ti\cD$ with $\cD$. This is the reason for the choice of terminology.
\\[.3em]
(ii)~Since every idempotent in the Karoubian envelope $\kar\cC$ of a
sovereign tensor category $\cC$ is also an idempotent in $\cC$, 
separability of $\cC$ implies separability of $\kar\cC$; owing to the 
functorial embedding $\cC\,{\to}\,\kar\cC$, the converse holds true, too.
Also, if $\cC$ is separable, then so is its dual $\ol\cC$.
\\
If $\cC$ and $\cD$ are sovereign tensor categories such that their 
product $\cC\Tic\cD$ (or $\cC\Ti\cD$) is separable, then  already 
$\cC$ and $\cD$ are separable.
\\
Furthermore, since, for $A$ an algebra in a sovereign tensor category 
$\cC$, every idempotent in $\cC_\AA$  is also an idempotent in $\cC$, 
separability of $\cC$ implies separability of $\cC_\AA$. By the same
argument, the category \calcal\ of local modules over a commutative \ssFA\ 
$A$ in a separable ribbon category $\cC$ is separable.
\\
Modular categories are in particular separable.

\bigskip

The proof of the stronger result involving modular tensor categories
relies also on the following 
\\[-2.3em]

\dtl{Lemma}{PPlemma}
Let $S,S'$ be two retracts of an object $U$ in a (not necessarily Karoubian)
separable sovereign tensor category $\cC$.
Suppose that the corresponding split idempotents satisfy
$P_SP_{S'}\eq P_{S'}P_S$ and $\tr_U( P_S ) \eq\tr_U( P_S P_{S'})\eq
\tr_U(P_{S'})$. Then $P_S \eq P_{S'}$ and $S \,{\cong}\, S'$ as retracts.

\medskip\noindent
Proof:\\
We write $S\eq (S,e,r)$ and $S'\eq (S',e',r')$, and consider the
morphisms $f\iN\Hom(S,S')$ and $g\iN\Hom(S',S)$ given by $f\,{:=}\,r'\cir e$ 
and $g\,{:=}\,r\cir e'$. Using the assumptions we see that $p\,{:=}\, 
g \cir f$ satisfies $p\cir p \eq r \cir P_{S'} \cir P_S \cir P_{S'} \cir e
\eq r\cir P_{S'} \cir e \eq p$, i.e.\ $p$ is an idempotent. Further we have
  \be
  \tr^{}_S\, p = \tr^{}_U( P_S P_{S'} ) = \tr^{}_U P_S = \dim(S)\,.  \ee
It follows that $\tr_S(\id_S{-}p)\eq 0$. By separability this implies that
$\id_S{-}p\eq 0$ so that $p\eq\id_S$. In the same way one shows that 
$f\cir g \eq\id_{S'}$. Thus $S$ and $S'$ are isomorphic as objects.
\\
{}From $\id_S\eq g\cir f \eq r \cir P_{S'} \cir e$ we deduce
(composing with $e$ from the left) that $e \eq P_{S'} \cir e \eq
    $\linebreak[0]$%
e' \cir f$ and (composing with $r$ from the right) that
$r \eq r \cir P_{S'} \eq g \cir r'$. The relation $e\eq e'\cir f$
implies that $S$ and $S'$ are isomorphic as subobjects, and
$P_S\eq e\cir r\eq e' \cir f \cir g \cir r'\eq P_{S'}$
shows that they are isomorphic even as retracts.
\qed

\bigskip

Having these ingredients at hand,%
  \foodnode{Recall also \convention s \ref{basicprops} and \ref{c:csplit}.}
we can formulate a much stronger result than the one of \Proposition 
\ref{prop:coset}: 
\\[-2.2em]

\dtl{Theorem}{thm:coset}
Let $\cQ$ be a (not necessarily Karoubian) ribbon category and $\cH$ 
a modular \tc\ (with \platz\ data $\ol\cH$, $T\,{\equiv}\,T_\cH$) such that 
the product $\QHh$ is separable, and let
$\A$ be a $\cQ$\haploid\ commutative \ssFA{} in the Karoubian product $\QH$.
\\[.2em]
(i)~\,The Frobenius algebra
  \be  L' := \lxt\OT{\Ao\!}  \ee
is haploid, commutative, symmetric and special, and there is an equivalence
  \be
  \kar\cQ \cong \EXt{\GHb}{L'}  \labl{eq:mod-coset-1}
of ribbon categories, with $\cG\eq\EXt\QH{\!\A}$.
\\[.4em]
(ii)~The Frobenius algebra $L'$ in $\GHb$ is even $\cG$\haploid.

\bigskip\noindent
Proof of (i):\\
1)~\,We start by checking that the conditions of \Proposition \ref{prop:coset} 
are fulfilled. Note that
  \be \bearll
  \dim_\koerper\,\Hom(\one_{\cQ}{\times}T_{\cH},L{\times}\one_{\Hb}) \!\!
  & = \dsty\sum_{k \in \II_{\cH}}
  \dim_\koerper\,\Hom(\one_{\cQ}{\times}U_k{\times}{\ol{U_k}},L{\times}
    \one_{\Hb}) \\{}\\[-.7em]
  & = \dim_\koerper\,\Hom(\one_{\cQ}{\times}\one_{\cH},L) = 1 \,,
  \eear\ee
since $L$ is in particular haploid, by \Remark \ref{rem:C-hapl}(i). Next we 
need to show that the algebras $\B\iN\Obj(\GHb)$ and $\Gama\iN\Obj(\kar\cQ)$ 
appearing in \Proposition \ref{prop:coset} 
have non-zero dimension. To see this, note that according to \Remark 
\ref{rem:C-hapl}(ii)
the categories $\kar\cQ$ and $\GHb$ are separable. Hence for any object $U$ in
one of these categories, the vanishing of $\dim(U)$ 
implies that $\tr(\id_U)\eq 0$ and thus $\id_U\eq 0$, so that $U$ is a zero object.
On the other hand, by \Remark \ref{prop:ssFA-unique}(vi), any Frobenius 
algebra has the tensor unit as a retract, and hence cannot be a zero object. 
\\
We can therefore apply \Proposition \ref{prop:coset}; 
in particular $L'\eq\B$ is haploid, commutative, symmetric and special.
To establish \erf{eq:mod-coset-1}, it remains to 
be shown that $\Gama \eq \lxt\AO\OT$ is trivial, $\Gama\,{\cong}\,\one_{\cQ}$.
\\[.3em]
2)~\,We regard $\QHc$ as a sub\cat\ of $\QH\eq\kar{(\QHc)}$ in the 
usual manner, and likewise for $\GHbc$. We start by noticing that the two 
\alg s $\Efu\Ao\OT \,{\cong}\,C_l(\TA)$ and $\OT$ are both retracts of 
$\TA\,{:=}\,(\OT)\oti(\Ao)$. The associated idempotents are
  %% [pic~32]
  \bea  \begin{picture}(280,101)(0,30)
  \put(0,0)   {\begin{picture}(0,0)(0,0)
              \scalebox{.38}{\includegraphics{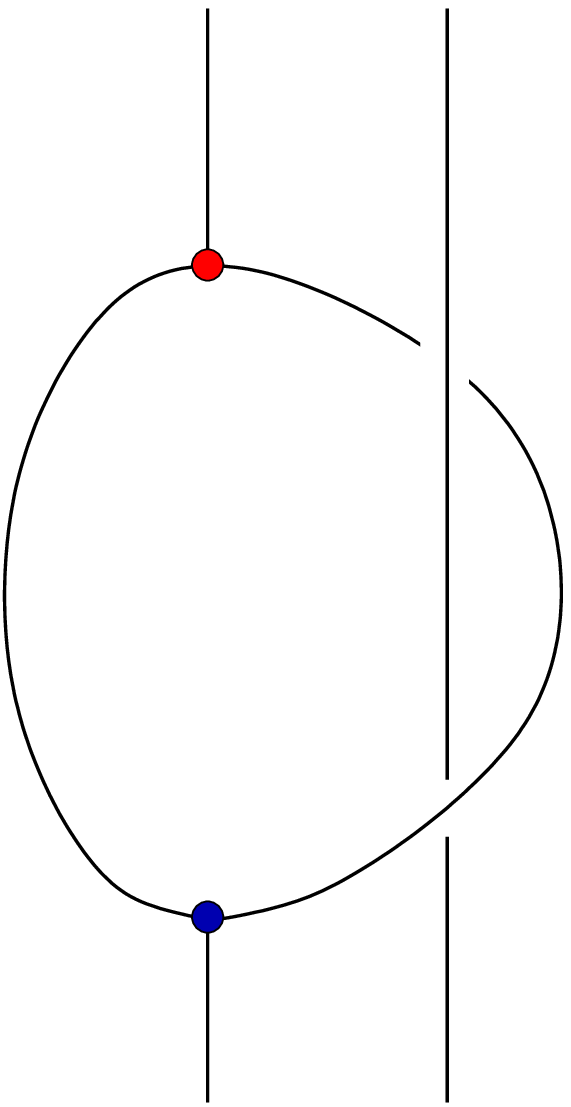}} \end{picture}}
  \put(252,0) {\begin{picture}(0,0)(0,0)
              \scalebox{.38}{\includegraphics{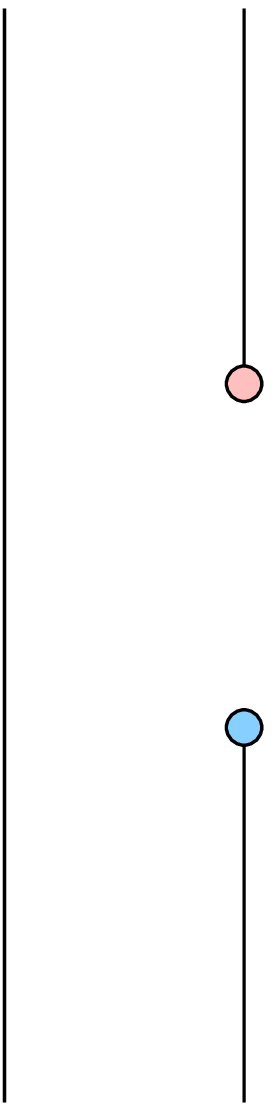}} \end{picture}}
  \put(-64,54.1)   {$P^{}_{\!C_l(\TA)}\; =$}
  \put(5.8,-9.9)   {\sse$\OT$}
  \put(6.8,125.5)  {\sse$\OT$}
  \put(42.5,-9.9)  {\sse$\Ao$}
  \put(43.3,125.5) {\sse$\Ao$}
  \put(98,54.1)    {and}
  \put(143,54.1)   {$P^{}_{\!\OT}\; =\ \dsty\frac1{\dim(\A)}$}
  \put(237.7,-9.9) {\sse$\OT$}
  \put(238.7,125.5){\sse$\OT$}
  \put(272.4,-9.9) {\sse$\Ao$}
  \put(272.9,125.5){\sse$\Ao$}
  \epicture21 \labl{pic32}
respectively. The idempotent $P^{}_{\!C_l(\TA)}$ is split by \convention\ 
\ref{c:csplit}. To see that $P^{}_{\!\OT}$ is split as well, consider $\OT$
as a retract of $F$, with embedding and restriction morphisms
$e\eq\id_{\OT}{\otimes}\eta_{\Ao}$ and
    % $\linebreak[0]$%
$r\eq\id_{\OT}{\otimes}\eps_{\Ao}{/}\dim(L)$, where in the definition of 
$e$ and $r$ the isomorphism $\OT\,{\cong}\,(\OT) \oti
(\one_{\cQ}{\times}\one_{\cH}{\times}\one_{\overline\cH})$ is implicit; 
clearly, $e \cir r\eq \id_{\OT}$ and $r \cir e\eq P^{}_{\!\OT}$.
\\[.2em]
Using the specialness of the \alg\ $T$, one
easily verifies that the idempotents \erf{pic32} satisfy
  \be  P^{}_{\!C_l(\TA)} \circ P^{}_{\!\OT}
  = P^{}_{\!\OT} = P^{}_{\!\OT} \circ P^{}_{\!C_l(\TA)} \,.  \ee
Their traces are computed as $\,\tr(P^{}_{\!\OT})\eq \dim(T)$ and as
  \be
  \tr(P^{}_{\!C_l(\TA)}) = s_{\OT,\Ao}^{\QHh} =
  \sum_{k\in\IH} s_{\one\times U_k,\A}^{\QH}\, s_{\ol U_k,\one}^{\ol\cH}
  \,, \labl{eq:trP1-0}
respectively, where the first equality holds by 
\Remark \ref{lem:[U]A-module-iii}, 
while in the second equality the explicit form \erf{TG} of $T$ is inserted.
\\[.3em]
3)~\,Next we use the fact that $\cH$ is modular and thus in particular
semisimple. Hence writing $\A\iN\Obj(\QH)$ as $\A\eq
(\A_\cQ{\times}\A_\cH;\pi)$ with suitable objects $\A_\cQ$ of $\cQ$ and 
$\A_\cH$ of $\cH$ and an idempotent $\pi\iN\End(\A_\cQ{\times}\A_\cH)$,
we know that $\A_\cH$ is a direct sum of simple objects $U_j$ of $\cH$,
with $j$ in the finite index set $\IH$, and as a consequence
  \be  \A \cong \bigoplus_{j\in\IH} \A_j\ti U_j  \ee
with suitable 
  % (not necessarily simple)
objects $L_j$ of $\cQ$. 
Inserting this decomposition into formula \erf{eq:trP1-0} we obtain
  \be
  \tr(P^{}_{\!C_l(\TA)}) = \sum_{j,k\in\IH} s_{\one,\A_j}^{\cQ}\,
  s_{U_k,U_j}^{\cH}\, s_{\ol U_k,\one}^{\ol\cH}
  = \sum_{j\in\IH} s_{\one,\A_j}^{\cQ} \sum_{k\in\IH}
  s_{U_k,U_j}^{\cH}\,s_{U_k,\one}^{\cH} \,.  \labl{eq:trP1-1}
By the identity \erf{ss=c}, modularity of $\cH$ also implies that the 
$k$-summation in the expression on the \rhs\ can be carried out, yielding  
$\delta_{j,\One}\sum_{k\in\IH}(s_{U_k,\one}^{\cH})^2_{} \eq \delta_{j,\One}
\dim(T)$, and hence $\tr(P^{}_{\!C_l(\TA)})\eq\dim(T)\,s_{\one,L_0}^{\cQ}$. 
Further, the hypothesis that $\A$ is $\cQ$\haploid\
means that $\A_\One\,{\cong}\,\oneQ$; thus we finally get
  \be
  \tr(P^{}_{\!C_l(\TA)}) = \dim(T)\,s_{\one,\one}^{\cQ} = \dim(T) \,.  
  \labl{eq:trP1}
It follows that $\tr(P^{}_{\!C_l(\TA)})\eq \tr(P^{}_{\!C_l(\TA)} 
{\circ}\, P^{}_{\!\OT})\eq \tr(P^{}_{\!\OT})$. By \Lemma \ref{PPlemma}
this implies, in turn, that the two idempotents \erf{pic32}
coincide, $P^{}_{\!C_l(\TA)}\eq P^{}_{\!\OT}$. We conclude that
  \be  \Efu\Ao\OT \cong \OT  \labl{TA-T}
as retracts of $\TA$.
\\[.2em]
It is also not difficult to check that the multiplication induced on
$\OT$ via its embedding in the \alg\ $\TA$ agrees with the one defined in 
\Lemma \ref{lem:TG}. The same holds for $\Efu\Ao\OT\,{\cong}\,C_l(\TA)$, as 
follows from \Proposition \ref{prop:tensor-center}. The isomorphism 
\erf{TA-T} therefore also holds as an isomorphism of algebras, and in 
fact even as an isomorphism of \ssFA s.
\\[.2em]
But the object $\one_{\cQ}{\times}T$ is the tensor unit in the \cat\
$\EXt{\QHh}{\OT} \,{\cong}\, \kar\cQ$, implying that $\lxt\AO\OT \,{\cong}\, 
\one$ as an object in $\kar\cQ$. The relation \erf{eq:mod-coset-1} 
now follows from \erf{QCB} with $L'\eq \B\eq \lxt\OT{\AO\!}$.
\\[.3em]
Proof of (ii):\\
It remains to be shown that the \alg\ $L'$ in $\GHb$ is 
$\cG$\haploid. We will establish that any object $M$ of $\cG$ with the
property that $M{\times}\one_{\cHb}$ is a retract
of $L'$, is itself a retract of $\one_{\cG}$. Since $\one_{\cG}$ is
simple, this implies that $M\,{\cong}\,\one_{\cG}$, and hence (ii).
\\[.3em]
Let us formulate these statements in terms of the category $\QHh$.
$L'$ is the algebra $\Efu\OT\Ao$, while $M$ is a local $L$-module in $\QH$. 
That $(M{\times}\one_{\Hb},e,r)$ is a retract of $L'$ in $\GHb$ thus means that
  \bea
  e \in \Hom_\Ao( M{\times}\oneHb , \lxt\OT\Ao )
  \qquad {\rm and} \\{}\\[-.7em]
  r \in \Hom_\Ao( \lxt\OT\Ao , M{\times}\oneHb )  \eear \labl{eq:M-retr}
as morphisms of $\QHh$. Now by the isomorphisms of \Proposition 
\ref{lem:[U]A-module-ii} and the reciprocity relation \erf{reciMV}, we have
  \be
  \Hom_\Ao(M{\times}\oneHb,\lxt\OT\Ao )
  \cong \Hom(\M{\times}\oneHb,\oneQ{\times}T) \,.  \ee
Using the explicit form of $T$ from formula \erf{TG}, this morphism space in 
$\QHh$ is, in turn, isomorphic to the space $\Hom(\M,\oneQ{\times}\oneH)$ of 
morphisms in $\QH$, and hence to $\Hom_L(M,L)$. Together with a similar 
argument for the second morphism space in \erf{eq:M-retr} we can conclude 
that there are bijections
  \bea
  f :\quad \Hom_\Ao(M{\times}\oneHb,\lxt\OT\Ao ) \stackrel\cong\longrightarrow
        \Hom_L(M,L) 
  \qquad {\rm and} \\{}\\[-.7em]
  g :\quad \Hom_\Ao(\lxt\OT\Ao,M{\times}\oneHb ) \stackrel\cong\longrightarrow
        \Hom_L(L,M) \,.  \eear \ee
Substituting the explicit form of these isomorphisms one can verify that for 
the morphisms $e$ and $r$ of \erf{eq:M-retr} we have $g(r) \cir f(e)\eq\id_M$. 
It follows that $(M,f(e),g(r))$ is a retract of $L$. Moreover, since $f(e)$ 
and $g(r)$ are morphisms of $L$-modules and $L$ is the tensor unit of the 
category $\cG$, this implies that $M$ is a retract of $\one_{\cG}$ in $\cG$.
\qed

\medskip

Combining \Theorem \ref{thm:coset} with \Proposition \ref{thm:mod} 
we arrive at the following statements about the \cat\ $\kar\cQ$:
\\[-2.1em]

\dtl{Corollary}{cor:modular}
For $\cQ$ a (not necessarily Karoubian) ribbon category and $\cH$ a modular 
tensor category such that the product $\QHh$ is separable, 
and $L$ a $\cQ$\haploid\ commutative symmetric special
Frobenius algebra in $\QH$, we have:
\\[.3em]
(i)~$\;$If $\EXt\QH{\!\A}$ is semisimple, then so is $\kar\cQ$.
\\[.3em]
(ii)~If  $\EXt\QH{\!\A}$ is a modular tensor category, then so is $\kar\cQ$.

\bigskip

Theorem \ref{thm:coset} allows us to construct the tensor category
$\cQ$ from the knowledge of the categories $\cG$ and $\cH$
and of the algebra $L'\,{=}\,\lxt\OT\Ao$ in $\GHb$.
For applications, e.g.\ in conformal quantum field theory, it turns out to 
be important to gain information about $L'$ by using as little information
about the category $\cQ$ as possible. The following result helps to
determine $L'$ as an object of $\GHb$ in case that $\cG$ is a {\em modular\/}
\tc\ (and hence, by corollary \ref{cor:modular}(ii), $\kar\cQ$ is a \mtc, 
too), so that in particular the set $\{M_\kappa\,{\mid}\,\kappa\iN\IG\}$ of
isomorphism classes of simple objects in $\cG$
(i.e.\ of simple local $\A$-modules in $\QH$) is finite.
\\[-2.1em]

\dtl{Lemma}{lem:B-obj}
Let $\cQ$, $\cH$ and $\A$ be as in \Theorem \ref{thm:coset}, 
and assume that $\cG\,{:=}\,\EXt{\QH}{\!\A}$ is modular.
Then as an object in $\GHb$ the algebra $L'\,{:=}\,\lxt\OT\Ao$ decomposes as
  \be
   L' \,\cong\, \bigoplus_{\kappa\in\IG} \bigoplus_{l\in\IH}
  \dim\llb\Hom^\QH(\M_\kappa,\one_\cQ\ti U_l) \lrb \, M_\kappa \ti \ol{U_l}
  \,. \ee

\medskip\noindent
Proof:\\
By \Theorem \ref{thm:coset}, $L'$ is a lift to $\GHb \,{\cong}\,\EXt\QHh\Ao$ of 
the algebra $\Efu\OT\Ao$, which is a local $\A{\times}\one_{\ol\cH}$\,-module.
Now owing to relation \erf{eq:CAxD1} every simple local
$\A{\times}\one_{\ol\cH}$\,-module is of the form $M \ti\ol{U_l}$, with $M$ a
simple local $\A$-module and $\ol{U_l}$ a simple object of $\ol\cH$. Invoking 
\Proposition \ref{lem:[U]A-module-ii} and the reciprocity relation 
\erf{reciMV}, it follows that the algebra $L'$ decomposes according to
  \be
  \Efu\OT\Ao \,\cong\,
  \bigoplus_{\kappa\in\IG} \bigoplus_{l\in\IH} \dim\llb\Hom^\QHh
  (\M_\kappa{\times}\ol{U_l}, \oneQ{\times}T)\lrb \, M_\kappa\ti\ol{U_l}  \ee
into simple local $\Ao$\,-modules. 
Moreover, the morphism spaces appearing here obey
  \be \begin{array}{ll}
  \Hom^\QHh (\dot M{\times}\ol{U_l}, \oneQ{\times}T) \!\!
  &\cong \Hom^\QHh (\dot M\ti\ol{U_l}, \oneQ \ti U_l \ti \ol{U_l} )
  \\[-.7em]{}\\
  &\cong \Hom^\QH(\dot M , \oneQ \ti U_l ) \,,
  \eear\ee
where the first isomorphism follows by inserting the explicit form of $T$
from \erf{TG} and observing that only the component $U_l \ti\ol{U_l}$
contributes.
\qed

\dtl{Remark}{rem:Dim-Coset}
If $\cG$, $\cQ$ and $\cH$ are modular, then from the observations in 
\Remarks \ref{rem:Dk-1}(i), \ref{rem1box}(iv) and \ref{rem:duals} 
one can easily determine the dimension of the algebra $L'$. Indeed, because of
$\cG\,{\cong}\,\EXt\QH{\!L}$ and $\cQ\,{\cong}\,\EXt{\GHb}{\!L'}$ we have 
  \be
  p^+(\cG) = \frac{p^+(\cQ)\,p^+(\cH)}{\dim^{\cQ\Ti\cH}(L)}
  \qquad {\rm and} \qquad
  p^+(\cQ) = \frac{p^+(\cG)\,p^-(\cH)}{\dim^{\cG\Ti\cHb}(L')}
  \,.  \ee
As a consequence,
  \be
  \dim^{\cQ\Ti\cH\!}(L)\, \dim^{\cG\Ti\cHb\!}(L')
    = {\rm Dim}(\cH)  \,.  \ee
This expresses the dimension of $L'$ in terms of those of $L$ and $\cH$.

%%%%%%%%%%%%%%%%%%%%%%%%%%%%%%%%%%%%%%%%%%%
\newpage

\appendix

\sect{Graphical calculus}

The computations in this paper are often presented in terms of a 
graphical calculus for ribbon categories, which was first advocated 
in \cite{joSt5}. To make these manipulations more easily accessible, we 
summarise in this appendix our conventions, and in particular recall the 
definition of various specific morphisms that are used in the main text.

\subsection{Morphisms} \label{apptab1}

In the following table we present the graphical notation for general 
morphisms of a tensor \cat, their composition and tensor product,
and for the embedding and restriction morphisms (see equation \erf{e-r}) 
of retracts. Also shown are the structural morphisms of a ribbon category:
the braiding, twist, and left and right dualities (see \Definition 
\ref{def:ribcat}), 
as well as the definition of the (left and right) dual of a general morphism:

\vskip 1.2em

\begin{tabular}{|cc|cc|cc|cc|}
                                                \hline\hline
%%%%%%%%%%%%%%%%%%%%%%%%%%%%%%%%%%%%%%%%%%%%%%%%%%%%%%%%%%%%%%%%%%%%%%
       %%%   general morphisms, their composition and tensor product

\mcll{
~~~$\id_U^{}=$
\begin{picture}(26,44)(0,22)       \apppicture{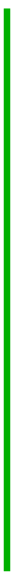}  {7}
\put(4.0,-8.8)   {\sse$U$}
\put(4.5,51.3)   {\sse$U$}
                 \end{picture}     }&
\mcll{
~~~$f       =$      
\begin{picture}(26,38)(0,22)       \apppicture{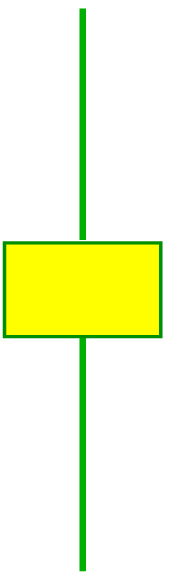}  {7}   
\put(10.3,-8.8)  {\sse$U$}
\put(11.2,51.3)  {\sse$V$}
\put(11.6,22.5)  {\tiny$f$}
                 \end{picture}     }&
\mcll{
~~$g\cir f =$      
\begin{picture}(36,38)(0,22)       \apppicture{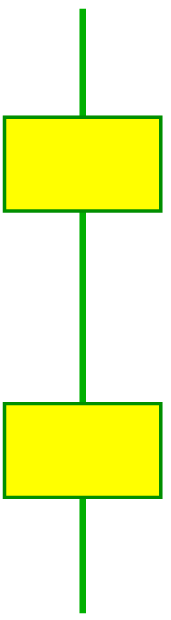}  {7}   
\put(10.1,54.3)  {\sse$W$}
\put(10.3,-8.8)  {\sse$U$}
\put(11.6,12.6)  {\tiny$f$}
\put(12.1,37.5)  {\tiny$g$}
\put(14.9,23.9)  {\tiny$V$}
                 \end{picture}     }&
\mcll{
~$f\Oti f'=$      
\begin{picture}(46,38)(0,22)       \apppicture{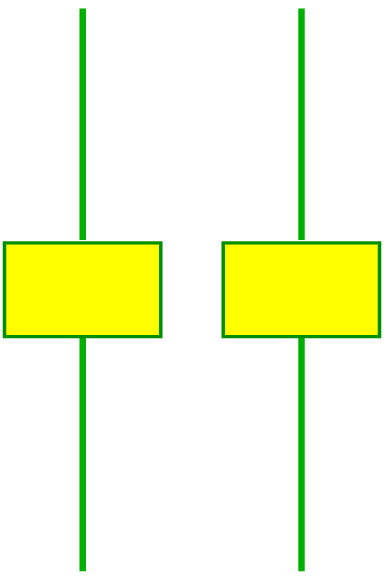}  {7}   
\put(10.3,-8.8)  {\sse$U$}
\put(11.2,51.3)  {\sse$V$}
\put(11.6,22.4)  {\tiny$f$}
\put(28.1,-8.8)  {\sse$U'$}
\put(29.0,51.3)  {\sse$V'$}
\put(29.6,22.4)  {\tiny$f'$}
                 \end{picture}     }
\\
\begin{picture}(0,29){}\end{picture} &&&&&&& \\ \hline

%%%%%%%%%%%%%%%%%%%%%%%%%%%%%%%%%%%%%%%%%%% 
                 %%%   embeddings and restrictions;

\mcll{
~~$e_{S\prec U}^{} =$  
\begin{picture}(42,38)(0,18)       \apppicture{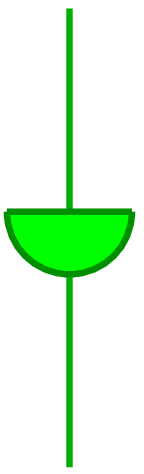}  {7}   
\put(9.2,-8.8)   {\sse$S$}
\put(9.8,43.3)   {\sse$U$}
                 \end{picture}     }&
\mcll{
~~$r_{\!U\succ S}^{} =$  
\begin{picture}(36,30)(0,18)       \apppicture{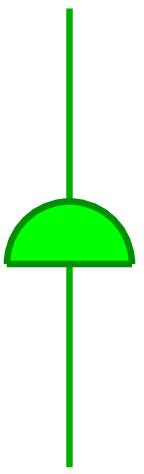}  {7}   
\put(9.2,-8.8)   {\sse$U$}
\put(9.8,43.3)   {\sse$S$}
                 \end{picture}     }
\\
\begin{picture}(0,26){}\end{picture} &&& \\ \hline\hline 

%%%%%%%%%%%%%%%%%%%%%%%%%%%%%%%%%%%%%%%%%%%
                 %%%   braiding, twist

\mcll{
~$c_{U,V}^{} =$
\begin{picture}(36,39)(0,18)       \apppicture{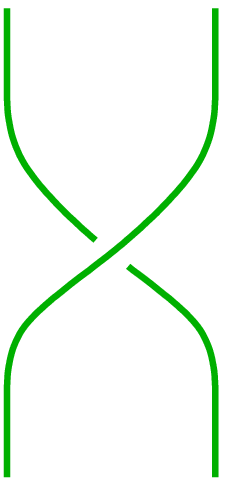}  {7}
\put(4.2,-8.8)   {\sse$U$}
\put(4.8,43.3)   {\sse$V$}
\put(21.2,-8.8)  {\sse$V$}
\put(22.9,43.3)  {\sse$U$}
                 \end{picture}     }&
\mcll{
~$c_{U,V}^{-1} =$
\begin{picture}(36,30)(0,18)       \apppicture{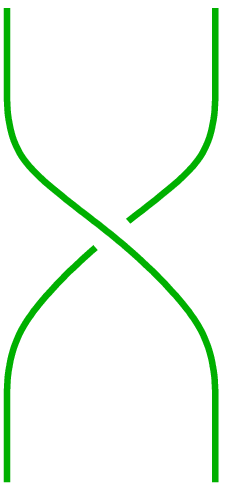}  {7}
\put(4.2,-8.8)   {\sse$V$}
\put(4.8,43.3)   {\sse$U$}
\put(21.2,-8.8)  {\sse$U$}
\put(22.9,43.3)  {\sse$V$}
                 \end{picture}     }&
\mcll{
~~~~$\theta_{U}^{}=$
\begin{picture}(26,30)(0,18)       \apppicture{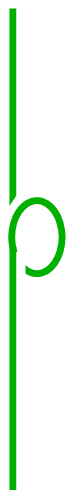}  {7}
\put(5.5,-8.8)   {\sse$U$}
\put(6.4,44.3)   {\sse$U$}
                 \end{picture}     }&
\mcll{
~~~$\theta_{U}^{-1}=$
\begin{picture}(26,30)(0,18)       \apppicture{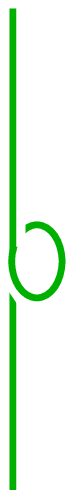}  {7}
\put(5.5,-8.8)   {\sse$U$}
\put(6.4,44.3)   {\sse$U$}
                 \end{picture}     }
\\
\begin{picture}(0,25){}\end{picture} &&&&&&& \\ \cline{1-8}

%%%%%%%%%%%%%%%%%%%%%%%%%%%%%%%%%%%%%%%%%%%
                 %%%   dualities

\mcll{
~$b_{U}^{} =$
\begin{picture}(36,34)(0,18)       \apppicture{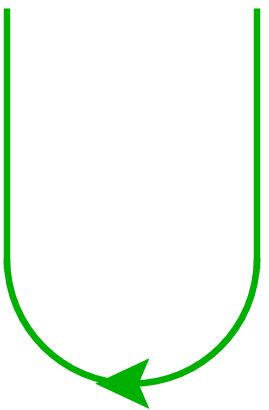}  {7}
\put(4.8,37.3)   {\sse$U$}
\put(24.9,37.3)  {\sse$U^\vee$}
                  \end{picture}    }&
\mcll{
~$d_{U}^{} =$
\begin{picture}(36,30)(0,11)       \apppicture{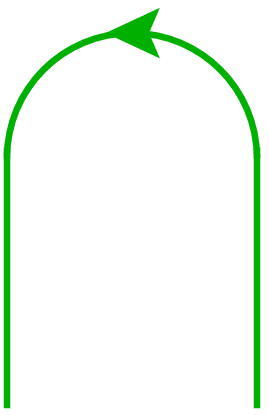}  {7}
\put(3.2,-8.8)   {\sse$U^\vee$}
\put(25.2,-8.8)  {\sse$U$}
                  \end{picture}    }&
\mcll{
~~$\tilde b_{U}^{}=$
\begin{picture}(36,30)(0,18)       \apppicture{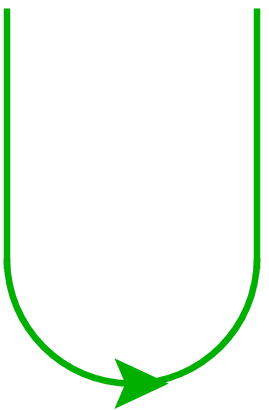}  {7}
\put(-4.4,37.3)  {\sse$
      {\phantom U}^\vee_{}\!{U}$}
\put(25.7,37.3)  {\sse$U$}
                  \end{picture}    }&
\mcll{
~$\tilde d_{U}^{}=$
\begin{picture}(36,30)(0,11)       \apppicture{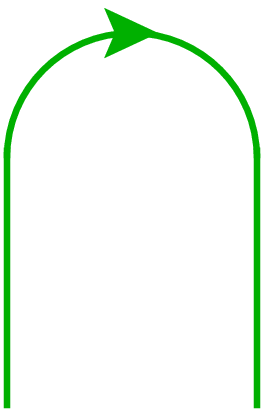}  {7}
\put(4.8,-8.8)   {\sse$U$}
\put(15.9,-8.8)  {\sse$
      {\phantom U}^\vee_{}\!{U}$}
                  \end{picture}    }
\\
\begin{picture}(0,19){}\end{picture} &&&&&&& \\ \cline{1-8}

%%%%%%%%%%%%%%%%%%%%%%%%%%%%%%%%%%%%%%%%%%%
                 %%%   dual morphisms

\mclo{\begin{picture}(61,0){}\end{picture}}&\mclll{
~$f^\vee_{} =$
\begin{picture}(36,53)(0,18)       \apppicture{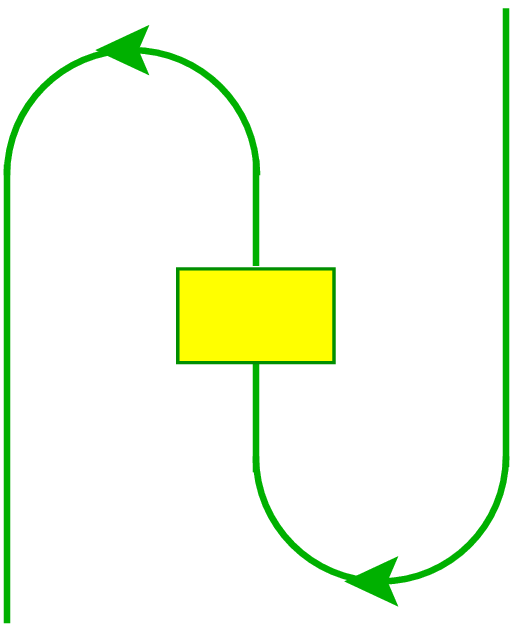}  {7}
\put(3.2,-8.8)   {\sse$V^\vee$}
\put(26.4,25.2)  {\tiny$f$}
\put(44.8,55.7)  {\sse$U^\vee$}
                  \end{picture}    }&
\mclo{\begin{picture}(61,0){}\end{picture}}&\mclll{
~${\phantom f}^{\vee\!}_{}\!{f} =$
\begin{picture}(36,30)(0,18)       \apppicture{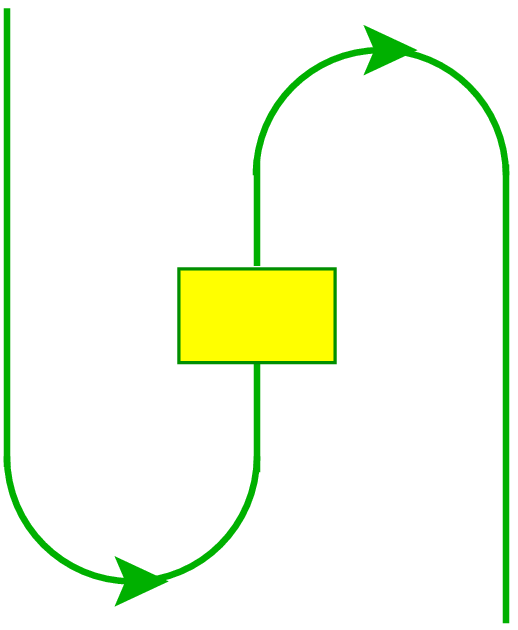}  {7}
\put(-5.8,55.7)  {\sse$
      {\phantom U}^\vee_{}\!{U}$}
\put(26.4,25.2)  {\tiny$f$}
\put(35.8,-8.8)  {\sse$
      {\phantom U}^\vee_{}\!{V}$}
                 \end{picture}     }
\\
\mclo{}&\mclll{\begin{picture}(0,27){}\end{picture}}&
\mclo{}&\mclll{\begin{picture}(0,27){}\end{picture}} \\
                                          \cline{2-4}\cline{6-8}
            \multicolumn8l{{}}\\[-1.08em] \cline{2-4}\cline{6-8}

\end{tabular}
%%%%%%%%%%%%%%%%%%%%%%%%%%%%%%%%%%%%%%%%%%%%%%%%%%%%%%%%%%%%%%%%%%%%%%
\newpage

\noindent
The next table lists the structural morphisms of a (co)algebra: the
product, unit, coproduct, and counit (see equations \erf{m-Delta} and
\erf{ass-coass});
the \rep\ morphism for a general left-module (see equation \erf{1m});
the \rep\ morphism for an induced left-module as well as the
right-\rep\ morphisms for $\alpha$-induced modules (see \erf{rr+-}):

\vskip 1.2em

\begin{tabular}{|c|c|c|c|}
\hline\hline
%%%%%%%%%%%%%%%%%%%%%%%%%%%%%%%%%%%%%%%%%%%%%%%%%%%%%%%%%%%%%%%%%%%%%%
                 %%%   product, unit, coproduct, counit

~$m =$
\begin{picture}(36,32)(0,18)       \apppicture{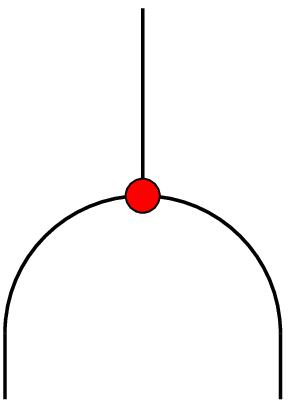}  {7}
\put(3.5,-8.8)   {\sse$A$}
\put(15.5,36.8)  {\sse$A$}
\put(26.5,-8.8)  {\sse$A$}
                 \end{picture}     &
~~~$\eta =$       
\begin{picture}(36,20)(0,12)        \apppicture{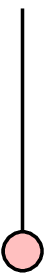}  {7}   
\put(5.8,25.8)   {\sse$A$}
                 \end{picture}     &
$\Delta =$       
\begin{picture}(36,30)(0,18)       \apppicture{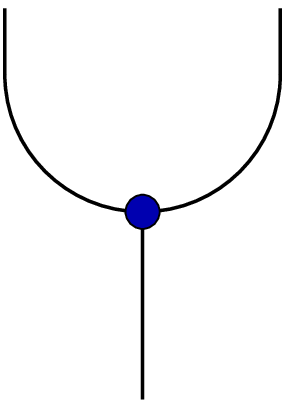}  {7}   
\put(3.8,36.8)   {\sse$A$}
\put(15.2,-8.8)  {\sse$A$}
\put(26.8,36.8)  {\sse$A$}
                 \end{picture}     &
~~~$\eps =$       
\begin{picture}(26,20)(0,6)        \apppicture{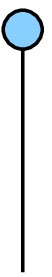}  {7}   
\put(5.1,-8.8)   {\sse$A$}
                  \end{picture}
\\ \begin{picture}(0,27){}\end{picture} &&& \\ \hline\hline

%%%%%%%%%%%%%%%%%%%%%%%%%%%%%%%%%%%%%%%%%%%
                 %%%   \rep\ morphism, induced \rep

~$\r_M^{} =$
\begin{picture}(36,45)(0,24)       \apppicture{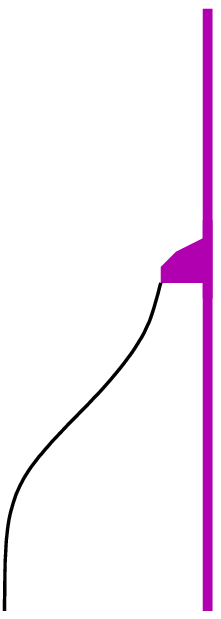}  {7}   
\put(3.5,-8.8)   {\sse$A$}
\put(19.5,-9.4)  {\sse$\M$}
\put(20.5,54.4)  {\sse$\M$}
                  \end{picture}    &
\mclll{
~~$\r_{\!\alpha_\AA^+(U)}^{\rm left} =
  \r_{\!\alpha_\AA^-(U)}^{\rm left} =
  \r_{\!\Ind_\AA(U)}^{} =\;$
\begin{picture}(46,45)(0,24)       \apppicture{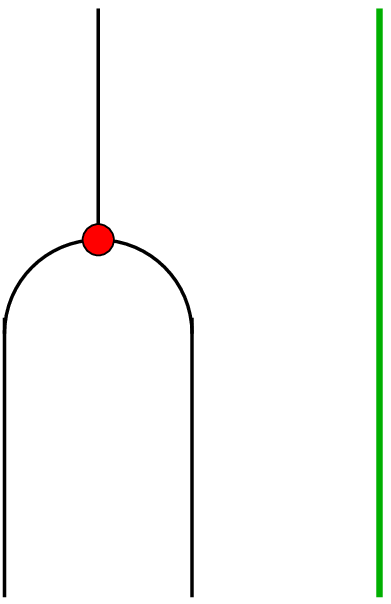}  {7}
\put(3.5,-8.8)   {\sse$A$}
\put(12.5,54.4)  {\sse$A$}
\put(19.5,-8.8)  {\sse$A$}
\put(35.5,-9.4)  {\sse$U$}
\put(36.2,54.4)  {\sse$U$}
                  \end{picture}    }
\\ \begin{picture}(0,31){}\end{picture} &\mclll{} \\ \cline{1-4}

%%%%%%%%%%%%%%%%%%%%%%%%%%%%%%%%%%%%%%%%%%%
                 %%%   right \rep\ morphisms for alpha ind 

\mcll{
~$\r_{\!\alpha_\AA^+(U)}^{\rm right} =$
\begin{picture}(46,45)(0,24)       \apppicture{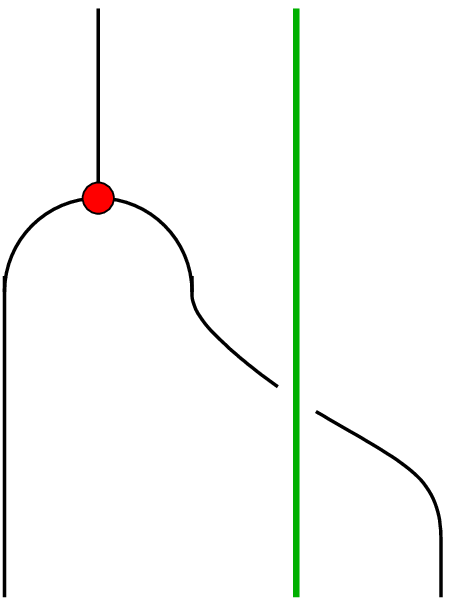}  {7}
\put(3.3,-8.8)   {\sse$A$}
\put(12.5,54.4)  {\sse$A$}
\put(28.5,-9.4)  {\sse$U$}
\put(29.2,54.4)  {\sse$U$}
\put(40.5,-8.8)  {\sse$A$}
                  \end{picture}    }&
\mcll{
~$\r_{\!\alpha_\AA^-(U)}^{\rm right} =$
\begin{picture}(46,45)(0,24)       \apppicture{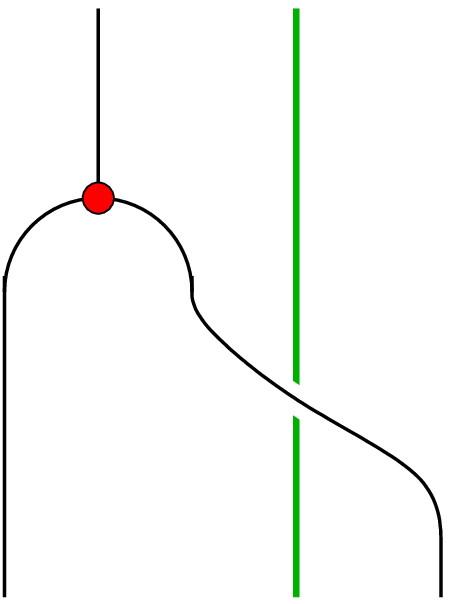}  {7}
\put(3.3,-8.8)   {\sse$A$}
\put(12.5,54.4)  {\sse$A$}
\put(28.5,-9.4)  {\sse$U$}
\put(29.2,54.4)  {\sse$U$}
\put(40.5,-8.8)  {\sse$A$}
                  \end{picture}    }
\\ \mcll{\begin{picture}(0,31){}\end{picture}} &\mcll{} \\ \hline\hline

\end{tabular}
%%%%%%%%%%%%%%%%%%%%%%%%%%%%%%%%%%%%%%%%%%%%%%%%%%%%%%%%%%%%%%%%%%%%%%
\newpage

In the following table we list some specific idempotents: the idempotents
$P^{l/r}_\AA(U)$ (see equation \erf{PU-def}) on which the left and right 
local induction are based and which appear in the \Definition \ref{csplit} 
of a \csplit\ 
Frobenius algebra; those appearing in the definition of the tensor product 
of local modules ($P_{M{\otimes}N}$, see formula \erf{eq:P-2-tensor}); 
and also the idempotents $Q_{r/l}(M_{l/r})$ defined in \erf{eq:pic19}, which 
appear in the functorial equivalences between $\cC_{C_l(A)}^{\sss\rmloc}$
and $\cC_{C_r(A)}^{\sss\rmloc}$.

\vskip 1.2em

\begin{tabular}{|ll|ll|}
\hline \hline
%%%%%%%%%%%%%%%%%%%%%%%%%%%%%%%%%%%%%%%%%%%%%%%%%%%%%%%%%%%%%%%%%%%%%%
                 %%%   specific idempotents: P_A

\mcll{
~$P^l_\AA(U)\;:=$
\begin{picture}(70,64)(0,44)       \apppicture{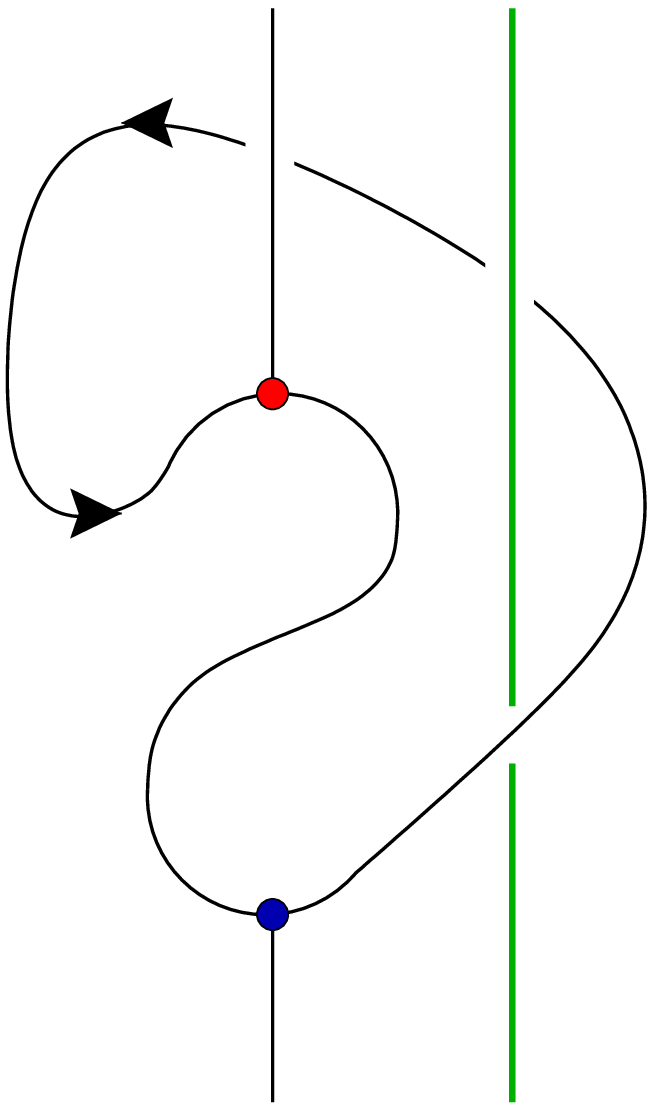}  {4}
\put(22.9,-8.8)  {\sse$A$}
\put(23.5,95.9)  {\sse$A$}
\put(43.5,-8.8)  {\sse$U$}
\put(43.9,95.9)  {\sse$U$}
                 \end{picture}     }&
\mcll{
~$P^r_\AA(U)\;:=$
\begin{picture}(70,60)(0,44)       \apppicture{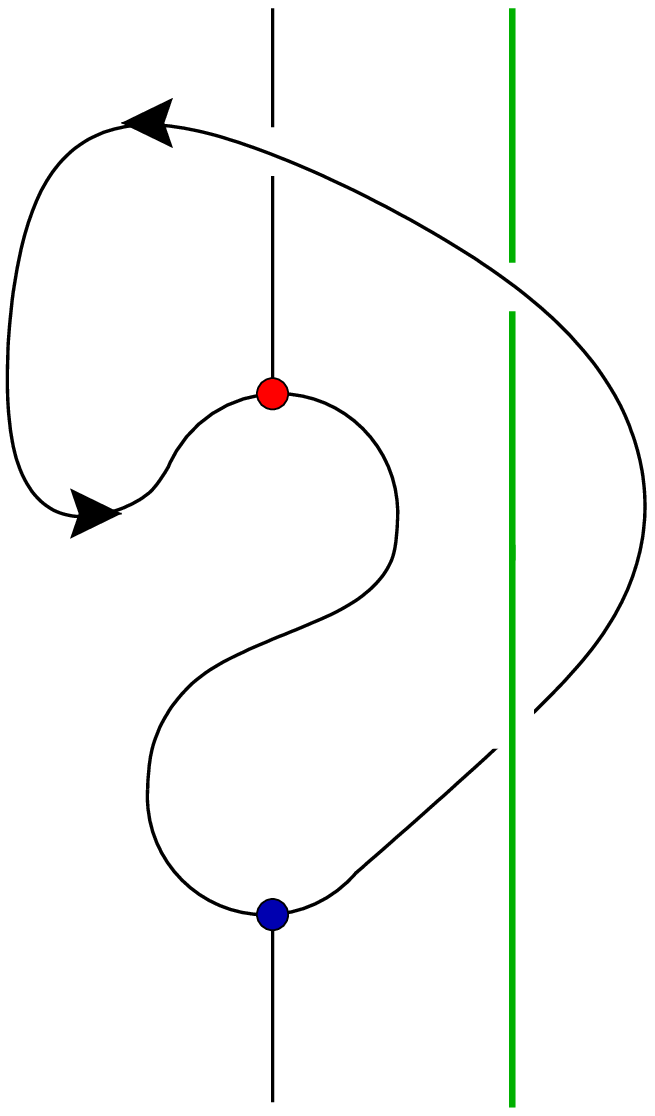}  {4}
\put(22.9,-8.8)  {\sse$A$}
\put(23.5,95.9)  {\sse$A$}
\put(43.5,-8.8)  {\sse$U$}
\put(43.9,95.9)  {\sse$U$}
                 \end{picture}     }
\\ \begin{picture}(0,52){}\end{picture} &&&\\ \cline{1-4}

%%%%%%%%%%%%%%%%%%%%%%%%%%%%%%%%%%%%%%%%%%%
                 %%%   specific idempotents: P_\otimes

\mclll{
~$P_{M{\otimes}N}\;:=$ 
\begin{picture}(135,54)(0,22)      \apppicture{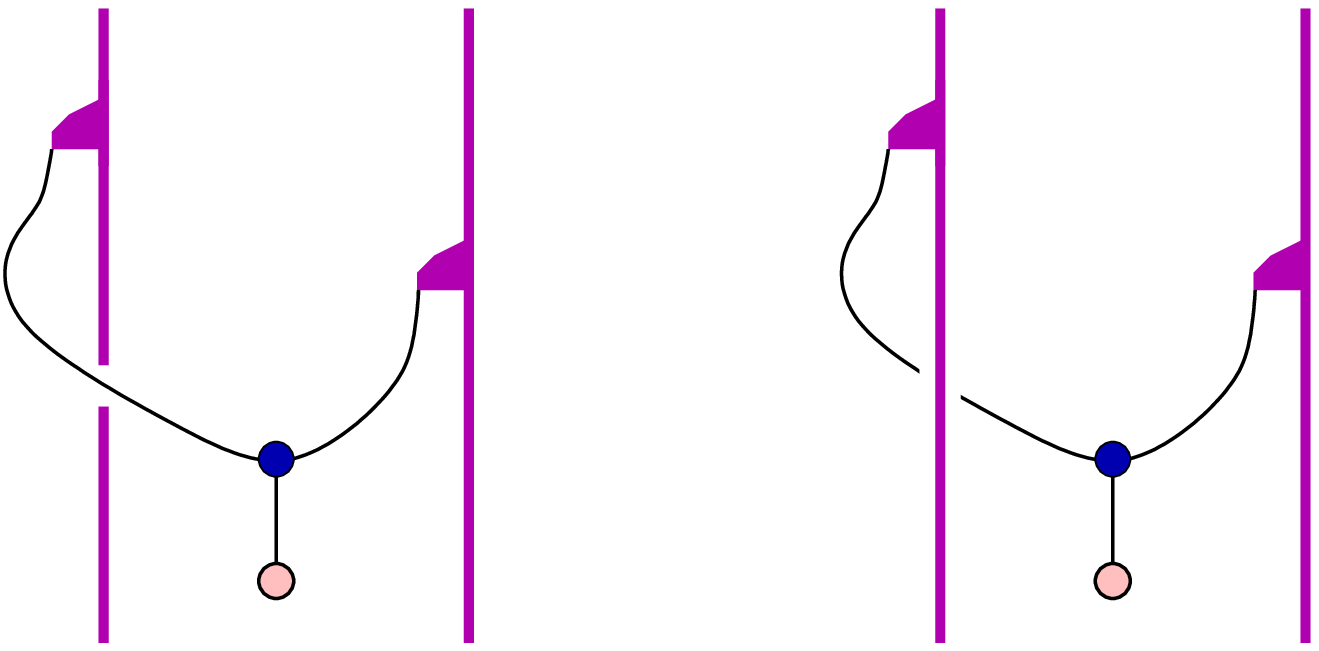}  {10}
\put(13.7,-9.6)  {\sse$\M$}
\put(14.3,57.4)  {\sse$\M$}
\put(45.1,-9.6)  {\sse$\dot N$}
\put(45.5,57.4)  {\sse$\dot N$}
\put(62.5,22)    {\small$=$}
\put(84.4,-9.6)  {\sse$\M$}
\put(84.9,57.4)  {\sse$\M$}
\put(115.8,-9.6) {\sse$\dot N$}
\put(116.2,57.4) {\sse$\dot N$}
                 \end{picture}
                 }
\\ \mclll{ \begin{picture}(0,34){}\end{picture} }\\
\cline{1-4}

%%%%%%%%%%%%%%%%%%%%%%%%%%%%%%%%%%%%%%%%%%% 
                 %%%   specific idempotents: Q

\mcll{
~$Q_r(M_l)\;:=$ 
\begin{picture}(70,60)(0,40)       \apppicture{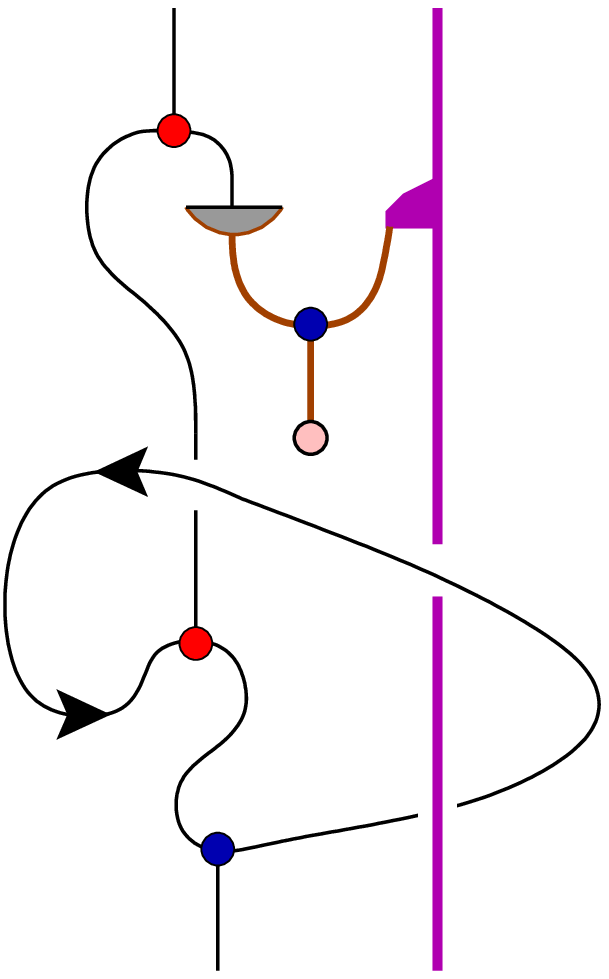}  {6}   
\put(17.2,84.5)  {\sse$A$}
\put(20.5,-9.2)  {\sse$A$}
\put(33.2,49.4)  {\tiny$C_{\!l}^{}$}
\put(37.3,-9.6)  {\sse$\M_l$}
\put(38.3,84.5)  {\sse$\M_l$} 
                 \end{picture}     }&
\mcll{
~$Q_l(M_r)\;:=$ 
\begin{picture}(70,52)(0,40)       \apppicture{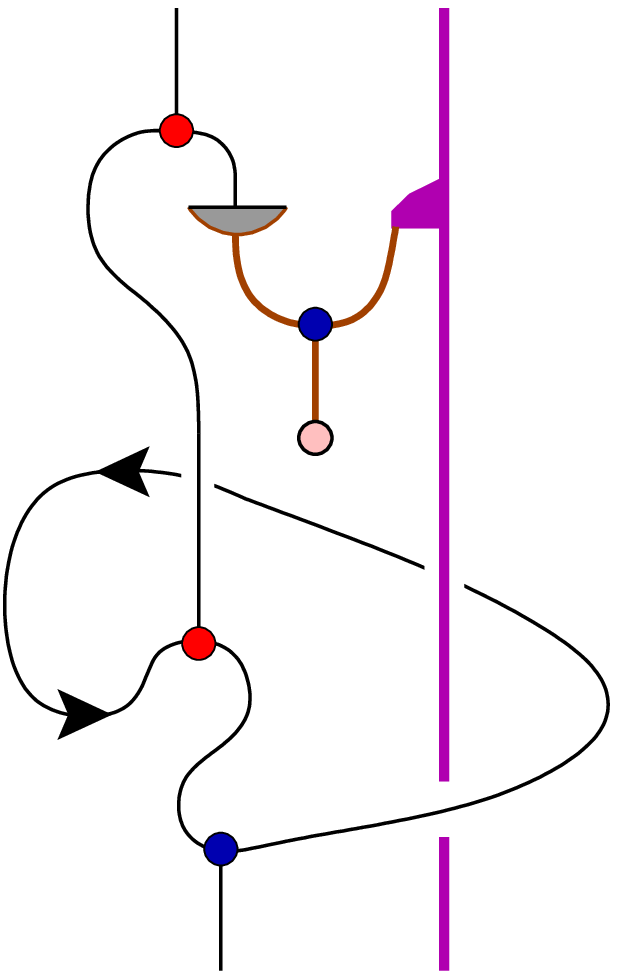}  {6}   
\put(17.2,84.5)  {\sse$A$}
\put(20.5,-9.2)  {\sse$A$}
\put(33.2,49.4)  {\tiny$C_{\!r}^{}$}
\put(37.8,-9.6)  {\sse$\M_r$}
\put(38.5,84.5)  {\sse$\M_r$} 
                 \end{picture}     }
\\ \begin{picture}(0,50){}\end{picture} &&&\\ \hline\hline

%%%%%%%%%%%%%%%%%%%%%%%%%%%%%%%%%%%%%%%%%%%

\end{tabular}

%%%%%%%%%%%%%%%%%%%%%%%%%%%%%%%%%%%%%%%%%%%%%%%%%%%%%%%%%%%%%%%%%%%%%%%
\newpage

\subsection{Defining properties} \label{apptab2}

We now present the defining properties of some of the morphisms displayed in 
\Section \ref{apptab1}.

\medskip

We start with the axioms of a ribbon \cat: the defining properties of
dualities; the functoriality and tensoriality of the braiding;
the functoriality of the twist, and the compatibility of the twist with 
duality and with braiding, see equation \erf{DTB}:

\vskip 1.2em

\noindent
\begin{tabular}{|c|c|cc|c|}
\hline\hline
%%%%%%%%%%%%%%%%%%%%%%%%%%%%%%%%%%%%%%%%%%%%%%%%%%%%%%%%%%%%%%%%%%%%%% 
                 %%%  axioms for dualities
 
\begin{picture}(97,56)(0,23)       \apppicture{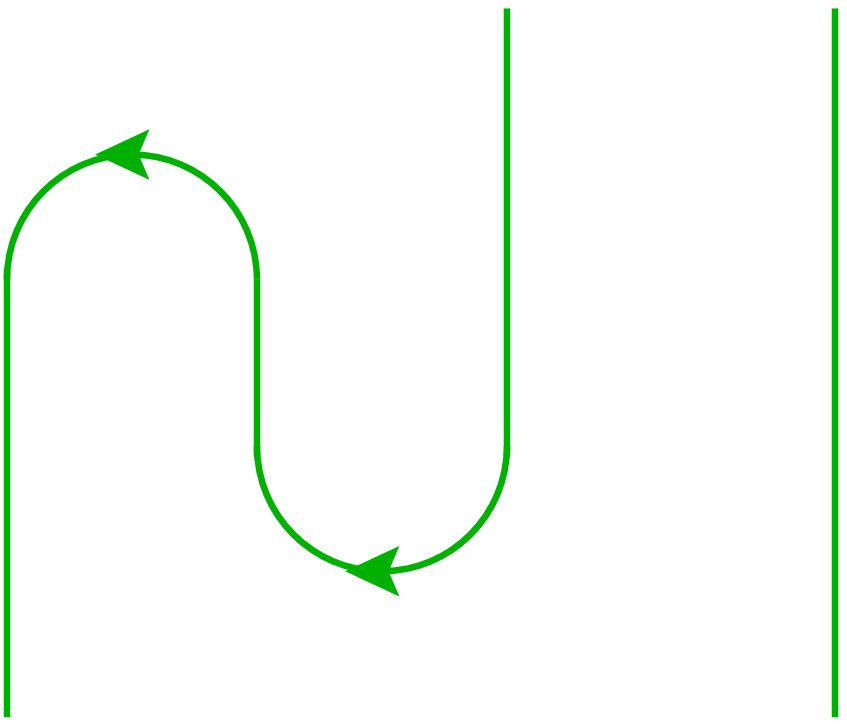}  {14}   
\put(10.2,-9.4)  {\sse$U^\vee$} 
\put(51.8,64.2)  {\sse$U^\vee$}
\put(65.5,25)    {\small$=$}
\put(79.2,-9.4)  {\sse$U^\vee$} 
\put(79.8,64.2)  {\sse$U^\vee$}
                 \end{picture}     &
\begin{picture}(97,49)(0,23)       \apppicture{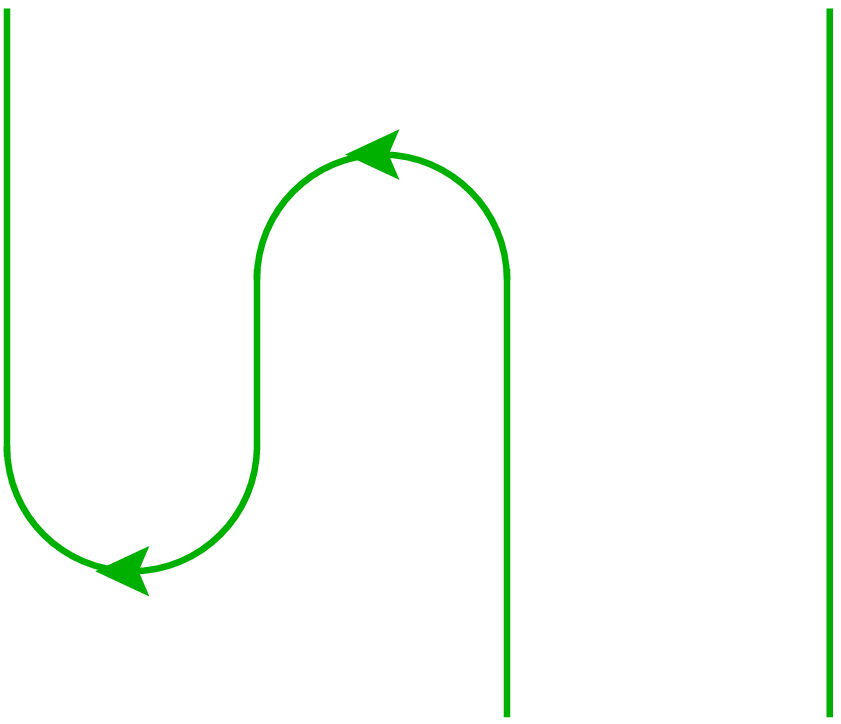}  {14}   
\put(12.2,64.2)  {\sse$U$} 
\put(52.8,-8.8)  {\sse$U$}
\put(65.5,25)    {\small$=$}
\put(80.4,-8.8)  {\sse$U$}
\put(81.1,64.2)  {\sse$U$} 
                 \end{picture}     &
\mcll{
\begin{picture}(97,49)(0,23)       \apppicture{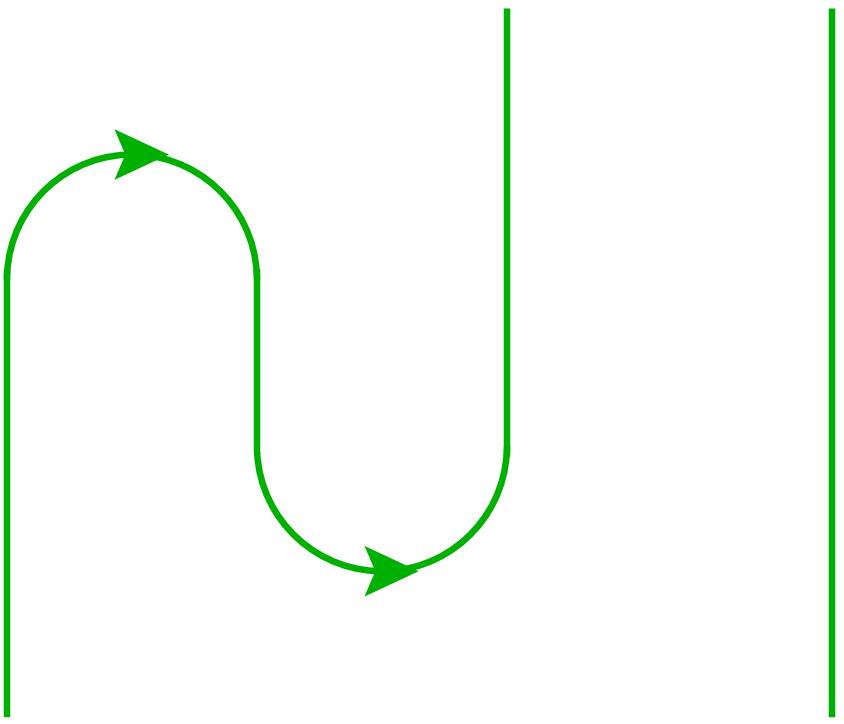}  {14}   
\put(10.8,-8.8)  {\sse$U$} 
\put(53.6,64.2)  {\sse$U$}
\put(65.5,25)    {\small$=$}
\put(80.4,-8.8)  {\sse$U$}
\put(81.1,64.2)  {\sse$U$} 
                 \end{picture}    }&
\begin{picture}(97,49)(0,23)       \apppicture{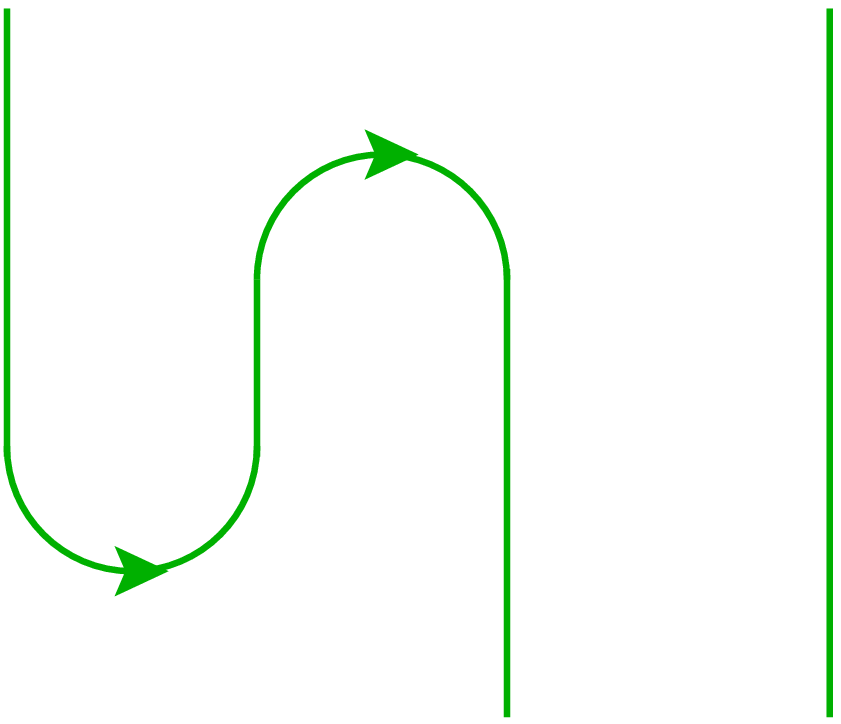}  {14}   
\put(3.2,64.2)   {\sse$
      {\phantom U}^\vee_{}\!{U}$} 
\put(43.8,-9.4)  {\sse$
      {\phantom U}^\vee_{}\!{U}$} 
\put(65.5,25)    {\small$=$}
\put(70.5,-9.4)  {\sse$
      {\phantom U}^\vee_{}\!{U}$} 
\put(70.8,64.2)  {\sse$
      {\phantom U}^\vee_{}\!{U}$} 
                 \end{picture}
\\
\begin{picture}(0,31){}\end{picture} &&&& \\ \cline{1-5}

%%%%%%%%%%%%%%%%%%%%%%%%%%%%%%%%%%%%%%%%%%%
                 %%%  functoriality of the braiding
                 %%%  tensoriality of the braiding

\begin{picture}(97,49)(0,18)       \apppicture{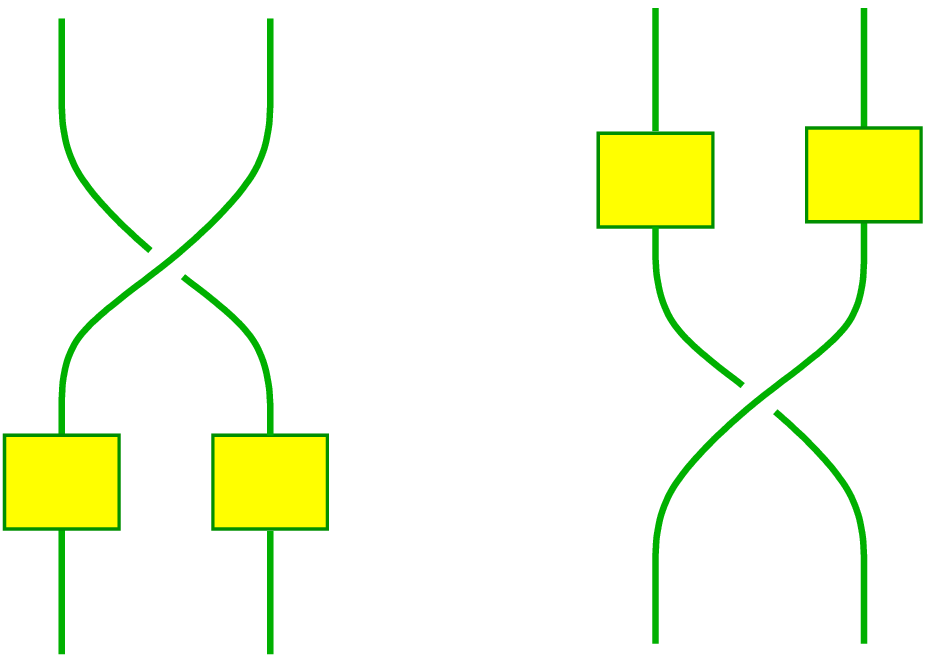}  {10}
\put(11.8,-8.5)  {\sse$U$}
\put(11.8,58.1)  {\sse$W$}
\put(12.5,13.5)  {\tiny$f$}
\put(31,14)      {\tiny$g$}
\put(28.8,-8.5)  {\sse$V$}
\put(29.2,58.1)  {\sse$X$}
\put(45,25)      {\small$=$}
\put(61.3,-8.5)  {\sse$U$}
\put(61.3,58.1)  {\sse$W$}
\put(63,39.7)    {\tiny$g$}
\put(80,39.3)    {\tiny$f$}
\put(78.8,-8.5)  {\sse$V$}
\put(79.2,58.1)  {\sse$X$}
                 \end{picture}     &
\mcll{
\begin{picture}(125,62)(0,28)      \apppicture{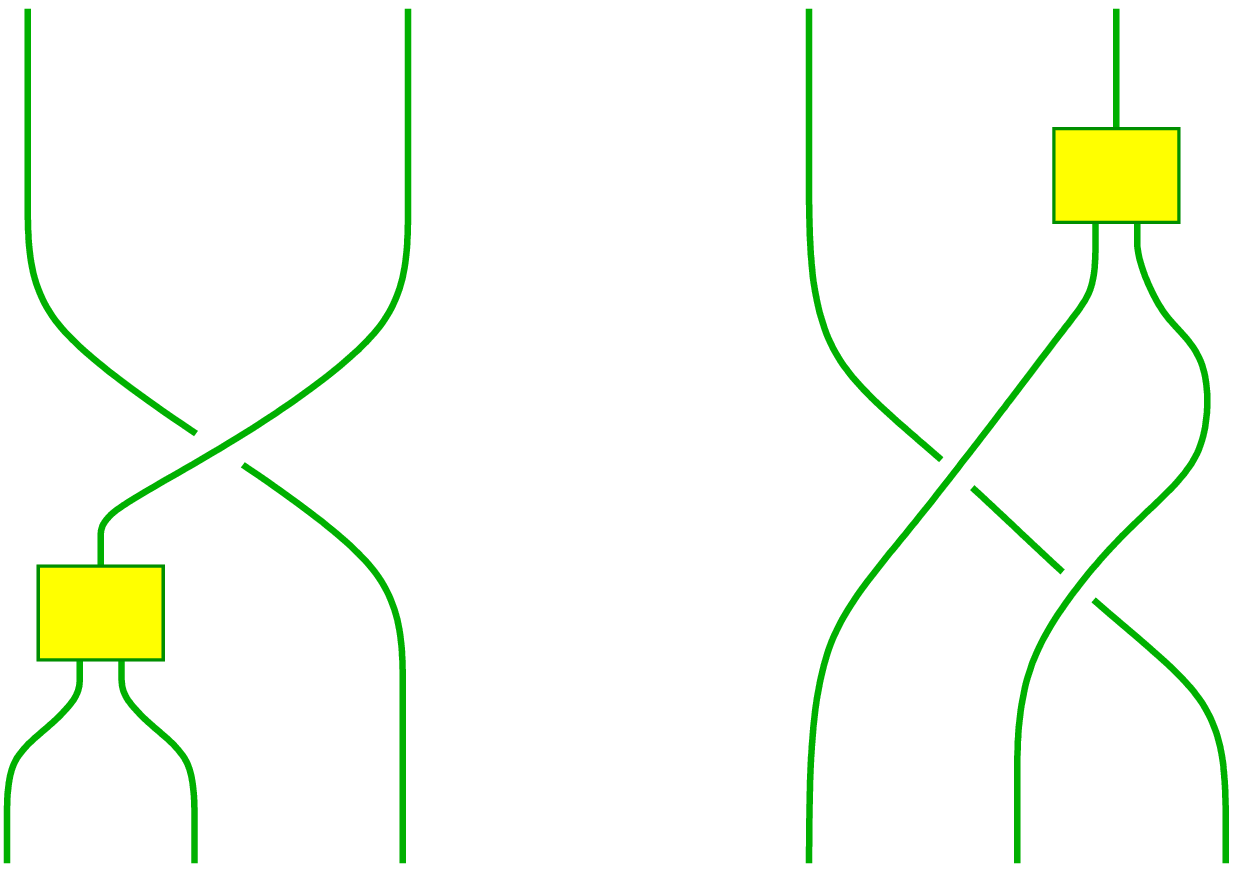}  {10}
\put(7.4,-8.5)   {\sse$U$}
\put(8.4,76.4)   {\sse$W$}
\put(16.3,19.9)  {\tiny$f$}
\put(22.8,-8.5)  {\sse$V$}
\put(34.4,76.4)  {\sse$U\Oti V$}
\put(39.4,-8.5)  {\sse$W$}
\put(58.5,31)    {\small$=$}
\put(74.1,-8.5)  {\sse$U$}
\put(74.4,76.4)  {\sse$W$}
\put(91.5,-8.5)  {\sse$V$}
\put(93.5,76.4)  {\sse$U\Oti V$}
\put(101.3,56.9) {\tiny$f$}
\put(107.7,-8.5) {\sse$W$}
                 \end{picture}     }
\\
\begin{picture}(0,36){}\end{picture} &\mcll{~} \\ \cline{1-4}

%%%%%%%%%%%%%%%%%%%%%%%%%%%%%%%%%%%%%%%%%%%
                 %%%  functoriality of the twist
                 %%%  twist-dual
                 %%%  twist-braid

\begin{picture}(63,36)(0,-7)       \apppicture{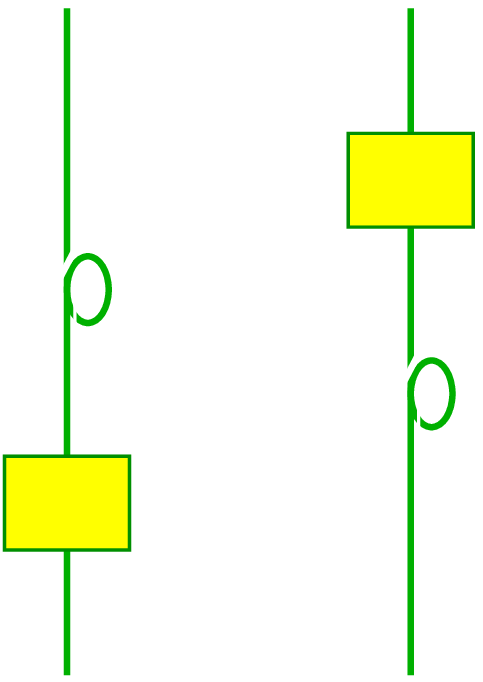}  {10}
\put(12.4,-8.5)  {\sse$U$}
\put(13.2,59.9)  {\sse$U$}
\put(13.5,13.6)  {\tiny$f$}
\put(26,26)      {\small$=$}
\put(41.4,-8.5)  {\sse$U$}
\put(42.2,59.9)  {\sse$U$}
\put(42.5,40.6)  {\tiny$f$}
                 \end{picture}     &
\begin{picture}(63,39)(0,-13)       \apppicture{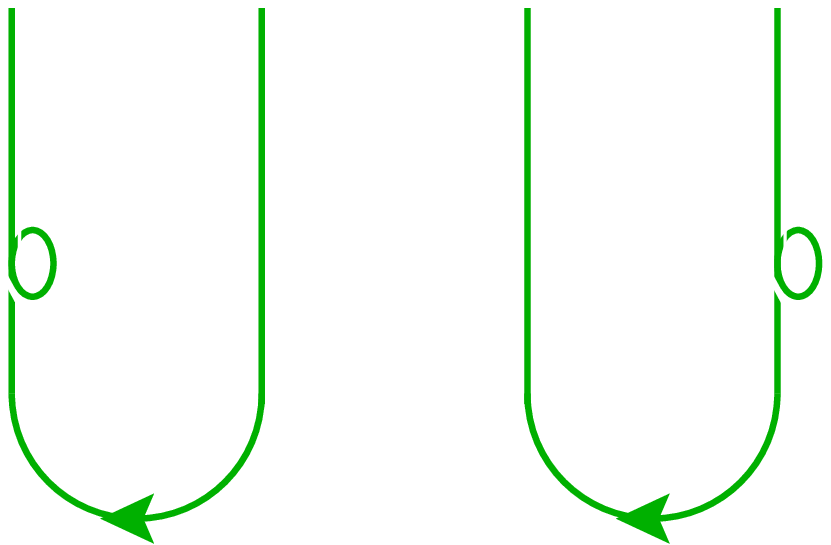} {-2}
\put(-2.8,49.3)  {\sse$U$}
\put(18.5,49.3)  {\sse$U$}
\put(30.5,18)    {\small$=$}
\put(40.8,49.3)  {\sse$U$}
\put(61.8,49.3)  {\sse$U$}
                 \end{picture}     &
\mcll{
\begin{picture}(63,94)(0,12)       \apppicture{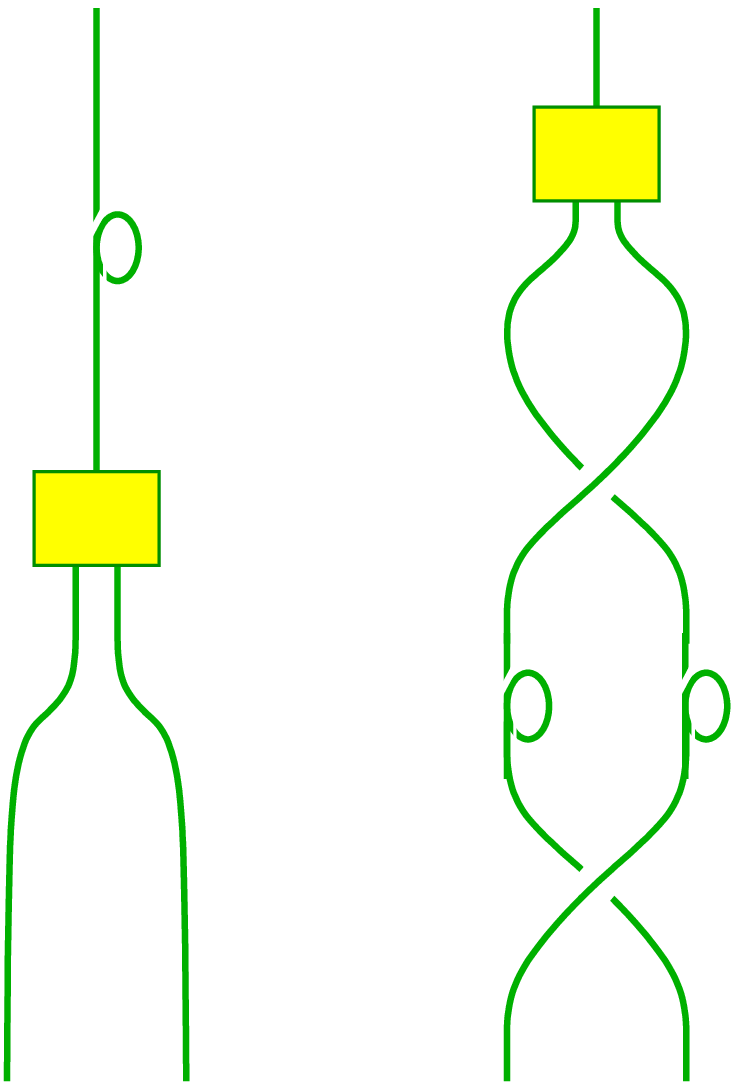}  {18}
\put(15.0,-8.5)  {\sse$U$}
\put(16.2,94.9)  {\sse$U\Oti V$}
\put(23.6,46.1)  {\tiny$f$}
\put(29.7,-8.5)  {\sse$V$}
\put(41,42)      {\small$=$}
\put(57.1,-8.5)  {\sse$U$}
\put(58.2,94.9)  {\sse$U\Oti V$}
\put(65.6,76.6)  {\tiny$f$}
\put(71.9,-8.5)  {\sse$V$}
                 \end{picture}     }
\\
\begin{picture}(0,19){}\end{picture} &&& \\ \cline{1-4}
                      \mclO{{}} \\[-1.08em] \cline{1-4}

\end{tabular}
%%%%%%%%%%%%%%%%%%%%%%%%%%%%%%%%%%%%%%%%%%%%%%%%%%%%%%%%%%%%%%%%%%%%%%
\newpage

Next we display the axioms of a \ssFA\ $A$: associativity of the product, the 
unit property, coassociativity of the coproduct, and the counit property, 
see equations \erf{m-Delta} and \erf{ass-coass}; the Frobenius property, 
the two specialness properties (with the normalisation $\beta_\AA\eq1$) 
and the symmetry property, see \Definition \ref{symm-frob-spec}. 
Finally we show the defining properties of the left and right centers $C_{l/r} 
\eq C_{l/r}(A)$ (see equation \erf{Cl-Cr-defprop}) as well as the two defining 
properties of a (left) \rep, and the defining property of a local (left) 
\rep, see equations \erf{1m} and \erf{eq:def-loc}.

 \vfill

\noindent
\begin{tabular}{|l|ll|l|l|}
\hline \hline
%%%%%%%%%%%%%%%%%%%%%%%%%%%%%%%%%%%%%%%%%%%%%%%%%%%%%%%%%%%%%%%%%%%%%%
                 %%%  algebra, coalgebra

\begin{picture}(99,41)(0,22)       \apppicture{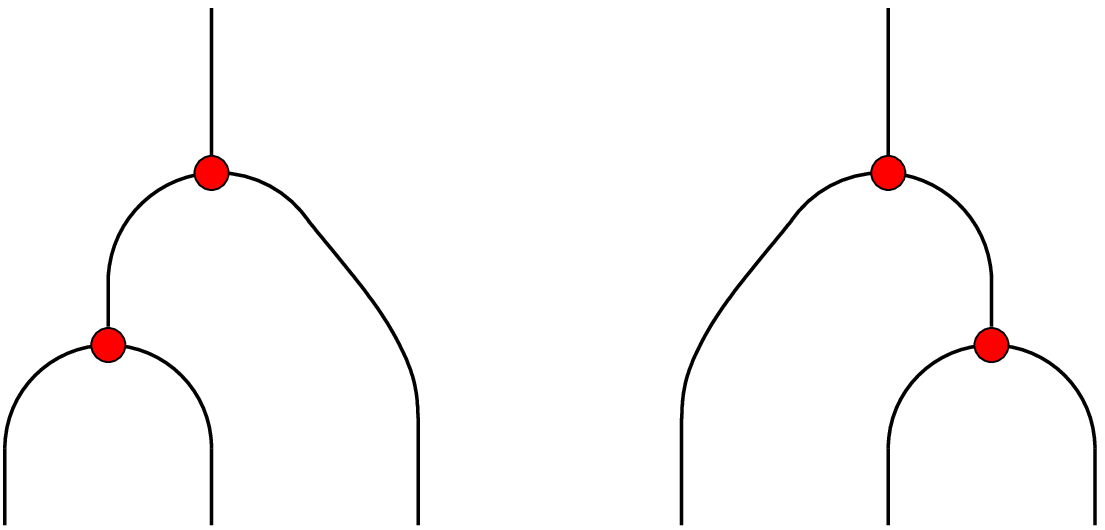}  {3}
\put(-.5,-9.2)   {\sse$A$}
\put(17.1,-9.2)  {\sse$A$}
\put(17.7,47.4)  {\sse$A$}
\put(34.1,-9.2)  {\sse$A$}
\put(44.8,18)    {\small$=$}
\put(56.3,-9.2)  {\sse$A$}
\put(73.9,-9.2)  {\sse$A$}
\put(74.7,47.4)  {\sse$A$}
\put(91.1,-9.2)  {\sse$A$}
                 \end{picture}     &
\mcll{
\begin{picture}(99,29)(0,22)       \apppicture{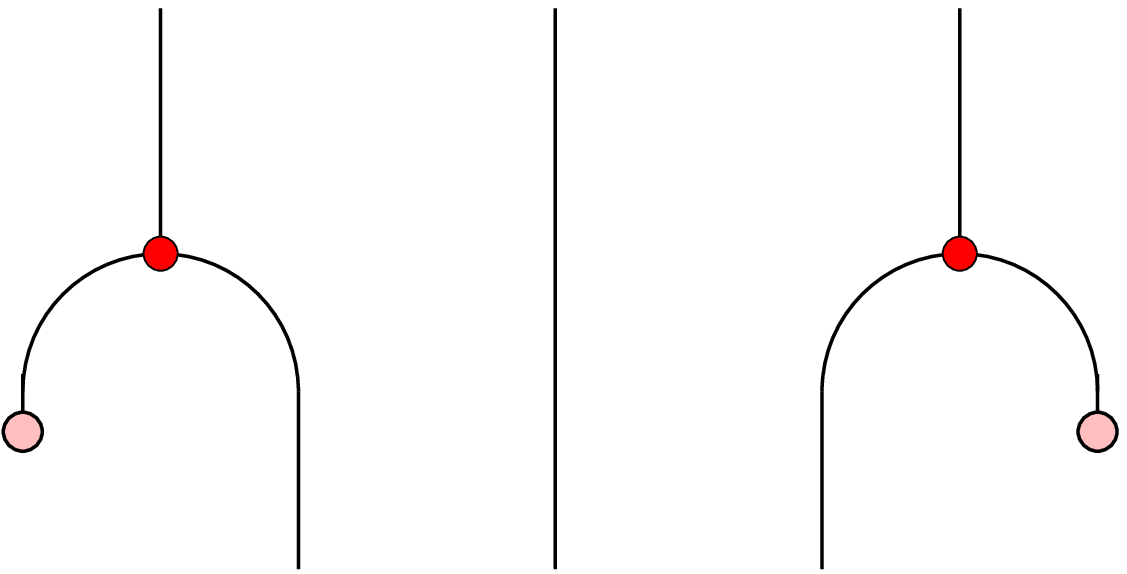}  {3}   
\put(13.3,50.9)  {\sse$A$}
\put(23.3,-9.2)  {\sse$A$}
\put(33.5,18)    {\small$=$}
\put(46.1,-9.2)  {\sse$A$}
\put(46.6,50.9)  {\sse$A$}
\put(56.5,18)    {\small$=$}
\put(68.2,-9.2)  {\sse$A$}
\put(80.1,50.9)  {\sse$A$}
                 \end{picture}     }&
\begin{picture}(99,29)(0,22)       \apppicture{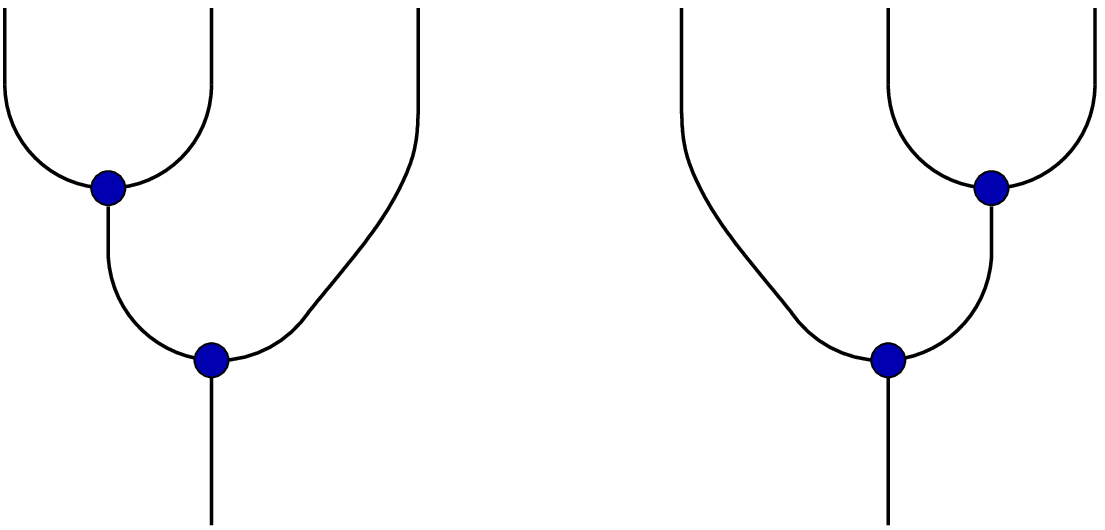}  {3}   
\put(-.2,47.7)   {\sse$A$}
\put(16.6,-9.2)  {\sse$A$}
\put(17.3,47.4)  {\sse$A$}
\put(34.5,47.4)  {\sse$A$}
\put(43.5,18)    {\small$=$}
\put(57.0,47.4)  {\sse$A$}
\put(73.4,-9.2)  {\sse$A$}
\put(74.3,47.4)  {\sse$A$}
\put(91.3,47.4)  {\sse$A$}
                 \end{picture}     &
\begin{picture}(99,29)(0,22)       \apppicture{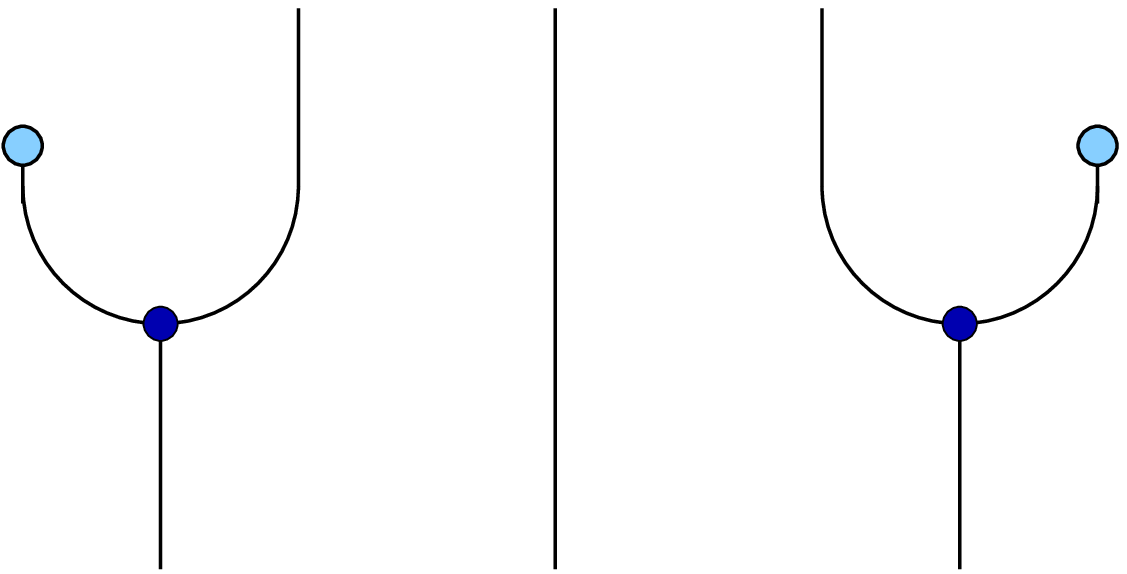}  {3}   
\put(12.9,-9.2)  {\sse$A$}
\put(24.1,50.9)  {\sse$A$}
\put(34.5,18)    {\small$=$}
\put(46.1,-9.2)  {\sse$A$}
\put(46.6,50.9)  {\sse$A$}
\put(56.5,18)    {\small$=$}
\put(68.9,50.9)  {\sse$A$}
\put(79.7,-9.2)  {\sse$A$}
                 \end{picture}     
\\
\begin{picture}(0,30){}\end{picture} &&&& \\ \cline{1-5}

%%%%%%%%%%%%%%%%%%%%%%%%%%%%%%%%%%%%%%%%%%%
                 %%%  Frobenius, special

\mclll{ 
\begin{picture}(97,46)(0,28)       \apppicture{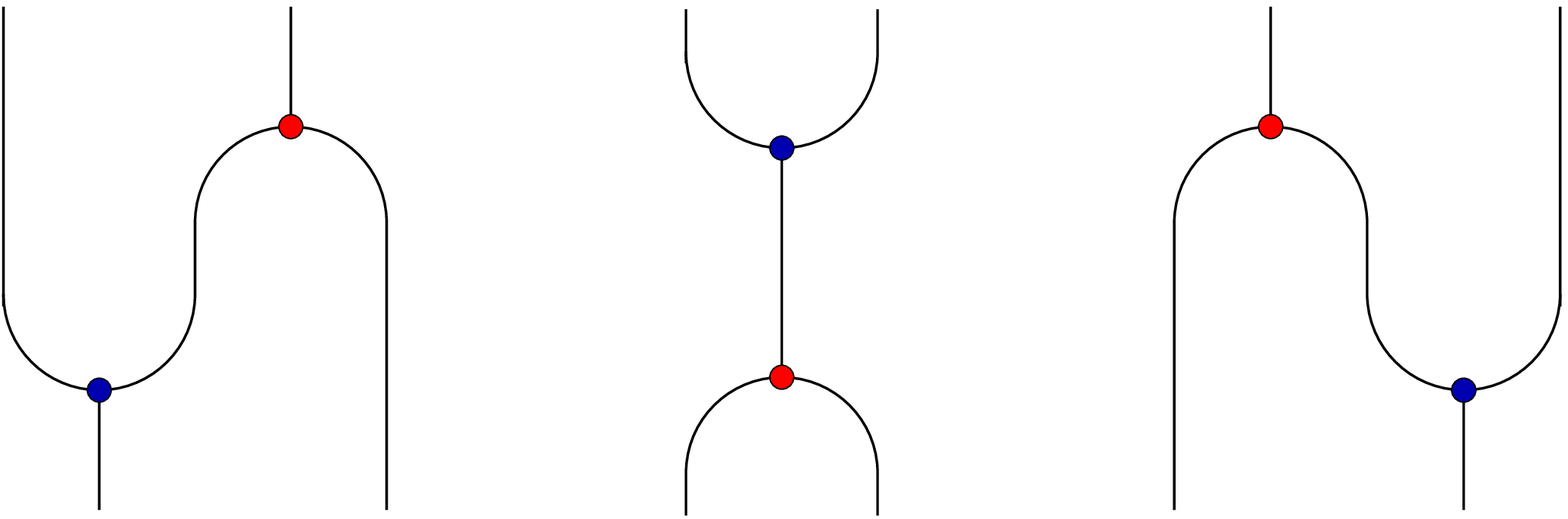}  {22}
\put(19.3,59.9)  {\sse$A$}
\put(28.7,-9.2)  {\sse$A$}
\put(51.2,59.9)  {\sse$A$}
\put(60.6,-9.2)  {\sse$A$}
\put(74.5,25)    {\small$=$}
\put(92.7,-9.2)  {\sse$A$}
\put(93.4,59.9)  {\sse$A$}
\put(113.7,-9.2) {\sse$A$}
\put(114.4,59.9) {\sse$A$}
\put(126.5,25)   {\small$=$}
\put(146.5,-9.2) {\sse$A$}
\put(157.4,59.9) {\sse$A$}
\put(178.5,-9.2) {\sse$A$}
\put(189.4,59.9) {\sse$A$}
                 \end{picture}     }&
\begin{picture}(99,40)(0,20)       \apppicture{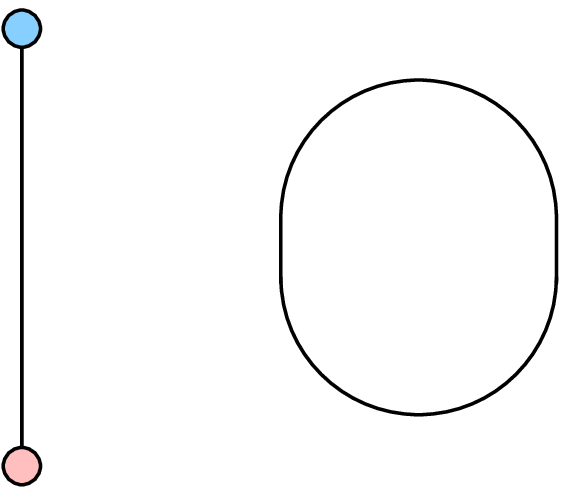}  {31}
\put(39.5,17)    {\small$=$}
\put(76.6,8.6)   {\sse$A$}
                 \end{picture}     &
\begin{picture}(99,47)(0,28)       \apppicture{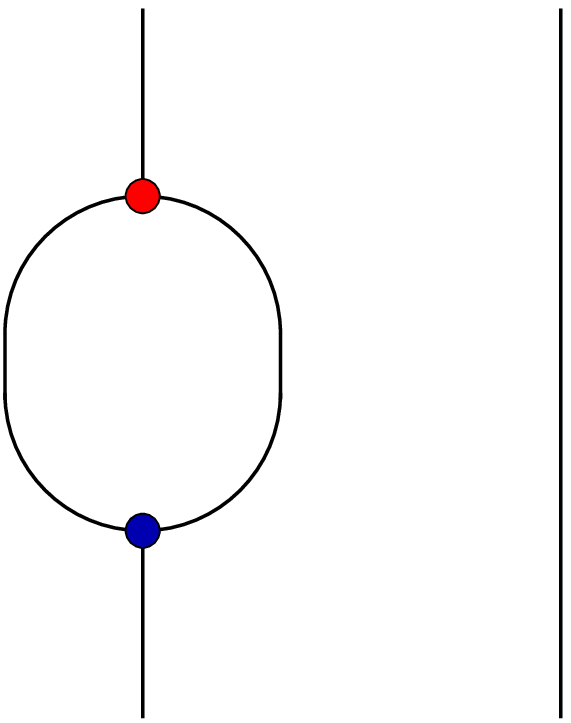}  {22}
\put(30.1,-9.2)  {\sse$A$}
\put(30.6,63.4)  {\sse$A$}
\put(52.5,28)    {\small$=$}
\put(64.9,-9.2)  {\sse$A$}
\put(65.4,63.4)  {\sse$A$}
                 \end{picture}
\\
\mclll{ \begin{picture}(0,36){}\end{picture} } && \\ \cline{1-5} 

%%%%%%%%%%%%%%%%%%%%%%%%%%%%%%%%%%%%%%%%%%%
                 %%%  symmetric %; rep properties

\mcll{
\begin{picture}(155,50)(0,28)      \apppicture{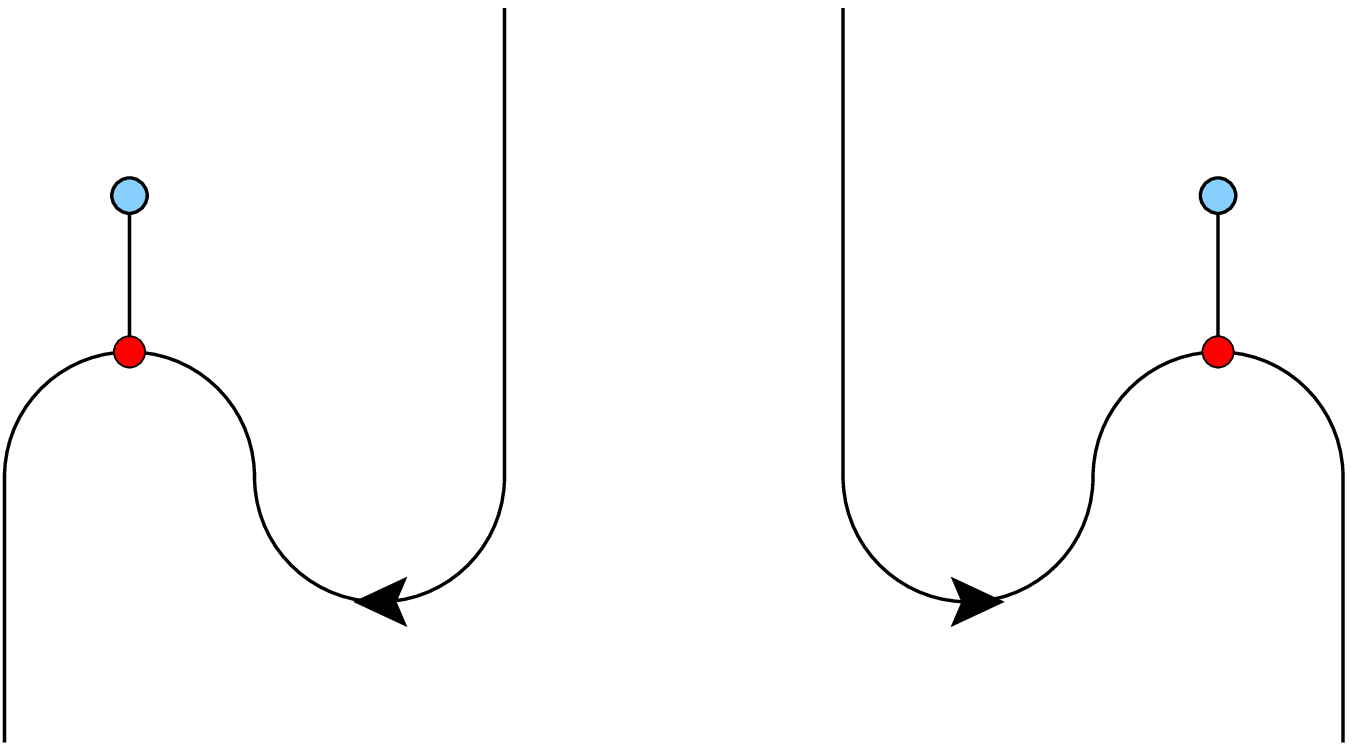}  {22}
\put(18.5,-9.2)  {\sse$A$}
\put(60.7,65.5)  {\sse$A$}
\put(74.9,29)    {\small$=$}
\put(89.1,65.5)  {\sse$A$}
\put(130.5,-9.2) {\sse$A$}
                 \end{picture}     }%&
%                 \multicolumn1{||l|}{
%\begin{picture}(90,48)(0,28)       \apppicture{68}  {6}
%\put(74.5,25)    {\small$=$}
%                 \end{picture}     }&
%\begin{picture}(90,48)(0,28)       \apppicture{69}  {6}
%\put(54.5,25)    {\small$=$}
%                 \end{picture}
\\
\mcll{ \begin{picture}(0,35){}\end{picture} } 
                   \\ \cline{1-2} \mclO{{}}\\[-1.08em] \cline{1-5}

%%%%%%%%%%%%%%%%%%%%%%%%%%%%%%%%%%%%%%%%%%%
                 %%%  centers

\begin{picture}(105,62)(0,33)      \apppicture{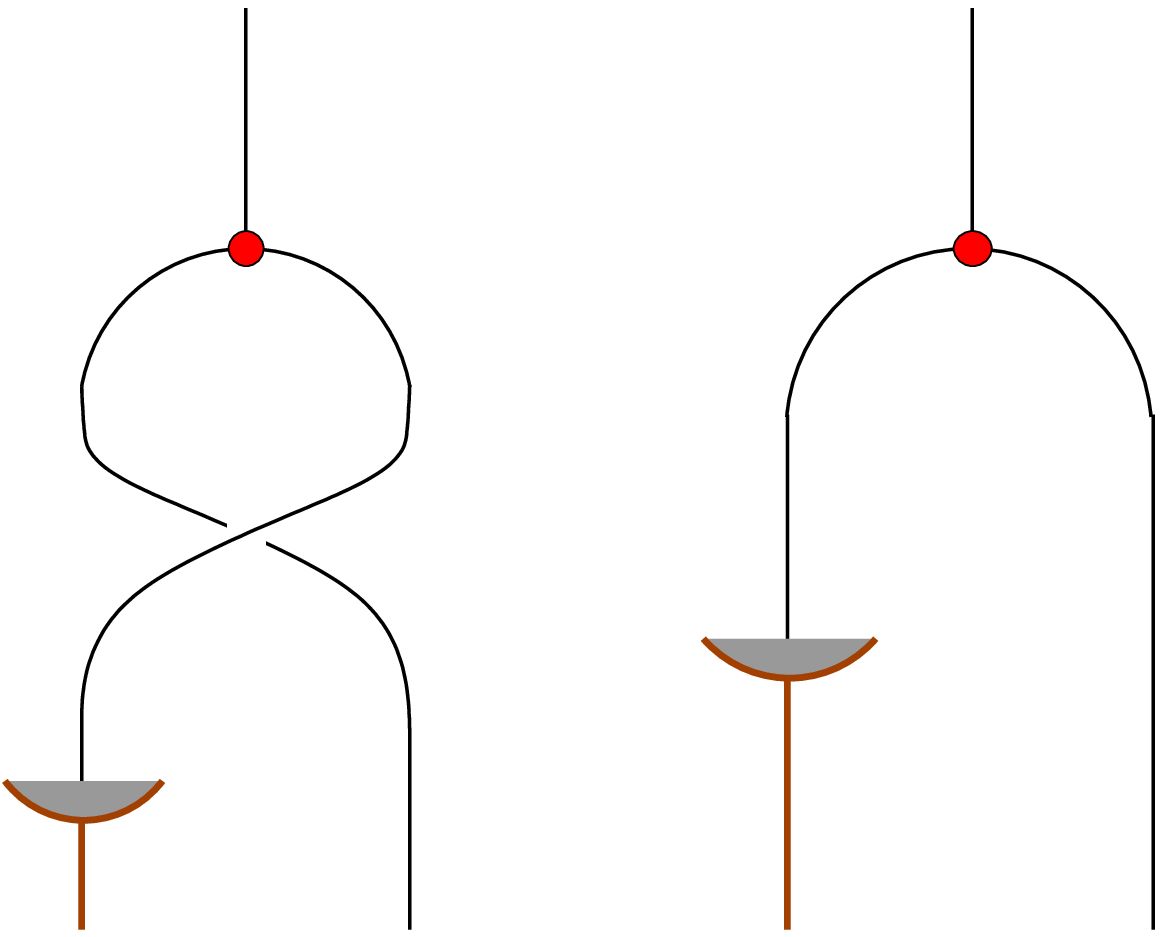}  {0}
\put(2.2,-9.2)   {\sse$C_l$}
\put(17.4,82.3)  {\sse$A$}
\put(30.3,-9.2)  {\sse$A$}
\put(45.2,35)    {\small$=$}
\put(62.2,-9.2)  {\sse$C_l$}
\put(78.5,82.3)  {\sse$A$}
\put(92.8,-9.2)  {\sse$A$}
                 \end{picture}     &
\mcll{
\begin{picture}(105,59)(0,33)      \apppicture{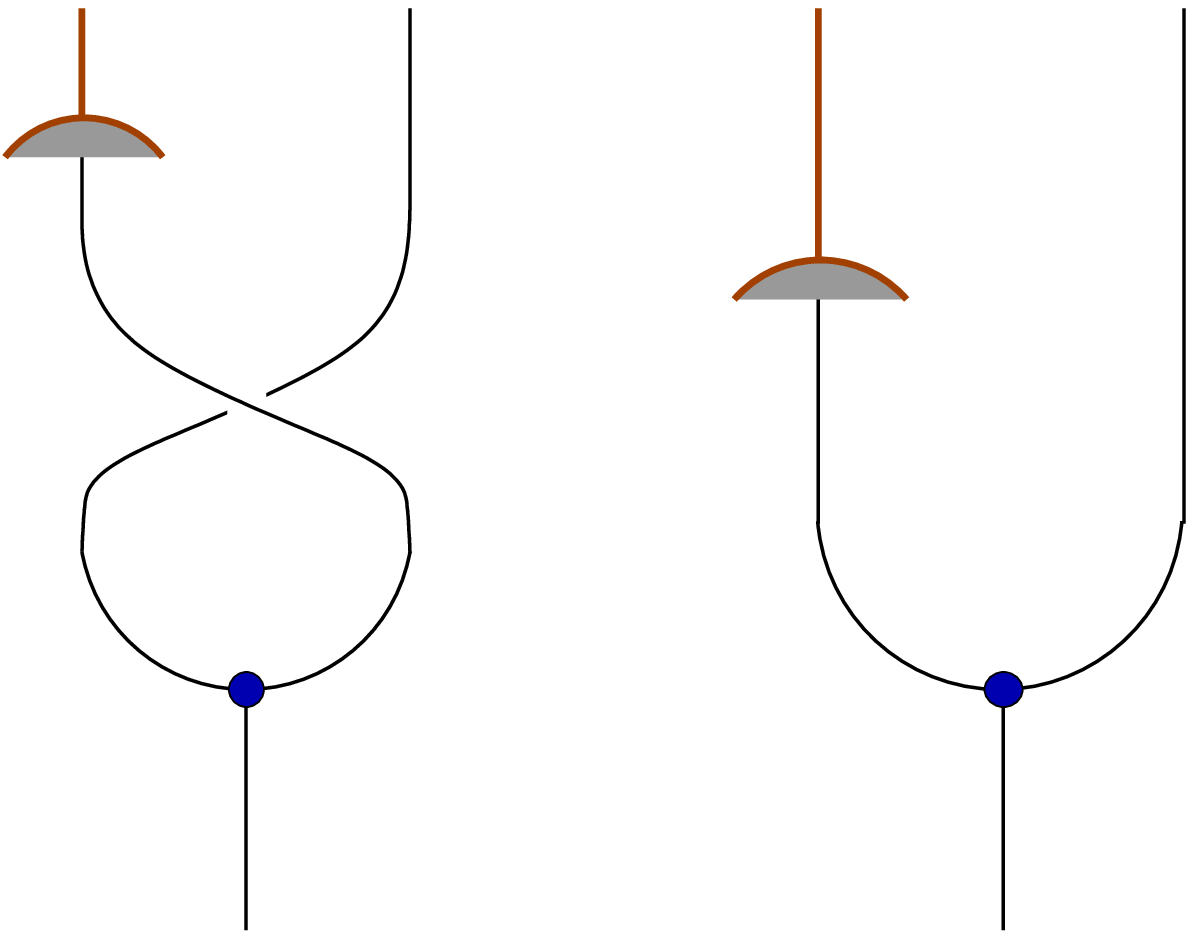}  {0}
\put(3.2,82.3)   {\sse$C_l$}
\put(16.5,-9.2)  {\sse$A$}
\put(31.3,82.3)  {\sse$A$}
\put(45.5,35)    {\small$=$}
\put(65.5,82.3)  {\sse$C_l$}
\put(80.2,-9.2)  {\sse$A$}
\put(96.3,82.3)  {\sse$A$}
                 \end{picture}     }&
\begin{picture}(105,59)(0,33)      \apppicture{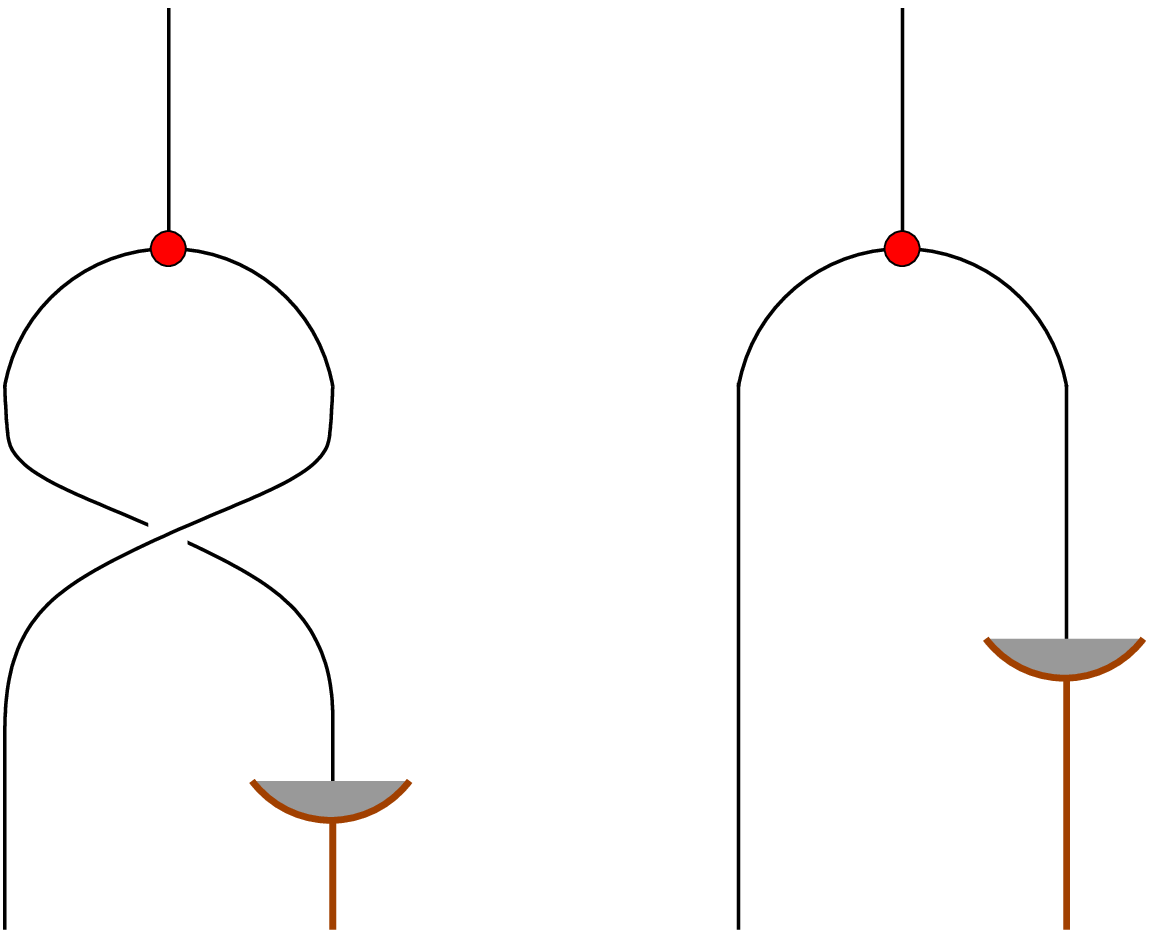}  {5}
\put(1.8,-9.2)   {\sse$A$}
\put(16.1,82.3)  {\sse$A$}
\put(29.1,-9.2)  {\sse$C_r$}
\put(45.5,35)    {\small$=$}
\put(62.7,-9.2)  {\sse$A$}
\put(77.5,82.3)  {\sse$A$}
\put(90.2,-9.2)  {\sse$C_r$}
                 \end{picture}     &
\begin{picture}(105,59)(0,33)      \apppicture{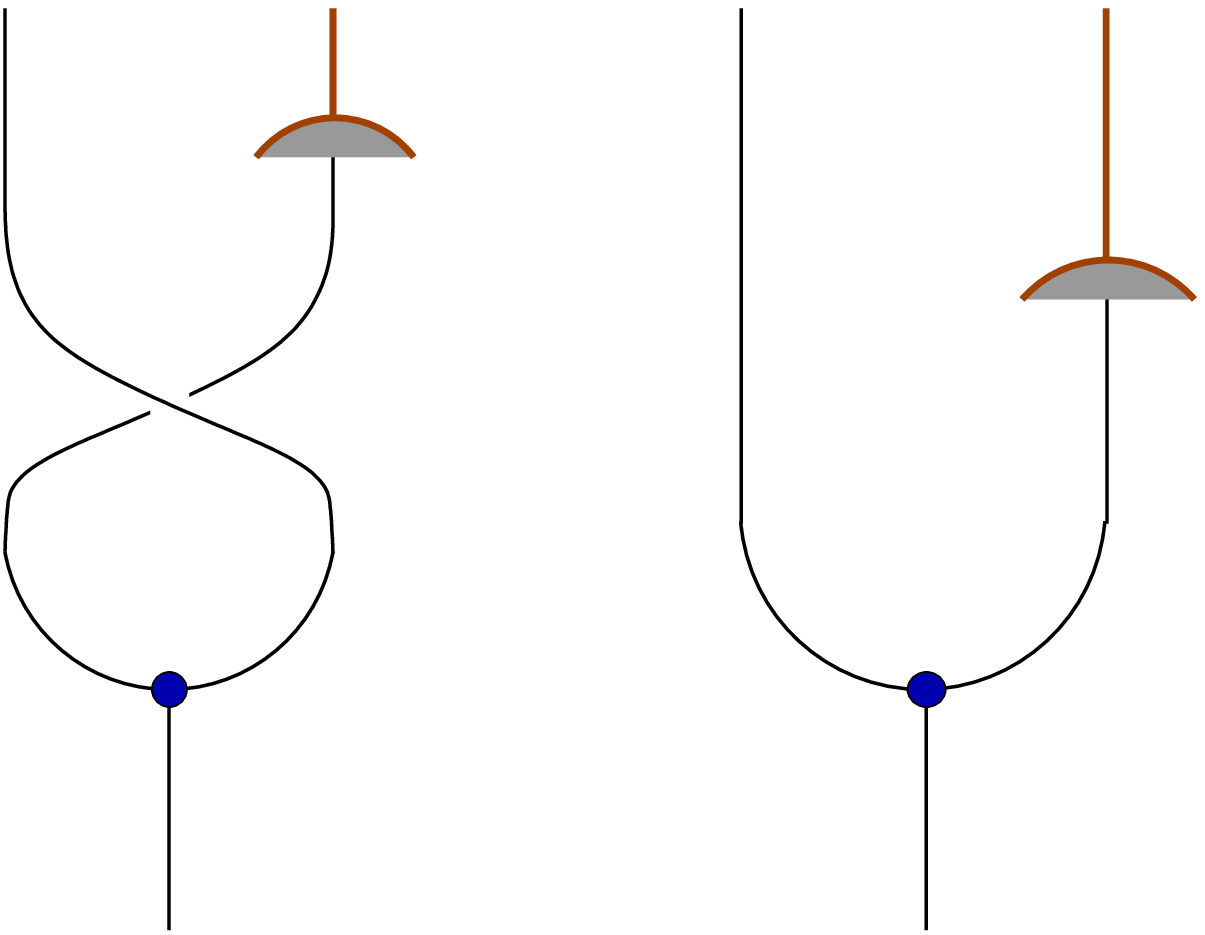}  {5}
\put(2.1,82.3)   {\sse$A$}
\put(15.2,-9.2)  {\sse$A$}
\put(29.1,82.3)  {\sse$C_r$}
\put(45.5,35)    {\small$=$}
\put(63.9,82.3)  {\sse$A$}
\put(79.2,-9.2)  {\sse$A$}
\put(93.7,82.3)  {\sse$C_r$}
                 \end{picture}      
\\
\begin{picture}(0,41){}\end{picture} & \mcll{{}} && \\ \hline\hline

%%%%%%%%%%%%%%%%%%%%%%%%%%%%%%%%%%%%%%%%%%%
                 %%%   rep properties, local rep

\begin{picture}(90,48)(0,28)       \apppicture{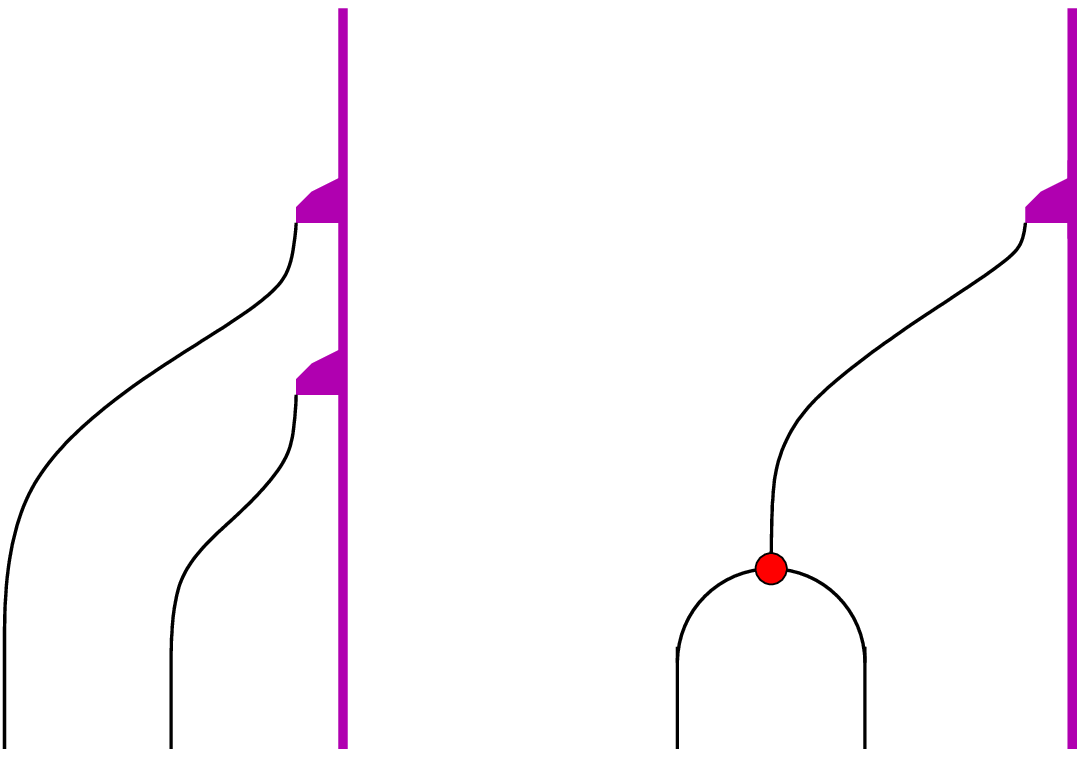}  {6}
\put(2.3,-8.5)   {\sse$A$}
\put(16.3,-8.5)  {\sse$A$}
\put(29.4,-9.8)  {\sse$\M$}
\put(30.5,66.2)  {\sse$\M$}
\put(48.5,28)    {\small$=$}
\put(58.7,-8.5)  {\sse$A$}
\put(74.2,-8.5)  {\sse$A$}
\put(90.9,-9.8)  {\sse$\M$}
\put(93.1,66.2)  {\sse$\M$}
                 \end{picture}     &
\mcll{
\begin{picture}(90,48)(0,28)       \apppicture{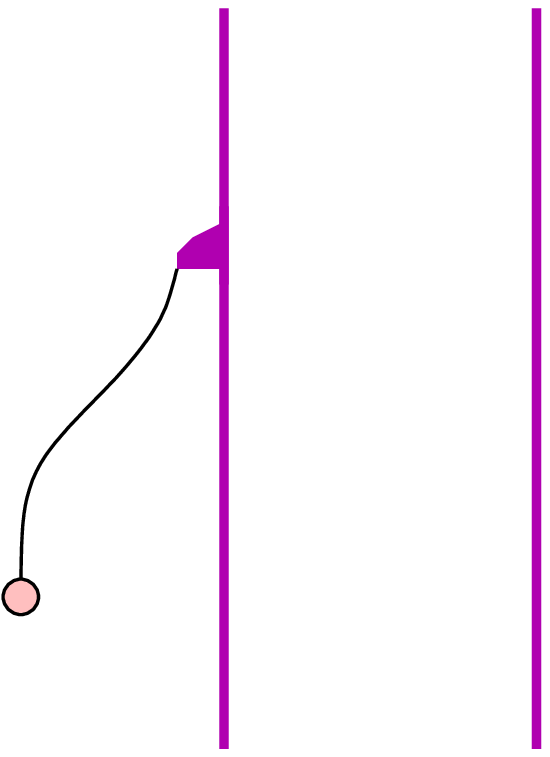}  {22}
\put(35.6,-9.8)  {\sse$\M$}
\put(36.8,66.2)  {\sse$\M$}
\put(49.5,28)    {\small$=$}
\put(61.6,-9.8)  {\sse$\M$}
\put(62.8,66.2)  {\sse$\M$}
                 \end{picture}     }&
\begin{picture}(90,64)(0,48)       \apppicture{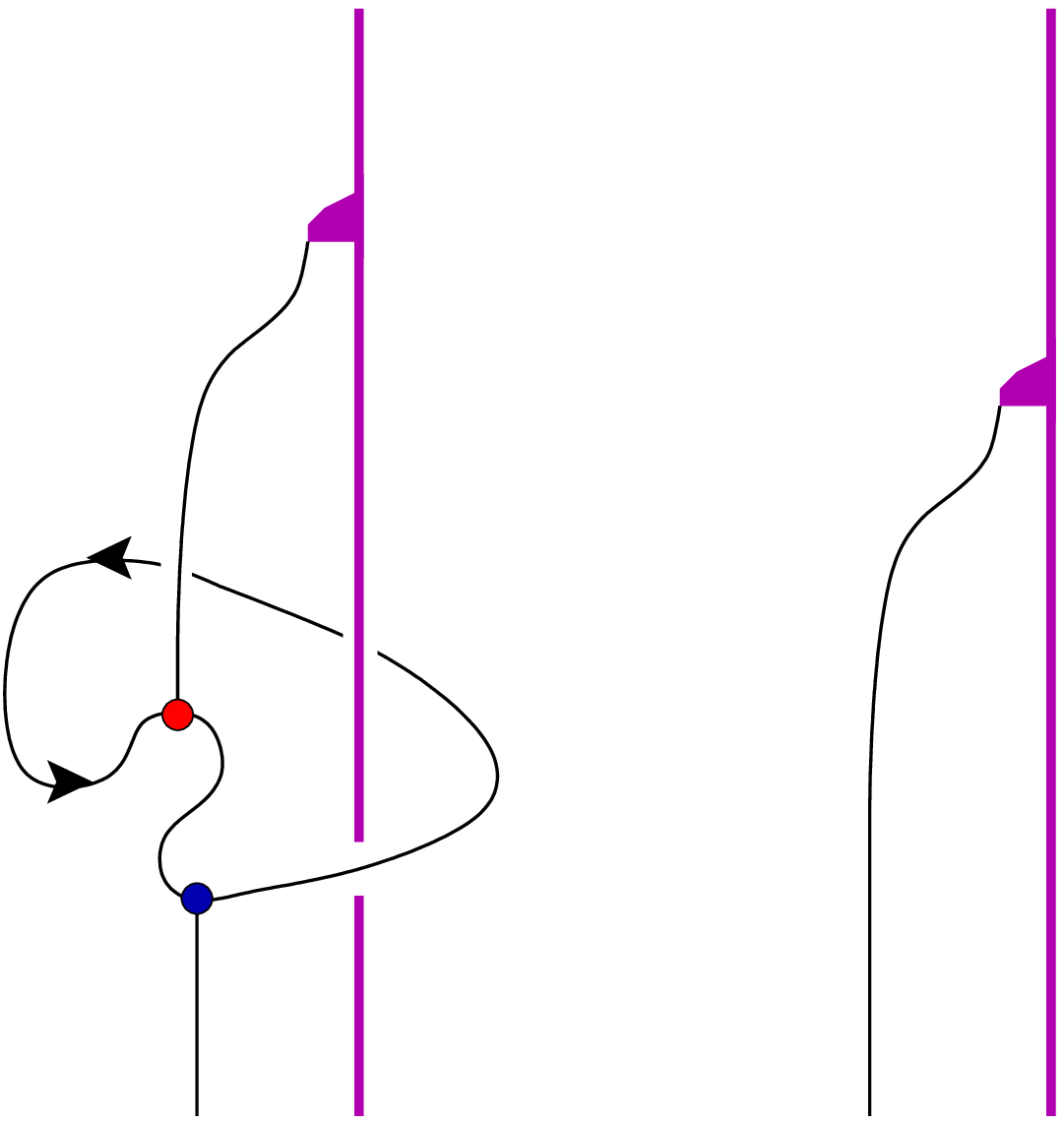}  {6}
\put(18.7,-8.5)  {\sse$A$}
\put(31.2,-9.8)  {\sse$\M$}
\put(32.4,98.5)  {\sse$\M$}
\put(57.5,38)    {\small$=$}
\put(76.3,-8.5)  {\sse$A$}
\put(90.2,-9.8)  {\sse$\M$}
\put(91.4,98.5)  {\sse$\M$}
                 \end{picture}
\\
\begin{picture}(0,55){}\end{picture} & \mcll{{}} & \\ \cline{1-4}
                                 \mclO{{}}\\[-1.08em] \cline{1-4}

%%%%%%%%%%%%%%%%%%%%%%%%%%%%%%%%%%%%%%%%%%%%%%%%%%%%%%%%%%%%%%%%%%%%%%%%%%%%%

\end{tabular}

%%%%%%%%%%%%%%%%%%%%%%%%%%%%%%%%%%%%%%%%%%%

\newcommand\wb{\,\linebreak[0]} \def\wB {$\,$\wb}
 \newcommand\Bi[1]    {\bibitem{#1}}
 \renewcommand\J[6]   {{\sl #6\/}, {#1} {#2} ({#3}) {#4--#5} }
 \newcommand\K[7]     {{\sl #7\/}, {#1} {#2} ({#3}) {#4--#5}}
 \newcommand\BOOK[4]  {{\sl #1\/} ({#2}, {#3} {#4})}
 \newcommand\inBO[8]{{\sl #8\/}, in:\ {\sl #1}, {#2}\ ({#3}, {#4} {#5}), p.\ {#6--#7}}
 \def\A     {Algebra }
 \def\dim   {dimension}
 \def\jf    {J.\ Fuchs}
 \def\adma  {Adv.\wb Math.}
 \def\anma  {Ann.\wb Math.}
 \def\anop  {Ann.\wb Phys.}
 \def\cocm  {Com\-mun.\wb Con\-temp.\wb Math.}
 \def\coma  {Con\-temp.\wb Math.}
 \def\comp  {Com\-mun.\wb Math.\wb Phys.}
 \def\cpma  {Com\-pos.\wb Math.}
 \def\duke  {Duke\wB Math.\wb J.}
 \def\fiic  {Fields\wB Institute\wB Commun.}
 \def\gafa  {Geom.\wB and\wB Funct.\wb Anal.}
 \def\ijmp  {Int.\wb J.\wb Mod.\wb Phys.\ A}
 \def\jams  {J.\wb Amer.\wb Math.\wb Soc.}
 \def\joal  {J.\wB Al\-ge\-bra}
 \def\jomp  {J.\wb Math.\wb Phys.}
 \def\josp  {J.\wb Stat.\wb Phys.}
 \def\jpaa  {J.\wB Pure\wB Appl.\wb Alg.}
 \def\maan  {Math.\wb Annal.}
 \def\mpla  {Mod.\wb Phys.\wb Lett.\ A}
 \def\npbp  {Nucl.\wb Phys.\ B (Proc.\wb Suppl.)}
 \def\nupb  {Nucl.\wb Phys.\ B}
 \def\phlb  {Phys.\wb Lett.\ B}
 \def\phrd  {Phys.\wb Rev.\ D}
 \def\rims  {Publ.\wB RIMS}
 \def\rvmp  {Rev.\wb Math.\wb Phys.}
 \def\slnm  {Sprin\-ger\wB Lecture\wB Notes\wB in\wB Mathematics}
 \newcommand\Slnm[1] {{\rm[\slnm\ #1]}}
 \def\AMS    {{American Mathematical Society}}
 \def\IPC    {{International Press Company}}
 \def\PL     {{Plenum Press}}
 \def\PUP    {{Princeton University Press}}
 \def\SV     {{Sprin\-ger Ver\-lag}}
 \def\WS     {{World Scientific}}
 \def\Be     {{Berlin}}
 \def\PR     {{Providence}}
 \def\pR     {{Princeton}}
 \def\Si     {{Singapore}}
 \def\NY     {{New York}}

\end{document}